\documentclass[12pt,a4paper,makeidx]{book}
\usepackage[T1]{fontenc}
\usepackage[english,french]{babel}
\usepackage{amssymb,amsbsy,amsmath,amsfonts,amssymb,amsthm,amscd,stmaryrd}
\usepackage{graphicx,latexsym,amsthm}
\usepackage{mathrsfs}
\usepackage{pst-all}
\usepackage{color}
\usepackage{fancyhdr}

\usepackage{makeidx}
\makeindex

\numberwithin{equation}{chapter}
\newcommand\Retire[1]{\relax}

\newtheorem{definition}{{\bf D\'efinition}}[chapter]
\newtheorem{theorem}[definition]{{\bf Th\'eor\`eme}}
\newtheorem{corollary}[definition]{{\bf Corollaire}}
\newtheorem{proposition}[definition]{\noindent {\bf Proposition}}
\newtheorem{lemma}[definition]{\noindent {\bf Lemme}}

\newtheorem{observation}[definition]{\noindent {\bf Observation}}
\newtheorem{claim}{\noindent {\bf Affirmation}}
\newtheorem{question}{\noindent {\bf Question}}%[chapter]

\newtheorem{example}[definition]{\noindent {\bf Exemple.}}

\newtheorem{remark}[definition]{\noindent {\bf Remarque}}
\newtheorem{propertie}[definition]{\noindent {\bf Propri\'et\'e}}
\newtheorem{properties}[definition]{\noindent {\bf Propri\'et\'es}}
\newtheorem{remarks}[definition]{\noindent {\bf Remarques}}
\newtheorem{conjecture}{\noindent {\bf Conjecture}}
\newtheorem{problem}{\noindent {\bf Probl\`eme}}
\newtheorem{problems}{\noindent {\bf Probl\`emes}}
\newtheorem{consequence}{\noindent{\bf Cons\'equence}}
\newtheorem{fact}{\noindent{\bf Fait}}

\def\endproof{\hfill {\kern 6pt\penalty 500
\raise -0pt\hbox{\vrule \vbox to5pt {\hrule width 5pt
\vfill\hrule}\vrule}}}

\newcommand{\qedclaim}{\hfill $\diamond$}
\newenvironment{proofclaim}{\noindent{\bf \emph{Preuve.}}}{\qedclaim\medbreak}

\newcommand{\clearemptydoublepage}{\newpage\thispagestyle{empty}\cleardoublepage}

\newcommand{\NN}{{\mathbb N}}
\usepackage{setspace}

\begin{document}

\thispagestyle{empty}

%\vspace{150mm}
\begin{center}
\bigskip

\bigskip

\bigskip

\bigskip

\bigskip

\bigskip

\bigskip

\bigskip

\bigskip

\bigskip

\bigskip

\bigskip

\bigskip

\bigskip

\bigskip

\bigskip

\bigskip
\end{center}

\bigskip

\bigskip

\bigskip

\bigskip

\bigskip

\bigskip

\bigskip

\begin{center}
%\vspace{4cm}

        \LARGE{\bf{Sur l'\'enum\'eration de structures discr\`etes, une

             approche par la théorie des relations}}
    \end{center}

\smallskip

\begin{center}
 \large{\textbf{Djamila OUDRAR}}
 \end{center}

 \medskip
\begin{center} 
 D\'epartement de Recherche Op\'erationnelle, Facult\'e de math\'ematiques, 
 
 Universit\'e des Sciences et de Technologie Houari Boumediene

 \medskip

 28 septembre 2015
 \end{center}

\clearemptydoublepage

\newpage
\thispagestyle{empty}

\vspace{5cm}
\begin{center}
\begin{minipage}{12cm}
This work is our doctoral thesis defended on september 28, 2015 at the faculty of mathematics of the university of sciences and technology Houari Boumediene (USTHB) at Algiers, under the supervision of professors Moncef ABBAS of USTHB and Maurice POUZET of UCB (Lyon 1). The jury of defense consisted of professors Hacène BELBACHIR (USTHB) as president, Isma BOUCHEMAKH (USTHB), Youssef BOUDABBOUS (King Saud University), Bachir SADI (UMMTO, Algeria), Nicolas Marc THI\'ERY (University of Paris-sud, France) and Robert WOODROW (University of Calgary, Canada) as referees.
\end{minipage}
\end{center}

\vspace{3cm}

\begin{center}
\begin{minipage}{12cm}
Ceci est le travail réalisé dans le cadre de la thèse de doctorat soutenue le 28 septembre 2015 à la faculté de mathématiques de l'université des sciences et de la technologie Houari Boumediene (USTHB) à Alger, sous la direction des professeurs Moncef ABBAS (USTHB) et Maurice POUZET (Universit\'e Claude-Bernard (Lyon1), France et Universit\'e de Calgary, Canada). Le jury de soutenance se composait du professeur Hacène BELBACHIR (USTHB) en tant que pr\'esident et des professeurs Isma BOUCHEMAKH (USTHB), Youssef BOUDABBOUS (King Saud University), Bachir SADI (UMMTO, Algerie), Nicolas Marc THI\'ERY (Universit\'e Paris-sud, France) et Robert WOODROW (Universit\'e de Calgary, Canada) en tant qu'examinateurs.
\end{minipage}
\end{center}

\newpage
\thispagestyle{empty}
\thispagestyle{empty}

\noindent{\bf Abstract :}

\vspace{1mm}

Theory of relations is the framework of this thesis. It is about enumeration of finite structures. Let $\mathscr C$ be a class of finite combinatorial structures,
    the \emph{profile} of $\mathscr C$ is the function $\varphi_{\mathscr C}$ which count, for every $n$, the number of members of $\mathscr{C}$ defined on $n$ elements, isomorphic structures been
identified. The generating function for $\mathscr C$ is $\mathcal H_{\mathscr C}(X):=\sum_{n\geqq 0}\varphi_{\mathscr C}(n)X^n$. Many results about the behavior of the function $\varphi_{\mathscr C}$ have been obtained. Albert and Atkinson have shown that the generating series of the profile of some classes of permutations are algebraic. This result fits in the frame of theory of relations, we show how its conclusion extends to classes of ordered binary structures using the notions of theory of relations. This is the subject of the first part of this thesis.

   \vspace{3mm}

   The second part is concerned with the notion of minimality. An hereditary class of finite structures is minimal if it is infinite and every proper hereditary subclass is finite.
    We consider hereditary classes of finite structures, containing infinitely many indecomposable structures but every proper subclass contains only finitely many such structures, we call them ind-minimal classes. We show, in particular, that
   these classes are wqo ages and their number is the continuum,
   we give some examples of ind-minimal ages of graphs with their profiles and generating functions.

   \vspace{3mm}

   The last part is motivated by the surprising phenomenon of the \emph{jump} observed in the behavior of the profile of hereditary classes of finite structures. We start with the following notion due to Pouzet and Thi\'ery. A
   \emph{monomorphic decomposition} of a relational structure  $\mathcal R$ is a partition of its domain  $V(\mathcal R)$ into a family of sets $(V_x)_{x\in X}$ such that the restrictions of $\mathcal R$ to two finite subsets $A$ and $A'$ of $V(\mathcal R)$ are isomorphic provided that the traces  $A\cap V_x$ and  $A'\cap V_x$ have the same size for each $x\in X$.
   We show that the profile of an hereditary classe made of ordered structures which have finite monomorphic decomposition is a polynomial. We also show that the class $\mathscr D$ of ordered structures which do not have a finite monomorphic decomposition contains a finite subset $\mathfrak A$ such that every member of $\mathscr D$ embeds some member of $\mathfrak A$. In the case of ordered binary structures, the profile of every member of $\mathfrak A$ is at least exponential. We deduce that if the profile of a hereditary class of finite ordered binary structures is not bounded by a polynomial then it is at least exponential. This result generalizes the result obtained by Balogh, Bollob\'{a}s and  Morris (2006) for ordered graphs.

\vspace{8mm}

 \bigbreak

\newpage
\thispagestyle{empty}
\thispagestyle{empty}
\noindent{\bf R\'esum\'e :}

\vspace{1mm}

Ce travail s'inscrit dans le cadre de la th\'eorie des relations. Il porte sur l'\'enum\'eration des structures finies. Etant donn\'ee une classe $\mathscr C$ de structures combinatoires finies,

   le \emph{profil} de $\mathscr C$ est la fonction $\varphi_{\mathscr C}$ qui, \`a chaque entier $n$, associe le nombre de structures \`a $n$ \'el\'ements appartenant \`a $\mathscr C$, les structures isomorphes \'etant identifi\'ees. La s\'erie g\'en\'eratrice de $\mathscr C$ est la fonction $\mathcal H_{\mathscr C}(X):=\sum_{n\geqq 0}\varphi_{\mathscr C}(n)X^n$. De nombreux  r\'esultats sur le comportement de $\varphi_{\mathscr C}$ ont \'et\'e obtenus. Albert et Atkinson ont montr\'e que les s\'eries g\'en\'eratrices de certaines classes de permutations sont alg\'ebriques. Nous montrons comment ce r\'esultat se g\'en\'eralise aux classes de structures binaires ordonn\'ees, en mettant l'accent sur les notions de la th\'eorie des relations. Ceci fait l'objet de la premi\`ere partie de notre travail.

   \vspace{3mm}

   Nous \'etudions ensuite la notion de minimalit\'e. Une classe minimale \'etant une classe infinie de structures finies dont les sous-classes h\'er\'editaires propres sont finies.
    Nous consid\'erons des classes h\'er\'editaires de structures finies, contenant une infinit\'e de structures ind\'ecomposables dont les sous-classes propres n'en contiennent qu'un nombre fini que nous appelons classes ind-minimales. Nous montrons en particulier que %obtenons des r�sultats structurels, exemple % puis, nous donnons quelques exemples de classes minimales de graphes.
   ces classes sont des \^ages belordonn\'es et sont en nombre contin\^upotent,
    puis nous donnons des exemples d'\^ages ind-minimaux de graphes dont nous calculons les profils et les fonctions g\'en\'eratrices.

   \vspace{3mm}

   Enfin, nous nous int\'eressons au ph\'enom\`ene de \emph{saut} observ\'e dans le comportement des profils de classes h\'er\'editaires de structures finies. Nous partons de la notion de
   \emph{d\'ecomposition monomorphe} d'une structure relationnelle.
   Nous montrons que le profil d'une classe h\'er\'editaire form\'ee de structures ordonn\'ees qui ont une d\'ecomposition monomorphe finie est un polyn\^ome. Nous montrons \'egalement que la classe $\mathscr D$ des structures ordonn\'ees qui n'ont pas de d\'ecomposition monomorphe finie contient un ensemble fini $\mathfrak A$ tel que tout \'el\'ement  de  $\mathscr D$ abrite un \'el\'ement de $\mathfrak A$. Dans le cas binaire, tout membre de $\mathfrak A$ a un profil au moins exponentiel. Il en r\'esulte que, si  le profil d'une classe h\'er\'editaire de structures binaires ordonn\'ees  n'est pas born\'e par un polyn\^ome, il  est au moins exponentiel.  Un r\'esultat g\'en\'eralisant une classification  obtenue par   Balogh, Bollob\'{a}s et  Morris en 2006  pour les graphes ordonn\'es.

\vspace{8mm}

\noindent{\bf Mots-cl\'es} : Ensemble ordonn\'e, belordre, structure relationnelle, profil, ind\'ecomposabilit\'e, graphe, tournoi, permutation, d\'ecomposition monomorphe, encha\^{i}nabilit\'e.
 \bigbreak

\begin{center}
\hbox{\raisebox{0.4em}{\vrule depth 0pt height 0.4pt width 16cm}}

\bigbreak

\end{center}
\clearemptydoublepage

\thispagestyle{empty}
\addcontentsline{toc}{chapter}{Remerciements}
\begin{center}{\bf \it \LARGE Remerciements}
\end{center}

\vspace{0.3in}

{\it Ce travail n'aurait pas pu se faire sans le concours de beaucoup de personnes auxquelles je tiens \`a rendre un hommage mérité.

\vspace{3mm}

Tout d'abord, je dédie ce travail à mon mari Mourad, mes deux anges adorés Ghilès et Sara et toute ma famille, sans leurs sacrifices cette thèse n'aurait jamais vu le jour.

\vspace{3mm}

Je remercie vivement  mon directeur de th\`ese Moncef ABBAS pour sa gentillesse, son aide, ses conseils et ses encouragements durant ces longues ann\'ees de collaboration.

\vspace{3mm}

J'exprime toute ma reconnaissance et ma profonde gratitude \`a mon co-directeur de thèse Maurice POUZET pour m'avoir initiée à ce domaine, tellement passionnant qu'est la théorie des relations, pour son accueil chaleureux durant chacun de  mes déplacements à Lyon, son aide, sa patience, sa disponibilité et pour m'avoir fait profiter de ses larges connaissances. 

\vspace{3mm}

Un grand merci pour Hacène BELBACHIR pour m'avoir fait l'honneur de présider mon jury de thèse, pour sa gentillesse, ses conseils et son aide. Je lui suis reconnaissante de m'avoir intégrée dans son équipe de l'accord programme CMEP-TASSILI, ce qui m'a permis de bénéficier de nombreux stages à l'institut Camille Jordan de l'université Claude Bernard de Lyon, le personnel de cet institut qui a facilité mon travail, m'a acueillie et installée dans de bonnes conditions mérite tous mes remerciements.  Je remercie par la m\^eme occasion le professeur Jean Gabriel LUQUE le responsable de cet accord du côté français pour avoir facilité mes déplacements et pour avoir répondu présent à chaque fois que je l'ai sollicité. 

\vspace{3mm}

Mes vifs remerciements s'adressent aux autres membres de mon jury de thèse les professeurs Isma BOUCHEMAKH de l'USTHB, Youssef BOUDABBOUS de King Saud University (Arabie Saoudite), Bachir SADI de l'UMMTO, Nicolas Marc THI\'ERY de l'Université de Paris Sud (France) et Robert WOODROW de l'Université de Calgary (Canada) pour avoir accepté de donner de leur temps si pr\'ecieux pour examiner la présente thèse. Je suis particulièrement reconnaissante envers Nicolas Marc THI\'ERY et Robert WOODROW pour leurs remarques et critiques qu'ils m'ont adressées et qui ont contribuées à améliorer la rédaction de ce document.

\vspace{3mm}

Je n'oublierai pas de remercier mes collègues de la faculté de mathématiques et tous mes amis, certains ne sont plus à l'université, pour leur soutien, je ne pourrai pas tous les citer mais je suis sûre qu'ils se reconnaitront. Qu'ils m'excusent de ne citer que les deux personnes avec lesquelles j'ai eu à partager et celle avec laquelle je partage encore le même bureau: Nacéra, Yamina et Karima. Merci pour l'ambiance chaleureuse qui y a toujours régné, je remercie en particulier Nacéra pour ses critiques et pour sa disponibilité.

\vspace{3mm}

Je ne vais pas terminer cette page sans exprimer ma reconnaissance aux responsables de la scolarité et du département de Recherche Opérationnelle pour m'avoir déchargée de certaines de mes obligations pédagogiques pour me permettre de mener à terme mon travail.}

\clearemptydoublepage

\tableofcontents
\clearemptydoublepage

\addcontentsline{toc}{chapter}{Introduction}
\markboth{\slshape{Introduction g\'en\'erale}} {\slshape{Introduction g\'en\'erale}}
\chapter*{Introduction générale}
%\addcontentsline{toc}{chapter}{Introduction}
%\markboth{\slshape{Introduction g\'en\'erale}} {\slshape{Introduction g\'en\'erale}}

Ce travail porte sur une fonction énumératrice: le \emph{profil}. %Voici le contexte.
%
%porte sur l'\'enum\'eration des structures finies dans le cadre de la th\'eorie des relations. Voici le contexte.
Soit $\mathscr C$ une classe de structures combinatoires finies.
%(une structure relationnelle \'etant form\'ee d'un ensemble $V$ sur lequel on d\'efinit un ensemble de relations $(\rho_i)_{i\in I}$ o\`u, pour chaque $i\in I$, $\rho_i$ est d'arit\'e  $m_i$, ce qui signifie que $\rho_i \subseteq V^{m_i}$),
   Le \emph{profil} de $\mathscr C$ est la fonction $\varphi_{\mathscr C}$ qui, à chaque entier $n$, associe le nombre de structures \`a $n$ \'el\'ements appartenant \`a $\mathscr C$, celles-ci étant comptées à l'isomorphie près. La s\'erie g\'en\'eratrice de $\mathscr C$ est la fonction $\mathcal H_{\mathscr C}(X):=\sum_{n\geqq 0}\varphi_{\mathscr C}(n)X^n$.
Cette fonction a des propriétés remarquables lorsque la classe $\mathscr C$ est héréditaire,  c'est à dire qu'elle contient toutes les sous-structures de ses propres structures.

\vspace{2mm}

    Le cadre théorique des profils est la théorie des relations. Soit $n$ un entier non négatif, une \emph{relation $n$-aire} sur un ensemble $E$ est un sous-ensemble $\rho$ de $E^n$, l'ensemble des $n$-uples d'éléments de $E$; au besoin, $\rho$ est identifiée à sa fonction caractéristique. L'entier $n$ est l'\emph{arité} de $\rho$ et $E$ est sa \emph{base} ou son \emph{domaine}. Une \emph{structure relationnelle} est une paire $\mathcal R:=(E,(\rho_i)_{i\in I})$ formée de relations $\rho_i$ d'arité $n_i$ et de base $E$. La famille $\mu:=(n_i)_{i\in I}$ est la \emph{signature} de $\mathcal R$.
   Si $\mathcal R$ est une structure relationnelle, l'\emph{âge} de $\mathcal R$ est la collection $\mathcal A(\mathcal R)$ des restrictions finies de $\mathcal R$ considérées à l'isomorphie près.  Un âge étant une classe héréditaire, le profil de $\mathcal R$ est le profil de son âge.

\vspace{2mm}

   De tr\`es nombreux travaux portent sur le comportement des profils  de classes héréditaires $\mathscr C$ de structures relationnelles finies. %particulièrement lorsque $\mathscr C$ est une classe h\'er\'editaire.
   Beaucoup de r\'esultats ont \'et\'e obtenus, certains sont g\'en\'eraux,  (Pouzet 1971, 1978 et 2006) voir \cite{pouzet06}, d'autres portent sur des structures relationnelles particuli\`eres telles que les graphes \cite{B-B-M(07),B-B-S-S}, les tournois \cite{BBM07,Bou-Pouz} et les graphes ordonn\'es \cite{B-B-M(06)}, pour plus de détails voir la synthèse de Klazar 2010 \cite{klazar}.

\vspace{2mm}

   Les nombreux travaux sur l'\'enum\'eration des classes héréditaires de permutations qui ont prolif\'er\'e ces quinze derni\`eres ann\'ees, motiv\'es par la conjecture de Stanley-Wilf, r\'esolue par Marcus et Tard\"os (2004) \cite{Mar-Tar}, rentrent \'egalement dans le cadre de l'\'enum\'eration des structures relationnelles finies. En effet, comme sugg\'er\'e par P. J. Cameron (2002) \cite{cameron}, les permutations peuvent \^etre vues comme des structures relationnelles particuli\`eres: les bicha\^{i}nes (couple de deux ordres totaux d\'efinis sur un m\^eme ensemble).

\vspace{2mm}

   Ce qui retient particulièrement l'attention dans les résultats de ces travaux est que les profils de ces classes sont loin d'être arbitraires, leurs taux de croissances comportent des sauts. Typiquement, la croissance du profil est soit polynomiale soit plus grande que tout polynôme (\cite {pouzet.tr.1978} pour les \^ages, \cite{pouzet06} pour une synthèse). Pour certaines classes de structures, si le profil n'est pas polynomial, il est ou bien au moins exponentiel (exemple pour les tournois \cite{BBM07,Bou-Pouz}, graphes ordonnés et hypergraphes \cite{B-B-M(06),B-B-M06/2,klazar08} et les permutations \cite{K-K}) ou bien au moins comme la croissance de la fonction partition d'entiers (exemple pour les graphes \cite {B-B-S-S}). Pour plus de détails, voir la synthèse de Klazar \cite{klazar}.

   \bigskip%\vspace{4mm}

   %Plus précisément, il existe souvent une famille $\mathscr F$ de fonctions $f$ de $\mathbb N$ dans $\mathbb N$ et une autre fonction $F$ de $\mathbb N$ dans $\mathbb N$, avec $F(n)$ beaucoup plus grand que $f(n)$ pour tout $f$ dans $\mathscr F$, telles que si le profil $\varphi_{\mathscr C}$ d'une classe $\mathscr C$ est plus grand que $f(n)$ pour tout entier $n$ et tout $f$ de $\mathscr F$ alors il est également plus grand que $F(n)$ pour tout $n$ dans $\mathbb N$. On dit que le profil fait un \emph{saut} de $\mathscr F$ à $F$.\\

    Ce travail comporte trois parties.  Dans la première partie, nous explicitons le lien existant entre les bichaînes et les permutations et traduisons des r\'esultats sur les permutations et les classes de permutations, en utilisant les concepts et les outils de la th\'eorie des relations: abritement, belordre, belordre héréditaire, ind\'ecomposabilit\'e,.... En outre, nous illustrons le rôle du belordre et de l'indécomposabilité dans les résultats d'énumération.

   \medskip

 %\vspace{1mm}

  \noindent Nous nous intéressons, en particulier, aux résultats d'Albert et Atkinson (2005) \cite{A-A}. Ils ont montr\'e qu'une classe h\'er\'editaire de permutations, contenant un nombre fini de permutations simples, est alg\'ebrique et ont calcul\'e les fonctions g\'en\'eratrices de plusieurs classes h\'er\'editaires de permutations. Nous \'etendons leurs r\'esultats aux classes de structures relationnelles binaires ordonn\'ees (classes dont les structures sont form\'ees d'un ensemble fini muni d'un ordre total et d'un nombre fini fixe de relations binaires) en montrant (Corallaire \ref{cor:wqoalgebraic} et Th\'eor\`eme \ref{theo:algebraic}) que si une classe $\mathscr C$, form\'ee de telles structures, est h\'er\'editaire, contient un nombre fini de structures ind\'ecomposables et est close par sommes alors elle est h\'er\'editairement alg\'ebrique (en ce sens que la fonction g\'en\'eratrice de toute sous-classe h\'er\'editaire est alg\'ebrique). Nous illustrons ce r\'esultat en construisant une classe de structures binaires ordonn\'ees dont les ind\'ecomposables sont de tailles au plus deux (Section \ref{section:exemple}); classe entièrement caractérisée par des bornes (structures minimales interdites) dans le cas des bi-relations ordonnées (structures form\'ees d'un ordre total et d'une relation binaire d\'efinis sur un m\^eme ensemble). Cette classe est algébrique quelque soit le nombre, fini et fixe, de relations binaires qui composent ses structures. Ce r\'esultat est signifiant dans le cas des bi-relations ordonn\'ees  et croise les r\'ecents travaux sur les $d$-permutations  (suite de $d$ permutations dont la premi\`ere est l'identit\'e) d'Asinowski et Mansour (2010) \cite{Asin-Mans}.

  \medskip

 \noindent  Nous ne savons pas encore s'il est possible d'étendre ce résultat au cas où la classe $\mathscr C$  contient un nombre infini de structures indécomposables mais nous pensons y parvenir en imposant des restrictions \`a la classe $\mathcal I$ des ind\'ecomposables de $\mathscr C$. Nous pensons et conjecturons qu'il suffit que $\mathcal I$ soit h\'er\'editaire (dans la classe des ind\'ecomposables), h\'er\'editairement belordonn\'ee (en ce sens que la classe $\mathcal I.P$ des structures de $\mathcal I$ \'etiquet\'ees par un ensemble ordonn\'e $P$ est belordonn\'ee d\`es que $P$ est belordonn\'e) et h\'er\'editairement alg\'ebrique. Un exemple, vérifiant notre conjecture, est donn\'e en prenant pour $\mathcal I$ la classe des bicha\^{i}nes critiques de Schmerl et Trotter (1993) \cite{S-T}. Imposer \`a $\mathcal I$ d'\^etre belordonn\'ee ne suffit pas. Pour illustrer ce fait, nous avons \'etudi\'e quelques exemples d'\^ages de bicha\^{i}nes infinies %(l'\^age d'une structure $\mathcal R$ \'etant la classe de toutes les sous structures finies de $\mathcal R$)
   dont la classe des ind\'ecomposables est belordonn\'ee mais non h\'er\'editairement belordonn\'ee. L'\^age obtenu est rationnel et belordonn\'e, sa cl\^oture par sommes est alg\'ebrique mais non h\'er\'editairement alg\'ebrique et non belordonn\'ee.\\
 Nous étendons \'egalement un des r\'esultats d'Albert, Atkinson et Vatter (2010) \cite{A-A-V}, qui stipule qu'une classe h\'er\'editaire de permutations  qui est h\'er\'editairement rationnelle est belordonn\'ee, en affaiblissant l'hypoth\`ese. Nous montrons que si une classe h\'er\'editaire $\mathscr C$ de structures est h\'er\'editairement alg\'ebrique alors elle est belordonn\'ee.%\\
% Nous étendons \'egalement un des r\'esultats d'Albert, Atkinson et Vatter (2010) \cite{A-A-V}, qui stipule qu'une classe h\'er\'editaire de permutations  qui est h\'er\'editairement rationnelle est belordonn\'ee, en affaiblissant l'hypoth\`ese. Nous montrons que si une classe h\'er\'editaire $\mathscr C$ de structures est h\'er\'editairement alg\'ebrique alors elle est belordonn\'ee.

 \bigskip%\vspace{4mm}

 Dans la deuxième partie, nous nous intéressons %aux classes héréditaires de structures binaires contenant un nombre infini de structures indécomposables mais pas "trop". En d'autre termes, nous nous intéressons
 à la notion de \emph{minimalité} \cite{P-S}. Un ensemble ordonn\'e est dit \emph{minimal} s'il est infini et toute section initiale propre est finie. L'exemple le plus simple est celui de la cha\^{i}ne des entiers naturels. Si nous consid\'erons les structures finies, une classe héréditaire $\mathscr C$ de telles structures, %relationnelles finies, %de signature $\mu$,
 considérées à l'isomorphie près, est dite \emph{minimale} si elle est infinie et toute sous-classe h\'er\'editaire propre de $\mathscr C$ est finie. Ces classes ont \'et\'e, implicitement, caract\'eris\'ees par Fraïssé et peuvent-\^etre d\'ecrites: %Si $\mathscr C$ est une classe h\'er\'editaire minimale, elle est
 ce sont les \^ages de structures infinies encha\^{i}nables (une structure binaire $\mathcal R:=(E,(\rho_i)_{i\in I})$  est dite \emph{encha\^{i}nable} s'il existe un ordre linéaire, $\leq$, sur $E$ tel que, pour tout $i\in I$,
$x\rho_i y\Leftrightarrow x'\rho_i y'$ pour tous $x,y, x',y'$ vérifiant $x\leq y\Leftrightarrow x'\leq y'$). Si la signature $\mu$ est finie, ces \^ages sont en nombre fini, chacun de ces \^ages n'a qu'un nombre fini de bornes, un  r\'esultat important d\^u \`a Frasnay (1965), et toute classe héréditaire de l'ensemble $\Omega_{\mu}$ des structures relationnelles de signature $\mu$, contient une sous-classe minimale. En particulier, si $\mathscr C$ est une classe h\'er\'editaire minimale de graphes non dirig\'es sans boucles, $\mathscr C$ est l'\^age du graphe infini qui est complet ou  vide. Cet \^age est ainsi une cha\^{i}ne pour l'abritement et son profil est donc constant et vaut $1$.

\medskip

 \noindent Cette notion est particulièrement int\'eressante si $\mathscr C$ est form\'ee de structures binaires ind\'ecomposables et est h\'er\'editaire (dans la classe des ind\'ecomposables).
En effet, il existe une correspondance biunivoque entre les classes h\'er\'editaires  minimales de structures binaires ind\'ecomposables et les classes h\'er\'editaires de structures binaires finies qui contiennent une infinit\'e de structures ind\'ecomposables mais dont toute sous-classe h\'er\'editaire propre n'en contient qu'un nombre fini.  Par exemple, \`a la classe des chemins finis correspond l'\^age d'un chemin infini. Ces dernières classes sont dites \emph{ind-minimale}.

\medskip

 Nous consid\'erons des classes ind-minimales, c'est à dire des classes h\'er\'editaires de structures binaires finies dont les sous-ensembles d'ind\'ecomposables qu'elles contiennent sont minimaux (dans la classes des structures ind\'ecomposables), ou encore des classes contenant une infinit\'e de structures ind\'ecomposables dont les sous-classes propres n'en contiennent qu'un nombre fini. Contrairement aux classes minimales, le profil d'une telle classe n'est pas n\'ecessairement constant.  Nous obtenons certains résultats généraux que nous affinons aux cas des graphes. Nous montrons en particulier que les classes ind-minimales sont les \^ages de structures ind\'ecomposables infinies, qu'elles sont belordonn\'ees et que, contrairement aux classes minimales, elles %les classes ind-minimales
sont en nombre contin\^upotent.  Par exemple, en transformant les mots de Sturm en chemins dirig\'es, nous montrons qu'il y a un nombre continûpotent de  classes ind-minimales form\'ees de graphes dirig\'es. %Nous montrons qu'il y'a un nombre contin\^upotent d'\^ages ind-minimaux de structures binaires.
%Nous montrons en particulier que ces classes sont des \^ages belordonn\'es et sont en nombre contin\^upotent. %donnons quelques exemples de classes minimales de graphes.
%Une structure relationnelle binaire infinie $\mathcal{R}$ est \emph{minimale ind\'ecomposable} si $\mathcal{R}$ s'abrite dans toute structure infinie ind\'ecomposable $\mathcal{R}'$  qui s'abrite dans $\mathcal{R}$ ($\mathcal{R}'$ s'abrite dans $\mathcal{R}$ si $\mathcal{R}'$ est isomorphe à une restriction de $\mathcal{R}$).

\vspace{1mm}

Dans \cite{P-Z}, plusieurs exemples de graphes indécomposables sans clique infinie ou indépendant (stable) infini sont donnés. Nous explorons leurs \^ages  ainsi que ceux de huit autres graphes indécomposables ayant chacun une clique infinie et un indépendant infini.  Nous  déterminons les profils et les fonctions génératrices de ces \^ages et montrons qu'ils sont ind-minimaux.

\bigskip%\vspace{4mm}

 Dans la troisième partie, nous présentons une approche structurelle de résultats de sauts dans le comportement des profils des classes héréditaires de structures finies. % est consacrée au phénomène de \emph{saut} dans la croissance des profils des classes héréditaires de structures relationnelles finies.
 Nous partons du r\'esultat de Kaiser et Klazar \cite{K-K} qui stipule que le profil d'une classe h\'er\'editaire de permutations est soit born\'e sup\'erieurement par un polyn\^ome et dans ce cas c'est un polyn\^ome, soit born\'e inf\'erieurement par une exponentielle. Ce r\'esultat a \'et\'e g\'en\'eralis\'e au cas des graphes (non dirig\'es) ordonn\'es par Balogh, Bollob\`{a}s et Morris \cite{B-B-M06/2} qui montrent que la croissance des profils des classes h\'er\'editaires de graphes ordonn\'es pr\'esente le m\^eme saut: d'une croissance polynomiale \`a une croissance exponentielle.
 Dans cette partie, nous retrouvons, d'une nouvelle mani\`ere, le r\'esultat de Balogh et collaborateurs et nous l'\'etendons au cas des structures binaires ordonn\'ees de signature finie. Pour se faire, nous utilisons une technique bas\'ee sur le th\'eor\`eme de Ramsey, donn\'e en termes de structures invariantes et les notions de d\'ecomposition monomorphe due \`a Pouzet et Thi\'ery \cite{P-T-2013} et de presque multi-encha\^{i}nabilit\'e due \`a Pouzet (voir \cite{pouzet06}). La m\^eme technique a \'et\'e utilisée dans \cite{mont-pou} pour les bicha\^{i}nes et \cite{Bou-Pouz} pour les tournois. Nous donnons, ci-dessous, un aperçu de la d\'emarche suivie ainsi qu'un bref résumé des r\'esultats obtenus.

Le point de d\'epart est la notion de d\'ecomposition monomorphe due \`a Pouzet et Thi\'ery \cite{P-T-2013}:
une \emph{d\'ecomposition monomorphe} d'une structure relationnelle  $\mathcal R$ est  une partition de son domaine  $V(\mathcal R)$ en une famille de parties $(V_x)_{x\in X}$ telles que les  restrictions de $\mathcal R$ \`a deux parties  finies $A$ et  $A'$ de $V(\mathcal R)$ sont isomorphes pourvu  que les traces  $A\cap V_x$ et  $A'\cap V_x$ aient m\^eme cardinalit\'e pour tout  $x\in X$. En imposant des conditions aux blocs $(V_x)$ de la partition, nous donnons d'autres d\'ecompositions en parties monomorphes. Nous d\'ecrivons en tout quatre types de d\'ecompositions; trois d'entre elles sont valables pour des structures relationnelles quelconques et une pour les structures ordonn\'ees. Ces d\'ecompositions sont toutes \'equivalentes dans le sens où, si une structure relationnelle poss\`ede une d\'ecomposition d'un certain type ayant un nombre fini de blocs, alors elle poss\`ede une d\'ecomposition d'un autre type ayant un nombre fini de blocs. Nous \'etudions les relations existant entre ces diff\'erentes d\'ecompositions, puis nous d\'ecrivons une relation d'\'equivalence qui permet de retrouver la d\'ecomposition monomorphe et nous l'étudions sur des structures particulières (graphes non dirig\'es, tournois, graphes dirig\'es, bicha\^{i}nes et structures binaires ordonn\'ees).

\vspace{1mm}

\noindent Nous \'etudions les propri\'et\'es des structures ayant une d\'ecomposition monomorphe finie et des classes form\'ees par ces structures. Nous montrons (Proposition \ref{prop-wqo-ages} et Corollaire \ref{cor:wqo-monomorphy}) que, si une structure relationnelle de signature finie poss\`ede une d\'ecomposition monomorphe ayant un nombre fini de blocs, alors son \^age est h\'er\'editairement belordonn\'e et a un nombre fini de bornes.

\noindent Nous considérons des classes héréditaires de structures relationnelles ordonnées. Celles telles que toute structure qu'elles contiennent poss\`ede une d\'ecomposition monomorphe finie ont, d'apr\`es \cite{P-T-2005, P-T-2013}, des profils polynomialement bornés; en g\'en\'eral leurs profils sont des quasi-polyn\^omes. Nous montrons (Th\'eor\`eme \ref{thm: polynomial-interval}) que, dans le cas ordonn\'e, les profils de ces classes sont en fait des polynômes.

\vspace{1mm}

\noindent Nous nous int\'eressons ensuite aux structures relationnelles ordonn\'ees, de signature finie, qui ne poss\`edent pas de d\'ecomposition monomorphe finie. Nous montrons, en utilisant le th\'eor\`eme de Ramsey, que toute structure ordonn\'ee qui ne poss\`ede pas de d\'ecomposition monomorphe finie abrite une structure presque multi-encha\^{i}nable dont la d\'ecomposition monomorphe n'est pas finie. Soit $\mathscr D_{\mu}$ la classe form\'ee des structures ordonn\'ees de signature $\mu$ finie qui ne poss\`edent pas de d\'ecomposition monomorphe finie; alors $\mathscr D_{\mu}$ contient un ensemble fini $\mathfrak A$ form\'e de structures presque multi-encha\^{i}nables incomparables tel que toute structure de $\mathscr D_{\mu}$ abrite un \'el\'ement de $\mathfrak A$.

Nous montrons \'egalement (Lemme \ref{lem:reduction}) que, si une classe h\'er\'editaire $\mathscr C$  form\'ee de structures relationnelles ordonn\'ees finies de signature $\mu$ finie contient, pour tout entier $\ell$, une structure finie qui ne poss\`ede pas de d\'ecomposition monomorphe ayant au plus $\ell$ blocs, alors elle contient une classe h\'er\'editaire $\mathscr A$, ayant la m\^eme propri\'et\'e, qui est minimale pour l'inclusion. La classe $\mathscr A$ est l'\^age d'une structure presque multi-encha\^{i}nable appartenant \`a l'ensemble $\mathfrak A$. Ceci nous conduit au th\'eor\`eme de dichotomie pour les classes de structures ordonn\'ees (Th\'eor\`eme \ref{theo:dichotomie-ordon}): \'etant donn\'ee une classe h\'er\'editaire $\mathscr C$ de structures relationnelles ordonn\'ees finies, de signature $\mu$ finie, ou bien il existe une borne sur le nombre de composantes monomorphes de chaque membre de $\mathscr C$ et le profil de $\mathscr C$ est un polyn\^ome, ou bien $\mathscr C$ contient l'\^age d'une structure presque multi-encha\^{i}nable appartenant \`a l'ensemble fini $\mathfrak A$.

\vspace{1mm}

Nous pensons et conjecturons que, dans ce cas ordonn\'e, les structures de l'ensemble $\mathfrak A$ ont des profils exponentiels et nous montrons que ce r\'esultat est vrai dans le cas particulier o\`u $\mathscr D_{\mu}$ est form\'e de structures binaires ordonn\'ees. Ceci entraine le r\'esultat principal de cette partie, un r\'esultat de dichotomie dans le cas particulier des structures binaires ordonn\'ees, donn\'e par le  Th\'eor\`eme \ref{theo:dichotomie}: le profil d'une structure binaire ordonn\'ee de type $k$ (structure ayant $k$ relations dont la premi\`ere est un ordre total) est soit polynomial soit born\'e inf\'erieurement par une exponentielle. Ce r\'esultat s'\'etend aux classes h\'er\'editaires form\'ees de telles structures. Il en r\'esulte (Th\'eor\`eme \ref{theo:dichotomie-classe}) que si  le profil d'une classe h\'er\'editaire de structures binaires ordonn\'ees  n'est pas born\'e par un polyn\^ome, il  est au moins exponentiel.  %R\'esultat faisant partie d'une classification  obtenue par   Balogh, Bollob\'{a}s et  Morris en 2006 \cite{B-B-M06/2} pour les graphes ordonn\'es.

%\medskip%
\vspace{1mm}

\noindent La preuve de ce r\'esultat passe par la description, très précise, des structures de l'ensemble fini $\mathfrak A$ dans le cas où les structures sont des graphes dirig\'es ordonn\'es. Cet ensemble est form\'e de mille deux cent quarante six graphes dirig\'es ordonn\'es tel que tout graphe dirig\'e ordonn\'e qui n'a pas de d\'ecomposition monomorphe fini abrite un \'el\'ement de $\mathfrak A$.
En outre,  pour chaque $\mathcal R\in \mathfrak A$, le profil de l'\^age $\mathcal A(\mathcal R)$  de $\mathcal R$  est au moins exponentiel (le profil est born\'e inf\'erieurement par la fonction de Fibonacci).

\vspace{1mm}

\bigskip%\vspace{3mm}

 Les résultats de ces trois parties ont été  exposés à différents  colloques internationaux, DIMACOS11 (2011-Mohammadia-Maroc), ISOR11 (2011-USTHB) et JGA  (2011-Lyon) pour les résultats de la première partie, RAMA08 (2012-Tipaza) pour ceux de la deuxième et ICGT (2014-Grenoble) pour ceux de la troisième. Les résultats de la première partie  ont donné lieu à un article soumis et accepté pour publication au "journal of multiple-valued logic and soft computing (MVLSC)", ceux des deux autres parties feront l'objet de publications ultérieures.\\
%Les résultats de cette partie constituent le chapitre 2.\\

%Les résultats de cette partie ont été exposés au colloque international RAMA08 (2012-Tipaza).\\

%Les résultat de cette partie ont été exposés au colloque international ICGT 2014 (Grenoble) et ferons l'objet de publications ultérieures.\\

 La thèse est organisée en trois parties, huit chapitres et deux annexes. La première partie comprend les chapitres 2, 3 et 4, la deuxième comprend les chapitres 5 et 6 et la troisième les chapitres 7 et 8. %sont consacrés aux résultats décrits dans les parties 1, 2 et 3 ci-dessus respectivement.
 Le chapitre 1 est consacré aux définitions et notions de base nécessaires à la compréhension du reste du document. L'annexe A complète un r\'esultat donn\'e dans la section \ref{subsection:classe S1}. En effet, dans cette section, nous donnons la s\'erie g\'en\'eratrice de la classe $\mathscr S^{re}_1$ des structures binaires ordonn\'ees s\'eparables r\'eflexives de type $1$, il se trouve que cette s\'erie est la s\'erie g\'en\'eratrice de la classe des $3-$permutations s\'eparables et celle des partitions guillotine en dimension $4$, la correspondance bijective entre ces deux classes a \'et\'e \'etabli dans \cite{Asin-Mans}. Nous construisons, dans l'annexe A,  une bijection entre les structures binaires ordonn\'ees s\'eparables r\'eflexives de type $1$ et les arbres binaires \'etiquet\'es par l'ensemble $\{1,2,3,4\}$ sachant que ces derniers sont en correspondance bijective avec les partitions guillotine en dimension $4$ \cite{Ack-Ba-Pin-Rom}. Ceci entraine que l'ensemble $\mathscr S^{re}_1$ et l'ensemble des partitions guillotine en dimension $4$ sont isomorphes. Cela entraine \'egalement que $\mathscr S^{re}_1$ est isomorphe \`a l'ensemble des $3-$permutations s\'eparables.

 \vspace{1mm}

 \noindent Dans l'annexe B, nous faisons quelques rappels sur des classes de permutations particuli\`eres, \`a savoir les classes grille-g\'eom\'etriques dont les propri\'et\'es ont permis de montrer que la cl\^oture par sommes de la classe des bicha\^{i}nes critiques est h\'er\'editairement alg\'ebrique. Ainsi, cette classe de bicha\^{i}nes v\'erifie notre conjecture pour les classes héréditaires de structures binaires ordonnées (Conjecture \ref{conjec} en page \pageref{conjec}).

 \vspace{2mm}

  En outre, pour faciliter la lecture du document et permettre de retrouver les différentes notations et définitions de concepts introduits, une table des notations est incluse en page \pageref{notation} et un index à la fin du document.  

\clearemptydoublepage

\chapter{Notions et outils de base}\label{chap:generalite}

Dans ce chapitre, nous introduisons les notions et les objets de base nécessaires pour la compréhension du document. Il n'a pas la pr\'etention d'\^etre exhaustif. Les notions plus sp\'ecifiques aux diff\'erents chapitres sont introduites le moment opportun.

\noindent Bien que la terminologie utilis\'ee soit, essentiellement, celle de Fra\"{i}ss\'e \cite{fraisse}, les r\'ef\'erences \cite{Flajolet, pouzet.81, pouzet06, pouzet12,  stanley, stanley2, trotter, Van-Wil} ont \'egalement \'et\'e utilis\'ees.

Nous utilisons les notations standards de la th\'eorie des ensembles. Si $E$ est un ensemble, $\vert E\vert$ d\'esigne sa cardinalit\'e. Si $E$ est fini et $\vert E\vert=p$ nous disons que $E$ est un $p$-ensemble\index{ensemble!$p$-ensemble}. Si $E$ est un $p$-ensemble, $F\subseteq E$ et $\vert F\vert=n$,  nous disons que $F$ est une $n$-partie de $E$. Nous notons  $\mathscr P(E)$ l'ensemble des parties\index{ensemble!des parties} de $E,\;[E]^n$ l'ensemble des   $n$-parties de l'ensemble $E$ et $E^n$ l'ensemble des $n$-uples d'\'el\'ements de $E$. Etant donn\'e un entier $n\in\mathbb N^{\star}$, nous notons par $[0,n[$ l'ensemble $\{0,\dots,n-1\}$ et par $[0,n]$ l'ensemble $\{0,\dots,n\}$.

\section{Relations binaires, ordres, graphes et ensembles ordonn\'es}\label{subsec:rela.bin}
    \subsection{D\'efinitions et notations}
 Soit $E$ un ensemble. Une \emph{relation binaire sur $E$},\index{relation!binaire} est un ensemble $\rho$ de couples de $E$. L'ensemble $E$ est la \emph{base} ou le \emph{domaine} de $\rho$. Pour tout ensemble $E$, l'ensemble $E^2$ est une relation binaire, nous convenons que l'ensemble vide $\varnothing$ est \'egalement une relation binaire et si $\rho$ est l'ensemble $\Delta_E$ des couples $(x,x)$ pour $x$ dans $E$ alors $\rho$ est la relation d'\'egalit\'e sur $E$.

 \vspace{1mm}

 Etant donn\'ee une relation binaire $\rho$ de base $E$ et deux \'el\'ements $x$ et $y$ de $E$,
si $(x,y)\in\rho$ nous disons que $x$ et $y$ sont en relation et notons indiff\'eremment $x\rho y$ ou $(x,y)\in \rho$, nous notons $x(non\rho)y$, $x\neg\rho y$ ou $(x,y)\notin\rho$ le fait contraire.
Si $\rho$ est une relation binaire sur $E$, sa relation \emph{compl\'ementaire}\index{relation!compl\'ementaire} %\footnote{C'est la n\'egation d\'efinie plus haut. Pour les relations binaires on a pr\'ef\'er\'e l'appellation classique. }
dans $E$ est la relation $\rho^c:=\{(x,y)\in E^2:\; (x,y)\notin \rho\}=E^2\setminus \rho$. Sa relation \emph{inverse}\index{relation!inverse} ou \emph{duale}\index{relation!duale} est la relation binaire $\rho^{-1}:=\{(x,y)\in E^2:\;(y,x)\in \rho\}$. Une relation qui est \'egale à sa relation duale est dite \emph{autoduale}\index{relation!autoduale}.

\vspace{1mm}

Si $\rho$ est une relation binaire sur $E$, une partie de $\rho$ est une \emph{sous-relation} de $\rho$ et si $A$ est une partie de $E$ la \emph{restriction}\index{restriction} de $\rho$ \`a $A$, dite \'egalement \emph{relation induite par} $\rho$ sur $A$, est la relation $\rho_{\restriction_A}:=\rho \cap A^2.$

\bigskip

Une relation binaire $\rho$ sur un ensemble $E$ est:
\begin{itemize}
\item \emph{r\'eflexive}\index{relation!r\'eflexive} si $x\rho x$ pour tout $x\in E$;
\item \emph{irr\'eflexive}\index{relation!irr\'eflexive} si $x\neg \rho x$ pour tout $x\in E$;
\item \emph{sym\'etrique}\index{relation!sym\'etrique} si $x\rho y$ implique $y\rho x$ pour tous $x, y\in E$;
\item \emph{antisym\'etrique}\index{relation!antisym\'etrique} si $x\rho y$ et $y\rho x$ impliquent $x=y$ pour tous $x, y\in E$;
\item \emph{totale}\index{relation!totale} si $x\neq y$ implique $x\rho y$ ou $y\rho x$ pour tous $x, y\in E$;
%\item \emph{compl\`ete}\index{relation!compl\`ete} si $x\rho y$ et $y\rho x$ pour tous $x, y\in E$, $x\neq y$;
\item \emph{transitive}\index{relation!transitive} si $x\rho y$ et $y\rho z$ impliquent $x\rho z$ pour tous $x, y, z\in E$.
\end{itemize}

Nous appelerons relation \emph{compl\`ete}\index{relation!compl\`ete} toute relation $\rho$ de base $E$ qui est sym\'etrique et totale, c'est \`a dire une relation pour laquelle  $x\rho y$ et $y\rho x$ pour tous $x, y\in E$, $x\neq y$.

\subsubsection{Graphe}
Un \emph{graphe\index{graphe!dirig\'e} dirig\'e}\footnote{\emph{graphe simple orient\'e} dans la terminologie de Claude Berge \cite{berge}.} est un couple $G:=(V,\rho)$ form\'e d'un ensemble $V$ et d'une relation binaire $\rho$ sur $V$. L'ensemble $V$ est souvent not\'ee $V(G)$, ses \'el\'ements sont appel\'es \emph{sommets}\index{graphe!sommet d'un -} de $G$ et sa cardinalit\'e $\vert V\vert$ est appel\'ee \emph{ordre}\index{graphe!ordre d'un -} de $G$. L'ensemble $\rho$, qui est une partie de $V^2$, est not\'e $E(G)$ et ses \'el\'ements sont appel\'es \emph{arcs}\index{graphe!arc d'un -} de $G$; les couples $(x,x)$ de $\rho$ sont les \emph{boucles} de $G$. L'inverse ou dual du graphe\index{graphe!inverse} dirig\'e\index{graphe!dirig\'e} $G$ est le graphe dirig\'e $G^{-1}:=(V,\rho^{-1})$ et son compl\'ementaire\index{graphe!compl\'ementaire} est $G^c:=(V,\rho^c).$ Si $A\subseteq V$, alors la restriction de $G$ \`a $A$, donn\'ee par $G_{\restriction_A}:=(A,\rho_{\restriction_A})$ est le \emph{sous-graphe\index{graphe!sous-graphe induit} dirig\'e induit} par $G$ sur $A$.

\vspace{1mm}

Un \emph{chemin}\index{chemin} d'un graphe dirig\'e $G$ %ce terme s'utilise dans un graphe orient\'e\index{graphe!orient\'e}, pour une cha\^{i}ne o\`u tous les arcs sont orient\'es dans le m\^eme sens.
 est une suite $P=(v_1,v_2,\ldots,v_k)$ de sommets distincts de $V(G)$ tels que, pour tout $i\in\{1,\ldots,k-1\}$, $(v_i, v_{i+1})\in E(G)$. Si $(v_i, v_j)\not\in E(G)$ pour tous $i,j$ tels que $j\neq i+1$, ceci d\'efinie sur l'ensemble $\{v_1,v_2,\ldots,v_k\}$ une relation appel\'ee la relation de \emph{cons\'ecutivit\'e}.\index{relation!de cons\'ecutivit\'e}\\

Un \emph{tournoi}\index{tournoi} est un graphe dirig\'e $G:=(V,\rho)$ pour lequel la relation $\rho$ est irr\'eflexive, antisym\'etrique et totale. Si $\rho$ est, de plus, transitive, $G$ est un \emph{tournoi transitif} appel\'e aussi \emph{tournoi acyclique}.\index{tournoi!transitif}\index{tournoi!acyclique}\\

Un \emph{graphe}\footnote{\emph{graphe simple non orient\'e sans boucle} dans la terminologie de C. Berge.} est un couple $G:=(V, E)$ dans lequel  $E$ est une partie de $[V]^2$ dont les \'el\'ements sont appel\'es \emph{ar\^etes}\index{graphe!ar\^ete d'un -} de $G$; les sommets $x$ et $y$ de $G$ sont dits \emph{adjacents} %, fait not\'e $x\sim y$,
si $\{x,y\}$ est une ar\^ete de $G$.

\noindent Nous pouvons identifier l'ensemble des paires $E$ d'\'el\'ements de $V$ \`a une relation binaire $\rho$ irr\'eflexive et sym\'etrique. %et inversement, c'est \`a dire
Ainsi $\{x,y\}\in E$ si et seulement si $(x,y)\in\rho$ et $(y,x)\in\rho$. Dans la suite, nous consid\'erons que toute relation irr\'eflexive et sym\'etrique d\'efinie sur un ensemble $V$ est un graphe et, inversement, tout graphe dont l'ensemble des sommets est $V$ est une relation irr\'eflexive et sym\'etrique sur $V$.

\noindent Si $G:=(V, E)$ est un graphe, son compl\'ementaire est le graphe $G^c:=(V, [V]^2\setminus E)$  not\'e aussi $\overline{G}$.
%\smallskip

\noindent Un graphe $G$ est \emph{complet}\index{graphe!complet}  si deux sommets distincts $x,y$ de $G$ sont toujours adjacents. Nous notons $K_n$ le graphe complet \`a $n$ sommets. Un graphe $G$ est un \emph{ind\'ependant}\index{graphe!ind\'ependant} %ou \emph{vide}\index{graphe!vide}
si deux sommets $x, y$ de $G$ sont toujours non adjacents. Nous notons $I_n$ l'ind\'ependant \`a $n$ sommets.

\noindent Soit un graphe $G$; une partie $A$ de $V(G)$ est une \emph{clique}\index{clique} si le sous-graphe induit sur $A$ est complet. $A$ est un \emph{stable}\index{stable} si deux sommets $x, y$ de $A$ sont toujours non adjacents, ou encore si le compl\'ementaire du sous-graphe induit sur $A$ est une clique\index{clique}.\\

Un \emph{chemin}\index{chemin}\footnote{\emph{chaîne \'el\'ementaire} dans la terminologie de C. Berge. %ce terme s'utilise dans un graphe orient\'e\index{graphe!orient\'e}, pour une cha\^{i}ne o\`u tous les arcs sont orient\'es dans le m\^eme sens.
} d'un graphe $G$ est une suite $P=(v_1,v_2,\ldots,v_k)$ de sommets distincts de $G$ tels que, pour tout $i\in\{1,\ldots,k-1\}$, $v_i v_{i+1}$ est une ar\^ete de $G$. Les sommets $v_1$  et $v_k$ sont les extr\'emit\'es de $P$. La \emph{longueur}\index{chemin!longueur d'un -} d'un chemin est le nombre de ses ar\^etes. Ainsi la longueur de $P=(v_1,v_2,\ldots,v_k)$ est $k-1$. Toute ar\^ete reliant deux sommets non cons\'ecutifs de $P$ est appel\'ee corde\index{corde}. Nous notons par $P_n$ un chemin sans corde \`a $n$ sommets.

    \subsubsection{Pr\'eordre}
Une relation binaire $\rho$ sur un ensemble $E$ est un \emph{pr\'eordre}\index{relation!de pr\'eordre} si elle est r\'eflexive et transitive. L'ensemble $E$ muni d'un pr\'eordre est un \emph{ensemble pr\'eordonn\'e},\index{ensemble!pr\'eordonn\'e} souvent  appel\'e  \emph{qoset}, abr\'eviation de l'appellation en anglais ``\emph{quasi-ordered set}''.  Dans ce cas, au lieu de $(x,y)\in \rho$ nous disons que $x$ est plus petit que $y$ et nous notons $x\leq y (mod\;\rho)$, ou que $y$ est plus grand que $x$, not\'e $y\geq x (mod\;\rho)$, nous notons $x\nleq y (mod\;\rho)$ le fait contraire (ou $y\ngeq x (mod\;\rho)$). Deux \'el\'ements distincts $x$ et $y$ de $E$ sont dits \emph{comparables}\index{element@\'el\'ements comparables} par $\rho$ si $x\leq y (mod\;\rho)$ ou $y\leq x (mod\;\rho)$, dans le cas contraire ils sont dits \emph{incomparables},\index{element@\'el\'ements incomparables} fait not\'e $x\parallel_{\rho} y$. Nous disons que $x$ est strictement plus petit que $y$ et notons $x<y$ si $x\leq y (mod\;\rho)$ et $y\nleq x (mod\;\rho)$. %Si $x$ et $y$ sont deux \'el\'ements distincts de $E$,
Nous disons que $x$ est \emph{couvert} par $y$ lorsque $x<y (mod~\rho)$ et il n'existe aucun \'el\'ement $z$ de $E$ tel que $x<z(mod~\rho)$ et $z<y (mod~\rho)$. Si $\rho$ est un pr\'eordre sur $E$, la relation duale $\rho^{-1}$ et la restriction $\rho_A$ \`a n'importe quelle partie $A$ de $E$ sont \'egalement des pr\'eordres sur $E$.  Deux vari\'et\'es particuli\`eres de pr\'eordres sont les \'equivalences et les ordres.\\

 Une \emph{relation d'\'equivalence} $\rho$\index{relation!d'\'equivalence} sur $E$ est un pr\'eordre sym\'etrique. Au lieu de $(x,y)\in \rho$ nous disons que $x$ est \'equivalent \`a $y$, not\'e $x\equiv y (mod\;\rho)$, nous notons $x\not \equiv  y(mod\;\rho)$ le fait contraire, %(on utilise souvent le symbole $\equiv$ pour d\'esigner une relation d'\'equivalence).
Si $\rho$ est une relation d'\'equivalence sur $V$ alors pour chaque $x\in E$, l'ensemble
$\overline{x}:=\{y\in E,\; x\equiv y (mod\;\rho)\}$ est la classe\index{classe!d'\'equivalence} (d'\'equivalence) de $x$. L'ensemble des classes d'\'equivalence, not\'e $E/{\rho}$, est appel\'e le \emph{quotient} de $E$ par la relation d'\'equivalence $\rho$. L'application $\varphi$ qui \`a chaque $x\in E$, associe $\overline{x}$ est la \emph{surjection canonique} de $E$ sur $E/{\rho}$. Si $X\in E/{\rho}$,  tout \'el\'ement de $\varphi^{-1}(X)$ est appel\'e un \emph{repr\'esentant de $X$ dans $E$}.

 \noindent Si $\rho$ est une relation d'\'equivalence sur $E$ alors $E/{\rho}$ est une partition de $E$. Nous rappelons qu'une partition d'un ensemble non vide $E$ est une famille $\mathscr P$ de sous-ensembles non vides $P\subseteq E$, deux \`a deux disjoints de $E$ telle que $\underset{P\in\mathscr P}\cup P=E$.

 \noindent R\'eciproquement, \`a toute partition $\mathcal A:=(A_j)_{j\in J}$ de $E$ correspond une unique relation d'\'equivalence $\rho$ sur $E$ telle que $\mathcal A=E/{\rho}$.
 \bigskip

 Un \emph{ordre partiel}\index{relation!d'ordre partiel}, ou simplement un \emph{ordre}\index{relation!d'ordre} $\rho$ sur $E$, est un pr\'eordre antisym\'etrique; l'ensemble $E$ muni de cet ordre $\rho$ est un \emph{ensemble ordonn\'e}\index{ensemble!ordonn\'e}, souvent appel\'e \emph{poset} en anglais, une abr\'eviation de %l'appellation en anglais
 ``\emph{partially ordered set}''. Lorsqu'il n'y a pas de confusion, nous ne ferons pas de diff\'erence entre ordre et ensemble ordonn\'e en utilisant la lettre, disons $P$, pour d\'esigner les deux.

 Si $\rho$ est un ordre sur $E$, il y va de m\^eme de sa relation inverse\index{relation!inverse} $\rho^{-1}$ et de sa restriction $\rho_A$ \`a n'importe quelle partie $A$ de $E.$  Nous d\'esignons souvent une relation d'ordre\index{relation!d'ordre} par le symbole "$\leq_{\rho}$" ou simplement par "$\leq$" s'il n'y a pas risque de confusion.%Dans ce cas la relation d'\'equivalence g\'en\'er\'ee par $\rho$ est la relation d'\'egalit\'e. Par cons\'equent, l'in\'egalit\'e strict $x<y(mod\;\rho)$ signifie que $x\leq y \; \text{et}\; x\neq y(mod\;\rho)$. Une autre cons\'equence: \'etant donn\'es $x, y\in E$, on a quatre cas possibles: soit $x=y\; \text{ou bien}\; x<y\;\text{ou bien}\; x>y\;\text{ou bien}\; x|y (mod\;\rho)$.\\
\vspace{1mm}

\noindent \textbf{Exemples d'ensembles ordonn\'es.}\index{ensemble!ordonn\'e} - Consid\'erons un ensemble $A$, tout sous-ensemble de $\mathscr P(A)$ est partiellement ordonn\'e\index{ensemble!partiellement ordonn\'e} par inclusion.\\$-$ La relation $\rho$ d\'efinie sur $\mathbb N$ (ensemble des entiers $\geq 0$ ou entiers naturels) par: $x\rho y$ si et seulement si $x$ divise $y$ pour tous $x,y\in \mathbb N$ est un ordre partiel.\index{ordre!partiel}\\

\noindent Une relation d'ordre\index{relation!d'ordre} $\leq$ sur $E$ est un ordre \emph{total} ou \emph{lin\'eaire},\index{ordre!lin\'eaire}  %appel\'e aussi \emph{cha\^{i}ne},\index{chaine@cha\^{i}ne!d'un ordre}
si deux \'el\'ements de $E$ sont toujours comparables; c'est \`a dire que $\leq$ v\'erifie la condition:
 $$ x\leq y\; \text{ou } y\leq x\;\text{ pour tous } x, y \in E.$$
Dans ce cas, les \'el\'ements de $E$ sont usuellement repr\'esent\'es sur une ligne (horizontale ou verticale), d'o\`u le mot ``lin\'eaire''.
L'ensemble $E$ muni de cet ordre lin\'eaire est appel\'e \emph{cha\^{i}ne}.\index{chaine@cha\^{i}ne!d'un ordre}

\textbf{Exemple.} La cha\^{i}ne naturelle $\omega$ de base $\mathbb N$ d\'efinie avec les comparaisons usuelles (l'ordre naturel sur $\mathbb N$) est un exemple fondamental de cha\^{i}nes infinies.\index{chaine@cha\^{i}ne!infinie}\\

\noindent La relation d'\'egalit\'e sur $E$ est une relation d'ordre\index{relation!d'ordre} sur $E$. L'ensemble $E$ muni de cet ordre est une \emph{anticha\^{i}ne}.\index{antichaine@anticha\^{i}ne}\\ %si deux \'el\'ements distincts de $E$ sont toujours incomparables.\\
Tout ordre sur un singleton %ensemble \`a un \'el\'ement
est \`a la fois une cha\^{i}ne et une anticha\^{i}ne, il est trivial.
\bigskip

\noindent Soit $\leq$ une relation d'ordre sur un ensemble $E$. La relation binaire $<$ sur $E$ d\'efinie par $x<y$ (ou de mani\`ere \'equivalente $y>x$) si et seulement si $x\leq y$ mais $x\ngeq y$  qui se lit $x$ \emph{strictement inf\'erieur} \`a $y$ ou $y$ \emph{strictement sup\'erieur} \`a $x$, est un \emph{ordre strict}\index{ordre!strict} (l'ordre strict associ\'e \`a $\leq$). La relation $<$ n'est pas une relation d'ordre, c'est une relation irr\'eflexive antisym\'etrique et transitive. Plus clairement, si $\rho$ est une relation d'ordre sur $E$, la relation $\delta:=\rho\setminus\Delta_E$, o\`u $\Delta_E:=\{(x,x),\;x\in E\}$, est l'ordre strict associ\'e \`a $\rho$. R\'eciproquement, si $\delta$ est un ordre strict sur un ensemble $E$, alors elle est antisym\'etrique et la relation $\rho:=\delta\cup\Delta_E$ est un ordre sur $E$.\\ Si $\leq$ est un ordre total sur $E$ alors l'ordre strict $<$ sur $E$ est un tournoi transitif. R\'eciproquement, \`a tout tournoi transitif correspond un ordre total.\index{tournoi!transitif}\\

Etant donn\'e un ensemble ordonn\'e $(E,\leq)$, nous pouvons lui associer une repr\'esentation sur le plan o\`u les \'el\'ements de l'ensemble fini $E$ sont repr\'esent\'es par des points tels que si $x$ est couvert par $y$ dans $E$ alors le point repr\'esentant $x$ se trouve en ``bas'' du point repr\'esentant $y$, ces deux points \'etant reli\'es par un trait ou un segment. Les boucles ne sont pas repr\'esent\'ees ainsi que les relations pouvant se d\'eduire par transitivit\'e. Cette repr\'esentation est appel\'ee \emph{diagramme de Hasse}.\index{diagramme de Hasse}

  \subsection{Homomorphisme et isomorphisme}
  Soient $\rho$ et $\rho'$ deux relations binaires de base $E$ et $E'$ respectivement. Une application $f$ de $E$ dans $E'$ est un \emph{homomorphisme}\index{homomorphisme} de $(E,\rho)$ dans $(E',\rho')$ si pour tous $x,y\in E$: $$(x,y)\in\rho\Rightarrow (f(x),f(y))\in\rho'.$$
  Si $E=E'$ et $\rho=\rho'$ l'homomorphisme $f$ est un \emph{endomorphisme}.\index{endomorphisme}
 L'application  identit\'e $1_E:E\rightarrow E$ est un endomorphisme.

 Dans le cas où $(E,\rho)$ et $(E',\rho')$ sont deux ensembles ordonn\'es, un homomorphisme de $(E,\rho)$ sur $(E',\rho')$ est appel\'e \emph{application croissante}.\index{application croissante}
 %Soient $P$ et $Q$ deux ensembles ordonn\'es. Une \emph{application croissante de $P$ dans $Q$} est toute application $f$ de $P$ dans $Q$ telle que pour tous $x,y\in P$:
 %$$x\leq y (mod\; P)\Rightarrow f(x)\leq f(y) (mod\; Q).$$
%Pour tout ensemble ordonn\'e $P$, l'application  identit\'e $1_P:P\rightarrow P$ est une application croissante.\\
 Soient $\rho$, $\rho'$ et $\rho''$ des relations binaires sur $E$, $E'$ et $E''$ respectivement. Si $f:E\rightarrow E'$ et $g:E'\rightarrow E''$ sont deux homomorphismes alors leur compos\'e $g\circ f:E\rightarrow E''$ est un homomorphisme. \\

Soient $\rho$ et $\rho'$ deux relations binaires de base $E$ et $E'$ respectivement. Une application $f$ de $E$ dans $E'$ est un \emph{isomorphisme} si
\begin{enumerate}
\item $f$ est bijective,
\item $(x,y)\in\rho$ si et seulement si $(f(x),f(y))\in\rho'$ pour tous $x,y\in E.$
\end{enumerate}

\vspace{1mm}

Donc un homomorphisme bijectif de $(E,\rho)$ dans $(E',\rho')$ est un isomorphisme de $(E,\rho)$ sur $(E',\rho')$ si sa bijection r\'eciproque est un homomorphisme de $(E',\rho')$ dans $(E,\rho)$.

\vspace{2mm}

Soient $P$ et $Q$ deux ensembles ordonn\'es. %relations binaires d\'efinies sur $E$ et $E'$ respectivement.
 Une application $f:P\rightarrow Q$ est un \emph{isomorphisme} si et seulement si elle est croissante, bijective et son inverse est croissante.

Il est \`a noter que l'inverse d'une application croissante bijective de $P$ sur $Q$ n'est pas forc\'ement croissante, elle l'est toutefois si ou bien $P$ est une cha\^{i}ne ou bien $P=Q$ et $P$ fini.

\vspace{1mm}

Une application $f:P\rightarrow Q$ est un \emph{plongement}\index{plongement} si et seulement si
\begin{enumerate}
\item $f$ est injective,
\item $x\leq y$ dans $P$ si et seulement si $f(x)\leq f(y)$ dans $Q$ pour tous $x,y\in P.$
\end{enumerate}

\vspace{1mm}

Donc une application $f$ est un plongement de $P$ dans $Q$ si et seulement si $f$ est un isomorphisme de $P$ sur une restriction de $Q$.

\subsection{Constructions et propri\'et\'es des ensembles ordonn\'es}
\subsubsection{Ensemble ordonn\'e et graphe}
Soit $P:=(E,\leq)$ un ordre ou ensemble ordonn\'e. Le \emph{graphe de comparabilit\'e}\index{graphe!de comparabilit\'e} de $P$ est le graphe, not\'e $Comp(P)$, ayant pour sommets les \'el\'ements de l'ensemble $E$ et pour ar\^etes les paires d'\'el\'ements distincts de $E$ qui sont comparables par $\leq$. Le \emph{graphe d'incomparabilit\'e}\index{graphe!d'incomparabilit\'e} de $P$, not\'e $Incomp(P)$ est le compl\'ementaire de $Comp(P)$, il a pour ar\^etes\index{arete@ar\^ete} les paires d'\'el\'ements incomparables de $P$. Ainsi, si $P$ est un ensemble ordonn\'e, les cliques\index{clique} et les stables de $Comp(P)$ sont respectivement les cha\^{i}nes et les anticha\^{i}nes de $P$. De m\^eme, les cliques et les stables de $Incomp(P)$ sont respectivement les anticha\^{i}nes et les cha\^{i}nes de $P$.

   % \subsection{Op\'erations sur les ensembles ordonn\'es:}
  % \subsection{Constructions et d\'ecompositions de graphes et d'ordres:}

  \subsubsection{Somme directe, somme ordinale et produit d'ensembles ordonn\'es}\label{subsec:somme et produit}
Soient $P$ et $Q$ deux ensembles ordonn\'es, nous supposons leurs ensembles de sommets disjoints. La \emph{somme directe}\index{somme!directe} de $P$ et $Q$, not\'ee $P\oplus Q$,  est l'ensemble ordonn\'e d\'efini sur la r\'eunion des deux ensembles de sommets $P\cup Q$ de la mani\`ere suivante:
$$x\leq_{P\oplus Q}y \text{ si }(( x,y\in P \text{ et } x\leq_P y) \text{ ou bien }( x,y\in Q \text{ et }x\leq_Q y)).$$
La somme directe  $P\oplus Q$ se repr\'esente par les deux diagrammes de Hasse de $P$ et $Q$ plac\'es l'un \`a c\^ot\'e de l'autre sans aucune relation entre les sommets des deux ordres.\\

\begin{figure}[h]
\centering
\input{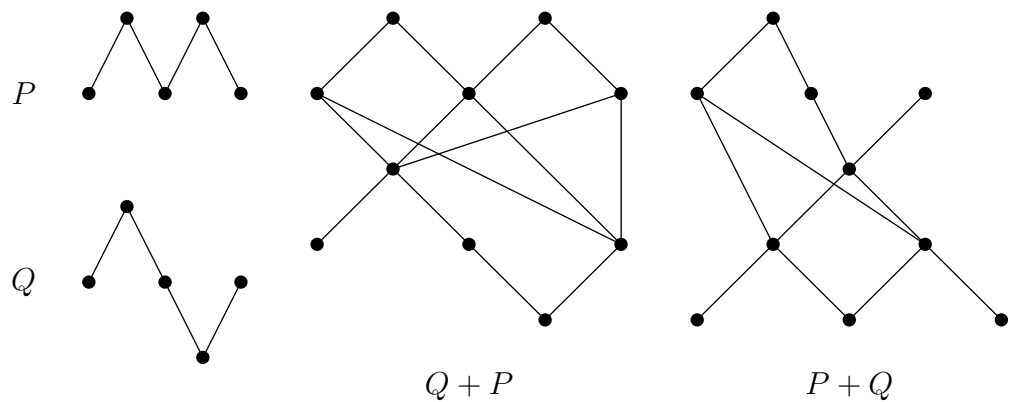}
\caption{\label{somme-ordinal}Diagramme de Hasse de la somme ordinale de deux ordres}
\end{figure}

La \emph{somme ordinale}\index{somme!ordinale} ou \emph{somme lexicographique}\index{somme!lexicographique} de $P$ et $Q$, not\'ee $P+Q$, est l'ordre d\'efini sur la r\'eunion des deux ensembles de sommets par (voir \figurename~\ref{somme-ordinal}):
$$x\leq_{P+Q}y \text{ si }(( x,y\in P \text{ et } x\leq_P y) \text{ ou bien }( x,y\in Q \text{ et }x\leq_Q y) \text{ ou bien }(x\in P \text{ et }y\in Q)).$$

Le \emph{produit direct}\index{produit direct} de $P$ et $Q$, not\'e $P\times Q$, est l'ordre d\'efini sur le produit $P\times Q$ des deux ensembles de sommets par:
$$(x,y)\leq_{P\times Q}(x',y') \text{ si } x\leq_P x' \text{ et } y\leq_Q y'.$$
En particulier, si $(P_i)_{1\leq i\leq n}$ est une suite finie d'ensembles ordonn\'es, le produit direct de ces ensembles est l'ensemble $\underset{1\leq i\leq n}\prod P_i$ muni de l'ordre produit donn\'e par: $(x_1,\cdots,x_n)\leq(y_1,\cdots,y_n)$ dans $\underset{i=1}{\overset{n}\prod}P_i$ si et seulement si $x_i\leq y_i$ dans $P_i$ pour tout $1\leq i\leq n$.

%    \subsubsection{Quelques \'el\'ements et sous-ensembles particuliers d'un ensemble ordonn\'e.}

    \subsection{El\'ements particuliers d'un ensemble ordonn\'e}
   Soit $\mathcal P$ un ensemble ordonn\'e. Un \'el\'ement $x$ de $\mathcal P$ est \emph{minimal}\index{element@\'el\'ement!minimal} (respectivement \emph{maximal})\index{element@\'el\'ement!maximal} dans $\mathcal P$ s'il n'existe pas d'\'el\'ement qui lui est strictement inf\'erieur (respectivement strictement sup\'erieur), c'est \`a dire si $y\in\mathcal P$ tel que $y\leq x$ (respectivement $y\geq x$) implique $y=x$. L'\'el\'ement $x$ est le \emph{plus petit \'el\'ement}\index{element@\'el\'ement!le plus petit -} ou le \emph{minimum} (respectivement le \emph{plus grand \'el\'ement} ou le \emph{maximum})\index{element@\'el\'ement!le plus grand -} de $\mathcal P$ si $y\geq x$ (respectivement $y\leq x$) pour tout $y$ dans $\mathcal P$.

   \vspace{1mm}

  Etant donn\'e un ensemble ordonn\'e $\mathcal P$, il n'existe pas forc\'ement un plus grand \'el\'ement (resp. un plus petit \'el\'ement) de $\mathcal P$.
  \begin{theorem}
  Si un ensemble ordonn\'e $\mathcal P$ est fini et non vide il contient au moins un \'el\'ement maximal et un \'el\'ement minimal.
  \end{theorem}

  Dans le cas infini, ce th\'eor\`eme n'est pas toujours vrai.
\medskip

 Soient $\mathcal P$ un ensemble ordonn\'e, $A$ une partie de $\mathcal P$ et $m$ un \'el\'ement de $\mathcal P$. Nous disons que $m$ minore (resp. majore) $A$ si $m\leq a$ (resp. $a\leq m$) pour tout $a\in A$. Nous noterons $A^-$ (resp. $A^+$) l'ensemble des \'el\'ements de $\mathcal P$ qui minorent (resp. majorent) $A$.\\ Lorsque $A^-$ (resp. $A^+$) poss\`ede un plus grand \'el\'ement (resp. un plus petit \'el\'ement) nous l'appelons la \emph{borne inf\'erieure} (resp. la \emph{borne sup\'erieure}) de $A$, nous la notons $inf A$ ou $\bigwedge A$ (resp. $sup A$ ou $\bigvee A$).

 \vspace{2mm}

 L'ensemble ordonn\'e $\mathcal P$ est dit \emph{inductif}\index{ensemble!inductif} si toute cha\^{i}ne non vide de  $\mathcal P$ poss\`ede un majorant.\\
 Dans ce cas, le \emph{lemme de Zorn} garantit l'existence d'un \'el\'ement maximal\index{element@\'el\'ement!maximal}. La validit\'e du lemme de Zorn est \'equivalente \`a l'axiome du choix.\\

\textbf{Lemme de Zorn:} {\em Tout ensemble inductif poss\`ede au moins un \'el\'ement maximal.}\index{Lemme de Zorn}

 \vspace{2mm}

   Un \emph{treillis}\index{treillis} ("lattice" en anglais) est un ensemble ordonn\'e $\mathcal P$ dont toute paire d'\'el\'ements (et donc toute partie finie) admet une borne sup\'erieure et une borne inf\'erieure. Si toute partie de $\mathcal P$ poss\`ede une borne inf\'erieure et une borne sup\'erieure, c'est un \emph{treillis complet}\index{treillis!complet} ("complete lattice"). Si $\mathcal P$ est un treillis (resp. treillis complet), une partie $T$ de $\mathcal P$ v\'erifiant $\bigvee\{a, b\}\in T$ et $\bigwedge\{a, b\}\in T$ pour tous $a,b\in T$ est un \emph{sous-treillis} (resp. sous-treillis complet) de $\mathcal P$.

             \subsection{Sections initiales et id\'eaux d'un ensemble ordonn\'e}
    Soit $\mathcal P$ un ensemble ordonn\'e. Une \emph{section initiale}\index{section initiale} (dite aussi \emph{segment initial})\index{segment initial} est toute partie $I$ de $\mathcal P$ telle que:
$$x\in I \text{ et } y\leq x\;\text{impliquent } y\in I.$$
Nous noterons $\underline{I}(\mathcal P)$ l'ensemble, ordonn\'e par inclusion, des sections initiales de $\mathcal P$.
L'ensemble $\mathcal P$ et l'ensemble vide $\varnothing$, sont des sections initiales de $\mathcal P$. L'intersection et la r\'eunion d'une famille quelconque de sections initiales sont des sections initiales.\\
Le compl\'ementaire d'une section initiale est appel\'e \emph{section finale}\index{section finale}, c'est toute partie $F$ de $\mathcal P$ v\'erifiant:
$$x\in F \text{ et } x\leq y\;\text{impliquent } y\in F.$$
Nous noterons $\underline{F}(P)$ l'ensemble, ordonn\'e par inclusion, des sections finales de $\mathcal P$. \\

Etant donn\'ee une partie $A$ de $\mathcal P$, l'ensemble $$\downarrow A:=\{x\in\mathcal P: x\leq y\; \text{pour au moins un } y\in A\}=\underset{a\in A}\bigcup\{a\}^-$$ (respectivement $\uparrow A:=\{x\in\mathcal P: x\geq y \;\text{pour au moins un } y\in A\}=\underset{a\in A}\bigcup\{a\}^+$) est la plus petite section initiale (respectivement section finale) contenant $A$. Cette section est \emph{engendr\'ee} par $A$. Une section initiale ou finale est dite \emph{finiment engendr\'ee}\index{section initiale!finiment engendr\'ee} lorsqu'elle est engendr\'ee par une partie $A$  finie. Si la partie $A$ est un singleton, disons $A=\{a\}$, nous notons $\downarrow a$ au lieu de $\downarrow A$ ce segment initial et nous disons qu'il est \emph{principal}.\index{section initiale!principale}\\
Si une section initiale $I$ est finiment engendr\'ee alors elle est engendr\'ee par ses \'el\'ements maximaux.

    \subsubsection{Partie filtrante et id\'eal d'un ensemble ordonn\'e}
Une partie $A$ d'un ensemble ordonn\'e $\mathcal P$ est dite \emph{filtrante\index{ordre!partie filtrante d'un -} sup\'erieurement} si deux \'el\'ements quelconques $x$ et $y$ de $A$ admettent toujours un majorant commun dans $A$ (c'est \`a dire qu'il existe un $z$ dans $A$ tel que $x\leq z$ et $y\leq z$).\\ Un \emph{id\'eal}\index{ideal@id\'eal} $J$ de $\mathcal P$ est une section initiale non vide  de $\mathcal P$ qui est filtrante sup\'erieurement. Par exemple, tout segment initial non vide d'une cha\^{i}ne est un id\'eal, ou encore tout segment initial engendr\'e par un singleton $\{a\}, ~a\in \mathcal P.$\\

\noindent Nous rappelons les propri\'et\'es suivantes d'un ensemble ordonn\'e  $\mathcal P$ :
\begin{properties}\label{propriete-ideaux}
\begin{enumerate}
\item Une section initiale $J$ de $\mathcal P$ est un id\'eal si et seulement si $J\neq\varnothing$ et, pour toutes sections initiales $I_1$ et $I_2$, si $J=I_1\cup I_2$ alors $J=I_1$ ou $J=I_2$ (ou, ce qui revient au m\^eme, si $J\subseteq I_1\cup I_2$ alors $J\subseteq I_1$ ou $J\subseteq I_2$).
\item  Une union d'un ensemble totalement ordonn\'e pour l'inclusion, ou simplement filtrant pour l'inclusion, d'id\'eaux est un id\'eal.
\item Tout id\'eal $J$ de $\mathcal P$ est contenu dans un id\'eal maximal\index{ideal@id\'eal!maximal} (pour l'inclusion) de $\mathcal P$.
\end{enumerate}
\end{properties}
%%%%%%%%%%%%%%%%%%%%%%%%%%%%%%%%%%%%%%%%%%%%%%%%%%%%%%%%%%%%%%%%%%%%%%%%%%%%%%%%%%%%%%%%%%%%%%%%%%%

        \subsection{Ordres bien fond\'es et ordres belordonn\'es}\label{sec:belordre}

        Un ordre $P$ sur $E$ est un \emph{bonordre}\index{bonordre} et $E$ est \emph{bien ordonn\'e}\index{ensemble!bien ordonn\'e} si toute partie non vide de $E$ a un minimum pour $P$. Toute cha\^{i}ne finie\index{chaine@cha\^{i}ne!finie} est un bonordre. L'ordre naturel sur $\mathbb N$ est un bonordre, par contre son ordre dual n'est pas un bonordre.

        \subsubsection{Ordres bien fond\'es}\label{subsec:bienfondé}
Un ensemble ordonn\'e $\mathcal P$ est \emph{bien fond\'e}\index{ordre!bien fond\'e} si toute partie non vide $A$ de $\mathcal P$ a au moins un \'el\'ement minimal. Il revient au m\^eme %(avec  l'axiome des choix d\'ependants)
de dire qu'il n'existe pas de suite infinie strictement d\'ecroissante d'\'el\'ements de $\mathcal P$, soit $x_{0} > x_{1} > \ldots > x_{n} >\ldots$.
Notons qu'une cha\^{i}ne bien fond\'ee est un bonordre. Tout ordre fini est bien fond\'e et toute restriction d'un ordre bien fond\'e est bien fond\'e.\\
%Les types d'isomorphie des ensembles bien ordonn\'es\index{ensemble!bien ordonn\'e} sont les {\it ordinaux}, outils indispensables pour l'induction.
Cette notion \'etend aux ensembles ordonn\'es la notion de bonordre.

 \subsubsection{D\'ecomposition en niveaux d'un ensemble ordonn\'e}\index{decomposition@d\'ecomposition!en niveaux}
Etant donn\'e un ensemble ordonn\'e bien fond\'e $\mathcal P$, posons $Min(\mathcal P):=\{a\in\mathcal P: a\text{ minimal dans }\mathcal P\}$. Posons \'egalement $P_0:=Min(\mathcal P)$, $P_1:=Min(\mathcal P\setminus P_0)$, $P_i:=Min(\mathcal P\setminus\underset{j<i}\cup P_j)$. Les sous-ensembles $P_i$ sont des anticha\^{i}nes et sont deux \`a deux disjoints. Ces ensembles sont appel\'es \emph{niveaux}\index{ordre!niveau d'un -}. Par exemple, si $\mathcal P$ est une cha\^{i}ne, chaque niveau est r\'eduit \`a un \'el\'ement. Si $\mathcal P$ est fini, alors il existe un entier $m$ tel que $P_m=\varnothing$. Le plus petit entier $m$ v\'erifiant ceci est la \emph{hauteur}\index{ordre!hauteur d'un -} de $\mathcal P$, not\'ee $h(\mathcal P)$. Exemple, si $\mathcal P$ est une cha\^{i}ne \`a $m$ \'el\'ements, la hauteur de $\mathcal P$ est $h(\mathcal P)=m$.

  Nous pouvons d\'efinir la hauteur de $\mathcal P$ en utilisant l'ensemble des cha\^{i}nes de $\mathcal P$. En effet, l'ensemble des cha\^{i}nes de $\mathcal P$ est partiellement ordonn\'e par inclusion et ses \'el\'ements maximaux sont appel\'es \emph{cha\^{i}nes maximales}\index{chaine@cha\^{i}ne!maximale}. %Une cha\^{i}ne $C$ est une cha\^{i}ne maximum si aucune autre cha\^{i}ne ne contient plus d'\'el\'ements que $C$.
 Si $\mathcal P$ est fini, la hauteur de $\mathcal P$ est le maximum du nombre d'\'el\'ements des cha\^{i}nes de $\mathcal P$.

 \subsubsection{Ordres belordonn\'es}
Un ensemble ordonn\'e $\mathcal P$ est \emph{belordonn\'e}\index{ordre!belordonn\'e} (en anglais \emph{well-quasi-ordered}, ou \emph{wqo}), si toute partie non vide $A$ de $\mathcal P$ contient un nombre fini, non nul, d'\'el\'ements minimaux\index{element@\'el\'ement!minimal}.
Par exemple, tout ensemble fini est belordonn\'e, tout bonordre est un belordre et toute restriction d'un belordre est un belordre.\\
La notion de belordre permet d'\'etendre \`a des ensembles ordonn\'es infinis les propri\'et\'es des ensembles ordonn\'es finis.

    \subsubsection{Principaux r\'esultats et propri\'et\'es}
Le th\'eor\`eme suivant, d\^u \`a Higman  1952 \cite{higman52} (voir aussi \cite{fraisse}) donne des conditions \'equivalentes d'un belordre.

\begin{theorem}\label{theo:higman-equivalence}
Soit $\mathcal P$ un ensemble ordonn\'e. Les propri\'et\'es suivantes sont \'equivalentes:
\begin{enumerate}
\item[(i)] $\mathcal P$ est belordonn\'e ;
\item[(ii)] $\mathcal P$ est bien fond\'e\index{ensemble!bien fond\'e} sans anticha\^{\i}ne infinie;
\item[(iii)] Toute suite infinie $x_{0}, \ldots,
x_{n}, \ldots$, d'\'el\'ements de $\mathcal P$, contient une sous suite infinie croissante $x_{i_{0}} \leq
x_{i_{1}} \leq \ldots x_{i_{n}} \leq \ldots$;
\item[(iv)] Toute suite infinie $x_{0}, \ldots,
x_{n}, \ldots$, d'\'el\'ements de $\mathcal P$, contient un "couple" croissant $x_{i} \leq x_{j}$ avec $i<j$;
\item[(v)] Toute section finale de $\mathcal P$ est finiment engendr\'ee ;
\item[(vi)] L'ensemble $\underline{F}(\mathcal P)$ des sections finales de $\mathcal P$ muni de
l'ordre oppos\'e de l'inclusion est bien fond\'e\index{ensemble!bien fond\'e};
\item[(vii)] L'ensemble $\underline{I}(\mathcal P)$ des sections initiales de $\mathcal P$ muni de
l'ordre d'inclusion est bien fond\'e.
\end{enumerate}
\end{theorem}

\begin{proposition}
Soit $P$ un ensemble ordonn\'e. Si $P=A\cup B$ alors $P$ est belordonn\'e d\`es que $A$ et $B$ le sont.
\end{proposition}

Voici quelques propri\'et\'es du belordre\index{ensemble!belordonn\'e}:

\begin{properties}\label{propr:belordre}
\begin{enumerate}
 \item[1)] Soient $P$ et $Q$ deux ensembles ordonn\'es et $f$ une application de $P$ dans $Q$.
    \begin{enumerate}
    \item[$a$)]  Si $f$ est un plongement et $Q$ un belordre alors $P$ est un belordre.
     \item[$b$)] Si $f$ est \emph{croissante} et surjective et $P$ un belordre alors $Q$ est un belordre.
    \end{enumerate}
 \item[2)] Si deux ordres $P$ et $Q$ sont belordonn\'es, il en est de m\^eme pour leur somme directe $P\oplus Q$ et leur produit direct $P\times Q$.
 \end{enumerate}
 \end{properties}

  \begin{corollary}
  Soit $(P_i)_{i=0,\cdots,n-1}$ une suite finie d'ensembles ordonn\'es, le produit $\underset{i=0}{\overset{n-1}\prod}P_i$ muni de l'ordre produit est belordonn\'e si et seulement si chaque $P_i$ est belordonn\'e. En particulier, si $P$ est belordonn\'e alors $P^k$ est belordonn\'e pour tout entier $k$. %Fraissé p 127
  \end{corollary}

        \subsubsection{Th\'eor\`eme d'Higman}
  Soit $E$ un ensemble. Un \emph{mot}\index{mot} fini sur $E$ est une suite finie $u=a_0a_1\cdots a_{k-1}$ d'\'el\'ements de $E$. Soit $E^{\star}$ l'ensemble des mots finis sur $E$. L'ensemble $E$ \'etant ordonn\'e, ordonnons $E^{\star}$ ainsi:
\'etant donn\'es deux \'el\'ements $x=(x_i)_{i=0,\cdots,n-1}$ et $y=(y_i)_{i=0,\cdots,m-1}$ de $E^{\star}$, notons $x\leq^{\star}y$ lorsqu'il existe une injection croissante $h:[0,n[\rightarrow [0,m[$ telle que pour tout $i=0,\cdots,n-1$ on ait $x_i\leq y_{h(i)}$. Cet ordre s'appelle  l'ordre de Higman.

 Nous avons le r\'esultat suivant:

\begin{theorem}\label{theo:higman}(Higman 52)\\
$E^{\star}$ muni de l'ordre $\leq^{\star}$ est belordonn\'e si et seulement si $E$ est belordonn\'e.
\end{theorem}

Nous rappelons \'egalement cette  autre propri\'et\'e, due \`a Erd\"{o}s et Tarski (43), %est donn\'ee par le th\'eor\`eme suivant:

\begin{theorem}\label{theo:erdos-tarski} Erd\"os-Tarski (1943).\\
Un ensemble belordonn\'e\index{ensemble!belordonn\'e} est une union finie d'id\'eaux.
\end{theorem}

    \subsubsection{Th\'eor\`eme de Ramsey}
Nous rappelons \'egalement le th\'eor\`eme de partition de Ramsey\index{Ramsey} qui est d'un usage fr\'equent dans l'\'etude du belordre;
\begin{theorem}\label{thm:ramsey}
Les paires d'entiers \'etant r\'eparties en deux classes  disjointes, il existe une partie infinie d'entiers dont toutes les paires sont contenues dans une m\^eme classe.
\end{theorem}
\vspace{1mm}

Ce th\'eor\`eme se g\'en\'eralise au cas de sous-ensemble à $m$ \'el\'ements r\'epartis en $k$ classes disjointes:
\begin{theorem}
Les parties \`a $m$ \'el\'ements de l'ensemble des entiers \'etant r\'eparties en $k$ ($k$ fini), classes disjointes (ces classes sont \'egalement appel\'ees des couleurs), il existe une partie $E$ infinie d'entiers telle que tous les sous-ensembles \`a $m$ \'el\'ements de $E$ ont la m\^eme couleur.
\end{theorem}

\vspace{2mm}

Le th\'eor\`eme de  partition de Ramsey remonte à 1930, il a inspir\'e toute une th\'eorie que l'on appelle "th\'eorie de Ramsey".
Pour plus de d\'etails sur cette th\'eorie voir par exemple \cite{Grah-Rot-Spen,land-robert}.
%$$$$$$$$$$$$$$$$$$$$$$$$$$$$$$$$$$$$$$$$$$$$$$$$$$$$$$$$$$$$$$$$$$$$$$$$$$$$$$$$$$$$$$$$$$

\section{Structures relationnelles}
        \subsection{D\'efinitions et notations}

Soit $n$ un entier positif. Une \emph{relation $n$-aire\index{relation!$n$-aire} sur un ensemble $E$} est un sous-ensemble $\rho$ de $E^n$. L'ensemble $E$ est la \emph{base}\index{relation!base d'une -} ou le \emph{domaine}\index{relation!domaine d'une -} de $\rho$, not\'e aussi $dom(\rho)$. L'entier $n$ est appel\'e l'\emph{arit\'e}\index{arite@arit\'e} de $\rho$. Si $n=1,2,3$ la relation $\rho$ est dite relation \emph{unaire}\index{relation!unaire}, \emph{binaire}\index{relation!binaire}, \emph{ternaire}\index{relation!ternaire}.
Pour tout ensemble $E$ et tout entier positif $n$, l'ensemble $E^n$ est une relation $n$-aire; l'ensemble vide $\varnothing$ est \'egalement une relation $n$-aire, appel\'ee \emph{relation vide.}\index{relation!vide} Etant donn\'ee une relation $n$-aire $\rho$  sur l'ensemble $E$, sa relation  \emph{compl\'ementaire}\index{relation!compl\'ementaire} $\rho^c$ est la relation $\rho^c:= E^n \setminus \rho$, sa \emph{fonction caract\'eristique}\index{fonction!caract\'eristique} est l'application  $f_{\rho}:E^n\rightarrow \{0,1\}$ telle que $\rho:=f_{\rho}^{-1}(\{1\})$. Au besoin, la relation $\rho$ est identifi\'ee à sa fonction caract\'eristique.\\

Une \emph{structure relationnelle}\index{structure relationnelle} est une paire $\mathcal R:=(E,(\rho_i)_{i\in I})$ form\'ee d'un ensemble $E$, appel\'e \emph{base}\index{structure relationnelle!base d'une -} ou \emph{domaine}\index{structure relationnelle!domaine d'une -} de $\mathcal R$ et d'une famille $(\rho_i)_{i\in I}$ de relations $n_i$-aires $\rho_i$ sur $E$. Chaque relation $\rho_i$ est appel\'ee \emph{composante}\index{structure relationnelle!composante d'une -} de la structure relationnelle $\mathcal R$. La famille $\mu:=(n_i)_{i\in I}$ est appel\'ee \emph{signature}\index{structure relationnelle!signature d'une -} de $\mathcal{R}$. La structure relationnelle $\mathcal R$ \'etant donn\'ee, nous noterons $V(\mathcal R)$ son domaine.

 Toute relation $\rho$ sur un ensemble $E$ donne lieu \`a une structure relationnelle $\mathcal R:=(E,\rho)$. Si $\vert I \vert=2,3,4$ la structure $\mathcal R$ est appel\'ee, respectivement, \emph{birelation}\index{birelation}, %%(ou \emph{une $2$-structure})%%
\emph{trirelation}\index{trirelation}, \emph{quadrirelation}\index{quadrirelation}. Quand $I$ est
fini nous disons plut\^ot que $\mathcal R$ est une \emph{multirelation}\index{multirelation}.
%Si toutes les composantes de la structure $\mathcal R$ sont des relations binaires\index{relation!binaire}, $\mathcal R$ est dite \emph{structure relationnelle binaire}\index{structure relationnelle!binaire} ou en bref \emph{structure binaire}, on l'appelle aussi \emph{une $2$-structure.}
Nous d\'esignons par $\Omega_{\mu}$ la classe\index{classe} des structures relationnelle finies\index{structure relationnelle!finie} de signature $\mu$.\\

 Une  structure relationnelle $\mathcal R:=(E,(\rho_i)_{i\in I})$ est \emph{ordonn\'ee}\index{structure relationnelle!ordonn\'ee} si une de ses relations $\rho_i$ est un ordre lin\'eaire; nous conviendrons alors que c'est la premi\`ere relation $\rho_1$. Autrement dit, une structure relationnelle est ordonn\'ee si elle s'\'ecrit  $\mathcal R:=(E,\leq,(\rho_j)_{j\in J})$ o\`u "$\leq$" est un ordre lin\'eaire sur $E$ et les $\rho_j$ sont des relations $n_j$-aire sur $E.$

 Comme exemples de structures ordonn\'ees nous avons les cha\^{i}nes ($J=\varnothing$), les multicha\^{i}nes ($J$ fini et $\rho_j$ est une cha\^{i}ne pour tout $j\in J$) et les graphes ordonn\'es ($\vert J\vert=1$ et $\rho_j$ est une relation binaire). \\

 Soit $\mathcal R:=(E,(\rho_i)_{i\in I})$ une structure relationnelle de signature $\mu:=(n_i)_{i\in I}$. La \emph{sous-structure induite}\index{structure relationnelle!sous-structure induite} par $\mathcal R$ sur un sous-ensemble $A$ de $E$, appel\'ee simplement \emph{restriction}\index{structure relationnelle!restriction d'une -} de $\mathcal R$ \`a $A$, est la structure relationnelle $\mathcal R_{\restriction_A}:=(A,({\rho_i}_{\restriction_A})_{i\in I})$ o\`u ${\rho_i}_{\restriction_A}:={\rho_i}\cap A^{n_i}$. Elle est parfois not\'ee $\mathcal R_A$ et, si $x\in E$, la restriction de $\mathcal R$ \`a $E\setminus\{x\}$ est not\'ee $\mathcal R_{-x}$.  Si $\mathcal R_{\restriction_A}$ est une restriction de $\mathcal R$ alors $\mathcal R$ est une \emph{extension}\index{structure relationnelle!extension d'une -} de $\mathcal R_{\restriction_A}$.

%%%%%%%%%%%%%%%%%%%%%%%%%%%%%%%%%%%%%%%%%%%%%%%%%%%%%%%%%%%%%%%%%%%ùù

    \subsection{Isomorphisme, abritement et \'equimorphisme}\label{sec:isom-abrit-equim}
 Soient $\mathcal{R}:= (E,(\rho_i)_{i\in I})$ et $\mathcal{R'}:= (E',(\rho'_i)_{i\in I})$ deux structures relationnelles de m\^eme signature $\mu:=(n_i)_{i\in I}.$

 Une application  $f:E\rightarrow E'$ est un \emph{isomorphisme}\index{isomorphisme} \emph{de $\mathcal{R}$ sur $\mathcal{R'}$} si
 \begin{enumerate}
 \item $f$ est bijective,
 \item $(x_1,\ldots,x_{n_i})\in \rho_i$  si et seulement si   $(f(x_1),\ldots,f(x_{n_i}))\in \rho'_i$ pour tout
$(x_1,\ldots,x_{n_i})\in E^{n_i}$, $i\in I.$
\end{enumerate}

Si $\mathcal R$ et $\mathcal R'$ sont \'egales le mot isomorphisme est remplac\'e par \emph{automorphisme}.\index{automorphisme}

\vspace{1mm}

La structure relationnelle $\mathcal{R}$ est \emph{isomorphe}\index{structure relationnelle!isomorphe \`a} \`a $\mathcal{R'}$ et nous notons $\mathcal R \cong \mathcal R'$,
s'il existe un isomorphisme de $\mathcal{R}$ sur $\mathcal{R'}$. A toute structure relationnelle $\mathcal R$, nous pouvons associer un \emph{type d'isomorphie}\index{type!d'isomorphie} $\tau (\mathcal R)$, de telle sorte que $\mathcal R \cong \mathcal R'$ si et seulement si
$\tau (\mathcal R)=\tau (\mathcal R')$. Si $\vert dom(\mathcal R)\vert=n$, $\tau (\mathcal R)$ peut-\^etre consid\'er\'e comme la structure relationnelle isomorphe \`a $\mathcal R$ et d\'efinie sur l'ensemble $\{0,\dots,n-1\}$.

La collection des types d'isomorphie des structures relationnelles finies de signature $\mu$ (c'est \`a dire des relations de $\Omega_{\mu}$)  est un ensemble que nous notons $\textit{T}_{\mu}$.\\

Une application $f$ d'un sous-ensemble $A$ de $E$ sur un sous-ensemble $A'$ de $E'$ est un \emph{isomorphisme local}\index{isomorphisme!local} ou un \emph{abritement local}\index{abritement!local} de $\mathcal R$ dans $\mathcal R'$ si $f$ est un isomorphisme de $\mathcal R_{\restriction_A}$ sur $\mathcal R'_{\restriction_{A'}}$. Dans ce cas, $A$ et $A'$ sont respectivement le \emph{domaine}\index{isomorphisme!domaine d'un - local} et le \emph{co-domaine}\index{isomorphisme!co-domaine d'un - local} de $f$.  En particulier, si $f$ est un isomorphisme de $\mathcal{R}$ sur $\mathcal{R'}$, alors toute restriction de $f$ est un isomorphisme local. Dans le cas o\`u $m$ est le maximum de l'arit\'e des relations $(\rho_i)$ et $(\rho'_i)$, une condition suffisante pour qu'une bijection $f$ d\'efinie sur $A$ soit un isomorphisme local de $\mathcal R$ dans $\mathcal R'$ est que la restriction de $f$ \`a tout sous-ensemble $X$ de $A$ de cardinalit\'e au plus $m$ est un isomorphisme local de $\mathcal R$ dans $\mathcal R'$ \cite{fraisse}.\\

  Un isomorphisme $f$ de $E$ sur un sous-ensemble $A'$ de $E'$ est un \emph{abritement}\index{abritement} de $\mathcal R$ dans $\mathcal R'$. Nous disons que $\mathcal{R}$ s'\emph{abrite} dans  $\mathcal{R'}$ ou que $\mathcal R$ est \emph{plus petite} que $\mathcal R'$ pour l'abritement, fait not\'e $\mathcal{R}\leq \mathcal{R'}$, ou que $\mathcal R'$ \emph{abrite} $\mathcal R$ ou est \emph{plus grande} que $\mathcal R$, fait not\'e $\mathcal R'\geq \mathcal R$, si et seulement si il existe une restriction de $\mathcal{R'}$ isomorphe \`a $\mathcal{R}$. Comme dans un tel cas, toute relation $\mathcal S$ isomorphe \`a $\mathcal R$ s'abrite dans toute relation $\mathcal S'$ isomorphe \`a $\mathcal R'$, on dit encore que le type d'isomorphie de $\mathcal R$ s'abrite dans le type d'isomorphie de $\mathcal R'$, ce que nous notons $\tau(\mathcal R)\leq\tau(\mathcal R')$.\\
 Nous disons que $\mathcal{R}$ s'\emph{abrite strictement} dans  $\mathcal{R'}$ ou est \emph{strictement plus petite} que $\mathcal R'$ ou que $\mathcal R'$ \emph{abrite strictement} $\mathcal R$ ou est \emph{strictement plus grande} que $\mathcal R$ et nous notons $\mathcal R<\mathcal R'$ ou $\mathcal R'>\mathcal R$ si et seulement si $\mathcal R\leq \mathcal R'$ mais $\mathcal R'\nleq \mathcal R.$\\

 Nous disons que $\mathcal R$ est \emph{\'equimorphe}\index{structure relationnelle!\'equimorphe \`a} avec $\mathcal R'$ %, not\'ee $\mathcal R \sim\mathcal R'$
 si et seulement si $\mathcal R\leq\mathcal R'$ et $\mathcal R'\leq\mathcal R$. Si les domaines de $\mathcal R$ et $\mathcal R'$ sont finis, l'\'equimorphisme \'equivaut \`a l'isomorphisme.\\

La relation d'abritement\index{relation!d'abritement} est r\'eflexive et transitive, elle d\'efinit un pr\'eordre sur la classe $\Omega_{\mu}$  des structures relationnelles de signature $\mu$ et induit un pr\'eordre sur la classe $\textit{T}_{\mu}$ de leurs types d'isomorphie. Sur l'ensemble $\textit{T}_{\mu}$ c'est un ordre.\\

L'\'etude des structures et classes de structures au moyen de ce pr\'eordre est apparue tr\`es t\^ot dans la litt\'erature. Au d\'ebut du vingti\`eme si\'ecle avec Cantor, Hausdorff et Sierpinski  pour l'\'etude des ordres lin\'eaires, puis dans Fra\"{\i}ss\'e qui a plac\'e ce pr\'eordre et les classes h\'er\'editaires, que nous introduirons ci-dessous, au centre de la th\'eorie des relations \cite{fraisse}.
%%%%%%%%%%%%%%%%%%%%%%%%%%%%%%%%%%%%%%%%%%%%%%%%%%%%%%%%%%%%%%%%%%%%%%%

\section{Classe h\'er\'editaire, \^age, obstruction et borne}\label{sec:clas-heredi-age-borne}
    \subsection{Classe h\'er\'editaire, classe filtrante}

Une classe\index{classe!h\'er\'editaire} $\mathscr{C}$ de structures relationnelles de m\^eme signature $\mu$ est dite \emph{h\'er\'editaire pour l'abritement} ou \emph{ferm\'ee pour l'abritement}\index{abritement!classe ferm\'ee pour l'-} %\footnote{Dans la terminologie des posets, une classe h\'er\'editaire est appel\'ee \emph{section initiale}\index{section!initiale} ou \emph{segment initial}.}
si $\mathcal{R} \in \mathscr{C}$ et $\mathcal{S}\leq \mathcal{R}$ impliquent $\mathcal{S}\in \mathscr{C}$. Dans la terminologie des ensembles ordonn\'es, c'est une section initiale de la classe des structures relationnelles de signature $\mu$.\\ La classe $\mathscr{C}$ est dite \emph{filtrante pour l'abritement}\index{classe!filtrante} si pour toutes structures relationnelles $\mathcal{R}, \mathcal{S}\in \mathscr{C}$ il existe une structure relationnelle $\mathcal{T}\in\mathscr{C}$ telle que $\mathcal{R}\leq \mathcal{T}$ et $\mathcal{S}\leq\mathcal{T}$.

 \noindent Etant donn\'ee une relation $\mathcal R$ et un entier positif $p$, l'ensemble des restrictions de $\mathcal R$ aux sous-ensembles de cardinalit\'es inf\'erieures ou \'egales \`a $p$ est une classe h\'er\'editaire; cependant, elle n'est pas toujours filtrante; par exemple si $\mathcal R$ est une cha\^{\i}ne cette classe\index{classe!filtrante} est filtrante, par contre si $\mathcal R$ est la relation de \emph{cons\'ecutivit\'e}\index{relation!de cons\'ecutivit\'e} sur les entiers (c'est \`a dire $\mathcal R=\{(i,i+1)/~i\in \mathbb N\}$) cette classe n'est filtrante que pour $p=0$ et $p=1$.\\
 L'ensemble des restrictions de $\mathcal R$ aux sous-ensembles finis de cardinalit\'es sup\'erieures ou \'egales \`a $p$ est filtrant, par contre ce n'est pas une classe h\'er\'editaire. %Le compl\'ementaire d'une classe h\'er\'editaire est appel\'e \emph{section finale}\index{section!finale}.\\

\subsection{\^Age}

 Soit $\mathcal R$ une structure relationnelle. L'ensemble des restrictions finies de $\mathcal R$, consid\'er\'ees \`a l'isomorphie pr\`es, est un ensemble appel\'e \emph{\^age}\index{age@\^age} de $\mathcal R$. Nous disons que $\mathcal R$ \emph{repr\'esente cet \^age}. De mani\`ere plus rigoureuse,
 %on appelle
 l'\emph{\^age} de la structure relationnelle $\mathcal R$ est l'ensemble $\mathcal A(\mathcal R)$ des types d'isomorphie des restrictions de $\mathcal R$ aux parties finies de son domaine. Une relation $\mathcal S$ est dite \emph{moins \^ag\'ee} que la relation $\mathcal R$ lorsque l'\^age de $\mathcal S$ est inclus dans l'\^age de $\mathcal R.$ Deux relations $\mathcal R$ et $\mathcal S$ ont le m\^eme \^age si chacune des deux relations est moins \^ag\'ee que l'autre.\\

La caract\'erisation suivante des \^ages est due \`a Fra\"{\i}ss\'e \cite{fraisse}, dans le cas o\`u la signature\index{signature!finie} $\mu$ est finie.
\begin{theorem}\label{theo:age=ideal}
Un ensemble non vide $\mathscr A$ de $\textit{T}_{\mu}$ est l'\^age\index{age@\^age} d'une relation si et seulement si c'est un id\'eal de $\textit{T}_{\mu}$ ordonn\'e par l'abritement.
\end{theorem}
Dans le cas o\`u $\mu$ n'est pas finie, un id\'eal n'est pas forc\'ement un \^age (voir exemple en \cite{Dolho-Pou-Sa-Sau}).\\

\noindent Etant donn\'e un ensemble d'\^ages, totalement ordonn\'e par inclusion,
l'union et l'intersection d'\^ages de cet ensemble sont des \^ages.
 Etant donn\'e un ensemble d'\^ages, filtrant pour l'inclusion, une union d'\^ages de cet ensemble est un \^age.  Cependant, en g\'en\'eral, l'union ou l'intersection de deux \^ages n'est pas n\'ecessairement un \^age.

\vspace{1mm}

Une autre caract\'erisation des \^ages qui d\'ecoule du Th\'eor\`eme \ref{theo:age=ideal} et qui est une propri\'et\'e classique des id\'eaux dans un ensemble ordonn\'e (voir assertion $1$ des Propri\'et\'es \ref{propriete-ideaux}) est:
\begin{theorem}
Un ensemble non vide $\mathcal A$ de relations finies de m\^{e}me arit\'e, qui est ferm\'e pour l'abritement, est un \^age si et seulement si $\mathcal A$ n'est pas l'union de deux sous-ensembles distincts de $\mathcal A$ et ferm\'es pour l'abritement.
\end{theorem}

\subsection{Borne et obstruction}

Soit $\mathscr{C}$ une classe h\'er\'editaire de $\Omega_{\mu}$. Il ressort de la d\'efinition d'une  classe h\'er\'editaire
que si une structure $\mathcal B \notin \mathscr C$ aucun \'el\'ement de $\mathscr C$ n'abritera $\mathcal B$, nous disons, dans ce cas, que $\mathcal B$ est une \emph{obstruction}\index{obstruction} de $\mathscr C$.
%Dans la suite, on s'int\'eresse aux classes de structures relationnelles finies et on note par $\Omega_{\mu}$ la classe des structures relationnelles finies de signature $\mu.$

Soit $\mathscr B$ un sous-ensemble de $\Omega_\mu$, posons
$\text{\it Forb}(\mathscr{B}):=\{\mathcal R \in \Omega_{\mu }: \mathcal B \nleq \mathcal R,~~\forall\mathcal B \in \mathscr{B}\}$ la sous-classe des \'el\'ements de $\Omega_\mu$ qui n'abritent aucun \'el\'ement de $\mathscr B$. Clairement, $\text{\it Forb}(\mathscr{B})$ est une classe h\'er\'editaire et toute sous-classe h\'er\'editaire $\mathscr C$ de $\Omega_\mu$ a cette forme. En fait, si nous consid\'erons les structures de $\Omega_\mu$ \`a l'isomorphie pr\`es, il existe une correspondance bijective entre les classes h\'er\'editaires et les familles d'anticha\^{\i}nes de $\Omega_{\mu}$; en effet, \`a une classe h\'er\'editaire $\mathscr C$ de $\Omega_\mu$, nous associons l'anticha\^{i}ne $\mathcal F_{\mathscr C}:=min(\Omega_\mu\setminus \mathscr C)$ form\'ee par les \'el\'ements minimaux pour l'abritement de $\Omega_\mu\setminus \mathscr C$ (ces \'el\'ements minimaux existent car l'ensemble $\Omega_{\mu}$ est bien fond\'e pour la relation d'abritement). Si $A$ est une anticha\^{i}ne de $\Omega_\mu$, nous lui associons la classe h\'er\'editaire  $\text{\it Forb}(A).$ \\

Ce fait, d\^u \`a Fra\"{\i}ss\'{e}, est bas\'e sur la notion de borne: une \emph{borne}\index{borne} d'une sous-classe h\'er\'editaire $\mathscr{C}$ de $\Omega_\mu$ est toute structure finie $\mathcal B$ n'appartenant pas \`a $\mathscr{C}$  telle que toute structure $\mathcal B'$ qui s'abrite strictement dans $\mathcal B$ appartient \`a $\mathscr{C}$. Il est clair que toute obstruction finie de $\mathscr{C}$ contient une borne. Si $\mathscr B(\mathscr{C})$ d\'esigne la collection des bornes de $\mathscr{C}$, consid\'er\'ees \`a l'isomorphie pr\`es, alors $\mathscr{C}=Forb(\mathscr B(\mathscr C))$, on note aussi
$\mathscr{C}=\Omega_\mu <\mathscr B(\mathscr C)>.$ Notons qu'une borne est n\'ecessairement non vide.

 \vspace{1mm}

 \noindent\textbf{Exemple de borne}: si $\mathscr{C}$ est une classe h\'er\'editaire de graphes finis, une borne de $\mathscr{C}$ est tout graphe $G=(V,E)$ tel que $G\notin \mathscr C$ et $G_{\restriction_{V\setminus\{x\}}}\in \mathscr C$ pour tout $x\in V.$ Si $\mathscr{B}=\{P_4\}$, alors $Forb(\mathscr B)$ est la classe\index{graphe!sans $P_4$} des graphes sans $P_4$ ou graphes s\'eries-parall\`eles\index{graphe!s\'eries-parall\`eles} appel\'es \'egalement cographes\index{cographe} ($P_4$ est le chemin d'ordre quatre, c'est \`a dire le graphe $G$ dont les sommets sont  $V(G):=\{x_1,x_2,x_3,x_4\}$ et tel que $x_i$ est adjacent \`a $x_{i+1}$ pour tout $i=1,2,3$).

 \vspace{1mm}

Cette d\'efinition de borne s'\'etend aux types d'isomorphie. Ainsi, une borne d'une classe h\'er\'editaire $\mathscr C$ de $\textit{T}_{\mu}$ est le type d'isomorphie d'une structure finie minimale parmi celles n'appartenant pas \`a $\mathscr C$.
%%%%%%%%%%%%%%%%%%%%%%%%%%%%%%%%%%%%%%%%%%%%%%%%%%%%%%%%%%%%%%%%%%%%%%%%%%%%%%%%%%%%%%%%%%%%%%%%%%%%%%%%%

\subsection{Classes h\'er\'editaires et id\'eaux}

Le r\'esultat suivant relie les classes h\'er\'editaires et les id\'eaux.
\begin{theorem}\label{theo:classehered-ideal}
Toute classe h\'er\'editaire infinie de $\Omega_{\mu}$, avec $\mu$ fini, contient un id\'eal infini\index{ideal@id\'eal!infini}.
\end{theorem}
Dans le cas o\`u $\mu$ n'est pas fini, nous avons le contre exemple suivant. Pour chaque $n\in\mathbb N$ soit $\mathcal R_n:=(\{0\},(u_k)_{k\in\mathbb N})$ où $u_k$ est la  relation unaire telle que $u_k(0)=1$ si $k=n$ et $u_k(0)=0$ pour $k\neq n$. Soit $\mathscr C$ la classe form\'ee des structures $\mathcal R_n$ et de la structure de base vide. \\

Le r\'esultat donn\'e par le Th\'eor\`eme \ref{theo:classehered-ideal} est d\^u \`a Pouzet, mais n'ayant pas \'et\'e publi\'e, nous reprenons sa preuve ici, elle n\'ecessite la proposition suivante:
\begin{proposition}\label{prop:classe-ideal}
Soit $\mathscr C$ une classe h\'er\'editaire de $\Omega_{\mu}$, ($\mu$ fini). Si $\mathscr C$ contient une anticha\^{i}ne infinie alors $\mathscr C$ contient un id\'eal infini $\mathscr I$ qui est belordonn\'e et qui a un nombre infini de bornes.
\end{proposition}

\begin{proof}
Si $\mathscr C$ poss\`ede une anticha\^{i}ne infinie, disons $A$, alors la section finale $(\uparrow A)\cap\mathscr C$ n'est pas finiment engendr\'ee. L'ensemble des sections finales non finiment engendr\'ees de $\mathscr C$, ordonn\'ees par l'inclusion, est un ensemble inductif, %toute chaine non vide a un majorant.
donc, d'apr\`es le Lemme de Zorn, il existe une section finale $\mathscr F_0$, non finiment engendr\'ee, maximale pour l'inclusion. Posons $\mathscr I=\mathscr C\setminus\mathscr F_0$, le compl\'ementaire de $\mathscr F_0$ dans $\mathscr C$, $\mathscr I$ est donc belordonn\'e. Posons $\mathscr B$ l'ensemble des \'el\'ements minimaux de $\mathscr F_0$. A l'isomorphie pr\`es, $\mathscr B$ est l'ensemble des bornes\index{ensemble!des bornes} de $\mathscr I$ et $\mathscr B$ est infini (c'est une anticha\^{i}ne infinie\index{antichaine@anticha\^{i}ne!infinie} et pour tout $\mathcal B\in\mathscr B$, $\mathcal B_{-x}\in\mathscr I$).\\
$\mathscr I$ est infini, car s'il \'etait fini, il poss\`ederait une structure $\mathcal R$ de cardinalit\'e maximum et donc les bornes de $\mathscr I$ seraient de tailles au plus $\vert\mathcal R\vert+1$ et seraient donc en nombre fini, car $\mu$ est fini.

\noindent Reste \`a montrer que $\mathscr I$ est une classe filtrante. Pour cela, montrons d'abord que tout \'el\'ement $\mathcal R$ de $\mathscr I$ est major\'e par presque toutes ses bornes. En effet, soit $\mathscr B_{\mathcal R}\subseteq\mathscr B$ l'ensemble des bornes contenant $\mathcal R$. %Supposons $\mathscr B_{\mathcal R}$ fini pour tout $\mathcal R\in\mathscr I$ et soit $\mathscr B_{\mathcal R_0}$ l'ensemble de cardinalit\'e maximum. Posons $\vert\mathscr B_{\mathcal R_0}\vert=k$ et consid\'erons une partition $\mathscr P=\{\mathscr B_0,\mathscr B_1,\ldots,\mathscr B_n,\ldots\}$ de $\mathscr B$ en ensembles de cardinalit\'e $k$, cette partition est \'evidemment infinie. Consid\'erons les ensembles $(\downarrow\mathscr B_i\setminus\mathscr B_i)_{i\geq 0}$. Ces ensembles sont contenus dans $\mathscr I$ et chaque ensemble contient un \'el\'ement $\mathcal R_i$ qui n'appartient pas aux autres ensembles. Ces \'el\'ements $(\mathcal R_i)_i$ forment une anticha\^{i}ne infinie de $\mathscr I$ ce qui est une contradiction puisque $\mathscr I$ est wqo.\\
Consid\'erons la classe $\mathscr F_1:=\mathscr F_0\cup\uparrow\mathcal R$. Cette classe est une section finale qui contient $\mathscr F_0$. Comme $\mathscr F_0$ est maximale parmi les classes finales non finiment engendr\'ees, alors $\mathscr F_1$ est finiment engendr\'ee. Or $Min(\mathscr F_1)=\{\mathcal R\}\cup(\mathscr B\setminus\uparrow\mathcal R)$, ce qui prouve que $\mathscr B_{\mathcal R}$ qui est \'egal à $\mathscr B\cap\uparrow\mathcal R$ est infini.

Soient maintenant $\mathcal R,~\mathcal R'\in\mathscr I$, posons $\vert\mathcal R\vert=n$ et $\vert\mathcal R'\vert=n'$. D'apr\`es ce qui pr\'ec\`ede,
$\mathscr B\setminus(\uparrow\mathcal R\cap\uparrow\mathcal R')$ est fini, c'est à dire $\mathscr B\cap(\uparrow\mathcal R\cap\uparrow\mathcal R')=\mathscr B_{\mathcal R}\cap\mathscr B_{\mathcal R'}$ est infini. Comme  $\mu$ fini, cet ensemble  poss\`ede des \'el\'ements de tailles arbitrairement grandes, en particulier de tailles sup\'erieures \`a $n+n'+1$. Donc il existe une borne $\mathcal S$ telle que $\mathcal R<\mathcal S$ et $\mathcal R'<\mathcal S$. Donc il existe $\mathcal S'<\mathcal S$ qui abrite $\mathcal R$ et $\mathcal R'$. Comme $\mathcal S$ est une borne, $\mathcal S'\in \mathscr I$ %D'apr\`es ce qui pr\'ec\`ede, il existe une infinit\'e de bornes qui majorent $\mathcal R$ et $\mathcal R'$, , d'o\`u $\mathcal R\cup\mathcal R'\in\mathscr I$
et par suite $\mathscr I$ est filtrante.
\end{proof}

\vspace{2mm}

\textbf{\emph{Preuve du Th\'eor\`eme \ref{theo:classehered-ideal}.}}  Soit $\mathscr C$ une classe h\'er\'editaire infinie de $\Omega_{\mu}$, nous avons deux cas:

\textbf{Cas 1:} Si $\mathscr C$ ne poss\`ede pas d'anticha\^{i}ne infinie, dans ce cas $\mathscr C$ est belordonn\'ee et est donc une union finie d'id\'eaux (Th\'eor\`eme \ref{theo:erdos-tarski} d'Erd\"{o}s-Tarski). Comme $\mathscr C$ est infinie, l'un des id\'eaux est infini.

\textbf{Cas 2:} Si $\mathscr C$ poss\`ede une anticha\^{i}ne infinie. Le r\'esultat d\'ecoule alors de la Proposition \ref{prop:classe-ideal}.
\hfill $\Box$
%%%%%%%%%%%%%%%%%%%%%%%%%%%%%%%%%%%%%%%%%%%%%%%%%%%%%%%%%%%%%%%%%%%%%%%%%%%%%%%%%%%%%%%%%%

\subsection{Classes h\'er\'editairement belordonn\'ees}\label{subsection:heredit-belordonne}
Soit $\mathscr{C}$ une sous-classe de $\Omega_\mu$ et soit $\mathcal{A}$ un poset. Posons $$\mathscr{C}.\mathcal{A}:=\{(\mathcal{R},f):\mathcal{R}\in \mathscr{C}, f:\mathcal{R}\rightarrow \mathcal{A}\}$$ et pour deux \'el\'ements $(\mathcal{R},f)$ et $(\mathcal{R'},f')$ de $\mathscr{C}.\mathcal{A}$,
$(\mathcal{R},f)\leq (\mathcal{R'},f')$ s'il existe un abritement $h:\mathcal{R}\rightarrow \mathcal{R'}$ tel que  $f(x)\leq f'(h(x))$ pour tout  $x\in
V(\mathcal{R})$.\\ %On rappelle que $\mathcal{A}$ est
%\textit{belordonn\'e (wqo)} si $\mathcal{A}$ ne contient pas d'anticha\^{i}ne infinie ou une cha\^{i}ne strictement d\'ecroissante.
Nous disons que $\mathscr{C}$ est \textit{h\'er\'editairement belordonn\'e}\index{classe!h\'er\'editairement belordonn\'ee} si $\mathscr{C}.\mathcal{A}$
est belordonn\'e pour tout belordre  $\mathcal{A}$.

\vspace{2mm}

Comme exemple, la classe des cha\^{i}nes finies est h\'er\'editairement belordonn\'ee. En effet, si $Ch$ est la classe\index{classe!des cha\^{i}nes finies} des cha\^{i}nes finies, $Ch.\mathcal A$ s'identifie \`a l'ensemble $\mathcal A^*$ des mots finis sur l'alphabet $\mathcal A$ muni de l'ordre de Higman. %(voir l'assertion $4.$ des Propri\'et\'es \ref{propr:belordre}).
Le fait que $Ch$ est h\'er\'editairement belordonn\'ee est le ''fameux'' r\'esultat  d\^u \`a Higman \cite{higman52} (voir Th\'eor\`eme \ref{theo:higman}).% dans la section \ref{sec:belordre}).

\vspace{2mm}

Voici quelques propri\'et\'es des classes h\'er\'editairement belordonn\'ees:

\begin{propertie}
Toute classe h\'er\'editairement belordonn\'ee  est belordonn\'ee.
\end{propertie}

\begin{proof}
Il suffit de prendre $\mathcal A=\{1\}$.
\end{proof}

\begin{propertie}\label{pro:R wqo}
Si $\mathscr{C}=\{\mathcal R\}$ alors $\mathscr C$ est h\'er\'editairement belordonn\'ee si $\mathcal R$ est finie.
\end{propertie}

\begin{proof}
Soit $\mathcal A$ un belordre, il s'agit de montrer que $\mathscr{C}.\mathcal{A}$ est belordonn\'e. %Comme $\mathscr C$ est finie, nous avons:
%$\mathscr{C}.\mathcal{A}=\underset{\mathcal R\in\mathscr C}{\cup}\{\mathcal R\}.\mathcal A.$

%Nous avons $\mathscr{C}.\mathcal{A}=\mathcal R.\mathcal A.$
%Puisque une r\'eunion finie de belordres est un belordre, ceci revient \`a montrer que $\mathcal R.\mathcal A$ est belordonn\'e pour $\mathcal R\in\Omega_{\mu}.$
 $\mathcal R:=(E,(\rho_i)_{i\in I})$ o\`u $E=\{a_1,\cdots,a_k\}.$
Comme $\mathcal A$ est belordonn\'e alors $\mathcal A^k$ est belordonn\'e. Posons $f_i:~\mathcal A^k\rightarrow\mathcal A$, pour $1\leq i\leq k$, la $i$\`eme projection de $\mathcal A^k$ dans $\mathcal A$ et consid\'erons l'application
$$\begin{array}{lrcl}
g:&\mathcal A^k& \rightarrow& \mathcal R.\mathcal A\\
& p& \mapsto& (\mathcal R,\widetilde{f})
\end{array}
\text{ o\`u }\begin{array}{lrcl}
\widetilde{f}:&E & \rightarrow& \mathcal A\\
& a_i& \mapsto& f_i(p)
\end{array}$$

L'application $g$ est surjective, en effet pour $(\mathcal R,\widetilde{f})\in \mathcal R.\mathcal A$, il existe $p\in\mathcal A^k$ tel que $g(p)=(\mathcal R,\widetilde{f})$, il suffit de prendre $p=(\widetilde{f}(a_1),\cdots,\widetilde{f}(a_k)).$ L'application $g$ est croissante par construction. Comme une image surjective par une application croissante d'un belordre est un belordre, on d\'eduit que $\mathcal R.\mathcal A$ est belordonn\'e. Ce qui entra\^{\i}ne que $\mathscr C$ est h\'er\'editairement belordonn\'ee.
\end{proof}
\bigskip

Si $\mathcal R$ est infinie, ceci n'a pas forc\'ement lieu. Si $\mathcal R$ est une anticha\^{\i}ne infinie, pour que $\mathcal R.\mathcal A$ soit belordonn\'e il faut et il suffit que $\underline{I}(\mathcal A)$ soit un belordre. %consid\'erer les sections initiales principales de $\mathcal A$

\vspace{1mm}

Une autre propri\'et\'e qui d\'ecoule du fait que toute union finie de belordres est un belordre %(voir section \ref{sec:belordre})
est % inclure la propriete dans la section belordre)
\begin{propertie}\label{pro:unionhere-wqo}
Toute union finie de classes h\'er\'editairement belordonn\'ees  est  h\'er\'editairement  belordonn\'e.
\end{propertie}

Et donc:

\begin{propertie}
Si $\mathscr C\subseteq\Omega_{\mu}$ est finie alors $\mathscr C$ est h\'er\'editairement belordonn\'ee.
\end{propertie}

\begin{proof}
Soit $\mathcal A$ un belordre, il s'agit de montrer que $\mathscr{C}.\mathcal{A}$ est belordonn\'e. Comme $\mathscr C$ est finie, nous avons: $\mathscr{C}.\mathcal{A}=\underset{\mathcal R\in\mathscr C}{\cup}\{\mathcal R\}.\mathcal A.$
Puisque une r\'eunion finie de belordres est un belordre, ceci revient \`a montrer que $\mathcal R.\mathcal A$ est belordonn\'e pour $\mathcal R\in\Omega_{\mu}.$ Ce qui est vrai d'apr\`es la Propri\'et\'e \ref{pro:R wqo}.
\end{proof}
% (inclure les references pour cette partie)
\medskip

Une question, rest\'ee ouverte depuis les ann\'ees $70$, est de savoir si cette condition est vraie lorsqu'elle est v\'erifi\'ee pour une anticha\^{i}ne \`a deux \'el\'ements $\mathcal{A}$.

\vspace{2mm}

 Rappelons \'egalement les r\'esultats suivants dus \`a Pouzet \cite{pouzet-belordre72, pouzet06}:

\begin{theorem}\label{theo:pouzet-wqo}
Si une  sous-classe $\mathscr C$ de $\Omega_\mu$ (avec $\mu$ fini) est h\'er\'editairement belordonn\'ee, alors $\downarrow\mathscr C$, la plus petite sous-classe h\'er\'editaire de $\Omega_\mu$ contenant $\mathscr C$, est h\'er\'editairement belordonn\'ee.
\end{theorem}

%\noindent On rappelle aussi le fait suivant d\^u \`a Pouzet 1972, cf\cite{pouzet72}:
 \begin{theorem}\label{theo:pouzet-borne}
 Si la signature $\mu$ est finie, une sous-classe de $\Omega_\mu$ qui est h\'er\'editaire et h\'er\'editairement belordonn\'ee a un nombre fini de bornes.
 \end{theorem}
%

%%%%%%%%%%%%%%%%%%%%%%%%%%%%%%%%%%%%%%%%%%%%%%%%%%%%%%%%%%%%%%%%%%%%%%%%%%%%%%%%%%%%%%%%%%%%%%%

\section{Profil et fonction g\'en\'eratrice}
\subsection{Profil de classe h\'er\'editaire}
Soit $\mathscr C$ une classe h\'er\'editaire. Le \emph{profil}\index{profil}\footnote{Dans \cite{B-B-S-S}, cette fonction est appel\'ee ``\emph{unlabelled speed}''. Il est \'egalement question de la fonction ``\emph{labelled speed}'' qui compte, pour tout entier $n$, le nombre de toutes les structures d\'efinies sur $n$ \'el\'ements.} de $\mathscr C$ est la fonction $\varphi _{\mathscr{C}}:\mathbb{N\longrightarrow N}$ qui compte, pour tout entier $n$, le nombre de structures d\'efinies sur $n$ \'el\'ements, appartenant \`a $\mathscr{C}$, les structures isomorphes \'etant identifi\'ees. Si on pose $\textit T_{\mathscr C}$ l'ensemble des types d'isomorphie des structures de $\mathscr C$ et, pour tout entier $n$,
$\mathscr C_n=\{\mathcal R \in \mathscr{C}: \vert dom(\mathcal R)\vert=n\}$ alors, $\varphi _{\mathscr{C}}(n):=\vert \textit T_{\mathscr C_n}\vert$. Notons que $\varphi _{\mathscr{C}}(0)=1$.\\

 La notion de profil est une sp\'ecialisation aux classes h\'er\'editaires d'une notion de base de la combinatoire \'enum\'erative le \emph{comptage} ou \emph{l'\'enum\'eration}. En fait, un probl\`eme de base de la combinatoire \'enum\'erative est le d\'enombrement des \'el\'ements d'un ensemble fini, c'est \`a dire, que l'on dispose d'un ensemble $S$ d'objets, sur lequel est d\'efini un param\`etre $p:S\rightarrow\mathbb N$, appel\'e \emph{la taille}, et l'on sait que le nombre $a_n$ d'\'el\'ements de $S$ de taille $n$ est fini. La fonction de comptage $f$ est alors donn\'ee par $f(n)=a_n$. Les r\'esultats g\'en\'eraux sur ces fonctions de comptage sont de deux types:

\vspace{1mm}

  - Des r\'esultats exacts dans lesquels les fonctions $f$ sont d\'efinies explicitement, par exemple, des fonctions polynomiales, ou des fonctions d\'efinies par des relations de r\'ecurrence.

\vspace{1mm}

   - Des r\'esultats asymptotiques, pour lesquels $f$ est d\'efinie par des \'equivalences asymptotiques ou des in\'egalit\'es asymptotiques, par exemple des fonctions avec une croissance au plus exponentielle ou au moins polynomiales.

\vspace{1mm}

   Pour plus de d\'etails sur les probl\`emes d'\'enum\'eration, consulter \cite{Flajolet, stanley, Van-Wil}.\\

 Pour \'enum\'erer des objets, il est parfois n\'ecessaire d'utiliser des outils tels que le principe d'inclusion-exclusion ou les codages.
 Le \emph{principe d'inclusion-exclusion}\index{principe d'inclusion-exclusion} est l'un des outils fondamentaux de la combinatoire \'enum\'erative. Il s'agit d'un type de raisonnement familier dans les math\'ematiques \'el\'ementaires.
 Son principe, pour effectuer un d\'ecompte exact de la cardinalit\'e d'un ensemble $S$, consiste \`a commencer par un ensemble plus large, puis apporter des corrections en \'eliminant les \'el\'ements "ind\'esirables".
 L'id\'ee est tr\`es simple, c'est essentiellement une g\'en\'eralisation de l'observation, toute simple, suivante:  si $A$, $B$ et $C$ sont trois ensembles alors  $$\vert A\cup B\vert=\vert A\vert+\vert B\vert-\vert A\cap B\vert$$ et
 $$\vert A\cup B\cup C\vert=\vert A\vert+\vert B\vert+\vert C\vert-\vert A\cap B\vert-\vert A\cap C\vert-\vert B\cap C\vert+\vert A\cap B\cap C\vert.$$ Notons que l'on commence par \emph{inclure} puis on exclut les termes que nous avons compt\'e plus d'une fois, en l'occurence ici,  $\vert A\cap B\vert$, $\vert A\cap C\vert$ et $\vert B\cap C\vert$, mais cette exclusion peut \'eliminer compl\'etement certains \'el\'ements, en l'occurence ici,  $\vert A\cap B\cap C\vert$, il faudra alors les inclure encore une fois (voir \cite{Anderson, Charal, Van-Wil} pour plus de d\'etails sur ce principe). Ce principe, d'apr\`es \cite{Van-Wil}, est ancien, il remonte au moins au $19^{eme}$ si\`ecle, il est apparu dans le papier de Da Silva (1854) et plus tard dans le papier de Sylvester (1883) (voir \cite{Van-Wil}).\\

 Le principe du \emph{codage}\index{codage} est d'\'etablir une bijection (respectant la taille) entre l'ensemble $S$ des objets que l'on veut \'enum\'erer et un ensemble $S'$ form\'e d'objets de structure plus simple, exemple (mots ou suites, arbres,...), des objets qu'une technique standard permettrait d'\'enum\'erer plus facilement (voir \cite{Flajolet} pour plus de d\'etails, des exemples  et pour d'autres r\'ef\'erences).

 \subsection{Fonction g\'en\'eratrice}
  Tr\`es souvent, l'\'evaluation de la fonction d'\'enum\'eration se fait au moyen de la fonction g\'en\'eratrice. Si $\mathscr{C}$ est une classe h\'er\'editaire, la \textit{fonction g\'en\'eratrice}\index{fonction!g\'en\'eratrice} de $\mathscr{C}$ est $$\mathcal {H}_{\mathscr{C}}(X):=\sum_{n\geq 0}\varphi _{\mathscr{C}}(n)X^{n}.$$
  La fonction g\'en\'eratrice donn\'ee ici est ordinaire. Il existe d'autres types de fonctions g\'en\'eratrices. Par exemple, la fonction g\'en\'eratrice exponentielle de $f(n)$ donn\'ee par: $\underset{n\geq 0}\sum f(n)\dfrac{x^n}{n!}$. Les fonctions g\'en\'eratrices ordinaires et exponentielles sont les plus couramment utilis\'ees.

  \vspace{2mm}

En fait, il est plus commode de d\'ecrire la suite de nombre $\varphi(n)$ par sa fonction g\'en\'eratrice qui est une s\'erie formelle, car certaines op\'erations sur les objets combinatoires correspondent \`a des op\'erations simples (somme, produit, d\'erivation,...) sur les s\'eries g\'en\'eratrices associ\'ees.\\

Comme pour la fonction profil, la fonction g\'en\'eratrice peut, dans certains cas, \^etre d\'efinie explicitement. Dans d'autres cas, il n'est pas ais\'e de la d\'efinir. Il est alors courant d'en donner une description g\'en\'erale en donnant la classe \`a laquelle elle appartient. On rappelle ici deux familles de s\'eries dont l'une est incluse dans l'autre. \\

Consid\'erons une s\'erie formelle\index{serie@s\'erie formelle} $F(x):=\underset{n\geq 0}{\sum }a_{n}.x^{n}$ \`a coefficients dans un corps $\mathbb K.$ La d\'efinition est donn\'ee ici de mani\`ere g\'en\'erale, toutes les s\'eries consid\'er\'ees dans ce document sont à coefficients entiers donc $\mathbb K=\mathbb Q,$ (voir preuve du Lemme \ref{lem:wqo}).

    \subsubsection{Les s\'eries rationnelles}
    La s\'erie $F(x)$ est \emph{rationnelle},\index{serie@s\'erie!rationnelle} dans $\mathbb K$, s'il existe deux polyn\^omes $P(x)$ et $Q(x)$, à coefficients dans $\mathbb K$, avec $Q(0)\neq 0$ tels que $$F(x)=\dfrac{P(x)}{Q(x)}.$$
    De fa\c{c}on \'equivalente, si les coefficients $a_n$ satisfont une r\'ecurrence lin\'eaire \`a coefficients constants: il existe un entier $k$, des nombres $\alpha_1,\ldots,\alpha_k$ dans $\mathbb K$ et un entier $n_0$ tels que, pour $n\geq n_0$,
    $$a_n+\alpha_1 a_{n-1}+\ldots+\alpha_k a_{n-k}=0.$$
    \subsubsection{Les s\'eries alg\'ebriques:} La s\'erie $F(x)$ est \emph{alg\'ebrique}\index{serie@s\'erie!alg\'ebrique} s'il existe un polyn\^ome $Q(x,y)$ dans $\mathbb{K}[x,y]$ tel que $Q(x,F(x))=0$. Si le polyn\^ome $Q$ est de degr\'e $1$ en $y$, la s\'erie est alors rationnelle, donc les s\'eries rationnelles sont des s\'eries alg\'ebriques.

    \vspace{2mm}

Les fonctions alg\'ebriques et les s\'eries alg\'ebriques sont, tout simplement, d\'efinies comme les solutions des \'equations ou des syst\`emes d'\'equations polynomiales.
\vspace{1mm}

 Nous rappelons que la somme et le produit de deux fonctions ou s\'eries alg\'ebriques sont alg\'ebriques, il en est de m\^eme pour l'inverse d'une fonction ou s\'erie alg\'ebrique non nulle en z\'ero ($f$ alg\'ebrique et $f(0)\neq 0$ alors $\dfrac{1}{f}$ est alg\'ebrique).
\vspace{1mm}

Si $\mathscr C$ est une classe de structures combinatoires finies, nous disons que $\mathscr C$ est alg\'ebrique (resp. rationnelle) si sa s\'erie g\'en\'eratrice est alg\'ebrique (resp. rationnelle).\\

%\vspace{2mm}

    En g\'en\'eral, dire qu'une s\'erie $f$ est alg\'ebrique ne permet pas de la d\'efinir explicitement, mais fournit des informations sur le comportement de ses coefficients. Nous savons, par exemple, que les coefficients d'une s\'erie alg\'ebrique croissent au plus de mani\`ere exponentielle voir \cite{Flajolet, stanley2} et qu'ils v\'erifient une r\'ecurrence \`a coefficients polynomiaux, en d'autres termes, les coefficients $a_n$ satisfont \cite{klazar, stanley2}:$$p_k(n)a_{n+k}+p_{k-1}(n)a_{n+k-1}+\ldots+p_0(n)a_n=0,$$ avec les coefficients $p_i\in\mathbb K[x]$ non tous nuls.

 \vspace{2mm}

     Nous rappelons \'egalement que la suite $a_n,~n\in\mathbb N$ est \emph{quasi-polynomiale}\index{serie@s\'erie!quasi-polynomiale} si on a, pour tout $n$, $$a_n=b_k(n)n^k+\ldots+b_1(n)n+b_0(n),$$ o\`u les $b_i,~0\leq i\leq k$ sont des fonctions p\'eriodiques. De mani\`ere \'equivalente, si $a_n$ est le coefficient de $x^n$ dans le d\'eveloppement en s\'erie de l'expression, $$\dfrac{p(x)}{(1-x)(1-x^2)\ldots(1-x^l)}$$ pour un $l$ dans $\mathbb N$ et un polyn\^ome $p\in \mathbb K[x]$.

\vspace{2mm}

    L'\'etude des s\'eries n'\'etant pas l'objet de ce document, nous ne donnerons pas plus de d\'etails, pour le lecteur int\'eress\'e voir par exemple \cite{Bous-melou, Charal, Flajolet, stanley2, Van-Wil}. % Des r\'esultats g\'en\'eraux existent sur le comportement asymptotique des coefficients d'une s\'erie alg\'ebrique

\subsection{Comportement des profils des classes h\'er\'editaires}

\subsubsection{Ordre de croissance du profil}\index{profil!croissance du -}
Etant donn\'ees deux fonctions $\varphi$ et $\psi$ de variable et valeurs enti\`eres, la fonction $\psi$ \emph{cro\^{i}t au moins aussi vite que} $\varphi$ (que la fonction $\psi$ soit croissante ou non) lorsqu'il existe un entier $a>0$ tel que $\varphi(n)\leq a\psi(n)$ pour $n$ assez grand.

\vspace{2mm}

Les fonctions $\varphi$ et $\psi$ ont \emph{m\^eme ordre de croissance} lorsque chacune cro\^{i}t au moins aussi vite que l'autre, par exemple si une classe est finie, son profil a le m\^eme ordre de croissance que la fonction nulle.

\vspace{1mm}

Si une fonction $\varphi$ a le m\^eme ordre de croissance que la fonction $\psi=n^k$ pour un entier non n\'egatif $k$, ou encore s'il existe deux r\'eels positifs $a$ et $b$ tels que $an^k\leq \varphi\leq bn^k$ pour $n$ suffisamment large, alors $\varphi$ a une \emph{croissance polynomiale\index{profil!\`a croissance polynomiale} de puissance} $k$. Les cas particuliers de la croissance polynomiale sont:

\vspace{1mm}

Si $k=0$ alors $\varphi$ est \emph{born\'ee} et ne s'annule pas (exemple le profil de la classe form\'ee de cha\^{i}nes finies).\index{profil!born\'e}

\vspace{1mm}

Si $k=1$ alors $\varphi$ a une \emph{croissance lin\'eaire}.\index{profil!\`a croissance lin\'eaire}

\vspace{1mm}

Si $k=2$ alors $\varphi$ a une \emph{croissance quadratique} ou \emph{parabolique}.\index{profil!\`a croissance quadratique}

\vspace{2mm}

La fonction $\varphi$ cro\^{i}t plus vite que tout polyn\^ome si et seulement si $\underset{n\rightarrow\infty}{lim}\frac{\varphi(n)}{n^k}=+\infty$ pour tout entier non n\'egatif $k$.

\vspace{1mm}

Si une fonction $\varphi$ a m\^eme ordre de croissance que la fonction $c^n$ pour une constante $c>0$ alors $\varphi$ a une \emph{croissance exponentielle}.\index{profil!\`a croissance exponentielle}

\subsubsection{Saut dans la croissance des profils}

Durant ces derni\`eres ann\'ees, il y a eu beaucoup de r\'esultats d'\'enum\'eration sur des classes h\'er\'editaires de structures concr\`etes telles que les graphes, les graphes ordonn\'es, les tournois et les permutations en connection avec la conjecture de Stanley-Wilf, sur laquelle nous reviendrons dans le chapitre \ref{chap:permut-bichaine}.

Tous les r\'esultats obtenus montrent, essentiellement, qu'il existe des \emph{sauts} dans les croissances des profils de ces classes h\'er\'editaires. Un \emph{saut} dans la croissance du profil d'une classe d'objets combinatoires $\mathscr C$ signifie que la croissance de $\varphi_{\mathscr C}$ ne se fait pas de mani\`ere continue, il y a des ``trous'', le profil peut, par exemple, passer d'une croissance polynomiale \`a une croissance exponentielle, ce qui veut dire qu'il n'existe pas de classe de tels objets ayant un profil v\'erifiant $n^k<\varphi_{\mathscr C}<c^n$ pour $n$ suffisamment large et pour tout $k>0$ et $c>1$.

Par exemple, dans le cas des graphes, Balogh, Bollob\'{a}s, Saks et S\'{o}s (2009) \cite{B-B-S-S} ont montr\'e:
\begin{theorem}\label{thm:graph-balogh}
Si $\mathscr C$ est une classe h\'er\'editaire de graphes finis alors une seule des assertions suivantes se produit.
\begin{enumerate}
\item Pour tout $n$, $\varphi_{\mathscr C}(n)=cn^k+O(n^{k-1})$, pour un entier $k>0$ et $c\in \mathbb Q$, $c>0$.
\item Pour $n$ suffisamment large, $\varphi_{\mathscr C}(n)\geq \mathfrak p(n),$ o\`u $\mathfrak p(n)$ est la fonction partition d'entier.
\end{enumerate}
\end{theorem}
\noindent Nous rappelons, que pour un entier positif $n$, une partition de $n$ est la donn\'ee d'une suite finie d\'ecroissante $n_1\geq n_2\geq \cdots\geq n_l$ d'entiers positifs appel\'es \emph{parties} telle que $n_1+n_2+\cdots+n_l=n.$ La fonction partition $\mathfrak p$ compte, pour tout entier $n$, le nombre $\mathfrak p(n)$ de partitions de $n$. La fonction g\'en\'eratrice associ\'ee \`a $\mathfrak p(n)$ est $\underset{n\geq 1}\sum \mathfrak p(n)x^n:=\underset{n\geq 1}\prod\dfrac{1}{1-x^n}$.

La valeur asymptotique de ce nombre est donn\'ee par un ''fameux'' r\'esultat d'Hardy et Ramanujan (1918) qui est: $\mathfrak p(n)\sim\frac{1}{4n\sqrt{3}}exp\Big(\pi\sqrt{\frac{2n}{3}}\Big)$.\\

Dans \cite{B-B-M(06), B-B-S-S, B-B-M(07)} il est montr\'e  que dans les cas des tournois, des graphes ordonn\'es et des permutations le saut de croissance du profil se fait d'une croissance polynomiale \`a une croissance exponentielle, g\'en\'eralisant ainsi les r\'esultats de Kaiser et Klazar (2004) sur les permutations sur lesquelles nous reviendrons dans le chapitre \ref{chap:permut-bichaine}.\\

\noindent Le comportement des profils des classes h\'er\'editaires particuli\`eres, les  "\^ages" de Fra\"{i}ss\'e et le lien avec les classes belordonn\'ees  ont \'et\'e consid\'er\'es par Pouzet au d\'ebut des ann\'ees soixante-dix (voir\cite{pouzet-belordre72} et \cite{pouzet06} pour une synth\`ese). Le cas des graphes, tournois et autres structures combinatoires a \'et\'e  \'elucid\'e plus r\'ecemment. Pour plus de d\'etails consulter l'article de Klazar \cite{klazar}. %Les autres structures combinatoires ont \'et\'e largement \'etudi\'ees dans \cite{klazar,B-B-S-S,B-B-M(07),B-B-M(06)}.
Nous \'etudierons le ph\'enom\`ene du saut du profil sur les structures relationnelles ordonn\'ees dans le chapitre \ref{chap:monomorphe}.
%Les notions de profil et de fonctions g\'en\'eratrice sont des sp\'ecialisations aux classes h\'er\'editaires des notions de base sur l'\'enum\'eration. Les nombreux travaux sur l'\'enum\'eration des classes de permutations rentrent dans le cadre de l'\'enum\'eration des structures relationnelles finies. En effet, comme sugg\'er\'e par P. J. Cameron (2002), les permutations peuvent \^etre vues comme des structures relationnelles particuli\`eres: des bicha\^{i}nes (couple de deux ordres totaux d\'efinis sur un m\^eme ensemble). Nous introduirons ces notions ainsi que la relation entre les permutations et les bicha\^{i}nes dans le chapitre 2.

    \subsection{Profil d'une structure relationnelle}\label{subsec:profil}
Soit $\mathcal R$ une structure relationnelle de base $E$ et de signature $\mu$. Nous d\'esignons \'egalement par \emph{profil\index{profil!d'une structure} de $\mathcal R$} la fonction $\varphi _{\mathcal{R}}$ qui compte, pour tout entier $n$, le nombre $\varphi _{\mathcal{R}}(n)$ de types d'isomorphie des sous-structures\index{structure relationnelle!sous-structure} de $\mathcal R$ induites sur les $n$-parties de $E$. Cette notion est due \`a Pouzet 1972 (voir \cite{pouzet06} pour plus de r\'esultats).\\ Notons que $\varphi_{\mathcal{R}}(0)=1$ et si $E$ est un $p$-ensemble, o\`u $p$ est un entier positif, alors $\varphi_{\mathcal{R}}(p)=1$ et $\varphi_{\mathcal{R}}(q)=0$ pour tout $q>p$.

\begin{example} 1) Si $\mathcal{R}$ est une cha\^{i}ne et $E$ un $n$-ensemble alors $\varphi_{\mathcal{R}}(k)=1$ pour $k\leq n$.\\ 2) Si $\mathcal R$ est une relation unaire d\'efinie sur un ensemble infini $E$ avec $\mathcal R=F\subset E$ o\`u $F$ est une $p$-partie de $E$, $p$ entier, alors le profil de $\mathcal R$ cro\^{i}t de $\varphi_{\mathcal{R}}(0)=1$ \`a $\varphi_{\mathcal{R}}(p)=p+1$ et reste stationnaire.
\end{example}

Il est clair que cette fonction ne d\'epend que de l'\^age $\mathcal A(\mathcal R)$ de $\mathcal R$. Si la signature $\mu$ de $\mathcal R$ est finie, le nombre de structures relationnelles de signature $\mu$ d\'efinies sur un $n$-ensemble \'etant fini, la valeur de $\varphi_{\mathcal{R}}(p)$, pour tout entier $p$, est n\'ecessairement un entier (donc finie). Dans le cas ou $\mu$ n'est pas fini, $\varphi_{\mathcal{R}}(n)$ peut-\^etre infini, mais nous nous int\'eressons pas \`a ce cas ici.\\
Les r\'esultats suivants, dus \`a Pouzet (voir \cite{fraisse} et\cite{pouzet06}), nous renseignent sur l'\'evolution de la fonction profil:

\begin{propertie}
Soit $n<\vert E\vert$. Alors $$\varphi_{\mathcal{R}}(n)\leq (n+1)\varphi_{\mathcal{R}}(n+1).$$
\end{propertie}

%Qui peut-\^etre am\'elior\'ee par:
\begin{theorem}
 Si $\mathcal R$ est une structure relationnelle de base infinie,\index{structure relationnelle!de base infinie} alors son profil $\varphi_{\mathcal R}$ est non-d\'ecroissant.
\end{theorem}

De plus;
\begin{theorem}
Soient $p, q$ deux entiers et $\mathcal R$ une structure relationnelle de base $E$ de cardinalit\'e au moins $2p+q$. Alors $$\varphi_{\mathcal{R}}(p)\leq \varphi_{\mathcal{R}}(p+q).$$
\end{theorem}

Nous rappelons en particulier ce th\'eor\`eme d\^u \`a Pouzet \cite{pouzet.tr.1978} (voir aussi \cite{P-T-2013}):% qui donne une relation entre le profil \index{profil} d'une relation et les bornes\index{borne} de son \^age

\begin{theorem}\label{theo:pouzet-polyborne}
L'\^age d'une structure relationnelle (de signature finie) dont le profil est born\'e par un polyn\^ome poss\`ede un nombre fini de bornes.% with polynomially bounded profile can be defined by finitely many bounds.
\end{theorem}

\subsubsection{Exemples de fonction profil}\label{par:exemple}

Beaucoup de fonctions usuelles sont des profils de relations, voici quelques exemples (voir \cite{pouzet06}):

\vspace{2mm}

\begin{enumerate}
\item \emph{Le coefficient binomial} $\binom{n+k}{k}$.
 Soit  $\mathcal R:=(\mathbb Q, \leq, u_1,\dots,u_k)$ o\`u $\leq $ est l'ordre naturel sur l'ensemble $\mathbb Q$  des nombres rationels,
  $u_1,\dots,u_k$ sont des relations unaires qui divisent $\mathbb Q$ en $k+1$ intervalles. Alors
$\varphi_{\mathcal R}(n) = {n+k \choose k}$.

\vspace{1mm}

\item \emph{La fonction exponentielle} $n \hookrightarrow k^{n}$.  Soit $\mathcal R:=(\mathbb Q, \leq , u_1, \dots, u_k)$, o\`u
  $u_{1}, \dots, u_{k}$ d\'esignent toujours $k$ relations unaires, qui divisent, cette fois, $\mathbb Q$ en $k$ ``couleurs''
  de telle sorte qu'entre deux nombres rationnels toutes les couleurs apparaissent. Alors $\varphi_{\mathcal R}(n) = k^n$.

\vspace{1mm}

\item \emph{La factorielle} $n \hookrightarrow n!$. Soit $\mathcal R:= (\mathbb Q, \leq,\leq')$, o\`u $\leq'$ est un autre ordre lin\'eaire sur $\mathbb Q$ d\'efini de telle sorte que les restrictions finies induisent toutes les paires possibles d'ordres lin\'eaires sur un ensemble fini.
    %(eg take for $\leq'$ an order with the same type as the natural order on the set $\N$ of non-negative integers).
 Alors $\varphi_{\mathcal R}(n) = n!$

\vspace{1mm}

\item  {\it La fonction partition d'entier} qui compte le nombre $\mathfrak p(n)$ de partitions de l'entier $n$. Soit $\mathcal R:= (\mathbb N, \rho)$ le chemin infini d'entiers (le graphe de comparabilit\'e de la relation de cons\'ecutivit\'e sur $\mathbb N$) dont les ar\^etes sont les paires  $\{n,n+1\}$ pour $n\in\mathbb N$. Alors $\varphi_{\mathcal R}(n)=\mathfrak p(n)$. La d\'etermination de sa croissance asymptotique et les difficult\'es rencontr\'ees pour prouver que $\mathfrak p(n)\simeq\frac{1}{4n\sqrt{3}}e^{\pi\sqrt{\frac{2n}{3}}}$ (Hardy et Ramanujan,  1918) montrent les difficult\'es rencontr\'ees dans l'\'etude des profils en g\'en\'eral.
\end{enumerate}

\vspace{3mm}

Les graphes fournissent des exemples simples de profils  \cite{pouzet06} (voir \cite{P-T-2013} pour d'autres exemples). Soit $G$ un graphe et $\varphi_{G}$
%Si on voit un graphe $G$ comme relation binaire irr\'eflexive et sym\'etrique,
son profil. % $\varphi_{G}$ n'est que la fonction qui pour
%chaque entier $n$ compte le nombre $\varphi_{G}(n)$ de sous-graphes induits sur les parties $V(G)$ ayant $n$ sommets, ceux-ci \'etant compt\'es \`a l'isomorphiepr\`es.

\vspace{2mm}

\begin{enumerate}

\item  $\varphi_{G}$ est constant, \'egal \`a $1$, pour tout $n\leq\vert V(G)\vert$ si et seulement si
$G$ est une clique ou un ind\'ependant.

\vspace{1mm}

\item  Si $G$ est la somme directe d'un graphe fini et d'une clique infinie ou d'un
ind\'ependant infini alors son profil est born\'e. En fait, le profil
$\varphi_{G}$  d'un graphe $G$ est born\'e si et seulement si $G$ est \emph{presque constant} (c'est \`a dire qu'il existe un
sous-ensemble
fini $F_{G}$ de sommets tel que si  deux paires  de sommets ont la m\^eme intersection sur $F_{G}$  alors toutes deux sont des
ar\^etes ou aucune des deux ne l'est).

\vspace{1mm}

\item  Si $G$ est la somme directe d'une infinit\'e de cliques \`a deux \'el\'ements ou bien la
somme directe de deux cliques ayant chacune  une infinit\'e d'\'el\'ements alors $\varphi_{G}(n)=\lfloor \frac{n}{2}
\rfloor +1$ et sa fonction g\'en\'eratrice est $\mathcal H_{G}(x)=\frac{1}{(1-x)(1-x^2)}$.

\item Si $G$ est la somme directe d'une clique infinie et d'un ind\'ependant infini alors $\varphi_{G}(n)=n$ pour $n\geq 1$ et $\mathcal H_{G}(x)=1+\frac{x}{(1-x)^2}=\frac{1-x+x^2}{(1-x)^2}=\frac{1+x^3}{(1-x)(1-x^2)}$.

\vspace{1mm}

\item  Si $G$ est la somme directe d'une infinit\'e de cliques \`a $k$ \'el\'ements ou bien la somme directe de $k$ cliques infinies alors %voir flajolet p43
     $\varphi_{G}(n)=\mathfrak p_{k}(n) \simeq
\frac{n^{k-1}} {(k-1)!k!}$ et $\mathcal H_{G}(x)=\frac{1}{(1-x)\dots(1-x^k)}$.

$\mathfrak p_k(n)$ est le nombre de partitions de $n$ en au plus $k$ parties, ce nombre est \'egal au nombre de partitions de $n$ en parties de tailles au plus $k$, il est donn\'e par $\mathfrak p_k(n)\simeq\frac{n^{k-1}}{(k-1)!k!}$.

\vspace{1mm}

\item  Si $G$ est la somme directe d'une infinit\'e de cliques infinies ou bien un chemin infini alors $\varphi_{G}(n)=\mathfrak p(n)$ o\`u $\mathfrak p$ est la
fonction partition d'entiers.
\end{enumerate}

\clearemptydoublepage

%$$$$$$$$$$$$$$$$$$$$$$$$$$$$$$$$$$$$$$$$$$$$$$$$$$$$$$$$$$$$$$$$$$$$$$$$$$$
\part{Structures ind\'ecomposables et profil alg\'ebrique}\label{part:profil algebrique}
%$$$$$$$$$$$$$$$$$$$$$$$$$$$$$$$$$$$$$$$$$$$$$$$$$$$$$$$$$$$$$$$$$$$$$$$$$$$

\chapter{Permutations et bichaines}\label{chap:permut-bichaine}

Ce chapitre comprend deux %cinq
sections. La premi\`ere est consacr\'ee aux permutations. Celles-ci ayant fait l'objet de nombreuses recherches ces derni\`eres ann\'ees, beaucoup de r\'esultats existent, certains sont des red\'ecouvertes de r\'esultats plus anciens sur d'autres structures, d'autres sont r\'ecents. Nous ne donnerons ici qu'un bref aper\c{c}u des r\'esultats qui sont directement li\'es \`a notre sujet. Pour plus de d\'etails, consulter les diff\'erentes r\'ef\'erences qui traitent du domaine cit\'ees dans la bibliographie, en particulier l'excellent ouvrage de \emph{S. Kitaev} \cite{Kit}, ou celui de \emph{M. B\'{o}na} \cite{bona}, qui reprennent une bonne partie des r\'esultats sur les permutations, parus ces derni\`eres ann\'ees. La deuxi\`eme section est consacr\'ee aux bicha\^{i}nes, des structures relationnelles form\'ees de deux ordres totaux d\'efinis sur un m\^eme ensemble. Suivant \emph{P.J. Cameron}, nous explicitons le lien existant entre ces derni\`eres et les permutations, ce qui justifie que l'on consid\`ere les permutations comme des structures relationnelles.

%%%%%%  Section 3  %%%%%%%%%%%%%%%%%%%%%%%%%%%%%%%%%%%%%%%%%%%%%%%%%%%%%%%%%%%%%%%%%%%%%%%%%%%%

%**************************************************************

%************************************************section********************************
\section{Permutations}\label{sec:permutation}

\begin{figure}[!hbp]
\centering
\input{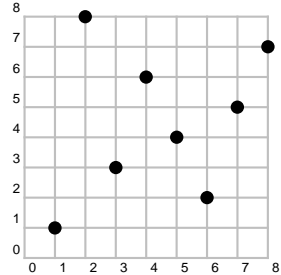}
\caption{\label{permutation-representation}Repr\'esentation graphique de la permutation $~\sigma:=18364257$.}
\end{figure}

Pour un entier $n\geq 1$, nous posons $[n]:=\{1,\dots, n\}$,  $[0]=\varnothing$. Une \emph{permutation}\index{permutation} $\sigma$ de $[n]$, $(n\geq 1)$, est une bijection de $[n]$ dans lui m\^eme, pour $n=0$ il n'y a qu'une seule permutation, la permutation vide ou nulle\footnote{Bien que la permutation nulle ne soit pas consid\'er\'ee dans les diff\'erents articles sur les permutations, on a pr\'ef\'er\'e la d\'efinir ici pour rester en phase avec la d\'efinition g\'en\'erale des structures.}. Si $\sigma$ est une permutation de $[n]$ alors $\sigma$ est de \emph{longueur}\index{permutation!longueur d'une -} $n$ et cette longueur est not\'ee $l(\sigma)$ ou $\vert\sigma\vert$. Nous notons $\sigma^{-1}$ la permutation inverse de $\sigma$. Une permutation $\sigma$ de $[n]$ peut \^etre repr\'esent\'ee de diff\'erentes mani\`eres:
 \begin{enumerate}
\item De mani\`ere lin\'eaire par la suite de ses valeurs, c'est \`a dire, $\sigma:=\sigma_1\sigma_2\dots \sigma_n,~$ o\`u $\sigma_i=\sigma (i)$, exemple $~\sigma:=18364257$.
\item Par un tableau \`a deux lignes $~\sigma:=\left(
\begin{array}{cccc}
1 & 2 & \cdots  & n \\
\sigma (1) & \sigma (2) & \cdots  & \sigma (n)%
\end{array}%
\right), $
exemple\\ $~\sigma:=\left(
\begin{array}{cccccccc}
1 & 2 & 3 & 4 & 5 & 6 & 7& 8\\
1&8&3&6&4&2&5&7%
\end{array}
\right),$
\item Par un produit de cycles, exemple $~\sigma:=(1)(287546)(3)$, les cycles de longueur $1$ n'\'etant souvent pas repr\'esent\'es.
\item Graphiquement comme le montre la \figurename~\ref{permutation-representation}.
\end{enumerate}

%*******************figure image1***********

Nous adopterons, dans toute la suite, la repr\'esentation lin\'eaire et, au besoin, la repr\'esentation graphique.\\

Il est bien connu qu'à toute permutation $\pi$ est associ\'e un graphe $G_{\pi}$, appel\'e \emph{graphe de permutations}.\index{graphe!de permutations} % et un ordre $P_{\pi}$, appel\'e \emph{ordre de permutation}.
 Si $\pi$ est une permutation de $[n]$ alors $G_{\pi}$ est le graphe ayant pour sommets l'ensemble $[n]$ et pour ar\^etes les paires de sommets  $\{i,j\}$ telles que l'ordre de $i$ et $j$ est invers\'e par $\pi$, en d'autres termes, $$\{i,j\}\in E(G_{\pi})\Leftrightarrow (i-j)(\pi_i^{-1}-\pi_j^{-1})<0.$$
Le graphe $G_{\pi}$ est un graphe d'incomparabilit\'e, les ordres qui lui sont associ\'es sont appel\'es \emph{ordres de permutations}.\index{ordre!de permutations}

\begin{figure}[t]
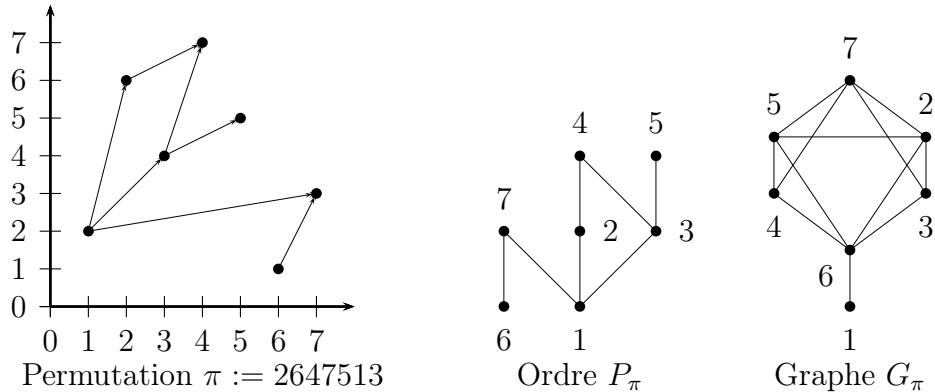

\centering
\input{figure2}\qquad\input{figure2a}\quad\input{figure2b}
\caption{\label{permut-ordre}Graphe $G_{\pi}$ et diagramme de Hasse de l'ordre $P_{\pi}$ pour $\pi:=2647513$.}
\end{figure}

Nous pouvons associer à $\pi$, de mani\`ere naturelle, un ordre $P_{\pi}$ d\'efini sur $[n]$ de la mani\`ere suivante. Consid\'erons la repr\'esentation graphique de $\pi$. Ordonnons "naturellement" l'ensemble $\{(i,\pi_i):i\in[n]\}$, c'est à dire $(i,\pi_i)\leq(j,\pi_j)$ si $i\leq j$ et $\pi_i\leq\pi_j$. Identifions chaque $i$ au couple $(i,\pi_i)$, on a alors un ordre sur $[n]$. Le graphe d'incomparabilit\'e de cet ordre est le graphe $G_{\pi}$. Les ordres $P_{\pi}$ et $P_{\pi^{-1}}$ sont isomorphes.\\%on n'a pas la m\^eme num\'erotation des sommets
\textbf{Exemple.} Consid\'erons la permutation $\pi=2647513$ de longueur 7. Sa repr\'esentation graphique, le diagramme de Hasse de l'ordre $P_{\pi}$ et le graphe $G_{\pi}$ sont donn\'es dans la \figurename~\ref{permut-ordre}.\\

%*******************figure2***********

%Pour plus de d\'etails sur les graphes et les ordres de permutations voir \cite{golumbic}

\subsection{Ordre sur les permutations}
 D\'esignons par $\mathfrak S_n$ l'ensemble des permutations de $[n]$ et posons $\mathfrak S:=\underset{n\in \mathbb  N}\bigcup\mathfrak S_n.$ Un ordre a \'et\'e d\'efini sur $\mathfrak S$ de la mani\`ere suivante.\\ Nous disons que deux suites finies de m\^eme longueur $n\geq 1$, $\alpha=\alpha_1\ldots\alpha_n$ et $\beta=\beta_1\ldots\beta_n$ ont la \emph{m\^eme forme} si, pour tous $i, j$, nous avons $\alpha_i<\alpha_j$ si et seulement si $\beta_i<\beta_j$. Ainsi, toute suite finie de nombres r\'eels distincts a la m\^eme forme qu'une unique permutation ayant pour longueur le nombre d'\'el\'ements de la suite.\\ Nous disons que la permutation $\pi$ de $[n]$ \emph{contient} la permutation $\sigma$ de $[k],$ et nous notons $\sigma \leq \pi,$ si $\pi$ a une sous-suite de longueur $k$ qui a la m\^eme forme que $\sigma$, cette sous-suite est appel\'ee \emph{copie de $\sigma$}\index{permutation!copie d'une -}. Plus pr\'ecis\'ement,
 $\sigma \leq \pi$ s'il existe des entiers $1\leq x_1 <\dots<x_k\leq n$ tels que pour $1\leq i,j\leq k,$
\begin{equation}
\sigma_i< \sigma_j~~\text{si et seulement si}~~ \pi_{x_i}< \pi_{x_j}.\label{perm-ordre}
\end{equation}

\begin{figure}[t]
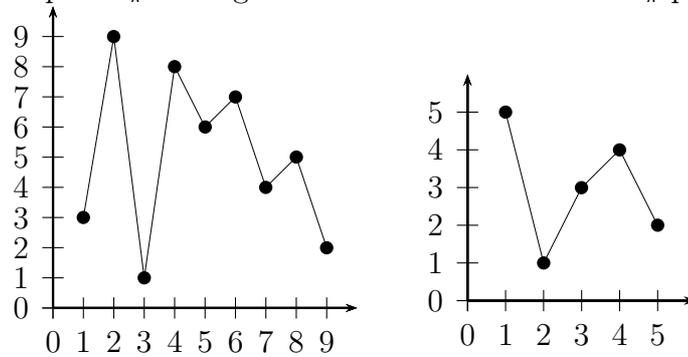

\centering
\input{figure1}\quad \input{figure1a}
\caption{\label{permutation-forme}Forme d'une permutation ($\pi:=391867452$ contient $\sigma:=51342$).}
\end{figure}

 En d'autres termes, $\sigma \leq \pi$ s'il existe des entiers $1\leq x_1 <\dots<x_k\leq n$ tels que la normalisation de $\pi_{x_1}\pi_{x_2}\cdots\pi_{x_k}$ sur $[k]$ donne $\sigma$.
Par exemple $\pi:=391867452$ contient $\sigma:=51342,$ comme nous pouvons le voir en consid\'erant la sous-suite $91672$ ($=\pi(2),\pi(3),\pi(5),\pi(6),\pi(9)$) (voir \figurename~\ref{permutation-forme}). Si $\sigma$ n'est pas contenue dans $\pi$, fait not\'e par $\sigma\nleq\pi$, nous disons que $\pi$ \emph{\'evite} $\sigma.$
L'une des caract\'eristiques de cet ordre est l'implication suivante: $~\sigma < \pi \Rightarrow l(\sigma)<l(\pi).$\\
%*******************figure figure1***********

Muni de cet ordre, l'ensemble ordonn\'e $\mathfrak S$ est bien fond\'e. En effet, $\mathfrak S$ est ordonn\'ee par niveaux de telle sorte que $\sigma\in\mathfrak S$ a le niveau $n$ si et seulement si $l(\sigma)=n$, donc il suffit de prendre, pour tout sous-ensemble $A$ de $\mathfrak S$, un \'el\'ement de longueur minimum. Par contre,  $\mathfrak S$ n'est pas belordonn\'e, en effet nous savons, depuis longtemps, que $\mathfrak S$ contient des anticha\^{i}nes infinies (voir \cite{Pratt, Tarjan} ou plus r\'ecemment \cite{Spi-Bona}).\\

 Un sous-ensemble $\mathscr C$ de $\mathfrak S$ est une classe \emph{h\'er\'editaire} si $\sigma <\pi \in\mathscr C$ implique $\sigma \in \mathscr C$. %Sa fonction \'enum\'eratrice, que l'on appelle
  Son profil est $\varphi _{\mathscr{C}}(n):=\vert \mathscr C\cap \mathfrak S_n\vert$, c'est \`a dire, la fonction qui compte, pour tout entier $n$, le nombre de permutations de $\mathscr C$ de niveau ou de longueur $n$. Par exemple, la fonction
 $\varphi_{\mathscr{C}}(n)=0$ pour $n\in \mathbb N$  est le profil de la classe vide $\mathscr C=\varnothing$, alors que la classe $\mathfrak S$ a pour profil la fonction\footnote{Cette formule serait connue depuis, au moins, cinquante si\`ecles \cite{Flajolet}.}
 $\varphi_{\mathfrak S}(n)=n!.$\\

 Mis \`a part ces deux exemples triviaux,  quelle est la valeur de  $\varphi_{\mathscr{C}}(n)$ si $\varnothing\neq\mathscr C \neq \mathfrak S ?$\\
 La conjecture de  \emph{Stanley-Wilf}\index{Stanley-Wilf} assure que cette valeur est born\'ee par une exponentielle. Cette conjecture est la suivante:\\

 \textbf{Conjecture de Stanley-Wilf.} \emph{Pour toute permutation $\sigma$, il existe une constante $c>0$ telle que $\varphi_{\text{\it Forb}(\sigma)}(n)<c^n$ pour tout entier $n\geq 1.$}\\

  o\`u $\text{\it Forb}(\sigma):=\{\pi\in\mathfrak S:\sigma\nleq\pi\}$ pour tout $\sigma \in \mathfrak S$.\\

 Il est tr\`es difficile de trouver une r\'ef\'erence exacte de la conjecture de \emph{Stanley-Wilf}. D'apr\`es R. Arratia \cite{Arratia}, elle aurait \'et\'e formul\'ee lors d'une communication orale en 1990 et que d'apr\`es une communication de Wilf (1999) cette conjecture, \`a l'origine, avait la formulation suivante:\\

  \textbf{Conjecture de Stanley-Wilf (autre version).}  \emph{Pour toute permutation $\sigma$, la limite de  $\varphi_{Forb(\sigma)}(n)^{1/n}$ existe et est finie.}\\

  Arratia, dans son papier \cite{Arratia}, a montr\'e que ces deux versions de la conjecture sont \'equivalentes.\\

  La conjecture de \emph{Stanley-Wilf} a suscit\'e beaucoup de travaux sur, notamment, la fonction profil et son comportement asymptotique, voir par exemple les articles \cite{A-A, A-A-K, Mar-Tar, Mur-Vat} pour les classes de permutations, \cite{ B-B-M(06),B-B-M(07), B-B-S-S, Bou-Pouz, mont-pou} pour les autres structures pour ne citer que ceux-l\`a, il y a \'egalement les articles de Klazar \cite{klazar, klazar08} qui r\'esument les r\'esultats sur l'\'enum\'eration parus ces derni\`eres ann\'ees.\\

  La validit\'e de la conjecture de Stanley-Wilf a \'et\'e confirm\'ee par \emph{Marcus et Tard\"{o}s}\footnote{Pour certaines classes particuli\`eres de permutations, la conjecture a \'et\'e \'etablie avant $2004$ voir \cite{Alo-fried, bona1, bona2, bona3}.} en 2004. Dans la m\^eme p\'eriode, \emph{Kaiser et Klazar} (2003),  ont reformul\'e la conjecture de \emph{Stanley-Wilf} sous la forme suivante:\\

  \textbf{Conjecture de Stanley-Wilf reformul\'ee.} \emph{Si $\mathscr C$ est une sous-classe h\'er\'editaire propre de permutations, alors il existe une constante $c>0$ telle que $\varphi_{\mathscr{C}}(n)< c^n$ pour tout entier $n\geq 1$.}\\

 En effet, si $\mathscr C$ est une classe h\'er\'editaire de permutations telle que $\sigma\notin\mathscr C$ alors $\varphi_{\mathscr C}(n)\leq \varphi_{\text{\it Forb}(\sigma)}(n)$ pour tout $n\geq 1$. Nous rappelons aussi que toute classe h\'er\'editaire de permutations $\mathscr C$ a la forme $\text{\it Forb}(F)$ pour un sous-ensemble $F\subset\mathfrak S$ (voir section \ref{sec:clas-heredi-age-borne}).\\

  Ainsi reformul\'ee, cette conjecture assure, qu'except\'e pour la classe $\mathfrak S$, il n'y a pas d'autre profil superexponentiel. A noter que
  $\varphi_{\mathfrak S}$ cro\^{i}t plus vite que toute fonction exponentielle, en effet, pour toute constante $c>0$, nous avons $\dfrac{c^{n+1}}{c^n}=c,~\forall n$ alors que $\dfrac{(n+1)!}{n!}=n+1.$\\

\emph{Kaiser et Klazar} \cite{K-K} ont montr\'e qu'il existe des ''trous`` dans le comportement des profils des classes h\'er\'editaires de $\mathfrak S.$ Il ont montr\'e que si $\mathscr C$ est h\'er\'editaire, alors, soit $\varphi_{\mathscr{C}}$ est born\'ee par un polyn\^ome ($\varphi_{\mathscr{C}}(n)\leq n^c$ pour $n\geq 1$ et une constante $c>0$) et dans ce cas c'est un polyn\^ome, ou bien elle est minor\'ee par le nombre de Fibonacci,\index{nombre de Fibonacci} $\varphi_{\mathscr{C}}(n)\geq F_n$ pour tout $n\geq 1$.\\ \emph{Le nombre de Fibonacci} %\footnote{Le nombre de Fibonacci donn\'e ici est d\'ecal\'e par rapport \`a la d\'efinition classique qui est: $F_0=F_1=1,~F_n=F_{n-1}+F_{n-2},~n\geq 2.$}
$F_n$ est d\'efinie par :
$$F_1=1,~F_2=2,~~F_n=F_{n-1}+F_{n-2},~n\geq 3;$$
Par r\'ecurrence, on a $F_n\leq 2^{n-1}$ pour tout $n\geq 1$.\\

Par la suite, ils ont introduit les \emph{nombres de Fibonacci g\'en\'eralis\'es}\index{nombre de Fibonacci!g\'en\'eralis\'e} $F_{n,k}$, en montrant que si $\varphi_{\mathscr C}(n)<2^{n-1}$ pour un entier $n$ alors il existe un unique entier $k\geq 1$ et une constante $c>0$ tels que $F_{n,k}\leq \varphi_{\mathscr{C}}(n)\leq n^c F_{n,k}$ pour tout $n\geq 1$. Le \emph{nombre de Fibonacci g\'en\'eralis\'e} $F_{n,k}$ est le coefficient de $x^n$ dans le d\'eveloppement en s\'erie de l'expression $$\dfrac{1}{1-x-x^2-\ldots-x^k}.$$
 En particulier:
$F_{n,1}=1$ pour tout $n\geq 1$ et $F_{n,2}=F_n$. De mani\`ere plus g\'en\'erale,
$$\left\{
\begin{array}{l}
F_{n,k}=0~ \text{ pour } n<0,\\
F_{0,k}=1,\\
F_{n, k}=F_{n-1, k}+F_{n-2, k}+\ldots+F_{n-k,k}~\text{ pour } n>0.\\
\end{array}%
\right. $$

Ils montrent \'egalement que $F_{n,k}\leq 2^{n-1}$ et $F_{n,n}= 2^{n-1}$ pour un $k\geq 1$ fix\'e et $n\geq 1$.
Le th\'eor\`eme de Kaiser et Klazar\index{Kaiser et Klazar} est le suivant:

\begin{theorem}\label{theo:kaiser-klazar}
 Si $\mathscr C$ est une classe h\'er\'editaire de permutations, alors un seul parmi les quatre cas suivants se produit.
\begin{enumerate}
\item Pour $n$ suffisament large  la fonction $\varphi_{\mathscr{C}}(n)$ est \'eventuellement constante.
\item Il existe des entiers $a_0,\ldots,a_k$, $~k\geq 1$ et $a_k>0$, tels que $\varphi_{\mathscr{C}}(n)=a_0 \dbinom{n}{0}+\ldots+a_k \dbinom{n}{k}$ pour un $n$ suffisamment large. De plus, $\varphi_{\mathscr{C}}(n)\geq n$ pour tout $n.$
\item Il existe une constante $c>0$ et un unique entier $k\geq 2$, tels que $F_{n,k}\leq \varphi_{\mathscr{C}}(n)\leq n^c. F_{n,k}$ pour tout $n$.
\item On a $\varphi_{\mathscr{C}}(n)\geq 2^{n-1}$ pour tout $n$.
\end{enumerate}
\end{theorem}

Dans les cas $1$ \`a $3$ du th\'eor\`eme \ref{theo:kaiser-klazar}, la fonction g\'en\'eratrice est rationnelle \cite{klazar}. \emph{Albert et Atkinson} (2005) \cite{A-A} ont donn\'e des exemples de classes h\'er\'editaires dont les fonctions g\'en\'eratrices sont  alg\'ebriques. Nous rappellerons leur r\'esultat dans le paragraphe suivant, puis nous donnerons une g\'en\'eralisation au cas des structures binaires ordonn\'ees.

    \subsection{Permutations simples et th\'eor\`eme d'Albert-Atkinson}
\begin{figure}[!hbp]
\centering
\input{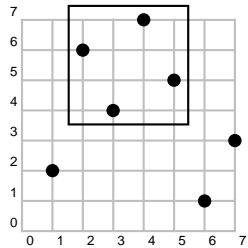}
\caption{\label{permutation-intervalle}Intervalle d'une permutation.}
\end{figure}

Soit un entier positif $n$ et soient $i, j\in [n]$ avec $i<j$. L'ensemble $[i,j]=\{i,i+1,\cdots,j\}$ est appel\'e un intervalle de $[n].$\\
Consid\'erons la permutation $\pi=2647513$ de $[7]$. L'intervalle $[2,5]$, c'est \`a dire l'ensemble des \'el\'ements cons\'ecutifs $2,3,4,5$ est envoy\'e par $\pi$ sur l'intervalle $[4,7]$ (voir \figurename~\ref{permutation-intervalle}). Le segment $[2,5]$ est appel\'e \emph{intervalle de $\pi$}.\index{intervalle!d'une permutation}

%************************image2*******************

En d'autres termes, si $\pi=a_1a_2\ldots a_n$ est une permutation de $[n]$, un \emph{intervalle}\index{permutation!intervalle d'une -} de $\pi$ est un ensemble d'indices cons\'ecutifs $I=[i,j]\subseteq [n]$ pour $1\leq i\leq j\leq n$, c'est \`a dire un intervalle de $[n]$, tel que l'ensemble des valeurs $\pi(I)=\{\pi(i):i\in I\}=\{a_i, a_{i+1},\cdots, a_j\}$ constitue un intervalle de $[n].$

Pour toute permutation $\pi$ de $[n]$, les singletons et l'ensemble $[n]$ sont des intervalles, ils sont dits \emph{triviaux}.\index{intervalle!trivial} Si $\pi$ n'a que des intervalles triviaux  elle est dite \emph{simple}.\index{permutation!simple} Exemple si $n\leq 2$ toutes les permutations sont simples. La permutation $\pi=58317462$ est simple. Les permutations simples de petites longueurs sont $1, 12, 21, 2413, 3142.$ Soit $S_n$ le nombre de permutations simples de $[n]$. Les valeurs de $S_n$ pour $n=1$ \`a $8$ sont: $1,2,0,2,6,46,338,2926$ (s\'erie A111111 dans \cite{Sloane}).
Asymptotiquement, $S_n$ \'equivaut \`a $\dfrac{n!}{e^2}$
\cite{A-A-K, Cor-Lou-Pem, nosaki}.\\% Ce r\'esultat signifie, qu'asymptotiquement, presque toutes les permutations sont simples.\\

Alors que les intervalles de permutations ont des applications en biomath\'ematiques, bioinformatique particuli\`erement pour les algorithmes g\'en\'etiques (voir \cite{Cor-Lou-Pem} pour explications et d'autres r\'ef\'erences), les permutations simples sont les \'el\'ements de base \`a partir desquels sont construites toutes les autres permutations, raison pour laquelle elles ont fait l'objet de nombreuses recherches ces derni\`eres ann\'ees.\\

Etant donn\'ees une permutation $\sigma$ de $[m]$ et des permutations $\alpha_1,\dots,\alpha_m$, l'\emph{inflation}\index{permutation!inflation d'une -} de $\sigma$ par $\alpha_1,\dots,\alpha_m$, not\'ee $\sigma[\alpha_1,\dots,\alpha_m]$ est la permutation obtenue en rempla\c{c}ant chaque terme $\sigma_i$ par un intervalle qui a la m\^eme forme que $\alpha_i.$ Exemple $2413[1,132,321,12]=479832156.$ L'inflation n'est autre que la somme lexicographique d\'efinie dans la section \ref{subsec:som-lex et decomp}, c'est ce terme qui est utilis\'e dans les articles qui traitent des permutations.

Inversement, une \emph{d\'eflation}\index{permutation!d\'eflation d'une -} de $\pi$ est toute expression de $\pi$ comme une inflation, c'est à dire $\pi=\sigma[\pi_1,\dots,\pi_k]$. La permutation $\sigma$ est alors appel\'ee \emph{quotient}\index{permutation!quotient d'une -} de $\pi.$ Le r\'esultat suivant est d\^u \`a \emph{Albert et Atkinson}\index{Albert et Atkinson} \cite{A-A}.

\begin{proposition}(Albert-Atkinson \cite{A-A}).
 Toute permutation peut-\^etre \'ecrite comme l'inflation d'une unique permutation simple.
Si $\pi$ peut s'\'ecrire comme $\sigma=\alpha[\sigma_1,\dots,\sigma_m]$ o\`u $\alpha$ est simple et $m\geq 4$ alors les $\alpha_i$ sont uniques.
\end{proposition}

Ce r\'esultat, qui donne la d\'ecomposition d'une permutation, est une sp\'ecialisation, au cas des permutations, des r\'esultats de d\'ecomposition des structures sur lesquels
nous reviendrons et dont l'un des pionniers est \emph{Gallai} \cite{gallai, Maf-Prei}.

\vspace{1mm}

\emph{Albert et Atkinson} montrent dans leur article \cite{A-A} que la connaissance des permutations simples dans une classe est souvent la clef pour comprendre sa structure et aborder le probl\`eme de son \'enum\'eration. Le th\'eor\`eme d'Albert et Atkinson est le suivant:

\begin{theorem}\label{theo:al-at}
  Si $\mathscr{C}$ est une classe h\'er\'editaire de permutations contenant un nombre fini de permutations simples alors la s\'erie g\'en\'eratrice  de $\mathscr{C}$, donn\'ee par
$\underset{n\geqq 1}\sum\varphi _{\mathscr{C}}(n)X^{n}$ est alg\'ebrique.
\end{theorem}

Nous donnerons, dans le chapitre \ref{sect:str.rela.bin}, une g\'en\'eralisation de ce th\'eor\`eme aux structures relationnelles binaires ordonn\'ees\index{structure relationnelle!binaire ordonn\'ee} (structure dont toutes les relations sont binaires l'une d'entre elles \'etant un ordre total).\\

Comme illustration du Th\'eor\`eme \ref{theo:al-at}, mentionnons que la classe des permutations qui ne  sont pas au dessus de  $2413$ et $3142$ ne contient aucune permutation simple non triviale (ces permutations sont appel\'ees \emph{permutations s\'eparables}).\index{permutation!s\'eparable}
%La s\'erie g\'en\'eratrice de cette classe est
%$\varphi(x)=\dfrac{1-x-\sqrt{1-6x+x^2}}{2}$ (voir \cite{A-A-V}).
Nous reviendrons sur la d\'efinition des permutations s\'eparables dans le paragraphe \ref{subsc:perm-separable}

    \subsubsection{Permutations exceptionnelles\index{permutation!exceptionnelle}} Etant donn\'ee une permutation simple $\pi$ on pourrait se demander quelles permutations simples sont contenues dans $\pi$?
     En particulier, y a t-il un terme qui pourrait \^etre supprim\'e de $\pi$ pour obtenir une suite qui aurait la m\^eme forme
     qu'une permutation simple\index{permutation!simple}? Ceci n'est pas vrai dans tous les cas. %il existe des permutations telles que la suppression de n'importe quel terme ne donne pas une permutation de m\^eme forme qu'une permutation simple.
     \emph{Albert et Atkinson}
     \cite{A-A}, adaptant, au cas des permutations, les r\'esultats de \emph{Schmerl et Trotter}\index{Schmerl et Trotter} \cite{S-T} sur les structures critiques (voir D\'efinition \ref{def:critique}),  montrent que la suppression d'un ou deux termes suffit. C'est ce qui est stipul\'e dans le th\'eor\`eme suivant, qui est une sp\'ecialisation au cas des permutations, d'un r\'esultat plus g\'en\'eral sur toute structure relationnelle dont les relations sont binaires et irr\'eflexives \cite{S-T}.

\begin{theorem}\cite{A-A, Br}
Toute permutation simple de longueur $n\geq 2$ contient une permutation simple de longueur $n-1$ ou $n-2$.
\end{theorem}

Dans la majorit\'e des cas, la suppression d'un point suffit. Si la suppression de n'importe quel terme de $\pi$ ne donne pas une suite de m\^eme forme qu'une permutation simple, alors $\pi$ est dite \emph{exceptionnelle}\index{permutation!exceptionnelle}.
%\footnote{Les permutations exceptionnelles correspondent aux structures critiques de Schmerl et Trotter \cite{S-T}.}
  \emph{Schmerl et Trotter}\index{Schmerl et Trotter} \cite{S-T} appellent de telles structures \emph{critiques}\index{structure binaire!critique}, ils pr\'esentent une caract\'erisation compl\`ete dans le cas des structures relationnelles binaires irr\'eflexives, des ensembles ordonn\'es, des graphes et des tournois. \emph{Albert et Atkinson}\index{Albert et Atkinson} \cite{A-A} ont donn\'e les formes des permutations exceptionnelles qui sont, en terme de permutations de $1,\dots,2m$ pour $m\geq 2:$\\

\noindent $(i)~~2.4.6....2m.1.3.5....2m-1.$\\
$(ii)~~2m-1.2m-3....1.2m.2m-2....2.$\\
$(iii)~~m+1.1.m+2.2....2m.m.$\\
$(iv)~~m.2m.m-1.2m-1....1.m+1.$

%%%%%%%%%%%%%%%%%%%%%%%%%%%%%%%%%%%%%%%%%%%%%%%%%%%%%%%%%%%%%%%%%%%%

%%%%%%%%%%%%%%%%%%%%%%%%%%%%%%%%%%%%%%%%%%%%%%%%%%%%%%%%%%%%%%%%%%%%%%%%

    \subsection{Permutations s\'eparables}\label{subsc:perm-separable}

Les \emph{permutations s\'eparables}\index{permutation!s\'eparable} ont \'et\'e d\'efinies dans \cite{A-A-V}. Ce sont des permutations qui peuvent-\^etre construites \`a partir de la permutation $1$ (permutation \`a un \'el\'ement) en appliquant, de mani\`ere it\'erative, les deux op\'erations de sommation suivantes d\'efinies, respectivement, sur les permutations $\pi$ de longueur $m$ et $\sigma$ de longueur $n$ par:\\
\begin{equation}
(\pi +\sigma)(i) =\left\{
\begin{array}{lll}
\pi (i) & & \text{si} \;1\leq i\leq m, \\
\sigma (i-m)+m & & \text{si} \;m+1\leq i\leq m+n,%
\end{array}%
\right.\label{equ:permut-separable1}
\end{equation}

\begin{equation}
(\pi \oplus \sigma)(i) =\left\{
\begin{array}{lll}
\pi (i)+n & & \text{si} \;1\leq i\leq m, \\
\sigma (i-m) & & \text{si} \;m+1\leq i\leq m+n,%
\end{array}%
\right. ~~~~~~~ \label{equ:permut-separable2}
\end{equation}

\begin{figure}[h]
\centering
\input{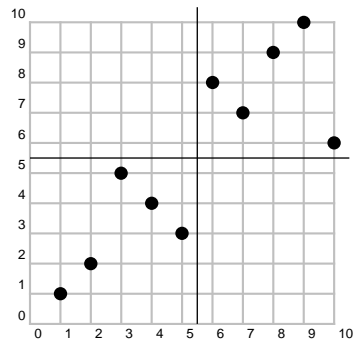}
\caption{\label{permutation-separable}Repr\'esentation d'une permutation s\'eparable.}
\end{figure}

 L'ensemble des permutations s\'eparables est form\'e par l'ensemble des permutations v\'erifiant \eqref{equ:permut-separable1} et \eqref{equ:permut-separable2} auquel nous incluons la permutation nulle. Du point de vu de la repr\'esentation graphique des permutations, une permutation non nulle est s\'eparable si elle est soit de longueur $1$ soit de longueur sup\'erieure \`a $1$ et il est possible de partitionner son graphe, par une ligne horizontale et une autre verticale, en quatre parties de telle sorte que deux parties oppos\'ees soient vides et les deux autres non vides et contiennent les graphes de permutations s\'eparables de plus petites longueurs, comme illustr\'e dans la \figurename~\ref{permutation-separable}.

 Bien que le terme "\emph{permutation s\'eparable}" date seulement du travail de \emph{Bose, Buss et Lubiw} (1998)\cite{Bose},  %%%a chercher)))),
ces permutations seraient apparues dans le travail d'\emph{Avis et Newborn} (1981) (voir \cite{A-A-V}). Les permutations s\'eparables sont l'analogue (pour les permutations) de deux autres classes d'objets qui ont \'et\'e largement \'etudi\'ees: les \emph{cographes} (ou les graphes sans $P_4$) et les ordres \emph{s\'eries-parall\`eles}\index{ordre!s\'eries-parall\`eles} (ou les ordres sans $N$). Un r\'esultat, qui d\'ecoule des caract\'erisations de ces deux classes sus-cit\'ees, caract\'erise les permutations s\'eparables.

\begin{proposition}(\cite{A-A-V})

Une permutation $\pi$ est s\'eparable si et seulement si elle ne contient ni la permutation $2413$ ni la permutation $3142$.
\end{proposition}

Ce qui signifie que la classe des permutations s\'eparables est, tout simplement, la classe h\'er\'editaire $Forb(\{2413,3142\}).$ Dans leur article \cite{A-A-V}, \emph{Albert, Atkinson et Vatter}\index{Albert, Atkinson et Vatter} donnent la fonction g\'en\'eratrice de cette classe:
\begin{proposition}
La fonction g\'en\'eratrice de la classe des permutations s\'eparables est:
$$\dfrac{1-x-\sqrt{1-6x+x^2}}{2}$$
et le nombre de permutations s\'eparables de longueur $n$ est le $(n-1)^{\text{\`eme}}$ nombre de Schr\"{o}der\index{nombre de Schr\"{o}der}.% large Schroder number.
\end{proposition}
La classe des permutations s\'eparables est alg\'ebrique, sa fonction g\'en\'eratrice $\varphi$ v\'erifie l'\'equation: $f=x+2f^2/(1+f)$ \cite{A-A-V}.

\vspace{1mm}

Nous rappelons que le $n^{\text{\`eme}}$ nombre de Schr\"{o}der $S_n$ est le nombre de chemins, sur une grille, allant du point $(0,0)$ au point $(n,n)$ en utilisant les \'etapes ou vecteurs $(0,1), (1,0)$ ou $(1,1)$ sans jamais traverser la droite $y=x$. Il est donn\'e par $S_n=\underset{i=0}{\overset{n}\sum}\binom{2n-i}{i}C_{n-i}$ o\`u $C_n=\dfrac{1}{n+1}\binom{2n}{n}$ est le $n^{\text{i\`eme}}$ nombre de Catalan (nombre d'arbres binaires ayant $n$ noeuds). La fonction g\'en\'eratrice des nombres de Schr\"{o}der v\'erifie $f=1+xf+xf^2$ (s\'erie A006318 \cite{Sloane}).

\vspace{2mm}

Les permutations s\'eparables ont \'et\'e g\'en\'eralis\'ees aux $d$-permutations s\'eparables \cite{Asin-Mans}. Une $d$-permutation s\'eparable est une suite de $d$ permutations de m\^eme longueur $n$ dont la premi\`ere est la permutation identit\'e $12\ldots n.$ Pour $d=2$, on obtient la permutation s\'eparable.

%%%%%%%%%%%%%%%%%%%%%%%%%%%%%%%%%%%%%%%%%%%%%%%%%%%%%%%%%%%%%%%%%%%%%%%%%%%%%%%%
%$$$$$$$$$$$$$$$$$$$$$$$$$$$$$$$$$$$$$$$$$$$$$$$$$$$$$$$$$$$$$$$$$$$$$$$$$$$$$
%£$$$$$$$$$$$$$$$$$$$$$$$$$$$$$$$$$$$$$$$$$$$$$$$$$$$$$$$$$$$$$$$$$$$$$$$$$$$$$
        %       Les bichaînes

%$$$$$$$$$$$$$$$$$$$$$$$$$$$$$$$$$$$$$$$$$$$$$$$$$$$$$$$$$$$$$$$$$$$$$$$$$$$$$$$

\section{Bicha\^{i}nes}\label{sec:bichaine}

Dans cette section nous codons les bicha\^{i}nes par des permutations  et traduisons, en termes de la th\'eorie des relations, les r\'esultats sur les permutations.\\

Une \emph{bicha\^{i}ne}\index{bichaine@bicha\^{i}ne} est une paire $\mathcal B:=(E,(L_1,L_2))$ o\`u $L_1$ et $L_2$ sont deux ordres totaux sur l'ensemble $E$, nous la noterons par la suite $\mathcal B:=(E,L_1,L_2)$. %L'ensemble $E$ est la \emph{base} ou le \emph{domaine} de $\mathcal B$.
Les notions d'isomorphisme, d'isomorphisme local et d'abritement d\'efinies dans la section \ref{sec:isom-abrit-equim} se d\'eduisent de fa\c{c}on naturelle sur les bicha\^{i}nes.

\vspace{1mm}

Soient $\mathcal B:=(E,L_1,L_2)$ et $\mathcal B':=(E',L'_1,L'_2)$ deux bicha\^{i}nes. Un \emph{isomorphisme}\index{isomorphisme!de bicha\^{i}nes} de $\mathcal B$ sur $\mathcal B'$ est une bijection $f:E\rightarrow E'$ telle que
$$(x,y)\in L_i \;~~\text{si et seulement si}\;~~ (f(x),f(y))\in L'_i,\;~~\forall x,y\in E\; \text{et}\; i\in \{1,2\}.$$

\noindent Si $A\subseteq E$, la \emph{restriction} de $\mathcal B$ \`a $A$ est la bicha\^{i}ne $\mathcal B_{\restriction_A}:=(A,{L_1}_{\restriction_A}, {L_2}_{\restriction_A})$ o\`u ${L_i}_{\restriction_A}:={L_i}\cap A^2$ (si $L_i$ est donn\'e sous forme de paires ordonn\'ees).

\vspace{1mm}

Un \emph{isomorphisme local}\index{isomorphisme!local} de $\mathcal B$ dans $\mathcal B'$ est tout isomorphisme d'une restriction $\mathcal B_{\restriction_A}$ de $\mathcal B$ \`a un sous-ensemble $A$ de $E$ sur une restriction $\mathcal {B'}_{\restriction_{A'}}$ de $\mathcal B'$ \`a un sous-ensemble $A'$ de $E'.$

\vspace{1mm}

Une bicha\^{i}ne $\mathcal B$ s'\emph{abrite} dans une bicha\^{i}ne $\mathcal B'$, nous notons $\mathcal B\leq \mathcal B'$, si $\mathcal B$ est isomorphe \`a une restriction de $\mathcal B'$.
Deux bicha\^{i}nes $\mathcal B$ et $\mathcal B'$ telles que $\mathcal B\leq \mathcal B'$ et $\mathcal B'\leq \mathcal B$ sont \emph{\'equimorphes}.\index{bichaine@bicha\^{i}ne!\'equimorphe} Si leurs domaines sont finis elles sont \emph{isomorphes}.\index{bichaine@bicha\^{i}ne!isomorphe }
Cette relation d'abritement d\'efinit un pr\'eordre sur la classe des bicha\^{i}nes. Sur la classe des types d'isomorphie c'est un ordre.
Dans ce cas particulier et comme sugg\'er\'e par \emph{Cameron}\index{Cameron} \cite{cameron}, les types d'isomorphie peuvent-\^etre cod\'es par des permutations, comme nous le verrons dans ce qui suit.

%************************************************
\subsection{Bicha\^{i}nes et permutations}
Soit $n\geq 1$ et soit $\sigma$ une permutation de $[n]$. Associons la bicha\^{i}ne $\mathcal C_{\sigma}:= ([n], \leq, \leq_{\sigma})$ o\`u "$\leq$" est l'ordre naturel sur $[n]$
et "$\leq_{\sigma}$" est l'ordre induit par $\sigma$ sur $[n]$ et d\'efini par:
\begin{equation}
i\leq_{\sigma} j\;~~\text{si et seulement si}\;~~\sigma_i \leq \sigma_j,\;~~\forall i,j \in [n].\label{eq:ordr-bichaine}
\end{equation}

\begin{example}\label{exp:ordre-permut}
Si $\pi =2647513$ alors, l'ordre lin\'eaire  $\leq_{\pi}$ sur $[7]=\{1,...,7\}$ est donn\'e par:
$$ 6< _{\pi }1< _{\pi }7< _{\pi }3<_{\pi }5< _{\pi }2< _{\pi}4.$$
De mani\`ere plus pr\'ecise, l'ordre  $\leq _{\pi }$ est donn\'e par:
$$\pi _{1}^{-1}< _{\pi }\pi _{2}^{-1}< _{\pi
}\pi _{3}^{-1}< _{\pi }\pi _{4}^{-1}< _{\pi }\pi
_{5}^{-1}< _{\pi }\pi _{6}^{-1}< _{\pi }\pi _{7}^{-1}.$$
Cette permutation et l'ordre qui lui est associ\'e sont repr\'esent\'es dans la \figurename~\ref{permut-ordre} en page \pageref{permut-ordre}. Les deux ordres $\leq$ et $\leq_{\pi}$ s'obtiennent naturellement à partir de cette repr\'esentation, en \'enum\'erant les points $(i,\pi_i)$ (identifi\'e à $i$) en parcourant la repr\'esentation de gauche à droite pour le premier ordre (l'ordre naturel $\leq$) et de bas en haut pour le deuxi\`eme ordre ($\leq_{\pi}$). L'ordre $P_{\pi}$ de la \figurename~\ref{permut-ordre} est l'intersection de ces deux ordres lin\'eaires.
\end{example}

\emph{Cameron}, dans son article \cite{cameron}, consid\`ere une permutation $\sigma$ de $[n]$ comme une paire d'ordres totaux sur $[n]$, le premier \'etant l'ordre naturel et le second l'ordre $\sigma_1<\sigma_2<\dots<\sigma_n.$ La bicha\^{i}ne obtenue ainsi est celle associ\'ee \`a $\sigma^{-1}$ avec l'ordre donn\'e par l'\'equation \eqref{eq:ordr-bichaine}.

\begin{lemma}\label{lem:isom-bichaine-type}
Si $\mathcal B:=(E,L_1,L_2)$ est une bicha\^{i}ne finie\footnote{Si $E=\varnothing$ alors $\mathcal B$ est la bicha\^{i}ne vide, son type d'isomorphie est repr\'esent\'e par la permutation nulle.} alors $\mathcal B$ est isomorphe \`a la bicha\^{i}ne $\mathcal C_{\sigma}$ pour un unique $\sigma$ sur $[\vert E\vert~].$
\end{lemma}

\begin{proof}
Observons dans un premier temps que si $L$ est un ordre lin\'eaire sur un $n$-ensemble $E$ ($n\in \mathbb N$), alors il existe un unique isomorphisme $l$ de $([n],\leq)$ sur $(E,L)$, ce qui revient \`a dire qu'il existe une fa\c{c}on unique d'\'ecrire les \'el\'ements de $E$ comme:\\ $l(1)<_L\cdots <_L l(n)$ suivant l'ordre $L$.

Posons $\vert E\vert=n$ et soient $l_1$ et $l_2$ les isomorphismes associ\'es \`a $L_1$ et $L_2$ respectivement. Posons $\sigma={l_1}^{-1}\circ l_2$, alors l'application $l_1:[n]\rightarrow E$ est un isomorphisme de $\mathcal C_{\sigma}$ sur $\mathcal B$.\\
En effet, en remarquant que $l_1$ transforme l'ordre naturel "$\leq$" sur $[n]$ en $L_1$, il suffit de montrer que $l_1$ transforme "$\leq_{\sigma}$" en $L_2$, c'est \`a dire:$$\sigma_i\leq \sigma_j\Leftrightarrow l_1(\sigma_i)\leq l_1(\sigma_j).$$ Comme $l_1\circ \sigma =l_2,$ nous avons l'\'equivalence recherch\'ee.
\end{proof}
\bigskip

\noindent Comme corollaire de ce lemme nous avons:

\begin{corollary}
L'application $\Phi:\sigma\mapsto \mathcal C_{\sigma}$ est une bijection de $\mathfrak S$ dans l'ensemble des types d'isomorphie des bicha\^{i}nes.
\end{corollary}

Nous avons donc un moyen de repr\'esenter les types d'isomorphie des bicha\^{i}nes. Comment s'\'etend la relation d'abritement?\\

Soient $\sigma$ et $\sigma'$ deux permutations sur $[n]$ et $[n']$ respectivement. Nous avons le r\'esultat suivant:

\begin{lemma}\label{lem:perm-bichaine}
 $\sigma \leq \sigma'$ si et seulement si $\mathcal C_{\sigma}\leq \mathcal C_{\sigma'}.$
\end{lemma}

\begin{proof}
Pour la condition n\'ecessaire, supposons $\sigma \leq \sigma'$ et posons $A=\{x_1,x_2,\cdots,x_n\}$ l'ensemble des indices donn\'es par la relation \eqref{perm-ordre}, on a $A\subseteq [n']$. Consid\'erons l'application,

$
%\begin{array}{c}
l:[n]\rightarrow A\text{ d\'efinie par }
l(i):=x_i \text{ pour tout }i\in [n].%
%\end{array}%
$

Nous avons bien $$i\leq_{\sigma}j \text{ si et seulement si }x_i\leq_{\sigma'}x_j,$$ donc  $l$ est un isomorphisme de $\mathcal C_{\sigma}$ sur ${\mathcal C_{\sigma'}}_{\restriction_A}$. D'o\`u $\mathcal C_{\sigma}\leq \mathcal C_{\sigma'}.$\\

Pour la condition suffisante, si $\mathcal C_{\sigma}\leq \mathcal C_{\sigma'}$ alors $\mathcal C_{\sigma}$ est isomorphe \`a ${\mathcal C_{\sigma'}}_{\restriction_B}$ pour un $n$-ensemble $B$ de $[n']$. Notons par $f$ cet isomorphisme, alors pour tous $i,j\in [n]$ nous avons:
\begin{enumerate}
\item $i\leq j$ si et seulement si $f(i)\leq f(j)$,
\item $i\leq_{\sigma}j$ si et seulement si $f(i)\leq_{\sigma'}f(j)$
\end{enumerate}
Posons $x_i=f(i)$ pour tout $i\in [n]$, les $x_i$ v\'erifient la relation \eqref{perm-ordre}, donc $\sigma \leq \sigma'$.
\end{proof}
\bigskip

Ce simple fait, donn\'e par le Lemme \ref{lem:perm-bichaine}, permet d'\'etudier la classe des permutations au moyen de la th\'eorie des relations.
En particulier, les classes h\'er\'editaires de permutations correspondent aux  classes h\'er\'editaires de bicha\^{i}nes  et comme nous le verrons ci-dessous, les permutations simples correspondent aux bicha\^{i}nes ind\'ecomposables.
%%%%%%%%%%%%%%%%%%%%%%%%%%%%%%%%%%%%%%%%%%%%%%%%%%%%%%%%%%%%%%%%%%%%%%%%%%%%%%

\subsection{Bicha\^{i}ne ind\'ecomposable}

Soit $L$ un ordre lin\'eaire sur un ensemble $E$ et soient $x, y\in E$ avec $x<_L y$. L'ensemble $[x,y]_L=\{z\in E/ x\leq_L z\leq_L y\}$ est un intervalle pour $L$.

\begin{definition}\label{def:inter-bichaine}
Soit $\mathcal B:=(E,L_1,L_2)$ une bicha\^{i}ne. Un sous-ensemble $A$ de $E$ est un \emph{intervalle}\index{bichaine@bicha\^{i}ne!intervalle d'une -} pour $\mathcal B$ (\'egalement appel\'e \emph{autonome}) si pour tout $i\in \{1,2\}$:
\begin{equation}(a,b)\in L_i\Leftrightarrow (a',b)\in L_i,\text{ pour tout } a, a'\in A\text{ et } b\notin A.\label{eq:bich-indec}\end{equation}
\end{definition}

L'ensemble vide $\varnothing$, les singletons $(\{x\}_{x\in E})$ et l'ensemble $E$ sont des intervalles, ils sont dits \emph{triviaux}.\index{intervalle!trivial}

Il ressort de la D\'efinition \ref{def:inter-bichaine} et de l'\'equation \eqref{eq:bich-indec} que pour une bicha\^{i}ne $\mathcal B:=(E,L_1,L_2)$, un sous-ensemble propre non vide $A$ de $E$ est un intervalle si et seulement si c'est un intervalle commun pour $L_1$ et $L_2$.

\begin{definition}
Une bicha\^{i}ne $\mathcal B$ est dite \emph{ind\'ecomposable}\index{bichaine@bicha\^{i}ne!ind\'ecomposable} si elle ne poss\`ede pas d'intervalle non trivial.
\end{definition}

Donc, la bicha\^{i}ne $\mathcal B:=(E,L_1,L_2)$ est ind\'ecomposable si et seulement si $L_1$ et $L_2$ ne poss\`edent pas d'intervalle en commun.

\vspace{1mm}

Soit $\mathcal B:=(E,L_1,L_2)$ une bicha\^{i}ne avec $\vert E\vert=n$ et soit $\mathcal C_{\sigma}:=([n],\leq,\leq_{\sigma})$ la bicha\^{i}ne associ\'ee \`a la permutation $\sigma$ de $[n]$ qui d\'efinit le type d'isomorphie de $\mathcal B$. Il est facile de voir, du fait de l'isomorphisme entre $\mathcal B$ et $\mathcal C_{\sigma}$,  que tout intervalle commun \`a  $L_1$ et $L_2$ induit un intervalle commun \`a   "$\leq$" et "$\leq_{\sigma}$" et vice-versa, d'o\`u:

\begin{lemma}
Une bicha\^{i}ne $\mathcal B$ est ind\'ecomposable si et seulement si la permutation $\sigma$ qui d\'efinit son type d'isomorphie est simple.
\end{lemma}

\begin{proof}
D'apr\`es ce qui pr\'ec\`ede, il suffit de montrer que $\sigma$ est simple si et seulement si $\mathcal C_{\sigma}$ est ind\'ecomposable.\\
Si $I$ est un intervalle de $\mathcal C_{\sigma}$, alors $I$ est un intervalle de $[n]$ (car c'est un intervalle pour $\leq$). Posons alors  $I=\{i,i+1,\dots,j-1,j\}$. D'apr\'es l'\'equation \eqref{eq:ordr-bichaine}, l'ordre $\leq_{\sigma}$ est donn\'e par: $\sigma_1^{-1}<_{\sigma}\sigma_2^{-1} <_{\sigma}\dots<_{\sigma} \sigma_n^{-1}$. Comme $I$ est un intervalle pour $\leq_{\sigma}$ alors il existe $k<l$ avec $l-k=j-i$ tels que $\{\sigma_k^{-1}, \sigma_{k+1}^{-1},\dots, \sigma_l^{-1}\}=I$. Donc $\sigma(I)=\{\sigma_i,\dots, \sigma_j\}=\{k,k+1,\dots,l\}$ qui est un intervalle de $[n]$. D'o\`u $I$ est un intervalle de $\sigma.$\\ Inversement, si $I$ est un intervalle de $\sigma$, alors $I$ et $\sigma(I)$ sont des intervalles de $[n]$. Donc $I$ est un intervalle de $\leq$ et il existe $m<k$ tels que $\{m,m+1,\dots,k\}=\sigma(I)$ avec $\vert I\vert=k-m+1$. Donc, $\sigma^{-1}\{m,m+1,\dots,k\}=\{\sigma_m^{-1},\sigma_{m+1}^{-1},\dots,\sigma_k^{-1}\}=I$, ce qui implique que $I$ est un intervalle pour $\leq_{\sigma}.$ Il s'ensuit que $I$ est un intervalle pour $\mathcal C_{\sigma}.$
\qedhere\end{proof}

\bigskip
La notion d'ind\'ecomposabilit\'e est plut\^ot ancienne. La notion d'intervalle remonte \`a Fra\"{\i}ss\'e\index{Fraisse@Fra\"{\i}ss\'e} \cite{fraisse2}, voir aussi \cite{fraisse3}. Un r\'esultat fondamental
sur la d\'ecomposition des structures binaires en intervalles a \'et\'e obtenu par Gallai \cite{gallai}
(voir \cite{ehren-H-Ros} pour d'autres extensions). Il n'est donc pas surprenant que certains  r\'esultats sur les permutations simples \'etaient d\'ej\`a connus (par exemple leur \'evaluation asymptotique \cite{nosaki}).\medskip

Avec cette correspondance bijective entre les permutations et les types d'isomorphie des bicha\^{i}nes, le Th\'eor\`eme \ref{theo:al-at} d'Albert et Atkinson devient:\index{Albert et Atkinson}

\begin{theorem}
Si $\mathscr C$ est une classe h\'er\'editaire de bicha\^{i}nes finies contenant un nombre fini de bicha\^{i}nes ind\'ecomposables alors la s\'erie g\'en\'eratrice de $\mathscr C$ est alg\'ebrique.
\end{theorem}

Nous \'etablirons dans le chapitre \ref{sect:str.rela.bin}  une extension de ce r\'esultat aux structures binaires ordonn\'ees.

%**************************************************************************
\subsection{Bicha\^{i}ne s\'eparable}\label{susec:bichaineseparable}

Les bicha\^{i}nes \emph{s\'eparables}\index{bichaine@bicha\^{i}ne!s\'eparable} sont des bicha\^{i}nes qui peuvent-\^etre construites \`a partir de la bicha\^{i}ne %vide ($\emptyset$)???????(il faudra peut-etre supprimer le vide) et la bicha\^{i}ne
unit\'e $1$ (d\'efinie sur un ensemble \`a un \'el\'ement) en appliquant, de mani\`ere it\'erative, les deux op\'erations de sommation suivantes:

Si $\mathcal B:=(E,L_1,L_2)$ et $\mathcal B':=(E',L'_1,L'_2)$ sont deux bicha\^{i}nes d\'efinies sur deux ensembles disjoints $E$ et $E'$, de cardinalit\'es $n\geq 1$ et $m\geq 1$ respectivement (si les deux ensembles ne sont pas disjoints, nous pouvons recopier les relations sur des ensembles disjoints et effectuer les op\'erations), nous posons:
$$\begin{array}{l}
\mathcal B+\mathcal B':=(E\cup E',L_1+L'_1,L_2+L'_2)\qquad\qquad\\
\mathcal B\oplus\mathcal B':=(E\cup E',L_1+L'_1,L'_2+L_2)%
\end{array}%
$$
o\`u $L_1+L'_1$ d\'esigne la somme lexicographique de $L_1$ et $L'_1$, c'est \`a dire:

$$x\leq_{(L_1+L'_1)}y\;\text{si}\left\{
\begin{array}{l}
x,y\in E\;~ \text{et}\; x\leq_{L_1}y\\
\text{ou bien}\\
x,y\in E'\;~ \text{et}\; x\leq_{L'_1}y\\
\text{ou bien}\\
x\in E\;\text{et}\; y\in E'%
\end{array}%
\right. $$
Nous incluons \`a cet ensemble la bicha\^{i}ne vide.\\

Si $\sigma$ et $\sigma'$ sont les permutations de $[n]$ et $[m]$ respectivement, qui repr\'esentent les types d'isomorphie des bicha\^{i}nes $\mathcal B$ et $\mathcal B'$ respectivement alors:

\begin{proposition}
 Les permutations $\sigma+ \sigma'$ et $\sigma\oplus\sigma',$ d\'efinies par les relations \eqref{equ:permut-separable1} et \eqref{equ:permut-separable2} repr\'esentent les types d'isomorphie des bicha\^{i}nes $\mathcal B+\mathcal B'$ et $\mathcal B\oplus\mathcal B'$ respectivement\footnote{La bichaîne nulle est repr\'esent\'ee par la permutation nulle.}.
\end{proposition}

\begin{proof}
D'apr\`es le Lemme \ref{lem:isom-bichaine-type}, les bicha\^{i}nes $\mathcal B$ et $\mathcal B'$ sont isomorphes aux bicha\^{i}nes $\mathcal C_{\sigma}$ et $\mathcal C_{\sigma'}$ respectivement. Donc $\mathcal B+\mathcal B'$ est isomorphe \`a $\mathcal C_{\sigma}+\mathcal C_{\sigma'}$ et $\mathcal B\oplus\mathcal B'$ est isomorphe \`a $\mathcal C_{\sigma}\oplus\mathcal C_{\sigma'}$ (en prenant soin, auparavant, de recopier les relations de $\mathcal C_{\sigma'}$ sur l'ensemble $[m]+n$, le translat\'e de l'ensemble $[m]$ obtenu en rajoutant la valeur $n$ \`a chaque \'el\'ement de l'ensemble $[m]$). Nous devons donc montrer que $\mathcal C_{\sigma}+\mathcal C_{\sigma'}=\mathcal C_{\sigma+\sigma'}$ et $\mathcal C_{\sigma}\oplus\mathcal C_{\sigma'}=\mathcal C_{\sigma\oplus\sigma'}.$\\
Pour la premi\`ere \'egalit\'e, il est clair que la somme lexicographique de l'ordre naturel sur $[n]$ et de l'ordre naturel sur $[m]+n$ donne l'ordre naturel sur $[n+m].$ D'un autre c\^ot\'e, la somme lexicographique $\leq_{\sigma}+\leq_{\sigma'}$ donne l'ordre:
$$\sigma_1^{-1}<\sigma_2^{-1}<\dots<\sigma_n^{-1}<{\sigma'}_1^{-1}+n<{\sigma'}_2^{-1}+n<\dots<{\sigma'}_n^{-1}+n.$$
Comme, d'apr\`es la relation \eqref{equ:permut-separable1}, $\sigma_i^{-1}=(\sigma+\sigma')_i^{-1}$ pour tout $1\leq i\leq n$ et ${\sigma'}_i^{-1}+n=(\sigma+\sigma')_{i-n}^{-1}+n$ pour tout $n+1\leq i\leq n+m$, alors nous avons $\leq_{\sigma}+\leq_{\sigma'}=\leq_{\sigma+\sigma'}.$ D'o\`u le r\'esultat.\\
La deuxi\`eme \'egalit\'e se d\'emontre de la m\^eme mani\`ere.
\end{proof}
\bigskip

Le r\'esultat suivant d\'ecoule alors de fa\c{c}on naturelle:

\begin{corollary}
Une bicha\^{i}ne finie $\mathcal B$ est s\'eparable si et seulement si la permutation qui repr\'esente son type d'isomorphie est s\'eparable.
\end{corollary}

\noindent Les r\'esultats d'\'enum\'eration donn\'es dans \cite{A-A-V} nous permettent d'\'enoncer ce qui suit:

\begin{proposition}
La fonction g\'en\'eratrice de la classe des bicha\^{i}nes s\'eparables est $$\dfrac{1-x-\sqrt{1-6x+x^2}}{2}$$
et le nombre de bicha\^{i}nes s\'eparables \`a $n$ \'el\'ements (\`a l'isomorphisme pr\`es) est le $n^{\text{i\`eme}}$ nombre de Schr\"{o}der.% \footnote{}
\end{proposition}

\clearemptydoublepage

\chapter{Structures relationnelles binaires ordonn\'ees}\label{sect:str.rela.bin}

Dans ce chapitre, nous g\'en\'eralisons des r\'esultats donn\'es dans \cite{A-A} et \cite{A-A-V} au cas des structures binaires ordonn\'ees.

\vspace{1mm}

Une  \emph{structure relationnelle binaire},\index{structure relationnelle!binaire} ou \emph{structure binaire}\index{structure binaire} ou encore une \emph{$2$-structure}\index{$2$-structure} est une structure $\mathcal R:=(E,(\rho_i)_{i\in I})$ o\`u  les $\rho_i$ sont des relations binaires sur $E.$
Soit un entier $k\geq 1$. Nous appelons \emph{structure binaire de \emph{type}\index{structure binaire!de type $k$} $k$}, ou simplement \emph{structure de type $k$} s'il n'y a pas de confusion, une structure binaire comportant $k$ relations binaires, c'est \`a dire lorsque $\vert I\vert=k$.
  Nous d\'esignons par $\Omega_k$ l'ensemble des structures binaires finies de type $k$.

\vspace{1mm}

Une  structure relationnelle binaire est \emph{ordonn\'ee}\index{structure binaire!ordonn\'ee} si une de ses relations $\rho_i$ est un ordre lin\'eaire. Nous conviendrons que c'est la premi\`ere relation $\rho_1$. Donc, une structure relationnelle binaire est ordonn\'ee si elle s'\'ecrit  $\mathcal R:=(E,\leq,(\rho_j)_{j\in J})$ o\`u "$\leq$" est un ordre lin\'eaire sur $E$ et les $\rho_j$ sont des relations binaires sur $E.$ Nous appellerons \emph{structure binaire ordonn\'ee de type $k$},\index{structure binaire!ordonn\'ee de type $k$} ou \emph{structure ordonn\'ee de type $k$} s'il n'y a pas de confusion,  une structure binaire ordonn\'ee form\'ee d'un ordre total et de $\vert J\vert=k$ relations binaires.\\
Les exemples basiques de ces structures sont les cha\^{i}nes, ou structures ordonn\'ees de type $0$ (pour $J=\varnothing$), les bicha\^{i}nes, ou structures ordonn\'ees de type $1$ (pour $J=\{1\}$ et $\rho_1$ est un ordre lin\'eaire) et les multicha\^{i}nes ($J$ fini et $\rho_j$ est un ordre lin\'eaire pour tout $j\in J$). Nous d\'esignons par $\Theta_k$ la collection des structures binaires ordonn\'ees finies de type $k$. Il est clair que $\Theta_k\subseteq\Omega_{k+1}$.

    \section{Ind\'ecomposabilit\'e et somme lexicographique}\label{subsec:som-lex et decomp}
    \subsection{Intervalle}
\begin{definition}\label{def:intervalle}
Soit  $\mathcal{R}:=(E,(\rho_i)_{i\in I})$ une structure binaire\index{structure binaire}. Un sous-ensemble $A$ de $E$ est un \emph{intervalle}\index{intervalle} de $\mathcal{R}$ si pour tout $i\in I$:
 $$(x\rho_i a \Leftrightarrow x\rho_i a')  \; \text{et} \; (a\rho_i x \Leftrightarrow a'\rho_i x) \; \text{pour tous} \; a,a'\in A \;\text{et}\;
x\notin A.$$
\end{definition}

En d'autres termes, $A$ est un intervalle de $\mathcal R$ si un \'el\'ement $x\notin A$ "voit" tous les \'el\'ements de $A$ de la m\^eme fa\c{c}on et deux \'el\'ements quelconques $a_1$ et $a_2$ de $A$ "voient" tout \'el\'ement $x\notin A$ de la m\^eme fa\c{c}on.

\vspace{1mm}

L'appellation \emph{intervalle} n'est pas la seule utilis\'ee dans la litt\'erature, on y trouve l'appellation \emph{module}, ou \emph{ensemble homog\`ene}, en anglais \emph{autonomous}, \emph{partitive set} ou \emph{clan}.

\vspace{1mm}

L'ensemble vide $\varnothing$, les singletons $\{x\}_{x\in E}$ et l'ensemble $E$ sont des intervalles et sont dits \emph{triviaux}.
Si $\mathcal{R}$ ne poss\`ede pas d'intervalle non trivial elle est dite \emph{ind\'ecomposable}\index{structure binaire!ind\'ecomposable}. Avec cette d\'efinition toute structure \`a deux \'el\'ements est ind\'ecomposable. Une structure ind\'ecomposable  d'au moins trois \'el\'ements est \emph{premi\`ere}\index{structure binaire!premi\`ere} (\emph{"prime"} en anglais).\\ %non triviale

Voici les propri\'et\'es les plus connues des intervalles:
\begin{enumerate}
\item Tout intervalle de $\mathcal{R}:=(E,(\rho_i)_{i\in I})$ est un intervalle de $\mathcal{R}^c:=(E,({\rho_i}^c)_{i\in I})$ et de $\mathcal{R}^{-1}:=(E,({\rho_i}^{-1})_{i\in I}).$
\item Si $A$ est un intervalle de $\mathcal{R}:=(E,(\rho_i)_{i\in I})$ et $F\subseteq E$ alors $F\cap A$ est un intervalle de $\mathcal R_{\restriction _F}.$
\item Si $A$ est un intervalle de $\mathcal{R}$ et $B\subseteq A$ alors $B$ est un intervalle de $\mathcal{R}$ si et seulement si $B$ est un intervalle de $\mathcal{R}_{\restriction_A}.$
\item Si $A$ et $B$ sont deux intervalles de $\mathcal{R}$ alors $A\cap B$ est un intervalle de $\mathcal{R}$; de plus si $A\cap B\neq\varnothing$ alors $A\cup B$ est un intervalle de $\mathcal{R}$ et si $A\setminus B\neq\varnothing$ alors $B\setminus A$ est un intervalle de $\mathcal R.$
%\item Si $\mathcal{R}:=(E,(\rho_i)_{i\in I})$ et $\mathcal{R'}:=(E,(\rho'_i)_{i\in I})$ sont deux structures binaires de m\^emes domaine et arit\'e telles que $\mathcal{R}\subseteq\mathcal{R'}$ (ce qui signifie que $\rho_i\subseteq\rho'_i$ pour tout $i\in I$) alors tout intervalle de $\mathcal{R}$ est un intervalle de $\mathcal{R'}.$ (\cite{ehren}) (a verifier, ca me semble incorrecte)
\end{enumerate}

Nous d\'esignons par $Ind(\Omega_k)$ (resp. $Ind(\Theta_k)$) la classe des membres ind\'ecomposables de $\Omega_k$ (resp. de $\Theta_k$).\\

Nous rappelons ici les r\'esultats de \emph{Schmerl et Trotter}\index{Schmerl et Trotter} \cite{S-T} sur des structures ind\'ecomposables particuli\`eres, les structures \emph{critiques}\index{structure binaire!critique}, notion introduite en (1993):

 \begin{definition}\label{def:critique}
 Une structure binaire $\mathcal{R}:=(E,(\rho_i)_{i\in I})$ est \textit{critique} si elle est premi\`ere et  $\mathcal{R}_{\restriction_{E\setminus \{x\}}}$ est non ind\'ecomposable pour tout $x\in E$.
 \end{definition}
  Notons que, d'apr\`es cette d\'efinition,  $E$ a au moins quatre \'el\'ements.
% Cette notion a \'et\'e introduite par Schmerl et Trotter \cite{S-T} en (1993).
Parmi les r\'esultats donn\'es dans leur article, nous avons:% le th\'eor\`eme suivant et son corollaire:

 \begin{theorem}\label{theo:schmerl-trotter}\cite{S-T}
Toute structure ind\'ecomposable de taille $n\geq 3$ abrite, au moins, une structure ind\'ecomposable de taille $n-1$ ou $n-2$.
\end{theorem}

\begin{theorem}\label{theo:indec-schmer-trott}\cite{S-T}
  Soit $\mathcal{R}:=(E,(\rho_i)_{i\in I})$ une structures binaire ind\'ecomposable d'ordre $n\geq 7$. Alors il existe deux \'el\'ements distincts $c, d\in E$ tels que $\mathcal{R}_{\restriction_{E\setminus \{c,d\}}}$ est ind\'ecomposable.
\end{theorem}

\begin{corollary}\label{cor:indec-schmer-trott}\cite{S-T}
 Supposons que  $\mathcal{R}$ soit une structure ind\'ecomposable d'ordre $n$ non critique et supposons  $5\leq m\leq n$.  Alors $\mathcal{R}$ poss\`ede une  sous-structure ind\'ecomposable d'ordre $m$.
\end{corollary}

 Nous disons qu'une classe de $Ind(\Omega_k)$ est h\'er\'editaire dans $Ind(\Omega_k)$ si elle contient tout membre de $Ind(\Omega_k)$ qui s'abrite dans l'un de ses membres.
Un r\'esultat analogue au Th\'eor\`eme \ref{theo:classehered-ideal} pour les classes de $Ind(\Omega_k)$ est le suivant:

\begin{theorem}\label{theo:ind-ideal}
Toute classe h\'er\'editaire infinie de $Ind(\Omega_k)$, contient un id\'eal infini.
\end{theorem}
\begin{proof}
Soit $\mathscr I$ une classe h\'er\'editaire infinie de $Ind(\Omega_k)$.\\
\quad\textbf{Cas 1:} $\mathscr I$ ne poss\`ede pas d'anticha\^{i}ne infinie, donc $\mathscr I$ est belordonn\'e et est donc une union fini d'id\'eaux (Th\'eor\`eme \ref{theo:erdos-tarski}). $\mathscr I$ \'etant infini, l'un des id\'eaux est infini.\\
\quad\textbf{Cas 2:} $\mathscr I$ poss\`ede une anticha\^{i}ne infinie. Nous d\'emontrons de la m\^eme fa\c{c}on que pour la Proposition \ref{prop:classe-ideal} que $\mathscr I$ poss\`ede une section initiale non vide et belordonn\'ee $\mathcal J$ ayant un nombre infini de bornes. Pour montrer que $\mathcal J$ est infini, nous utilisons le Th\'eor\`eme \ref{theo:schmerl-trotter} de Schmerl et Trotter. Si $\mathcal J$ \'etait fini, la taille de ses \'el\'ements serait born\'ee. Soit $k$ cette borne. D'apr\`es Schmerl et Trotter, les bornes de $\mathcal J$ seraient de tailles $k+1$ ou $k+2$ et seraient donc en nombre fini.\\
Si $\mathcal J$ est filtrante alors $\mathcal J$ est un id\'eal, sinon, $\mathcal J$, \'etant belordonn\'ee, est une union finie d'id\'eaux et l'un d'eux est infini.
\end{proof}

  \subsection{Composition et r\'eduction}\label{compo-reduction}

Soient $\mathcal R:=(E,(\rho_i)_{i\in [k]})\in \Omega_k$ et $v\in E$. L'op\'eration de \emph{composition de $\mathcal R$ par $\mathcal S$ suivant $v$}\index{operation@op\'eration de composition}, o\`u $\mathcal S\in \Omega_k$ et $dom(\mathcal S)\cap E=\varnothing$, consiste en la substitution\footnote{A noter que dans cette op\'eration, le fait que $(v,v)$ appartienne ou non \`a $\rho_i,$ pour un $i\in I$, n'appara\^{i}t pas dans la structure r\'esultante.} du sommet $v$ de $E$ par la relation $\mathcal S$ dans $\mathcal R$. Dans la structure r\'esultante, not\'ee $\mathcal R[v:\mathcal S]$,  qui est une structure $\mathcal R':=(E',(\rho'_i)_{i\in [k]})$ o\`u $E'=(E\setminus\{v\})\cup dom(\mathcal S)$, tout \'el\'ement de $E\setminus \{v\}$ est reli\'e \`a tous les \'el\'ements de $dom(\mathcal S)$ de la m\^eme fa\c{c}on qu'il est reli\'e \`a $v$ dans $\mathcal R$, c'est \`a dire, pout tout $i\in [k]$:
$$\forall u\in E\setminus \{v\},\; \left\{\begin{array}{l}
(u,v)\in \rho_i\Leftrightarrow (u,x)\in \rho'_i\\
\text{ et }\\
 (v,u)\in \rho_i\Leftrightarrow (x,u)\in \rho'_i \\
 \end{array}\right. \; \forall x\in dom(\mathcal S).$$

Avec cette d\'efinition, un sous-ensemble $A$ de $E$ est un intervalle de $\mathcal R$ si $\mathcal R$ est de la forme $\mathcal R'[v:\mathcal S]$ o\`u $\mathcal S=\mathcal R_{\restriction_A}.$
\vspace{1mm}

\noindent Une propri\'et\'e de l'op\'eration de composition est la suivante: si $w\in dom(\mathcal S)$ alors
$$\mathcal R[v:\mathcal S][w:\mathcal S']=\mathcal R[v:\mathcal S[w:\mathcal S']].$$

Soit $\mathcal R:=(E,(\rho_i)_{i\in [k]})\in\Omega_k$ poss\'edant un intervalle non trivial $A$. L'op\'eration de \emph{r\'eduction}\index{operation@op\'eration de r\'eduction} de $\mathcal R$ suivant $A$ consiste  en la contraction de l'ensemble $A$ en un sommet $x_A$ ($x_A\notin E$), la structure r\'esultante, not\'ee $red_A[\mathcal R]$ est une structure $\mathcal R''= (E'',(\rho''_i)_{i\in [k]})$ d\'efinie par:
\begin{itemize}
\item $E''=(E\setminus A)\cup \{x_A\}$ avec $x_A\notin E$.
\item pour tout $i\in [k]$, la relation $\rho''_i$ est donn\'ee par:
  \begin{itemize}
  \item $(x,y)\in \rho''_i\Leftrightarrow(x,y)\in \rho_i,~~\; \forall x,y\in E\setminus A,$
      \item $(x,x_A)\in \rho''_i\Leftrightarrow(x,a)\in \rho_i,~~\; \forall a\in A\;\text{et}\; x\in E\setminus A,$
          \item $(x_A,x)\in \rho''_i\Leftrightarrow(a,x)\in \rho_i,~~\; \forall a\in A \;\text{et}\; x\in E\setminus A,$
          \item $(x_A,x_A)\in\rho''_i\Leftrightarrow \rho_i$ est r\'eflexive $~~\forall i\in [k].$
      \end{itemize}
  \end{itemize}
\bigskip

Il existe un lien entre les deux op\'erations de composition et de r\'eduction. Sur les structures o\`u toutes les relations sont r\'eflexives ou irr\'eflexives elles sont inverses l'une de l'autre, en d'autres termes:
\begin{itemize}
\item Si $\mathcal R'=red_A[\mathcal R]$ alors $\mathcal R=\mathcal R'[x_A:\mathcal R_{\restriction_A}],$
\item Si toutes les relations de $\mathcal R$ et $\mathcal S$ sont r\'eflexives (ou irr\'eflexives) alors:
     $$\mathcal R'=\mathcal R[v:\mathcal S]\Rightarrow \mathcal R=red_{dom(S)}[\mathcal R'].$$
\end{itemize}

  \subsection{Somme lexicographique et d\'ecomposition de $2$-structures}
\begin{definition}\label{def:somme lexicographique}
   Soient $\mathcal{R}:=(E,(\rho_i)_{i\in [k]})$ une structure binaire de $\Omega_k$ et $\mathfrak{F}:=(\mathcal{S}_x)_{x\in E}$ une  famille de structures binaires
$\mathcal{S}_x:=(E_x,({\rho_i}^x)_{i\in [k]})$, index\'ee par les
\'el\'ements de $E$. Les ensembles $E$ et $E_x$ sont suppos\'es non-vides.
La \emph{somme lexicographique}\index{somme!lexicographique} de $\mathfrak{F}$ suivant $\mathcal{R}$ ou tout simplement la $\mathcal R$-\emph{somme} de $\mathfrak{F}$, not\'ee par $\underset{x\in \mathcal{R}}{\oplus }\mathcal{S}_x$ (ou encore si $E=\{x_1,\cdots,x_n\},\mathcal{S}_{x_1}\underset{\mathcal R}{\oplus} \cdots \underset{\mathcal R}{\oplus}\mathcal{S}_{x_n}$) est la structure binaire $\mathcal{T}$ obtenue par la composition de $\mathcal R$ par $\mathcal S_x$ suivant $x$, pour tout $x\in E$. Autrement dit, en substituant \`a chaque \'el\'ement $x\in E$  la structure $\mathcal{S}_x$. Plus pr\'ecis\'ement $\mathcal{T}=(Z,(\tau_i)_{i\in [k]})$ o\`u  $Z:=\{(x,y):x\in E,~y\in E_x\}$ et pour tout $i\in [k]$,
$(x,y)\tau_i (x',y')$  si
\begin{itemize}
    \item ou bien $x\neq x'$ et $x\rho_i x'$
    \item ou bien  $x=x'$ et $y{\rho_i}^x y'$.
    \end{itemize}
\end{definition}
\medskip

\noindent De mani\`ere triviale, si on remplace chaque $\mathcal{S}_x$ de la D\'efinition \ref{def:somme lexicographique} par une structure  binaire isomorphe $\mathcal S_x'$, alors $\underset{x\in \mathcal{R}}{\oplus }\mathcal S_x'$ est isomorphe \`a
 $\underset{x\in \mathcal{R}}{\oplus }\mathcal{S}_x.$
 Donc, on peut supposer que les domaines des structures $\mathcal{S}_x$ sont deux \`a deux disjoints. Dans ce cas on peut modifier la d\'efinition \ref{def:somme lexicographique} ci-dessus, en posant
$Z:=\underset{x\in E}{\cup}{E_x}$ et pour deux \'el\'ements
 $z\in E_x$ et $z'\in E_{x'}$ et tout $i\in [k]$, nous avons $z{\tau_i} z'$ si
 \begin{itemize}
\item  ou bien $x\neq x' \;\text{et}\; x\rho_i x'$,
\item ou bien  $x=x'\; \text{et}\; z{\rho_i}^x z'.$
\end{itemize}
\medskip

 Avec cette d\'efinition, chaque ensemble $E_x$ est un intervalle de la somme $\mathcal T.$ Par ailleurs, soit $Z/_{\equiv}$ le quotient de $Z$ constitu\'e par les blocs de cette partition\index{partition} en intervalles.\\ Soit $p:Z\rightarrow Z/_{\equiv}$
la surjection canonique de $Z$ dans $Z/_{\equiv}$ et soit $\mathcal R'$ l'image de $\mathcal T$ par $p$ (ce qui signifie que $\mathcal R'=(Z/_{\equiv},(\rho'_i)_{i\in [k]})$ o\`u
$\rho'_i=\{(p(x_1),p(x_2)):(x_1,x_2)\in \rho_i \}$).\\
 Si on identifie chaque bloc $E_x$ \`a l'\'el\'ement $x,$ alors $\mathcal R$ et $\mathcal R'$ co\"{\i}ncident sur les paires d'\'el\'ements distincts. Elles co\"{\i}ncident si on ne consid\`ere que des  relations r\'eflexives (ou irr\'eflexives).
Inversement, si  $\mathcal{T}:=(Z,(\tau_i)_{i\in [k]})$ est une structure binaire et $(E_x)_{x\in E}$ une partition de $Z$ en intervalles non vides de
$\mathcal T,$ alors $\mathcal T$ est une somme lexicographique de
$(\mathcal T_{\restriction_{E_x}})_{x\in E}$ suivant le quotient $Z/_{\equiv}.$
%En d'autres termes:
De mani\`ere plus formelle:

\begin{definition}\label{def:decomposition}
Soit $\mathcal{R}:=(E,(\rho_i)_{i\in [k]})$ une structure binaire et $\mathscr P:=\{E_1,\cdots, E_q\}$ avec $q>1$ une partition\index{partition!en intervalles} de $E$ en intervalles non vides de $\mathcal R$.
La \emph{d\'ecomposition}\index{structure binaire!d\'ecomposition d'une -} de $\mathcal R$ suivant $\mathscr P$ est la donn\'ee de la structure $\mathcal S$ et de la famille des structures
$\mathscr H=\{\mathcal R_1,\cdots, \mathcal R_q\}$ telles que $\mathcal R$ est la somme lexicographique de $\mathscr H$ suivant $\mathcal S$, autrement dit, $\mathcal R=\mathcal R_1\underset{\mathcal S}{\oplus}\cdots\underset{\mathcal S}{\oplus}\mathcal R_q$ o\`u, pour tout $i\in \{1,\cdots,q\}$ on a
$\mathcal R_i=\mathcal R_{\restriction_{E_i}}$ et $\mathcal S=red_{E_1,\cdots, E_q}[\mathcal R]$ est la r\'eduction de $\mathcal R$ suivant $E_1,\cdots, E_q$ simultan\'ement. $\mathcal S$ est appel\'ee \emph{structure quotient}\index{structure binaire!quotient} de $\mathcal R$ par $\mathscr P$ et est not\'ee
$\mathcal S=\mathcal R/{\mathscr P}$.
\end{definition}

\begin{observation}
Toute structure relationnelle binaire $\mathcal{R}:=(E,(\rho_i)_{i\in [k]})$ admet au moins deux d\'ecompositions, dites \emph{triviales}:\index{structure binaire!d\'ecomposition triviale d'une -} la d\'ecomposition suivant la partition $\mathscr P:=\{\{x\},\; x\in E\}$ (dans laquelle $\mathcal S=\mathcal R$ et $\vert E_i\vert=1$ pour tout $1\leq i\leq \vert E\vert$) et celle suivant $\mathscr P':=\{E\}$ (dans laquelle $\mathscr H=\{\mathcal R\}$ et $\mathcal S$ est une structure \`a un \'el\'ement). Ces partitions sont dites \emph{partitions triviales}.\index{partition!triviale}
\end{observation}

D'apr\`es ce qui pr\'ec\`ede nous avons:
\begin{observation}
 Il existe une correspondance bijective entre les d\'ecompositions des structures binaires en sommes  lexicographiques et les  partitions de leurs  domaines en intervalles.
\end{observation}

Disons qu'une somme lexicographique $\underset{x\in \mathcal{R}}{\oplus }\mathcal{S}_x$ est triviale si $\vert E \vert=1$ ou  $\vert E_x \vert=1$ pour tout $x\in E$, sinon elle est \emph{non triviale}; aussi, une structure binaire est \emph{somme-ind\'ecomposable}\index{structure binaire!somme-ind\'ecomposable} si elle n'est pas isomorphe \`a une somme lexicographique non triviale. Nous avons imm\'ediatement:

%Ce qui nous am\`ene \`a l'observation suivante:
\begin{observation}\label{obs:unique-decomp}
Une structure binaire $\mathcal R$ est ind\'ecomposable si et seulement si elle poss\`ede uniquement deux d\'ecompositions si et seulement si elle est somme-ind\'ecomposable.
\end{observation}

\noindent Une propri\'et\'e  importante de ces d\'ecompositions est la suivante:

\begin{propertie}
 L'ensemble des  partitions de $E$, ordonn\'e par affinement, est un treillis dont l'ensemble des  partitions de $E$ en intervalles de $\mathcal R$ est un sous-treillis.
\end{propertie}
Nous rappelons qu'une partition $\mathscr P$ est plus fine qu'une partition $\mathscr P'$ si tout \'el\'ement de $\mathscr P'$ est une union d'\'el\'ements de $\mathscr P$, c'est \`a dire que les parties de  $\mathscr P'$ sont fractionn\'ees dans $\mathscr P$ en plus petites parties. La partition la plus fine est la partition triviale form\'ee par des singletons et la moins fine est celle form\'ee de l'ensemble $E$.

    \subsection{Partition forte et th\'eor\`eme de d\'ecomposition}

Parmi toutes les d\'ecompositions d'une structure binaire finie $\mathcal R$, il existe une qui est canonique. Cette d\'ecomposition est due \`a Gallai \cite{gallai, Maf-Prei}. Nous la d\'esignons par \emph{d\'ecomposition de Gallai},\index{structure binaire!d\'ecomposition de Gallai d'une -} elle est d\'efinie ci-dessous.

\begin{definition}
Deux ensembles $A$ et $B$ se rencontrent si $A\cap B\neq \varnothing$; ils \emph{se chevauchent} si $A\cap B\neq \varnothing, \; A\setminus B\neq \varnothing$ et $B\setminus A\neq \varnothing.$
\end{definition}

\begin{definition}\label{def:intervalle fort,maximal}
Soit $\mathcal R$ une structure binaire. Un intervalle $A$ de $\mathcal R$ est dit \emph{fort}\index{intervalle!fort} s'il est non vide et ne chevauche aucun autre intervalle. Autrement dit, $A$ est fort si:
\begin{itemize}
\item $A\neq\varnothing$.
\item Pour tout intervalle $B\neq A$ on a $B\cap A\neq \varnothing\Rightarrow B\subset A$ ou $A\subset B$.
\end{itemize}
 Un intervalle fort $A$ est dit \emph{maximal}\index{intervalle!fort!maximal} s'il est maximal pour l'inclusion parmi les intervalles forts.
\end{definition}

\noindent Il est clair que les intervalles triviaux non vides sont forts.

\begin{lemma}\label{lem:part-forte}
Pour toute structure binaire  finie $\mathcal R$, il existe une unique partition de son domaine, $dom(\mathcal R)$, en intervalles forts maximaux.% qu'on appellera \emph{partition forte.}\index{partition!forte}
\end{lemma}

\begin{proof}
Posons $E=dom(\mathcal R)$. Si $\vert E\vert =1$, le r\'esultat est \'evident. Supposons que $E$ poss\`ede au moins deux \'el\'ements. L'existence de cette partition est due au fait que tout \'el\'ement $x\in E$ est inclus dans un intervalle fort propre, en l'occurrence $\{x\}$. Donc, $x\in A_x$, o\`u $A_x$ est un intervalle fort maximal parmi les intervalles forts propres contenant $x$. Aussi, si $A_x\cap A_y\neq\varnothing$ alors $A_x=A_y$.\\ Pour l'unicit\'e, supposons que $\mathcal R$ poss\`ede deux telles partitions $\mathscr P:=(A_j)_{j\in J}$ et $\mathscr P':=(B_k)_{k\in K}$. Il existe alors deux indices $j$ et $k$ tels que $A_j\neq B_k$ et $A_j\cap B_k\neq\varnothing$. Nous avons les deux cas suivants:\\1)- Soit $A_j\subset B_k$ (ou bien $B_k\subset A_j$) ce qui constitue une contradiction avec la maximalit\'e de $A_j$ ou de $B_k$.\\2)- Soit $A_j$ et $B_k$ se chevauchent, ce qui contredit le fait que $A_j$ et $B_k$ sont forts.\\
Il s'ensuit qu'une telle partition est unique.
\end{proof}

\bigskip

Nous d\'esignons cette partition par \emph{partition de Gallai}.
\vspace{1mm}

\noindent La \emph{d\'ecomposition de Gallai} de $\mathcal R$ est sa d\'ecomposition suivant la \emph{partition de Gallai}. Le quotient de $\mathcal R$ induit par cette partition est le \emph{quotient de Gallai}.\\

C'est sur la partition de Gallai  que se base le th\'eor\`eme de d\'ecomposition des $2$-structures finies (Th\'eor\`eme \ref{theo:decomposition1}). Ce th\'eor\`eme (voir \cite{C-D, gallai, Ke} pour les graphes, \cite{ehren-H-Ros} pour les $2$-structures) que nous donnons ici sous une autre forme (Th\'eor\`eme \ref{theo:decomposition}) est, semble t-il, \`a maintes fois red\'ecouvert, depuis le th\'eor\`eme de d\'ecomposition des graphes, appel\'e th\'eor\`eme de d\'ecomposition modulaire, d\^u \`a \emph{Gallai}\index{Gallai} (1967) \cite{gallai}, \'etendu par Kelly (1985)\index{Kelly} \cite{Ke} au cas des graphes infinis. Pour le cas des graphes dirig\'es, les g\'en\'eralisations et d'autres r\'ef\'erences voir \cite{C-D, ehren-H-Ros, Montg}.

\begin{remark}
Dans \cite{ehren-H-Ros}, la $2$-structure ou structure binaire n'a pas la m\^eme d\'efinition que celle que nous avons adopt\'ee dans ce document, elle consiste en un couple $g:=(D,R)$ o\`u $D$ est un ensemble fini non vide et $R$ une relation d'\'equivalence sur $D_2$, $R\subseteq D_2\times D_2$ avec $D_2=\{(x,y)/ x\in D, y\in D, x\neq y\}$. Les $2$-structures, \`a notre sens, sont des $2$-structures au sens de \cite{ehren-H-Ros} mais avec plus d'informations; en effet, soit $\mathcal{S}:=(E,(\rho_i)_{i\in I})$ une $2$-structure telle que nous l'avons d\'efinie, avec toutes les relations irr\'eflexives. Soit $R$ une relation sur $E_2$ d\'efinie par
$$(x_1,y_1)R(x_2,y_2) \Leftrightarrow \left\{\begin{array}{l}
(x_1,y_1)\in\rho_i\Leftrightarrow (x_2,y_2)\in\rho_i\\
\text{et}\\
(y_1,x_1)\in\rho_i\Leftrightarrow (y_2,x_2)\in\rho_i\end{array}\right.\text{ pour tout } i\in I$$
 $R$ est une relation d'\'equivalence sur $E_2$ et $g:=(E,R)$ est une $2$-structure, au sens de \cite{ehren-H-Ros}. Cette d\'efinition rejoint plus celle du squelette (skeleton) d'une structure binaire d\'efini par \emph{Schmerl et Trotter} \cite{S-T}.\index{Schmerl et Trotter}
\end{remark}

\begin{theorem}\label{theo:decomposition1}
Soit $\mathcal R$ une $2$-structure finie ayant au moins deux \'el\'ements et $\mathscr P_{max}$ sa partition de Gallai. Alors le quotient $\mathcal R/{\mathscr P_{max}}$ poss\`ede au moins deux \'el\'ements et a une seule des formes suivantes:
\begin{enumerate}
\item $\mathcal R/{\mathscr P_{max}}$ est ind\'ecomposable,
\item Chacune des relations qui le composent a une seule des formes suivantes:
    \begin{enumerate}
    \item c'est une relation vide.
    \item c'est une relation compl\`ete.
    \item c'est un ordre lin\'eaire. De plus, si deux relations ont cette forme, elles sont \'egales ou duales.
    \end{enumerate}
\end{enumerate}
\end{theorem}

\begin{remark}\label{rem:dec-enchai}
Si une relation $\mathcal R$ est vide ou compl\`ete sur un ensemble $X$, alors toute partie de $X$ est un intervalle de $\mathcal R$. Si $\mathcal R$ est un ordre lin\'eaire, ses intervalles sont de la forme $[x,y]:=\{z\in X/ x\leq z\leq y\}$ pour tout $x, y\in X,~~x\leq y.$ Les relations vides, compl\`etes et les ordres lin\'eaires ne poss\`edent pas d'intervalle fort non trivial. Nous pouvons les d\'ecomposer de plusieurs mani\`eres, notament en somme lexicographique de deux relations index\'ee par une relation \`a deux \'el\'ements. Nous imposerons, ci-dessous, une seule d\'ecomposition dans  ces cas l\`a pour les structures binaires ordonn\'ees.
\end{remark}

Avec la Remarque \ref{rem:dec-enchai}, l'Observation \ref{obs:unique-decomp} et le Lemme \ref{lem:part-forte} nous avons un corollaire imm\'ediat du Th\'eor\`eme \ref{theo:decomposition1}:

\begin{corollary}\label{cor:sumlex}
 Soit $\mathcal{R}$ une structure binaire finie ayant au moins deux \'elements. Alors  $\mathcal{R}$ est isomorphe \`a une somme lexicographique   $\underset{x\in \mathcal{S}}{\oplus }\mathcal{R}_x$ o\`u  $\mathcal{S}$ est ind\'ecomposable ayant au moins deux \'elements. De plus, lorsque $\mathcal S$ a au moins trois \'elements, la partition de $\mathcal R$ en intervalles est unique.
\end{corollary}

Gr\^{a}ce \`{a} l'unicit\'e de la d\'ecomposition de Gallai (Lemme \ref{lem:part-forte}), nous arrivons au th\'eor\`eme de d\'ecomposition suivant (voir par exemple \cite{C-D} pour des extensions aux structures  infinies).
\medskip

\noindent Disons qu'une structure binaire $\mathcal R:=(E,(\rho_i)_{i\in I})$ est  \emph{encha\^{\i}nable}\index{structure binaire!encha\^{i}nable} %\footnote{Nous reviendrons sur la notion d'encha\^{i}nabilit\'e dans le  chapitre \ref{chap:monomorphe}.}
s'il existe un ordre lin\'{e}aire $\leq$ sur $E$ v\'erifiant, pour tout $i\in I$,
$x\rho_i y\Leftrightarrow x'\rho_i y'$ pour tous  $x,y, x',y'$ tels que
$x\leq y\Leftrightarrow x'\leq y'$. Si $\mathcal R$ est r\'eflexive, ceci implique que toute relation $\rho_i$ est soit la relation d'\'egalit\'e
$\triangle_E$, la relation compl\`ete $E\times E$ ou un ordre lin\'eaire;
 de plus, si $\rho_i$ et $\rho_j$, pour $i\neq j$, sont deux ordres lin\'eaires, ils co\"{\i}ncident ou sont oppos\'es. Notons que si de plus $\mathcal R$ est ordonn\'ee alors les relations  $\rho_i$ qui sont des ordres lin\'eaires sont \'egales ou oppos\'ees \`{a} l'ordre donn\'e.

\begin{theorem}\label{theo:decomposition}
Soit $\mathcal R$ une structure binaire finie ayant au moins deux \'el\'ements alors $\mathcal R$ est une somme  lexicographique $\underset{x\in \mathcal{S}}{\oplus}{\mathcal R}_x$ o\`{u} $\mathcal S$ est soit ind\'ecomposable avec au moins trois \'el\'ements soit une structure binaire encha\^{\i}nable ayant au moins deux \'el\'ements et les ensembles $E_x=dom({\mathcal R}_x)$ sont des intervalles forts maximaux de $\mathcal R$.
\end{theorem}

Nous \'etendons la notion de somme lexicographique aux collections de structures binaires non vides. Etant donn\'ees une structure binaire non vide  $\mathcal{R}$ et des classes  $\mathscr{A}_x$ de strucures binaires non vides pour tout $x\in V(\mathcal{R})$, notons
 $\underset{x\in \mathcal{R}}{ \oplus}\mathscr{A}_x$ la classe de toutes les structures binaires de la forme
$\underset{x\in \mathcal{R}}{\oplus }\mathcal{S}_x$ avec $\mathcal{S}_x \in \mathscr{A}_x$ pour tout $x\in V(\mathcal R).$ Si $\mathscr{A}_x=\mathscr A$ pour tout $x\in V(\mathcal R)$  cette classe sera not\'ee  $\underset{ \mathcal{R}}{ \oplus}\mathscr{A}$. Si  $\mathcal R:=(\{0,1\},\leq,(\rho_i)_{i\in J})$,
 avec $0<1$, $\mathscr{A}_0:=\mathscr A$  et  $\mathscr{A}_1 := \mathscr B$, nous posons   $\underset{x\in \mathcal{R}}{ \oplus}\mathscr{A}_x=\mathscr A\underset{\mathcal R}{\oplus}\mathscr {B}$.

\noindent Aussi, si $\mathscr A$ et $\mathscr B$ sont deux classes de structures binaires, nous posons
$$\underset{\mathscr{A}}{\oplus}\mathscr{B}:=\{\underset{x\in \mathcal{R}}{\oplus }\mathcal{S}_x:~~\mathcal{R}\in \mathscr{A},~~\mathcal{S}_x \in \mathscr{B} \; \text{pour tout} \; x\in dom(\mathcal{R}) \}.$$

Nous disons qu'une collection  $\mathscr C$ de structures binaires est \emph{ferm\'ee par sommes} si $\underset{\mathscr{C}}{\oplus}{\mathscr{C}}\subseteq{\mathscr{C}}$.
La \emph{cl\^oture par sommes}  de $\mathscr{C}$, est la plus petite classe, $cl(\mathscr{C})$, ferm\'ee par sommes contenant $\mathscr{C}$. Si nous posons  $\mathscr C_0:=\mathscr C$ et $\mathscr C_{n+1}:=\underset{\mathscr{C}}{\oplus}\mathscr C_n$, alors $cl(\mathscr{C})=\underset{n=0}{\overset{\infty}{\cup}} \mathscr{C}_{n}$ (notons que $\underset{\mathscr C_{n}}{\oplus}
\mathscr C_n \subseteq \underset{\mathscr C}{\oplus}
\mathscr C_{2n}$).
Si $\mathscr C$ est form\'ee de structures r\'eflexives et contient la structure \`a un \'el\'ement,  alors $cl(\mathscr C)=\underset{cl(\mathscr{C})}{\oplus}{cl(\mathscr{C})}$ (en effet, toute structure $\mathcal R\in cl(\mathscr C)$ est une somme lexicographique suivant $\mathcal R$ de copies de la structure \`a un \'el\'ement). Si de plus $\mathscr C$ contient la structure vide alors $cl(\mathscr C)$ est h\'er\'editaire.

    \subsection{Classe h\'er\'editaire et structures ind\'ecomposables}

  Si $\mathscr C$ est une sous-classe de $\Omega_k$, nous notons par $Ind(\mathscr C)$ la collection de ses membres ind\'ecomposables et si $\mathcal{R}$ est une structure binaire, nous noterons par $Ind(\mathcal{R})$ la collection de ses sous-structures finies qui sont ind\'ecomposables. Par exemple, si $\mathcal{R}$ est un cographe\index{cographe} ou un ordre s\'eries-parall\`{e}les\index{ordre!s\'eries-parall\`eles} alors les membres de $Ind(\mathcal{R})$ ont au plus deux \'el\'ements.

\begin{theorem}\label{theo:closed}
Soit $\mathscr{D}$ une classe h\'er\'editaire de $Ind(\Omega_k)$. Posons $\sum \mathscr{D}:= \{\mathcal{R} \in {\Omega_k} : Ind(\mathcal{R}) \subseteq\mathscr{D}\}$.
 Si tous les membres de $\mathscr D$ sont r\'eflexifs\footnote{Nous avons le m\^eme r\'esultat si nous supposons tous les membres de $\mathscr D$ irr\'eflexifs ou si $\mathscr D$ contient toutes les structures \`a un \'el\'ement de $\Omega_k.$} alors $\sum \mathscr{D}=cl(\mathscr{D}).$
\end{theorem}

\begin{proof}

L'inclusion $\sum \mathscr{D}\supseteq {cl(\mathscr{D})}$ a lieu sous l'hypoth\`ese que tous les membres de $\mathscr{D}$ sont r\'eflexifs. Inversement, si $\mathcal{R}\notin cl(\mathscr{D})$ alors ou bien $\mathcal{R}$ est ind\'ecomposable auquel cas $\mathcal{R}\notin \mathscr{D}$ ou bien $\mathcal{R}$ ne peut pas s'exprimer comme une somme lexicographique de structures de $\mathscr{D}$ donc $\mathcal{R}\notin \sum \mathscr{D}.$
\end{proof}
\medskip

Le r\'esultat suivant donne la relation qui existe entre les bornes de $\mathscr D$ dans $Ind(\Omega_k)$ et les bornes de $\sum \mathscr{D}$ dans $\Omega_k$.

\begin{proposition}\label{prop:borneindec}
Si $\mathscr B$ est l'ensemble des bornes de $\mathscr D$ dans $Ind(\Omega_k)$ et $\mathscr B'$ l'ensemble des bornes de $\sum \mathscr{D}$ dans $\Omega_k$ alors $\mathscr B=\mathscr B'$.
\end{proposition}

\begin{proof}
Soit $\mathcal B\in\mathscr B$, donc $\mathcal B$ est ind\'ecomposable et $\mathcal B\notin\mathscr D.$ Comme $\mathcal B$ est une borne de $\mathscr D$ alors $\forall \mathcal B'<\mathcal B$ et $\mathcal B'$ ind\'ecomposable on a $\mathcal B'\in\mathscr D$. On a $\mathcal B\notin\sum \mathscr{D}=cl(\mathscr{D})$ mais $\mathcal B_{-x}\in\sum \mathscr{D}=cl(\mathscr{D}),~\forall x\in\mathcal B$, en effet, si $\mathcal B_{-x}$ est ind\'ecomposable alors $\mathcal B_{-x}\in\mathscr D\subset\sum\mathscr D$ et si $\mathcal B_{-x}$ n'est pas ind\'ecomposable alors tout $\mathcal B'$ ind\'ecomposable tel que $\mathcal B' <\mathcal B_{-x}$ se trouve dans $\mathscr D$, c'est \`a dire que $\mathcal B_{-x}\in\sum \mathscr{D}$. Il s'ensuit que  $\mathcal B\in \mathscr B'.$

Inversement, soit $\mathcal B\in\mathscr B'$, c'est \`a dire que $\mathcal B\notin\sum \mathscr{D}$, donc $\mathcal B$ abrite un ind\'ecomposable $\mathcal B'\notin\mathscr D$ d'où $\mathcal B'\notin\sum \mathscr{D}$ ce qui implique que $\mathcal B'=\mathcal B$ car $\mathcal B$ est une borne. Mais $\mathcal B'\in\mathscr B$, car $\mathcal B'$ est ind\'ecomposable, donc $\mathcal B\in\mathscr B$. Il s'ensuit que $\mathscr B=\mathscr B'.$
\end{proof}
\medskip

Comme cons\'equence de ce r\'esultat nous avons:
\begin{consequence}
Toutes les bornes de $\sum \mathscr{D}$ sont ind\'ecomposables.
\end{consequence}

Inversement, nous avons:
\begin{lemma}\label{lem:ind-sum}
Si $\mathscr C$ est une classe de $\Omega_k$ ayant toutes ses bornes ind\'ecomposables alors $\mathscr C=\sum Ind(\mathscr C).$
\end{lemma}
\begin{proof}
L'inclusion $\mathscr C\subseteq\sum Ind(\mathscr C)$ est \'evidente et a lieu sans conditions. Supposons que l'inclusion inverse n'a pas lieu et soit $\mathcal R$ un \'el\'ement de $\sum Ind(\mathscr C)\setminus\mathscr C.$ Puisque $\mathcal R\notin \mathscr C$ alors $\mathcal R$ abrite une borne, disons $\mathcal B$, de $\mathscr C$, mais $\mathcal R\in\sum Ind(\mathscr C)$ et $\mathcal B$ est ind\'ecomposable par hypoth\`ese, donc $\mathcal B\in Ind(\mathscr C)\subset\mathscr C,$ ce qui constitue une contradiction.
\end{proof}
%\vspace{1mm}

\subsection{Rôle du belordre}

Nous avons  le r\'esultat suivant:
\begin{proposition}\label{prop:wqo}
Si une classe h\'er\'editaire $\mathscr{D}$ de $Ind(\Omega_k)$ est h\'er\'editairement belordonn\'ee alors  $\sum \mathscr{D}$ est h\'er\'editairement belordonn\'ee et $\sum \mathscr{D}$ a un nombre fini de bornes.
\end{proposition}

\begin{proof}
La seconde partie de la proposition vient du r\'esultat  d\^u \`a Pouzet (Th\'eor\`eme \ref{theo:pouzet-borne}). %Dans notre cas de  structures binaires  la preuve est directe (voir la preuve de Pouzet???).
La premi\`ere partie utilise les propri\'et\'es des ordres belordonn\'es et d\'ecoule du th\'eor\`eme de Higman (1952) sur  les alg\`ebres pr\'eordonn\'ees par divisibilit\'e  \cite{higman52}. Sans rappeler le r\'esultat nous donnons une preuve directe.  Soit $\mathcal A$ un ordre qui est belordonn\'e, montrons que $(\sum \mathscr D). \mathcal A$ est belordonn\'e. Si  $(\sum \mathscr D). \mathcal A$ n'est pas belordonn\'e, alors d'apr\`es l'un des r\'esultats pr\'eliminaires d'Higman\index{Higman} (voir  Th\'eor\`eme \ref{theo:higman-equivalence}), il contient une section finale non finiment engendr\'ee. D'apr\'es le lemme de Zorn, il existe une section finale, disons $\mathcal F$, maximale pour l'inclusion parmi les sections finales ayant cette propri\'et\'e (car l'ensemble des sections finales ordonn\'e par inclusion est un ensemble inductif).

 Posons $\mathcal I:= (\sum \mathscr D). \mathcal A \setminus \mathcal F$ le compl\'ementaire de $\mathcal F$ dans $(\sum \mathscr D). \mathcal A.$ L'ensemble $\mathcal I$ est donc belordonn\'e.
Soit $\mathscr R:=(\mathcal R_0,f_0),\cdots, (\mathcal R_n,f_n), \cdots$ une anticha\^{i}ne infinie d'\'el\'ements  minimaux de $\mathcal F.$ Comme  $\mathscr D. \mathcal A$ est belordonn\'e car $\mathscr D$ est h\'er\'editairement belordonn\'e, on peut supposer qu'aucun \'el\'ement de cette anticha\^{i}ne ne soit dans  $\mathscr D. \mathcal A.$ %car s'il existe des elements dans $\mathscr D. \mathcal A.$ ils sont en nombre fini.
D'apr\'es le Corollaire \ref{cor:sumlex} et le th\'eor\`eme \ref{theo:decomposition}, pour tout $i\geq 0$ il existe une structure ind\'ecomposable $\mathcal S_i$ et des structures non vides $(\mathcal R_{ix})_{x\in S_i}$ telles que $\mathcal R_i= \underset{x\in \mathcal S_i}\oplus \mathcal R_{ix}$.  Puisque  $(\mathcal R_{ix},{f_i}_{\restriction_{\mathcal R_{ix}}})$ s'abrite strictement dans $(\mathcal R_i,f_i)$ on a
$(\mathcal R_{ ix},{f_i}_{\restriction_{\mathcal R_{ix}}})\in \mathcal I$ pour tout $i\geq 0$ et $x\in \mathcal S_i$. Puisque $\mathcal I$ est belordonn\'e et
$\mathscr D$ est h\'er\'editairement belordonn\'e alors
$\mathscr D.\mathcal I$ est belordonn\'e, donc la suite infinie $(\mathcal S_0,g_0),\cdots,(\mathcal S_i,g_i),\cdots$ de  $\mathscr D.\mathcal I$,
o\`u $g_i(x):=(\mathcal R_{ix},{f_i}_{\restriction_{\mathcal R_{ix}}})$,  contient une paire croissante $(\mathcal S_p,g_p)\leq (\mathcal S_q,g_q)$ pour deux entiers  $p<q$ (Th\'eor\`eme \ref{theo:higman-equivalence}). Ce qui signifie qu'il existe un abritement $h:\mathcal S_p \rightarrow \mathcal S_q$ tel que  $g_p(x)\leq g_q(h(x))$ pour tout $x\in S_p$, c'est \`a dire
$(\mathcal R_{px},{f_p}_{\restriction_{\mathcal R_{px}}})\leq (\mathcal R_{qh(x)},{f_q}_{\restriction_{\mathcal R_{qh(x)}}})$ pour tout $x\in S_p$. Il s'ensuit que   $(\mathcal R_p,f_p)=\underset{x\in \mathcal S_p}\oplus (\mathcal R_{px},{f_p}_{\restriction_{\mathcal R_{px}}})\leq \underset{x\in \mathcal S_q}\oplus
(\mathcal R_{qx},{f_q}_{\restriction_{\mathcal R_{qx}}})=(\mathcal R_q,f_q)$
ce qui contredit le fait que $\mathscr R$ est une  anticha\^{i}ne.

Donc $(\sum \mathscr D). \mathcal A$ est belordonn\'e et par suite $(\sum \mathscr D)$ est h\'er\'editairement belordonn\'e.
\end{proof}

\bigskip

 La Proposition \ref{prop:wqo} est  particuli\`erement vraie si $\mathscr{D}$ est finie. Si $\mathscr{D}$ est la classe $Ind_{\ell}(\Omega_k)$ de structures ind\'ecomposables de $\Omega_k$ de tailles au plus $\ell$ alors d'apr\`es un r\'esultat de Schmerl et Trotter \cite{S-T} (voir Th\'eor\`eme \ref{theo:schmerl-trotter}), les bornes de $\sum Ind_{\ell}(\Omega_k)$ sont de tailles au plus $\ell+2$.  Lorsque $\mathscr{D}$ est form\'ee de bicha\^{i}nes, la Proposition \ref{prop:wqo} a \'et\'e  obtenue par Albert et Atkinson \cite {A-A}.\\

Un corollaire imm\'ediat est le suivant:

\begin{corollary}\label{cor:wqo}
Si une classe h\'er\'editaire de $\Omega_k$ contient un nombre fini de membres ind\'ecomposables alors elle est belordonn\'ee et a un nombre fini de bornes.
\end{corollary}

%*****************************************************************************

    \section{Structures binaires ordonn\'ees}
    Dans toute la suite, nous consid\'erons uniquement les structures ordonn\'ees compos\'ees de relations binaires r\'eflexives. Soit $\Gamma_{k}$ la classe de telles structures finies de type $k$. %Soit $\mathscr A$ une sous-classe, pour $i\in \mathbb{N}$ soit $\mathscr A_{(i)}$, resp. $\mathscr A_{(\geq i)}$, les sous-classes form\'ees de ses membres qui ont $i$ \'el\'ements, resp. au moins $i$ \'el\'ements.

 \vspace{1mm}

Toutes les structures  consid\'er\'ees ici sont des structures qui  repr\'esentent des types d'isomorphie. De mani\`ere g\'en\'erale, si $\mathcal R:=(E,(\rho_i)_{i\in I})$ est une structure o\`u $E$ est un $n$-ensemble, alors son type d'isomorphie peut-\^etre repr\'esent\'e par une structure $\mathcal R'$ isomorphe \`a $\mathcal R$ et d\'efinie sur $\{0,1,\cdots,n-1\}$ et si $\mathcal R$ est ordonn\'ee, l'ordre sur $\mathcal R'$ est l'ordre naturel.

\vspace{1mm}

\noindent Si, dans le %corollaire \ref{cor:sumlex} a au moins deux \'el\'ements, (ou si $\mathcal S$ dans le
 Th\'eor\`eme \ref{theo:decomposition}, $\mathcal S$ est encha\^{\i}nable, alors, comme soulign\'e dans la Remarque \ref{rem:dec-enchai}, la d\'ecomposition n'est pas n\'ecessairement unique, un fait qui est \'egalement d\^{u} \`{a} la notion d'intervalle fort. Nous nous int\'eressons \`a ce cas l\`a pour des structures binaires ordonn\'ees. Pour nos besoins, nous introduisons la notion suivante.\\

 D\'esignons par $\textbf{1}$ la structure \`a un \'el\'ement. Soit $\tau$ une structure de $\Gamma_k$ ayant deux \'el\'ements. Une structure ordonn\'ee $\mathcal S$ est dite $\tau$-\emph{ind\'ecomposable} si elle ne peut pas se d\'ecomposer en une somme lexicographique index\'ee par $\tau.$
Si $\mathscr A$ est une classe de structures de $\Gamma_k$, on d\'esigne par $\mathscr A(\tau)$ l'ensemble des membres de $\mathscr A$ qui sont $\tau$-ind\'ecomposables.

\begin{lemma}\label{lem:decomposition}
Soit  $\mathcal S:=(\{0,1\},\leq,(\rho_i)_{i\in [k]})$ avec  $0<1$, une structure ordonn\'ee de type $k$ ayant deux \'el\'ements. Si $\mathcal R$ est une somme lexicographique $\underset{x\in \mathcal{S}}{\oplus }\mathcal{R}_x=\mathcal R_0\underset{\mathcal{S}}{\oplus}\mathcal R_1$ et $\mathcal R_0$ est $\mathcal S$-ind\'ecomposable, alors la partition
$\mathcal R_0,\mathcal R_1$ est unique.
\end{lemma}

\begin{proof}
Evident. Si $\mathcal R_0$ n'est pas $\mathcal S$-ind\'ecomposable nous aurons une autre partition issue de la d\'ecomposition de $\mathcal R_0$ suivant $\mathcal S.$ Mais puisque $\mathcal R_0$ est $\mathcal S$-ind\'ecomposable, il n'existe pas d'autre d\'ecomposition dans laquelle cette hypoth\`ese est v\'erifi\'ee.
\end{proof}
\medskip

Soit $\mathscr C$ une sous-classe de $\Gamma_k$. Pour $i\in \mathbb{N}$ posons $\mathscr C_{(i)}$, (resp. $\mathscr C_{(\geq i)}$), les sous-classes form\'ees des membres de $\mathscr C$ qui ont $i$ \'el\'ements, (resp. au moins $i$ \'el\'ements).
%\medskip

\begin{lemma}\label{lem:union}
Soit $\mathscr D$ une classe de structures ind\'ecomposables non vides de $\Gamma_k$ contenant la structure \`{a} un \'el\'ement ${\bf 1}$. Soit $\mathscr A=\sum\mathscr D$ la cl\^oture par sommes de $\mathscr D$ et pour tout $\mathcal S\in  \mathscr D_{(2)}$, notons par $\mathscr A (\mathcal S)$ la sous-classe des membres
 $\mathcal S$-ind\'ecomposables de $\mathscr A$.\\ %En \'ecrivant chaque  $\mathcal S\in  \mathscr D_{(2)}$ comme $\mathcal S:=(\{0,1\},\leq,(\rho_i)_{i\in J})$ avec $0<1$, nous avons le r\'esultat suivant:
  Posons $\mathscr A_{\mathcal S}:= \underset{ \mathcal{S}}{ \oplus}\mathscr{A}$ si $\mathcal S\in  \mathscr D_{( \geq3)}$ et $\mathscr A_{\mathcal S}:=\mathscr A(\mathcal S)\underset{\mathcal S}{\oplus}\mathscr A$ si $\mathcal S\in  \mathscr D_{(2)}$ et \'ecrivons chaque $\mathcal S\in  \mathscr D_{(2)}$ comme  $\mathcal S:=(\{0,1\},\leq,(\rho_i)_{i\in [k]})$ avec $0<1$. Alors:

\begin{equation} \label{eq:classe}
\mathscr A=\{\bf 1\}\cup \underset{\mathcal S\in \mathscr D_{(\geq 2)}}{\bigcup}\mathscr A_{\mathcal S}\end{equation}
et
\begin{equation}\label{eq:classeindec}
\mathscr A(\mathcal S)=\mathscr A \setminus \mathscr A_{\mathcal S}
\end{equation}
pour tout $\mathcal S \in \mathscr D_{(2)}$.

De plus, tous les ensembles de l'\'equations \eqref{eq:classe} sont deux \`{a} deux disjoints.
\end{lemma}

\begin{proof}
D\'esignons par $(1)$ (respectivement par $(2)$) le terme de gauche (respectivement de droite) de l'\'equation \eqref{eq:classe}. L'inclusion $(2)\subseteq (1)$ est \'evidente car $\mathscr A$ est ferm\'e par sommes d'apr\`es le Th\'eor\`eme \ref{theo:closed}.
Pour montrer l'inclusion $(1)\subseteq (2)$, soit $\mathcal R$ dans $(1)$, si $\mathcal R$ a un \'el\'ement alors elle est dans $(2)$, sinon, d'apr\`es le  Th\'eor\`eme \ref{theo:decomposition}, $\mathcal R$ est une somme  lexicographique $\underset{x\in \mathcal S}\oplus \mathcal R_x$ o\`{u} $\mathcal S$ est soit ind\'ecomposable avec au moins trois  \'el\'ements soit une structures binaire encha\^{\i}nable avec au moins deux \'el\'ements et pour chaque
$\mathcal R_x\in \mathscr A$, le domaine $V(\mathcal R_x)$  est un intervalle fort de $\mathcal R$ pour tout $x\in V(\mathcal S)$.  Dans le premier cas, $\mathcal R$ est dans $\underset{\mathcal R'\in \mathscr D_{(\geq  3)}} \bigcup \underset{x\in \mathcal R'}\oplus \mathscr A_x$ donc dans $(2).$
Dans le second cas,  $\mathcal S$ est encha\^{\i}nable avec $n$ \'el\'ements ($n\geq 2$) nous pouvons alors poser
$\mathcal S:=(\{0,1,\cdots,n-1\},\leq,(\rho_i)_{i\in [k]})$ avec $0<1<\cdots<n-1$.
Posons $\mathcal S':=\mathcal S_{ \restriction_{\{0,1\}}}$, $\mathcal S'':=\mathcal S_{ \restriction_{\{1,\cdots,n-1\}}}$, $\mathcal R'_0:=\mathcal R_0$ et
$\mathcal R'_1:= \underset{x\in \mathcal S''}\oplus \mathcal R_x$. Nous avons \'evidemment $\mathcal R= \mathcal R'_0 \underset{\mathcal S'}\oplus \mathcal R'_1$.\\
$\mathcal R_0$  est $\mathcal S'$-ind\'ecomposable car  $V(\mathcal R_0)$ est un intervalle fort de $\mathcal R$, ; %puisque $\mathcal R'_0:=\mathcal R_0$,
ainsi $\mathcal R$ appartient \`{a} $\underset{\mathcal S\in \mathscr D_{(2)}}{\bigcup}(\mathscr A_0(\mathcal S)\underset{\mathcal S}{\oplus}\mathscr A_1)$ qui est un sous-ensemble de $(2).$

 Le fait que ces ensembles soient deux \`{a} deux disjoints vient du Corollaire \ref{cor:sumlex} et du Lemme \ref{lem:decomposition}.

 L'\'equation \eqref{eq:classeindec} est \'evidente; les membres $\mathcal S$-ind\'ecomposables de $\mathscr A$ sont ceux qui ne peuvent pas s'\'ecrire comme des  $\mathcal S$-sommes.
\end{proof}

\bigskip

Soit $\{\mathcal S_1,\cdots,\mathcal S_p\}$  une \'enum\'eration des structures de $\mathscr D_{(2)}$ (structures \`{a} $2$ \'el\'ements de $\mathscr D$).  Soit $\mathcal H$, resp. $\mathcal H_i$, resp.  $\mathcal K$, les s\'eries g\'en\'eratrices de $\mathscr{A}$, resp. $\mathscr {A}(\mathcal S_i)$, resp. $\mathscr{D}_{(\geq 3)}$ et soit $\mathcal {K(H)}$ la s\'erie obtenue \`a partir de $\mathcal K$ en substituant \`a l'inconnue $x$ la s\'erie $\mathcal H$.

\vspace{1mm}

A partir du Lemme \ref{lem:union} ci-dessus, nous d\'eduisons

\begin{lemma}\label{lem:seriesclasses}
\begin{equation}\label{eq:serieclasse}
 (p-1)\mathcal {H}^2+(x-1+\mathcal K(\mathcal{H})){\mathcal {H}}+x+\mathcal K(\mathcal{H})=0.
 \end{equation}

 \begin{equation}\label{eq:serieclasseindec}
\mathcal H_i=\frac{\mathcal H}{1+\mathcal H}\; ~~\text{pour tout}\;  i\in \{1,\cdots,p\}.
\end{equation}
\end{lemma}

\begin{proof}
 Montrons d'abord l'\'equation \eqref{eq:serieclasseindec}. Posons $\mathcal H_{\mathcal S}$ la s\'erie g\'en\'eratrice de $\mathscr A_{\mathcal S}$ pour $\mathcal S\in \mathscr{D}_{(\geq 2)}$. Soit $\mathcal S\in \mathscr{D}_{( 2)}$. D'apr\`es les notations du Lemme \ref{lem:union},  $\mathscr A_{\mathcal S}=\mathscr A(\mathcal S)\underset{\mathcal S}{\oplus}\mathscr {A}$, nous avons alors $\mathcal H_{\mathcal S}=\mathcal H_i.\mathcal H$. A partir de l'\'equation (\ref{eq:classeindec}) nous d\'eduisons que $\mathcal H_i= \mathcal H-\mathcal H_i.\mathcal H$ pour tout $i\in \{1, \dots , p\}$. Or les coefficients de $\mathcal H$ sont non n\'egatifs, donc la s\'erie  $1+ \mathcal H$ est inversible, ainsi $\mathcal H_i=\dfrac{\mathcal H}{1+\mathcal H}$. D'où l'\'equation \eqref{eq:serieclasseindec}.

Montrons maintenant l'\'equation \eqref{eq:serieclasse}. Soit $\mathcal S\in  \mathscr{D}_{(n)}$ pour $n\geq 3$ alors, toujours d'apr\`es le Lemme \ref{lem:union}, $\mathscr A_{\mathcal S}=\underset{\mathcal S}{\oplus}\mathscr {A}$ et donc $\mathcal H_{\mathcal S}=\mathcal H^n$. Nous d\'eduisons alors que la s\'erie g\'en\'eratrice de $\underset{\mathcal S\in \mathscr D_{(\geq 3)}}{\bigcup}\mathscr A_{\mathcal S}$ est \'egale \`a $\mathcal K(\mathcal{H})$. D'apr\`es ce qui pr\'ec\`ede  la s\'erie g\'en\'eratrice de
$\underset{\mathcal S\in \mathscr D_{( 2)}}{\bigcup}\mathscr A_{\mathcal S}$ est \'egale \`a $\underset{i=1}{\overset{p}{\sum }}\mathcal H_i.\mathcal H$, donc, en utilisant l'\'equation \eqref{eq:serieclasseindec}, elle vaut $p\dfrac{\mathcal H^2}{1+\mathcal H}$.

En substituant ces r\'esultats dans l'\'equation \eqref{eq:classe} nous obtenons:
\begin{equation}\label{equation:3}\mathcal H=x+p\frac{\mathcal H^2}{1+\mathcal H}+\mathcal K(\mathcal{H}).
\end{equation}
Un calcul direct permet de trouver l'\'equation \eqref{eq:serieclasse}.

\noindent (Remarquons que nous pouvons faire explicitement les calculs \`a partir du Lemme \ref{lem:union}, l'\'equation \eqref{eq:classe} donne:
\begin{equation}
\begin{array}{rcl}
%\label{eq}
\underset{n\geq 1}\sum \varphi_{\mathscr A}(n)x^n&=&x+\underset{i=1}{\overset{p}{\sum }}\underset{n\geq 2}{\sum }\left[ \overset%
{n-1}{\underset{l=1}{\sum }}\varphi_{\mathscr A(\mathcal S_i)}(l).\varphi_{\mathscr A}(n-l)\right]x^n\\
 && +\: \underset{n\geq 3}{\sum }\left[\overset{n}{\underset{l=3}{\sum }}\;
 \underset{n_1+\cdots+n_l=n}{\sum }\varphi_{\mathscr D_{(\geq 3)}}(l)\varphi_{\mathscr A}(n_1)\cdots \varphi_{\mathscr A}(n_l)\right]x^n.
\end{array}
\label{eq}
\end{equation}
En utilisant les notations donn\'ees ci-dessus, nous obtenons l'\'equation \eqref{equation:3}).

\end{proof}

%*******************************subsection********************************************

\section{Classes alg\'ebriques et h\'er\'editairement alg\'ebriques}
Disons qu'une classe de structures finies est alg\'ebrique si sa s\'erie g\'en\'eratrice est alg\'ebrique.

\noindent Un corollaire imm\'ediat du Lemme \ref{lem:seriesclasses} (analogue du corollaire 13 de \cite{A-A}) est le suivant:

\begin{corollary}\label{cor:wqoalgebraic}
Soit $\mathscr D$ une classe de structures non vides et ind\'ecomposables contenant la structure \`{a} un \'el\'ement ${\bf 1}$. Si $\mathscr D$ est alg\'ebrique alors sa cl\^oture par sommes ainsi que chaque sous-classe form\'ee des membres $\mathcal S$-ind\'ecomposables de cette cl\^oture, avec $\mathcal S\in \mathscr D_{(2)}$, sont alg\'ebriques.
\end{corollary}

 Nous disons qu'une classe $\mathscr{C}$ de structures relationnelles  est \emph{h\'er\'editairement rationnelle}\index{classe!h\'er\'editairement rationnelle}, resp. \emph{h\'er\'editairement alg\'ebrique}\index{classe!h\'er\'editairement alg\'ebrique} si la fonction g\'en\'eratrice\index{fonction!g\'en\'eratrice} de toute sous-classe h\'er\'editaire de $\mathscr{C}$ est rationnelle, resp. alg\'ebrique.
Albert, Atkinson et Vatter\index{Albert, Atkinson et Vatter} \cite{A-A-V} ont montr\'e que les classes h\'er\'editairement rationnelles de permutations sont  belordonn\'ees. Ce fait peut-\^etre \'etendu  aux classes h\'er\'editairement  alg\'ebriques.

\begin{lemma} \label{lem:wqo}
Une classe h\'er\'editaire $\mathscr{C}$ qui est h\'er\'editairement alg\'ebrique est belordonn\'ee.
\end{lemma}

\begin{proof} Si  $\mathscr{C}$ contient une anticha\^{i}ne infinie, il existe un nombre non d\'enombrable de sous-classes  h\'er\'editaires de  $\mathscr{C}$ et en fait une cha\^{i}ne non d\'enombrable de sous-classes; ces classes fournissent un nombre non d\'enombrable de s\'eries g\'en\'eratrices. Certaines de ces s\'eries ne peuvent pas \^etre alg\'ebriques. En effet, d'apr\`es C. Retenauer\index{Retenauer} \cite{C-R}, une s\'erie g\'en\'eratrice avec des coefficients rationnels qui est alg\'ebrique sur $\mathbb C$ est alg\'ebrique sur $\mathbb Q$. Comme les s\'eries g\'en\'eratrices consid\'er\'ees ont des coefficients entiers, elles sont alg\'ebriques sur $\mathbb Q$, donc le nombre de ces s\'eries est au plus d\'enombrable.
\end{proof}
\bigskip

Si $\mathscr C$ et $\mathscr C'$ sont deux classes h\'er\'editaires, alors leurs s\'eries g\'en\'eratrices v\'erifient l'identit\'e
$\mathcal H_{\mathscr C\cup \mathscr C'}= \mathcal H_{\mathscr C}+\mathcal H_{\mathscr C'}-\mathcal H_{\mathscr C\cap \mathscr C'}$. A partie de cette \'equation simple nous avons:

\begin{lemma}\label{lemme:union here-rati}
L'union de deux classes h\'er\'editairement rationnelles (resp. h\'er\'editairement alg\'ebriques) est h\'er\'editairement rationnelle (resp. h\'er\'editairement alg\'ebrique).
\end{lemma}

De ce lemme nous pouvons d\'eduire
\begin{corollary}
Une classe $\mathscr C$ non h\'er\'editairement rationnelle ou non h\'er\'editairement alg\'ebrique qui est minimale pour cette propri\'et\'e est %n\'ecessairement
un id\'eal.
\end{corollary}

\begin{proof}
D'apr\`es le Lemme \ref{lemme:union here-rati}, $\mathscr C$ ne peut pas \^etre une r\'eunion de deux sous-classes h\'er\'editaires propres, c'est alors un id\'eal, donc un \^age.
\end{proof}

%qu'une classe non h\'er\'editairement rationnelle ou non h\'er\'editairement alg\'ebrique qui est minimale pour cette propri\'et\'e est n\'ecessairement un id\'eal. %c'est d\^u \`a la minimalit\'e, car sinon, elle serait une union d'id\'eaux, qui sont tous h\'er\'editairement alg\'ebriques resp. rationnels
\vspace{1mm}

Nous \'enon\c{c}ons le th\'eor\`eme principal de cette partie, qui g\'en\'eralise le th\'eor\`eme d'Albert et Atkinson \cite{A-A}:

\begin{theorem}\label{theo:algebraic}
 Soit $k$ un entier. Si une classe h\'er\'editaire $\mathscr C$ de $\Gamma _{k}$ contient uniquement un nombre fini de membres ind\'ecomposables alors elle est alg\'ebrique.
\end{theorem}

\begin{proof}
Pour la preuve de ce th\'eor\`eme, nous suivons essentiellement les \'etapes de la preuve d'Albert-Atkinson.
Nous donnons une preuve inductive sur les sous-classes h\'er\'editaires de $\mathscr C$. Mais pour cela, nous avons besoin de montrer plus, notamment que $\mathscr C$ et chaque $\mathscr C(\mathcal S)$ pour $\mathcal S\in \mathscr D_{(2)}$, o\`u $\mathscr D:= Ind(\mathscr C)$,  sont alg\'ebriques (ceci est la seule diff\'erence avec la preuve d'Albert-Atkinson). Pour \'eviter des complications qui ne sont pas n\'ecessaires,  nous retirons la structure relationnelle vide de $\Gamma_k$, ce qui signifie que nous supposons que $\mathscr C$ est form\'ee de structures non vides.

Posons $\mathscr A:=\sum Ind(\mathscr{C})$. Nous avons alors deux cas:

\noindent\textbf{1)-} Si  $\mathscr C= \mathscr A$ alors d'apr\'es le Corollaire \ref{cor:wqoalgebraic},  $\mathscr C$ et
 chaque $\mathscr C(\mathcal S)$ pour $\mathcal S\in \mathscr D_{(2)}$  %o\`u $\mathscr D:= Ind(\mathscr C)$,
 sont alg\'ebriques, le r\'esultat est donc d\'emontr\'e.

\vspace{1mm}

\noindent\textbf{2)-} Si $\mathscr C\not = \mathscr A$, nous pouvons supposer que pour chaque sous-classe h\'er\'editaire propre $\mathscr C'$ de $\mathscr C$; $\mathscr C'$ et chaque $\mathscr C'(\mathcal S)$ pour $\mathcal S\in \mathscr D'_{(2)}$,  o\`u $\mathscr D':= Ind(\mathscr C')$,
sont alg\'ebriques. En effet, si ceci n'est pas v\'erifi\'e, sachant que  $\mathscr C$ est belordonn\'ee (Corollaire \ref{cor:wqo}) elle contient une sous-classe h\'er\'editaire  minimale qui ne v\'erifie pas cette propri\'et\'e et nous pouvons remplacer $\mathscr C$ par cette sous-classe.

\vspace{1mm}
%\noindent La derni\`ere \'etape de la preuve du Th\'eor\`eme \ref{theo:algebraic} est comme suit.

\noindent Soit  $\mathcal S\in \mathscr D_{(\geq 2)}$.  Si $V(\mathcal S)$ a deux \'el\'ements  $0$ et $1$, avec $0<1$, nous posons
 $\mathscr C_{\mathcal S}:=(\mathscr C(\mathcal S)\underset {\mathcal S }\oplus\mathscr C)\cap \mathscr C$. Si $V(\mathcal S)$ a au mois trois \'el\'ements, nous posons $\mathscr C_{\mathcal S}:= (\underset{ \mathcal{S}}{\oplus}\mathscr{C})\cap \mathscr C$.  Comme dans le Lemme \ref{lem:union} nous avons
\begin{equation}\label{eq:global}
\mathscr C= \{\bf1\} \cup \underset {\mathcal S\in \mathscr D_{(\geq 2)}}\bigcup \mathscr C_{\mathcal S}.
\end{equation}
et
\begin{equation}\label{eq:classeindecbis}
\mathscr C(\mathcal S)=\mathscr C\setminus\mathscr C_{\mathcal S}\; \text{ pour tout }\;  \mathcal S \in \mathscr D_{( 2)}.
\end{equation}

Soient $\mathcal H$, $\mathcal H_{\mathscr C_{(2)}}$ et $\mathcal H_{\mathscr C_{(\geq 3)}}$   les  s\'eries g\'en\'eratrices de
$\mathscr{C}$, %et $\mathscr D_{(\geq 3)}$
$\mathscr C_{(2)}:=\underset{\mathcal S\in \mathscr D_{(2)}} \bigcup {\mathscr C_\mathcal S}$ et   $\mathscr C_{(\geq 3)}:=\underset{\mathcal S\in \mathscr D_{(\geq 3)}} \bigcup {\mathscr C_\mathcal S}$ respectivement.  Notons par $\mathcal H_{\mathscr C(\mathcal S)}$ et $\mathcal H_{\mathscr C_{\mathcal S}}$ les s\'eries g\'en\'eratrices de $\mathscr C(\mathcal S)$ et $\mathscr C_{\mathcal S}$ respectivement. %$Ind(\mathfrak C)_{(\geq 3)}$ respectively.

Nous avons:
\begin{equation}\label{eq:global1}
\mathcal  H= x + \mathcal H_{\mathscr C_{(2)}}+\mathcal H_{\mathscr C_{(\geq 3)}}\end{equation}
et

\begin{equation}\label{eq:classeindecter}
\mathcal H_{\mathscr C(\mathcal S)}=\mathcal H - \mathcal H_{\mathscr C_{\mathcal S}} \; \text {pour tout}\; \mathcal S\in   \mathscr D_{(2)}.\end{equation}

 Nous d\'eduisons que  $\mathcal  H$  et $\mathcal H_{\mathscr C(\mathcal S)}$ sont alg\'ebriques pour tout  $\mathcal S \in \mathscr D_{(2)}$, \`a partir des Affirmations suivantes que nous d\'emontrons plus loin:

\begin{claim}\label{claim:claim1}
La s\'erie g\'en\'eratrice $\mathcal H_{\mathscr C_{(\geq 3)}}$ %de $\mathscr B:=\underset{\mathcal R\in \mathscr D_{(\geq 3)}} \bigcup {\mathscr C_\mathcal R}$
est un polyn\^ome en la s\'erie g\'en\'eratrice $\mathcal H$ de $\mathscr C$  dont les coefficients sont des s\'eries alg\'ebriques.
\end{claim}

 \begin{claim}\label{claim:claim2}
 Pour tout $\mathcal S\in \mathscr D_{(2)}$, la s\'erie g\'en\'eratrice $\mathcal H_{\mathscr C(\mathcal S)}$ de  ${\mathscr C(\mathcal S)}$ est soit un polyn\^ome en la s\'erie g\'en\'eratrice $\mathcal H$ de la forme
 \begin{equation} \label {eq:almostfinal2}
\mathcal H_{\mathscr C(\mathcal S)}=\frac{(1-\alpha)\mathcal H-\delta}{1+\beta};
\end{equation}
 dont les coefficients sont des s\'eries alg\'ebriques, soit une fraction rationnelle de la forme %quotient de deux polyn\^omes en la s\'erie g\'en\'eratrice $\mathcal H_{\mathscr C}$.
 \begin{equation} \label {eq:almostfinal1}
\mathcal H_{\mathscr C(\mathcal S)}=\frac{\mathcal H}{1+ \mathcal H}.
\end{equation}
 \end{claim}

En substituant dans la formule \eqref{eq:global1} les valeurs de $\mathcal H_{\mathscr C_{(\geq 3)}}$ et $\mathcal H_{\mathscr C_{( 2)}}$ par celles donn\'ees dans l'Affirmation \ref{claim:claim1} et l'Affirmation \ref{claim:claim2} nous obtenons une \'equation polynomiale $C=0$ en la variable
$\mathcal{H}$ dont les coefficients sont des s\'eries alg\'ebriques. Le polyn\^ome $C$ n'est pas identiquement nul. En effet, $C$ est la somme de deux  polyn\^omes $A=a_0+a_1\mathcal H+a_2\mathcal H^2$ et $B=b_0+b_1\mathcal H+\cdots+b_k\mathcal H^k$ dont les coefficients sont des s\'eries alg\'ebriques (en fait, $B=\mathcal H_{\mathscr C_{(\geq 3)}}(1+\mathcal H)$).  Les valuations de $A$ et $B$, en tant que s\'eries en $x$, sont distinctes (on rappelle que la valuation d'un polyn\^ome non nul est le plus petit degr\'e de ses mon\^omes). En effet, la  valuation de $A$ est $1$ (notons que $a_0=x+\delta$ où $\delta$ est soit \'egal \`a zero soit une  s\'erie algebrique de valuation au moins $2$).  Ainsi, si $B\neq 0$ (lorsque $\mathscr D_{(\geq 3)}$ est non vide) sa valuation est au moins $3$. %$A+B$ is non zero, hence $\mathcal H_{\mathfrak C}$ is algebraic. If $B\neq 0$, the valuation of $B$, as a series in $x$ is $3$ at least $3$. The valuation of $A$, as polynomial in $x$ is at least $1$.
 Puisque $A$ et $B$ n'ont pas la même valuation, alors $A+B$ n'est pas identiquement nul.   Etant une racine d'un polyn\^ome non nul, $\mathcal H$ est alg\'ebrique. Avec ce r\'esultat et l'Affirmation \ref{claim:claim2}, $\mathcal H_{\mathscr C(\mathcal S)}$ est alg\'ebrique.
Avec ceci, la preuve du th\'eor\`eme est compl\`ete.
 \end{proof}

\vspace{2mm}

Pour la preuve de nos deux Affirmations pr\'ec\'edentes, nous avons besoin des lemmes suivants (respectivement Lemme 15 et Lemme 18 dans \cite{A-A}).

\begin{lemma}\label{lem:intersection}
 Soit $\mathcal R$ une structure binaire ordonn\'ee ind\'ecomposable et $\mathscr A:=(\mathscr A_x)_{x\in \mathcal R}$, $\mathscr B:=(\mathscr B_x)_{x\in \mathcal R}$ deux suites de sous-classes de structures binaires ordonn\'ees index\'ees par les \'el\'ements de $V(\mathcal R)$, alors
  \begin{enumerate}
  \item Si $V(\mathcal R)$ a au moins trois \'el\'ements alors
    $$(\underset{x\in \mathcal R}{\oplus}{\mathscr A_x})\cap(\underset{x\in \mathcal R}{\oplus}{\mathscr B_x})=\underset{x\in \mathcal R}{\oplus}
    ({\mathscr A_x\cap \mathscr B_x}).$$
\item Si $\mathcal R:=(\{0,1\},\leq,(\rho_i)_{i\in J})$ avec $0<1$ alors
    $$(\mathscr A_0(\mathcal R)\underset{\mathcal R}{\oplus}\mathscr A_1)\cap (\mathscr B_0(\mathcal R)\underset{\mathcal R}{\oplus}\mathscr B_1)=(\mathscr A_0(\mathcal R)\cap \mathscr B_0(\mathcal R))\underset{\mathcal R}{\oplus}(\mathscr A_1\cap \mathscr B_1).$$
\end{enumerate}
\end{lemma}

\begin{proof}
La preuve vient directement de l'unicit\'e de la d\'ecomposition des structures (Corollaire \ref{cor:sumlex}, Th\'eor\`eme \ref{theo:decomposition} et Lemme \ref{lem:decomposition}).
\end{proof}

\bigskip

Pour le second lemme, nous introduisons les notions suivantes:

Soit  $\mathscr C$ une classe de structures finies et $\overline {\mathcal B} := \mathcal B_1,\mathcal B_2,...,\mathcal B_l$ une suite de structures finies, nous posons
 $\mathscr C <\overline {\mathcal B}>:=\mathscr C < \mathcal B_1,\mathcal B_2,...,\mathcal B_l>:=Forb(\{\mathcal B_1,\mathcal B_2,...,\mathcal B_l\})\cap \mathscr C.$

\begin{definition}
 Si $\mathscr C$ est h\'er\'editaire, une sous-classe h\'er\'editaire propre $\mathscr C'$ de $\mathscr C$ est dite \emph{forte} si toute borne de  $\mathscr C'$ dans $\mathscr C$ s'abrite dans au moins une borne de $\mathscr C.$
 \end{definition}
\noindent Notons que l'intersection de sous-classes fortes est forte.

\medskip

 Soit $\mathscr A:=(\mathscr A_x)_{x\in \mathcal S}$, o\`u $\mathcal S$ est ind\'ecomposable avec au moins trois \'el\'ements. Une \emph{d\'ecomposition d'une structure binaire} $\mathcal B$ \emph{suivant} $\mathscr A$ est une application $h:\mathcal B\rightarrow\mathcal S$ telle que
$\mathcal B=\underset{x\in \mathcal S_{ \restriction_{Im(h)}}}{\oplus}{\mathcal B}_
{\restriction_{h^{-1}(x)}}$ et $\mathcal B_{ \restriction_{h^{-1}(x)}} \in \mathscr A_x$ pour tout $x\in Im(h).$ Donc, chaque $\mathcal B_ {\restriction_{h^{-1}(x)}}$ est un intervalle de $\mathcal B.$ D\'esignons par  $H_{\mathcal B}$ l'ensemble de telles d\'ecompositions de $\mathcal B$.

\begin{lemma}\label{lem:egality}
Soient $\mathcal S$ une structure ordonn\'ee ind\'ecomposable, $\mathscr A:=(\mathscr A_x)_{x\in \mathcal S}$ une famille de sous-classes de structures binaires ordonn\'ees index\'ee par les \'el\'ements de $\mathcal S$,  $\overline {\mathcal B} := \mathcal B_1,\mathcal B_2,...,\mathcal B_l$ une suite de structures finies et $\mathscr C:=(\underset{x\in \mathcal S}{\oplus}{\mathscr A_x})<\overline {\mathcal B}>$. Si $\mathcal S$ a au moins trois \'el\'ements alors $\mathscr C$ est une r\'eunion d'ensembles de la forme
$\underset{x\in \mathcal S}{\oplus}\mathscr D_x$ o\`u chaque $\mathscr D_x$ est soit $\mathscr A_x<\overline {\mathcal B}>$ soit l'une de ses sous-classes fortes.
\end{lemma}

\begin{proof}
Nous d\'emontrons le r\'esultat pour $l=1$. Pour cela nous montrons que:
\begin{equation} \label{eq:eqstrong}
\mathscr C=\underset{h\in H_{\mathcal B}}{\bigcap}\;\underset{x\in Im(h)}{\bigcup}\underset{y\in \mathcal S}{\oplus}{\mathscr A^{(x)}_y}
\end{equation}
 o\`u $\mathscr A^{(x)}_x:=\mathscr A_x<\mathcal B_{\restriction_{h^{-1}(x)}}>~$ et $~\mathscr A^{(x)}_y:=\mathscr A_y~$ pour $y\neq x$.

D\'esignons par $(1)$ (respectivement par $(2)$) le membre de gauche (respectivement de droite) de l'\'equation \eqref{eq:eqstrong}. L'inclusion $(1)\subseteq (2)$ a lieu sans aucune hypoth\`ese. En effet, soit $\mathcal T$ dans $(1)$. Nous montrons que $\mathcal T$ est dans $(2)$. Si $h$ est une d\'ecomposition de $\mathcal B$, nous voulons trouver un $x\in Im(h)$ tel que $\mathcal T\in \underset{y\in \mathcal S}{\oplus}\mathscr A^{(x)}_y.$
Puisque $\mathcal T$ est dans $(1)$, elle  a une d\'ecomposition suivant $\mathcal S.$  Soit $h\in H_{\mathcal B}$, puisque $\mathcal B \nleq \mathcal T$, il existe $x\in \mathcal S_{\restriction_{Im(h)}}$ tel que $\mathcal B_{\restriction_{h^{-1}(x)}}\nleq\mathcal T$ (car si $\mathcal B \leq \mathcal T=\underset{x\in \mathcal S}{\oplus}\mathcal T_x$ o\`u $\mathcal T_x\in\mathscr A_x,~\forall x\in\mathcal S$ alors il existe $h\in H_{\mathcal B}$ tel que $\mathcal B_{\restriction_{h^{-1}(x)}}\leq\mathscr A_x$ pour tout $x\in Im(h)$), donc
$\mathcal T\in \underset{y\in \mathcal S}{\oplus}\mathscr A^{(x)}_y.$

 L'inclusion $(2)\subseteq (1)$ a lieu sous l'hypoth\`ese que la structure dans $(2)$ a une unique d\'ecomposition suivant $\mathcal S$ %du a l'unicite de la decomposition des structures ordonnee
et qu'elle est ordonn\'ee (ce qui signifie que $\mathcal S$ est rigide, c'est \`a dire que $\mathcal S$
n'a aucun automorphisme distinct de l'identit\'e).
Soit $\mathcal T$ dans $(2)$, alors pour tout $h\in H_{\mathcal B}$ il existe $x_h\in Im (h)$ tel que $\mathcal T\in \underset{y\in \mathcal S}{\oplus}{\mathscr A^{(x_h)}_y}$ donc $\mathcal T\in \underset{h\in H_{\mathcal B}}{\bigcap}\underset
{y\in \mathcal S}{\oplus}{\mathscr A^{(x_h)}_y}.$
Nous avons $\mathcal T\in \underset{y\in \mathcal S}{\oplus}{\mathscr A_y}$ car $\mathscr A^{(x_h)}_y\subseteq \mathscr A_y$ pour tout $h.$ Donc $\mathcal T=\underset{y\in \mathcal S}{\oplus}\mathcal T_y.$
 Nous avons $\mathcal B\nleq \mathcal T$. En effet, supposons que  $\mathcal B\leq \mathcal T$ et
soit $f$ un abritement de $\mathcal B$ dans $\mathcal T$ et  $h:=pof$, o\`u $p$ est la projection de $\mathcal T$ dans $\mathcal S,$ nous devons avoir $\mathcal B_{\restriction_{h^{-1}(x)}}\leq \mathcal T_x$ pour tout $x\in Im(h)$ ce qui constitue une contradiction avec le fait que
$\mathcal T\in \underset{h\in H}{\bigcap}\underset{y\in \mathcal S}{\oplus}{\mathscr A^{(x_h)}_y}.$ Donc $\mathcal T$ est dans $(1).$

\noindent En utilisant la distributivit\'e de l'intersection par rapport \`a la r\'eunion, nous pouvons \'ecrire $(2)$ comme une r\'eunion de termes, dont chacun est l'intersection de termes de la forme $\underset{y\in \mathcal S}{\oplus}\mathscr A_y<\mathcal B_x>,$ o\`u $\mathcal B_x$ est un intervalle de $\mathcal B$ tel que, il existe une d\'ecomposition $h$ de $\mathcal B$ et  $\mathcal B_x=\mathcal B_{ \restriction_{h^{-1}(x)}}.$\\
 Ces intersections, par le Lemme \ref{lem:intersection} et le fait que, parmi toutes les d\'ecompositions de $\mathcal B$ nous avons celles qui envoient  $\mathcal B$ vers un \'el\'ement $x$ de $V(\mathcal S)$ (c'est \`a dire celles pour lesquelles $h$ est constante), ont la forme $\underset{x\in \mathcal S}{\oplus}{\mathscr D_x}$ o\`u chaque $\mathscr D_x$ est de la forme $\mathscr A_x<\mathcal B,\cdots>$ (où les structures qui se trouvent apr\`es $\mathcal B$ -s'il en existe- sont des intervalles de $\mathcal B.$)\\ Par cons\'equent, $\mathscr D_x$ est soit la classe $\mathscr A_x<\mathcal B>$ soit l'une de ses sous-classes fortes.
Le cas $l>1$ s'ensuit par r\'ecurrence.

\end{proof}

\vspace{2mm}

Dans le cas o\`u $\mathcal S$ a deux \'el\'ements le Lemme \ref{lem:egality} devient:

\begin{lemma}\label{lem:twoelements}
Si $\mathcal S:=(\{0,1\},\leq,(\rho_i)_{i\in J})$ avec $0<1$ alors  $(\mathscr A(\mathcal S)\underset{\mathcal S}{\oplus}\mathscr A)<\overline{\mathcal B}>$ est une union de classes de la forme $(\mathscr {A'}(\mathcal S)<\overline{\mathcal B}>)\underset{\mathcal S}{\oplus}(\mathscr {A''}<\overline{\mathcal B}>)$, o\`u
$\mathscr A'<\overline{\mathcal B}>$ et $\mathscr A''<\overline{\mathcal B}>$ sont soit \'egales \`a $\mathscr A<\overline{\mathcal B}>$ soit \`a une sous-classe forte de
$\mathscr A<\overline{\mathcal B}>$.
\end{lemma}

\begin{proof} Comme pour le Lemme \ref{lem:egality}, nous supposons $l=1$.
L'\'egalit\'e \eqref{eq:eqstrong} devient dans ce cas:
\begin{equation} \label{eq:eqbinary}
\big(\mathscr A(\mathcal S)\underset{\mathcal S}{\oplus}\mathscr A\big)<\mathcal B>=\underset{h\in H_{\mathcal B}}{\bigcap} \left[\big(\mathscr A(\mathcal S)<\mathcal B_{\restriction_{h^{-1}(0)}}>
\underset{\mathcal S}{\oplus}\mathscr A\big)\bigcup \big(\mathscr A(\mathcal S)\underset{\mathcal S}{\oplus}\mathscr A<\mathcal B_{\restriction_{h^{-1}(1)}}>\big)\right]
\end{equation}
Le cas $l>1$ se d\'eduit par r\'ecurrence.
\end{proof}

\vspace{3mm}

\noindent{\bf Preuve de l'Affirmation \ref{claim:claim1}.}
%\begin{proofclaim}
 Puisque l'ensemble $\mathscr  D:=Ind (\mathscr C)$ est fini, il est h\'er\'editairement belordonn\'e, donc $\mathscr  A=\sum \mathscr{D}$ est belordonn\'e (Proposition \ref {prop:wqo}).
   Comme $\mathscr C$ est une sous-classe h\'er\'editaire propre de $\mathscr  A$, nous avons $\mathscr C=\mathscr A < \overline{\mathcal B}>:=\mathscr A<\mathcal B_1,\mathcal B_2,...,\mathcal B_l>$ pour une famille finie
$\overline{\mathcal B}:=\mathcal B_1,\mathcal B_2,...,\mathcal B_l$ de $\mathscr A$.

Soit $\mathcal S\in \mathscr D_{(\geq 3)}$,  le Lemme \ref{lem:egality} assure que
$\mathscr C_{\mathcal S}$ est une r\'eunion de classes, non n\'ecessairement disjointes, de la forme $\underset{x\in \mathcal S}{\oplus}\mathscr C_x$ o\`u chaque $\mathscr C_x$ est soit $\mathscr C$ soit l'une de ses sous-classes fortes.
La s\'erie g\'en\'eratrice de $\underset{x\in \mathcal S}{\oplus}\mathscr C_x$ est un mon\^ome en la s\'erie g\'en\'eratrice $\mathcal H$ de  $\mathscr C$ dont le coefficient est le produit des s\'eries g\'en\'eratrices des sous-classes fortes propres de $\mathscr C$ et par cons\'equent ce coefficient est une s\'erie alg\'ebrique (d'apr\`es l'hypoth\`ese faite dans la preuve du th\'eor\`eme \ref{theo:algebraic}). En utilisant le principe d'inclusion-exclusion, nous obtenons que la s\'erie g\'en\'eratrice $\mathcal H_{\mathscr C_{\mathcal S}}$ de $\mathscr C_{\mathcal S}$ est un polyn\^ome en la s\'erie g\'en\'eratrice $\mathcal H$  dont les coefficients sont des s\'eries alg\'ebriques. Puisque les classe $\mathscr C_{\mathcal S}$ sont deux \`a deux disjointes, la s\'erie g\'en\'eratrice $\mathcal H_{\mathscr C_{(\geq 3)}}$ %de $\mathscr B:=\underset{\mathcal R\in \mathscr D_{(\geq 3)}} \bigcup {\mathscr C_\mathcal R}$
est aussi un polyn\^ome en la s\'erie g\'en\'eratrice $\mathcal H$  dont les coefficients sont des s\'eries alg\'ebriques.
%\end{proofclaim}

 \hfill $\Box$

\vspace{2mm}

\noindent{\bf Preuve de l'Affirmation \ref{claim:claim2}.}
 %\begin{proofclaim}
 Soit $\mathcal S\in \mathscr D_{(2)}$,  le Lemme \ref{lem:twoelements} assure que $\mathscr C_{\mathcal S}$ est une r\'eunion, non n\'ecessairement disjointes, de classes de la forme ${\mathscr C'(\mathcal S)}\underset{\mathcal S}{\oplus}\mathscr C''$, o\`u chacune des classes $\mathscr C'$ et $\mathscr C''$ est soit \'egale \`a  $\mathscr C$ soit \'egale \`a l'une des sous-classes fortes de $\mathscr C$.
Les s\'eries g\'en\'eratrices de ces classes sont ou bien de la forme $\mathcal H_{\mathscr C(\mathcal S)}\mathcal H$ (si toutes les bornes $\mathcal B_i$ sont $\mathcal S$-ind\'ecomposables) ou bien de la forme $\alpha \mathcal H$ ou $\beta \mathcal H_{\mathscr C(\mathcal S)}$, avec $\alpha$ et $\beta$  des s\'eries alg\'ebriques (s\'eries des sous-classes fortes de $\mathscr C$) de valuation au moins 2 en $x$. En utilisant le principe d'inclusion-exclusion,  la s\'erie g\'en\'eratrice $\mathcal H_{\mathscr C_{\mathcal S}}$ obtenue est soit de la forme
$\mathcal H_{\mathscr C(\mathcal S)}\mathcal H$ (dans le cas o\`u toutes les bornes $\mathcal B_i$ sont $\mathcal S$-ind\'ecomposables) soit de la forme $\alpha \mathcal H+ \beta \mathcal H_{\mathscr C(\mathcal S)}+\delta$, avec $\alpha, \beta$ et $\delta$  des s\'eries alg\'ebriques,  $\alpha$ et $\beta$ \'etant de valuations au moins $2$. En particulier  $\mathcal H_{\mathscr C_{\mathcal S}}$ a la forme $\alpha_{\mathcal S} \mathcal H_{\mathscr C(\mathcal S)}+ \beta_{\mathcal S}$ o\`u $\alpha_{\mathcal S}$ et $\beta_{\mathcal S}$ sont des polyn\^{o}mes en $\mathcal H$ de degr\'e au plus $1$ avec des s\'eries alg\'ebriques comme coefficients (en particulier, $\alpha_{\mathcal S}$ est de valuation au moins 2, en $x$, s'il est de degr\'e z\'ero en $\mathcal H$).\\

En utilisant l'\'equation \eqref{eq:classeindecbis} alors ou bien
\begin{equation} \label {eq:almostfinal1}
\mathcal H_{\mathscr C(\mathcal S)}=\frac{\mathcal H}{1+ \mathcal H};
\end{equation}
dans le cas o\`u toutes les bornes $\mathcal B_i$ sont $\mathcal S$-ind\'ecomposables, ou bien

\begin{equation} \label {eq:almostfinal2}
\mathcal H_{\mathscr C(\mathcal S)}=\frac{(1-\alpha)\mathcal H-\delta}{1+\beta};
\end{equation}
si au moins une des bornes $\mathcal B_i$ n'est pas  $\mathcal S$-ind\'ecomposable.
%\end{proofclaim}

\hfill $\Box$
\vspace{2mm}

La conclusion du Th\'eor\`eme \ref{theo:algebraic} ci-dessus n'est pas v\'erifi\'ee pour des structures qui ne sont pas n\'ecessairement ordonn\'ees. %does not hold with structures which are not necessarily ordered.
\vspace{1mm}

%\begin{remark}\label{rem:cograph}
\begin{example}\label{exemp:cograph}

%La conclusion du th\'eor\`eme \ref{theo:algebraic} ci-dessus n'a pas lieu pour des structures qui ne sont pas n\'ecessairement ordonn\'ees.
Soit $K_{\infty,\infty}$ la somme directe d'une infinit\'e de copies du graphe complet sur un ensemble infini. Comme il est facile de le voir, le profil de $Age(K_{\infty,\infty})$ est la fonction partition d'entier (voir la section \ref{subsec:profil} du chapitre \ref{chap:generalite}). Sa s\'erie g\'en\'eratrice est donn\'ee par $\underset{n\geq 1}\prod\frac{1}{1-x^n}$, elle n'est pas alg\'ebrique. Cependant, $Age(K_{\infty,\infty})$ ne contient aucun membre ind\'ecomposable ayant plus de deux \'el\'ements. Plus g\'en\'eralement, notons que la classe $Forb(P_4)$, des cographes finis ne contenant pas de sous-graphe induit isomorphe au chemin $P_4$, ne contient aucun cographe ind\'ecomposable ayant plus de deux sommets. Cette classe contient $Age(K_{\infty,\infty})$ et donc n'est pas h\'er\'editairement alg\'ebrique. Les cographes finis sont les graphes de comparabilit\'e des ordres s\'eries parall\`eles, qui sont \`a leurs tours les ordres d'intersection associ\'es aux bicha\^{\i}nes s\'eparables. D'apr\`es le th\'eor\`eme d'Albert-Atkinson, la classe de ces  bicha\^{\i}nes est h\'er\'editairement alg\'ebrique. Ceci nous dit que l'alg\'ebricit\'e n'est pas n\'ecessairement pr\'eserv\'ee par la transformation d'une classe en une autre via un processus comme ci-dessus
(les processus de ce type sont les op\'erateurs libres de Fra\"{\i}ss\'e \cite{fraisse}).

\end{example}

\clearemptydoublepage

\chapter{Exemple et conjecture}\label{chap:exemple-conjecture}

\section{Exemple: les structures binaires ordonn\'ees s\'eparables}\label{section:exemple}
Nous construisons dans cette section, pour tout entier $k\geq 1$, une classe de structures binaires ordonn\'ees de type $k$ qui est h\'er\'editairement alg\'ebrique.
Consid\'erons l'ensemble $\mathscr D_{(2)}^k$ des structures binaires ordonn\'ees de type $k\geq 1$, r\'eflexives\index{structure binaire!ordonn\'ee s\'eparable} d\'efinies sur deux \'el\'ements ($\mathscr D_{(2)}^k\subseteq\Gamma_k$) (nous pouvons, \'egalement, supposer que les $k$ relations sont irr\'eflexives, cela ne changera rien \`a la suite (voir note sur les sommes  lexicographiques en page \pageref{compo-reduction})).  Ces structures, \`a l'isomorphie pr\`es, ont la forme %\footnote{Ces structures sont celles qui repr\'esentent les types d'isomorphie.}
$\mathcal S:=(\{0,1\},\leq,\rho_1,\cdots,\rho_k)$ o\`u $\leq$ est l'ordre naturel ($0<1$).

\vspace{2mm}

 Soit $\mathcal S\in \mathscr D_{(2)}^k$, consid\'erons $\mathcal S$ comme une op\'eration qui associe \`a tout couple $(\mathcal R_1, \mathcal R_2)$ d'\'el\'ements de $\Theta_k$ l'\'el\'ement $\mathcal S(\mathcal R_1, \mathcal R_2):=\mathcal R_1 \underset{\mathcal S}\oplus \mathcal R_2$.\\
 D\'esignons par $\mathscr S_k$ la classe des structures binaires ordonn\'ees, engendr\'ee par les structures \`{a} un \'el\'ement de $\Theta_k$ et les op\'erations $\mathcal S$ pour $\mathcal S \in \mathscr D_{(2)}^k$ et par $\mathscr T_k^s$ la classe des types d'isomorphie %\footnote{Une structure de $\mathscr T_k^s$ d'ordre $n$ est une structure d\'efinie sur $\{0,1,\cdots,n-1\}$ telle que $0<1<\cdots<n-1$ et est isomorphe \`a une structure de $\mathscr S_k$.}
 des structures de $\mathscr S_k$.  \\Nous appellerons les structures de $\mathscr S_k$ des \emph{structures binaires ordonn\'ees s\'eparables de type $k.$} Le terme ''\emph{s\'eparable}'' a \'et\'e emprunt\'e \`{a} la terminologie des permutations.

\vspace{2mm}

 La classe des structures binaires ordonn\'ees s\'eparables est analogue aux classes d'objets largement \'etudi\'ees qui sont les classes des co-graphes (ou graphes sans $P_4$), des ordres s\'erie-parall\`eles (ou ordres sans N), des permutations s\'eparables \cite{A-A-V} et des $d$-permutations s\'eparables \cite{Asin-Mans}.

Cette classe partage des propri\'et\'es avec les classes pr\'ec\'edemment cit\'ees, notament la taille de leurs structures ind\'ecomposables (simples pour les permutations) qui ne d\'epasse pas deux.

\vspace{1mm}

Si $k=1$ et $\rho_1$ un ordre total (il existe, \`a l'isomorphie pr\`es, deux telles structures dans $\mathscr D_{(2)}^1$ qui sont $\mathcal S_1:=(\{0,1\},\leq,\rho_1)$ avec $0<1$ et $0<_{\rho_1}1$  %\footnote{$\mathcal S_1:=(\{0,1\},\leq,\rho_1)$ avec $0<1$ et $0<_{\rho_1}1$.}
et $\mathcal S_2:=(\{0,1\},\leq,\rho_1')$ avec $0<1$ et $1<_{\rho_1'}0$), la sous-classe $\mathscr B$ de $\mathscr S_1$, induite par la structure r\'eflexive \`a un \'el\'ement et ces deux op\'erations ($\mathcal S_1$ et $\mathcal S_2$), est celle des bicha\^{i}nes s\'eparables (voir section \ref{susec:bichaineseparable}) et la classe $\mathscr B^t$, des types d'isomorphie des structures de $\mathscr B$, est celle des permutations s\'eparables (voir section \ref{subsc:perm-separable}).

\vspace{2mm}

Notons par $\mathscr S_k^{re}$ (resp. $\mathscr S^{ir}_k$) la sous-classe form\'ee des structures r\'eflexives (resp. irr\'eflexives) de $\mathscr S_k$ et par $\mathscr T^{re}_k$ (resp. $\mathscr T^{ir}_k$) la classe des types d'isomorphie des structures de $\mathscr S_k^{re}$ (resp. $\mathscr S^{ir}_k$).

\begin{fact}\label{fact:bijection}
 Il existe une correspondance bijective entre les deux sous-classes $\mathscr S^{re}_k$ et $\mathscr S^{ir}_k$ (et donc entre $\mathscr T^{re}_k$ et $\mathscr T^{ir}_k$).
 \end{fact}

 Nous avons les r\'esultats suivants:

\begin{lemma}
Pour tout $k\geq 1$, les classe $\mathscr S_k$, $\mathscr S^{re}_k$ et $\mathscr S^{ir}_k$ sont des classes h\'er\'editaires pour la relation d'abritement dans $\Theta_k$ et les classes $\mathscr T^s_k$, $\mathscr T^{re}_k$ et $\mathscr T^{ir}_k$ sont h\'er\'editaires dans la classe des types d'isomorphie des structures de $\Theta_k$, pour la relation d'abritement.
\end{lemma}

\begin{proof} Il suffit de remarquer que toute structure qui s'abrite dans une structure s\'eparable et elle m\^eme s\'eparable et qu'une structure r\'eflexive (resp. irr\'eflexive) ne s'abrite que dans une structure r\'eflexive (resp. irr\'eflexive).
\end{proof}

\begin{proposition}\label{prop:relationseparable}
Pour tout $k\geq 1$, les classes $\mathscr S_k$, $\mathscr S^{re}_k$ et $\mathscr S^{ir}_k$ v\'erifient:
\begin{equation}\label{eq:separable}
\underset{r \in \mathscr D_{(2)}^k}{\bigcup}\mathscr S {\underset{r}\oplus} \mathscr S\subseteq \mathscr S.
\end{equation}
 et les classes $\mathscr S^{re}_k$ et $\mathscr S^{ir}_k$ sont des classes h\'er\'editaires minimales parmi celles v\'erifiant cette relation, dans $\Theta_k.$
\end{proposition}

\begin{proof}
La premi\`ere partie vient directement de la d\'efinition. Pour la minimalit\'e il suffit de remarquer que si une classe h\'er\'editaire v\'erifie la relation \eqref{eq:separable} et contient la structure r\'eflexive (resp. irr\'eflexive) \`a un \'el\'ement  alors elle contient la classe $\mathscr S^{re}_k$ (resp. $\mathscr S^{ir}_k$).
\end{proof}

%***************************************************************************
\subsection{Ind\'ecomposabilit\'e et structures s\'eparables}

Nous disons qu'une structure $\mathcal R$ est $2$-\emph{d\'ecomposable}\index{structure relationnelle!$2$-d\'ecomposable} si elle s'\'ecrit comme somme lexicographique de deux structures index\'ee par une structure \`a deux \'el\'ements. Nous disons qu'elle est $2$- \emph{ind\'ecomposable}\index{structure relationnelle!$2$-ind\'ecomposable} dans le cas contraire. Si $\mathcal R$ est $2$-d\'ecomposable, nous \'ecrivons simplement (sauf n\'ecessit\'e de pr\'ecision) $\mathcal R= \mathcal R_1 \oplus \mathcal R_2.$ Nous avons le r\'esultat suivant:

\begin{proposition}\label{prop:2-decoposable}
Une structure binaire ordonn\'ee de type $k$ est s\'eparable si et seulement si toute structure de taille sup\'erieure ou \'egale \`a $2$ qu'elle abrite est $2$-d\'ecomposable.
\end{proposition}

\begin{proof}
\begin{enumerate}
\item $\Rightarrow$) \'evident (par d\'efinition d'une structure s\'eparable).
\item $\Leftarrow$) D'apr\`es la relation \eqref{eq:separable} de la Proposition \ref{prop:relationseparable}, si une structure binaire ordonn\'ee $\mathcal R$ de type $k$  n'est pas s\'eparable, elle n'est pas dans $\mathscr S_k$ donc n'est pas dans $\underset{r}{\bigcup}\mathscr S_k {\underset{r}\oplus} \mathscr S_k$ et par suite abrite au moins une structure $2$-ind\'ecomposable.
\end{enumerate}
\end{proof}

De la proposition \ref{prop:2-decoposable} nous avons:

\begin{corollary}\label{cor:separable}
Une structure binaire ordonn\'ee $\mathcal R$ est s\'eparable si et seulement si les structures de $Ind(\mathcal R)$ sont de tailles au plus deux.
\end{corollary}
Si nous posons $\mathscr D_{(\leq2)}^k=\mathscr D_{(2)}^k\cup\{\textbf{1}\}$ o\`u $\textbf{1}$ est la structure \`a un \'el\'ement de $\Gamma_k$ et si $\mathscr {D'}_{(\leq 2)}^k$ d\'esigne l'ensemble des structures binaires ordonn\'ees de type $k\geq 1$ irr\'eflexives d\'efinies sur au plus deux \'el\'ements et $\overline{\mathscr D}_{(\leq2)}^k$ l'ensemble des structures binaires ordonn\'ees de type $k$ d\'efinies sur au plus deux \'el\'ements, alors
\begin{corollary}\label{cor:sumseparable}
$$1)\quad\mathscr S^{re}_k=\sum\mathscr D_{(\leq2)}^k.$$
$$2)\quad\mathscr S^{ir}_k=\sum\mathscr {D'}_{(\leq2)}^k.$$
$$3)\quad\mathscr S_k=\sum\overline{\mathscr D}_{(\leq2)}^k.$$
\end{corollary}
\medskip

D'apr\`es la Cons\'equence 1 de la Proposition \ref{prop:borneindec}, nous avons:
\begin{corollary}
Les bornes des classes $\mathscr S^{re}_k$, $\mathscr S^{ir}_k$ et $\mathscr S_k$ sont ind\'ecomposables.
\end{corollary}

D'apr\`es le Th\'eor\`eme \ref{theo:closed} et la note qui lui est associ\'ee, nous avons imm\'ediatement

\begin{corollary}\label{cor:close-separable}
$$1)\quad\mathscr S^{re}_k=cl(\mathscr D_{(\leq2)}^k).$$
 $$2)\quad\mathscr S^{ir}_k=cl(\mathscr {D'}_{(\leq2)}^k).$$
  $$3)\quad\mathscr S_k=cl(\overline{\mathscr D}_{(\leq2)}^k).$$
\end{corollary}
Dans la section suivante, nous donnons, dans le cas des structures binaires ordonn\'ees s\'eparables de type $1$, une caract\'erisation des bornes.
%************************************************************************************
\subsection{Etude de la classe $\mathscr S_1$}\label{subsection:classe S1}

Soit $\rho$ une relation binaire d\'efinie sur un ensemble $E$.
Pour simplifier la pr\'esentation, nous identifions $\rho$ \`a sa fonction caract\'eristique,\index{fonction!caract\'eristique} ce qui signifie que l'on consid\`ere $\rho$ comme une application
 $\rho:E^2\rightarrow \{0,1\}$.
Posons pour tout $x, y \in E$:
\begin{equation}\label{eq:valeur relation}
d(x,y):=(\rho(x,y),\rho(y,x)).
\end{equation}

\noindent D\'esignons par $A$ et $B$ deux ensembles de structures de type $1$ %et notons $T_A$  l'ensemble de ces structures,
%(resp. de \emph{type $B$}%, d'ensemble not\'ee $T_B$
de la forme $\mathcal S:=(E,\leq,\rho)$, telles que:

$\mathcal S\in A\Leftrightarrow\left\{
\begin{array}{l}
E=\{0,1,2\} \\
0<1<2 \\
d(0,1)\neq d(0,2) \\
d(1,2)\neq d(0,2)
\end{array}
\right. \qquad   \mathcal S \in B\Leftrightarrow\left\{
\begin{array}{l}
E=\{0,1,2,3\} \\
0<1<2<3 \\
d(0,1)= d(0,3)= d(2,3) \\
d(1,2)= d(1,3)= d(0,2)\\
d (0,1) \neq d (0,2)%
\end{array}%
\right. $
\medskip

 Posons $\mathcal F=A\cup B$. %Un d\'ecompte simple montre que nous avons $36$ structures r\'eflexives de \emph{type} $A$ et $12$ structures r\'eflexives de \emph{type} $B$.
Nous avons:

\begin{fact}\label{fact:separable}
Une structure binaire ordonn\'ee de type $1$ d\'efinie sur trois \'el\'ements est ind\'ecomposable si et seulement si elle est isomorphe \`a une structure de $A$.
\end{fact}

\begin{proof}
Il est facile de v\'erifier que toutes les structures de  $A$ sont ind\'ecomposables. Inversement, soit $\mathcal S$ une structure binaire ordonn\'ee \`a trois \'el\'ements  ind\'ecomposable. $\mathcal S$ est isomorphe \`a une structure $\mathcal S'$ d\'efinie sur $\{0,1,2\}$ tels que $0<1<2$. Donc les ensembles $\{0,1\}$ et $\{1,2\}$ ne sont pas des intervalles de $\mathcal S'$, ce qui imlique que $d(0,2)\neq d(1,2)$ et $d(0,1)\neq d(0,2)$. Donc $\mathcal S'\in A$.
\end{proof}

\begin{remarks}\label{rem:separable}

\begin{enumerate}
\item Il y a exactement deux structures dans $A$ pour lesquelles  $\rho$ est un tournoi (c'est \`a dire $d(i,j)\in \{(1,0),(0,1)\}$, pour tout $i,j\in E$); ce sont celles pour lesquelles $\rho$ est un cycle de longueur $3$.
\item Il n'existe aucune structure de $A$ pour laquelle $\rho$ est un ordre total  et dans $B$ il y a deux bicha\^{i}nes \`a quatre \'el\'ements correspondant aux permutations $2413$ et $3142.$
\item Les structures de $B$ sont ind\'ecomposables.
\end{enumerate}
\end{remarks}

\begin{theorem}\label{theo:separable}
Une structure binaire ordonn\'ee de type $1$ est s\'eparable si et seulement si elle n'abrite aucune structure de $\mathcal F.$ En d'autres termes $\mathscr S_1=\text{\it Forb}(\mathcal F)$.
\end{theorem}

Avant de donner la preuve de ce th\'eor\`eme, montrons le r\'esultat suivant:

\begin{lemma}\label{lem:separable}
Toute structure binaire ordonn\'ee de type $1$ ind\'ecomposable d'ordre quatre abrite au moins une des structures de $\mathcal F.$
\end{lemma}

\begin{proof}
Il suffit de faire la preuve pour une structure repr\'esentant un type d'isomorphie. Soit alors $\mathcal R:=(\{0,1,2,3\},\leq,\rho)$ une structure ind\'ecomposable telle que $0<1<2<3$. Les valeurs de $d(0,1), d(0,2)$ et $d(0,3)$ ne peuvent \^etre \'egales car sinon $\{1,2,3\}$ serait un intervalle. Donc, ou bien ces trois valeurs sont toutes diff\'erentes, ou bien deux d'entre elles sont \'egales et diff\'erentes de la troisi\`eme. Nous avons alors les cas suivants:\newline
\textbf{Cas $1$:} Supposons $d(0,1)=d(0,2)\neq d(0,3)$ alors:%avec ce choix $\{2,3,4\}$ ne pas \^etre un intervalle, reste \`a \'eliminer $\{1,2\}$ et $\{2,3\}$ et \{1,2,3\}$.
  \begin{itemize}
  \item soit $d(1,3)\neq d(0,3)$ et dans ce cas $\mathcal R_{\restriction_{\{0,1,3\}}}\in \mathcal F,$
   \item soit $d(1,3)= d(0,3)$ et dans ce cas nous avons forc\'ement $d(2,3)\neq d(0,3)$, car sinon $\{0,1,2\}$ serait un intervalle de $\mathcal R$ ce qui contredit le fait que $\mathcal R$ est ind\'ecomposable, d'o\`u $\mathcal R_{\restriction_{\{0,2,3\}}}\in \mathcal F.$
   \end{itemize}

\textbf{Cas $2$:} Si $d(0,1)\neq d(0,2)= d(0,3).$ Ce cas est identique au premier.\newline

\textbf{Cas $3$:} Si $d(0,1)= d(0,3)\neq d(0,2)$ alors nous avons les cas suivants:
    \begin{enumerate}
  \item soit $d(1,2)\neq d(0,2)$ et dans ce cas $\mathcal R_{\restriction_{\{0,1,2\}}}\in \mathcal F,$
   \item soit $d(1,2)= d(0,2)$ et dans ce cas nous avons forc\'ement $d(1,3)\neq d(0,3)$ car sinon $\{0,1\}$ serait un intervalle. Deux cas se pr\'esentent:
        \begin{enumerate}
            \item soit $d(2,3)\neq d(0,3)$ auquel cas $\mathcal R_{\restriction_{\{0,2,3\}}}\in \mathcal F,$
            \item soit $d(2,3)= d(0,3)$ et dans ce cas nous avons:
                \begin{itemize}
                    \item ou bien $d(1,2)\neq d(1,3)$ et nous avons d\'ej\`a $d(1,3)\neq d(0,3)= d(2,3)$, donc $\mathcal R_{\restriction_{\{1,2,3\}}}\in \mathcal F,$
                    \item ou bien $d(1,2)= d(1,3)$, ce qui donne une structure de  $B$, donc dans $\mathcal F$.
                \end{itemize}
        \end{enumerate}
   \end{enumerate}
\textbf{Cas $4$:} Si $d(0,1), d(0,2)$ et $d(0,3)$ ont des valeurs toutes diff\'erentes alors:
   \begin{enumerate}
    \item si $d(1,2) \neq d(0,2)$ alors $\mathcal R_{\restriction_{\{0,1,2\}}}\in \mathcal F,$
    \item si, $d(1,2) = d(0,2)$ et dans ce cas nous avons:
        \begin{enumerate}
        \item si $d(2,3) \neq d(0,3)$ alors $\mathcal R_{\restriction_{\{0,2,3\}}}\in \mathcal F,$
        \item si, $d(2,3) = d(0,3)$ et dans ce cas nous avons, forc\'ement, $d(1,3)\neq d(0,3)$, car sinon $\{0,1,2\}$ serait un intervalle. Nous avons alors $\mathcal R_{\restriction_{\{0,1,3\}}}\in \mathcal F.$
        \end{enumerate}
   \end{enumerate}
   Donc toute structure ind\'ecomposable d'ordre quatre abrite, au moins, une structure de $\mathcal F.$
\end{proof}
\bigskip

\noindent{\bf \emph{Preuve du Th\'eor\`eme \ref{theo:separable}}.}
\begin{enumerate}
\item[$\Rightarrow$)] Si $\mathcal R$ est s\'eparable alors, d'apr\`es le Corollaire \ref{cor:separable}, elle n'abrite pas de structure ind\'ecomposable de taille sup\'erieure \`a deux et donc ne peut pas abriter une structure de $\mathcal F$ (d'apr\`es le Fait \ref{fact:separable} et le point $(3)$ des Remarques \ref{rem:separable}).
\item[$\Leftarrow$)] Si $\mathcal R$ n'est pas s\'eparable alors, toujours d'apr\`es le Corollaire \ref{cor:separable} elle abrite une structure ind\'ecomposable $\mathcal S$ de taille $n\geq 3$. Montrons par r\'ecurrence sur $n$ que toute structure ind\'ecomposable $\mathcal S$ d\'efinie sur un $n$-ensemble abrite une structure de $\mathcal F$.
    \begin{itemize}
    \item Si $n=3$ alors, d'apr\`es le Fait \ref{fact:separable}, $\mathcal S\in\mathcal F$ et si $n=4$ alors d'apr\`es le Lemme \ref{lem:separable},  $\mathcal S$ abrite une structure de $\mathcal F.$
    \item Supposons que toute structure ind\'ecomposable d'ordre $(n-1)$ ($n\geq 4$) abrite une structure de $\mathcal F$. Soit $\mathcal S$ une structure ind\'ecomposable d'ordre  $n\geq 5$ alors d'apr\`es le Th\'or\`eme \ref{theo:schmerl-trotter} de Schmerl et Trotter, %\footnote{Schmerl et Trotter ont montr\'e, dans \cite{S-T}, que toute structure ind\'ecomposable de taille $n\geq 3$ abrite, au moins, une structure ind\'ecomposable de taille $n-1$ ou $n-2$.} \cite{S-T},
         $\mathcal S$ abrite une structure ind\'ecomposable d'ordre $(n-1)$ ou $(n-2)$. Ce qui termine la preuve du th\'eor\`eme.
    \end{itemize}
\end{enumerate}
\hfill $\Box$

%*****************************************************
\subsubsection{Profil de la classe $\mathscr S^{re}_1$ (resp. $\mathscr S^{ir}_1$)}
D'apr\`es le Fait \ref{fact:bijection}, les classes $\mathscr S^{re}_1$ et $\mathscr S^{ir}_1$ ont m\^eme profil qui est celui des classes $\mathscr T^{re}_1$ et $\mathscr T^{ir}_1.$

D\'esignons par $S^r_1$ la fonction g\'en\'eratrice de la classe $\mathscr S^{re}_1$ (donc de $\mathscr S^{ir}_1$).
 Nous avons les r\'esultats suivants:

\begin{proposition}\label{prop:generatriceseparable}
La fonction g\'en\'eratrice $S^r_1$ v\'erifie:
\begin{equation}\label{eq:6}
3(S^r_1)^2+(x-1)S^r_1+x=0.
\end{equation}
 et est donn\'ee par:
 \begin{equation}\label{eq:7}
 S^r_1=\dfrac{1-x-\sqrt{1-14x+x^{2}}}{6}.
 \end{equation}
\end{proposition}

\begin{proof}
$\mathscr D^1_{(2)}$ poss\`ede, \`a l'isomorphie pr\`es, quatre structures. Ecrivons chaque $r$ de $\mathscr D^1_{(2)}$ comme $r:=(\{0,1\},\leq, \rho)$ avec $0<1$, et
posons $\mathscr D^1_{(2)}=\{r_1,r_2,r_3,r_4\}.$ D'apr\`es le lemme \ref{lem:union} et le corollaire \ref{cor:close-separable} nous avons:
$$\mathscr S^{re}_1=\{\bf{1}\}\cup \underset{1\leq i\leq 4}\bigcup (\mathscr S^{re}_1(r_i)\underset{r_i}\oplus \mathscr S^{re}_1).$$
et pour chaque $i\in \{1,2,3,4\}$ nous avons:
$$\mathscr S^{re}_1(r_i)=\{\bf 1\}\cup\underset{j\neq i} \bigcup (\mathscr S^{re}_1(r_j)\underset{r_j}\oplus \mathscr S^{re}_1).$$
o\`u $\textbf{1}$ d\'esigne la structure \`a un \'el\'ement et $\mathscr S^{re}_1(r_i)$ la classe des structure de $\mathscr S^{re}_1$ qui sont $r_i$-ind\'ecomposables.
En passant aux s\'eries g\'en\'eratrices nous  obtenons:
\begin{equation}\label{eq:gensepar1}
S^{r}_1=x+S^{r}_1\underset{1\leq i\leq 4}\sum S'_i,
\end{equation}
et
\begin{equation}\label{eq:gensepar2}
S'_i=x+S^{r}_1\underset{j\neq i}{\sum}S'_j \text{ pour tout }i\in\{1,2,3,4\}.
\end{equation}
o\`u $S'_i$ est la s\'erie g\'en\'eratrice de $\mathscr S^{re}_1(r_i).$ En \'eliminant les $S'_i,~i\in\{1,2,3,4\}$ de l'\'equation \eqref{eq:gensepar1}, nous obtenons l'\'equation \eqref{eq:6}. La r\'esolution de cette \'equation donne la fonction \eqref{eq:7}. L'expansion de cette fonction donne la s\'erie dont les premiers termes sont
$$S^{r}_1=x+4x^2+28x^3+244x^4+2380x^5+\cdots.$$
\end{proof}

Cette s\'erie se trouve dans ''On line encyclop\'edia of integer sequences'' sous le num\'ero $A103211$ \cite{Sloane}.

\vspace{2mm}

Dans \cite{Ack-Ba-Pin-Rom}, les auteurs ont \'etudi\'e les partitions guillotine en dimension $d$ et ont donn\'e la s\'erie g\'en\'eratrice du nombre de ces partitions, elle v\'erifie $f=1+xf+(d-1)xf^2,~(*)$. Ils ont \'egalement \'etabli une correspondance bijective entre les partitions guillotine en dimension $d$ et les arbres binaires \'etiquet\'es par l'ensemble $\{1,\cdots,d\}$.

\vspace{2mm}

Dans \cite{Asin-Mans}, les auteurs ont donn\'e la s\'erie g\'en\'eratrice des $d$-permutations s\'eparables, elle v\'erifie $f=1+xf+(2^{d-1}-1)xf^2,~~(**)$ et ont \'etabli une correspondance bijective entre les $d$-permutations s\'eparables et les partitions guillotine en dimension $2^{d-1}$. Il se trouve que notre s\'erie $S^{r}_1$ v\'erifie $(*)$ pour $d=4$ et $(**)$ pour $d=3$.

\vspace{2mm}

Pour ne pas allourdir cette partie avec des rappels sur les partitions guillotine qui ne rentrent pas dans le cadre de ce travail, nous renvoyons ces rappels \`a l'annexe \ref{codage} où nous construisons \'egalement,  une bijection entre les structures binaires ordonn\'ees s\'eparables de type $1$ et les arbres binaires \'etiquet\'es par l'ensemble $\{1,2,3,4\}$, ceci entraine que l'ensemble $\mathscr S^{re}_1$ et l'ensembles des partitions guillotine en dimension $4$ sont isomorphes. Cela entraine \'egalement que $\mathscr S^{re}_1$ est isomorphe \`a l'ensemble des $3-$permutations.

\begin{corollary}\label{cor:algebriqueS1}
Les classes $\mathscr S^{re}_1$ et $\mathscr S^{ir}_1$ sont alg\'ebriques.
\end{corollary}
\medskip
Nous donnons pour $n=1,2,3,4$ les formes de ces structures voir \tablename~\ref{tableseparable} (les boucles ne sont pas repr\'esent\'ees).
\begin{itemize}
\item Pour $n=1$ nous avons une seule structure.
\item Pour $n=2$ nous avons quatre structures qui sont $r_i,\; i\in \{1,2,3,4\}.$
\item Pour $n=3$, les structures s'obtiennent en faisant la somme, suivant l'une des structures $r_i$, de la structure $\textbf{1}$ et d'une structure \`{a} deux \'el\'ements $r_j$. Nous avons deux cas (voir \tablename~\ref{tableseparable}).
        Dans le cas $(1)$ nous avons $4$ structures diff\'erentes ($r_j$) pour chaque $i=1,\cdots,4$, donc $4\times 4=16$ structures diff\'erentes et dans le cas $(2)$ nous devons exclure pour chaque $i$ la structure $r_j=r_i$ (car elle est isomorphe \`a une structure trouv\'ee dans le cas $(1)$) ce qui donne $3\times 4$ structures diff\'erentes. Au total, nous obtenons les $28$ structures d'ordre $3.$
\item Pour $n=4$, Les structures s'obtiennent en faisant la somme, suivant l'une des structures $r_i$, de la structure $\textbf{1}$ avec une structure \`{a} trois \'el\'ements, cas $(1)$ et $(2)$ dans la \tablename~\ref{tableseparable}, ou de deux structures \`a deux \'el\'ements cas $(3)$.
       Dans le cas $(1)$ pour chaque $r_i$, nous avons $28$ structures $r_j$ ce qui donne $4\times 28$ structures diff\'erentes, dans le cas $(2)$ nous avons pour chaque $r_i$, $28$ structures $r_j$ mais nous devons exclure celles qui se d\'ecomposent suivant $r_i$ qui sont au nombre de $7$ ce qui fait $21$ structures pour chaque $r_i$, donc $4\times 21$ structures. Dans le cas $(3)$ la structure $r_j$ doit \^etre diff\'erente de $r_i$ ce qui donne $3\times 4=12$ structures diff\'erentes pour chaque $r_i$ donc $4\times12$ structures. Au total, nous obtenons les $244$ structures d'ordre $4.$
\end{itemize}

\begin{table}[!hbp]
\begin{tabular}[c]{|c|c|c|c|c|}
\hline
\input{fig1}&\input{fig2}&\input{fig3}&\input{fig4}&\\
\hline
$d(0,1)=(1,0)$& $d(0,1)=(0,1)$ & $d(0,1)=(1,1)$& $d(0,1)=(0,0)$&\\
\hline
\multicolumn{5}{|c|}
{Structures \`a deux \'el\'ements}\\
\hline
\input{imag1}&\quad \input{imag2}&\input{imag1a}&\input{imag2a}&\input{imag3a}\\
\hline
Cas $(1)$& Cas $(2)$&Cas $(1)$& Cas $(2)$&Cas $(3)$\\
\hline
\multicolumn{2}{|c|}
{Structures d'ordre $3$}& \multicolumn{3}{|c|}
{Structures d'ordre $4$}\\
\hline
\end{tabular}
\caption{\label{tableseparable}Construction des structures s\'eparables d'ordres $3$ et $4.$}
\end{table}

%********************************************************************
\subsubsection{Profil de la classe $\mathscr S_1$}
%Dans le cas o\`u la relation $\rho$ n'est pas r\'eflexive ou irr\'eflexive alors si $\mathscr S_1$ d\'esigne la classe de ces structures:
D\'esignons par $S_1$ la s\'erie g\'en\'eratrice de $\mathscr S_1$. Nous avons:
\begin{proposition}
La s\'erie  $S_1$ v\'erifie:
\begin{equation}\label{eq:8}
3{S_1}^2+(2x-1)S_1+2x=0.
\end{equation}
 et est donn\'ee par:
 \begin{equation}\label{eq:9}
S_1=\dfrac{1-2x-\sqrt{1-28x+4x^{2}}}{6}.
\end{equation}
\end{proposition}

\begin{proof}
La preuve de cette proposition se fait de la m\^eme fa\c{c}on que celle de la Proposition \ref{prop:generatriceseparable} \`{a} la diff\'erence que dans ce cas nous avons deux structures \`{a} un \'el\'ement.
\end{proof}
\medskip

 L'expansion de cette fonction en s\'erie donne:
$$S_1=2x+16x^2+224x^3+3904x^4+\cdots.$$

\bigskip
Nous pouvons, comme dans le cas pr\'ec\'edent, donner, pour $n=1,2,3,4,$ les formes des structures (ces structures se construisent de la m\^eme fa\c{c}on que les pr\'ec\'edentes, mais leurs nombres augmentent car elles ne sont pas r\'eflexives ou irr\'eflexives).
\begin{itemize}
\item Pour $n=1$ nous avons deux structures.
\item Pour $n=2$ nous avons, pour chaque $r_i$ quatre structures qui coincident avec $r_i$ sur le couple $(0,1)$ donc il y a $16$ structures \`{a} deux \'el\'ements.
\item Pour $n=3$, les structures s'obtiennent en faisant la somme, suivant l'une des structures $r_i$, d'une structure \`{a} un \'el\'ement et d'une structure \`{a} deux \'el\'ements.
    Dans le cas $(1)$ nous avons $2$ structures \`a un \'el\'ement et $16$ structures diff\'erentes ($r_j$) pour chaque $r_i,~1\leq i\leq 4$, donc $2\times16\times 4=128$ structures diff\'erentes et dans le cas $(2)$ nous devons exclure pour chaque $i$ la structure $r_j=r_i$ (car elle est isomorphe \`a une structure trouv\'ee dans le cas $(1)$), donc le nombre de structures $r_j$ est $16-4$ ce qui donne $12\times2\times 4$ structures diff\'erentes. Au total, nous obtenons les $224$ structures d'ordre $3.$
\item Pour $n=4$, Les structures s'obtiennent en faisant la somme, suivant l'une des structures $r_i$, d'une structure \`a un \'el\'ement et d'une  structure \`{a} trois \'el\'ements, cas $(1)$ et $(2)$, ou de deux structures \`a deux \'el\'ements, cas $(3)$.
    Dans le cas $(1)$ pour chaque $r_i$, nous avons $224$ structures \`a trois \'el\'ements $r_j$ et deux structures \`a un \'el\'ement, ce qui donne $4\times224\times 2$ structures diff\'erentes. Dans le cas $(2)$ nous avons pour chaque $r_i$, deux structures \`a un \'el\'ement et $224$ structures $r_j$ mais nous devons exclure celles qui se d\'ecomposent suivant $r_i$ qui sont, d'apr\`es le cas pr\'ec\'edent,  au nombre de $32+24$, ce qui fait $168$ structures pour chaque $r_i$, donc  $4\times 168\times2$ structures. Dans le cas $(3)$ nous avons, pour chaque $r_i$,  $16$ structures pour $r_k$, la structure $r_j$ ne doit pas se d\'ecomposer suivant $r_i$ ce qui donne $3\times 4=12$ structures pour $r_j$, donc $4\times 12\times 16$ structures diff\'erentes. Au total, nous obtenons les $3904$ structures d'ordre $4.$
\end{itemize}

%***************************************************************************************
\subsection{Profil de $\mathscr S_k$ pour $k\geq 2$}
Rappelons que $\mathscr S_k$ pour $k\geq 2$, est l'ensemble des structures binaires ordonn\'ees s\'eparables de type $k$. D\'esignons par $S_k$ sa s\'erie g\'en\'eratrice.

\begin{proposition}\label{prop:gener-sepa-k}
La s\'erie g\'en\'eratrice $S_k$ v\'erifie:
\begin{equation}\label{eq:sep-k1}
(4^{k}-1){S_k}^2+(2^{k}x-1)S_k+2^{k}x=0.
\end{equation}
 et est donn\'ee par:
 \begin{equation}\label{eq:sep-k2}
 S_k=\dfrac{1-2^{k}x-\sqrt{1-2^{k+1}(2.4^{k}-1)x+2^{2k}x^{2}}}{2(4^{k}-1)}.
 \end{equation}
\end{proposition}

\begin{proof}
Dans ce cas $\mathscr D^k_{(2)}$ comporte $4^{k}$ structures et nous avons $2^{k}$ structures \`a un \'el\'ement.
D'apr\`es le Lemme \ref{lem:union}, le Corollaire \ref{cor:sumseparable} et le Corollaire \ref{cor:close-separable}, nous avons:
$$\mathscr S_k=\underset{1\leq i\leq 2^{k}}\bigcup\{1_i\}\cup \underset{1\leq i\leq 4^{k}}\bigcup (\mathscr S_k(r_i)\underset{r_i}\oplus \mathscr S_k).$$
et pour chaque $i,~~1\leq i\leq 4^{k}$ nous avons:
$$\mathscr S_k(r_i)=\underset{1\leq i\leq 2^{k}}\bigcup\{1_i\}\cup\underset{j\neq i}\bigcup (\mathscr S_k(r_j)\underset{r_j}\oplus \mathscr S_k).$$
o\`u $1_i$ d\'esigne une structure \`a un \'el\'ement et $\mathscr S_k(r_i)$ la classe des structure de $\mathscr S_k$ qui sont $r_i$-ind\'ecomposables.
En passant aux s\'eries g\'en\'eratrices nous  obtenons:
\begin{equation}\label{eq:gen-separ-k1}
S_k=2^{k}x+S_k\underset{i=1}{\overset{4^{k}}{\sum}}S^i_k,
\end{equation}
et
\begin{equation}\label{eq:gen-separ-k2}
S^i_k=2^{k}x+S_k\underset{j\neq i}{\sum}S^j_k.
\end{equation}
o\`u $S^i_k$ est la s\'erie g\'en\'eratrice de $\mathscr S_k(r_i).$ En \'eliminant les $S^i_k,~ 1\leq i\leq 4^{k}$ de l'\'equation \eqref{eq:gen-separ-k1}, nous obtenons l'\'equation \eqref{eq:sep-k1}. La r\'esolution de cette \'equation donne la fonction \eqref{eq:sep-k2}.
\end{proof}

%$$$$$$$$$$$$$$$$$$$$$$$$$$$$$$$$$$$$$$$$$$$$$$$$$$$$$$$$$$$$$$$$$$$$$$$$$$$

        %   Conjecture

%$$$$$$$$$$$$$$$$$$$$$$$$$$$$$$$$$$$$$$$$$$$$$$$$$$$$$$$$$$$$$$$$$$$$$$$$$$$$$

\section{Conjecture et questions}\label{sec:conjecture}
Dans leur papier \cite{A-A}, Albert et Atkinson ont fait remarquer qu'il existe des ensembles infinis de permutations simples dont la cl\^oture par sommes est alg\'ebrique, mais il s'av\`ere que certaines sous-classes h\'er\'editaires ne sont pas n\'ecessairement alg\'ebriques comme c'est le cas de la collection des oscillations d\'ecroissantes (voir la fin de cette section). Dans le but d'\'etendre leur preuve \`a d'autres classes, ils  se sont demand\'es s'il existe \emph{un ensemble infini de permutations simples dont la cl\^oture par sommes est belordonn\'ee}. Comme nous le verrons dans la proposition \ref{criticalwqo} ci-dessous, l'ensemble des permutations exceptionnelles v\'erifie cette propri\'et\'e. En fait, cet ensemble est h\'er\'editairement belordonn\'e. Nous pensons que cette notion de belordre h\'er\'editaire est le concept qui permet d'\'etendre le th\'eor\`eme d'Albert et Atkinson.

\vspace{1mm}

Les permutations exceptionnelles correspondent aux bicha\^{\i}nes qui sont \emph{critiques} au sens de Schmerl et Trotter (voir D\'efinition \ref{def:critique}). %Rappelons qu'une structure binaire $\mathcal R$ de base $E$ est dite \emph{critique} si $\mathcal R$ est ind\'ecomposable mais $\mathcal R\restriction_{E\setminus  \{x\}}$ n'est pas ind\'ecomposable pour tout $x\in E$.
Schmerl et Trotter  \cite{S-T} ont donn\'e une description des ordres critiques. Ils forment deux classes infinies:
$\mathscr{P}:=\{\mathcal P_n: n\in \mathbb{N}\}$ et $\mathscr{P'}:=\{\mathcal P'_n: n\in \mathbb{N}\}$ o\`{u}
$\mathcal P_n:=(V_n,\leq_n)$, $V_n:=\{0, \dots, n-1\}\times \{0,1\}$,  $(x,i)<_n(y,j)$ si  $i<j \; \text{et} \; x\leq y$;
$\mathcal P'_n:=(V_n,\leq'_n)$   et $(x,i)<'_n (y,j)$ si $ j\leq i \; \text{et} \; x< y$.\medskip

\noindent Ces ordres sont de dimension deux. Ce qui signifie qu'ils sont intersections de deux ordres lin\'eaires qui sont: $L_{n,1}:=(0,0)<(0,1)<\cdots<(i,0)<(i,1)\cdots <(n-1,0)<(n-1,1)$ et $L_{n,2}:=(n-1,0) <\cdots<(n-i,0) <\cdots <(0, 0)< (n-1,1)<\cdots  <(n-i,1)\cdots <(0, 1)$ pour $\mathcal P_n$ et $L'_{n,1}:=L_{n,1}$ et $L'_{n,2}:=(L_{n,2})^{-1}$  pour $\mathcal P'_n.$\medskip

 Il est bien connu  qu'un ordre $\mathcal P:=(V, L)$ de dimension deux qui est ind\'ecomposable a une unique r\'ealisation (c'est \`{a} dire qu'il existe une unique paire  $\{L_1,L_2\}$ d'ordres lin\'eaires dont l'intersection est l'ordre $L$ de $\mathcal P$) (voir \cite{gallai}). Donc,
il y a au plus deux bicha\^{\i}nes,  $(V, L_1,L_2)$ et $(V, L_2,L_1)$ telles que $L_1\cap L_2=L$. Les ordres critiques  d\'ecrits ci-dessus donnent quatre sortes de bicha\^{\i}nes critiques,\index{bichaine@bicha\^{i}ne!critique}
$(V_n,L_{n,1},L_{n,2})$, $(V_n,L_{n,2},L_{n,1})$, $(V_n,L_{n,1},(L_{n,2})^{-1})$ et $(V_n,(L_{n,2})^{-1},L_{n,1})$. %Ces bicha\^{i}nes sont critiques.
En effet, une bicha\^{\i}ne est ind\'ecomposable si et seulement si l'ordre intersection est ind\'ecomposable  (\cite{R-Z} pour les bicha\^{i}nes finies et \cite {zaguia} pour les bicha\^{i}nes infinies).  Les types d'isomorphie de ces bicha\^{\i}nes sont d\'ecrits dans l'article d'Albert et Atkinson \cite{A-A} en terme de permutations de $1, \dots, 2m$ pour $m\geq 2$:\\
\qquad $(i)$\label{prmut-critique}~~$2.4.6....2m.1.3.5....2m-1.$\\
\qquad $(ii)$~~$2m-1.2m-3....1.2m.2m-2....2.$\\
\qquad $(iii)$~~$m+1.1.m+2.2....2m.m.$\\
\qquad $(iv)$~~$m.2m.m-1.2m-1....1.m+1.$\\

Par exemple, le type d'isomorphie de la bicha\^{\i}ne  $(V_m, L_{m,1}, L_{m,2})$ est la permutation donn\'ee en $(iv)$, alors que le type d'isomorphie de $(V_m, L_{m,2}, L_{m,1})$ est son inverse, donn\'e en  $(ii)$
(\'enum\'erer les \'el\'ements de $V_m$ dans la suite $1, \dots, 2m$,  suivant l'ordre $L_{m,1}$, puis r\'eordonner cette suite suivant l'ordre $L_{m,2}$ (voir Exemple \ref{exp:ordre-permut} en page \pageref{exp:ordre-permut}); ceci donne la suite
$\sigma^{-1}:= \sigma^{-1}(1),\dots, \sigma^{-1}(2m)$; d'apr\`es notre d\'efinition, le type de $(V_m, L_{m,1}, L_{m,2})$ est la permutation $\sigma$, celle-ci est celle donn\'ee en $(iv)$).
Pour $m=2$, les permutations donn\'ees en $(i)$ et $(iv)$ co\"{\i}ncident avec  $2413$ alors que celles donn\'ees en $(ii)$ et $(iii)$ co\"{\i}ncident avec $3142$; pour des valeurs de $m$ plus grandes, ces permutations sont toutes diff\'erentes.\medskip

Les quatre classes de bicha\^{\i}nes ind\'ecomposables  sont obtenues \`a partir de  $$\mathscr B:=\{(V_n, L_{n,1}, L_{n,2}): n\in \mathbb{N}\}$$  en permutant les deux ordres dans chaque bicha\^{\i}ne ou bien en inversant le premier ordre ou bien en inversant le second ordre. Par cons\'equent la structure d'ordre (pour l'abritement) de ces classes est la m\^{e}me et  reste la m\^{e}me si on \'etiquette les \'el\'ements des ces bicha\^{\i}nes.

\begin{proposition}\label{criticalwqo}
 La classe des bicha\^{\i}nes critiques est h\'er\'editairement belordonn\'ee.
\end{proposition}

\begin{proof}
Comme nous l'avons mentionn\'e ci-dessus, la classe des bicha\^{i}nes est l'union de quatre classes,  donc pour montrer qu'elle est h\'er\'editairement belordonn\'ee, il suffit de montrer que chacune de ces classe est h\'er\'editairement belordonn\'ee. D'apr\`es l'observation ci-dessus, il suffit de montrer qu'une seule parmi elle l'est, par exemple $\mathscr B$.\\ Soit $\mathcal A$ un ensemble belordonn\'e. Nous devons montrer que $\mathscr B. \mathcal A$ est belordonn\'e. Pour cela, posons $\mathcal B:=\mathcal A^2$, o\`{u} $\mathcal A^2:=\{e:\{0,1\}\rightarrow \mathcal A\}$ est l'ensemble des applications de $\{0,1\}$ dans $\mathcal A$, et ordonnons $\mathcal B$ suivant ses composantes.
Soit $\mathcal B^*$ l'ensemble de tous les mots sur l'alphabet ordonn\'e $\mathcal B$ (c'est \`a dire de toutes les suites d'\'el\'ements de $\mathcal B$). Nous d\'efinissons une application $F$, qui pr\'eserve l'ordre, de $\mathcal B^*$ sur $\mathscr B. \mathcal A$. Ceci est suffisant. En effet, $\mathcal B$ est  belordonn\'e comme produit de deux ensembles belordonn\'e donc, d'apr\`es le th\'eor\`eme d'Higman sur les mots d\'efinis sur un alphabet ordonn\'e \cite{higman52} (voir Th\'eor\`eme \ref{theo:higman}), $\mathcal B^*$ est belordonn\'e. Puisque $\mathscr B. \mathcal A$ est l'image d'un belordre par une application croissante il est belordonn\'e (voir Propri\'et\'es \ref{propr:belordre}.). \\Nous d\'efinissons l'application $F$ comme suit.
 Soit $w:=w(0)w(1)\cdots w(n-1)\in \mathcal B^*$. Posons $F(w):=(\mathcal R, f_w)\in \mathscr B. \mathcal A$ o\`{u} $\mathcal R:= (V_n,  L_{n,1}, L_{n,2})$ et $f_w(i,j):=w(i)(j)$ pour $j\in \{0,1\}$.  Observons d'abord que $w\leq w'$ dans $\mathcal B^*$ implique $F(w)\leq F(w')$ dans $\mathscr B. \mathcal A$. En effet, si $w\leq w'$ il existe un abritement $h$ de la cha\^{\i}ne $0<\cdots<n-1$ dans la cha\^{\i}ne  $0<\cdots<n'-1$ tel que $w(i)\leq w'(h(i))$  pour tout $i<n$.
  Soit alors $\overline h:\{0,\dots, n-1\}\times \{0,1\}\rightarrow \{0,\dots, n'-1\}\times \{0,1\}$ d\'efini en posant $\overline{h}(i,j):=(h(i), j)$. Comme il est facile de le v\'erifier,  $\overline h$ est un abritement de $F(w)$ dans  $F(w')$,
car pour $x=(i,j)\in \{0,\dots, n-1\}\times \{0,1\}$ nous avons
$$f_w(x)=f_w(i,j)=w(i)(j)\leq w'(h(i))(j)=f_{w'}(h(i),j)=f_{w'}(\overline{h}(i,j))=f_{w'}(\overline{h}(x)).$$
Notons ensuite que $F$ est surjective. En effet, si $(\mathcal R, f)\in \mathscr B. \mathcal A$ avec  $\mathcal R:= (V_n,  L_{n,1}, L_{n,2})$ alors le mot  $w:=w(0)w(1)\cdots w(n-1)$ avec $w(i)(j):=f(i,j)$ est tel que  $F(w)=(\mathcal R, f)$.

\end{proof}

\vspace{2mm}

\noindent Avec la Proposition \ref{prop:wqo} et le Th\'eor\`eme \ref{theo:closed}, nous avons:

\begin{corollary}
 La cl\^oture par sommes de la classe des bicha\^{\i}nes critiques est h\'er\'editairement belordonn\'ee et a un nombre fini de bornes.
 \end{corollary}

\noindent Dans \cite{A-A} il est mentionn\'e que cette classe a un nombre fini de bornes.
La s\'erie g\'en\'eratrice de la classe des bicha\^{\i}nes critiques est rationnelle (la classe est couverte par quatre cha\^{i}nes). D'apr\`es le Corollaire \ref{cor:wqoalgebraic} (ou corollaire 13 de \cite {A-A}) sa cl\^oture par sommes est alg\'ebrique.

\vspace{2mm}

%\begin{question}
%\begin{enumerate}
%\item

Il est naturel de se demander si la cl\^oture par sommes de la classe des bicha\^{i}nes critiques\index{bichaine@bicha\^{i}ne!critique} est h\'er\'editairement alg\'ebrique.

\vspace{2mm}

Comme indiqu\'e par V.Vatter \cite{vatter4},\index{Vatter} la r\'eponse \`a cette question est positive. Ceci d\'ecoule d'un r\'esultat plus g\'en\'eral sur "les classes grille-g\'eom\'etriques" (\emph{geometric grid classes} en anglais)\index{classe!grille-g\'eom\'etrique} de permutations et le fait que ces classes renferment les classes de permutations exceptionnelles. En effet, comme il a \'et\'e montr\'e dans  \cite{albert-vatter1} et \cite{albert-vatter2}:

\begin{theorem}\label{theo:classegrill}\cite{albert-vatter1}

Toute classe grille-g\'eom\'etrique est h\'er\'editairement rationnelle et a un nombre fini de bornes.
 \end{theorem}

\begin{theorem}\label{theo:grillgeoalgebrique}\cite{albert-vatter2}

 La cl\^oture par sommes de toute classe  g\'eom\'etriquement en grille est h\'er\'editairement alg\'ebrique.
\end{theorem}

Une classe g\'eom\'etriquement en grille est une classe contenue dans une classe grille-g\'eom\'etrique. Pour ne pas allourdir cette partie, nous renvoyons \`a l'annexe \ref{chapitre:grille-geo} les  d\'efinitions et les rappels de r\'esultats concernant les classes grille-g\'eom\'etriques.

\vspace{2mm}

Le fait que la cl\^oture par sommes de la classe des bicha\^{i}nes critiques soit h\'er\'editairement alg\'ebrique serait \'egalement une cons\'equence de la conjecture pour les classes h\'er\'editaires de structures binaires ordonn\'ees que nous formulons ci-dessous. %Nous conjecturons que la r\'eponse est positive. Ceci sera une cons\'equence de la conjecture pour les classes h\'er\'editaires de structures binaires ordonn\'ees que nous formulons ci-dessous.
En effet, la classe des bicha\^{\i}nes critiques satisfait les hypoth\`eses de cette conjecture.

\begin{conjecture}\label{conjec}
Si $\mathscr D$ est une classe h\'er\'editaire de structures binaires ordonn\'ees ind\'ecomposables qui est h\'er\'editairement belordonn\'ee et h\'er\'editairement alg\'ebrique, alors sa cl\^oture par sommes est h\'er\'editairement alg\'ebrique.
\end{conjecture}

 L'hypoth\`ese que $\mathscr D$ est belordonn\'e ne suffit pas dans la Conjecture \ref{conjec}.

\vspace{2mm}

\noindent En effet, soit $\mathcal P_{\mathbb{Z}}$ le chemin doublement infini dont l'ensemble des sommets est $\mathbb{Z}$ et l'ensemble des ar\^etes est  $E:=\{\{n,m\}\in \mathbb{Z}\times \mathbb{Z}: \vert n-m\vert=1\}$. L'ensemble des ar\^etes $E$ a deux orientations transitives,
e.g. $P:=\{(n,m)\in \mathbb{Z}\times \mathbb{Z}: \vert n-m\vert=1 \; \text{et}\; n\;  \text {est pair}\}$ et son inverse ou duals $P^{-1}$, nous les appelons zigzags.\index{zigzags} Comme ordre, $P$ est l'intersection des ordres lin\'eaires:\\ $L_1:=\cdots<2n<2n-1<2(n+1)<2n+1<\cdots$ et $L_2:=\cdots<2(n+1)<2n+3<2n<2n+1<\cdots.$

\vspace{2mm}

Soit $\mathcal B:= (\mathbb{Z}, L_1,L_2)$ et
$\mathscr D_{osc}:=Ind(\mathcal B)$.

 \begin{lemma}  $\mathscr D_{osc}$ est belordonn\'e mais non h\'er\'editairement belordonn\'e.
 \end{lemma}

 \begin{proof}
 Les membres de $\mathscr D_{osc}$ de taille $n$ sont obtenus en restreignant $\mathcal B$ aux intervalles de taille $n\not =3$ de la cha\^{i}ne $(\mathbb{Z}, \leq)$ (observons d'abord que le graphe $\mathcal P_{\mathbb{Z}}$ est ind\'ecomposable ainsi que toutes ses restrictions aux intervalles de tailles diff\'erentes de $3$ de la cha\^{i}ne $(\mathbb{Z}, \leq)$ et par suite il n'y a pas d'autres restrictions ind\'ecomposables; puis,
utilisons le fait que l'ind\'ecomposabilit\'e d'un graphe de comparabilit\'e d\'ecoule de l'ind\'ecomposabilit\'e de ses orientations \cite{Ke}, et que l'ind\'ecomposabilit\'e d'un ordre de dimension deux d\'ecoule de l'ind\'ecomposabilit\'e des bicha\^{i}nes associ\'ees \`a l'ordre \cite{zaguia}).
A l'isomorphie pr\`es, il y a deux bicha\^{i}nes ind\'ecomposables de taille $n\not=3$:
$\mathcal B_n:=\mathcal B_{\restriction \{0, \dots, n-1\}}$ et $\mathcal {B}_n^{-1}:=\mathcal B_{\restriction \{0, \dots, n-1\}}^{-1}$ o\`{u}
$\mathcal B^{-1}:= (\mathbb{Z}, L_1^{-1},L_2^{-1})$.
Ces deux bicha\^{i}nes  abritent tous les membres de $\mathscr D_{osc}$ ayant des tailles inf\'erieures \`{a} $n$. Etant couvert par deux cha\^{i}nes,
$\mathscr D_{osc}$ est belordonn\'e. Pour montrer que $\mathscr D_{osc}$ n'est pas h\'er\'editairement belordonn\'e, nous pouvons associer \`{a} tout membre ind\'ecomposable de $\mathscr D_{osc}$  le graphe de comparabilit\'e de l'intersection des deux ordres et observer que cette association pr\'eserve la relation d'abritement, m\^{e}me si on rajoute des \'{e}tiquettes.  La classes des graphes obtenue par cette association consiste en les chemins de tailles diff\'erentes de $3$. Elle n'est pas h\'er\'editairement belordonn\'ee. En fait, comme nous pouvons le voir, si une classe $\mathscr G$ de graphes  contient une infinit\'e de chemins de tailles distinctes, alors $\mathscr G.\underline{2}$  n'est pas belordonn\'e, o\`u $\underline{2}$ est l'anticha\^{i}ne \`a deux \'el\'ements $0$ et $1$. En effet, si on \'{e}tiquette les extr\'emit\'es de chaque chemin par $1$ et les autres sommets par $0$,  nous obtenons une anticha\^{i}ne infinie. Donc $\mathscr D_{osc}.\underline{2}$ n'est pas belordonn\'e.
\end{proof}

\begin{lemma}\label{lem:D-her-alg}
$\mathscr D_{osc}$ est h\'er\'editairement alg\'ebrique.
\end{lemma}

\begin{proof}
Comme l'ensemble $\mathscr D_{osc}$ est couvert par deux cha\^{i}nes, toute sous-classe h\'er\'editaire propre de $\mathscr D_{osc}$ est finie, donc alg\'ebrique. La s\'erie g\'en\'eratrice de $\mathscr D_{osc}$ est rationnelle, sa fonction g\'en\'eratrice est donn\'ee par: $$\dfrac{x+x^2-2x^3+2x^4}{1-x}.$$ D'o\`u l'on d\'eduit que $\mathscr D_{osc}$ est h\'er\'editairement alg\'ebrique.
\end{proof}

\begin{lemma}\label{lem:Dnonwqo}
$\sum \mathscr D_{osc}$ n'est pas belordonn\'e.
\end{lemma}

\begin{proof}
C'est d\^u au fait que l'on peut abriter l'ensemble ordonn\'e $\mathscr D_{osc}.\underline{2}$ dans
$\sum \mathscr D_{osc}$ par une application qui pr\'{e}serve l'ordre. Un argument simple consiste \`{a} observer d'abord que la famille $(G_n)_{n\in \mathbb{N}}$, o\`{u} $G_n$ est le graphe obtenu \`{a} partir du chemin \`{a} $n$ sommets  $P_n$ en rempla\c{c}ant ses sommets extr\'emit\'es par un ind\'ependant \`{a} deux sommets, est une anticha\^{\i}ne, puis que ces graphes sont des graphes de comparabilit\'e associ\'es aux membres de $\sum\mathscr D_{osc}$.
\end{proof}

\begin{lemma}
$\sum \mathscr D_{osc}$ est alg\'ebrique mais non h\'er\'editairement alg\'ebrique.
\end{lemma}

\begin{proof}
$\mathscr D_{osc}$ \'etant alg\'ebrique (Lemme \ref{lem:D-her-alg}), d'apr\`es le Corollaire \ref{cor:wqoalgebraic}, sa cl\^oture par sommes, $\sum \mathscr D_{osc}$, est alg\'ebrique, en fait la fonction g\'en\'eratrice $D$ de $\sum \mathscr D_{osc}$ v\'erifie: $$2D^5+2D^4-D^3+(2-x)D^2-D+x=0.$$
Mais d'apr\`es le Lemme \ref{lem:Dnonwqo}, $\sum \mathscr D_{osc}$ n'est pas belordonn\'e, donc n'est pas h\'er\'editairement alg\'ebrique d'apr\`es le Lemme \ref{lem:wqo}. \end{proof}
\bigskip

Les permutations correspondant aux membres de $\mathscr D_{osc}$ s'appellent \emph{oscillations d\'ecroissantes} (\emph{decreasing oscillations} en anglais) voir \cite{brignall-al}. Ces permutations ont les formes suivantes pour $n\geq 4$:\\%Elles ont fait l'objet de nombreuses \'etudes: ajouter des references

\noindent Si $n$ est pair on a:

\qquad $\sigma_1=(n-2)n(n-4)(n-1)(n-6)(n-3)(n-8)\ldots 1.3$ et

\qquad $\sigma_2=(n-1)(n-3)n(n-5)(n-2)(n-7)(n-4)\ldots 4.2$,
avec $\sigma_1=\sigma^{-1}_2,$

\noindent et si $n$ est impair, alors

\qquad $\sigma'_1=(n-1)(n-3)n(n-5)(n-2)(n-7)(n-4)\ldots 1.3$ et

\qquad $\sigma'_2=(n-2)n(n-4)(n-1)(n-6)(n-3)(n-8)\ldots 4.2$,
avec $\sigma'_1=\sigma'^{-1}_1,~\sigma'_2=\sigma'^{-1}_2.$

\vspace{1mm}

Par exemple, pour $n=7$ on a $\sigma'_1=6472513$ et $\sigma'_2=5736142$ et pour $n=8$ on a $\sigma_1=68472513$ et $\sigma_2=75836142.$\\

 Nous rappelons la d\'efinition de l'oscillation d\'ecroissante donn\'ee dans \cite{brignall-al}. La \emph{suite oscillante croissante} est la suite infinie donn\'ee par: $4, 1, 6, 3,8,5,\ldots,2k+2, 2k-1,\ldots$. Une \emph{oscillation croissante} est toute permutation simple qui est contenue dans la suite oscillante croissante. Une \emph{oscillation d\'ecroissante} est la permutation duale d'une oscillation croissante, c'est \`a dire, si $\sigma=a_1a_2\ldots a_n$, son dual est $\sigma^*=a_n\ldots a_2 a_1.$\\

La cl\^oture par abritement $\downarrow\mathscr D_{osc}$ de $\mathscr D_{osc}$ est l'\^{a}ge, $Age (\mathcal B)$, de $\mathcal B$; cet ensemble est rationnel et a quatre obstructions,  sa s\'erie g\'en\'eratrice est $$\dfrac{1-x}{1-2x-x^{3}},$$ la fonction g\'en\'eratrice \'etant la s\'equence A052980 de \cite {Sloane}, commen\c{c}ant par 1, 1, 2, 5, 11, 24.

\vspace{2mm}

\noindent En effet, dans \cite{brignall-al}, il a \'et\'e montr\'e que la classe des permutations contenues dans les oscillations croissantes est $Forb(\{321,2341,3412,4123\}).$ Cette classe a,  pour s\'erie la s\'equence $A052980$ \cite{Sloane} et pour fonction g\'en\'eratrice $$\dfrac{1-x}{1-2x-x^{3}}.$$ Comme, par d\'efinition, une permutation $\sigma$ est une oscillation d\'ecroissante si et seulement si son dual $\sigma^*$ est une oscillation croissante, il est normal que ces deux classes aient la m\^eme fonction profil. Nous pouvons \'egalement en d\'eduire que
la classe des permutations contenues dans les oscillations d\'ecroissantes est $Forb(\{123,1432,2143,3214\}).$

\subsubsection{Questions}
Est-il vrai que:
\begin{enumerate}
\item une classe h\'er\'editaire $\mathscr D$ de structures binaires ordonn\'ees ind\'ecomposables est h\'er\'editairement belordonn\'ee lorsque sa cl\^oture par sommes est h\'er\'editairement alg\'ebrique? Dans \cite{vatter4}, Vatter conjecture qu'une classe de permutations est h\'er\'editairement alg\'ebrique si et seulement si elle est belordonn\'ee.

\item  la s\'erie g\'en\'eratrice d'une classe h\'er\'editaire de structures relationnelles est rationnelle lorsque le profil de cette classe est born\'e par un polyn\^{o}me? Ce r\'esultat \'etant vrai pour les graphes \cite{B-B-S-S} et les tournois \cite{Bou-Pouz}.

\item le profil d'une classe h\'er\'editaire belordonn\'ee de structures relationnelle est born\'e sup\'erieurement par une exponentielle?

\end{enumerate}         
\clearemptydoublepage

%$$$$$$$$$$$$$$$$$$$$$$$$$$$$$$$$$$$$$$$$$$$$$$$$$$$$$$$$$$$$$$$$$$$$$$$$$$$
\part{Minimalité}\label{part:minimalite}
%$$$$$$$$$$$$$$$$$$$$$$$$$$$$$$$$$$$$$$$$$$$$$$$$$$$$$$$$$$$$$$$$$$$$$$$$$$$

\chapter{Notion de minimalit\'e et \^ages de graphes infinis}\label{chap:minimalite}
\section{Introduction}

Dans ce chapitre nous nous int\'eressons au r\^ole que peut jouer la notion de minimalit\'e de la th\'eorie des ensembles ordonn\'es, dans la th\'eorie des relations.
Un ensemble ordonn\'e $P$ est dit \emph{minimal}\index{ensemble!ordonn\'e!minimal} s'il est infini et tout segment initial propre de $P$ est fini.

 L'exemple le plus simple est celui de la cha\^{i}ne $\omega$ des entiers naturels.  Si nous consid\'erons les structures finies, une classe h\'er\'editaire $\mathscr C$ de telles structures, %relationnelles finies, %de signature $\mu$,
 consid\'er\'ees à l'isomorphie près, est dite \emph{minimale}\index{classe!minimale} si elle est infinie et toute sous-classe h\'er\'editaire propre de $\mathscr C$ est finie. Ces classes minimales ont \'et\'e caract\'eris\'ees par Fra\"{\i}ss\'e, elles ont un profil constant \'egal \`a $1$ et si la signature est finie, elles en nombre fini.

Cette notion est particulièrement int\'eressante si $\mathscr C$ est form\'ee de structures ind\'ecomposables et est h\'er\'editaire (dans la classe des ind\'ecomposables). Aussi, nous disons qu'une classe $\mathscr C$ de structures binaires ordonn\'ees est \emph{ind-minimale}\index{classe!ind-minimale} si $\mathscr C$ contient un nombre infini de structures ind\'ecomposables mais toute sous-classe h\'er\'editaire propre de $\mathscr C$ n'en contient qu'un nombre

Nous consid\'erons des classes ind-minimales, c'est à dire des classes h\'er\'editaires de structures finies dont les sous-ensembles d'ind\'ecomposables qu'elles contiennent sont minimaux (dans la classes des structures ind\'ecomposables), ou encore des classes contenant une infinit\'e de structures ind\'ecomposables dont les sous-classes propres n'en contiennent qu'un nombre fini. Contrairement aux classes minimales, le profil d'une telle classe n'est pas n\'ecessairement constant. Nous obtenons certains r\'esultats g\'en\'eraux %que nous affinons aux cas des graphes. Nous montrons
en particulier nous montrons que ces classes sont des \^ages belordonn\'es et qu'elles sont en nombre contin\^upotent.

\section{Classes minimales}
\subsection{Posets minimaux}
    \begin{definition}
    Un ensemble ordonn\'e $P$ est dit \emph{minimal}\index{ensemble!ordonn\'e!minimal} ou de \emph{type minimal}\index{ensemble!ordonn\'e!de type minimal} s'il est infini et tout segment initial propre de $P$ est fini.
    \end{definition}

Ces posets minimaux ont \'et\'e caract\'eris\'es par Pouzet-Sauer dans \cite{P-S}:

%Nous donnons ci-dessous quelques propri\'et\'es \'equivalentes pour une classe minimale, r\'esultat d\^u \`a Pouzet-Sauer \cite{P-S}.

\begin{theorem} \label{minimalposet} Soit $P$ un ensemble ordonn\'e infini. Alors, les propri\'et\'es suivantes sont \'equivalentes:
\begin{enumerate}
\item Tout segment initial propre de $P$ est fini.
\item $P$ est belordonn\'e et tous les id\'eaux distincts de $P$ sont principaux.%\footnote{Un id\'eal $\mathcal I$ est principal s'il est engend\'e par $\{x\}$ pour un $x\in\mathcal I$, c'est \`a dire $\mathcal I=\downarrow\{x\}$.};
\item $P$ n'a pas d'anticha\^{i}ne infinie et tous les id\'eaux distincts de $P$ sont finis.
\item Toute extension lin\'eaire de $P$ est de type $\omega$.
\item $P$ a des nivaux finis de hauteur $\omega$ et pour tout $n<\omega$ il existe  $m<\omega$ tel que tout \'el\'ement de hauteur au plus $n$ est au dessous de tout \'el\'ement de hauteur au moins  $m$.
\item $P$ n'abrite aucun des posets suivants: une anticha\^{i}ne infinie; une cha\^{i}ne de type $\omega ^{\star}$ (le dual de $\omega$); une cha\^{i}ne de type $\omega +1$; la somme directe $\omega\oplus 1$ d'une cha\^{i}ne de type $\omega$ et d'une cha\^{i}ne \`a un \'el\'ement.
 \end{enumerate}
\end{theorem}
 Une mani\`ere d'obtenir des posets minimaux est donn\'ee par ce corollaire \cite{P-S}
\begin{corollary}
Soient $P$ un ensemble ordonn\'e et $n$ un entier. $P$ est l'intersection de $n$ ordres lin\'eaires de type $\omega$ si et seulement si $P$ est l'intersection de $n$ ordres lin\'eaires et $P$ est de type minimal.
\end{corollary}
Nous avons le r\'esultat suivant:
%
%An easy way of obtaining posets with minimal type is given by the following corollary
%\begin {corollary}
%Let $n$ be an integer and $P$ be  a poset.  The order on $P$ is the intersection of $n$ linear
%orders of  order type $\omega$  if and only if $P$ is the intersection of $n$ linear orders  and $P$ has minimal type.
%\end {corollary}
%
\begin{lemma} \label{lem:contains minimal}
Tout ensemble ordonn\'e $P$ qui se d\'ecompose en un nombre d\'enombrable de niveaux  finis contient un segment initial qui est minimal.
\end{lemma}

\begin{proof}
Il suffit d'appliquer le lemme de Zorn \`a l'ensemble $\mathcal J$ des segments initiaux infinis de $P$ ordonn\'e par l'inverse de l'inclusion. Pour cela, nous montrons que $\mathcal J$ est ferm\'e pour les intersections de cha\^{i}nes non vides. %donc il est inductif
 En effet, soit $\mathcal C$ une cha\^{i}ne non vide (par rapport \`a l'inclusion) d'\'el\'ements de $\mathcal J$. Posons $I:= \cap \mathcal C$. Soit $P_n$ le niveau $n$ de $P$ et posons $\mathcal {C}_n:=  \{C_n=C\cap P_n: C\in \mathcal C\}$. Les membres $C_n$ de $\mathcal{C}_n$ sont finis, non vides et totalement ordonn\'es par l'inclusion. Donc,  $I_n:=\cap C_n$ est non vide. Puisque $I=\cup \{ I_n: n\in \mathbb N\}$, $I\in \mathcal J$.
\end{proof}

\subsection{Classes minimales de $\Omega_{\mu}$}
\begin{definition}
Une classe h\'er\'editaire $\mathscr C$ de structures de $\Omega_{\mu}$, consid\'er\'ees \`a l'isomorphie pr\`es, est dite \emph{minimale}\index{classe!minimale} ou de \emph{type minimale}\index{type!minimal} si elle est infinie et toute sous-classe h\'er\'editaire propre de $\mathscr C$ est finie. %Aussi, une classes h\'er\'editaire de structures binaires de $Ind(\Omega_k)$ est dite \emph{minimale} si elle est infinie et toute sous-classe h\'er\'editaire propre est finie.
\end{definition}

Etant donn\'e une structure relationnelle $\mathcal R$ de $\Omega_{\mu}$, son \^age, $Age(\mathcal R)$, est un ensemble ordonn\'e qui est naturellement minimal s'il est infini et totalement ordonn\'e. En effet, dans ce cas, il est ordonn\'e comme la cha\^{i}ne $\omega$ des entiers naturels. Par exemple, si $\mathcal R$ est une cha\^{i}ne, son \^age est minimal. \\L'inverse (et plus encore) est vrai comme le montre le th\'eor\`eme suivant
 essentiellement d\^u \`a Fra\"{\i}ss\'{e} (voir son livre \cite{fraisse} pour une illustration) qui donne la structure d'une classe minimale (voir aussi \cite{pouzet.81, pouzet06}).

  \begin{theorem}
Une classe h\'er\'editaire $\mathscr C$ de structures relationnelles de signature $\mu$ (consid\'er\'ees \`a l'isomorphie pr\`es) est minimale si et seulement si c'est l'\^age d'une structure relationnelle infinie monomorphe.
\end{theorem}

D'apr\`es \emph{Fra\"{\i}ss\'{e}} qui a introduit la  notion de monomorphie en 1954, une structure relationnelle $\mathcal R$ pour laquelle $\varphi_{\mathcal R}(n)=1$ pour tout $n\leq \vert V(\mathcal R)\vert$ est \emph{monomorphe}\footnote{Voir cette notion plus en d\'etail dans le chapitre \ref{chap:monomorphe}.}.

Il s'ensuit imm\'ediatement

\begin{lemma}\label{lem:min-ideal}
Toute classe minimale de $\Omega_{\mu}$ est un id\'eal de $\Omega_{\mu}$.
\end{lemma}
\vspace{2mm}

\vspace{2mm}

Si la signature $\mu$ est finie, le nombre de structures infinies monomorphes (encha\^{i}nables) est fini, les classes minimales de $\Omega_{\mu}$ sont donc en nombre fini. En particulier, si $\mathscr C$ est une classe h\'er\'editaire minimale de graphes,  $\mathscr C$ est l'\^age du graphe infini qui est complet ou vide. Il y a donc deux classes minimales de graphes, l'une form\'ee des graphes finis vides, l'autre de graphes finis complets.

\begin{lemma}
Si l'arit\'e $\mu$ est finie, toute classe h\'er\'editaire $\mathscr C$ de $\Omega_{\mu}$  contient une sous-classe minimale.
\end{lemma}
\begin{proof}
 Ce r\'esultat est une cons\'equence du Th\'eor\`eme \ref{theo:classehered-ideal}, en effet, la preuve de ce dernier montre que $\mathscr C$ contient une section initiale $\mathscr I$ qui est infinie et belordonn\'ee. La collection des sections initiales d'un belordre \'etant bien fond\'ee (Th\'eor\`eme \ref{theo:higman-equivalence}), prendre $\mathscr I'\subseteq\mathscr I$ qui est minimale parmi les sections initiales infinies.
\end{proof}

            \section{Classe ind-minimale}

Cette notion de classe minimale est particuli\`erement int\'eressante si $\mathscr C$ est form\'ee de structures binaires ind\'ecomposables et est h\'er\'editaire (dans la classe des ind\'ecomposables).

L'int\'er\^et des classes h\'er\'editaires minimales  form\'ees de structures ind\'ecomposables apparait \'egalement dans le lemme suivant:

\begin{lemma}
Dans la classe des structures binaires finies de type $k$ ($k$ fini), il existe une correspondance biunivoque entre les classes h\'er\'editaires  minimales de structures ind\'ecomposables et les classes h\'er\'editaires qui contiennent une infinit\'e de structures ind\'ecomposables mais dont toute sous-classe h\'er\'editaire propre n'en contient qu'un nombre fini.
\end{lemma}

\begin{proof}
A une classe h\'er\'editaire  minimale $\mathcal I$ de structures finies ind\'ecomposables on associe la classe h\'er\'editaire $\downarrow\mathcal I$. Si $\mathscr C$ est une classe h\'er\'editaire de structures finies qui contient une infinit\'e de structures ind\'ecomposables mais dont toute sous-classe h\'er\'editaire propre n'en contient qu'un nombre fini on associe la classe $Ind(\mathscr C)$.
\end{proof}

\vspace{2mm}

Par exemple, \`a la classe des chemins finis correspond l'\^age d'un chemin infini. Ces derni\`eres classes seront dites \emph{ind-minimale}.

\begin{definition}
Une classe $\mathscr C$ de $\Omega_k$ est dite \emph{ind-minimale}\index{classe!ind-minimale} si $\mathscr C\cap Ind(\Omega_k)$ est infinie et pour toute sous-classe h\'er\'editaire propre $\mathscr C'$ de $\mathscr C$, la sous-classe $\mathscr C'\cap Ind(\Omega_k)$ est finie.
\end{definition}

Nous avons, de mani\`ere naturelle le lemme suivant:

\begin{lemma}\label{lem:equiv-min}
 Soit $\mathscr C$ une classe h\'er\'editaire de $Ind(\Omega_k)$, alors les trois propri\'et\'es suivantes sont \'equivalentes:
 \begin{enumerate}
 \item[$(i)$] $\mathscr C$ est minimale dans $Ind(\Omega_k)$.
  \item[$(ii)$] $\downarrow\mathscr C$ est ind-minimale.
  \item[$(iii)$] $\mathscr C$ est infinie et toute sous-classe propre $\mathscr C'$ de  $\downarrow\mathscr C$ a un nombre fini d'ind\'ecomposables tous de taille finie.
\end{enumerate}
\end{lemma}
\begin{proof}
 Si $\mathscr C$ est minimale dans $Ind(\Omega_k)$ alors toute sous-classe propre de $\downarrow\mathscr C$ contient une sous-classe propre de $\mathscr C$ qui est finie, donc $\downarrow\mathscr C$ est ind-minimale, d'o\`u $(i)$ implique $(ii)$. Avec $(ii)$ comme hypoth\`ese, $\mathscr C$ est infinie.    Soit $\mathscr C'$ une sous-classe propre de $\downarrow\mathscr C$, donc $\mathscr C'$ a un nombre fini d'ind\'ecomposables qui sont de tailles finies (car $\mathscr C$ \'etant h\'er\'editaire dans $Ind(\Omega_k)$, leur nombre serait infini, d'apr\`es Schmerl et Trotter, si ce n'\'etait pas le cas), %Donc, le nombre de structures pouvant s'abriter dans chacun de ces ind\'ecomposables est fini, donc $\mathscr C'$ fini,
 d'o\`u $(iii)$. Supposons $(iii)$ alors $\mathscr C$ est infinie et toute sous-classe propre de $\mathscr C$ induit une sous-classe propre de $\downarrow\mathscr C$ qui a un nombre fini d'ind\'ecomposables, d'o\`u $(i)$.%Inversement, toute sous-classe propre de $\mathscr C$ induit une sous-classe propre de $\downarrow\mathscr C$ qui a un nombre fini d'indecomposable, d'o\`u le r\'esultat.
\end{proof}

\vspace{2mm}

Nous avons  \'egalement le r\'esultat suivant:
\begin{theorem}\label{minimal}
Soit $\mathscr D$ une classe h\'er\'editaire infinie de structures ind\'ecomposables finies.  Toute sous-classe de $\mathscr D$ contient une sous-classe minimale si et seulement si $\mathscr D$ ne contient qu'un nombre fini de membres de taille $1$ et $2$. %contenant un nombre fini de membres de taille $1$ ou $2$.
%Alors $\mathscr D$ contient une sous-classe h\'er\'editaire minimale.
\end{theorem}

La preuve du Th\'eor\`eme \ref{minimal} repose sur le lemme suivant:

\begin{lemma}\label{lem:niveaufini}
Tout niveau de $\mathscr D$ est fini.
\end{lemma}

L'arit\'e des structures de $\mathscr D$ n'est pas fini, à priori (si elle est finie, le Lemme \ref{lem:niveaufini} devient \'evident). Il existe beaucoup d'exemples int\'eressants de structures relationnelles dont l'arit\'e n'est pas finie, nous pouvons citer des exemples qui viennent des espaces m\'etrique (voir \cite{Dolho-Pou-Sa-Sau, Dolho-Laf-Pou-Sau}).

Consid\'erons un espace m\'etrique $\mathbb M:=(M,d)$ où $d$ est une distance sur l'ensemble $M$. Nous pouvons associer une structure relationnnelle binaire $\mathcal R_{\mathbb M}:=(M, (\rho_{r})_{r\in \mathbb Q^+})$ avec $\rho_{r}(x,y):=1$ si et seulement si $d(x,y)\leq r$. Avec cette d\'efinition, $$d(x,y)=inf\{r\in\mathbb Q^+/ \rho_r(x,y)=1\}.$$
Ainsi, l'espace m\'etrique $\mathbb M$ peut-\^etre retrouv\'e \`a partir de la structure $\mathcal R_{\mathbb M}$. Il s'ensuit \'egalement que pour deux espaces m\'etriques $\mathbb M:=(M,d)$ et $\mathbb M':=(M',d')$, une application $f: M\rightarrow M'$ est une isom\'etrie de $\mathbb M$ dans $\mathbb M'$ si et seulement si $f$ est un abritement de $\mathcal R_{\mathbb M}$ dans $\mathcal R_{\mathbb M'}$ \cite{Dolho-Pou-Sa-Sau}. Nous rappelons que $f$ est une \emph{isom\'etrie} de $\mathbb M$ dans $\mathbb M'$ si $d'(f(x),f(y))=d(x,y)$ pour tous $x,y\in M$.

\vspace{2mm}

 Le Th\'eor\`eme \ref{minimal} est bas\'e sur le r\'esultat de Schmerl et Trotter \cite{S-T}, \`a partir duquel d\'ecoule la relation ci-dessous entre la hauteur de toute structure $\mathcal{R}$ de $Ind(\Omega_{k})$ et son ordre $\vert \mathcal{R} \vert$.
%           $$h(\mathcal{R}) \leq \vert \mathcal{R} \vert  \leq 2(h(\mathcal{R})-1).$$
Dans leur papier,  Schmerl et Trotter ont donn\'e des exemples de structures critiques ind\'ecomposables dans les classes des graphes, des ensembles ordonn\'es, des tournois, des graphes orient\'es et des structures relationnelles binaires (voir D\'efinition \ref{def:critique} pour une structure critique). L'ensemble des structures critiques dans chacune de ces classes est une union finie de cha\^{i}nes. Consid\'erons la classe h\'er\'editaire $\mathscr D$. D\'ecomposons-l\`a en niveaux, dans le niveau $i$, pour $i\leq 2$,
sont les structures d'ordres un et deux. Nous avons une relation entre la hauteur\index{hauteur} de toute structure $\mathcal{R}$ de  $\mathscr D$, not\'ee, $h(\mathcal{R})$ et son ordre $\vert \mathcal{R} \vert$:
\begin{enumerate}
\item Si $\mathscr{D}$ est compos\'ee de structures critiques\index{structure relationnelle!critique} auxquelles nous rajoutons les structures ind\'ecomposables d'ordres un et deux, alors, \`a partir de l'ordre $n>2$, nous avons uniquement des structures ind\'ecomposables d'ordres pairs.
    Nous obtenons la relation suivante:
      $$h(\mathcal{R}) \leq \vert \mathcal{R} \vert  \leq 2(h(\mathcal{R})-1).$$
    En fait, si $\mathcal{R}$ a un ou deux \'el\'ements, alors $\vert \mathcal{R} \vert=h(\mathcal{R})$ et si $\vert \mathcal{R} \vert>2$, nous avons
     $\vert \mathcal{R} \vert =2(h(\mathcal{R})-1).$
\item si $\mathscr{D}$ est compos\'ee de structures ind\'ecomposables non critiques, alors, d'apr\`es le Corollaire \ref{cor:indec-schmer-trott} nous avons des membres de tout ordre (except\'e, peut-\^etre, pour l'ordre trois si $\mathscr D$ est une classe de graphes ou de bicha\^{i}nes), nous obtenons la relation suivante:
      $$h(\mathcal{R}) \leq \vert \mathcal{R} \vert  \leq h(\mathcal{R})+1.$$
\end{enumerate}
De ce fait d\'ecoule le Lemme \ref{lem:niveaufini}. % si un niveau n'est pas fini, les structures ne peuvent pas s'abriter dans les niveaux sup_'erieur qui sont finis car critiques.

\vspace{2mm}

\textbf{\emph{Preuve du Th\'eor\`eme \ref{minimal}.}}  D'apr\`es le Lemme \ref{lem:contains minimal}, il suffit de montrer que $\mathscr D$ se d\'ecompose en une infinit\'e de niveaux finis si et seulement si $\mathscr D$ ne contient qu'un nombre fini de membres de taille $1$ et $2$.

La condition n\'ecessaire est \'evidente, pour la condition suffisante,
puisque les membres de $\mathscr D$ sont finis, cet ensemble est bien fond\'e, donc il se d\'ecompose en niveaux. Pour $i\leq 2$, le $i$-\`eme niveau de $\mathscr D$ comporte les structures relationnelles \`a  $i$ \'el\'ements qui sont en nombre fini par hypoth\`ese. Si $i\geq 3$, le r\'esultat vient du Lemme \ref{lem:niveaufini}.

\hfill $\Box$

\begin{theorem}\label{theo:equiv-min}
Une classe h\'er\'editaire $\mathscr C$ de $\Omega_k$ est ind-minimale si et seulement si
    \begin{enumerate}
    \item $\mathscr C\cap Ind(\Omega_k)$ est minimale dans $Ind(\Omega_k)$.
    \item $\mathscr C:=\downarrow(\mathscr C\cap Ind(\Omega_k))$.
    \end{enumerate}
\end{theorem}

\begin{proof}
Condition n\'ecessaire: supposons que $\mathscr C$ soit ind-minimale. Alors $\mathscr C\cap Ind(\Omega_k)$ est infini. Soit $\mathcal I$ une sous-classe h\'er\'editaire  propre de $\mathscr C\cap Ind(\Omega_k)$. D'un c\^ot\'e $\downarrow \mathcal I$ est une sous-classe propre de $\mathscr C$, donc $(\downarrow\mathcal I)\cap Ind(\Omega_k)$ est finie car $\mathscr C$ est ind-minimale. Mais $(\downarrow\mathcal I)\cap Ind(\Omega_k)=\mathcal I$, il s'ensuit que $\mathscr C\cap Ind(\Omega_k)$ est minimale dans $Ind(\Omega_k)$. D'un autre c\^ot\'e, nous avons $\downarrow(\mathscr C\cap Ind(\Omega_k))\subseteq\mathscr C$, donc si $\mathscr C\neq\downarrow(\mathscr C\cap Ind(\Omega_k))$, alors $\downarrow(\mathscr C\cap Ind(\Omega_k))$ est une sous-classe propre de $\mathscr C$ ayant un nombre infini d'ind\'ecomposables, ce qui contredit son ind-minimalit\'e. Donc $\mathscr C:=\downarrow (\mathscr C\cap Ind(\Omega_k))$.

Condition suffisante: avec les assertions $1.$ et $2.$ du Th\'eor\`eme, $\mathscr C$ est infinie. Si $\mathscr C$ n'est pas ind-minimale, alors il existe une sous-classe propre $\mathscr C'$ de $\mathscr C$ telle que $\mathscr C'\cap Ind(\Omega_k)$ est infinie. Or $\mathscr C'\cap Ind(\Omega_k)$ est une sous-classe propre de $\mathscr C\cap Ind(\Omega_k)$, il s'ensuit que $\mathscr C\cap Ind(\Omega_k)$ n'est pas minimale dans $Ind(\Omega_k)$, ce qui contredit l'hypoth\`ese. %(pour le deuxieme point, je ne sais pas s'il est suffisant??)
\end{proof}

\vspace{3mm}

Une propri\'et\'e g\'en\'erale des classes minimales de $Ind(\Omega_k)$ est le r\'esultat suivant.

\begin{theorem} \label{thm:main1}
Soit $\mathscr C$ une classe minimale de $Ind(\Omega_k)$. Alors
 \begin{enumerate}
\item $\downarrow\mathscr C$ est belordonn\'ee pour l'ordre d'abritement.
\item Il existe une structure ind\'ecomposable infinie $\mathcal{R}$ telle que $Age (\mathcal{R})= \downarrow \mathscr C$.
\end{enumerate}
\end{theorem}

\begin{proof}
Montrons la premi\`ere partie du th\'eor\`eme. La classe $\mathscr C$ \'etant minimale dans $Ind(\Omega_k)$, elle est belordonn\'ee d'apr\`es l'implication $(2)$ du Th\'eor\`eme \ref{minimalposet}. Pour voir que $\downarrow\mathscr C$ est belordonn\'ee, il faut et il suffit de voir que pour chaque structure $\mathcal S\in\downarrow\mathscr C$, la sous-classe $\mathscr C':= \text{\it Forb}(\mathcal S)\cap\downarrow\mathscr C$ est belordonn\'ee. %car une antichaine de $\downarrow\mathscr C'$ est une antichaine de $\downarrow\mathscr C$ et une antichaine infinie $S_0, S_1,\dots,S_n,\dots$ de $\downarrow\mathscr C$ induit une antichaine infinie de $\downarrow\mathscr C'$ pour $\mathcal S=S_0$.
Pour cela nous utilisons la Proposition \ref{prop:wqo}.\\ Soit $\mathcal S\in\downarrow\mathscr C$. Comme $\mathscr C'\subsetneqq\downarrow\mathscr C$ et $\downarrow\mathscr C$ est ind-minimale d'apr\`es le Lemme \ref{lem:equiv-min}, alors $\mathscr C'$ ne contient qu'un nombre fini d'ind\'ecomposables, elle est donc belordonn\'ee.

Pour la seconde partie du th\'eor\`eme, $\mathscr C$ \'etant minimale, c'est un id\'eal (Th\'eor\`eme \ref{minimalposet}) %(Lemme \ref{lem:min-ideal})
de $Ind(\Omega_k)$, il s'ensuit que $\downarrow \mathscr C$ est un id\'eal de $\Omega_k$. Il existe alors une suite $\mathcal P=\mathcal R_0\leq \mathcal R_1\leq\ldots\leq\mathcal R_n\leq\ldots$ d'\'el\'ements de $\mathscr C$ tel que tout \'el\'ement de $\downarrow \mathscr C$ s'abrite dans un \'el\'ement de $\mathcal P$. La cha\^{i}ne $\mathcal P$ est, d'apr\`es le Th\'eor\`eme \ref{minimalposet}, de type $\omega$, donc la limite $\mathcal R_{\omega}$ de $\mathcal P$ a pour \^age $\downarrow \mathscr C$.
\end{proof}

\begin{remark}
Une classe ind-minimale n'est pas n\'ecessairement h\'er\'editairement belordonn\'ee. Si $P_{\infty}$ est le chemin infini sur $\mathbb N$, alors son \^age $Age(P_{\infty})$ est ind-minimal mais pas h\'er\'editairement belordonn\'e (voir preuve du Lemme \ref{lem:Dnonwqo}). %
\end{remark}

\vspace{3mm}

Nous avons vu dans la section pr\'ec\'edente que le nombre de classes minimales de $\Omega_{\mu}$ est fini, la situation est diff\'erente pour les classes ind-minimales, en effet nous avons le r\'esultat suivant:

\begin{theorem} \label{many} Il y a un nombre contin\^upotent %($2^{\aleph_0}$)
de classes h\'er\'editaires minimale form\'ees de graphes ind\'ecomposables qui sont dirig\'es sans boucles ou non dirig\'es avec des boucles.
\end{theorem}

%Remarquons d'abord qu'un graphe $G:=(V,E)$ non dirig\'e avec des boucles peut-\^etre vu comme une structure relationnelle $\mathcal R:=(V,\rho,u)$ où $\rho$ est une relation binaire irr\'eflexive et sym\'etrique et $u$ une relation unaire form\'ee des \'el\'ements $x\in V$ tels que $(x,x)\in E$.

\vspace{2mm}

 Ce r\'esultat vient de la dynamique symbolique. La notion de minimalit\'e apparait quand nous consid\'erons les mots sur un alphabet fini et que nous ordonnons ces mots par l'ordre des facteurs. En effet, les ensembles de mots finis qui sont minimaux sont les facteurs finis de mots infinis uniform\'ement r\'ecurrents. Si par exemple nous prenons des mots de Sturm de pentes diff\'erentes nous obtenons des ensembles de facteurs finis diff\'erents. Maintenant, si nous codons chaque mot de Sturm par la relation de cons\'ecutivit\'e sur $\NN$ et une relation unaire, les ind\'ecomposables finis de cette multirelation correspondent aux facteurs finis du mot de Sturm. Ainsi nous obtenons un nombre contin\^upotent d'\^ages ind-minimaux.

\vspace{1mm}

 Pr\'ecisons tout cela, rappelons quelques d\'efinitions et r\'esultats \cite{All-Sha, Lo}.\\

Soit $A$ un ensemble fini. %Une suite $\alpha=\alpha_0\alpha_1\dots $ telle que $\alpha_i\in A$ est un mot\index{mot}  sur  l'alphabet $A$.
Un mot $\alpha$ sur l'alphabet\index{alphabet} $A$ est une application \`a valeur dans $A$ dont le domaine est soit $\{0,1,\dots,n-1\}$, auquel cas le mot est fini\index{mot!fini} et sa \emph{longueur}\index{mot!longueur d'un -}, not\'ee $\vert \alpha\vert$, est  $n$ et le mot $\alpha$ s'\'ecrit simplement $\alpha_0\alpha_1\dots\alpha_{n-1}$, soit $\NN$, soit $\NN^-$, l'ensemble des entiers n\'egatifs, soit $\mathbb Z$ auxquels cas le mot est soit infini \`a droite, soit infini\index{mot!infini} \`a gauche, soit infini des deux côt\'es. Le mot de longueur nulle est le mot \emph{vide}.\index{mot!vide}

\vspace{2mm}

D\'esignons par $A^{<\NN}$, respectivement $A^{\NN}$, l'ensemble des mots finis,\index{mot!fini} resp. infinis\index{mot!infini}, de domaine $\NN$, de $A$ et  $A^{\leq \NN}:=A^{<\NN}\cup  A^{\NN}$.  Si $u\in A^{<\NN}$ et  $v\in A^{\leq \NN}$,  nous notons $u\leq v$ et disons que  $u$ est un \emph{facteur}\index{mot!facteur d'un -} de $v$ s'il existe $u_1\in A^{<\NN}$  et $u_2\in A^{\leq \NN}$ tels que $v=u_1\cdot u\cdot u_2$ (l'op\'eration repr\'esent\'ee par un point \'etant la concat\'enation\index{mot!concat\'enation de -} des mots, c'est à dire l'op\'eration qui à deux mots $u,v$, où $u$ est fini, associe le mot $u\cdot v=uv$). Donc un \emph{facteur} d'un mot $\alpha$ est toute sous-suite finie contig\"{u}e de $\alpha$. Si $u$ est un facteur de $\alpha$ tel que $\alpha=uv$ alors $u$ est dit \emph{pr\'efixe}\index{mot!pr\'efixe d'un -} de $\alpha$. D\'esignons par $Fac(u)$ l'ensemble des facteurs finis du mot $u$.

\vspace{2mm}

 La fonction \emph{shift}\index{fonction!shift} (ou \emph{d\'ecalage})\index{fonction!d\'ecalage} est la fonction $\sigma:A^{\NN}\rightarrow A^{\NN}$ d\'efinie, pour tout mot infini $u=u_0u_1\dots$, par $\sigma(u_0u_1\dots)=u_1u_2\dots$.\\ Si $u,v\in A^{\NN}$, alors $u\leq v$ si $u=\sigma^n(v)$ pour un entier $n$.

\vspace{1mm}

La relation $\leq$ est l'\emph{ordre des facteurs}.\index{ordre!des facteurs} Cette relation est un pr\'eordre sur $A^{\leq \NN}$ et un ordre sur $A^{<\NN}$.

\vspace{2mm}

Muni de la topologie produit, $A^{\NN}$ est un espace topologique compact d\'esign\'e par \emph{espace de Cantor} ($A$ \'etant muni de la topologie discr\`ete).
%Si $\alpha=uv$ alors $u$ est dit \emph{pr\'efixe} de $\alpha$.

\vspace{2mm}

Si $u$ est un mot fini non vide, le mot $u^{\omega}=uuu\ldots$ est dit \emph{p\'eriodique}.\index{mot!p\'eriodique} Un mot infini de la forme $vu^{\omega}$, où $v$ est fini, est dit \emph{\'eventuellement-p\'eriodique}.\index{mot!\'eventuellement-p\'eriodique} Un mot infini est dit \emph{ap\'eriodique}\index{mot!ap\'eriodique} s'il n'est pas \'eventuellement-p\'eriodique.

\vspace{1mm}

$u$ est p\'eriodique si et seulement si $u=\sigma^p(u)$ pour un $p>0$. Si deux mots s'abritent mutuellement alors ils sont p\'eriodiques, c'est le cas pour les mots $(01)^{\omega}$ et $(10)^{\omega}$.

\vspace{1mm}

Un mot $u$ est \'eventuellement p\'eriodique si et seulement si le nombre $Fac_n(u)$ de facteur de longueur $n$ de $u$ est born\'e pour tout entier $n$.
%$Fac_n(u)\leq \vert xv\vert$ où $u=xv^{\omega}$   thm 1.3.13 p19 \cite{Lo}

\vspace{1mm}

Un mot infini $\alpha$ est dit \emph{r\'ecurrent}\index{mot!r\'ecurrent} si tout facteur fini de $\alpha$ apparait une infinit\'e de fois. Un mot infini $\alpha$ est r\'ecurrent si et seulement si pour tout $u\in Fac(\alpha)$ il existe $v\in Fac(\alpha)$ tel que $u\cdot v\cdot u\in Fac(\alpha)$.

\vspace{1mm}

Un mot infini $\alpha$ est dit \emph{uniform\'ement r\'ecurrent}\index{mot!uniform\'ement r\'ecurrent} (dit aussi  \emph{presque-p\'eriodique})\index{mot!presque-p\'eriodique} si pour tout facteur fini $u$ de $\alpha$ il existe un entier $r$ tel que tout facteur de longueur $r$ contient $u$. Il est clair qu'un mot uniform\'ement r\'ecurrent est r\'ecurrent.

\vspace{1mm}

Il est clair que tout mot infini p\'eriodique est uniform\'ement r\'ecurrent et si un mot infini r\'ecurrent est \'eventuellement-p\'eriodique alors il est p\'eriodique.

\vspace{2mm}

Un sous-ensemble $X$ de $A^{\NN}$ est dit invariant pour le shift si $\sigma(X)\subseteq X$. Un sous-ensemble $X$, ferm\'e (topologiquement) et invariant pour le shift est dit \emph{minimal} s'il est non vide et minimal pour l'inclusion parmi les sous-ensembles ferm\'es et invariants pour le shift (si $T$ est ferm\'e invariant pour le shift et $T\subset X$ alors $T=\varnothing$ ou $T=X$).

En particulier, \'etant donn\'e un mot $x\in A^{\NN}$, le plus petit ensemble ferm\'e et invariant pour le shift contenant $x$ est l'ensemble $S(x)=\{y\in A^{\NN}/ Fac(y)\subseteq Fac(x)\}$ (\cite{Lo} page 26 ou \cite{All-Sha} Proposition 10.8.9 page 327).

\vspace{1mm}

Nous rappelons ce r\'esultat (\cite{Lo} Theorem 1.5.9 page 26 ou \cite{All-Sha} Theorem 10.9.4 page 331):

\begin{theorem}\label{theorem mot} %(\cite{Lo})
Soit $\alpha\in A^{\NN}$, les trois assertions suivantes sont \'equivalentes:
\begin{enumerate}
\item $\alpha$  est uniform\'ement r\'ecurrent,
\item $S(\alpha)$ est minimal,
\item $Fac(\beta)=Fac(\alpha)$ pour tout $\beta\in S(\alpha)$.
\end{enumerate}
\end{theorem}

Nous avons \'egalement le r\'esultat suivant (\cite{Lo} Theorem 1.5.11):
\begin{theorem}\label{theorem shift recur}
Tout ensemble non vide ferm\'e et invariant pour le shift contient un mot uniform\'ement r\'ecurrent.
\end{theorem}

%que tout ensemble non vide ferm\'e et invariant pour le shift contient un mot uniform\'ement r\'ecurrent ,
Ce dernier th\'eor\`eme  signifie que tout ensemble non vide ferm\'e et invariant pour le shift qui est minimal est de la forme $S(x)$ pour un mot uniform\'ement r\'ecurrent $x$.

\vspace{2mm}

 Soit $X$ un sous-ensemble non vide de $A^{\NN}$ qui est ferm\'e et invariant pour le shift, associons \`a $X$ l'ensemble $\mathcal A(X)$ des mots finis $u$ tel que $u$ est un pr\'efixe d'un mot de $X$. L'ensemble $X$ \'etant invariant pour le shift, l'ensemble $\mathcal A(X)$ est l'ensemble des facteurs finis des mots de $X$. Ordonnons $\mathcal A(X)$ par l'ordre des facteurs, nous avons alors:

\begin{theorem}\label{theorem motmini} (\cite{pouzet12})%thm6.1 p 333

L'ensemble ordonn\'e $\mathcal A(X)$ est de type minimal si et seulement si $X$ est un ensemble ferm\'e, invariant pour le shift et minimal.
\end{theorem}

\begin{proof}
Si $X$ n'est pas minimal, il contient un mot uniform\'ement r\'ecurrent $y$ (Th\'eor\`eme \ref{theorem shift recur}) avec $S(y)\neq X$. %  sous-ensemble propre $X'$ non vide, ferm\'e et invariant pour le shift.
Donc $\mathcal A(S(y))\varsubsetneq\mathcal A(X)$ et $\mathcal A(X)$ n'est pas de type minimal. % $\mathcal A(X)$ est infini car form\'e de pr\'efixes de mots infinis

Inversement, si $\mathcal A(X)$ n'est pas de type minimal alors il poss\`ede un segment initial propre $\mathcal I$ qui est infini (Th\'eor\`eme \ref{minimalposet}). Soit $F$ l'ensemble des \'el\'ements minimaux de $(\mathcal A(X)\setminus\mathcal I)\cup (A^{<\NN}\setminus\mathcal A(X))$, alors $\mathcal I=Forb(F)\cap A^{<\NN}$. Soit $Y=Forb(F)\cap A^{\NN}$ l'ensemble des mots infinis $z$ tel que aucun facteur de $z$ n'est dans $F$. L'ensemble $Y$ est un ensemble ferm\'e et invariant pour le shift (voir \cite{Lo} Proposition 1.5.1 page 23) et $\mathcal A(Y)=\mathcal I$, donc $Y$ est non vide et $Y\varsubsetneq X$ ce qui contredit la minimalit\'e de $X$. Il s'ensuit que $\mathcal A(X)$ est de type minimal.
\end{proof}

\vspace{2mm}

Tenant compte du fait que tout ensemble ferm\'e et invariant pour le shift qui est minimal est de la forme $S(x)$, pour un mot uniform\'ement r\'ecurrent $x$ et du fait que $\mathcal A(S(x))=Fac(x)$, nous pouvons d\'eduire \`a partir des Th\'eor\`emes \ref{theorem mot} et \ref{theorem motmini} le r\'esultat suivant:

\begin{theorem}\label{theo:Facteur-minimal}
Un mot $\alpha$ sur un alphabet fini $A$ est uniform\'ement r\'ecurrent si et seulement si $Fac(\alpha)$, muni de l'ordre des facteurs, est de type minimal.
\end{theorem}

En outre, le th\'eor\`eme suivant montre que toute section initiale de mots finis qui est de type minimal est l'ensemble des facteurs finis d'un mot uniform\'ement r\'ecurrent.

\begin{theorem} (\cite{pouzet12})

Toute section initiale infinie de $A^{<\NN}$, ordonn\'e par l'ordre des facteurs, contient l'ensemble des facteurs d'un mot uniform\'ement r\'ecurrent.
\end{theorem}

\begin{proof}
Soit $\mathcal I$ une section initiale infinie de $A^{<\NN}$. L'ensemble $A^{<\NN}$ \'etant bien fond\'e, $\mathcal I$ peut-\^etre d\'ecompos\'ee en niveaux, le niveau $n$ contenant les facteurs de longueur $n$ des mots de niveaux sup\'erieurs. L'alphabet $A$ \'etant fini, les niveaux de $\mathcal I$ sont finis. D'apr\'es le Lemme \ref{lem:contains minimal}, $\mathcal I$ contient une section initiale minimale $\mathcal J$. L'ensemble $A^{<\NN}$ \'etant ordonn\'e par l'ordre des facteurs, il existe un mot infini, disons $\alpha$, tel que $Fac(\alpha)\subseteq\mathcal J$. Comme $\mathcal J$ est de type minimal, $\mathcal J\subseteq Fac(\alpha)$ et donc $\alpha$ est uniform\'ement r\'ecurrent.

\end{proof}

\vspace{2mm}

 Posons $A=\{0,1\}$. Un mot sur l'alphabet $A$ est dit \emph{binaire}.\index{mot!binaire} Identifions chaque sous-ensemble de $\NN$ avec sa fonction caract\'eristique, c'est donc un mot binaire.

\vspace{1mm}

Soit $\alpha$ un mot binaire (fini ou infini). Associons la structure relationnelle $C_{\alpha}:=(E,c,u_{\alpha})$, où $E$ est soit $\NN$ si $\alpha$ est infini, ou l'ensemble $\{0,1,\dots,n-1\}$ si $\alpha$ est fini de longueur $n$, $c$ est la relation de cons\'ecutivit\'e sur $E$ et $u_{\alpha}$ une relation unaire sur $E$,  avec $u_{\alpha}={\alpha}^{-1}(1)=\{i\in E;~ \alpha_i=1\}$. Une relation unaire $u$ pouvant-\^etre consid\'er\'ee comme une relation binaire $b$ telle que $x\in u$ si et seulement si $(x,x)\in b$, la structure $C_{\alpha}$ est une structure binaire de type $2$.

\vspace{1mm}

 Remarquons que $C_{\alpha}$ est ind\'ecomposable, que $\beta\leq \alpha$  si et seulement si $C_{\beta}$ est une restriction de $C_{\alpha}$ \`a un intervalle de $\NN$ et que $Fac(\alpha)$ correspond \`a l'ensemble des restrictions ind\'ecomposables finies de $C_{\alpha}$. Ainsi, $Fac(\alpha)$ est minimal dans $A^{<\NN}$ muni de l'ordre des facteurs si et seulement si l'\^age de $C_{\alpha}$ est ind-minimal dans $\Omega_2$, ce qui entra\^{i}ne  im\'ediatement:

 \begin{theorem}\label{theo:mot-structure}
 Un mot binaire infini $\alpha$ est uniform\'ement r\'ecurrent si et seulement si l'\^age de la structure $C_{\alpha}$ est ind-minimal.
 \end{theorem}

Remarquons que deux mots infinis peuvent avoir le m\^eme ensemble de facteurs finis sans \^etre \'egaux, en fait si $u$ est uniform\'ement r\'ecurrent, le Th\'eor\`eme \ref{theorem mot} assure que $Fac(v)=Fac(u)$ pour tout $v\in S(u)$ et l'ensemble $S(u)$ est contin\^upotent si et seulement si $u$ n'est pas p\'eriodique (\cite{All-Sha} Theorem 10.8.12, page 328), il y a donc un nombre contin\^upotent de mots avec un m\^eme ensemble de facteurs.

\vspace{1mm}

 D'un autre c\^ot\'e et comme nous pouvons le voir ci-dessous, il y a \'egalement  un ensemble contin\^upotent  de mots uniform\'ement r\'ecurrents dont les ensembles de facteurs sont distincts. Un tel ensemble peut-\^etre obtenu en prenant des mots de Sturm de pentes distinctes.

\vspace{2mm}

Un mot sur un alphabet binaire\index{alphabet!binaire} $\{0,1\}$  est un mot de \emph{Sturm}\index{mot!de Sturm} si pour tout entier $n$, le nombre de ses facteurs distincts de longueur $n$ est $n+1$. Une propri\'et\'e caract\'eristique des mots de Sturm est que le nombre de $1$ dans deux facteurs finis $u$, $u'$ de m\^eme longueur diff\`ere d'au plus un.

\vspace{1mm}

Soit un mot binaire\index{mot!binaire} fini $x$, la \emph{hauteur}\index{mot!binaire!hauteur d'un -} de $x$ est le nombre $h(x)$ de lettres \'egale \`a $1$ dans $x$ et la \emph{pente}\index{mot!binaire!pente d'un -} de $x$ est la valeur (comprise dans l'intervalle $[0,1]$), $\pi(x)=\frac{h(x)}{\vert x\vert}$ (fr\'equence des lettres \'egales \`a $1$ dans $x$).

\vspace{1mm}

Pour un mot binaire infini $x$, la pente est donn\'ee par $\overline{\pi}(x)=\underset{n\rightarrow \infty}{lim}\pi(x_n)$ où $x_n$ est le pr\'efixe de longueur $n$ de $x$ (une suite $(u_n)_{n\geq 0}$ de mots finis de $A^{<\NN}$ converge vers un mot infini $x$ si tout pr\'efixe de $x$ est le pr\'efixe de tous, sauf peut-\^etre un nombre fini, les mots $u_n$, le mot $x$ est unique). La pente d'un mot de Sturm\index{mot!de Sturm} est irrationnelle.

\vspace{2mm}

 Soit $\alpha$ un nombre irrationnel compris entre $0$ et $1$  et $S_{\alpha}$ l'ensemble des mots de Sturm de pente $\alpha$. Les mots de Sturm sont uniform\'ements r\'ecurrents et ont les propri\'et\'es suivantes (\cite{Lo} Proposition 2.1.5 page 43, Proposition 2.1.18 page 53 et Proposition 2.1.25 page 56):

 \begin{proposition}\label{prop:sturm}

 \begin{enumerate}
 \item Les mots de Sturm sont ap\'eriodiques.
 \item L'ensemble $S_{\alpha}$ est un sous-ensemble minimal de  $A^{\NN}$ pour tout nombre irrationnel $\alpha$ de $(0,1)$.
 \end{enumerate}
 Soient $x$ et $y$ deux mots de Sturm,
\begin{enumerate}
\item Si $x$ et $y$ ont une m\^eme pente alors $Fac(x)=Fac(y)$.
\item Si $x$ et $y$ ont des pentes distinctes alors $Fac(x)\cap Fac(y)$ est fini.
\end{enumerate}
\end{proposition}

\vspace{1mm}

De la Proposition \ref{prop:sturm} nous d\'eduisons que deux ensemble $S_{\alpha}$ et $S_{\alpha'}$ sont distincts si et seulement si les valeurs de $\alpha$ et $\alpha'$ sont distinctes. Il s'ensuit %en prenant des valeurs de $\alpha$ irrationnelle dans $(0,1)$.

\begin{lemma}
Il existe un nombre contin\^upotent de mots uniform\'ement r\'ecurrents dont les ensembles de facteurs sont distincts.
\end{lemma}

\vspace{2mm}

\textbf{Preuve du Th\'eor\`eme \ref{many}.}
Posons $A=\{0,1\}$. Nous d\'efinissons une application $P$ de l'ensemble $A^{<\NN}$  des mots finis sur l'alphabet binaire $A$ dans l'ensemble des orientations de tous les chemins finis.
 %Soit $2:= \{0,1\}$;
 Posons $\overline 0:=1$ et $\overline 1:= 0$. Soit $n\in \NN$ et  $u:=u_0 \dots u_{n-1}\in A^{n}$, le \emph{conjugu\'e} de $u$ est $\overline u:= u_{n-1} \dots u_{0}$.  Soit $P_u$ l'orientation du chemin $P_{n}$ sur $n$ sommets d\'efinie comme suit: son ensemble de sommets est $\{0, \dots, n-1\}$ et son ensemble d'arc est form\'e des paires ordonn\'ees $(i,i+1)$, $i<n-1$, telles que $u_i=1$ et des paires ordonn\'ees $(i+1, i)$, $i<n-1$,  telles que $u_{i}=0$. Notons que, $P_{n}$ \'etant ind\'ecomposable,  $P_u$ est \'egalement ind\'ecomposable; de plus, toute restriction ind\'ecomposable de $P_u$ est isomorphe \`a $P_v$ où $v$ est un facteur de $u$ et si $v\in A^{n}$, $v$ ou $\overline v$ est un facteur de $u$ si et seulement si $P_v$ s'abrite dans $P_u$. Par cons\'equent,
  si $X$ est un segment initial de $A^{<\NN}$  pour l'ordre des facteurs, son image $P(X):=\{P_u: u\in X\}$ est une classe h\'er\'editaire de l'ensemble des relations binaires ind\'ecomposables $Ind(\Omega_1)$ (ces membres \'etant consid\'er\'es \`a l'isomorphie pr\`es). De plus $X$ est un poset minimal si et seulement si $P(X)$ est une classe h\'er\'editaire minimale de $Ind(\Omega_1)$. Puisque il y a $2^{\aleph_0}$ segments initiaux minimaux de $A^{<\NN}$ (exemple, les segments initiaux  correspondant aux facteurs finis des mots de Sturm pour des valeurs distinctes de la pente, voir Proposition \ref{prop:sturm}) nous avons $2^{\aleph_0}$ classes h\'er\'editaires minimales. Ceci d\'emontre la premi\`ere partie du Th\'eor\`eme \ref{many}.

  Pour la seconde  partie,  nous transformons les mots binaires en chemins (non dirig\'es) avec des boucles. A un mot $u:=u_0 \dots u_{n-1}\in A^{n}$ nous associons  $Q_u$, le chemin $P_{n}$ sur $n$ sommets d\'efini sur $\{0, \dots, n-1\}$ avec une boucle au sommet  $i$ si et seulement si  $u_i=1$. Comme ci-dessus,  un segment initial $X$ de $A^{<\NN}$ est minimal si et seulement si $Q(X):= \{Q_u: u\in X\}$ est une classe h\'er\'editaire minimale de $Ind(\Omega_1)$.
\hfill $\Box$

\vspace{3mm}

Comme nous l'avons vu ci-dessus d'apr\`es le  Th\'eor\`eme \ref{theo:mot-structure}, les mots uniform\'ement r\'ecurrents fournissent des exemples d'\^ages ind-minimaux. En effet,  en consid\'erant des mots uniform\'ement r\'ecurrents, par exemples des mots  de sturm pour des valeurs diff\'erentes de la pente, nous avons

\begin{theorem}\label{theo:age-sturm}
 Il existe un nombre contin\^upotent d'\^ages ind-minimaux de structures binaires de type $2$.
 \end{theorem}

    \subsection{Structures ind\'ecomposables minimales}
Une autre notion de minimalit\'e est donn\'ee par la d\'efinition suivante;
\begin{definition}
Une structure relationnelle binaire infinie $\mathcal{R}$ est \emph{minimale ind\'ecomposable}\index{structure relationnelle!minimale ind\'ecomposable} si $\mathcal{R}$ s'abrite dans toute structure infinie ind\'ecomposable $\mathcal{R}'$  qui s'abrite dans $\mathcal{R}$.
\end{definition}

Plusieurs exemples, dans le cas des graphes, sont donn\'es dans le papier de Pouzet et Zaguia \cite{P-Z}.

\begin{problem}
Si $\mathcal R$ est une structures binaire minimale ind\'ecomposable est-ce que son \^age est ind-minimal?
\end{problem}

Nous n'avons pas de r\'eponse pour cette question mais elle a motiv\'e la recherche sur les \^ages des structures minimales de Pouzet-Zaguia (voir chapitre \ref{sec:age de graphe}).

\vspace{2mm}

Notons qu'un \^age ind-minimal $\mathcal A$ n'est pas n\'ecessairement l'\^age d'une structure minimale ind\'ecomposable. Pour le voir, nous allons appliquer certaines propri\'et\'es des mots aux \^ages de relations.

\vspace{2mm}

Soit $\alpha$ un mot binaire infini et $C_{\alpha}$ la structure binaire qui lui est associ\'ee, alors

\begin{lemma}\label{lem:minimal-periodic}
$C_{\alpha}$ est une structure ind\'ecomposable minimale si et seulement si $\alpha$ est p\'eriodique.
\end{lemma}

\begin{proof}
 Si $\alpha$ est p\'eriodique, alors pour tout $\beta$ infini tel que $\beta\leq\alpha$ nous avons $\alpha\leq\beta$, ce qui entraine que $C_{\alpha}$ est minimale. Inversement si $C_{\alpha}$ est minimale alors elle s'abrite dans toute structure ind\'ecomposable infinie qu'elle abrite et donc ${\alpha}$ est p\'eriodique.
\end{proof}

\vspace{2mm}

Montrons alors que parmi les \^ages du Th\'eor\`eme \ref{theo:age-sturm}, ceux qui sont associ\'es \`a des mots non p\'eriodiques ne sont pas les \^ages de structures ind\'ecomposables minimales.

\vspace{2mm}

Soit $\alpha$  un mot binaire uniform\'ement r\'ecurrent ap\'eriodique alors la structure $C_{\alpha}$ n'est pas minimale (d'apr\`es le Lemme \ref{lem:minimal-periodic}) et %D'apr\`es le Th\'eor\`eme \ref{theorem mot},
$Fac(\alpha)$ est minimal dans $\{0,1\}^{<\NN}$ muni de l'ordre des facteurs (d'apr\`es le Th\'eor\`eme \ref{theo:Facteur-minimal}). Donc l'ensemble $Ind(C_{\alpha})$ est minimal dans $Ind(\Omega_2)$ (d'apr\`es le Th\'eor\`eme \ref{theo:mot-structure}). Posons $\mathcal A_{\alpha}=\downarrow Ind(C_{\alpha})=Age(C_{\alpha})$.

\vspace{2mm}

L'ensemble $\mathcal A_{\alpha}$ est un \^age ind-minimal de $\Omega_2$. Existe-t-il une structure ind\'ecomposable minimale $\mathcal R$  telle que $Age(\mathcal R)=\mathcal A_{\alpha}$?

\vspace{2mm}

Soit $\mathcal R$ une structure de $\Omega_2$ et soit $E$ sa base. Posons $\mathcal R:=(E,c',u)$ avec $c'$ (respectivement $u$) une relation binaire (respectivement unaire) sur $E$. Nous d\'efinissons une relation d'\'equivalence sur $E$ de la mani\`ere suivante. Deux \'el\'ements $x,y$ de $E$ sont dits \'equivalents s'ils sont \'egaux ou s'il existe une suite finie d'\'el\'ements distincts de $E$,
$x=x_1,\dots,x_{n}=y$ %(ou $y=x_1,\dots,x_{n}=x$)
telle que $(x_i,x_{i+1})\in c'$ ou $(x_{i+1},x_i)\in c'$ pour tout $i=1,\dots,n-1$.  $\mathcal R$ est donc une somme directe de ses restrictions aux classes d'\'equivalence (les classes d'\'equivalence sont les composantes connexes de la sym\'etris\'ee de la relation $c'$, c'est \`a dire, $c'\cup c'^{-1}$).

\vspace{1mm}

Si $Age(\mathcal R)=\mathcal A_{\alpha}$ alors la restriction de $\mathcal R$ \`a une classe d'\'equivalence $V$ est isomorphe à une structure $\mathcal R':=(E',c'',u')$ où $E'$ est soit $\{0,1,\dots,n-1\}$ si $n=\vert V\vert$, soit $\NN$, soit $\NN^-$ l'ensemble des entiers n\'egatifs soit $\mathbb Z$, $c''$ est la cons\'ecutivit\'e sur $E'$ et $u'$ est la relation unaire d\'efinie sur $E'$ et isomorphe \`a la restriction de $u$ \`a $V$. Chacune de ces restrictions (à une classe d'\'equivalence) est donc isomorphe \`a une structure de la forme $C_{\beta}$ où $\beta\leq \alpha$.

\vspace{1mm}

Si $\mathcal R$ est ind\'ecomposable, elle est form\'ee d'une seule classe d'\'equivalence. Si son domaine $E$ est comme $\NN$ ou $\NN^-$,  elle est de la forme $C_{u}$ où $u$ est un mot uniform\'ement r\'ecurrent non p\'eriodique (car autrement, $Fac_n(u)$ serait born\'e pour tout $n$, ce qui entrainerait la m\^eme conclusion pour $Fac_n(\alpha)$ et donc $\alpha$ serait p\'eriodique, ce qui n'est pas le cas). Donc $\mathcal R$ n'est pas minimale. Si le domaine $E$ de $\mathcal R$ est comme $\mathbb Z$, alors le mot infini $x$ qui lui correspond abrite tous les mots de $Fac(\alpha)$ et tout les facteurs finis de $x$ sont dans $Fac(\alpha)$, puisque l'\^age de $\mathcal R$ est $\mathcal A_{\alpha}$. Donc la restriction de $x$ \`a $\NN$ a le m\^eme ensemble de facteur que $x$ et la restriction de $\mathcal R$ \`a $\NN$ est de m\^eme \^age que $\mathcal R$ et correspond à un  mot uniform\'ement r\'ecurrent non p\'eriodique. Donc $\mathcal R$ n'est pas minimale.

\vspace{2mm}

Comme cons\'equence nous avons:

\begin{proposition}
Soit  $\mathcal A$ l'\^age associ\'e \`a un mot uniform\'ement r\'ecurrent $u$. Il existe une structure minimale ind\'ecomposable d\'enombrable $\mathcal R$ telle que $Age(\mathcal R)=\mathcal A$ si et seulement si $u$ est p\'eriodique.
\end{proposition}

%%%%%%%%%%%%%%%%%%%%%%%%%%%%%%%%%%%%%%%%%%%%%%%%%%%%%%%%%%%%%%%%%%%%%%%%%%%%%%%%%%%%%%%

\section{Minimalit\'e, presque in\'epuisabilit\'e et belordre h\'er\'editaire}

 Nous avons montr\'e l'existence d'un  nombre contin\^upotent d'\^ages ind-minimaux.  Nous allons voir qu'une famille particuli\`ere de tels \^ages est d\'enombrable.

\vspace{2mm}

Soit $\mathcal R$ une structure relationnelle d\'efinie sur $E$. Le \emph{noyau}\index{structure relationnelle!noyau d'une -} de $\mathcal R$, not\'e $Ker(\mathcal R)$ est l'ensemble des \'el\'ements $x\in E$ tel que l'\^age de $\mathcal R_{\restriction_{E\setminus\{x\}}}$ est strictement inclus dans l'\^age de $\mathcal R$.
Un expos\'e de cette notion est dans \cite{Pou-Sobr}.

\vspace{1mm}

Un fait crucial est que si $\mathcal R$ et $\mathcal R'$ ont m\^eme \^age alors les restrictions de $\mathcal R$ et $\mathcal R'$ \`a leurs noyaux sont isomorphes. En particulier, si le noyau de $\mathcal R$ est fini alors toute structure $\mathcal R'$ de m\^eme \^age que $\mathcal R$ a un noyau fini de m\^eme cardinalit\'e que le noyau de $\mathcal R$.

\vspace{1mm}

Etant donn\'e un \^age $\mathcal A$, nous appelons \emph{noyau}\index{age@\^age!noyau d'un -} de $\mathcal A$ et notons $Ker(\mathcal A)$, le noyau de n'importe quelle relation $\mathcal R$ d'\^age $\mathcal A$.

\vspace{2mm}

Une structure relationnelle de noyau vide (respectivement fini) est dite \emph{in\'epuisable}\index{structure relationnelle!in\'epuisable} (resp. \emph{presque in\'epuisable})\index{structure relationnelle!presque in\'epuisable} et son \^age est dit \emph{in\'epuisable}\index{age@\^age!in\'epuisable} (resp. \emph{presque in\'epuisable})\index{age@\^age!presque in\'epuisable}.

\vspace{2mm}

Nous Montrons:

\begin{theorem}\label{thm:noyaufini}
Le noyau d'un \^age ind-minimal est fini.
\end{theorem}

\begin{problem}
Est-ce que le noyau d'un \^age ind-minimal a au plus deux \'el\'ements?
\end{problem}

\begin{theorem}\label{thm:indmin-denombrable}
Les  \^ages ind-minimaux de noyau non vide sont en nombre tout au plus d\'enombrable.
\end{theorem}

\begin{question}: Y en a-t-il une infinite?\end{question}

Les deux resultats ci-dessus s'appuient sur la notion de belordre h\'er\'editaire.

\vspace{2mm}

Soit $\mathscr C$ une classe h\'er\'editaire de structures binaires et soit $n\in\NN$. Notons par $\mathscr C_{+n}$ la classe  des structures binaires $\mathcal R$ telles que $\mathcal R_{-F}:=\mathcal R_{\restriction_{dom(\mathcal R)\setminus F}}\in\mathscr C$ pour au moins une partie $F$ d'au plus $n$ \'el\'ements de $dom(\mathcal R)$.

\begin{proposition}
Soit $n\in\NN$. Si une classe h\'er\'editaire $\mathscr C$ de structures binaires est h\'er\'editairement belordonn\'ee alors la classe $\mathscr C_{+n}$ %, des structures $\mathcal R$ telles que $\mathcal R_{-F}=\mathcal R_{\restriction_{dom(\mathcal R)\setminus F}}\in\mathscr C$ pour au moins une partie $F$ d'au plus $n$ \'el\'ements de $dom(\mathcal R)$,
est h\'er\'editairement belordonn\'ee.
\end{proposition}

\begin{proof}
Remarquons que $\mathscr C\subseteq\mathscr C_{+n}\subseteq\mathscr C_{+(n+1)}$ pour tout $n\in\NN$ et que $\mathscr C_{+(n+1)}=(\mathscr C_{+n})_{+1}$. %et que $\mathscr C_{+n}\subseteq\mathscr C_{+(n+1)}$.
Donc, il suffit de traiter le cas $n=1$. %puis nous it\'erons sur les autres valeurs de $n$.

Etant donn\'ee une structure $\mathcal R:=(V,(\rho_i)_{i\in I})\in \mathscr C_{+1}$, nous choisissons $x_{\mathcal R}\in V$ tel que $\mathcal R_{-{x_{\mathcal R}}}\in\mathscr C$. Soit $A$ un belordre. \`A une structure $(\mathcal R,f)\in (\mathscr C_{+1}).A$, nous associons $(\mathcal S,(f_{x_{\mathcal R}},g))$  o\`u $\mathcal S= \mathcal R_{-{x_{\mathcal R}}}$, $f_{x_{\mathcal R}}$ est la restriction de $f$ \`a $\mathcal S$ et $g$ est l'application  qui \`a $y\in dom(\mathcal S)$ associe $(\rho_i(x_{\mathcal R},y),\rho_i(y,x_{\mathcal R}))_{i\in I}$ (o\`u $I$ est l'ensemble des indices des pr\'edicats binaires).

Si on pose $B= A \times 4^I$  alors $B$ est un belordre, donc $\mathscr C.B$ est belordonn\'e.
Un plongement de $(\mathcal S,(f_{x_{\mathcal R}},g))$ dans $(\mathcal S',(f_{x_{\mathcal R}}',g'))$
induit un plongement de $(\mathcal R,f)$ dans $(\mathcal R',f')$ pourvu que la valeur de $\mathcal R$
sur $x_{\mathcal R}$ co\"{\i}ncide avec la valeur de $\mathcal R'$ sur $x_{\mathcal R'}$. Comme il n'y a qu'un
nombre fini de possibilit\'es pour ces valeurs, alors $(\mathscr C_{+1}).A$ est belordonn\'ee.
\end{proof}

\begin{corollary}\label{cor:age-wqo}
Soit $\mathcal R$ une structure binaire. Soit $x\in Ker(\mathcal R)$. Si l'\^age de $\mathcal R_{-x}$ est h\'er\'editairement belordonn\'e alors l'\^age de $\mathcal R$ est h\'er\'editairement belordonn\'e.
\end{corollary}

%\vspace{2mm}

\begin{corollary}\label{cor:indmin wqo}
Un \^age ind-minimal  de noyau non vide est h\'er\'editairement belordonn\'e.
\end{corollary}

\begin{proof}
Soit $\mathcal A$ un \^age ind-minimal de noyau non vide et soit $\mathcal R$ une structure binaire telle que $Age(\mathcal R)=\mathcal A$. Soit $x\in Ker(\mathcal R)$. Alors  $Age(\mathcal R_{-x})\varsubsetneq \mathcal A$, donc $Age(\mathcal R_{-x})$ contient un nombre fini de structures ind\'ecomposables, il s'ensuit que $Age(\mathcal R_{-x})$ est h\'er\'editairement belordonn\'e d'apr\`es la Proposition \ref{prop:wqo}. La conclusion d\'ecoule alors du Corollaire \ref{cor:age-wqo}.
\end{proof}

\vspace{2mm}

Nous avons le th\'eor\`eme suivant d\^u \`a Pouzet (1978):
\begin{theorem} \cite{pouzet.tr.1978}

Un \^age h\'er\'editairement belordonn\'e est presque in\'epuisable.
\end{theorem}

De ce dernier th\'eor\`eme avec le Corollaire \ref{cor:indmin wqo} d\'ecoule le Th\'eor\`eme \ref{thm:noyaufini}.

\vspace{2mm}

\textbf{Preuve du Th\'eor\`eme \ref{thm:indmin-denombrable}.}
D'apr\`es le Th\'eor\`eme \ref{theo:pouzet-borne}, un \^age ind-minimal de noyau non vide a un nombre fini de bornes. Comme dans une arit\'e (finie) donn\'ee il y a  tout au plus un nombre d\'enombrable d'\^ages ayant un nombre fini de bornes, le th\'eor\`eme s'ensuit. %l'ensemble de tous les sous-ensembles finis de $\NN$ est d\'enombrable
\hfill $\Box$
%%%%%%%%%%%%%%%%%%%%%%%%%%%%%%%%%%%%%%%%%%%      En plus

\begin{theorem}\label{theorem:age indmin-minimal}
Soit $\mathcal A$ un \^age ind-minimal de noyau non vide.
Toute structure binaire ind\'ecomposable $\mathcal R$  d'\^age $\mathcal A$ abrite une structure minimale.
\end{theorem}

 La preuve s'appuie sur un r\'esultat de Christian Delhomm\'e \cite{Dolho} sur le belordre de classes de relations d\'enombrables dont l'\^age ne contient qu'un nombre fini de relations ind\'ecomposables.

\vspace{2mm}

On \'etend les notions de bonne fondation et de belordre  \`a des classes  de  structures  qui ne sont pas forc\'ement finies. Disons (voir section \ref{subsec:bienfondé}) qu'une classe de structures est \emph{bien fond\'ee} si elle ne contient pas de sous-suite infinie strictement d\'ecroissante $\mathcal R_0>\dots > \mathcal R_n\dots$, la relation $<$ \'etant  l'ordre strict associ\'e au pr\'eordre d'abritement. %ainsi $\mathcal R<\mathcal S$ veut dire $R\leq S$ mais $ S\not\leq R$.

\begin{lemma} Soit $\mathscr C$ une classe h\'er\'editaire de $\Omega_{\mu}$. Si la collection  des structures $\mathcal R$ d\'enombrables telles que  $Age(\mathcal R)\subseteq \mathscr C$ est bien fond\'ee alors:

 Toute structure ind\'ecomposable  $\mathcal R$ telle que $Age(\mathcal R)\subseteq \mathscr C$ a une restriction  $\mathcal R'$ qui est minimale
\end{lemma}

\begin{proof}
 La collection   $\mathfrak R$  des structures $\mathcal R'$ d\'enombrables qui s'abritent dans $\mathcal R$  et qui sont ind\'ecomposables est  non vide (elle contient $\mathcal R$).  Cette collection est bien fond\'ee. Elle a donc un \'el\'ement minimal. N'importe lequel de ces \'el\'ements minimaux donne une restriction ind\'ecomposable minimale de $\mathcal R$.
\end{proof}
\vspace{1mm}

%Notons que toute structure ind\'ecomposable qui s'abrite strictement dans la structure minimale $\mathcal R'$ est forc\'ement finie, d'o\`u l'on d\'eduit que l'\^age de la structure minimale $\mathcal R'$ est ind-minimal.
%NB De ce fait on ne peut deduire que l'age de $R'$ est ind-minimal.
%R' ne peut abriter strictement une structure ind\'ecomposable infinie car elle est minimale, donc toute structure qu'elle abrite strictement est finie.

\vspace{2mm}

Avec ce lemme,  la preuve du Th\'eor\`eme \ref{theorem:age indmin-minimal} s'ensuit une fois prouv\'e que  si un \^age  $Age(\mathcal R)$ est ind-minimal et son noyau fini,  la  collection  des  $\mathcal R'$ d\'enombrables telles que  $Age(\mathcal R')\subseteq Age (\mathcal R)$ est bien fond\'ee pour le pr\'eordre d'abritement.

\vspace{2mm}

Il s'av\`ere qu'une telle classe est belordonn\'ee. En fait, pour obtenir le belordre de cette classe il faut prouver plus. La notion clef et celle de meilleurordre invent\'ee par C.St.J. A. Nash Williams \cite{Nash-W1,Nash-W2}. Rappelons bri\`evement qu'un ensemble ordonn\'e $\mathcal P$  est \emph{meilleurordonn\'e}\index{ensemble!meilleurordonn\'e} si la classe $\omega_1.\mathcal P$ des suites ordinales d\'enombrables \`a valeurs dans $\mathcal P$ est belordonn\'ee pour l'abritement des suites, la comparaison de deux suites ordinales se faisant ainsi: la suite $f:\alpha \rightarrow \mathcal P$ s'abrite dans la suite $g:\beta \rightarrow \mathcal P$ s'il existe une injection croissante $h: \alpha \rightarrow \beta$  telle que $f(\gamma)\leq g(h(\gamma))$ pour tout $\gamma<\alpha$.

Une classe meilleurordonn\'ee est belordonn\'ee et donc bien fond\'ee. Le r\'esultat suivant donne la conclusion d\'esir\'ee.

\begin{proposition}\label{prop:extensiondelhomme}
Si  $\mathscr C$ est ind-minimale et non in\'epuisable alors la collection des  $\mathcal R$ d\'enombrables telles que $Age(\mathcal R)\subseteq \mathscr C$ est meilleurordonn\'ee par abritement.
\end{proposition}

Ce r\'esultat est une cons\'equence d'un r\'esultat g\'en\'eral d\^u \`a C. Delhomm\'e \cite{Dolho}.

\begin{theorem}\label{thm:delhomme}
Si une classe h\'er\'editaire $\mathscr C$ de  $\Omega_k~(k\in\NN)$ ne contient qu'un nombre fini d'ind\'ecomposables alors la collection des $\mathcal R$ d\'enombrables telles que $Age(\mathcal R)\subseteq \mathscr C$ est meilleurordonn\'ee par abritement
\end{theorem}

Ce r\'esulat, \'etend le r\'esultat de Thomass\'e disant que la classe des ordres s\'erie-parall\`eles d\'enombrables est meilleurordonn\'ee \cite{thomasse}.  Le r\'esultat de Thomass\'e \'etant une extension du fameux r\'esultat de Laver \cite{laver} 1971 disant que la classe des cha\^{i}nes d\'enombrables est meilleurordonn\'ee, r\'esultat r\'epondant positivement \`a une conjecture tr\`es connue de Fra\"{\i}ss\'e (1948) affirmant que la classe des cha\^{i}nes d\'enombrables est  belordonn\'ee. La notion de meilleur ordre est l'outil qui permet de prouver le caract\`ere belordonn\'e de classes d'objets infinis. Le r\'esultat de Delhomm\'e s'appuie sur un r\'esultat de Corominas montrant le meilleurordre de certaines classes d'arbres \'etiquet\'es.

\vspace{2mm}

Nous nous contentons de prouver la Proposition \ref{prop:extensiondelhomme} \`a partir du Th\'eor\`eme \ref{thm:delhomme}.

\vspace{2mm}

\textbf{Preuve de la Proposition \ref{prop:extensiondelhomme}.}
Soit $\mathscr C$ une classe ind-minimale et non in\'epuisable. Soit  $\overline{\mathscr C}$ la collection des structures $\mathcal S$ d\'enombrables d'\^age inclus dans $\mathscr C$.  Nous devons prouver que cette collection est meilleurordonn\'ee. Comme $\mathscr C$ est ind-minimal, c'est un \^age. Soit $\mathcal R:=(V,(\rho_i)_{i\in I})$ d'\^age $\mathscr C$ et soit $K(\mathcal R)$ son noyau. Ce noyau est fini, d'apr\`es le Th\'eor\`eme \ref{thm:noyaufini}. Soit $x\in K(\mathcal R)$. On consid\`ere les relations unaires sur $V\setminus\{x\}$, $u_i^+$ et $u_i^-$ d\'efinies par $u_i^+(y):= \rho_i(x,y)$ et $u_i^-(y):= \rho_i(y, x)$. Soit $\mathcal A$ l'\^age de $\mathcal R_{-x}$ et $\mathcal A'$ l'\^age de $\mathcal R_{-x}$ augment\'ee de ces unaires.  Comme l'\^age de $\mathcal R_{-x}$ ne contient qu'un nombre fini d'ind\'ecomposable, il en va de m\^eme de $\mathcal A'$. D'apr\`es le Th\'eor\`eme \ref{thm:delhomme}, la classe des structures d\'enombrables $\mathcal R'$ telles que $Age(\mathcal R')\subseteq \mathcal A'$ est meilleurordonn\'ee. La collection $\overline{\mathscr C}$ se d\'ecompose en deux parties $\mathscr C'$ et $\mathscr C''$ form\'ees respectivement des $\mathcal S$  dont l'\^age est inclus dans l'\^age de $\mathcal R_{-x}$ et des autres.  D'apr\`es le Th\'eor\`eme \ref{thm:delhomme} la premi\`ere classe est meilleurodonn\'ee. Comme la r\'eunion de deux meilleurordres est un meilleurordre, pour obtenir le meilleurordre de $\overline {\mathscr C}$ il nous suffit d'obtenir  que la classe $\mathscr C''$ est meilleurordonn\'ee.  Si $\mathcal S\in \mathscr C''$ alors on va trouver un \'el\'ement $y\in \mathcal S$ de sorte que la transformation de $y$ en $x$ s'\'etende a toute partie finie  de $\mathcal S$ en un isomorphisme local  \`a valeurs dans $\mathcal R$. %(c'est un peu technique, admettez le).
Ainsi on peut d\'efinir sur $\mathcal S_{-y} $ des relations unaires comme ci-dessus. L'\^age de la  structure $\overline{\mathcal S}$ ainsi obtenue est inclus dans  $\mathcal A'$. Donc, d'apr\`es le r\'esultat de Delhomm\'e la collection de ces $\overline {\mathcal S}$ est meilleurordonn\'ee.  Il en r\'esulte que la collection des $\mathcal S$ est meilleurordonn\'ee. %si deux structures ont meme age, il existe un isomorphisme entre leurs noyaux qui a la propri\'et\'e que pour chaque sous-ensemble fini du noyau il s'etend à tout sous-ensemble fini. P-S
\hfill  $\Box$

%******************************************************************************

\clearemptydoublepage

\chapter{Exemples de classes minimales de graphes ind\'ecomposables}\label{sec:age de graphe}

\section{Introduction}

Dans \cite{P-Z}, Pouzet et Zaguia ont d\'ecrit les graphes minimaux  ind\'ecomposables qui ne contiennent pas une clique infinie ou un stable infini. Ces graphes sont au nombre de huit. Quatre, que nous d\'ecrivons ci-dessous, sont sans clique infinie et les quatre autres sont leurs compl\'ementaires. Nous donnons \'egalement des exemples de graphes ayant une clique infinie et un ind\'ependant infini. Parmi tous ces graphes  figurent des graphes minimaux de noyau non vide. Nous pr\'esentons ainsi un ensemble $\mathcal N$ form\'e de huit graphes minimaux de noyau non vide. Nous conjecturons que cette  liste est compl\`ete.

\vspace{1mm}

Dans ce chapitre, nous explorons les \^ages de tous ces graphes. %des graphes minimaux ind\'ecomposables de \cite{P-Z}.
Nous montrons que tous ces \^ages sont ind-minimaux et nous calculons leurs profils ainsi que leurs fonctions g\'en\'eratrices. Les croissances des profils des \^ages de ces graphes sont soit polynomiales, et dans ce cas les profils sont quasi-polynomiaux, soit exponentielles. Except\'e une, toutes les fonctions g\'en\'eratrices sont rationnelles.

\section{Description des graphes}
Tous les graphes que nous consid\'erons ici sont non dirig\'es et sans boucles. Posons  $\mathcal{G}:=\{G_0, G_1, G_2, G_3\}$ et $\overline{\mathcal G}$ l'ensemble des compl\'ementaires des graphes de $\mathcal G$. Les graphes de $\mathcal G$ sont d\'efinis comme suit.

\begin{figure}[!hbp]
\centering
\input{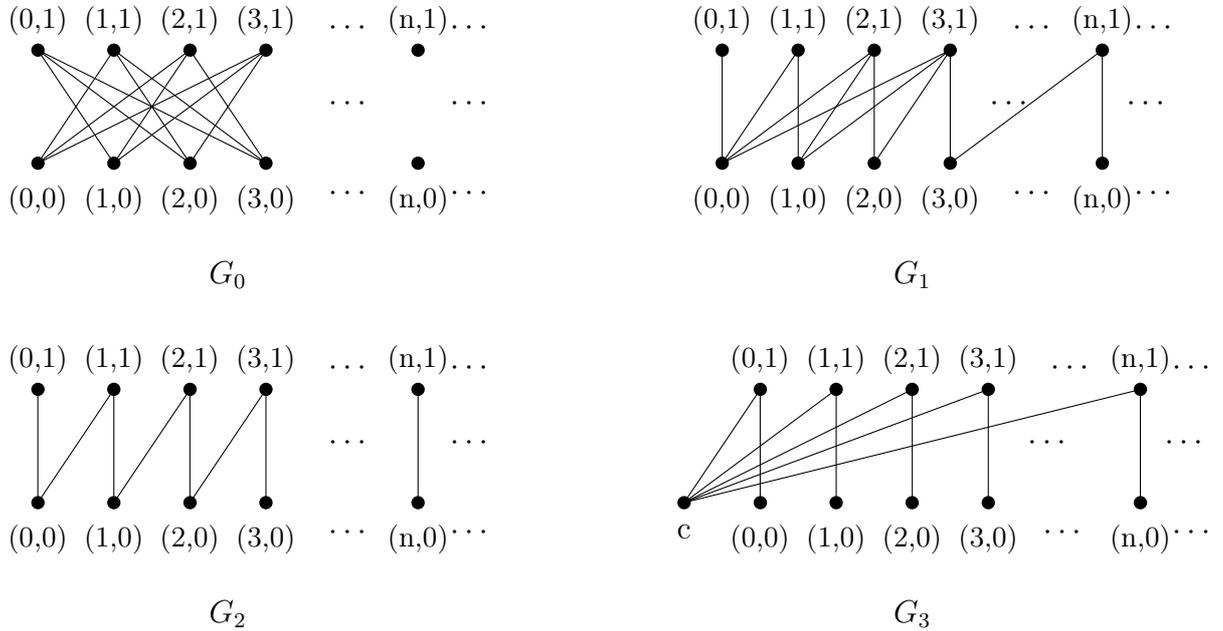}
\caption{\label{critique}Graphes ind\'ecomposables minimaux sans clique infinie.}
\end{figure}

 Tous ces graphes sont bipartis, tous sauf $G_3$ ont le m\^eme ensemble de sommets $\mathbb N \times \{0,1\}$ qui se d\'ecompose en deux ensembles ind\'ependants disjoints $A:=\mathbb N \times \{0\}$ et $B:=\mathbb N \times \{1\}$. Une paire $\{(i,0),(j,1)\}$ est une ar\^ete dans $G_0$ si $i\neq j$, une ar\^ete dans $G_1$ si $i\leq j$, une ar\^ete dans $G_2$ si $j=i$ ou $j=i+1$ et, finalement,  une ar\^ete dans $G_3$ si $j=i$. Pour $G_3$, un nouveau sommet $c$ adjacent \`a tout \'el\'ement de $B$ est rajout\'e \`a $A\cup B$.
Le graphe $G_1$ est le \emph{graphe biparti demi-complet}\index{graphe!biparti demi-complet}, le graphe critique de Schmerl et Trotter \cite{S-T}. Le graphe $G_2$ est le chemin infini, %\footnote{La cha\^{i}ne au sens des graphes.} $P_{\mathbb N}$,
tandis que le graphe $G_3$ est un arbre fait d'un nombre infini d\'enombrable d'ar\^etes disjointes reli\'ees \`a un unique sommet (not\'e $c$). $G_3$ est \'egalement appel\'e \emph{araign\'ee}.\index{araignee@araign\'ee}
Ces graphes sont repr\'esent\'es dans la \figurename~\ref{critique}.

\vspace{2mm}

\begin{figure}[t]
\centering
\input{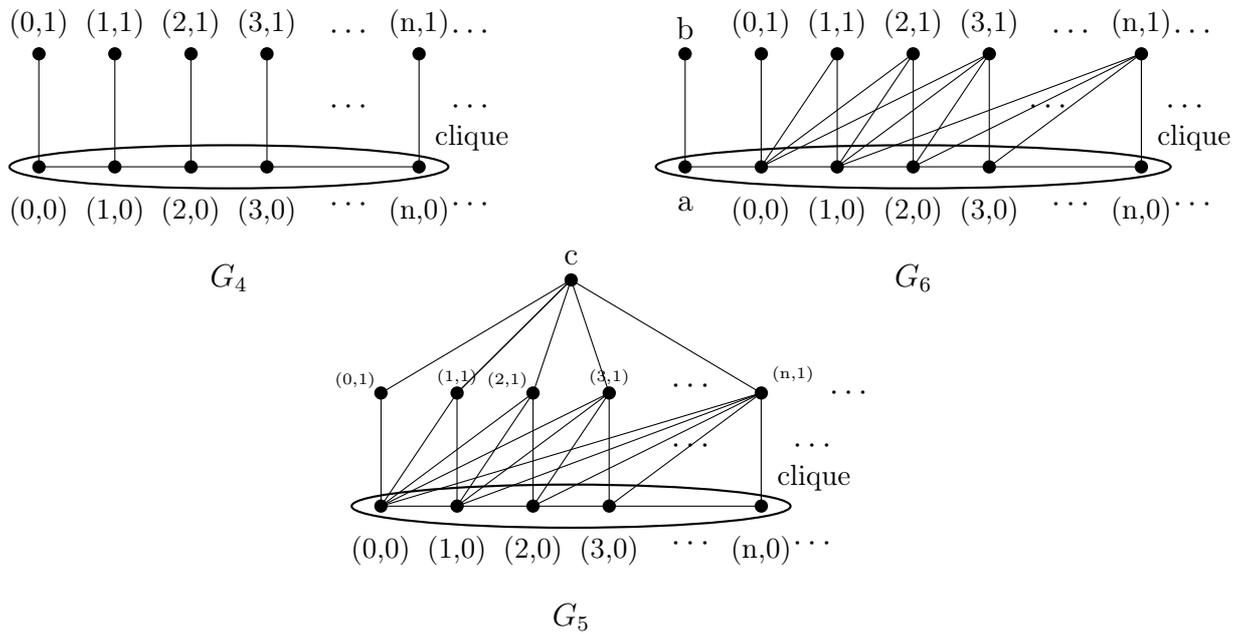}
\caption{\label{cliqueind}Les graphes $G_4$, $G_5$ et $G_6$.}
\end{figure}

\begin{figure}[!hbp]
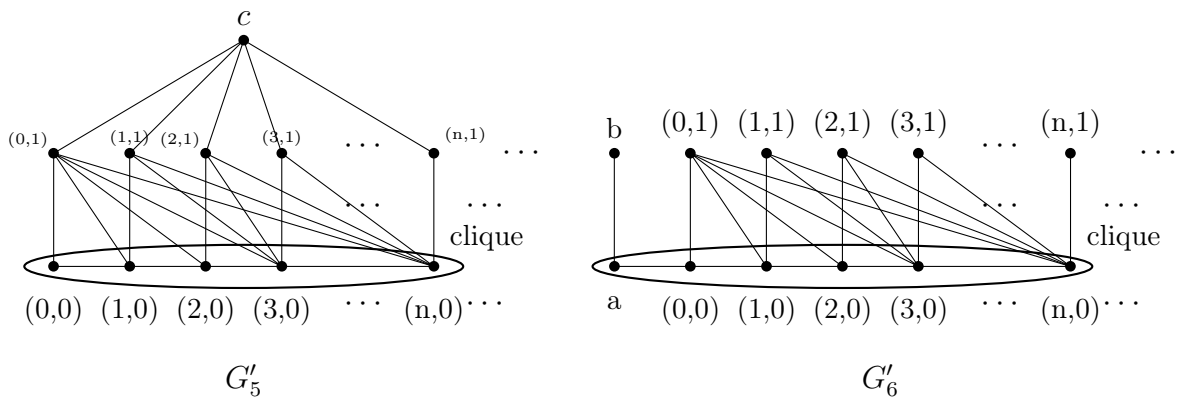

\centering
\input{image4a}\qquad\quad\input{image4b}
\caption{\label{cliqueind2}Les graphes  $G'_5$ et $G'_6$.}
\end{figure}

\noindent Pour le cas des graphes avec une clique\index{clique!infinie} infinie et un stable\index{stable!infini} infini, nous consid\'erons l'ensemble $\mathcal{G}':=\{G_4, G_5, G'_5, G_6, G'_6\}\cup\{\overline{G_4},\overline{G_6},\overline{G_6'}\}$, o\`u   $G_i:=(V_i, E_i)$, pour $i=4,5,6$ et $G_i':=(V_i',E_i')$, pour $i=5,6$, avec   $V_4:=\mathbb N \times \{0,1\}$, $V_5=V_5':=V_4\cup \{c\}$ et $V_6=V_6':=V_4\cup \{a,b\}$. Dans tous les graphes $G_i~(4\leq i\leq 6)$ et $G_i'~(5\leq i\leq 6)$, l'ensemble $V_4$ se d\'ecompose en une clique $\mathbb N \times \{0\}$ et un stable $\mathbb N \times \{1\}$.  Une paire $\{(i,0),(j,1)\}$ est une ar\^ete dans $G_4$ si $i= j$, une ar\^ete dans $G_5$ et $G_6$ si $i\leq j$ et une ar\^ete dans $G_5'$ et $G'_6$ si $i\geq j$. Pour $G_5$ et $G'_5$ nous rajoutons l'ensemble suivant d'ar\^etes $\{\{c,(i,1)\}, i\in \mathbb N\}$ et pour $G_6$ et $G'_6$ nous rajoutons toutes les ar\^etes de l'ensemble $\{\{a,(i,0)\}, \; i\in \mathbb N\}\cup \{\{a,b\}\}$.
Les graphes $G_i,~i=4,5,6$ sont repr\'esent\'es dans la \figurename~\ref{cliqueind} et les graphes $G'_i,~i=5,6$ dans la \figurename~\ref{cliqueind2}. %Ils ont les propri\'et\'es suivantes:

%Ces graphes ont les propri\'et\'es suivantes:

\vspace{2mm}

 Notons que le graphe $G_5$ (respectivement $G'_5$) s'abrite dans son graphe compl\'ementaire $\overline{G_5}$ (respectivement $\overline{G'_5}$). %et que le graphe $G'_5$ s'abrite dans son graphe compl\'ementaire $\overline{G'_5}$.
 En effet, si nous consid\'erons un abritement $f$ de la cha\^{i}ne $\omega$ des entiers non n\'egatifs $\mathbb N$ v\'erifiant $n< f(n)$ pour tout $n\in\mathbb N$ alors l'application $g$ d\'efinie par $g(n,0):=(f(n),1)$ et $g(n,1):=(f(n),0)$ est un abritement de $G_5$ dans $\overline{G_5}$. En remplaçant la chaîne $\omega$ ci-dessus par sa chaîne duale $\omega^{\star}$, l'application $g$ devient un abritement de $G_5'$ dans $\overline{G_5'}$.

\vspace{2mm}

Les graphes de $\mathcal G'$  sont ind\'ecomposables minimaux. Ces graphes sont d\'ecrits dans \cite{P-Z}.

\vspace{2mm}

Soit $\mathcal  L=\mathcal G\cup\overline{\mathcal G}\cup\mathcal G'$ la liste form\'ee de tous les graphes cit\'es. %$G_i,~(0\leq i\leq 6)$, des graphes $\overline{G_i},~(0\leq i\leq 4)$ et des graphes $G_5'$, $G_6'$, $\overline{G_6}$ et $\overline{G_6'}$.
Nous explorons, ci-dessous, les \^ages %$Age(G_i) (~0\leq i\leq 6)$
des graphes de $\mathcal L$ et calculons leurs profils ainsi que leurs fonctions g\'en\'eratrices.

Comme les \^ages d'un graphe $G$ et de son compl\'ementaire $\overline{G}$ sont isomorphes, les profils et fonctions g\'en\'eratrices de l'\^age de $G$ sont \'egalement ceux de l'\^ages de son compl\'ementaire. %de ces graphes. N'oublions pas \'egalement que les \^ages de $G_5$ et $G_5'$ sont isomorphes, c'est le cas aussi pour les \^ages de $G_6$ et $G_6'$. %$Age(\overline{G_i}) (~0\leq i\leq 6)$.
Remarquons \'egalement que les \^ages de $G_5$ et $G_5'$ sont isomorphes et qu'il en est de m\^eme pour les \^ages de $G_6$ et $G_6'$. Donc, nous \'etudierons uniquement les \^ages des graphes $G_i~(0\leq i\leq 6)$. Posons $\mathcal L'=\{G_i, ~ 0\leq i\leq 6\}$.

\begin{lemma}
Les \^ages des graphes de $\mathcal L$ sont ind-minimaux.
\end{lemma}

En effet, pour tous les graphes de $\mathcal L'$,  %$G_i,~ i=0,\cdots,6$  et leurs compl\'ementaires (sauf  $\overline{G_5}$),
les \^ages sont ind-minimaux car, comme nous le verrons ci-dessous, les classes de leurs ind\'ecomposables sont totalement ordonn\'ees pour $G_2, G_3, G_4$ et $G_6$. Pour $G_1$, les ind\'ecomposables de son \^age forment une cha\^{i}ne\index{chaine@cha\^{i}ne} pour tous les entiers pairs. Pour $G_0$ il y a un ind\'ecomposable pour chaque ordre impair  $n$, ($n\geq 5$) et deux pour tout ordre pair $n$ ($n\geq 6$) tel que tout membre d'ordre $n$ s'abrite dans tout membre d'ordre sup\'erieur \`a $n$.
Pour $G_5$ il y a deux ind\'ecomposables pour tout ordre $n\geq 5$, ils forment deux cha\^{i}nes (\`a partir de l'ordre $n\geq 5$) tel que tout membre d'ordre  $n$ s'abrite dans tout membre d'ordre sup\'erieur \`a $n$.

\vspace{2mm}

 Comme nous le verrons ci-dessous, le noyau, $Ker(G_3)$, de $G_3$ est r\'eduit a un singleton (l'\'el\'ement c), de m\^eme pour les noyaux de $G_5$ et de $G_5'$, tandis que $Ker(G_6)$ et $Ker(G_6')$  ont deux \'el\'ements chacun (les sommets $a$ et $b$).

\vspace{2mm}

Soit $\mathcal N$ l'ensemble des huit graphes de noyau non vide parmi les graphes de $\mathcal L$.
$$\mathcal N=\{G_3, \overline{G_3}, G_5, G'_5, G_6, \overline{G_6}, G_{6}', \overline{G_6'}\}.$$
%le graphe  $G_3$ et son compl\'ementaire, les graphes $G_5$ et $G'_5$, le  graphes $G_6$ et son compl\'ementaire, le graphe $G_{6}'$ et son compl\'ementaire.

Nous conjecturons que ces huit graphes sont les seuls graphes minimaux de noyau non vide. Plus pr\'ecis\'ement:

\begin{conjecture} (\cite{alzohairi-all})

Soit  $G$ un graphe infini ind\'ecomposable. Si la restriction de $G$ \`a un ensemble cofini de sommets ne contient pas de sous-graphe ind\'ecomposable infini alors $G$ abrite l'un des graphes de l'ensemble $\mathcal N$.
\end{conjecture}

Dans un graphe $G$, un ensemble de sommets $A\subseteq V(G)$ est cofini si $V(G)\setminus A$ est fini.%adjacent à un sommet de $A$.

\vspace{2mm}

Comme nous l'avons signal\'e plus haut, les noyaux des  huit graphes de $\mathcal N$ ont tous, au plus, deux \'el\'ements et leurs \^ages sont ind-minimaux.

\vspace{2mm}

Nous avons le r\'esultat suivant:

\begin{theorem}\label{pro:profiles}
Les  fonctions g\'en\'eratrices des graphes de $\mathcal L$ sont celles des graphes de $\mathcal L'$.
Pour tout $i~(0\leq i\leq 6)$, la fonction g\'en\'eratrice $F_{G_i}$ de l'\^age, $Age(G_i)$, du graphe $G_i$ et la croissance asymptotique de son profil $\varphi_i$ sont donn\'ees par:
\begin{enumerate}
\item $F_{G_0}(x)=\dfrac{1-x^2+x^3+2x^4-2x^5-x^6+x^7}{(1-x)(1-x^2)^2}~$ et $~\varphi_0(n)\simeq \dfrac{n^2}{8}$;
\item $F_{G_1}(x)=\dfrac{1-x-2x^2+x^3}{(1-2x)(1-2x^2)}~$ et $~\varphi_1(n)\simeq \dfrac{2^n}{4}$;
\item $F_{G_2}(x)=\underset{k=1}{\overset{\infty}\prod}\dfrac{1}{1-x^k}~$ et $~\varphi_2(n)\simeq \mathfrak{p}(n)$ la fonction partition d'entier;
\item $F_{G_3}(x)=\dfrac{1-x+x^3+x^4-x^6}{(1-x)^2(1-x^2)}~$ et $~\varphi_3(n)\simeq \dfrac{n^2}{4}$;
\item $F_{G_4}(x)= \dfrac{1-x+2x^3-x^5}{(1-x)^3(1+x)}~$ et $~\varphi_4(n)\simeq \dfrac{n^2}{4}$;
\item $F_{G_5}(x)=\dfrac{1-2x+x^2+3x^4+x^5+2x^6}{(1-x)(1-2x)}~$ et $~\varphi_5(n)\simeq 2^n$;
\item $F_{G_6}(x)=\dfrac{1-3x+3x^2-x^3+x^4}{(1-2x)(1-x)^2}~$ et $~\varphi_6(n)\simeq 3.2^{n-2}$.
\end{enumerate}
Ces fonctions g\'en\'eratrices sont toutes distinctes. A part $F_{G_2}$, elles sont toutes rationnelles. En outre, la croissance du profil est polynomiale de puissance $2$ pour les graphes $G_0$, $G_3$ et $G_4$ et exponentielle pour les graphes $G_1$, $G_5$ et $G_6$. %celle de $G_2$ est sous-exponentielle.
\end{theorem}

Notons que les profils des \^ages $Age(G_i)~(i=0,3,4)$  sont quasi-polynomiaux, ceux des \^ages $Age(G_i)~(i=1,5, 6)$ sont exponentiels et le profil de $Age(G_2)$ est sous-exponentiel.

\vspace{2mm}

Les assertions de la Proposition \ref{pro:profiles} seront d\'emontr\'ees s\'epar\'ement, chacune dans une des sections suivantes de ce chapitre.

%Pour la deuxi\`eme question, nous montrons que  les profils des \^ages $Age(G_i)~(i=0,3,4)$ sont polyn\^omiaux et que la fonction g\'en\'eratrice $F_{G_i}~(0\leq i\leq 6,~i\neq 2)$ est rationnelle.
%$$$$$$$$$$$$$$$$$$$$$$$$$$$$$$$$$$$$$$$$$$$$$$$$$$$$$$$$$$$$$$$$$$$$$$$
%$$$$$$$$$$$$$$$$$$$$$$$$$$$$$$$$$$$$$$$$$$$$$$$$$$$$$$$$$$$$$$$$$$$$$$$$

\section{Le graphe $G_0$}\label{subsec:graphe $G_0$}

Le graphe $G_0:=(V_0,E_0)$ est biparti, son ensemble de sommets $V_0:=\mathbb N \times \{0,1\}$ se d\'ecompose  en deux stables disjoints $A:=\mathbb N \times \{0\}$ et $B:=\mathbb N \times \{1\}$. Une paire $\{(i,0),(j,1)\}$ est une ar\^ete de  $G_0$ si $i\neq j$.

\medskip
\subsection{Repr\'esentations des sous-graphes de $G_0$}

D\'esignons par \emph{paire} dans $G_0$ tout sous-ensemble $\{(i,0),(i,1)\}$, pour $i\in \mathbb N$. Ces deux sommets ne sont pas li\'es dans $G_0$.

Soit $H=(V,E)$ un sous-graphe de $G_0$, posons $V_A=V\cap A$ et $V_B=V\cap B$. Nous pouvons %Observons que tout sous-graphe $S$ de $G_0$ peut-\^etre
repr\'esenter $H$ de la mani\`ere suivante:

Consid\'erons deux lignes parall\`eles horizontales $L_0, L_1$. Pla\c{c}ons les sommets de $V_A$ sur une des lignes et ceux de $V_B$ sur l'autre de telle sorte que toutes les paires de $G_0$ qui appartiennent \`a $H$ soient \`a gauche, les sommets d'une m\^eme paire \'etant align\'es verticalement. Les sommets restants sont plac\'es \`a droite de ces paires avec ceux sur $L_1$ compl\`etement \`a droite de ceux sur $L_0$ (voir  \figurename~\ref{representationG0}).

Observons que si nous permutons les places de $V_A$ et $V_B$ sur les lignes $L_0$ et $L_1$, le graphe obtenu sera isomorphe \`a $H$. Nous pouvons donc choisir de placer sur, disons $L_0$, l'ensemble ($V_A$ ou $V_B$) ayant le nombre minimum de sommets.

Avec cette repr\'esentation, nous pouvons associer, comme le montre la  \figurename~ \ref{representationG0}, \`a chaque membre $H$ de $Age(G_0)$, d'ordre $n$, un triplet d'entiers $(n_0,n_1,p)$, o\`u $p$ est le  nombre de paires et tel que:
\begin{equation}\label{eqG0}
 p\leq n_0\leq n_1~\text{et } n_0+n_1=n.
 \end{equation}

 \begin{figure}[h]
\centering
\input{image5}
\caption{\label{representationG0}Repr\'esentation des sous-graphes de $G_0$.}
\end{figure}

Inversement, \`a tout triplet $t:=(n_0,n_1,p)$ v\'erifiant les conditions \eqref{eqG0} ci-dessus, avec $n_0+n_1\geq 1$,  nous associons le type d'isomorphie $H_t$ de ${G_0}_{\restriction_X}$ o\`u $X=\{0,\cdots,n_0-1\}\times \{0\}\cup (\{0,\cdots,p-1\}\cup\{n_0,\cdots,n_0+n_1-p-1\})\times \{1\}.$ Il est \'evident que si $t=t'$  alors $H_t=H_{t'}$.

\begin{remark}\label{rem:G0}

\begin{itemize}
\item Pour $n=1$, il existe un seul sous-graphe, il est repr\'esent\'e de mani\`ere unique par le triplet $(0,1,0)$.
 \item Pour $n=2$ nous avons deux sous-graphes, la clique \`a deux sommets, dont l'unique repr\'esentation est le triplet $(1,1,0)$ et l'ind\'ependant \`a deux sommets qui peut-\^etre repr\'esent\'e par deux triplets diff\'erents $(0,2,0)$ et $(1,1,1)$.
\end{itemize}
\end{remark}

Soit $F=\{t=(n_0,n_1,p)\in\mathbb N^3/\; t \text{ v\'erifie \eqref{eqG0} et }p=n_0\leq 1\}\cup\{(2,2,2)\}$.

\begin{remark}\label{rem:F-G0}
Si $t\in F$ alors,
\begin{itemize}
\item ou bien $t=(0,n_1,0)$ et dans ce cas $H_t$ est un ind\'ependant d'ordre $n_1$.
\item ou bien $t=(1,n_1,1)$ et dans ce cas $H_t$ est la somme directe d'un sommet et d'une \'etoile\index{etoile@\'etoile} \`a $n_1-1$ extr\'emit\'es  $S_{1,n_1-1}$ (une \'etoile \`a $k$ extr\'emit\'e, not\'ee $S_{1,k}$, est un graphe \`a $k+1$ sommets $\{x_0,x_1,\cdots,x_k\}$ ayant pour ar\^etes l'ensemble $\{\{x_0,x_i\},\;i=1,\cdots,k\}$).
\item ou bien $t=(2,2,2)$ et dans ce cas $H_t$ est le graphe $2.K_2$, (voir \tablename~\ref{grapheF}).
\end{itemize}
\end{remark}
\begin{table}[h]
\hspace{2cm}\begin{tabular}[c]{|c|c|c|}
\hline
\input{imag4} &\input{imag5}& \input{imag6}\\
\hline
$\begin{array}{c}
H_t=A_n\\
t=(0,n,0)
\end{array}$& $\begin{array}{c}
H_t=1\oplus S_{1,k}\\
t=(1,k+1,1)
\end{array}$& $\begin{array}{c}
H_t=2K_2\\
t=(2,2,2)
\end{array}$\\
\hline
\end{tabular}
\caption{\label{grapheF}Les sous-graphes $H_t$ pour $t\in F$.}
\end{table}

\begin{lemma}\label{lem:G0}
Soit $t=(n_0,n_1,p)$ un triplet d'entiers v\'erifiant les conditions \eqref{eqG0} tel que $n_0+n_1\geq 3$ et soit $H_t$ le type d'isomorphie associ\'e, alors
 $$t\notin F \text{si et seulement si } H_t \text{ est connexe}.$$
\end{lemma}

\begin{proof}\\
$\Leftarrow)$ Vient de la Remarque \ref{rem:F-G0}.\\
$\Rightarrow)$ Si $t\notin F$ alors ou bien $p\neq n_0$, ou bien $p=n_0\geq 2$ et $n_1\neq 2$.\\
 Dans le premier cas, $p\neq n_0$, nous avons les situations suivantes:

    \qquad \textbf{a)} Si $p=0$ alors $H_t$ est isomorphe au graphe biparti complet\index{graphe!biparti complet} $K_{n_0,n_1}$ qui est connexe\index{graphe!connexe}.

     \qquad \textbf{b)} Si $p\neq 0$ alors $H_t={G_0}_{\restriction_{X_t}}$ avec $X_t=\{0,\cdots,n_0-1\}\times \{0\}\cup (\{0,\cdots,p-1\}\cup\{n_0,\cdots,n_0+n_1-p-1\})\times \{1\}.$ Comme tout sommet de $\{0,\cdots,n_0-1\}\times \{0\}$ est reli\'e \`a tout sommet de $\{n_0,\cdots,n_0+n_1-p-1\}\times \{1\}$ et tout sommet de $\{0,\cdots,p-1\}\times \{1\}$ est reli\'e \`a tout sommet de $\{p,\cdots,n_0-1\}\times \{0\}$, le graphe $H_t$ est connexe. \\%(voir represenation***) \\
Dans le deuxi\`eme cas, $p=n_0\geq 2$ et $n_1\neq 2$, nous avons \\%(voir figure**):

    \qquad \textbf{a)} Si $p=2$ alors $n_1>2$ et donc $H_t$ est connexe.

    \qquad \textbf{b)} Si $p\geq 3$ alors $H_t$ connexe quelque soit la valeur de $n_1\geq n_0$.
\end{proof}

\begin{lemma}\label{lem:abrit-G0}
Soient $t=(n_0,n_1,p)$ et $t'=(n'_0,n'_1,p')$ deux triplets d'entiers v\'erifiant les condition \eqref{eqG0} et tels que $n_0+n_1=n\geq 3$ et $n'_0+n'_1=n'\geq 3$, $H_t$ et $H_{t'}$ les types d'isomorphie des sous-graphes\index{graphe!sous-graphe} de $G_0$ associ\'es respectivement \`a $t$ et $t'$, alors on a
$$H_t\leq H_{t'}\Leftrightarrow p\leq p', n_0\leq n'_0 \text{ et } n_1\leq n'_1.$$
\end{lemma}

\begin{proof}
D'apr\`es la repr\'esentation des sous-graphes de $G_0$ donn\'ee ci-dessus, $H_t$ est le type d'isomorphie de ${G_0}_{\restriction_X}$ o\`u $X=\{0,\cdots,n_0-1\}\times \{0\}\cup (\{0,\cdots,p-1\}\cup\{n_0,\cdots,n_0+n_1-p-1\})\times \{1\}$ et $H_{t'}$ est le type d'isomorphie de ${G_0}_{\restriction_{X'}}$ o\`u $X'=\{0,\cdots,n'_0-1\}\times \{0\}\cup (\{0,\cdots,p'-1\}\cup\{n'_0,\cdots,n'_0+n'_1-p'-1\})\times \{1\}$. Il est \'evident que si $p\leq p'$, $n_0\leq n'_0$  et  $n_1\leq n'_1$ on a $H_t\leq H_{t'}$. Inversement, si $H_t\leq H_{t'}$, alors il existe un abritement de $H_t$ dans $H_{t'}$ qui envoie les paires de $H_t$ sur les paires de $H_{t'}$, donc $p\leq p'$. D'un autre c\^ot\'e, si $H_t$ et $H_{t'}$ sont connexes alors, \'etant des graphes bipartis\index{graphe!biparti}, les partitions % \footnote{C'est un r\'esultat classique de la th\'eorie des graphes, si un graphe biparti est connexe alors sa partition en stables est unique.}
de leurs sommets en deux stables sont uniques, se sont donc celles des ensembles $X$ et $X'$. Donc $n_0\leq n'_0$ et $n_1\leq n'_1$. Si $H_t$ et $H_{t'}$ ne sont pas connexes, alors d'apr\`es le Lemme \ref{lem:G0}, $t,t'\in F$ et d'apr\`es la Remarque \ref{rem:F-G0} et le fait que les sous-graphes sont d'ordre au moins trois, on a n\'ecessairement $n_0\leq n'_0$ et $n_1\leq n'_1$.
\end{proof}

\begin{lemma}
Si deux sous-graphes de $G_0$, d'ordre au moins trois, sont isomorphes alors ils sont repr\'esent\'es par un m\^eme triplet.
%A deux sous-graphes de $G_0$, d'ordre au moins trois, qui sont isomorphes correspond une m\^eme repr\'esentation et donc un m\^eme triplet.
\end{lemma}

\begin{proof}
Soient $H$ et $H'$ deux sous-graphes isomorphes, d'ordre $n\geq 3$, de $G_0$ et soient $t=(n_0,n_1,p)$ et $t'=(n'_0,n'_1,p')$ les triplets associ\'es respectivement. Montrons que $t=t'$.  $H$ et $H'$ sont les restrictions de $G_0$ aux sous-ensembles $V=A_1\times \{0\}\cup A_2\times \{1\}$ et $V'=A'_1\times \{0\}\cup A'_2\times \{1\}$ respectivement o\`u $A_1, A_2, A'_1, A'_2$ sont des sous-ensembles de $\mathbb N$ et $\vert A_1\cup A_2\vert= \vert A'_1\cup A'_2\vert=n.$ Comme $H$ et $H'$ sont isomorphes, ils ont le m\^eme nombre de paires, donc $p=p'$.\\
 S'ils sont connexes alors, \'etant des graphes bipartis, les partitions de leurs sommets en deux stables sont uniques. Donc, soit $\vert A_1\vert=\vert A'_1\vert$ et $\vert A_2\vert=\vert A'_2\vert$ ou bien $\vert A_1\vert=\vert A'_2\vert$ et $\vert A_2\vert=\vert A'_1\vert$. Avec notre convention de placer, en $L_0$, le stable de cardinalit\'e minimale, nous avons n\'ecessairement $n_0=n'_0$ et $n_1=n'_1$.\\
 Si $H$ et $H'$ ne sont pas connexes, alors d'apr\`es le Lemme \ref{lem:G0}, $t,t'\in F$ et d'apr\`es la Remarque \ref{rem:F-G0}, on a n\'ecessairement $t=t'$.
\end{proof}

\subsection{Les ind\'ecomposables de $Age(G_0)$}

A partir de la repr\'esentation des sous-graphes de $G_0$ et de leurs codages par des triplets d'entiers, nous pouvons d\'eterminer l'ensemble $Ind(G_0)$ des sous-graphes ind\'ecomposables de $G_0$. \\Soit $H$ un sous-graphe d'ordre $n\geq 4$ de $G_0$ induit par le sous-ensemble $V_H$ de sommets et soit $t_H=(n_0,n_1,p)$ le triplet qui lui est associ\'e. \\Posons $P=\{x\in V_H/ x \text{\ appartient \`a une paire de } H\}$, on a $\vert P\vert=2p$. Pour tout $x\in P$, notons $\overline x$ le sommet de $P$ tel que $\{x,\overline x\}$ est une paire.

\begin{lemma}\label{lem:int-G0}
Si $n_0-p\geq 2$ ou $n_1-p\geq 2$ alors $H$ poss\`ede un intervalle non trivial. La r\'eciproque n'est pas vraie.
\end{lemma}

\begin{proof}
Si $n_0-p\geq 2$, alors dans la partition en deux stables des sommets de $V_H$, en consid\'erant le stable de cardinalit\'e $n_0$, les sommets de ce stable qui n'appartiennent pas \`a $P$ forment un intervalle de $H$. Pour le cas $n_1-p\geq 2$, consid\'erer les sommets du stable de cardinalit\'e $n_1$ qui ne sont pas dans $P$. Pour la r\'eciproque, il suffit de consid\'erer le sous-graphe $2.K_2$ qui poss\`ede un intervalle et qui correspond au triplet $(2,2,2)$.
\end{proof}

\begin{lemma}\label{lem:int-P}
Soit $E\subseteq P$, avec $E\neq P$ et $\vert E\vert\geq 2$ , alors;
\begin{enumerate}
\item Si $p\geq 3$, l'ensemble $E$ n'est pas un intervalle.
\item Si $p=2$, $H$ poss\`ede un intervalle si et seulement si $n_0=n_1= p$.
\end{enumerate}
\end{lemma}

\begin{proof}
\begin{enumerate}
\item Si $p\geq 3$, alors $E$ n'est pas un intervalle, il suffit de consid\'erer un sommet $x$ de $P\setminus E$ tel que $\overline x\in E$, s'il en existe, sinon, $E$ est form\'e de paires, prendre alors n'importe quel \'el\'ement $x\in P\setminus E$.
\item Si $p=2$, alors,
    \begin{itemize}
    \item  $n_0= p$ et $n_1= p$ impliquent que $H$ est isomorphe \`a $2K_2$ qui poss\`ede un intervalle.
    \item Si $E$ est un intervalle alors n\'ecessairement $E=\{x,y\}$ tels que $(x,y)$ est une ar\^ete de $H$, donc on a forc\'ement $n_0= p$ et $n_1= p$.
    \end{itemize}
\end{enumerate}\end{proof}

\begin{lemma}\label{lem:ind-G0}
Soit $H$ un sous-graphe de $G_0$ d'ordre $n\geq 4$ et $t_H=(n_0,n_1,p)$ le triplet qui lui est associ\'e. Alors
$$H\in Ind(G_0) \text{ si et seulement si }t_H\in\{(k,k,k-1),(k,k+1,k); k\geq 2\}\cup\{(k,k,k); k\geq 3\}.$$
\end{lemma}

\begin{proof}
Si $H$ est ind\'ecomposable\index{graphe!ind\'ecomposable}, alors d'apr\`es le Lemme \ref{lem:int-G0}, $n_0-p\leq 1$ et $n_1-p\leq 1$, donc soit $n_0=n_1=p$ avec $p\geq 3$, d'apr\`es le Lemme \ref{lem:int-P}, soit $n_0=p$ et $n_1=p+1$ soit $n_0=n_1=p+1$ avec $n_0+n_1\geq 4$ et on retrouve l'ensemble donn\'e.\\ Inversement, si $t_H\in\{(k,k,k-1),(k,k+1,k); k\geq 2\}\cup\{(k,k,k); k\geq 3\}$, il est facile de voir que %alors d'apr\`es les lemmes \ref{lem:int-G0} et \ref{lem:int-P},
$H$ est ind\'ecomposable.
\end{proof}

\bigskip
D'apr\`es le Lemme \ref{lem:ind-G0}, le graphe $G_0$ a deux sous-graphes ind\'ecomposables de tout ordre pair $n:=2k, \; k\geq 3$ qui correspondent aux triplets $(k,k,k-1)$ et $(k,k,k)$ et un de tout ordre impair $n:=2k+1,\; k\geq 2$ qui est cod\'e par $(k,k+1,k)$. ces ind\'ecomposables sont repr\'esent\'es dans la \figurename~\ref{indecG0}. Pour $n=4$ il y a un seul ind\'ecomposable qui est isomorphe \`a $P_4$ et est cod\'e par $(2,2,1)$, il n'y a pas d'ind\'ecomposable de taille $3$,  Chaque ind\'ecomposable de taille $n$ s'abrite dans tout ind\'ecomposable de taille sup\'erieure, d'apr\`es le Lemme \ref{lem:abrit-G0}. Donc l'ensemble des ind\'ecomposables de $G_0$ forme une classe minimales de $Ind(\Omega_1)$. Donc $Age(G_0)$, l'\^age de $G_0$ (et aussi $Age(\overline {G_0})$, l'\^age du compl\'ementaire de $G_0$) est ind-minimal d'apr\`es  le Th\'eor\`eme \ref{theo:equiv-min}.

\begin{figure}[h]
\centering
\input{image6}
\caption{\label{indecG0}Les sous-graphes ind\'ecomposables d'ordres $n\geq 5$ de $G_0$.}
\end{figure}

\subsection{Profil de l'\^age de $G_0$}\label{subsec:profilG0}

Les premi\`eres valeurs du profil de $Age(G_0)$ pour $n=0,1,2,3,4,5,6,7$ sont respectivement $1,1$,$2,3$,$6,6$, $10,10$.

\begin{lemma}
La fonction profil de $G_0$ est quasi-polynomiale, elle est donn\'ee, pour tout $n\in \mathbb N$,  par:
%$$\varphi(2)=2,~~\varphi(n)=\underset{k=0}{\overset{\lfloor n/2\rfloor}{\sum}}(k+1),~n\neq 2.$$
%ou encore, en posant $l=\lfloor n/2\rfloor$;
$$\left\{\begin{array}{ll}\varphi_0(2)=2,&\\
\varphi_0(n)=\dfrac{n^2+6n+8}{8}&\text{si }n \text{ est pair }, ~n\neq 2,\\
\varphi_0(n)=\dfrac{n^2+4n+3}{8}&\text{si }n \text{ est impair }.
\end{array}\right.$$
En outre, la croissance de $\varphi_0(n)$ est polynomiale de puissance $2$,  $\varphi_0(n)\simeq\dfrac{n^2}{8}$.
\end{lemma}

\begin{proof}
A partir du codage des sous-graphes par des triplets et de l'\'enum\'eration de ces derniers, nous avons $$\left\{\begin{array}{l}\varphi_0(2)=2,\\
\varphi_0(n)=\underset{k=0}{\overset{\lfloor n/2\rfloor}{\sum}}(k+1)=\dfrac{({\lfloor n/2\rfloor}+1)(\lfloor n/2\rfloor+2)}{2},~\forall n\neq 2.\end{array}\right.$$
Apr\`es calcul, nous obtenons la fonction donn\'ee dans le lemme.

 Le cas $n=2$ n'ob\'eit pas \`a cette formule car, d'apr\`es la Remarque \ref{rem:G0}, l'ind\'ependant \`a deux \'el\'ements peut-\^etre cod\'e par deux triplets diff\'erents.

 \vspace{1mm}

 Il est clair que le profil $\varphi_0(n)$ est quasi-polynomial. En effet, pour tout $n\neq 2$, nous avons $\varphi_0(n)=a_2(n)n^2+a_1(n)n+a_0(n)$ o\`u $a_2(n)=\dfrac{1}{8}$ pour tout $n$, $a_0$ et $a_1$ sont p\'eriodiques de p\'eriode $2$ et valent: $a_0(2n)=1$,  $a_0(2n+1)=\dfrac{3}{8}$, $a_1(2n)=\dfrac{3}{4}$ et $a_1(2n+1)=\dfrac{1}{2}$ pour tout $n$.
\end{proof}

\begin{proposition}
La fonction g\'en\'eratrice de $Age(G_0)$ (et de $Age(\overline{G_0})$) est rationnelle\index{fonction!g\'en\'eratrice!rationnelle} et est donn\'ee par:
$$F_{G_0}(x)=\dfrac{1-x^2+x^3+2x^4-2x^5-x^6+x^7}{(1-x)(1-x^2)^2}.$$
\end{proposition}

\begin{proof}
 D'apr\`es le codage des \'el\'ements de $Age(G_0)$ par des triplets, la s\'erie g\'en\'eratrice peut-\^etre calcul\'ee comme suit:
$$F_{G_0}=\underset{
\begin{array}{c}
 n_0+n_1=n\\
p\leq n_0\leq n_1%
\end{array}%
}{\sum} \delta_{(n_0,n_1,p)}X^n$$
o\`u $\delta_{(n_0,n_1,p)}$ est la fonction caract\'eristique de l'ensemble $\{(n_0,n_1,p);n_0+n_1=n\;\text{et } \; p\leq n_0\leq n_1\}.$\\
En prenant toutes les partitions de $n$ en deux parts $n_0$ et $n_1$ avec $n_0\leq n_1$, nous obtenons, en posant $n_1=n_0+m$.\\

$\begin{array}{cl}
F_{G_0}&=\underset{
\begin{array}{c}
 2n_0+m=n\\
p\leq n_0%
\end{array}%
}{\sum} \delta_{(n_0,n_1,p)}X^{2n_0+m}=\underset{2n_0+m=n}{\sum} (n_0+1)X^{2n_0}.X^{m}\\
&=\underset{n_0\geq 0}{\sum} (n_0+1).X^{2n_0}.\underset{m\geq 0}{\sum} X^{m}%
\end{array}%
$

Qui, apr\`es calculs et corrections (bas\'ees sur le calcul des premi\`eres valeurs du profil), donne:
$$F_{G_0}(X)=\dfrac{1+X}{(1-X^2)^3}-X^2=\dfrac{1-X^2+X^3+2X^4-2X^5-X^6+X^7}{(1-X)(1-x^2)^2}.$$
\end{proof}
%

%***********************************************************************************************
%$$$$$$$$$$$$$$$$$$$$$$$$$$$$$$$$$$$$$$$$$$$$$$$$$$$$$$$$$$$$$$$$$$$$$$$

\section{Le graphe $G_1$}\label{subsec:graphe $G_1$}

Le graphe $G_1:=(V_1,E_1)$ est biparti, son ensemble de sommets $V_1:=\mathbb N \times \{0,1\}$ se d\'ecompose  en deux ind\'ependants disjoints $A:=\mathbb N \times \{0\}$ et $B:=\mathbb N \times \{1\}$. Une paire de sommets $\{(i,0),(j,1)\}$ est une ar\^ete de  $G_1$ si et seulement si $i\leq j$. Le graphe $G_1$ est le graphe biparti demi-complet de Schmerl et Trotter\index{Schmerl et Trotter} \cite{S-T}, c'est un graphe critique\index{graphe!critique}. Le graphe $G_1$ et son graphe compl\'ementaire sont des graphes de permutations\index{graphe!de permutations} (ils sont associ\'es \`a des ordres de dimension deux,\index{ordre!de dimension deux}  voir section \ref{sec:conjecture}).

\subsection{Repr\'esentation des sous-graphes de $G_1$}
Soit $H$ un sous-graphe d'ordre $n$ de $G_1$ induit par le sous-ensemble de sommets $V_H=X\cup Y$ avec $X\subseteq \mathbb N\times\{0\}$,  $Y\subseteq \mathbb N\times\{1\}$ et $\vert X\cup Y\vert=n$.

\begin{remarks}\label{rem:G1}
Nous avons, d'apr\`es la structure du graphe $G_1$,
\begin{enumerate}
\item $H$ est connexe si et seulement si
    \begin{itemize}
    \item $X\neq\varnothing$ et $Y\neq\varnothing$,
    \item $min\{i\in\mathbb N/(i,0)\in X\}\leq min\{j\in\mathbb N/ (j,1)\in Y\}$,
    \item $max\{i\in\mathbb N/(i,0)\in X\}\leq max\{j\in\mathbb N/ (j,1)\in Y\}$.
    \end{itemize}
\item Si $H$ n'est pas connexe, il poss\`ede, au plus, une composante connexe\index{composante!connexe} d'ordre sup\'erieur ou \'egal \`a deux.
\item Pour deux sommets quelconques $(i,0),~(i',0)$ de $X$ on a $$d_H((i,0))\geq d_H((i',0))\Leftrightarrow i\leq i'.$$
\item Pour deux sommets quelconques $(j,1),~(j',1)$ de $Y$ on a $$d_H((j,1))\leq d_H((j',1))\Leftrightarrow j\leq j',$$
o\`u $d_H(x)$ est le degr\'e du sommet $x$ dans $H$.
\end{enumerate}
\end{remarks}
\bigskip

\begin{figure}[h]
\centering
\input{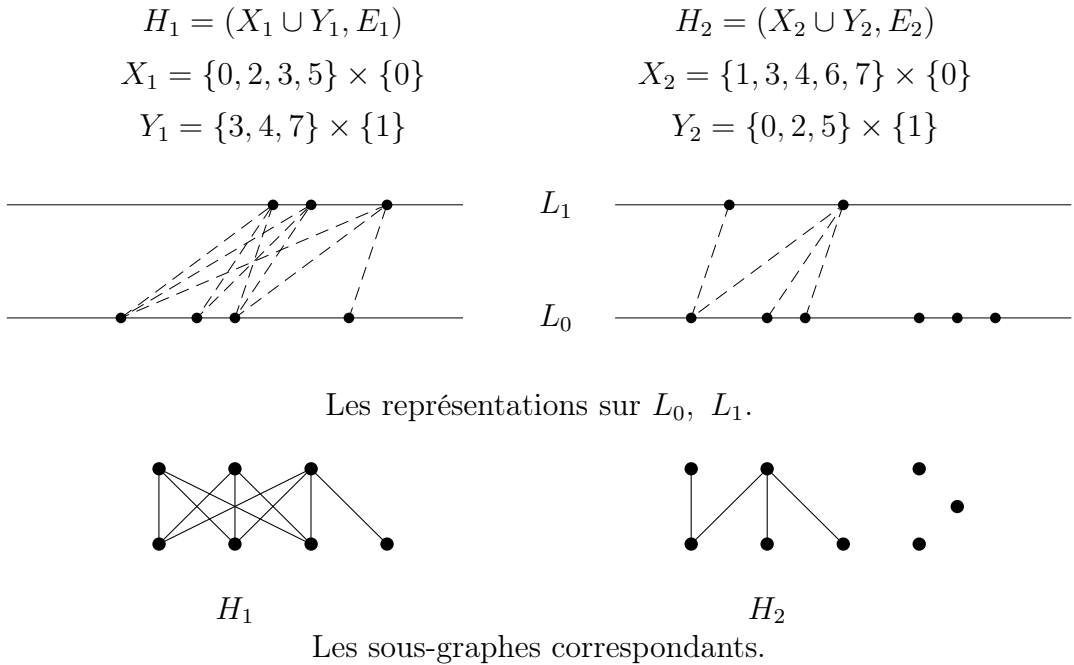}
\caption{\label{representationG1}Repr\'esentation des sous-graphes de $G_1$.}
\end{figure}

Nous associons au sous-graphe $H$ la repr\'esentation suivante.
Consid\'erons deux droites parall\`eles horizontales $L_0,~L_1$, alors
\begin{enumerate}
\item Si $H$ est connexe et $\vert X\vert\geq\vert Y\vert$ (respectivement $\vert Y\vert\geq\vert X\vert$), placer les sommets $(i,0)$ de $X$ (respectivement $(j,1)$ de $Y$) sur $L_0$ suivant l'ordre croissant des $i$ (respectivement d\'ecroissant des $j$) et les sommets $(j,1)$ de $Y$ (respectivement $(i,0)$ de $X$) sur $L_1$ suivant l'ordre croissant des $j$ (respectivement d\'ecroissant des $i$) de telle sorte que  tout $(j,1)\in Y$ se trouve à droite (respectivemnt à gauche) de tout $(i,0)\in X$  tel que $j\geq i$ et à gauche (respectivement à droite) de tout $(i,0)$  tel que $j<i$ %les sommets de $ N_H(x)\subseteq Y$ (respectivement\footnote{$N_H(x)$ d\'esigne l'ensemble des voisins de $x$ dans $H$.} $N_H(y)\subseteq X$), se trouvent \`a sa droite et ceux de $Y\setminus N_H(x)$ (respectivement $X\setminus N_H(y)$) se trouvent \`a sa gauche.
    (voir l'exemple $H_1$ sur la \figurename~\ref{representationG1})
\item Si $H$ n'est pas connexe alors, d'apr\`es la Remarque \ref{rem:G1}, nous avons
    \begin{itemize}
    \item ou bien $H$ est un ind\'ependant, dans ce cas %l'un des deux ensembles $X$ ou $Y$ est vide,
    nous pla\c{c}ons tous les sommets sur $L_0$.
    \item ou bien $H$ poss\`ede une composante connexe $H'$ d'ordre au moins deux, nous pla\c{c}ons les sommets de $H'$ comme dans le cas $H$ connexe, puis nous pla\c{c}ons les sommets restants (qui n'appartiennent pas \`a la composante connexe $H'$) sur $L_0$ \`a droite de tous les sommets de $H'$, (voir l'exemple du sous-graphe $H_2$ sur la \figurename~\ref{representationG1}).
    \end{itemize}
\end{enumerate}

\begin{remark}\label{rem:uniqueG1}
Si $\vert X\vert=\vert Y\vert$ la repr\'esentation ci-dessus n'est pas n\'ecessairement unique.
\end{remark}

\begin{figure}[t]
\centering
\input{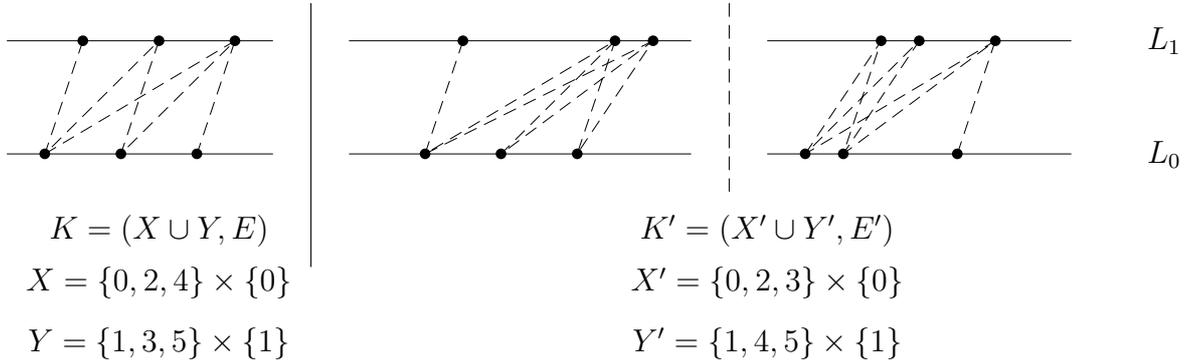}
\caption{\label{exempleG1}Exemples de repr\'esentation des sous-graphes de $G_1$.}
\end{figure}

\begin{example}
\begin{enumerate}
\item Si $X=\{(0,0),(2,0),(4,0)\}$ et $Y=\{(1,1),(3,1),(5,1)\}.$  Une seule repr\'esentation est associ\'ee \`a $H$, (voir le sous-graphe $K$ de la \figurename~\ref{exempleG1}).
\item Si $X=\{(0,0),(2,0),(3,0)\}$ et $Y=\{(1,1),(4,1),(5,1)\}.$  Deux repr\'esentations diff\'erentes sont associ\'ees \`a $H$, (voir le sous-graphes $K'$ de la \figurename~\ref{exempleG1}).
\end{enumerate}
\end{example}

\medskip

Nous associons \`a la repr\'esentation de $H$, ci-dessus, un mot ou une suite de $\{0,1\}^n$ obtenue en balayant les lignes $L_0,~L_1$, par une ligne verticale, de gauche \`a droite et en associant, \`a chaque point $v$ rencontr\'e, dans cet ordre, $0$ si $v$ se trouve sur $L_0$ et $1$ s'il est sur $L_1.$\\ Par exemple, pour
les sous-graphes de la \figurename~\ref{representationG1}, nous associons \`a $H_1$ le mot $S_1=0001101$ et \`a $H_2$ le mot $S_2=01001000$.

\medskip
Donc, \`a tout sous-graphe $H$ d'ordre $n$ de $G_1$, nous pouvons associer un mot $S_H=u_0u_1\cdots u_{n-1}$ o\`u $u_i\in\{0,1\}$ pour $i\in\{0,1,\cdots,n-1\}$, tel que $u_0=0$ et si $k$ est le plus grand indice pour lequel $u_k=1$ alors, $$l_0(S_H):=\vert\{i\in\{0,\cdots,k\}/u_i=0\}\vert\geq l_1(S_H):=\vert\{i\in\{0,\cdots,k\}/u_i=1\}\vert.$$

\medskip

Inversement, notons par $e$ le mot vide et soient $$\mathcal W:=\underset{n\in\mathbb N^*}\cup\{S=u_0u_1\cdots u_{n-1}\in\{0,1\}^n/ u_0=0,u_{n-1}=1\text{ et } l_0(S)\geq l_1(S)\}\cup\{e\}$$ et $$\mathcal O:=\{e,0,00,000,\cdots,0\cdots 0,\cdots\}$$
$\mathcal O$ est l'ensemble de tous les mots de longueur $n\in\mathbb N$ dont toutes les lettres sont \'egales \`a z\'eros, le mot de longueur $0$ \'etant \'egale \`a $e$. Posons $\mathcal S:=\{S_1S_2/ S_1\in \mathcal W,~S_2\in\mathcal O\}$.
Soit le mot $S\in\mathcal S$, alors
\begin{itemize}
\item si $S=e$, nous lui associons le sous-graphe vide,
\item sinon $S=u_0u_1\cdots u_{n-1},~ n\in\mathbb N^*$, nous lui associons le type d'isomorphie de ${G_1}_{\restriction_{X\cup Y}}$ o\`u $X=\{i\in\{0,\cdots,n-1\}/ u_i=0\}\times\{0\}$ et $Y=\{i\in\{0,\cdots,n-1\}/ u_i=1\}\times\{1\}.$
\end{itemize}

\noindent Par exemple pour $S=010101$ on a $X=\{0,2,4\}\times\{0\}$ et $Y=\{1,3,5\}\times\{1\}$ et pour $S=001101000$ on a $X=\{0,1,4,6,7,8\}\times\{0\}$ et $Y=\{2,3,5\}\times\{1\}.$\\

\noindent Soit $S=u_0u_1\cdots u_{n-1}\in\mathcal S$. Posons $\overline{0}=1,~\overline{1}=0$ et
$$\overline{S}=\left\{\begin{array}{ll}\overline{u}_{n-1} \cdots \overline{u}_1  \overline{u}_0& \text{ si } S\in\mathcal W\\
 \overline{S'}O& \text{ si } S=S'O \text{ avec }S'\in\mathcal W,~O\in\mathcal O\end{array}\right.$$

\begin{remark}
\begin{enumerate}
\item[1)] Pour tout $S\in\mathcal W$ nous avons:
    \begin{enumerate}
    \item[$a)$]  $l_0(S)>l_1(S)\Rightarrow\overline{S}\notin\mathcal S.$
    \item[$b)$]  $l_0(S)=l_1(S)\Leftrightarrow\overline{S}\in\mathcal W.$
    \end{enumerate}
\item[2)] Si $H$ est un sous-graphe connexe de $G_1$ alors $S_H\in\mathcal W$ et si $V(H)=X\cup Y$ avec $\vert X\vert=\vert Y\vert$ alors $\overline{\mathcal S_H}\in\mathcal W.$
\end{enumerate}
\end{remark}

\begin{lemma}\label{lem:2repreG1}
Soit $H=(X\cup Y,E)$ un sous-graphe connexe de $G_1$ d'ordre $n$ avec $\vert X\vert=\vert Y\vert$ et soit $S_H$ le mot correspondant. Si $S_H\neq \overline{S_H}$ alors $\overline{S_H}$ est un autre mot qui repr\'esente $H$.
\end{lemma}

\begin{proof}
Ceci est d\^u au fait que si $S_H$ correspond \`a la repr\'esentation de $H$ telle que les sommets de $X$ sont sur $L_0$ et ceux de $Y$ sur $L_1$ alors $\overline{S_H}$ est le mot associ\'e \`a la repr\'esentation qui place les sommets de $Y$ sur $L_0$ et ceux de $X$ sur $L_1$. En effet, si $S_H=u_0u_1\cdots u_{n-1}$ on associe \`a chaque terme $u_i$ pour $0\leq i\leq n-1$ le sommet $x^j_i$ avec $j=u_i$. Le sommet $x^j_i$ est le $(i+1)$-\`eme sommet de $H$ rencontr\'e en balayant les lignes $L_0,~L_1$ de gauche \`a droite, ce sommet est sur $L_0$ si $j=u_i=0$ et sur $L_1$ si $j=u_i=1$. Alors  $\overline{S_H}$ est le mot qui repr\'esente $f(H)$ o\`u $f$ est l'automorphisme\index{automorphisme} de $H$ donn\'e par:
$$\begin{array}{cccc}
f:& V_H&\rightarrow & V_H\\
  &x^j_i&\mapsto & x^{1-j}_{n-1-i}%
  \end{array}$$
  Il est clair que $f$ est bijective et que $\{x^j_i,x^{j'}_{i'}\}\in E$ si et seulement si $\{f(x^j_i),f(x^{j'}_{i'})\}\in E.$
  En effet; si $\{x^j_i,x^{j'}_{i'}\}\in E$ avec $i<i'$ alors n\'ecessairement $j=u_i=0$ et $j'=u_{i'}=1$ donc $f(x^j_i)=x^{1-j}_{n-1-i}$
  et $f(x^{j'}_{i'})=x^{1-j'}_{n-1-i'}$ avec $1-j=1$, $1-j'=0$ et $n-1-i>n-1-i'$.
  Or si $\overline{S_H}=\overline{u}_{n-1} \cdots \overline{u}_1 \overline{u}_0={u'}_0{u'}_1\cdots {u'}_{n-1}$
  alors $\overline{u}_i={u'}_{n-1-i}$ donc le sommet $x^{1-j}_{n-1-i}$ correspond au terme ${u'}_{n-1-i}=\overline{u_i}$ de
  $\overline{S_H}$ et $x^{1-j'}_{n-1-i'}$ correspond au terme ${u'}_{n-1-i'}=\overline{u_{i'}}$.
  L'implication inverse se d\'emontre de la m\^eme fa\c{c}on.
\end{proof}

\begin{example}
\begin{itemize}
\item Soit $H$ d\'efini par $X=\{0,2,3\}\times\{0\}$ et $Y=\{1,4,5\}\times\{1\}$. En pla\c{c}ant $X$ sur $L_0$ et $Y$ sur $L_1$ on a comme suite repr\'esentative $S=010011$ et en pla\c{c}ant $Y$ sur $L_0$ et $X$ sur $L_1$ on a comme suite repr\'esentative $S'=001101=\overline{S}.$
    \item Les sous-graphes de la \figurename~\ref{exempleG1} ont pour suites repr\'esentatives $S_1=010101=\overline{S}_1$ pour le sous-graphe $K$ et $S_2=010011$, $\overline{S}_2=001101$ pour le sous-graphe $K'$.
    \end{itemize}
\end{example}

\begin{lemma}\label{lem:suiteG1}
Soient $H$ et $H'$ deux sous-graphes de $G_1$ de m\^eme ordre $n$ et soient $S_H$ et $S_{H'}$ les suites correspondantes respectivement. % $S_H=SO$ et $S_{H'}=S'O'$ avec $S,~S'\in\mathcal W$ et $O,~O'\in\mathcal O$. Si $l_0(S)>l_1(S)$ et $l_0(S')>l_1(S')$
Alors:
\begin{center}
$H$ et  $H'$  sont isomorphes $\Leftrightarrow S_H=S_{H'}$  ou bien  $S_H=\overline{S_{H'}}.$%\footnote{On dira, dans ce cas, que $S_H$ et $\overline{S_{H'}}$ sont isomorphes.}
\end{center}
%Dans ce dernier cas ($S_H=\overline{S_{H'}}$), nous dirons que $S_H$ et $\overline{S_{H'}}$ sont isomorphes.
\end{lemma}

\begin{proof}\\ $\Leftarrow)$ \'evident par construction et d'apr\`es le Lemme \ref{lem:2repreG1}.\\ $\Rightarrow)$ Posons $S_H=SO$ et $S_{H'}=S'O'$ avec $S,~S'\in\mathcal W$ et $O,~O'\in\mathcal O$. Si $H$ et $H'$ sont isomorphes, ils ont le m\^eme nombre de composantes connexes de m\^eme ordre, donc $S$ et $S'$ sont de m\^eme longueur, il en est de m\^eme pour $O$ et $O'$.\\ Si $H$ et $H'$ sont des ind\'ependants d'ordre $n$, alors leur repr\'esentation est unique, dans ce cas $S=S'=e$ et $O=O'$ de taille $n$.\\ Si $H$ et $H'$ poss\`edent, chacun, une composante connexe d'ordre au moins $2$, alors, les bipartitions en stables de ces composantes sont uniques. Nous avons deux cas, ou bien, pour chacun des graphes, les deux stables ne sont pas de m\^eme cardinalit\'e (donc $l_0(S)>l_1(S)$ et $l_0(S')>l_1(S')$). $H$ et $H'$ \'etant isomorphes, les sommets qui co\"{\i}cident par cet isomorphisme se superposent et on a une seule fa\c{c}on de les repr\'esenter (le stable de cardinalit\'e maximum avec les sommets isol\'es sur $L_0$ et l'autre sur $L_1$) on a alors $S_H=S_{H'}$, ou bien les deux stables sont de m\^eme cardinalit\'e alors, en vertu du Lemme \ref{lem:2repreG1}, chacun des deux graphes poss\`ede deux suites repr\'esentatives (qui peuvent-\^etre \'egales). Soit $f$ un isomorphisme de $H$ sur $H'$. Si $f$ envoie les sommets de $H$ qui sont sur $L_0$ sur les sommets de $H'$ qui sont sur $L_0$ on a $S_H=S_{H'}$ et si $f$ envoie les sommet de $H$ qui sont sur $L_0$ sur les sommets de $H'$ qui sont sur $L_1$ on a $S_H=\overline{S_{H'}}.$
\end{proof}

%Lorsque $\overline{S}\in\mathcal S$ pour $S\in\mathcal S$, les suites $S$ et $\overline{S}$ sont dites isomorphes.

\subsection{Les ind\'ecomposables de $Age(G_1)$}
D'apr\`es la repr\'esentation des sous-graphes de $G_1$, il est \'evident qu'un sous-graphe $H$ poss\`ede un intervalle si et seulement si ou bien $H$ n'est pas connexe, ou bien dans la suite $S_H$ qui le repr\'esente, il existe deux entiers $j,r$ tels que $u_j=u_{j+1}=\cdots=u_{j+r}.$ Donc nous avons
\begin{fact}
Si $H$ est un sous-graphe  d'ordre $n\geq 4$ de $G_1$ et $S_H=u_0\cdots u_{n-1}$ la suite  qui le repr\'esente alors
$$H \text{ est ind\'ecomposable }\Leftrightarrow u_0=0,~ u_{n-1}=1 \text{ et } u_j\neq u_{j+1},\forall j,~0\leq j\leq n-2.$$
\end{fact}
Cette condition signifie que $n$ est pair. Ces sous-graphes ind\'ecomposables sont les graphes critiques de Schmerl et Trotter.
% est un ind\'ecomposable de $G_1$ si et seulement si la est telle que
%$u_0=0,~ u_{n-1}=1$ et $u_j\neq u_{j+1}$ pour tout $j,~0\leq j\leq n-2,$ ce qui signifie que $n$ est pair. %Ces sous-graphes sont des graphes critiques de Schmerl et Trotter.
En incluant les sous-graphes d'ordres $1$ et $2$, ces sous-graphes forment une cha\^{i}ne \`a partir de l'ordre $4$ telle que le sous-graphe $H$ qui se trouve au niveau $h$ pour $h\geq 3$ est d'ordre $n=2(h-1)$ et pour $h=1,2$ d'ordre $n=1,2$ respectivement. Donc $$h(H)\leq \vert H\vert\leq 2(h(H)-1).$$

 Il s'ensuit que l'ensemble des ind\'ecomposables est minimal dans la classe des graphes ind\'ecomposables. $Age(G_1)$ est donc ind-minimal d'apr\`es le Th\'eoreme \ref{theo:equiv-min}. %et donc minimal d'apr\`es le corollaire \ref{cor:ind-min-minimal}.

\subsection{Profil de l'\^age de $G_1$}

D'apr\`es ce qui pr\'ec\`ede, le profil de $G_1$ s'obtient, en \'enum\'erant, pour tout $n\in\mathbb N$, les suites de longueur $n$ appartenant \`a $\mathcal S$, en tenant compte du Lemme \ref{lem:suiteG1}.

\vspace{1mm}

Pour tout $n\in\mathbb N$, d\'esignons par $\varphi_1(n)$ le nombre de sous-graphes, non isomorphes, d'ordre $n$ de $G_1$ et par $m(l,n)$ le nombre de mots de longueur $n$, non isomorphes, de $\mathcal S$ ayant exactement $l$ termes \'egaux \`a $1$.

Nous avons:
\begin{itemize}
\item $m(l,n)=0$ pour $l>\dfrac{n}{2},~\forall n\in\mathbb N$.
\item $m(0,n)=1,~\forall n\in\mathbb N.$
\item $m(1,n)=n-1,~\forall n\in\mathbb N^*$.
\end{itemize}

\vspace{1mm}

\begin{lemma}
$$m(l,n)=\left\{\begin{array}{ll}
m(l,n-1)+\binom{n-2}{l-1}& \text{si } l<\dfrac{n}{2}\\
\dfrac{1}{2}(\binom{2l-2}{l-1}+2^{l-1})& \text{si } l=\dfrac{n}{2}%
\end{array}%
\right.$$
pour tout $n\in\mathbb N^*$ et $l\geq 1.$
\end{lemma}

\vspace{1mm}

\begin{proof}
Soit $n\in\mathbb N^*$. Posons\\
$m(l,n)=N_1+N_2$ o\`u\\
$N_1$ est le nombre de mots $S\in\mathcal W$ de longueur $n$ et\\
$N_2$ le nombre de mots $S'\in\mathcal S\setminus\mathcal W$ de longueur $n$. Nous distingons deux cas:

\textbf{Cas 1:} Si $n$ est impair alors n\'ecessairement $l<\dfrac{n}{2}$ et pour toute suite $S\in\mathcal S$ de longueur $n$, il est clair que  $\overline{S}\notin\mathcal S.$ Nous avons dans ce cas:\\
$N_1=\binom{n-2}{l-1}$ car $S=u_0\cdots u_{n-1}\in\mathcal W$ est telle que $u_0=0$ et $u_{n-1}=1$, reste \`a choisir $l-1$ termes de valeur $1$ parmi $n-2$ termes.\\
$N_2=m(l,n-1)$ car dans ce cas $u_0=u_{n-1}=0$, ce qui revient au nombre de suites, de longueur $n-1$ ayant $l$ termes \'egaux \`a $1$. D'o\`u\\
$m(l,n)=m(l,n-1)+\binom{n-2}{l-1}$.

\textbf{Cas 2:} Si $n$ est pair alors dans ce cas les suites $S$ de $\mathcal S\setminus\mathcal W$ sont telles que $l<\dfrac{n}{2}$ donc:\\
$N_2=m(l,n-1)$. Pour $N_1$ nous avons deux cas:
\begin{enumerate}
\item Si $l<\dfrac{n}{2}$ nous retrouvons la m\^eme valeur que pr\'ec\'edemment $N_1=\binom{n-2}{l-1}.$
\item Si $l=\dfrac{n}{2}$ alors les suites $S\in\mathcal W$ concern\'ees sont telles que $\overline{S}\in\mathcal W$. Posons\\
    $N_1={N'}_1+{N'}_2$ o\`u\\
${N'}_1$ est le nombre de suites $S\in\mathcal W$ de longueur $n$ telles que $S=\overline{S}$ et\\
${N'}_2$ le nombre de suites $S'\in\mathcal W$ de longueur $n$  telles que $S'\neq\overline{S'}$. Alors\\
$\bullet~ {N'}_1=\underset{k=0}{\overset{l-1}{\sum}}\binom{l-1}{k}=2^{l-1}$ car dans ce cas les suites $S$ consern\'ees sont de la forme $S=S_1 S_2$ avec $S_2=\overline{S}_1$ cela revient \`a choisir les termes de $S_1$ qui valent $1$.\\
\medskip
$\bullet ~{N'}_2=\dfrac{1}{2}(\binom{n-2}{l-1}-{N'}_1)=\dfrac{1}{2}(\binom{n-2}{l-1}-2^{l-1}).$  %nombre de toutes les suites \`a $n/2$ termes \'egaux \`a $1={N'}_1+2.{N'}_2$

 Donc pour r\'esumer, nous avons dans ce cas:\\
\quad - Si $l<\dfrac{n}{2}$ alors $m(l,n)=m(l,n-1)+\binom{n-2}{l-1};$\\
\quad - Si $l=\dfrac{n}{2}$ alors $m(l,n)=\dfrac{1}{2}(\binom{n-2}{l-1}-2^{l-1}).$
\end{enumerate}
Nous retrouvons alors la formule
$$m(l,n)=\left\{\begin{array}{ll}
m(l,n-1)+\binom{n-2}{l-1}& \text{si } l<\dfrac{n}{2}\\
\dfrac{1}{2}(\binom{2l-2}{l-1}+2^{l-1})& \text{si } l=\dfrac{n}{2}%
\end{array}%
\right.$$
\end{proof}

\vspace{2mm}

Les premi\`eres valeurs du profil de $Age(G_1)$ pour $n=0,1,2,3,4,5,6,7,8,9$ sont respectivement  $1,1,2,3,6,10,20,36,72,136.$ Pour les autres nous avons le lemme suivant:

\begin{lemma}\label{lem:profG1}
Le profil $\varphi_1$ de l'\^age de $G_1$ est exponentiel, il est donn\'e par:
$$\left\{\begin{array}{ll}
\varphi_1(0)=\varphi_1(1)=1,~\varphi_1(2)=2,&\\
\varphi_1(n)=\varphi_1(n-1)+2^{n-3}& \text{si }n \text{ impair } n\geq 3\\
\varphi_1(n)=\varphi_1(n-1)+2^{n-3}+2^{\frac{n-4}{2}}& \text{si }n \text{ pair } n\geq 4%
\end{array}\right.$$
En outre, $\varphi_1\simeq\dfrac{2^n}{4}$.
\end{lemma}

\begin{proof}
Pour $n=0,1,2$ le r\'esultat est imm\'ediat. Soit $n\geq 3$ nous distingons les deux cas:

$1)$ \textbf{Si $n$ est impair:} posons $n=2k+1$ avec $k\geq 1$ alors\\
$\begin{array}{lll}
\varphi_1(2k+1)&=\underset{l=0}{\overset{k}\sum}m(l,2k+1)&=m(0,2k+1)+\underset{l=1}{\overset{k}\sum}(m(l,2k)+\binom{2k-1}{l-1})\\
&&=\underset{l=0}{\overset{k}\sum}m(l,2k)+\underset{l=1}{\overset{k}\sum}\binom{2k-1}{l-1}\\
&&=\varphi_1(2k)+\underset{j=0}{\overset{k-1}\sum}\binom{2k-1}{j}%
\end{array}$

Ceci en utilisant les relations $m(0,n)=1,~\forall n\in\mathbb N$ et $\binom{n}{j}=\binom{n-j}{n}$. Nous avons alors dans ce cas,
$$\varphi_1(n)=\varphi_1(n-1)+2^{n-3},~\forall n\geq 3,~ n \text{ impair.}$$

$2)$ \textbf{Si $n$ est pair:} Posons $n=2k$ avec $k\geq 2$ nous avons alors\\
$\begin{array}{lll}
\varphi_1(2k)&=\underset{l=0}{\overset{k}\sum}m(l,2k)&=m(0,2k)+\underset{l=1}{\overset{k-1}\sum}m(l,2k)+m(k,2k)\\
&&=m(0,2k)+\underset{l=1}{\overset{k-1}\sum}(m(l,2k-1)+\binom{2k-2}{l-1})+\dfrac{1}{2}(\binom{2k-2}{k-1}+2^{k-1})\\
&&=\underset{l=0}{\overset{k-1}\sum}(m(l,2k-1)+\underset{j=0}{\overset{k-2}\sum}\binom{2k-2}{j}+\dfrac{1}{2}(\binom{2k-2}{k-1}+2^{k-1})\\
&&=\varphi_1(2k-1)+\underset{j=0}{\overset{k-2}\sum}\binom{2k-2}{j}+\dfrac{1}{2}\binom{2k-2}{k-1}+2^{k-2}%
\end{array}$

En utilisant le fait que $\underset{j=0}{\overset{k-2}\sum}\binom{2k-2}{j}=2^{2k-3}-\dfrac{1}{2}\binom{2k-2}{k-1}$ car $\binom{n}{j}=\binom{n-j}{n}$ et $\underset{j=0}{\overset{n}\sum}\binom{n}{j}=2^n,$ nous trouvons
 $$\varphi_1(2k)=\varphi_1(2k-1)+2^{2k-3}+2^{k-2}.$$
 D'o\`u
 $$\varphi_1(n)=\varphi_1(n-1)+2^{n-3}+2^{\frac{n-4}{2}},~\forall n\geq 4,~ n \text{ pair.}$$
 Pour la croissance asymptotique, elle a \'et\'e calcul\'ee sur Maple gr\^ace à la fonction g\'en\'eratrice donn\'ee dans la proposition qui suit.
\end{proof}

\begin{proposition}
La fonction g\'en\'eratrice\index{fonction!g\'en\'eratrice} de $Age(G_1)$, l'\^age de $G_1$, est rationnelle\index{fonction!g\'en\'eratrice!rationnelle} et est donn\'ee par:
$$F_{G_1}(x)=\dfrac{1-x-2x^2+x^3}{(1-2x)(1-2x^2)}.$$
\end{proposition}

\begin{proof}
Nous avons
$$\begin{array}{rl}
A=F_{G_1}(x)=\underset{n\geq 0}\sum\varphi(n)x^n=& \varphi(0)+\varphi(1)x+\varphi(2)x^2
        +\underset{\begin{array}{c}
        n\geq 4\\
        n \text{\small{ pair}}\end{array}}\sum\varphi(n)x^n+\underset{\begin{array}{c}
    n\geq 3\\
    n \text{\small{ impair}}\end{array}}\sum\varphi(n)x^n.\\
    =& \varphi(0)+\varphi(1)x+\varphi(2)x^2+\underset{\begin{array}{c}
        n\geq 4\\
        n \text{\small{ pair}}\end{array}}\sum(\varphi(n-1)+2^{n-3}+2^{\frac{n-4}{2}})x^n\\
        &\qquad +\underset{\begin{array}{c}
    n\geq 3\\
    n \text{\small{ impair}}\end{array}}\sum(\varphi(n-1)+2^{n-3})x^n.\\
   = &  \varphi(0)+\varphi(1)x+\varphi(2)x^2+\underset{n\geq 3}\sum\varphi(n-1)x^n+\underset{n\geq 3}\sum 2^{n-3}x^n\\
    & \qquad +\underset{\begin{array}{c}
        n\geq 4\\
        n \text{\small{ pair}}\end{array}}\sum 2^{\frac{n-4}{2}}x^n.\\
    \end{array}$$
    Comme $$\underset{n\geq 3}\sum\varphi(n-1)x^n=\underset{n\geq 0}\sum\varphi(n)x^{n+1}-\varphi(0)x-\varphi(1)x^2.$$
    Nous obtenons alors:
 $$\begin{array}{rl}
A=F_{G_1}(x)= & \varphi(0)+\varphi(1)x+\varphi(2)x^2+\underset{n\geq 0}\sum\varphi(n)x^{n+1}-\varphi(0)x-\varphi(1)x^2\\
&\qquad  +\underset{n\geq 3}\sum 2^{n-3}x^n+\underset{n\geq 2}\sum 2^{n-2}x^{2n}.\end{array}$$
Nous avons \\ $\underset{n\geq 3}\sum 2^{n-3}x^n=x^3\underset{n\geq 0}\sum (2x)^n=\dfrac{x^3}{1-2x}$ et \\
$\underset{n\geq 2}\sum 2^{n-2}x^{2n}=x^4\underset{n\geq 0}\sum (2x^2)^n=\dfrac{x^4}{1-2x^2}.$ D'o\`u l'on d\'eduit:
$$A=xA+\dfrac{x^3}{1-2x}+\dfrac{x^4}{1-2x^2}+1+x^2.$$ Ce qui donne
$$F_{G_1}(x)=\dfrac{1-x-2x^2+x^3}{(1-2x)(1-2x^2)}.$$
\end{proof}

%\bigskip
%**************************************************************************************
%$$$$$$$$$$$$$$$$$$$$$$$$$$$$$$$$$$$$$$$$$$$$$$$$$$$$$$$$$$$$$$$$$$$$$$$$$$$$$$$$

\section{Le graphe $G_2$}

Le graphe $G_2:=(V_2,E_2)$ est biparti, son ensemble de sommets $V_2:=\mathbb N \times \{0,1\}$ se d\'ecompose  en deux stables disjoints $A:=\mathbb N \times \{0\}$ et $B:=\mathbb N \times \{1\}$. Une paire $\{(i,0),(j,1)\}$ est une ar\^ete de  $G_1$ si et seulement si $j=i$ ou $j=i+1$. Le graphe $G_2$ est le chemin infinie\index{chemin!infini} (en tant que graphe). Le graphe $G_2$ et son graphe compl\'ementaire sont des graphes de permutations, $G_2$ est le graphe de comparabilit\'e du zigzags\index{zigzags} (introduit dans la section \ref{sec:conjecture} du chapitre \ref{chap:exemple-conjecture}), il est \'egalement le graphe d'incomparabilit\'e de l'ordre $P:=(\mathbb N\times\{0,1\}, L_1\cap L_2)$ o\`u $L_1:=(0,1)<(0,0)<(2,1)<(1,0)<\cdots<(n,1)<(n-1,0)<\cdots$ et $L_2:=(0,0)<(0,1)<(1,0)<(1,1)<\cdots<(n,0)<(n,1)<\cdots$.

\vspace{2mm}

Les sous-graphes ind\'ecomposables de $G_2$ sont, en plus des ind\'ecomposables d'ordres $1$ et $2$, les sous-graphes connexes de $G_2$ d'ordres $n> 3$. Il y'en a un %d'ordre $1$, deux d'ordre $2$ et un
pour chaque ordre $n$, $n\geq 4$. On y inclura l'ind\'ecomposable d'ordre $1$ et les deux ind\'ecomposables d'ordre $2$. Ces graphes ind\'ecomposables  forment une cha\^{i}ne telle que l'ind\'ecomposable d'ordre $n$ se trouve au niveau $n-1$ pour $n\geq 4$. Ceux d'ordres $1$ et $2$ se trouvent respectivement aux niveaux $1$ et $2$. L'ensemble de ces sous-graphes est minimal dans la classe des ind\'ecomposables. $Age(G_2)$ (et donc $Age(\overline{G_2})$ l'\^age de son compl\'ementaire) est ind-minimal (Th\'eor\`eme \ref{theo:equiv-min}). %donc minimal (corollaire \ref{cor:ind-min-minimal}).

\begin{lemma}
Le profil, $\varphi_2$, de $Age(G_2)$ est la fonction partition d'entiers\index{fonction!partition d'entier} $\mathfrak p(n)$. Sa fonction g\'en\'eratrice est donn\'ee par:
$$F_{G_2}(x)=\underset{k=1}{\overset{\infty}\prod}\dfrac{1}{1-x^k}.$$
\end{lemma}

\begin{proof}
Les sous-graphes d'ordre $n$ de $G_2$ peuvent-\^etre repr\'esent\'es par les partitions\index{partition!d'un entier} de $n$. En effet, pour identifier un sous-graphe de $G_2$, \`a l'isomorphie pr\`es, il faut et il suffit de d\'eterminer les tailles de ses composantes connexes, la somme de toutes ces tailles \'etant \'egale \`a l'ordre du sous-graphe. \end{proof}

%***************************************************************************************
%************************************************************************************
%\bigskip

\section{Le graphe  $G_3$}\label{subsec:graphe $G_3$}

Le graphe $G_3:=(V_3,E_3)$ (repr\'esent\'e dans la \figurename~\ref{critique}) est biparti, son ensemble de sommets $V_3:=(\mathbb N \times \{0,1\})\cup \{c\}$, se d\'ecompose en deux sous-ensembles disjoints $A_3:=(\mathbb N \times \{0\})\cup\{c\}$ et $B:=\mathbb N\times \{1\}$ qui forment des stables. Une paire $\{(i,0),(j,1)\}$ est une ar\^ete de $G_3$ si $i=j$ et nous ajoutons \`a $E_3$ toutes les paires $\{c,(i,1)\}$ pour $i\in \mathbb N$.

\vspace{1mm}

Remarquons que pour tout $x\in V_3\setminus\{c\}$, le graphe $G_3$ s'abrite dans ${G_3}_{\restriction_{V_3\setminus\{x\}}}$, donc $x\notin Ker(G_3)$ et $Age(G_3)\neq Age({G_3}_{\restriction_{V_3\setminus\{c\}}})$, donc $Ker(G_3)=\{c\}$.

\subsection{Repr\'esentation des sous-graphes de $G_3$}

Soit $H$ un sous-graphe de $G_3$, nous avons deux cas:

\textbf{Cas 1:} $H$ ne contient pas le sommet $c$ ($H$ est un sous-graphe de $G_3\setminus \{c\}$) alors $H$ est de la forme $p.K_2\oplus q$ o\`u $''\oplus''$ est la somme directe, $K_2$ est la clique \`a deux sommets et $q$ le stable d'ordre $q$. Donc, tout sous-graphe $H$ de $G_3\setminus \{c\}$ est enti\`erement d\'efini par son ordre $n$ et le nombre $p$ de copies de $K_2$ qu'il contient. Nous pouvons alors coder chaque type d'isomorphie des sous-graphes $H$, d'ordre $n$, de $G_3\setminus \{c\}$ par un couple d'entiers $(n,p)$ tel que $p$ repr\'esente le nombre de copies de $K_2$ contenus dans $H$ et $p\leq \dfrac{n}{2}$. Inversement, \`a tout couple d'entiens $(n,p)$ avec $p\leq \dfrac{n}{2}$, nous pouvons  associer le type d'isomorphie $T_{n,p}$ de ${G_3}_{\restriction_X}$ o\`u $X=\{0,\cdots, n-p\}\times \{0\}\cup \{0,\cdots,p-1\}\times\{1\}$.

\begin{observation}\label{obs:G3}
Deux sous-graphes de $G_3\setminus \{c\}$ sont isomorphes si et seulement si ils sont cod\'es par le m\^eme couple.
\end{observation}
\medskip

\textbf{Cas 2:} $H$ contient le sommet $c$ alors $H$ est caract\'eris\'e par son ordre $n$, le nombre $k$ de copies de $K_2$ qu'il contient et le degr\'e $p$ du sommet $c$ dans $H$. Nous associons alors, \`a tout sous-graphe $H$ de $G_3$ contenant le sommet $c$, le triplet d'entiers $t_H:=(n,k,p)$ tel que $p\geq k$ et $n\geq k+p+1,$ donc $k\leq\dfrac{n-1}{2}.$ \\Posons $\mathcal T:=\{t:=(n,k,p)\in\mathbb N^3/~p\geq k,~n\geq k+p+1\}.$

\vspace{1mm}

\begin{figure}
\centering
\input{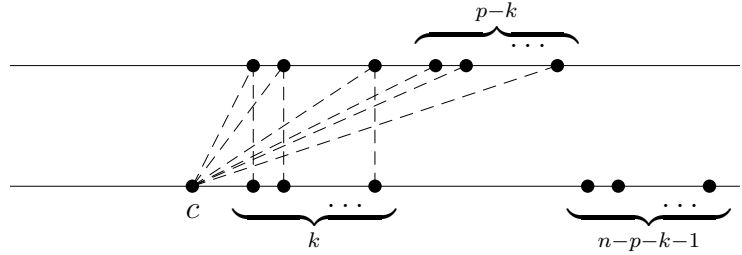}
\caption{\label{representationG3}Repr\'esentation des sous-graphes de $G_3$.}
\end{figure}

Inversement, \`a tout triplet $t:=(n,k,p)\in\mathcal T$, nous associons le type d'isomorphie $T_t$ de ${G_3}_{\restriction_Y}$ o\`u
$Y=\{c\} \cup (\{0,\cdots,k-1\}\cup\{p,\cdots,n-k-2\})\times \{0\}\cup \{0,\cdots,p-1\}\times\{1\}$ (voir \figurename~\ref{representationG3}).

\begin{remarks}\label{rem:G3}
\begin{enumerate}
\item    Un sous-graphe $H$ de $Age (G_3)$ qui contient le sommet $c$ peut-\^etre isomorphe \`a un sous-graphe de $G_3\setminus \{c\}$, en effet nous avons:
\begin{enumerate}
 \item Pour chaque entier $n\geq 2$, le sous-graphe isomorphe \`a $K_2\oplus (n-2)$ peut-\^etre repr\'esent\'e soit par le triplet $(n,0,1)$ soit par le couple $(n,1)$.
\item Pour tout entier $n$, le stable d'ordre $n$ peut-\^etre repr\'esent\'e soit par le triplet $(n,0,0)$ soit par le couple $(n,0)$.
\end{enumerate}
\item De m\^eme nous avons, pour tout entier $n\geq 3$, le sous-graphe isomorphe \`a $P_3\oplus (n-3)$ peut-\^etre repr\'esent\'e par l'un des deux triplets $(n,0,2)$ et $(n,1,1)$.
\end{enumerate}
\end{remarks}

Posons $$\mathcal T':=\{t:=(n,k,p)\in\mathcal T/n\geq 4,~ k\neq 0\}\cup\{t:=(n,0,p)\in\mathcal T/n\geq 4,~ p\geq 3\}.$$
Nous avons le lemme suivant:

\begin{lemma}\label{lemma1}
 Soient  $t:=(n,k,p)$ et $t':=(n',k',p')$ deux triplets de $\mathcal T$ et soient  $T_t$ et $T_{t'}$ les sous-graphes correspondants respectifs. Si $t,t'\in\mathcal T'$ alors
       $$T_t\text{ et } T_{t'} \text{ sont isomorphes si et seulement si } t=t'.$$
        De plus, $T_t$ et $T_{t'}$ ne sont pas des sous-graphes de  $G_3\setminus \{c\}$.
\end{lemma}

\begin{proof}
 La condition suffisante est \'evidente. Pour la condition n\'ecessaire,
supposons que $T_t$ et $T_{t'}$ sont isomorphes mais $t\neq t'$,  donc $n=n'$. Soit $f$ cet isomorphisme.

\textbf{Cas 1:} Si  $k\neq k'$, supposons que $k<k'$, nous avons les sous-cas suivant:\newline
 \qquad\textbf{a)} Si $k\geq 1$ alors $k'\geq 2$ ce qui signifie que dans $T_{t'}$ le sommet $c$ a un degr\'e au moins \'egal \`a $2$ et $k'$ sommets, distincts de $c$, sont de degr\'e $2$ (sont reli\'es \`a $c$), donc, $c$ est invariant par $f$, c'est \`a dire $f(c)=c$, il s'ensuit que $p=d_{T_t}(c)=d_{T_{t'}}(c)=p'$, où $d_T(c)$ repr\'esente le degr\'e de $c$ dans le graphe $T$. Mais dans $T_t$, nous avons $k$ sommets diff\'erents de $c$ dont le degr\'e est $2$. Comme $k'>k$ les deux graphes ne peuvent-\^etre isomorphes ce qui contredit l'hypoth\`ese.  \\
\qquad\textbf{b)} Si $k= 0$ alors $p\geq 3$, donc $k'\geq 1$ et $f(c)=c$ d'o\`u $p'\geq 3.$  %(il ne peut y avoir un autre sommet de degr\'e $3$).
Dans ce cas  $T_{t'}$ contient $P_4$ mais $T_t$ ne le contient pas. Il s'ensuit que $T_t$ et $T_{t'}$ ne sont pas isomorphes ce qui contredit l'hypoth\`ese.
\bigskip

\textbf{Cas 2:} Si $k=k'$ alors $p\neq p'$, donc si $f(c)=c$ nous ne pouvons avoir isomorphisme car $p=d_{T_t}(c)\neq p'=d_{T_{t'}}(c)$ %car $t\neq t'$
 et si $f(c)\neq c$ alors nous avons n\'ecessairement $d_{T_t}(c)=2$ et $k=0$. Donc $T_t$ est cod\'e par $(n,0,2)$ ce qui contredit l'hypoth\`ese.

\vspace{1mm}

Pour la deuxi\`eme partie du lemme, les sous-graphes $T_t$ et $T_{t'}$ ne peuvent s'abriter dans $G_3\setminus \{c\}$ car ils contiennent, chacun, un sommet de degr\'e sup\'erieur \`a $1$.
\end{proof}

\subsection{Les ind\'ecomposables de $Age(G_3)$}
Il est clair que tous les sous-graphes ind\'ecomposables, d'ordre $n\geq 4$, de  $G_3$, contiennent le sommet $c$, car tout sous graphe ind\'ecomposable d'ordre au moins quatre doit contenir le graphe $P_4$. On a alors:
\begin{lemma}
Soit $H$ un sous-graphe contenant le sommet $c$ et ayant au moins quatres sommets et $t_H=(n,k,p)$ le triplet qui lui est associ\'e. Alors
$$H \text{ est ind\'ecomposable }\Leftrightarrow k\geq 1,~p-k\leq 1 \text{ et } n=p+k+1.$$
\end{lemma}
\begin{proof}
Evident d'apr\`es la repr\'esentation \figurename~\ref{representationG3}.
\end{proof}

\bigskip

\begin{figure}
\centering
\input{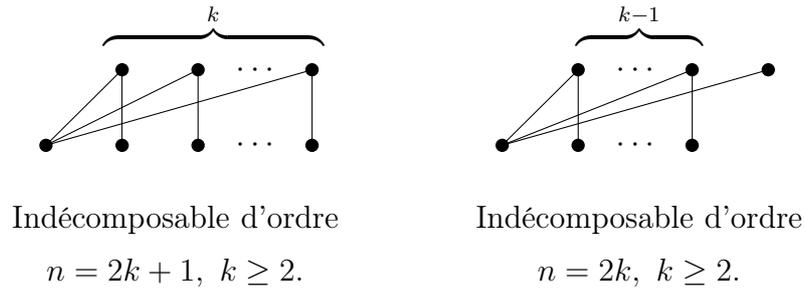}
\caption{\label{indecG3}Les sous-graphes ind\'ecomposables de $G_3$.}
\end{figure}

Il y'a un ind\'ecomposable pour chaque ordre $n\geq 4$, ils sont donn\'es dans la \figurename~\ref{indecG3}. Ces ind\'ecomposables  forment une cha\^{i}ne. En incluant, \`a cet ensemble, les ind\'ecomposables de tailles au plus $2$, le sous-graphe d'ordre $n$ se trouve au niveau $n-1$ de la cha\^{i}ne, pour $n\geq 4$. La classe des ind\'ecomposables de $Age(G_3)$ est minimale dans la classe des graphes ind\'ecomposables. Donc $Age(G_3)$ est ind-minimal.

\subsection{Profil de $Age(G_3)$}\label{subsec:profiG3}
Les premi\`eres valeurs du profil de $Age(G_3)$ pour $n=0,1,2,3,4,5,6,7$ sont respectivement $1$, $1$, $2$, $3$, $6$, $9$, $13, 17.$

\vspace{1mm}

D\'esignons par $\varphi_3$ la fonction profil de $Age(G_3)$.

\begin{lemma}\label{lem:profG3}
Le profil $\varphi_3$ est quasi-polynomial, il est donn\'e pour tout $n\in\mathbb N$ par:
$$\left\{\begin{array}{lll}
\varphi_3(0)=\varphi_3(1)=1,~\varphi_3(2)=2,& &\\
\varphi_3(n)=\dfrac{n^2+4n-8}{4}& &\text{si }n \text{ est pair } n\geq 4\\
\varphi_3(n)=\dfrac{n^2+4n-9}{4}& &\text{si }n \text{ est impair } n\geq 3%
\end{array}%
\right.$$
En outre, $\varphi_3\simeq\dfrac{n^2}{4}$
\end{lemma}

\begin{proof}
D'apr\`es ce qui pr\'ec\`ede, nous avons:
\begin{equation}
\forall n\in\mathbb N,~~\varphi_3(n)=\varphi'_3(n)+\varphi''_3(n)\label{eq:profilG3}
\end{equation}
o\`u $\varphi'_3(n)$ est le nombre de sous-graphes \`a $n$ sommets de $G_3$, non isomorphes,  qui s'abritent dans $G_3\setminus\{c\}$
et $\varphi''_3(n)$ est le nombre de sous-graphes \`a $n$ sommets de $G_3$, non isomorphes, qui ne s'abritent pas dans $G_3\setminus\{c\}.$\\
D'apr\`es l'Observation \ref{obs:G3} nous avons: $$\varphi'_3(n)=\vert\{(n,p)/ p\leq \dfrac{n}{2}, p\in\mathbb N\}\vert=\lfloor\frac{n}{2}\rfloor +1,\forall n\in\mathbb N.$$
D'apr\`es la Remarque \ref{rem:G3} et le Lemme \ref{lemma1}, le nombre $\varphi''_3(n)$ de types d'isomorphie  de sous-graphes  d'ordre $n$, qui contiennent le sommet $c$ et qui ne s'abritent pas dans $G_3\setminus\{c\}$ peut-\^etre trouv\'e en d\'enombrant tous les triplets $(n,k,p)\in\mathcal T'$.
Donc,
$$\varphi''_3(n)=\left\{\begin{array}{ll}
0& \text{ pour }n\leq 2\\
1& \text{ pour } n=3\\
\vert\{(k,p)\in\mathbb N^2/ (n,k,p)\in\mathcal T'\}\vert& \text{ pour } n\geq 4%
\end{array}%
\right.$$
Pour $n\geq 4$, nous d\'enombrons tous les couples $(k,p)$ tels que $k\neq 0$, $k\leq p\leq n-k-1$ et $k\leq\lfloor\frac{n-1}{2}\rfloor$ ou bien $k=0$ et $3\leq p\leq n-1$.

Ce nombre, apr\`es calcul est:
$$\varphi''_3(n)=\left\{
\begin{array}{ll}
0& \text{ pour }n\leq 2\\
(\frac{n}{2})^2+\frac{n}{2}-3 & \text{ pour }\; n\geq 4\text{ et }n\text{ pair} \\
(\frac{n-1}{2})^2+n-3 & \text{ pour }\; n\geq 3 \text{ et }n\text{ impair}%
\end{array}%
\right. $$

En remplaçant dans l'\'equation \eqref{eq:profilG3}, nous retrouvons la formule de $\varphi_3$ donn\'ee.

\vspace{1mm}

Il est clair que le profil de $Age(G_3)$ est asymptotiquement \'equivalent \`a la fonction $\dfrac{n^2}{4}$ et qu'il est quasi-polynomial. En effet, pour $n\geq 3$ nous avons
$$\varphi_3(n)=b_2(n)n^2+b_1(n)n+b_0(n),$$
où $b_0(n)=-\frac{9}{4}$  pour $n$ impair et $b_0(n)=-2$ pour $n$ pair, $b_1(n)=b_2(n)=\frac{1}{4}$ pour tout $n$.
\end{proof}

\begin{proposition}
La fonction g\'en\'eratrice de $Age(G_3)$ est rationnelle et est donn\'ee par:
$$F_{G_3}(x)=\dfrac{1-x+x^3+x^4-x^6}{(1-x)^2(1-x^2)}.$$
\end{proposition}

\begin{proof}
En utilisant le Lemme \ref{lem:profG3} on a:
$$\begin{array}{rl}
F_{G_3}(x)=\underset{n\geq 0}\sum\varphi_3(n)x^n=&\underset{n\geq 0}\sum\varphi_3(2n)x^{2n}+\underset{n\geq 0}\sum\varphi_3(2n+1)x^{2n+1}\\
=& \underset{n\geq 2}\sum(n^2+2n-2)x^{2n}+\underset{n\geq 1}\sum(n^2+3n-1)x^{2n+1}+1+x+2x^2
\end{array}$$
Apr\`es calcul nous obtenons la fonction donn\'ee.
\end{proof}

%********************************************************************************
%$$$$$$$$$$$$$$$$$$$$$$$$$$$$$$$$$$$$$$$$$$$$$$$$$$$$$$$$$$$$$$$$$$$$$$$$$$$$$$$$$$$$$

%\bigskip
\section{Le graphe  $G_4$}\label{subsec:graphe $G_4$}

Le graphe $G_4:=(V_4,E_4)$ est le graphe d\'efini par $V_4:=\mathbb N \times \{0,1\}$ qui se d\'ecompose en deux sous-ensembles, $A=\mathbb N \times \{0\}$ qui forme une clique et $B=\mathbb N \times \{1\}$ qui est un stable. Une paire de sommets $\{(i,0),(j,1)\}$ est une ar\^ete dans $G_4$ si $i=j$.

\subsection{Repr\'esentation des sous-graphes de $G_4$}
Un sous-graphe $H$ de $G_4$ est form\'e d'une clique \`a $n_0$ sommets, parmi lesquels $p$ sont reli\'es \`a des sommets pendants et d'un stable form\'e par des sommets isol\'es. $H$ poss\`ede, au plus, une composante connexe ayant plus d'un sommet.

\begin{figure}
\centering
\input{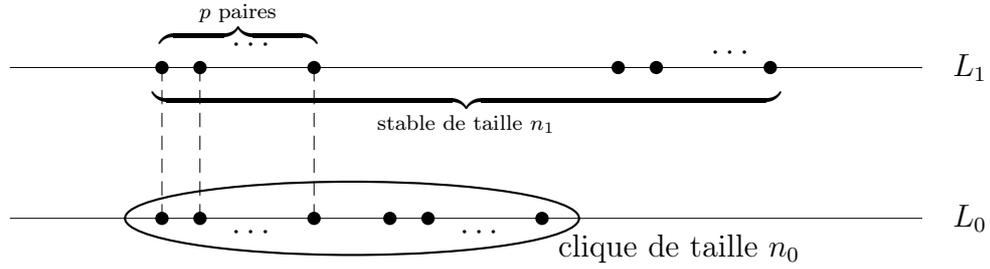}
\caption{\label{representationG4}Repr\'esentation des sous-graphes de $G_4$.}
\end{figure}

Remarquons que $G_4$ s'obtient en supprimant les ar\^etes $\{(i,1),(i',1)\}$ pour $i\neq i'$, au compl\'ementaire de $G_0$. Comme pour le graphe $G_0$, nous d\'esignons par \emph{paire} de $G_4$, le sous-ensemble de sommets $\{(i,0),(i,1)\}$ pour $i\in\mathbb N$. Soit $H=(X,F)$ un sous-graphe d'ordre $n$ de $G_4$, nous repr\'esentons les sommets de $X$ sur deux lignes parall\`eles horizontales $L_0$, $L_1$ tels que les sommets de $X\cap A$ sont sur $L_0$ et ceux de $X\cap B$ sur $L_1$ de telle sorte que toutes les paires de $G_4$ qui appartiennent \`a $H$ soient \`a gauche, les sommets d'une m\^eme paire \'etant align\'es verticalement. Les sommets restants sont plac\'es \`a droite de ces paires avec ceux sur $L_1$ compl\`etement \`a droite de ceux sur $L_0$ (\figurename~\ref{representationG4}).

\vspace{1mm}

Avec cette repr\'esentation, nous pouvons associer \`a chaque membre $H$ de $Age(G_4)$, d'ordre $n$, un triplet d'entiers $t_H:=(p,n_0,n_1)$, o\`u $p$ est le  nombre de paires de $H$, $n_0$ la taille de sa clique maximum et $n_1=n-n_0.$ %la taille de son stable maximum,
Ces entiers v\'erifient:
\begin{equation}\label{eq:equationG4}
 p\leq min(n_0, n_1)~\text{et } n_0+n_1=n.
 \end{equation}

Inversement, \`a tout triplet $t:=(p,n_0,n_1)$ v\'erifiant les conditions \eqref{eq:equationG4} ci-dessus, avec $n_0+n_1\geq 1$,  nous associons le type d'isomorphie $H_t$ de ${G_4}_{\restriction_X}$ o\`u $X=\{0,\cdots,n_0-1\}\times \{0\}\cup (\{0,\cdots,p-1\}\cup\{n_0,\cdots,n_0+n_1-p-1\})\times \{1\}.$ Il est \'evident que si $t=t'$  alors $H_t=H_{t'}$, mais si nous avons deux sous-graphes $H$ et $H'$ isomorphes nous n'avons pas n\'ecessairement $t_H=t_{H'}$ comme nous pouvons le voir dans la remarque qui suit:

\begin{remarks}\label{rem:G4}
\begin{enumerate}
\item Les sous-graphes associ\'es aux deux triplets $(0,1,n)$ et $(0,0,n+1)$ sont isomorphes pour tout $n\in\mathbb N$ (ils repr\'esentent l'ind\'ependant \`a $n+1$ sommets).
\item Les sous-graphes associ\'es aux deux triplets $(0,2,n)$ et $(1,1,n+1)$ sont isomorphes pour tout $n\in\mathbb N$ (ils repr\'esentent $K_2\oplus n$).
\item Si $H$ est un sous-graphe de $G_4$ et $t_H=(p,n_0,n_1)$ le triplet associ\'e alors  $H$ est connexe si et seulement si $n_1=p$.
\end{enumerate}
\end{remarks}

\begin{fact}
Si $t$ et $t'$ sont deux triplets v\'erifiant les conditions \eqref{eq:equationG4} avec $t$ et $t'$ n'appartenant pas \`a l'ensemble $\{(0,1,n),(0,2,n),~n\in\mathbb N\}$ alors
$$t=t'\Leftrightarrow H_t \text{ et } H_{t'}\text{ sont isomorphes.}$$
\end{fact}

\begin{proof}
La preuve est \'evidente d'apr\`es la structure des sous-graphes de $G_4$ et les Remarques \ref{rem:G4}.
\end{proof}

\subsection{Les ind\'ecomposables de $Age(G_4)$}
%Il est clair, d'ap\`es la repr\'esentation des sous-graphes de $G_4$ qu'un sous-graphe $H$ ind\'ecomposable, d'ordres au moins quatre, est tel que la taille de son stable maximum est \'egale au nombre de ses  paires et la taille de sa clique maximum comprise entre le nombre de paires et le nombre de paires plus un.

\begin{lemma}
Soit $H$ un sous-graphe ayant au moins quatre sommets et $t_H$ le triplet qui lui est associ\'e. Alors
$$H \text{ est ind\'ecomposable  }\Leftrightarrow t_H\in\{(p,p,p),(p,p+1,p),~p\geq 2\}.$$
\end{lemma}
\begin{proof}
Evident d'apr\`es la repr\'esentation \figurename~\ref{representationG4}.
\end{proof}

\bigskip

Il y'a un ind\'ecomposable pour chaque ordre $n\geq 4$. Ces ind\'ecomposables  forment une cha\^{i}ne. En incluant, \`a cet ensemble, les ind\'ecomposables d'ordre au plus $2$, le sous-graphe ind\'ecomposable d'ordre $n$ se trouve au niveau $n-1$ de la cha\^{i}ne, pour $n\geq 4$. La classe des ind\'ecomposables de $Age(G_4)$ est donc minimale dans $Ind(\Omega_1)$. Par cons\'equent,  $Age(G_4)$ est ind-minimal d'apr\`es le Th\'eor\`eme \ref{theo:equiv-min}. %il s'ensuit qu'il est minimal d'apr\`es le corollaire \ref{cor:ind-min-minimal}.

\subsection{Profil de $Age(G_4)$}\label{subsec:profilG4}
Les premi\`eres valeurs du profil pour $n=0,1,2,3,4,5,6,7,8,9$ sont respectivement $1$, $1$, $2$, $4$, $7$, $10$, $14$, $18$, $23,28$.
 Soit $\varphi_4$ la fonction profil de $Age(G_4)$. Nous avons le r\'esultat suivant:

\begin{lemma}\label{lem:lemmeG4}
Le profil $\varphi_4$ est quasi-polynomial, il est donn\'e par
$$\left\{\begin{array}{ll}
\varphi_4(0)=\varphi_4(1)=1, \varphi_4(2)=2,&\\
\varphi_4(n)=\dfrac{n^2+4n-5}{4}& \text{si }n\text{ est impair }n\geq 3\\
\varphi_4(n)=\dfrac{n^2+4n-4}{4}& \text{si }n\text{ est pair }n\geq 4%
\end{array}\right.$$
En outre $\varphi_4\simeq\dfrac{n^2}{4}$.
\end{lemma}

\begin{proof}
D'apr\`es ce qui pr\'ec\`ede, la valeur de $\varphi_4(n)$ pour tout entier $n\geq 3$, s'obtient en \'enum\'erant tous les triplets d'entiers $(p,k,l)$ v\'erifiant $p\leq min(k,l)$ et $k+l=n$ qui sont diff\'erents de $(0,1,n-1)$ et $(0,2,n-2)$. Nous obtenons alors:
$$\left\{\begin{array}{ll}
\varphi_4(0)=\varphi_4(1)=1, \varphi_4(2)=2,&\\
\varphi_4(n)=2\underset{k=0}{\overset{\frac{n-1}{2}}\sum}(k+1)-2=\dfrac{1}{4}(n-1)(n+5)& n\text{ impair }n\geq 3\\
\varphi_4(n)=2\underset{k=0}{\overset{\frac{n-2}{2}}\sum}(k+1)+(\dfrac{n}{2}+1)-2=\dfrac{1}{4}n(n+4)-1& n\text{ pair }n\geq 4%
\end{array}\right.$$
%La formule donn\'ee s'obtient apr\`es calcul. %on obtient la formule donn\'ee.
Il est clair que le profil de $\varphi_4$ est quasi-polynomial.
\end{proof}

\begin{proposition}
La fonction g\'en\'eratrice de $Age(G_4)$ est rationnelle et est donn\'ee par:
$$F_{G_4}(x)= \dfrac{1-x+2x^3-x^5}{(1-x)^3(1+x)}.$$
\end{proposition}

\begin{proof} La fonction s'obtient en utilisant le Lemme \ref{lem:lemmeG4} et en exprimant les profils par une r\'ecurrence. En effet nous avons:
$$\left\{\begin{array}{ll}
\varphi_4(0)=\varphi_4(1)=1, \varphi_4(2)=2,&\\
\varphi_4(n)=\varphi_4(n-1)+\frac{1}{2}(n+1)& \text{ si }n\text{ est impair }n\geq 3\\
\varphi_4(n)=\varphi_4(n-1)+\frac{1}{2}(n+2)& \text{ si }n\text{ est pair }n\geq 4%
\end{array}\right.$$
La fonction g\'en\'eratrice est alors donn\'ee par
$$\begin{array}{rl}
A=F_{G_4}(x)=\underset{n\geq 0}\sum\varphi_4(n)x^n=& \varphi_4(0)+\varphi_4(1)x+\varphi_4(2)x^2
        +\underset{\small{\begin{array}{c}
        n\geq 3\\
        n \text{\small{ impair}}\end{array}}}\sum\varphi_4(n)x^n+\underset{\small{\begin{array}{c}
    n\geq 4\\
    n \text{\small{ pair}}\end{array}}}\sum\varphi_4(n)x^n.\\
    =& \varphi_4(0)+\varphi_4(1)x+\varphi_4(2)x^2+\underset{\small{\begin{array}{c}
        n\geq 3\\
        n \text{\small{ impair}}\end{array}}}\sum(\varphi_4(n-1)+\dfrac{1}{2}(n+1))x^n\\
        &\qquad +\underset{\small{\begin{array}{c}
    n\geq 4\\
    n \text{\small{ pair}}\end{array}}}\sum(\varphi_4(n-1)+\dfrac{1}{2}(n+2))x^n.\\
   = &  \varphi_4(0)+\varphi_4(1)x+\varphi_4(2)x^2+\underset{n\geq 3}\sum\varphi_4(n-1)x^n+\dfrac{1}{2}\underset{n\geq 3}\sum (n+1)x^n\\
    & \qquad +\dfrac{1}{2}\underset{\small{\begin{array}{c}
        n\geq 4\\
        n \text{\small{ pair}}\end{array}}}\sum x^n.\\
    \end{array}$$
    Comme $$\underset{n\geq 3}\sum\varphi_4(n-1)x^n=\underset{n\geq 0}\sum\varphi_4(n)x^{n+1}-\varphi_4(0)x-\varphi_4(1)x^2.$$
    Alors:
 $$\begin{array}{rl}
A=F_{G_4}(x)= & \varphi_4(0)+\varphi_4(1)x+\varphi_4(2)x^2+\underset{n\geq 0}\sum\varphi_4(n)x^{n+1}-\varphi_4(0)x-\varphi_4(1)x^2\\
&\qquad  +\dfrac{1}{2}\underset{n\geq 3}\sum (n+1)x^n+\dfrac{1}{2}\underset{n\geq 2}\sum x^{2n}.\end{array}$$
Nous avons $$\underset{n\geq 3}\sum (n+1)x^n=x\underset{n\geq 3}\sum nx^{n-1}+\underset{n\geq 3}\sum x^n=\dfrac{-2x^3+3x^4}{(1-x)^2}$$
et
$$\underset{n\geq 2}\sum x^{2n}=x^4\underset{n\geq 0}\sum x^{2n}=\dfrac{x^4}{1-x^2}.$$
D'o\`u:
$$A=xA+\dfrac{-2x^3+3x^4}{2(1-x)^2}+\dfrac{x^4}{2(1-x^2)}+1+x^2.$$
Ce qui donne
$$F_{G_4}(x)=\dfrac{1-x+2x^3-x^5}{(1-x)^3(1+x)}.$$
\end{proof}

%*****************************************************************************
%$$$$$$$$$$$$$$$$$$$$$$$$$$$$$$$$$$$$$$$$$$$$$$$$$$$$$$$$$$$$$$$$$$$$$$$$$$$$
%$$$$$$$$$$$$$$$$$$$$$$$$$$$$$$$$$$$$$$$$$$$$$$$$$$$$$$$$$$$$$$$$$$$$$$$$$$$
%\bigskip

\section{Le graphe  $G_5$}\label{subsec:graphe $G_5$}
Les graphes $G_5$ et $G'_5$ ont m\^eme profil, nous \'etudions alors un seul d'eutre eux, ce sera le graphe $G_5$.

Le graphe $G_5:=(V_5,E_5)$ est le graphe d\'efini par $V_5:=(\mathbb N \times \{0,1\})\cup\{c\}$ tel que $V_5\setminus\{c\}$ se d\'ecompose en deux sous-ensembles, $A=\mathbb N \times \{0\}$ qui forme une clique et $B=\mathbb N \times \{1\}$ qui forme un stable. L'ensemble des ar\^etes contient, en plus des ar\^etes de la clique, les paires $\{(i,0),(j,1)\},~i,j\in\mathbb N$ telles que $i\leq j$ et les paires $\{c,(j,1)\}, ~j\in\mathbb N$.

\vspace{1mm}

Remarquons, comme pour le graphe $G_3$, que pour tout $x\in V_5\setminus\{c\}$, le graphe $G_5$ s'abrite dans ${G_5}_{\restriction_{V_5\setminus\{x\}}}$, donc $x\notin Ker(G_5)$ et $Age(G_5)\neq Age({G_5}_{\restriction_{V_5\setminus\{c\}}})$, donc $Ker(G_5)=\{c\}$.

\vspace{2mm}

Soit $H$ un sous-graphe d'ordre $n\geq 1$ de $G_5$ induit par un $n$-ensemble $V_H$ de sommets.

\begin{remarks}\label{rem:G5}
\begin{enumerate}
\item Si $c\notin V_H$ alors $H$ est connexe si et seulement si
    \begin{itemize}
\item ou bien $V_H\cap(\mathbb N\times\{1\})=\varnothing$ auquel cas $H$ est une clique,
\item ou bien $i_0=min\{i\in\mathbb N/(i,0)\in V_H\}\leq j_0=min\{j\in\mathbb N/(j,1)\in V_H\}$ et dans ce cas $(i_0,0)$ est un sommet de degr\'e maximum et $(j_0,1)$ un sommet de degr\'e minimum.
    \end{itemize}
\item Si $c\in V_H$ alors $H$ est connexe si et seulement si $V_H\cap(\mathbb N\times\{1\})\neq\varnothing$ et si $V_H\cap(\mathbb N\times\{0\})\neq\varnothing$ alors $min\{i\in\mathbb N/(i,0)\in V_H\}\leq max\{j\in\mathbb N/(j,1)\in V_H\}.$
\item Si $H$ n'est pas connexe nous avons
\begin{enumerate}
\item  Si $c\notin V_H$ alors $H$ contient, au plus, une composante connexe de taille au moins deux.
\item Si $c\in V_H$ alors $H$ poss\`ede exactement deux composantes connexes.
\end{enumerate}
\end{enumerate}
\end{remarks}

\begin{lemma}\label{lem:lemmeG5}
Un sous-graphe $H$ d'ordre $n\geq 4$ de $G_5$ s'abrite dans $G_5\setminus\{c\}$ si et seulement s'il ne contient aucun sous-graphe isomorphe \`a $P_4,~2K_2$ ou $C_4$ où $2K_2$ est la somme directe de deux copies de $K_2$ et $C_4$ et le cycle à quatre sommets.
\end{lemma}
\begin{proof}
Il est facile de v\'erifier que ces trois graphes ne s'abritent pas dans $G_5\setminus\{c\}$, d'o\`u la condition n\'ecessaire. Pour la condition suffisante, montrons que si $H$ ne s'abrite pas dans $G_5\setminus\{c\}$, il contient n\'ecessairement un des graphes cit\'es. Si $H$ ne s'abrite pas dans $G_5\setminus\{c\}$, donc $H$ contient le sommet $c$. Nous avons deux cas:
\medskip

\quad\textbf{Cas 1:} Si $H$ est connexe alors, d'apr\`es la Remarque \ref{rem:G5}, nous avons $V_H\cap(\mathbb N\times\{1\})\neq\varnothing$. Comme $H$ ne s'abrite pas dans $G_5\setminus\{c\}$, alors n\'ecessairement $V_H\cap(\mathbb N\times\{0\})\neq\varnothing,$ donc $i_0=min\{i\in\mathbb N/(i,0)\in V_H\}\leq j_0=max\{j\in\mathbb N/(j,1)\in V_H\}.$ Posons $x_{i_0}=(i_0,0)$ et $x_{j_0}=(j_0,1)$. Mais $n\geq 4$ donc il existe au moins un autre sommet $x$ dans $V_H$.  Ce sommet $x$ est
\begin{itemize}
\item soit de la forme $(i_x,0)$ avec $i_x> i_0$ et comme $H$ ne s'abrite pas dans $G_5\setminus\{c\}$, nous avons n\'ecessairement $i_x> j_0$  et dans ce cas les sommets $\{c,x_{i_0},x_{j_0},x\}$ forment un $P_4$. %si $i_x> j_0$.
\item soit de la forme $(j_x,1)$ avec $j_x< j_0$ et dans ce cas les sommets $\{c,x_{i_0},x_{j_0},x\}$ forment un $P_4$ si $j_x< i_0$ ou un $C_4$ sinon.
\end{itemize}
\medskip

\quad\textbf{Cas 2:} Si $H$ n'est pas connexe alors, toujours d'apr\`es la Remarque \ref{rem:G5}, $H$ poss\`ede exactement deux composantes connexes, l'une donn\'ee par $\{c\}\cup V_H\cap(\mathbb N\times\{1\})$ et l'autre par $V_H\cap(\mathbb N\times\{0\})$. Comme $H$ ne s'abrite pas dans $G_5\setminus\{c\}$, ces deux composantes ont, chacune, au moins deux sommets, d'o\`u le sous-graphe $2K_2$.

\end{proof}

\subsection{Repr\'esentation des sous-graphes de $G_5$}
Nous associons \`a $H$ la repr\'esentation suivante.
Comme pour les sous-graphes de $G_1$, nous plaçons les sommets de $V_H\setminus\{c\}$ sur deux lignes horizontales $L_0$ et $L_1$, telles que les sommets de $V_H\cap(\mathbb N\times\{0\})$ sont sur $L_0$ et ceux de $V_H\cap(\mathbb N\times\{1\})$ sur $L_1$ avec les sommets de la forme $(j,1)$ qui sont \`a droite de tous les sommets de la forme $(i,0)$ pour $i\leq j$ et \`a gauche de tous les sommets de la forme $(i,0)$ pour $i>j$. Un sommet qui est sur $L_0$ est reli\'e \`a tout sommet qui se trouve sur $L_0$ et \`a tout sommet qui se trouve sur $L_1$ mais \`a sa droite.  Un sommet sur $L_1$ est reli\'e \`a tout sommet qui se trouve sur $L_0$ mais \`a sa gauche. Si $c\in H$ nous le repr\'esentons par un point, au dessus de la ligne $L_1$, qui sera reli\'e \`a tous les sommets qui sont sur $L_1$.
\medskip

  Nous repr\'esentons $H$ par une suite de $\{0,1\}^{n-1}\times\{0,1,2\}$ obtenue en balayant les lignes $L_0$, $L_1$, de gauche \`a droite, par une ligne verticale et en associant \`a tout point $v$, rencontr\'e dans cet ordre, $0$ si $v$ se trouve sur $L_0$ et $1$ s'il est sur $L_1$. Cette suite obtenue est de longueur $n$ si $c\notin V_H$ et $n-1$ sinon, dans ce dernier cas nous rajoutons \`a la suite un $n^{\text{i\`eme}}$ terme de valeur $2$. Par exemple le sous-graphe $H=(V,E)$ avec $V=(\{0,2,3,4\}\times\{0\})\cup(\{0,1,4\}\times\{1\})$ est repr\'esent\'e par la suite $S=0110001$ et le sous-graphe $H'=(V',E')$ o\`u $V'=V\cup\{c\}$ est repr\'esent\'e par $S'=01100012$ (voir \figurename~\ref{exempleG5}).\\

Donc, \`a tout sous-graphe $H$ d'ordre $n$ de $G_5$ nous associons une suite $S_H=u_0u_1\cdots u_{n-1}$ telle que $u_i\in\{0,1\}$ pour tout $i,~0\leq i\leq n-2$ et $u_{n-1}\in\{0,1,2\}$.

Inversement, nous associons \`a la suite vide le sous-graphe vide et \`a toute suite $S=u_0u_1\cdots u_{n-1}$ pour $n\in\mathbb N^{\star}$ le type d'isomorphie de ${G_5}_{\restriction_X}$ avec $X=\{i\in\{0,\cdots, n-1\}/u_i=0\}\times\{0\}\cup\{i\in\{0,\cdots, n-1\}/u_i=1\}\times\{1\}$ si $u_{n-1}\neq 2$ et $X=\{i\in\{0,\cdots, n-1\}/u_i=0\}\times\{0\}\cup\{i\in\{0,\cdots, n-1\}/u_i=1\}\times\{1\}\cup\{c\}$ si $u_{n-1}=2$.
\begin{figure}
\centering
\input{imag5a}
\caption{\label{exempleG5}Exemple de repr\'esentation de sous-graphes de $G_5$.}
\end{figure}

\begin{example}\label{ex:G5}
\begin{enumerate}
\item Le sous-graphe d'ordre $1$ est repr\'esent\'e par les deux suites $0$ et $1$.
\item Les sous-graphes d'ordre $2$ sont repr\'esent\'es par $00$ et $01$ pour la clique et par $10$ et $11$ pour l'ind\'ependant.
\item Les cliques d'ordres $n\geq 3$ sont repr\'esent\'ees par $\underset{n}{\underbrace{00\cdots 0}}$ et $\underset{n-1}{\underbrace{00\cdots 0}} 1.$
\item Les ind\'ependants d'ordres $n\geq 3$ sont repr\'esent\'es par $\underset{n-1}{\underbrace{11\cdots 1}} 0$ et $\underset{n}{\underbrace{11\cdots 1}} .$
    \end{enumerate}
\end{example}

\begin{remark}\label{rem:remG5}
Les suites de $\{0,1\}^n$, pour tout $n\geq 1$ qui ne diff\`erent que du $n^{\text{i\`eme}}$ terme ($S=u_0u_1\cdots u_{n-2}0$ et $S'=u_0u_1\cdots u_{n-2}1$) sont associ\'ees \`a des sous-graphes isomorphes.
\end{remark}

Par exemple le sous-graphe $H$ donn\'e dans la \figurename~\ref{exempleG5} est isomorphe au sous-graphe $K=(U,F)$ avec $U=(\{0,2,3,4,5\}\times\{0\})\cup(\{0,1\}\times\{1\})$ qui est repr\'esent\'e par la suite $S''=0110000.$

\begin{lemma}\label{lem:sanscG5}
Soient $H$ et $H'$ deux sous-graphes de $G_5$ de m\^eme ordre $n\geq 1$ et soient $S=u_0u_1\cdots u_{n-1}$ et $S'=u'_0u'_1\cdots u'_{n-1}$ les suites associ\'ees respectivement. Si $c\notin H\cup H'$ (autrement dit, $u_{n-1}\neq 2$ et $u'_{n-1}\neq 2$) alors $$H \text{ et }H'\text{ sont isomorphes }\Leftrightarrow u_0u_1\cdots u_{n-2}=u'_0u'_1\cdots u'_{n-2}.$$
\end{lemma}

\begin{proof}
La condition suffisante d\'ecoule de la Remarque \ref{rem:remG5}. Pour la condition n\'ecessaire, supposons que $H$ et $H'$ soient isomorphes. Faisons un raisonnement par r\'ecurrence sur $n$.
\begin{itemize}
\item Si $n=1$, nous avons deux suites associ\'ees $0$ ou $1$ donc le r\'esultat est vrai.
\item Supposons le r\'esultat vrai \`a l'ordre $n\geq 2$ et soient $H$ et $H'$ deux sous-graphes isomorphes d'ordre $n+1$ ne contenant pas le sommet $c$ et $f$ un isomorphisme de $H$ sur $H'$. Nous avons les cas suivants:
   % \begin{enumerate}

\vspace{1mm}

   \quad \textbf{Cas 1:}
   Si $H$ est connexe, d'apr\`es la Remarque \ref{rem:G5}, $H$, et donc $H'$, est ou bien une clique dans ce cas les seules suites possibles sont $\underset{n}{\underbrace{00\cdots 0}}0$ ou  $\underset{n}{\underbrace{00\cdots 0}}1$ et le r\'esultat est v\'erifi\'e, ou bien $H$ n'est pas une clique et donc les degr\'es de ses sommets ne sont pas tous \'egaux. Soit alors $x_0$ un sommet de $H$ de degr\'e maximum. Consid\'erons les sous-graphes $K=H\setminus\{x_0\}$ et $K'=H'\setminus\{f(x_0)\}$. $K$ et $K'$ sont deux sous-graphes isomorphes d'ordre $n$, d'apr\`es l'hypoth\`ese de r\'ecurrence, les suites qui leurs sont associ\'es ne peuvent diff\'erer que du dernier terme consid\'erons les repr\'esentations de $K$ est $K'$ sur $L_0,~L_1$.  Les sommets $x_0$ et $f(x_0)$ \'etant de degr\'e maximum, dans les repr\'esentations de $H$ et $H'$, ils ne peuvent \^etre plac\'es que sur $L_0$ \`a gauche de tous les autres sommets (car $H$ est connexe). Ceci se traduit, dans les suites par un $0$ \`a la premi\`ere position. D'o\`u le r\'esultat.

\vspace{1mm}

      \quad\textbf{Cas 2:} Si $H$ n'est pas connexe, alors d'apr\`es la Remarque \ref{rem:G5},  $H$ (et donc $H'$), contient, au plus, une composante connexe d'ordre $k\geq 2$. Donc ou bien $H$ (et donc $H'$), est un ind\'ependant alors le r\'esultat est vrai (voir Exemple \ref{ex:G5}), ou bien $H$ n'est pas un ind\'ependant mais poss\`ede donc, au moins un sommet isol\'e $x_1$. On applique alors l'hypoth\`ese de r\'ecurrence aux sous-graphes $H\setminus\{x_1\}$ et $H'\setminus\{f(x_1)\}$ et la seule position possible pour $x_1$ est sur $L_1$ tout \`a fait \`a gauche. Le r\'esultat est alors v\'erifi\'e.
   % \end{enumerate}
\end{itemize}
\end{proof}
%\medskip

\begin{remark}
Les sous-graphes donn\'es par le Lemme \ref{lem:sanscG5}, ne sont plus isomorphes si on leurs rajoute le sommet $c$.
 \end{remark}
 En effet, les sous-graphes repr\'esent\'es par les suites $u_0u_1\cdots u_{n-3}02$ et $u_0u_1\cdots u_{n-3}12$ pour $n\geq 3$ ne sont pas isomorphes, car le sommet $c$ n'a pas le m\^eme degr\'e dans les deux sous-graphes. Comme exemple  reprenons les sous-graphes isomorphes $H$ et $K$ donn\'es dans la Remarque \ref{rem:remG5} et consid\'erons les sous-graphes $H'$ et $K'$ obtenus en rajoutant le sommet $c$ \`a $H$ et $K$ respectivement. $H'$ est le sous-graphe de la \figurename~\ref{exempleG5} et $K'$ est le m\^eme mais sans l'ar\^ete $\{x,c\}$.
\bigskip

\begin{figure}
\centering
\input{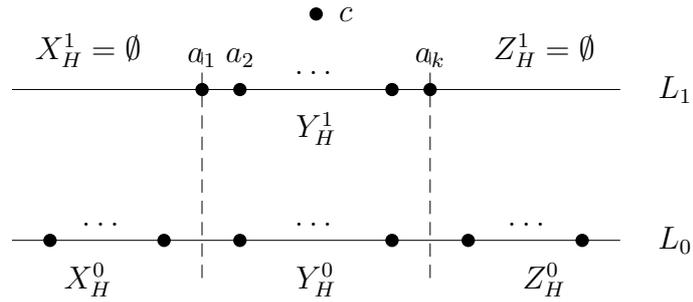}
\caption{\label{schemaG5}R\'epartition des sommets d'un sous-graphe $H$.}
\end{figure}

Soit $H=(V_H,E_H)$ un sous-graphe de $G_5$ contenant le sommet $c$. Consid\'erons la repr\'esentation de $H\setminus\{c\}$ sur $L_0,~L_1$, soient $X^i_H,~Y^i_H,~Z^i_H$ avec $i\in\{0,1\}$ d\'efinis par (voir \figurename~\ref{schemaG5}):

$Y^1_H=\{a_1,\cdots,a_k\}$ est l'ensemble de tous les sommets de $H$ qui sont sur $L_1$, donc $V_H\cap X^1_H=V_H\cap Z^1_H=\varnothing$,

$X^0_H$ est l'ensemble de tous les sommets qui sont sur $L_0$ et \`a gauche de $a_1$,

$Y^0_H$ est  l'ensemble des sommets sur $L_0$ qui sont \`a droite de $a_1$ et \`a gauche de $a_k$ et

$Z^0_H$ les sommets de $L_0$ qui sont \`a droite de $a_k$.

\begin{lemma}\label{lem:lemm1G5}
Soit $H$ un sous-graphe de $G_5$ d'ordre $n\geq 4$ contenant le sommet $c$. Le sous-graphe $H$ s'abrite dans $G_5\setminus\{c\}$ si et seulement si nous avons l'une des situations suivantes:
\begin{enumerate}
\item[(1)] $Y^1_H=\varnothing$.
\item[(2)] $\vert Y^1_H\vert=1$ (donc $\vert Y^0_H\vert=\varnothing$), $X^0_H\neq\varnothing$ et $Z^0_H=\varnothing.$
\item[(3)] $\vert Y^1_H\vert >1,~X^0_H=Y^0_H=\varnothing$ et $\vert Z^0_H\vert\leq 1.$
\end{enumerate}
\end{lemma}

\begin{proof}
La condition suffisante se v\'erifie de fa\c{c}on  imm\'ediate. Pour la condition n\'ecessaire elle utilise le Lemme \ref{lem:lemmeG5}.  Supposons que $H$ s'abrite dans $G_5\setminus\{c\}$ et soit la repr\'esentation de $H\setminus\{c\}$ sur $L_0,~L_1$ alors
\begin{enumerate}
\item Si la repr\'esentation n'utilise aucun point de $L_1$, donc $H$ est de la forme $K_{n-1}\oplus 1$,  s'abrite dans $G_5\setminus\{c\}$ et nous sommes  dans le cas $(1)$.
\item Si nous avons un seul point sur $L_1$, alors comme $H$ a au moins quatre sommets, nous avons n\'ecessairement $Z^0_H=\varnothing$, car sinon $H$ abriterait $P_4$ si $\vert Z^0_H\vert= 1$ ou $2K_2$ si $\vert Z^0_H\vert> 1$ et nous sommes alors dans le cas $(2)$ du lemme.
\item Si $\vert Y^1_H\vert >1$ alors n\'ecessairement $X^0_H=\varnothing$, car sinon $H$ abriterait $C_4$ et $Y^0_H=\varnothing$ car sinon $H$ abriterait $P_4$ et $\vert Z^0_H\vert\leq 1$ car sinon $H$ abriterait $2K_2$. Nous sommes dans le cas $(3)$ du lemme.
\end{enumerate}
\end{proof}

 \begin{corollary}\label{cor:corG5}
 Soient $H$ et $H'$ deux sous-graphes de $G_5$ de m\^eme ordre $n\geq 4$  et soient $S=u_0u_1\cdots u_{n-1}$ et $S'=u'_0u'_1\cdots u'_{n-1}$ les suites associ\'ees respectivement. Si $c\in H\setminus H'$ (autrement dit, $u_{n-1}= 2$ et $u'_{n-1}\neq 2$) alors $H$ et $H'$ sont isomorphes si et seulement si nous avons l'une des situations suivantes:
 \begin{enumerate}
 \item $S=\underset{n-1}{\underbrace{00\cdots 0}}2$ et $S'=1\underset{n-1}{\underbrace{00\cdots 0}}.$
 \item $S=\underset{n-2}{\underbrace{00\cdots 0}}12$ et $S'=01\underset{n-2}{\underbrace{00\cdots 0}}.$
 \item $S=\underset{n-1}{\underbrace{11\cdots 1}}2$ et $S'=0\underset{n-1}{\underbrace{11\cdots 1}}.$
 \item $S=\underset{n-2}{\underbrace{11\cdots 1}}02$ et $S'=10\underset{n-2}{\underbrace{11\cdots 1}}.$
 \end{enumerate}
 \end{corollary}

 \begin{proof}
 La condition suffisante se v\'erifie facilement. Pour la condition n\'ecessaire si $H$ est isomorphe \`a $H'$ nous avons l'un des cas $(1),~(2)$ ou $(3)$ du Lemme \ref{lem:lemm1G5} appliqu\'ees à $H$.
 \begin{enumerate}
 \item Dans le cas $(1)$, $S=\underset{n-1}{\underbrace{00\cdots 0}}2$ et $H$ est une clique avec un sommet isol\'e qui est aussi repr\'esent\'e par $1\underset{n-1}{\underbrace{00\cdots 0}}$.
 \item Dans le cas $(2)$, $S=\underset{n-2}{\underbrace{00\cdots 0}}12$ et $H$ est une clique avec un sommet pendant qui est \'egalement repr\'esent\'e par $01\underset{n-2}{\underbrace{00\cdots 0}}$.
 \item Dans le cas $(3)$, nous avons deux possibilit\'es, ou bien $S=\underset{n-1}{\underbrace{11\cdots 1}}2$ et dans ce cas $H$ est une \'etoile $S_{1,n-1}$ qui peut \'egalement \^etre repr\'esent\'ee par $0\underset{n-1}{\underbrace{11\cdots 1}}$, ou bien $S=\underset{n-2}{\underbrace{11\cdots 1}}02$ et dans ce cas $H$ est une \'etoile $S_{1,n-2}$ avec un sommet isol\'e qui est aussi repr\'esent\'e par $10\underset{n-2}{\underbrace{11\cdots 1}}$.
 \end{enumerate}
 \end{proof}

\begin{remarks}\label{rem:isoG5}
\begin{enumerate}
\item[(1)] Les sous-graphes $H$ et $H'$ d'ordre $n\geq 4$ contenant $c$ repr\'esent\'es par les suites $01\underset{n-3}{\underbrace{00\cdots 0}}2$ et $1\underset{n-3}{\underbrace{00\cdots 0}}12$ sont isomorphes mais $H\setminus\{c\}$ et $H'\setminus\{c\}$ ne le sont pas.
     \begin{center}
\psset{unit=1cm}
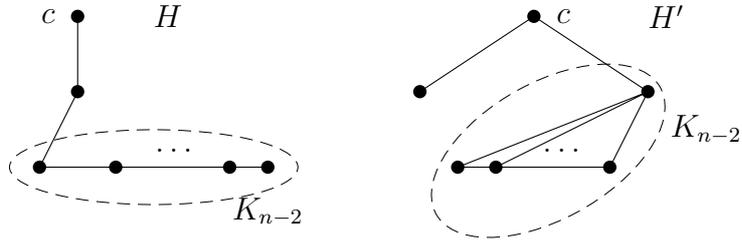
\begin{pspicture}(-8,-2)(8,1.5)
\psdots[dotsize=5pt](-4,-1)(-3,-1)(-1.5,-1)(-1,-1)(-3.5,0)(-3.5,1)
\psline[linewidth=0.4pt](-3.5,1)(-3.5,0)(-4,-1)(-3,-1)(-1.5,-1)(-1,-1)
\uput{0.3}[l](-3.5,1){$c$}
\uput{1}[r](-3.5,1){$H$}
\uput{0.2}[u](-2.2,-1){$\dots$}
\psellipse[linestyle=dashed,linewidth=0.2pt](-2.5,-1)(1.9,0.5)
\uput{0.4}[d](-1,-1){$K_{n-2}$}
\psdots[dotsize=5pt](1,0)(1.5,-1)(2,-1)(3.5,-1)(4,0)(2.5,1)
\psline[linewidth=0.4pt](1,0)(2.5,1)(4,0)(1.5,-1)(2,-1)(3.5,-1)(4,0)(2,-1)
\uput{0.3}[r](2.5,1){$c$}
\uput{1.5}[r](2.5,1){$H'$}
\uput{0.2}[u](2.9,-1){$\dots$}
\rput{30}(0.1,-1.7){\psellipse[linestyle=dashed,linewidth=0.2pt](2.7,-0.5)(1.7,0.9)}
\uput{0.2}[r](4.1,-0.5){$K_{n-2}$}
\end{pspicture}
\end{center}
\item[(2)] Les sous-graphes $H$ et $H'$ d'ordre $n\geq 5$ contenant $c$ repr\'esent\'es par les suites $0\underset{n-3}{\underbrace{11\cdots 1}}02$ et $10\underset{n-3}{\underbrace{11\cdots 1}}2$ sont isomorphes mais $H\setminus\{c\}$ et $H'\setminus\{c\}$ ne le sont pas.
  \begin{center}
\psset{unit=1cm}
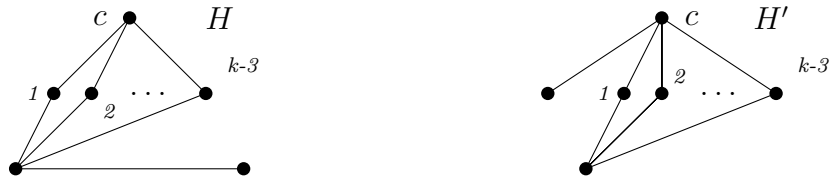
\begin{pspicture}(-8,-1.5)(8,1)
\psdots[dotsize=5pt](-5,-1)(-4.5,0)(-4,0)(-2.5,0)(-2,-1)(-3.5,1)
\psline[linewidth=0.4pt](-5,-1)(-2.5,0)(-3.5,1)(-4,0)(-5,-1)(-2,-1)
\psline[linewidth=0.4pt](-3.5,1)(-4.5,0)(-5,-1)
\uput{0.3}[l](-3.5,1){$c$}
\uput{1}[r](-3.5,1){$H$}
\uput{0.2}[l](-4.5,0){\scriptsize{1}}
\uput{0.2}[dr](-4,0){\scriptsize{2}}
\uput{0.4}[ur](-2.5,0){\scriptsize{k-3}}
\uput{0.5}[r](-4,0){$\dots$}
\psdots[dotsize=5pt](2,0)(2.5,-1)(3,0)(3.5,0)(5,0)(3.5,1)
\psline[linewidth=0.4pt](2,0)(3.5,1)(3.5,0)(2.5,-1)(3.5,0)(3.5,1)(5,0)(2.5,-1)
\psline[linewidth=0.4pt](3.5,1)(3,0)(2.5,-1)
\uput{0.3}[r](3.5,1){$c$}
\uput{1.2}[r](3.5,1){$H'$}
\uput{0.2}[l](3,0){\scriptsize{1}}
\uput{0.2}[ur](3.5,0){\scriptsize{2}}
\uput{0.4}[ur](5,0){\scriptsize{k-3}}
\uput{0.5}[r](3.5,0){$\dots$}
\end{pspicture}
\end{center}
\end{enumerate}
Tous les sous-graphes cit\'es dans $(1)$ et $(2)$ contiennent $P_4, ~2K_2$ ou $C_4$ et aucun des isomorphismes transformant $H$ en $H'$  ne laisse $c$ invariant.
\end{remarks}

\begin{lemma}\label{lem:lemme2G5}
Deux sous-graphes $H$ et $H'$ de m\^eme ordre $n\geq 4$ contenant le sommet $c$ et ne s'abritant pas dans $G_5\setminus\{c\}$ sont isomorphes avec $H\setminus\{c\}$ et $H'\setminus\{c\}$ non isomorphes si et seulement si $H$ et $H'$ sont donn\'es par l'un des points $(1)$ ou $(2)$ de la Remarque \ref{rem:isoG5}.
\end{lemma}

\begin{proof}
La condition suffisante est donn\'ee par la Remarque \ref{rem:isoG5}. Pour la condition n\'ecessaire, soient $H$ et $H'$ deux sous-graphes de m\^eme ordre $n\geq 4$, v\'erifiant les conditions du Lemme \ref{lem:lemme2G5} et soit $f$ un isomorphisme de $H$ sur $H'$. Le sommet $c$ ne peut pas \^etre invariant par $f$ car sinon $H\setminus\{c\}$ et $H'\setminus\{c\}$ seraient isomorphes, ce qui contrediraient l'hypoth\`ese. Donc $f(c)\neq c.$ Posons $c'=f(c)$. Consid\'erons les repr\'esentations de $H\setminus\{c\}$ et $H'\setminus\{c\}$ sur $L_0,~L_1$.
Comme $H$ et $H'$ ne s'abritent pas dans $G_5\setminus\{c\}$, d'apr\`es le Lemme \ref{lem:lemm1G5} nous avons $Y^1_H\neq\varnothing$ et $Y^1_{H'}\neq\varnothing$,
c'est \`a dire $d_H(c)\neq 0$ et $d_{H'}(c)\neq 0$. Nous avons les cas suivants qui sont repr\'esent\'es sur les tables \tablename~\ref{tablescas1-G5} et \tablename~\ref{tablescas2-G5}:
\begin{enumerate}
\item Si $d_H(c)=1$, c'est \`a dire $\vert Y^1_H\vert=1$, donc $d_{H'}(c')=1$ (les voisins de $c$ dans $H$ sont envoy\'es par $f$ sur les voisins de $c'$ dans $H'$). D'apr\`es le Lemme \ref{lem:lemm1G5} on a $X^0_H=\varnothing$ ou $Z^0_H\neq\varnothing$,   alors
    \begin{enumerate}
    \item Si $c'$ se trouve sur $L_1$ (donc $c'$ et $c$ sont reli\'es dans $H'$) alors $X^0_{H'}=\varnothing$ (car sinon on aurait $d_{H'}(c')>1$). Nous avons alors
        \begin{enumerate}
        \item Si $X^0_H=\varnothing$ donc $H$ n'est pas connexe et $H'$ \'egalement. Dans ce cas $H\setminus\{c\}$ et $H'\setminus\{c\}$ sont isomorphes ce qui contredit l'hypoth\`ese.
        \item Si $X^0_H\neq\varnothing$ alors $Z^0_H\neq\varnothing$, nous avons forc\'ement $\vert X^0_H\vert=1$, car sinon, $H$ et $H'$ ne seraient pas isomorphes, puisque l'unique voisin $x$ de $c$ dans $H$ est envoy\'e, par $f$ sur $c$ dans $H'$. Donc les voisins de $x$ dans $H$ (qui sont tous sur $L_0$ et forment donc une clique dans $H$) sont envoy\'es sur les voisins de $c$ dans $H'$ (qui sont tous sur $L_1$ et forment donc un stable).  Soit alors $y\neq c$ le voisin de $x$ dans $H$ et $y'=f(y)$ ($y\in X^0_H$ et $y'\in Y^1_{H'}$), tous les voisins de $y$, diff\'erents de $x$, sont dans $Z^0_H$ qui est non vide. Ceci se traduit, dans $H'$, par autant de points dans $Y^0_{H'}$. Ceci donne la situation $(1)$ de la Remarque \ref{rem:isoG5}.
        \end{enumerate}
    \item Si $c'$ se trouve sur $L_0$, soit $x'$ l'unique voisin de $c'$ dans $H'$ (car $d_{H'}(c')=1$), donc $x'=f(x)$ on a alors
            \begin{enumerate}
            \item Si $x'$ se trouve sur $L_0$ alors $x'$ est plac\'e \`a gauche de $c'$ et tous ses voisins, diff\'erents de $c'$, sont sur $L_1$ et \`a gauche de $c'$ (sinon ils seraient des voisins de $c'$ qui est de degr\'e $1$). Donc $x'$ ne peut-avoir plus d'un voisin sur $L_1$ car ces voisins (qui forment un stable) sont les images, par $f$, des voisins de $x$ dans $H$ qui eux se trouvent sur $L_0$ et forment une clique. Soit alors $y\neq c$ le voisin de $x$ dans $H$ et $y'=f(y)$. Comme $y'\in L_1$, donc $y'$ est un voisin de $c$ dans $H'$, ce qui donne dans $H$ un voisin de $y$ qui se trouve dans $Z^0_H$. Dans ce cas $H\setminus\{c\}$ et $H'\setminus\{c\}$ sont isomorphes ce qui contredit l'hypoth\`ese.
            \item Si $x'$ se trouve sur $L_1$ donc $c'$ est l'unique sommet de $H'$ qui se trouve sur $L_0$ (car sinon on aurait $d_{H'}(c')>1$). Il s'ensuit que $H$ et $H'$ sont d'ordre $3$ et par cons\'equent s'abritent dans $G_5\setminus\{c\}$, ce qui contredit l'hypoth\`ese.
            \end{enumerate}
    \end{enumerate}
\item Si $d_H(c)=p>1$, soit $\{x_1,\cdots, x_p\}$ l'ensemble des voisins de $c$ dans $H$, plac\'es sur $L_1$ dans l'ordre croissant de leurs indices. Soit $\{x'_1,\cdots, x'_p\}$ l'ensemble des voisins de $c'$ dans $H'$ avec $x'_i=f(x_i)$ pour $1\leq i\leq p.$ Alors nous avons les cas suivant:
    \begin{enumerate}
    \item Si $c'\in L_1$ alors forc\'ement $p=2$ car sinon on aurait un stable de $H$ qui serait envoy\'e sur une clique de $H'$. Alors, si $x_1$ ne poss\`ede pas d'autres voisins que $c$, nous obtenons le cas $(1)$ de la Remarque \ref{rem:isoG5} et si $x_1$ poss\`ede un autre voisin  alors $H$ et $H'$ sont tels que $H\setminus\{c\}$ et $H'\setminus\{c\}$ sont isomorphes ce qui contredit l'hypoth\`ese.
    \item Si $c'\in L_0$ alors $c'$ a, au plus, un voisin sur $L_0$ et il serait \`a droite de tous ses autres voisins sur $L_1$ (car les voisins de $c'$ forment un stable).
        \begin{enumerate}
        \item Si $c'$ a un voisin sur $L_0$, ce sera $x'_1$, alors les autres voisins $x'_2,\cdots, x'_p$ sont sur $L_1$ (\`a droite de $c'$ et \`a gauche de $x'_1$). Comme les sommets $x'_2,\cdots, x'_p$ sont reli\'es \`a $c$ dans $H'$, on en d\'eduit que, dans $H$, $\vert Y^0_H\vert=1$ et ce point de $Y^0_H$ est \`a gauche des sommets $x_2,\cdots, x_p$ et à droite de $x_1$. Il ne peut y avoir d'autres sommets. Ceci engendre la situation $(2)$ de la Remarque \ref{rem:isoG5}.
        \item $c'$ n'a pas de voisins dans $L_0$ ($c'$ est donc l'unique sommet sur $L_0$). Comme les voisins $x'_1,x'_2,\cdots, x'_p$ de $c'$ sont reli\'es \`a $c$ dans $H'$, donc, dans $H$, $\vert X^0_H\vert=1$. Soit $y$ le sommet de $X^0_H$, donc $f(y)=c$. Le sommet $y$ poss\`ede, au plus, un voisin sur $L_0$.
                \begin{itemize}
                \item Si $y$ a un voisin $z$ sur $L_0$ alors $z\in Z^0_H$ et $z'=f(z)$ se trouve sur $L_1$ et \`a gauche de $c'$ et nous obtenons le cas $(2)$ de la Remarque \ref{rem:isoG5}.
                \item Si $y$ n'a pas de voisin sur $L_0$ alors $H\setminus\{c\}$ et $H'\setminus\{c\}$ sont isomorphes ce qui contredit l'hypoth\`ese.
                \end{itemize}
        \end{enumerate}
    \end{enumerate}
\end{enumerate}
\end{proof}

\begin{table}[h]
\small
\begin{tabular}[c]{|c|c|}
\hline
\input{imag7}&\quad \input{imag7a}\\
\hline
Cas $1.(a).ii$& Cas $1.(b).i$\\
\hline
\end{tabular}
\caption{\label{tablescas1-G5}Les cas $1.(a).ii$ et $1.(b).i$ de la preuve du Lemme \ref{lem:lemme2G5}.}
\end{table}

\begin{table}[h]
\small
\begin{tabular}[c]{|c|c|}
\hline
\input{imag8}&\input{imag8a}\\
\hline
Cas $2.(a)$ pour $d(x_1)>1$& Cas $2.(a)$ pour $d(x_1)=1$\\
\hline
\multicolumn{2}{|c|}
{Les deux situations du cas $2.(a)$}\\
\hline
\input{imag9}&\input{imag9a}\\
\hline
Cas $2.(b).i$& Cas $2.(b).ii$\\
\hline
\multicolumn{2}{|c|}
{Les deux situations du cas $2.(b)$}\\
\hline
\end{tabular}
\caption{\label{tablescas2-G5}Les diff\'erents sous-cas du cas $2$ de la preuve du Lemme \ref{lem:lemme2G5}.}
\end{table}

\subsection{Les ind\'ecomposables de $Age(G_5)$}
Comme le sous-graphe $G_5\setminus\{c\}$ est sans $P_4$, donc s\'eries-parall\`eles, les ind\'ecomposables qu'il abrite sont de tailles au plus $2$. Il s'ensuit que tout sous-graphe ind\'ecomposable, d'ordre sup\'erieure \`a quatre, de $G_5$ contient le sommet $c$.

\medskip
Soit $H$ un sous graphe d'ordre $n\geq 4$ de $G_5$ et $S=u_0 u_1\cdots u_{n-1}$ la suite qui lui est associ\'ee.

\begin{fact}
 $H$ est ind\'ecomposable si et seulement si $u_{n-1}=2$ et $u_i\neq u_{i+1},~\forall i, ~0\leq i\leq n-3$ pour $n\geq 4.$
\end{fact}
\begin{proof}
Evident car des termes successifs de m\^eme valeur dans la suite engendrent un intervalle de $H$.
\end{proof}
\bigskip

 Il s'ensuit, d'apr\`es le Lemme \ref{lem:lemm1G5}, la Remarque \ref{rem:isoG5} et le Lemme \ref{lem:lemme2G5}, qu'il existe un ind\'ecomposable d'ordre $n=4$ et deux ind\'ecomposables d'ordre $n$ pour tout $n\geq 5$. Il sont donn\'es dans la \figurename~\ref{indecG5}. Chaque ind\'ecomposable d'ordre $n$ s'abrite dans tout ind\'ecomposable d'ordre sup\'erieure. Les sous-graphes ind\'ecomposables d'ordre $n$ se trouvent au niveau $n-1$, pour $n\geq 4$. La classe des ind\'ecomposables est donc minimale dans $Ind(\Omega_1)$ et $Age(G_5)$ est ind-minimal.% donc minimal.

\begin{figure}[h]
\centering
\input{imag10}
\caption{\label{indecG5}Les sous-graphes ind\'ecomposables de $G_5$.}
\end{figure}

\subsection{Profil de $Age(G_5)$}

\noindent Les premi\`eres valeurs du profil sont respectivement $1, 1, 2, 4, 11, 26, 58, 122$ pour $n=0$, $1$, $2$, $3$, $4$, $5$, $6$, $7$.

\bigskip
Soit $\varphi_5$ la fonction profil de $Age(G_5)$. Nous avons

\begin{lemma}\label{lem:proG5}
Le profil $\varphi_5$ est exponentiel et est donn\'e par
$$\left\{\begin{array}{l}
\varphi_5(0)=\varphi_5(1)=1,~\varphi_5(2)=2,~\varphi_5(3)=4,~\varphi_5(4)=11;\\
%\varphi_5(3)=4,~\varphi_5(4)=11;\\
\varphi_5(n)=2^n-6\qquad\quad \text{ pour }n\geq 5%
\end{array}%
\right.$$
En outre $\varphi_5\simeq 2^n$.
\end{lemma}

\begin{proof}
Le nombre $\varphi_5(n)$ de sous-graphes d'ordre $n$, non isomorphes, de $G_5$,  s'obtient en additionnant le nombre $\varphi_{-c}(n)$ de sous-graphes d'ordre $n$, non isomorphes, de $G_5\setminus\{c\}$ et le nombre $\varphi_c(n)$ de sous-graphes, non isomorphes, d'ordre $n$ contenant $c$ et ne s'abritant pas dans $G_5\setminus\{c\}$. Donc
$$\varphi_5(n)=\varphi_{-c}(n)+\varphi_c(n).$$
 Nous avons, d'apr\`es le Lemme \ref{lem:sanscG5}, $$\varphi_{-c}(n)=2^{n-1}\qquad \forall n\in\mathbb N^{\star}.$$ car ceci revient \`a \'enum\'erer toutes les suites distinctes de $\{0,1\}^{n-1}.$\\

 $\varphi_c(n)$ s'obtient en \'enum\'erant tous les sous-graphes contenant $c$, non isomorphes (voir Lemme \ref{lem:lemme2G5}) qui ne v\'erifient pas le Lemme \ref{lem:lemm1G5} ou toutes les suites de $\{0,1\}^{n-1}\times\{2\}$, qui ne v\'erifient pas le Corollaire \ref{cor:corG5}.\\
 Les sous-graphes d'ordres $n\geq 4$ v\'erifiant le Lemme \ref{lem:lemm1G5} sont au nombre de quatre pour tout $n\geq 4$, un pour chacun des points $(1)$ et $(2)$ du lemme et deux pour le point $(3)$. Les sous-graphes d'ordre $n\geq 4$ v\'erifiant le Lemme \ref{lem:lemme2G5} sont au nombre de un pour $n=4$ et deux pour tout $n\geq 5$, \`a signaler que tous les sous-graphes d'ordres $n\leq 3$ s'abritent dans $G_5\setminus\{c\}$.
 Ceci donne $$\varphi_c(n)=2^{n-1}-5=3,~ \text{ pour } n= 4\quad\text{ et }~~\varphi_c(n)=2^{n-1}-6,~ \forall n\geq 5.$$
 D'o\`u nous d\'eduisons la formule donn\'ee.
\end{proof}

\begin{proposition}
La fonction g\'en\'eratrice de $Age(G_5)$ est rationnelle et est donn\'ee par:
$$F_{G_5}(x)=\dfrac{1-2x+x^2+3x^4+x^5+2x^6}{(1-x)(1-2x)}.$$
\end{proposition}

\begin{proof}
S'obtient par un calcul simple.
\end{proof}

%********************************************************************************
%$$$$$$$$$$$$$$$$$$$$$$$$$$$$$$$$$$$$$$$$$$$$$$$$$$$$$$$$$$$$$$$$$$$$$$$$$$$$$$$$$
%$$$$$$$$$$$$$$$$$$$$$$$$$$$$$$$$$$$$$$$$$$$$$$$$$$$$$$$$$$$$$$$$$$$$$$$$$$$$
%\bigskip
\section{Le graphe  $G_6$}
Les graphes $G_6$ et $G'_6$ ont m\^eme profil, nous \'etudions alors un seul d'eutre eux, ce sera le graphe $G_6$.

Le graphe $G_6:=(V_6,E_6)$ est tel que $V_6:=(\mathbb N \times \{0,1\})\cup \{a,b\}$ ($a,b\notin \mathbb N \times \{0,1\}$) se d\'ecompose en deux sous-ensembles, une clique $(\mathbb N \times \{0\})\cup \{a\}$ et un stable $(\mathbb N \times \{1\})\cup \{b\}$. Une paire de sommets $\{(i,0),(j,1)\}$ est une ar\^ete de $G_6$ si $i\leq j$, l'ensemble des ar\^etes $E_6$ compte, en plus, la paire $\{a,b\}$.

\subsection{Repr\'esentation des sous-graphes de $G_6$}
Observons que $G_6\setminus\{a,b\}$ est isomorphe \`a $G_5\setminus \{c\}$. Donc, ils ont le m\^eme \^age.  Observons \'egalement que $Age(G_6\setminus\{a,b\})=Age(G_6\setminus\{a\})=Age(G_6\setminus\{b\})\neq Age(G_6)$. Donc, le noyau, $Ker(G_6)=\{a,b\}$ et
  $Age(G_6)$ est form\'e des types d'isomorphie des sous-graphes de $G_6\setminus\{a,b\}$ et des types d'isomorphie des sous-graphes de $G_6$ contenant les sommets $a$ et $b$ et ne s'abritant pas dans $G_6\setminus\{a,b\}$. Il convient de signaler que tous les sous-graphes d'ordres $n\leq 3$ de $G_6$ s'abritent dans $G_6\setminus\{a,b\}.$

\vspace{1mm}

Les sous-graphes d'ordre $n\geq 4$ de $G_6$ contenant les sommets $a$ et $b$ sont obtenus en rajoutant les sommets $a$ et $b$ \`a tous les sous-graphes d'ordre $n-2$ de $G_6\setminus\{a,b\}$. Ces derniers  sont repr\'esent\'es, comme nous l'avons vu ci-dessus pour le graphe $G_5$, par des suites de $\{0,1\}^{n-2}$, avec la m\^eme remarque que pour $G_5$, les suites qui diff\`erent uniquement par le dernier terme donnent le m\^eme type d'isomorphie dans  $G_6\setminus \{a,b\}$ mais lorsque nous rajoutons les sommets $a,b$, ils donnent deux types d'isomorphie diff\'erents.

\vspace{1mm}

Nous adopterons la m\^eme repr\'esentation par des suites de $\{0,1\}^n$ pour les sous-graphes de $G_6\setminus\{a,b\}$. Pour les sous-graphes contenant $a,~b$, nous rajoutons $ab$ \`a gauche du premier terme de la suite qui repr\'esente le sous-graphes sans les sommets $a,~b$. Par exemple, si $H$ est un sous-graphe de $G_6\setminus\{a,b\}$ repr\'esent\'e par la suite $S=00101$, le sous-graphe $H'$ obtenu en rajoutant \`a $H$ les sommets $a,~b$ aura pour suite $S'=ab00101$.% o\`u le sommet $a$ sera reli\'e \`a tous les sommets repr\'esent\'es par des $0$.

\begin{lemma}\label{lem:lemme-G6}
Un sous-graphe $H$ d'ordre $n\geq 4$ de $G_6$ s'abrite dans $G_6\setminus\{a,b\}$ si et seulement s'il ne contient pas un sous-graphe isomorphe \`a $P_4$.
\end{lemma}

\begin{proof}
Il est facile de v\'erifier que $P_4$ ne s'abrite pas dans $G_6\setminus\{a,b\}$, d'o\`u la condition n\'ecessaire. Pour la condition suffisante, montrons que si $H$ ne s'abrite pas dans $G_6\setminus\{a,b\}$, il contient n\'ecessairement $P_4$. Si $H$ ne s'abrite pas dans $G_6\setminus\{a,b\}$, donc $H$ contient les sommets $a,b$. Nous avons forc\'ement $V_H\cap(\mathbb N\times\{1\})\neq\varnothing$, $V_H\cap(\mathbb N\times\{0\})\neq\varnothing$ et
 $i_0=min\{i\in\mathbb N/(i,0)\in V_H\}\leq j$ pour au moins un $j$ tel que $(j,1)\in V_H$ car sinon, $H$ s'abriterait dans $G_6\setminus\{a,b\}$. Les sommets $a, b, (i_0,0)$ et $(j,1)$ forment un $P_4$.
\end{proof}

\vspace{2mm}

En consid\'erant la repr\'esentation donn\'ee par la \figurename~\ref{schemaG5} de la page \pageref{schemaG5}, nous avons

\begin{lemma}\label{lem:lemm1G6}
Un sous-graphe $H$ d'ordre $n\geq 4$ de $G_6$ contenant les sommets $a$ et $b$ s'abrite dans $G_6\setminus\{a,b\}$ si et seulement si nous avons l'une des situations suivantes:
\begin{enumerate}
\item[(1)] $Y^1_H=\varnothing$.
\item[(2)] $\vert Y^1_H\vert\geq 1$ et $X^0_H=Y^0_H=\varnothing$.
\end{enumerate}
\end{lemma}

\begin{proof}
La condition suffisante se v\'erifie de fa\c{c}on  imm\'ediate. Pour la condition n\'ecessaire elle utilise le Lemme \ref{lem:lemme-G6}.  Supposons que $H$ s'abrite dans $G_6\setminus\{a,b\}$ et soit la repr\'esentation de $H\setminus\{a,b\}$ sur $L_0,~L_1$ alors
\begin{enumerate}
\item Si la repr\'esentation n'utilise aucun point de $L_1$, alors $H$  s'abrite dans $G_6\setminus\{a,b\}$ et nous sommes  dans le cas $(1)$.
\item Si nous avons au moins un point sur $L_1$, alors nous avons n\'ecessairement $X^0_H=Y^0_H=\varnothing$, car sinon $H$ abriterait $P_4$.
\end{enumerate}
\end{proof}

\vspace{2mm}

\begin{lemma}\label{lem:lemmeG6}
 Soit $H$ un sous-graphe, d'ordre $n\geq 4$, contenant $a,~b$ et soit $S=ab u_0u_1\cdots u_{n-3}$ la suite qui lui est associ\'ee. Le sous-graphe $H$ s'abrite dans $G_6\setminus\{a,b\}$ si et seulement si $S$ est de la forme $ab\underset{i}{\underbrace{1\cdots 1}}\underset{n-2-i}{\underbrace{0\cdots 0}}$ avec $0\leq i\leq n-2$.
\end{lemma}

\begin{proof}
%Les trois cas de types isomorphiques obtenus pour $i=0,\; i=n-2$ et  $1\leq i \leq n-3$ apr\`es rajout de $\{a,b\}$ sont illustr\'es dans la repr\'esentation ci-dessous (voir*****).
Pour la condition suffisante, il suffit de remarquer que
\begin{enumerate}
\item Le sous-graphe (obtenu pour $i=0$) rep\'esent\'e par $ab\underset{n-2}{\underbrace{00 \cdots 0}}$ est isomorphe au sous-graphe repr\'esent\'e par $01\underset{n-2}{\underbrace{00\cdots 0}}$.
\item Le sous-graphe (obtenu pour $i=n-2$) repr\'esent\'e par $ab\underset{n-2}{\underbrace{11 \cdots 1}}$ est isomorphe au sous-graphe repr\'esent\'e  par $\underset{n-2}{\underbrace{11\cdots 1}}01$.
\item Le sous-graphe (obtenu pour $1\leq i\leq n-3$) repr\'esent\'e par $ab\underset{i}{\underbrace{11 \cdots 1}}\underset{n-2-i}{\underbrace{0\cdots 0}}$ est isomorphe au sous-graphe repr\'esent\'e  par $\underset{i}{\underbrace{11\cdots 1}}01\underset{n-2-i}{\underbrace{00\cdots 0}}$.
    \end{enumerate}
    Pour la condition n\'ecessaire, si $H$ s'abrite dans $G_6\setminus\{a,b\}$ alors nous avons l'une des conditions du Lemme \ref{lem:lemm1G6}.
    \begin{enumerate}
\item Si c'est la condition $(1)$ qui est v\'erifi\'ee, alors $S=ab\underset{n-2}{\underbrace{0\cdots 0}}$ qui est de la forme donn\'ee pour $i=0$.
\item Si c'est la condition $(2)$ qui est v\'erifi\'ee, alors $S=ab\underset{i}{\underbrace{1\cdots 1}}\underset{n-2-i}{\underbrace{0\cdots 0}}$, où $i=\vert Y^1_H\vert$.
    \end{enumerate}
\end{proof}

\begin{lemma}\label{lem:lemme2G6}
Soient $H$ et $H'$ deux sous-graphes de m\^eme ordre $n\geq 4$ contenant les sommets $a,~b$ et ne s'abritant pas dans $G_6\setminus\{a,b\}$ et soient $S=abu_0\cdots u_{n-3}$ et $S'=abu'_0\cdots u'_{n-3}$ les suites qui les repr\'esentent respectivement. Alors $$ H \text{ et }H'\text{ sont isomorphes }\Leftrightarrow S=S'.$$
\end{lemma}

\begin{proof}
La condition suffisante est \'evidente. Pour la condition n\'ecessaire, supposons $H$ et $H'$ isomorphes. Comme ils ne s'abritent pas dans $G_6\setminus\{a,b\}$, les suites $S$ et $S'$ n'ont pas la forme donn\'ee par le Lemme \ref{lem:lemmeG6} et tout isomorphisme qui envoie $H$ sur $H'$ laisse invariant les sommets $a$ et $b$. Donc $H\setminus\{a,b\}$ et $H'\setminus\{a,b\}$ sont isomorphes. Comme $G_6\setminus\{a,b\}$ et $G_5\setminus\{c\}$ sont isomorphes, d'apr\`es le Lemme \ref{lem:sanscG5}, nous avons $u_0\cdots u_{n-4}=u'_0\cdots u'_{n-4}.$  Si $u_{n-3}\neq u'_{n-3}$, $H$ et $H'$ ne sont pas isomorphes comme signal\'e ci-dessus (le sous-graphe pour lequel le dernier terme de la suite est $0$ abrite une clique de taille sup\'erieure à celle abrit\'ee par l'autre sous-graphe). Donc $u_{n-3}= u'_{n-3}$ et par suite $S=S'$.
\end{proof}

       \subsection{Les ind\'ecomposables de $Age(G_6)$}
\begin{figure}[h]
\centering
\input{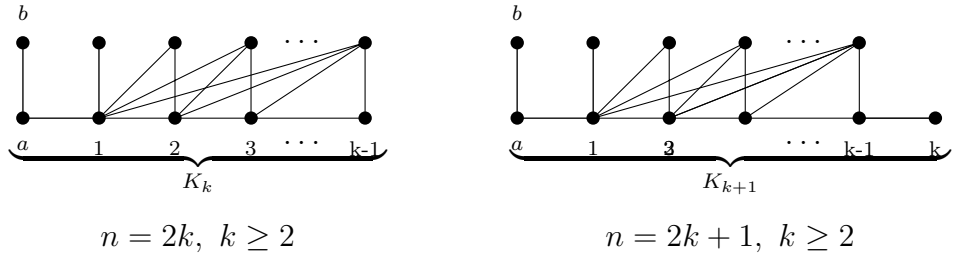}
\caption{\label{indecG6}Les sous-graphes ind\'ecomposables de $G_6$.}
\end{figure}

Comme pour le graphe $G_5$, tous les sous-graphes ind\'ecomposables de $G_6$ d'ordre $n\geq 4$ contiennent les sommets $a,~b$. Ce sont donc tous les sous-graphes repr\'esent\'es par des suites $S=abu_0\cdots u_{n-3}$ avec $u_0=0$ et $u_i\neq u_{i+1}$ pour tout $i,~0\leq i\leq n-4.$
Il s'ensuit donc qu'il existe un ind\'ecomposable d'ordre $n$ pour tout $n\geq 4$. Il sont donn\'es dans la \figurename~\ref{indecG6}. L'ind\'ecomposable d'ordre $n$ s'abrite dans tout ind\'ecomposable de taille sup\'erieure. Les sous-graphes ind\'ecomposables forment une cha\^{i}ne.  Le sous-graphe ind\'ecomposable d'ordre $n$ se trouve au niveau $n$ pour $n\leq 2$ et au niveau $n-1$, pour $n\geq 4$. La classe des ind\'ecomposables est minimale dans $Ind(\Omega_1)$ et $Age(G_6)$ est ind-minimal.

       \subsection{Profil de $Age(G_6)$}

Les premi\`eres valeurs du profil pour $n=0,1,2,3,4,5,6$ sont respectivement $1,1,2$, $4,9$, $20,43$.
Posons $\varphi_6$ la fonction profil de $Age(G_6)$. Nous avons

\begin{lemma}
Le profil $\varphi_6$ est exponentiel, il est donn\'e par:
$$\left\{
\begin{array}{ll}
\varphi_6(0)=\varphi_6(1)=1;&\\
\varphi_6(n)= 2^{n-1}+2^{n-2}-n+1&  \text{ pour }\; n\geq 2%
\end{array}%
\right. $$
En outre, $\varphi_6\simeq 3.2^{n-2}$.
\end{lemma}

\begin{proof}
Nous avons $\varphi_6(n)=\varphi_{-\{a,b\}}(n)+\varphi_{a,b}(n)$ pour tout $n\geq 1$, o\`u $\varphi_{-\{a,b\}}(n)$ repr\'esente le nombre de types d'isomorphie de sous-graphes d'ordre $n$ de $G_6\setminus\{a,b\}$ et $\varphi_{a,b}(n)$ le nombre de types d'isomorphie de sous-graphes d'ordre $n$ de $G_6$ contenant les sommets $a$ et $b$ et ne s'abritant pas dans $G_6\setminus\{a,b\}$, donc $\varphi_{a,b}(n)=0$ pour $n\leq 3.$
\smallskip

 Comme $G_6\setminus\{a,b\}$ est isomorphe \`a $G_5\setminus\{c\}$ nous avons $\varphi_{-\{a,b\}}(n)=2^{n-1}$ pour tout $n\geq 1.$

\noindent D'apr\`es les Lemmes \ref{lem:lemmeG6} et \ref{lem:lemme2G6} nous avons  $\varphi_{a,b}(n)=2^{n-2}-(n-1)$ pour tout $n\geq 4.$
En faisant la somme nous obtenons le r\'esultat donn\'e.
\end{proof}

\begin{proposition}
La fonction g\'en\'eratrice de $Age(G_6)$ est rationnelle et est donn\'ee par:
$$F_{G_6}(x)=\dfrac{1-3x+3x^2-x^3+x^4}{(1-2x)(1-x)^2}.$$
\end{proposition}

\begin{proof}
S'obtient par un calcul simple.
\end{proof} 
\clearemptydoublepage

%$$$$$$$$$$$$$$$$$$$$$$$$$$$$$$$$$$$$$$$$$$$$$$$$$$$$$$$$$$$$$$$$$$$$$$$$$$$
\part{D\'ecomposition monomorphe et profil}\label{part:monomorphie}
%$$$$$$$$$$$$$$$$$$$$$$$$$$$$$$$$$$$$$$$$$$$$$$$$$$$$$$$$$$$$$$$$$$$$$$$$$$$

\chapter*{Introduction}
\addcontentsline{toc}{chapter}{Introduction}

Nous nous int\'eressons dans cette partie au ph\'enom\`ene du saut\index{profil!saut} dans le comportement des profils des classes h\'er\'editaires, ph\'enom\`ene mis en \'evidence par \emph{Pouzet}\index{Pouzet} (voir \cite{pouzet.tr.1978}) il y a de cela plus d'une trentaine d'ann\'ees. Durant les dix derni\`eres ann\'ees, de nombreux articles ont \'et\'e consacr\'es \`a l'\'etude du profil des classes h\'er\'editaires de structures combinatoires telles que les graphes (dirig\'es ou non), les graphes ordonn\'es\index{graphe!ordonn\'e}, les hypergraphes, les ensembles ordonn\'es, les tournois, les permutations, $\ldots$. Tous les r\'esultats de ces \'etudes montrent qu'il existe des sauts dans le comportement des profils des classes h\'er\'editaires form\'ees pas ces structures. Typiquement, le profil est soit polynomial\index{profil!polynomial} ou cro\^{i}t plus vite que tout polyn\^ome et pour certaines classes de structures, il est au moins exponentiel comme dans le cas des tournois (\cite{B-B-M(07),Bou-Pouz}), des graphes ordonn\'es et hypergraphes ordonn\'es (\cite{B-B-M(06),klazar08}) et des permutations (voir Th\'eor\`eme \ref{theo:kaiser-klazar}, \cite{K-K}) ou a au moins la croissance de la fonction partition d'entier\index{fonction!partition d'entier} comme dans le cas des graphes (\cite{B-B-S-S}); pour plus de d\'etails voir \cite{klazar}.

\vspace{3mm}

Dans le chapitre \ref{chap:monomorphe}, nous consid\'erons des classes h\'er\'editaires de structures relationnelles ordonn\'ees. %et nous apportons, dans ce cas, des r\'eponses aux questions pos\'ees. %dans certains cas, en particulier dans le cas des structures ordonn\'ees o\`u
Nous d\'ecrivons les classes dont le profil est born\'e par un polyn\^ome et montrons que ce profil est en fait un polynôme. Nous identifions celles dont le profil n'est pas born\'e par un polyn\^ome et qui sont minimales pour l'inclusion. %et  montrons,

\vspace{1mm}

Au chapitre \ref{chap:graphe ordonne}, %que lorsque ces classes sont form\'ees de structures binaires, leurs profils sont exponentiels\index{profil!exponentiel}. % si  le profil d'une classe h\'er\'editaire de structures binaires ordonn\'ees  n'est pas born\'e par un polyn\^ome, il  est au moins exponentiel.  Un r\'esultat faisant partie d'une classification plus profonde obtenue par   Balogh, Bollob\'{a}s et  Morris en (2006) pour les graphes et les hypergraphes ordonn\'es. %que le saut du profil se fait d'une croissance polyn\^omiale \`a une croissance exponentielle.
nous \'etudions de mani\`ere plus d\'etaill\'ee le cas binaire. Nous montrons que le profil d'une structure binaire ordonn\'ee de type $k$ (ie, comportant outre la relation d'ordre $k$ autres relations binaires) est soit polynomial soit born\'e par une exponentielle. %structures binaires ordonn\'ees
 Pour ce faire, nous d\'ecrivons, dans le cas  des graphes dirig\'es ordonn\'es, tous les graphes dont les \^ages sont minimaux parmi ceux des structures sans d\'ecomposition monomorphe finie. Ces graphes sont au nombre de mille deux cent quarante six.
Les outils utilis\'es sont le th\'eor\`eme de \emph{Ramsey}\index{theoreme@th\'eor\`eme de Ramsey}, donn\'e en termes de structures invariantes, la notion de d\'ecomposition monomorphe d'une structure relationnelle introduite par \emph{Pouzet-Thi\'ery} \cite{P-T-2013}   qui \'etend %d\'ecoule de
la notion de monomorphie due \`a \emph{Fra\"{\i}ss\'e}\index{Fraisse@Fra\"{\i}ss\'e} \cite{fraisse} et la notion de presque multi-encha\^{i}nabilit\'e due \`a \emph{Pouzet} (voir \cite{pouzet06}).

\clearemptydoublepage

\chapter[Structures sans d\'ecomposition monomorphe finie]{Structures sans d\'ecomposition monomorphe finie. Application au profil des classes h\'er\'editaires}\label{chap:monomorphe}

\section{Encha\^{i}nabilit\'e et monomorphie}
    \subsection{D\'efinitions et rappel de r\'esultats}
      Soit $\mathcal R$ une structure relationnelle. Nous rappelons qu'un \emph{automorphisme local} de $\mathcal R$ est un isomorphisme entre deux restrictions de $\mathcal R$.
    \begin{definition}
   Soit $\mathcal R:= (E,(\rho_i)_{i\in I})$ une structure relationnelle et $\mathcal S$ une autre structure relationnelle ayant le m\^eme domaine que $\mathcal R$. La structure $\mathcal R$ est dite \emph{libre interpr\'etable}\index{structure relationnelle!libre interpr\'etable}  en $\mathcal S$, on dit aussi que $\mathcal S$ \emph{interpr\`ete librement} $\mathcal R$, si tout automorphisme local de $\mathcal S$ est un automorphisme local de $\mathcal R$.
     \end{definition}
     \vspace{1mm}

     Si $\mathcal S$ est une cha\^{i}ne,  $\mathcal R$ est dite \emph{encha\^{i}n\'ee}\index{structure relationnelle!encha\^{i}n\'ee} par ${\mathcal S}$ %\emph{${\mathcal S}$-encha\^{i}nable}
   ou que $\mathcal S$ \emph{encha\^{i}ne} $\mathcal R$. La structure $\mathcal R$ est dite \emph{encha\^{i}nable}\index{structure relationnelle!encha\^{i}nable} s'il existe une cha\^{i}ne $C$ sur la base de $\mathcal R$ qui encha\^{i}ne $\mathcal R$. %Par exemple, si $\mathcal R$ est une structure binaire, alors $\mathcal R$ est encha\^{i}n\'ee par $\mathcal S$ si, pour tout $i\in I$:
%$(x,y) \in \rho_i\Leftrightarrow (x',y')\in \rho_i$ pour tout $x,y,x',y'$ tels que: $(x,y)\in \mathcal S\Leftrightarrow(x',y')\in \mathcal S.$

\vspace{1mm}

 \textbf{Exemples.}  Si $\mathcal R$  et $\mathcal S$ sont deux cha\^{i}nes, alors $\mathcal S$ encha\^{i}ne $\mathcal R$ si et seulement si $\mathcal R$ co\"{\i}ncide avec $\mathcal S$ ou avec sa cha\^{i}ne duale.

\vspace{1mm}

%si $C$ est une cha\^{i}ne, l'ordre total strict obtenu \`a partir de $C$ et l'ordre dual\index{ordre!dual} de $C$ sont encha\^{i}n\'es par $C$.\\
\noindent Si $G$ est un graphe non dirig\'e, alors $G$ est encha\^{i}nable si et seulement si $G$ est un graphe complet ou un graphe vide (ie, sans ar\^ete).

\vspace{1mm}

\noindent Si $\mathcal R$ est une structure binaire, alors $\mathcal R$ est encha\^{i}n\'ee par $\mathcal S$ si, pour tout $i\in I$:
$(x,y) \in \rho_i\Leftrightarrow (x',y')\in \rho_i$ pour tout $x,y,x',y'$ tels que: $(x,y)\in \mathcal S\Leftrightarrow(x',y')\in \mathcal S.$

\vspace{1mm}

Voici quelques exemples pour des structures d'arit\'e sup\'erieure \`a deux.\\
Soit $C:= (E \leq)$ une cha\^{i}ne alors $C$ encha\^{i}ne la relation "\emph{entre}", l'\emph{"ordre circulaire"} associ\'es \`a $C$ et la relation "\emph{entre}" associ\'ee \`a l'ordre circulaire et d\'efinie sur les quadruplets de $E$. Ces trois relations sont:
 \begin{itemize}
\item La relation \emph{"entre" associ\'ee \`a  $C$} est la relation ternaire  $B_C:= (E, b_C)$  o\`u $b_C$ est l'ensemble des triplets $(x_1,x_2,x_3)$ tels que ou bien $x_1<x_2<x_3$ ou bien $x_3<x_2<x_1$.
\item L'\emph{ordre circulaire associ\'e \`a $C$} ou la relation \emph{ternaire cyclique} est la relation ternaire $T_C:= (E, t_C)$ où $t_C$ est l'ensemble des triplets
$(x_1,x_2,x_3)$ tels que $x_{\sigma(1)}<x_{\sigma(2)}<x_{\sigma(3)}$ pour une permutation circulaire $\sigma$ de $\{1,2,3\}$.
\item La relation \emph{"entre" associ\'ee \`a  l'ordre circulaire} est la relation quaternaire (d'arit\'e quatre) $D_C:=(E, d_C)$ où $d_C$ est l'ensemble des quadruplets
$(x_1,x_2,x_3, x_4)$ tels que ou bien $x_{\sigma(1)}<x_{\sigma(2)}<x_{\sigma(3)}<x_{\sigma(4)}$ ou bien $x_{\sigma(4)}<x_{\sigma(3)}<x_{\sigma(2)}<x_{\sigma(1)}$ pour une permutation circulaire $\sigma$ de $\{1,2,3, 4\}$.
 \end{itemize}

\bigskip

 Soit $n$ un entier positif, une structure relationnelle  $\mathcal R$ est dite \emph{$n$-monomorphe}\index{structure relationnelle!$n$-monomorphe} si les restrictions de $\mathcal R$ \`a deux $n$-ensembles quelconques sont isomorphes. Il est \`a signaler que toute relation est $0$-monomorphe, toute relation d\'efinie sur un $n$-ensemble est $n$-monomorphe, toute relation binaire r\'eflexive est $1$-monomorphe et tout tournoi est $2$-monomorphe.
\medskip

 Une structure relationnelle $\mathcal R$ est dite \emph{($\leq n$)-monomorphe}\index{structure relationnelle!$\leq n$-monomorphe} si elle est $p$-monomorphe pour tout $p\leq n$. La structure relationnelle $\mathcal R$ est dite \emph{monomorphe}\index{structure relationnelle!monomorphe} si elle est $n$-monomorphe pour tout entier $n$. Il est \`a signaler qu'une structure relationnelle de base $E$ est monomorphe d\`es qu'elle est $n$-monomorphe pour tout entier $n\leq \vert E\vert$.
\vspace{1mm}

  Si $\mathcal R$ est $n$-monomorphe alors toute restriction de $\mathcal R$ est $n$-monomorphe. Il en est de m\^eme si $\mathcal R$ est ($\leq n$)-monomorphe ou monomorphe.

 \medskip

  Une caract\'erisation simple de la $n$-monomorphie est donn\'ee par ce r\'esultat d\^u \`a Fra\"{\i}ss\'e.

\begin{lemma}\cite{fraisse}.

Une structure relationnelle $\mathcal R$ est $n$-monomorphe si et seulement si chaque restriction de $\mathcal R$ \`a un $(n+1)$-ensemble est $n$-monomorphe.
\end{lemma}

Nous rappelons \'egalement le r\'esultat suivant d\^u \`a Pouzet qui fait le lien entre les $n$-monomorphies.
\begin{lemma}\label{lem:r-monomorphie}\cite{fraisse,pouzet76}.

Toute structure relationnelle $\mathcal R$ qui est $n$-monomorphe et d\'efinie sur un $(n+p)$-ensemble est $r$-monomorphe pour tout $r\leq min(n,p)$. En particulier, si la base de $\mathcal R$ est infinie ou de cardinalit\'e au moins $2n-1$, alors $\mathcal R$ est $n$-monomorphe si et seulement si $\mathcal R$ est ($\leq n$)-monomorphe.
\end{lemma}

  %\subsubsection{Monomorphie et encha\^{i}nabilit\'e: rappel de r\'esultats et extensions}

 Comme exemple de relation monomorphe, une cha\^{i}ne, puisque deux cha\^{i}nes finies de m\^eme cardinalit\'e sont isomorphes. Plus g\'en\'eralement une relation encha\^{i}nable est monomorphe.

 La r\'eciproque n'est pas vraie; en effet, il existe des relations monomorphes non encha\^{i}nables, comme par exemple la relation binaire qui est un cycle \`a trois \'el\'ements.  Fra\"{\i}ss\'e (1954) a montr\'e que la r\'eciproque est vraie pour une structure infinie:

 \begin{theorem}\label{theo:fraisse-monomorphie}(voir \cite{fraisse,pouzet06}).

 Une structure relationnelle de base infinie est monomorphe si et seulement si elle est encha\^{i}nable.
 \end{theorem}

Suivant les traces de Fra\"{\i}ss\'e, Frasnay\index{Frasnay} (1965) \cite{frasnay 65} s'est int\'eress\'e au cas d'une structure de base finie et a montr\'e:

\begin{theorem}\cite{frasnay 65}

Pour tout entier $m\geq 1$ il existe un entier  $p$ tel que toute structure relationnelle $p$-monomorphe dont le maximum de l'arit\'e est au plus $m$ et le domaine de cardinalit\'e finie assez grande est encha\^{i}nable.
\end{theorem}

Le plus petit entier $p$ dans l'\'enonc\'e ci-dessus est le \emph{seuil de monomorphie $m$-aires}, not\'e $p(m)$. Ses valeurs ont \'et\'e donn\'ees par C. Franay en 1990 \cite{frasnay 90}.

\begin{theorem}
$p(1)=1$, $p(2)=3$, $p(m)=2m-2$ pour $m\geq 3$.
\end{theorem}

Pour un expos\'e d\'etaill\'e de ce r\'esultat (voir Fra\"{\i}ss\'e \cite{fraisse}, chapitre 13, particuli\`erement p 378).

\medskip

Une cons\'equence du Th\'eor\`eme \ref{theo:fraisse-monomorphie} est la suivante:
\begin{theorem}\cite{frasnay 65}

Toute structure relationnelle encha\^{i}nable infinie et d'arit\'e finie n'a qu'un nombre fini de bornes.
\end{theorem}

En effet, si $m$ est le maximum de l'arit\'e de la structure $\mathcal R$, le maximum de la cardinalit\'e des bornes est au plus le seuil de monomorphie $m$-aires (voir Franay 1973 \cite{frasnay 73}, Pouzet 1981 \cite{pouzet.81}).
%\vspace{2mm}

\subsection{G\'en\'eralisations}\label{subsec:generalisation-monomo}

Les notions d'encha\^{i}nabilit\'e et de  monomorphie ont \'et\'e g\'en\'eralis\'ees comme suit (\cite{pouzet06, P-T-2013}).

 Soit $\mathcal R$ une structure relationnelle de base $E$ et soit $F$ un sous-ensemble de  $E$. Nous disons que $\mathcal R$ est \emph{$F$-encha\^{i}nable}\index{structure relationnelle!$F$-encha\^{i}nable}  s'il existe une cha\^{i}ne $C$ sur $E\setminus F$ telle que tout isomorphisme local de $C$ \'etendu par l'identit\'e sur $F$ est un isomorphisme local de $\mathcal R$; nous disons que $\mathcal R$  est \emph{$F$-monomorphe}\index{structure relationnelle!$F$-monomorphe} si pour tout entier $n$ et toute paire de $n$-parties $A$ et $A'$ de $E\setminus F$, il existe un isomorphisme de  $\mathcal R_{\restriction_ A}$ sur  $\mathcal R_{\restriction_ {A'}}$ qui, \'etendu par l'identit\'e sur $F$, est un isomorphisme local de $\mathcal R$. Nous disons que $\mathcal R$ est \emph{presque-encha\^{i}nable}\index{structure relationnelle!presque-encha\^{i}nable}, respectivement \emph{presque-monomorphe}\index{structure relationnelle!presque-monomorphe},  s'il existe un ensemble fini $F$ pour lequel $\mathcal R$ est $F$-encha\^{i}nable,  respectivement $F$-monomorphe. Les r\'esultats mentionn\'es ci-dessus s'\'etendent \`a ces notions (voir \cite{fraisse, pouzet06}).

\subsection{Cas binaire}
%\subsubsection{$\textbf{n}$-monomorphie pour $n\leq 3$ et encha\^{i}nabilit\'e: cas des structures binaires}

$\textbf{1}$ - Toute relation unaire $1$-monomorphe est constante. De mani\`ere plus g\'en\'erale, toute relation unaire $n$-monomorphe d\'efinie sur au moins $n+1$ \'el\'ements est constante donc encha\^{i}nable \cite{frasnay 65}.

\medskip

$\textbf{2}$ - Soit $\rho$  une relation binaire de base $E$, alors

\quad $a)$ - $\rho$ est $1$-monomorphe si et seulement si $\rho$ est soit r\'eflexive soit irr\'eflexive.

\vspace{2mm}

\quad $b)$ - Si $\rho$ est $2$-monomorphe et $\vert E\vert\geq 3$ alors elle est $1$-monomorphe et soit sym\'etrique soit anti-sym\'etrique. Le cas r\'eflexif donne quatre types de relations: %(\`a l'isomorphie pr\`es),
la relation compl\`ete r\'eflexive, la relation d'\'egalit\'e\index{relation!d'\'egalit\'e} et une relation total (asym\'etrique) % $3$-cycle orient\'e
 r\'eflexive, qui peut-\^etre transitive (ordre total) ou non. Le cas irr\'eflexif donne quatre autres types de relations %(\`a l'isomorphie pr\`es)
 qui sont: la relation compl\`ete irr\'eflexive, la relation vide et un tournoi qui est transitif\index{tournoi!transitif} (un ordre strict) ou intransitif.

\vspace{2mm}

 \quad $c)$ - Le r\'esultat suivant,  d\^u \`a Frasnay\index{Frasnay} (\cite{frasnay 65}), \'etablit une relation entre la monomorphie et l'encha\^{i}nabilit\'e pour une relation binaire.

\begin{lemma}\label{lem:franay binaire}%\cite{frasnay 65}.
Si $\rho$ est une relation binaire, ($\leq 3$)-monomorphe d\'efinie sur un ensemble d'au moins quatre \'el\'ements alors $\rho$ est encha\^{i}nable.
\end{lemma}

\begin{proof}
Soit $\rho$ une relation binaire ($\leq 3$)-monomorphe sur un ensemble $E$ ayant au moins quatre \'el\'ements. $\rho$ \'etant $1$-monomorphe nous pouvons la supposer r\'eflexive d'apr\`es $2-a$). La $2$-monomorphie et le fait que $\vert E\vert\geq 4$ impliquent alors, d'apr\`es $2-b$), que $\rho$ est l'une des relations suivantes:
\begin{enumerate}
\item[(i)] la relation compl\`ete r\'eflexive ($\rho=\mathscr P(E)$),
\item[(ii)] la relation d'\'egalit\'e ($\rho=\Delta_E$),
\item[(iii)] la relation totale (asym\'etrique) r\'eflexive qui est soit transitive, soit non transitive.
\end{enumerate}
Dans les deux premiers cas, $\rho$ est encha\^{i}nable. Dans le cas (iii), si $\rho$ est transitive, c'est un ordre total, elle est donc encha\^{i}nable. Si $\rho$ n'est pas transitive, la $3$-monomorphie fait que, pour toute partie $X=\{a,b,c\}$ de $E$, la restriction de $\rho$ \`a $X$ contient un cycle de longueur trois, disons que $\rho$ est de la forme $$\{(a,a),(b,b),(c,c),(a,b),(b,c),(c,a)\}.$$
Comme $\vert E\vert\geq 4$ on obtient une contradiction avec le fait que $\rho$ est asym\'etrique; en effet, $\rho_{\restriction_{\{a,b,d\}}}$ contient le couple $(b,d)$ et $\rho_{\restriction_{\{b,c,d\}}}$ contient le couple $(d,b)$. Donc $\rho$ ne peut pas avoir cette forme. Il s'ensuit que $\rho$ est encha\^{i}nable.
\end{proof}

\medskip

Dans le Lemme \ref{lem:ench-monomorphe} ci-dessous, nous \'etendons le Lemme \ref{lem:franay binaire} \`a une structure relationnelle ordonn\'ee d'arit\'e finie.

\vspace{1mm}

Signalons que la seule relation $(\leq 3)$-monomorphe non encha\^{i}nable est le cycle\index{cycle} \`a trois \'el\'ements (r\'eflexif ou pas) et
qu'une relation binaire peut-\^etre $(\leq2)$-monomorphe de base infinie\index{base!infinie} mais non encha\^{i}nable (un tournoi poss\`edant un cycle de longueur trois).

\medskip

$\textbf{3}$ - Soit $\mathcal R:=(E,(\rho_i)_{i\in I})$ une structure relationnelle;

  \quad $a)$ - Si $\rho_i$ est binaire pour tout $i\in I$, alors  $\mathcal R$ est $1$-monomorphe si et seulement si pour tout $i\in I$, la relation $\rho_i$ est soit r\'eflexive ($\rho_i(x,x)=1$ pour tout $x\in E$) soit irr\'eflexive ($\rho_i(x,x)=0$ pour tout $x\in E$). \\ Si nous \'etendons la notion de r\'eflexivit\'e \`a une relation d'arit\'e quelconque en disant qu'une relation $\rho$ d'arit\'e $n$ est r\'eflexive (respectivement irr\'eflexive) si $\rho(x,\dots,x)=1$ (respectivement $\rho(x,\dots,x)=0$) pour tout $n$-uple $(x,\dots,x)\in E^n$, alors $\mathcal R$ est $1$-monomorphe si et seulement si pour tout $i\in I$, la relation $\rho_i$ est soit r\'eflexive soit irr\'eflexive.

\vspace{2mm}

   \quad $b)$ - Si $\mathcal R$ est $2$-monomorphe alors $(E,\rho_i)$ est $2$-monomorphe pour tout $i\in I$. La r\'eciproque n'est pas vraie, en effet, consid\'erons la structure binaire $\mathcal R:=(E,\rho_1,\rho_2)$ avec $E=\{x,y,z\}$, $\rho_1=\{(x,y),(y,z),(x,z)\}$ et $\rho_2=\{(x,y),(y,z),(z,x)\}$. Les relations $(E,\rho_1)$ et $(E,\rho_2)$ sont $2$-monomorphes mais $\mathcal R$ n'est pas $2$-monomorphe car ses restrictions aux ensembles $\{x,y\}$ et $\{x,z\}$ ne sont pas isomorphes.

\vspace{2mm}

   \quad $c)$ - Le Lemme \ref{lem:franay binaire} se g\'en\'eralise naturellement au cas d'une structure binaire. En effet, si $\mathcal R:=(E,(\rho_i)_{i\in I})$ est une structure binaire v\'erifiant les hypoth\`eses du Lemme \ref{lem:franay binaire}, alors $(E,\rho_i)$ est $(\leq 3)$-monomorphe pour tout $i\in I$. Donc, pour tout $i\in I$, si la relation $\rho_i$ est r\'eflexive, elle est donn\'ee par l'une des trois relations cit\'ees dans la preuve du Lemme \ref{lem:franay binaire} (relation compl\`ete r\'eflexive, relation d\'egalit\'e ou ordre total). Si $\rho_i$  est irr\'eflexive, elle donn\'ee par l'une des trois relations obtenues \`a partir des pr\'ec\'edentes en supprimant tous les couples de la formes $(a,a)$ pour $a\in E$ (ie, c'est une relation compl\`ete irr\'eflexive, une relation vide ou un tournoi transitif). Si chacune des relations $\rho_i$ est une relation compl\`ete r\'eflexive ou irr\'eflexive, une relation d'\'egalit\'e ou une relation vide alors $\mathcal R$ est naturellement encha\^{i}nable (elle est encha\^{i}n\'ee par toute cha\^{i}ne de base $E$). Si au moins une des relations $\rho_i$ est soit un ordre total soit un tournoi transitif (ie, un ordre strict transitif), alors la $2$-monomorphie de $\mathcal R$ implique que deux telles relations co\"{\i}ncident ou sont oppos\'ees sur toute paire de sommets distincts de $E$. Il s'ensuit que $\mathcal R$ est encha\^{i}nable.

    \subsection{Cas des structures relationnelles ordonn\'ees}

    Soit $\mathcal R:=(E,\leq,\rho_1,\dots,\rho_k)$ une structure relationnelle ordonn\'ee %(par exemple $\mathcal R:=(E,\leq,\rho_1,\dots,\rho_k)$)
    avec $\vert E\vert\geq 2$. Alors

\begin{lemma}\label{lem:reduction2monomorphe} Pour tout entier $p$ positif,
$\mathcal R$ est $p$-monomorphe si et seulement si toutes les relations $\mathcal R_i:=(E,\leq,\rho_i)$ pour $i\in\{1,\dots,k\}$ sont $p$-monomorphes.
\end{lemma}
\begin{proof}
Les structures \'etant ordonn\'ees, si deux restrictions d'une relation sont isomorphes l'isomorphisme est unique.
\end{proof}

\medskip

Si $\mathcal R$ est binaire et $(\leq2)$-monomorphe alors pour $i\in\{1,\dots,k\}$, la restriction de $\mathcal R_i:=(E,\leq,\rho_i)$ \`a un ensemble \`a deux \'el\'ements $\{x,y\}$ a l'une des formes donn\'ees dans la \tablename~\ref{exemple-relation}.

\begin{table}[!hbp]
\begin{center}
\begin{tabular}[c]{|c|c|c|c|}
\hline
\input{cas1}&\input{cas2}&\input{cas3}&\input{cas4}\\
\hline
\input{cas5}&\input{cas6}&\input{cas7}&\input{cas8}\\
\hline
\end{tabular}
\caption{\label{exemple-relation}Relations binaires $(\leq2)$-monomorphe \`a deux \'el\'ements.}
\end{center}
\end{table}

\vspace{1mm}

D'apr\`es le Lemme \ref{lem:reduction2monomorphe} et la liste des relations binaires $(\leq2)$-monomorphe \`a deux \'el\'ements (\tablename~\ref{exemple-relation}), nous avons

\begin{lemma}\label{lem:binaire2-monomorphe}
Une structure binaire ordonn\'ee $\mathcal R$ de type $k$ (form\'ee d'un ordre total et de $k$ relations binaires) est $(\leq 2)$-monomorphe si et seulement si chaque relation $\rho_i$, $i\in\{1,\dots,k\}$ est soit la relation compl\`ete r\'eflexive ($\rho_i=E\times E$) ou irr\'eflexive ($\rho_i=(E\times E)\setminus \Delta_E$), soit la relation d'\'egalit\'e ($\rho_i=\Delta_E$), soit un ordre total qui co\"{\i}ncide avec $\leq$ ou son dual, soit la relation vide ($\rho_i=\varnothing$) soit un tournoi transitif qui con\"{\i}cide avec l'ordre strict ($<$) ou son dual.
\end{lemma}

Ainsi nous avons exactement $8^k$ structures binaires ordonn\'ees de type $k$ qui sont $(\leq2)$-monomorphes. %encha\^{i}nables.
\medskip

 Le lemme suivant \'etend le Lemme \ref{lem:franay binaire} au cas d'une structure relationnelle ordonn\'ee. %Nous avons les r\'esultats suivants:

\begin{lemma}\label{lem:ench-monomorphe}
Soit $\mathcal R:=(E,\leq, (\rho_i)_{i\in I})$ une structure ordonn\'ee et
soit $m$ le maximum de l'arit\'e des relations $(\rho_i)_{i\in I}$, posons $C:=(E,\leq)$   alors $\mathcal R$ est encha\^{i}n\'ee par $C$ si et seulement si $\mathcal R$ est $(\leq m)$-monomorphe.
\end{lemma}

\begin{proof}
Toute structure relationnelle encha\^{i}nable \'etant monomorphe (Th\'eor\`eme \ref{theo:fraisse-monomorphie}), il suffit de prouver que la $(\leq m)$-monomorphie entra\^ {i}ne l'encha\^{i}nabilit\'e. %Si $\vert E\vert\leq m$, le r\'esultat est \'evident. Supposons $\vert E\vert\geq m+1$
Soit donc $f$ un automorphisme local de $C$. Notons $A$ et $A'$ son domaine et son codomaine. Pour montrer que $f$ est un automorphisme local de
$\mathcal R$ il suffit de montrer que chaque restriction de $f$ \`a une partie $X$ de cardinalit\'e au plus $m$ est un automorphisme local %\footnote{Voir rappel sur l'isomorphisme local donn\'e en Section \ref{sec:isom-abrit-equim} du Chapitre \ref{chap:generalite}.}
de $\mathcal R$.
Soit $X$ une telle partie. Comme $\mathcal R$ est $(\leq m)$-monomorphe, les restrictions $\mathcal R_{\restriction_X}$ et $\mathcal R_{\restriction_{f(X)}}$ sont isomorphes. Un tel isomorphisme, disons $\varphi$, est un isomorphisme de $C_{\restriction_X}$ sur $C_{\restriction_{f(X)}}$. Comme $C$ est une cha\^{i}ne, cet isomorphisme est unique, donc $f_{\restriction_X}=\varphi$. Ainsi $f$ est un isomorphisme de $\mathcal R_{\restriction_X}$ sur $\mathcal R_{\restriction_{f(X)}}$.
\end{proof}

\medskip

Les Lemmes \ref{lem:r-monomorphie} et \ref{lem:ench-monomorphe} entra\^{i}nent imm\'ediatement pour les structures binaires ordonn\'ees:

\begin{corollary}\label{cor:encha-monomorphe}
Si toutes les relations $\rho_i,~i\in I$, sont binaires, %$\mathcal R:=(V,\leq, (\rho_i)_{i\in I})$ une $2$-structure ordonn\'ee
alors $\mathcal R$ est encha\^{i}n\'ee par $C$  si et seulement si $\mathcal R$ est $2$-monomorphe.
\end{corollary}

\section{D\'ecompositions d'une structure relationnelle en parties monomorphes}

Dans cette section, nous consid\'erons des d\'ecompositions d'une structure relationnelle en parties monomorphes assujetties \`a certaines conditions.  Pour chaque condition, il y a  une d\'ecomposition qui est canonique en ce sens que toute autre d\'ecomposition satisfaisant la m\^eme condition est plus fine. Ces d\'ecompositions canoniques ne sont pas toutes identiques, mais elles ne diff\`erent que par les blocs finis.

\subsection{D\'ecomposition fortement monomorphe}\label{subsec:decomposition fortement monomorphe}

Cette d\'efinition vient naturellement de la g\'en\'eralisation de la notion de monomorphie donn\'ee dans le paragraphe \ref{subsec:generalisation-monomo}.

\begin{definition}\label{def:part-fortement monomorphe}
Soit $\mathcal R$ une structure relationnelle de base $E$. Un sous-ensemble $B$ de $E$ est un \emph{bloc fortement monomorphe}\index{bloc fortement monomorphe} de $\mathcal R$ si $\mathcal R$ est $E\setminus B$-monomorphe. Autrement dit, pour tout entier $n$ et pour toutes $n$-parties $A$ et $A'$ de $B$, il existe un isomorphisme de $\mathcal R_{\restriction_A}$ sur $\mathcal R_{\restriction_{A'}}$ qui, \'etendu par l'identit\'e sur $E\setminus B$, est un automorphisme local de $\mathcal R$.
\end{definition}

L'ensemble vide et les singletons de $E$ sont des blocs fortement monomorphes de $\mathcal R$.

\begin{lemma}\label{lem:union fortement-monomo}
L'union de deux blocs fortement monomorphes d'intersection non vide de $\mathcal R$ est un bloc fortement monomorphe de $\mathcal R$.
\end{lemma}

\begin{proof}
 Soient $X$ et $Y$ deux blocs fortement monomorphes de $\mathcal R$ avec $X\cap Y$ non vide. Montrons d'abord que si $A$ et $A'$ sont deux parties finies de m\^eme cardinalit\'e de $X\cup Y$ telles que leur diff\'erence sym\'etrique soit de cardinalit\'e 2 alors il existe un isomorphisme envoyant $A$ sur $A'$ qui, prolong\'e par l'identit\'e sur $E\setminus (X\cup Y)$ est un isomophisme local de $\mathcal R$. \\% (pour deux parties quelconques on passera de l'une a l'autre en enlevant et ajoutant un element).
 Soient $\{x, x'\}$ les deux \'el\'ements de la diff\'erence sym\'etrique. Ainsi pour $C:= A\cap A'$ nous pouvons supposer  $A= \{x\} \cup C$ et $A'= \{x'\}\cup C$.

  \textbf{Cas 1):} Supposons que $x,x'$ soient tous les deux dans  $X$ ou tous les deux dans $Y$. Supposons $x,x'\in X$. Comme $\mathcal R$ est $E \setminus X$-monomorphe et $\vert A\cap X\vert=\vert A'\cap X\vert$ alors les restrictions de $\mathcal R$ \`a $A$ et $A'$ sont isomorphes via l'identit\'e sur $E\setminus X$.

    \textbf{Cas 2):} Supposons que $x$ soit dans $X\setminus Y$ et $x'$ dans $Y\setminus X$. Nous avons alors deux cas.

    \quad \textbf{Cas 2.1:} $X\cap Y\not \subseteq C$. Soit $z\in  X\cap Y\setminus C$. Alors les restrictions de $\mathcal R$ \`a  $A'$ et $C\cup \{z\}$ sont isomorphes via l'identit\'e sur $E\setminus Y$. En remplaçant $A'$ par $C\cup \{z\}$, nous sommes dans le cas 1.

    \quad \textbf{Cas 2.2:} Si $X\cap Y\subseteq C$. Soit alors $z\in X\cap Y$. Les restrictions de $\mathcal R$ \`a  $C\setminus \{z\}\cup\{x,x'\}$ et $C\cup \{x'\}$ sont isomorphes via l'identit\'e sur $E\setminus X$; de m\^eme les restrictions de $\mathcal R$ \`a  $C\setminus \{z\}\cup\{x,x'\}$ et $C\cup \{x\}$ sont isomorphes via l'identit\'e sur $E\setminus Y$. En composant ces deux isomorphismes nous obtenons un isomorphisme entre les restrictions de $\mathcal R$ \`a $A$ et $A'$ qui prolong\'e par l'identit\'e sur $E\setminus (X\cup Y)$ est un isomorphisme local de $\mathcal R$.

 \vspace{1mm}

    Soient maintenant deux parties $A$ et $A'$, de m\^eme cardinalit\'e, de $X\cup Y$. Posons $K=A\setminus A'$ et $K'=A'\setminus A$. Si $\vert K\vert=0$, il n'y a rien \`a montrer. Si $\vert K\vert=1$, la preuve est donn\'ee ci-dessus. Supposons $\vert K\vert=k> 1$. Nous pouvons trouver une suite de $k+1$ sous-ensembles de $X\cup Y$, disons $A_0,\dots$, $A_i,\dots, A_k$ tels que $A_0=A$, $A_k=A'$ et la diff\'erence sym\'etrique de $A_i$ et $A_{i+1}$ est de cardinalit\'e $2$ ($A_{i+1}$ est obtenu en enlevant \`a $A_i$ un \'el\'ement de $A_i\cap K$ et en le remplaçant par un \'el\'ement de $K'\setminus A_i$). D'apr\`es ce qui pr\'ec\`ede, il existe, pour tout $i<k$, un isomorphisme envoyant $A_i$ sur $A_{i+1}$ qui, prolong\'e par l'identit\'e sur $E\setminus (X\cup Y)$ est un isomorphisme local de $\mathcal R$. Ceci induit un isomorphisme local de $\mathcal R$ qui envoie $A$ sur $A'$ et qui co\"{\i}ncide avec l'identit\'e sur $E\setminus (X\cup Y)$.
\end{proof}

\vspace{2mm}

Il s'ensuit que l'union d'une cha\^{i}ne croissante de blocs fortement monomorphes est un bloc fortement monomorphe.

\vspace{2mm}

Un bloc fortement monomorphe maximal pour l'inclusion est une \emph{composante fortement monomorphe}\index{composante!fortement monomorphe} de $\mathcal R$.

\begin{definition}
Une \emph{d\'ecomposition fortement monomorphe}\index{decomposition@d\'ecomposition!fortement monomorphe} de $\mathcal R$  est une partition $\mathscr P$ de $E$ en blocs fortement monomorphes.
\end{definition}

\begin{proposition}
Soit $\mathcal R$ une structure relationnelle de base $E$. Tout bloc fortement monomorphe de $\mathcal R$ est contenu dans une composante fortement monomorphe. Les composantes fortement monomorphes de $\mathcal R$ forment une d\'ecomposition fortement monomorphe de $\mathcal R$ dont toute autre d\'ecomposition fortement monomorphe est plus fine.
\end{proposition}

\begin{proof}
Soit un \'el\'ement $a\in E$. Comme $\{a\}$ est un bloc fortement monomorphe alors, d'apr\`es le Lemme \ref{lem:union fortement-monomo}, l'union $P(a)$ de tous les blocs fortement monomorphes de $\mathcal R$ contenant $a$ est un bloc fortement monomorphe. De plus, tout bloc fortement monomorphe maximal pour l'inclusion est de la forme $P(a)$.

Par ailleurs, deux composantes fortement monomorphes $P(a)$ et $P(b)$ sont soit \'egales soit disjointes, donc l'ensemble $\{P(a)\}_{a\in E}$ forme une partition de $E$ en composantes fortement monomorphes et si $P'$ est une autre partition de $E$ en parties fortement monomorphes alors chaque \'el\'ement de $P'$ est contenu dans une composante fortement monomorphe.
\end{proof}

\vspace{2mm}

La d\'ecomposition en composantes fortement monomorphes est dite \emph{d\'ecomposition fortement monomorphe canonique}.

\vspace{2mm}

Remarquons que l'intersections de blocs fortement monomorphes de $\mathcal R$ n'est pas un bloc fortement monomorphe de $\mathcal R$. Voici un contre- exemple:

\medskip

Consid\'erons la cha\^{i}ne $(E,\leq)$ o\`u $E=\{a,b,c,d\}$ et $a<b<c<d$. Soit $\mathcal R=(E,\rho)$ o\`u $\rho$ est l'ensemble des triplets $(x_1,x_2,x_3)$ d'\'el\'ements de $E$ tels que $x_1<x_2<x_3$ ou $x_{\sigma(1)}<x_{\sigma(2)}<x_{\sigma(3)}$ pour une permutation circulaire $\sigma$ de $\{1,2,3\}$. Pour tout \'el\'ement $x\in E$, l'ensemble $E\setminus\{x\}$ est un bloc fortement monomorphe; en effet, on peut d\'efinir un ordre total sur $E\setminus \{x\}$ en posant, pour tout $u,v\in E\setminus\{x\}$, $u\leq' v$ si $(x,u,v)\in\rho$. Tout isomorphisme local de cet ordre est un isomorphisme local de $\mathcal R_{\restriction_{E\setminus\{x\}}}$ et il le reste une fois prolong\'e par l'identit\'e sur $\{x\}$.

Consid\'erons les deux blocs fortement monomorphes $E\setminus\{a\}$ et $E\setminus\{c\}$, leur intersection est l'ensemble $\{b,d\}$ qui n'est pas un bloc fortement monomorphe; en effet, il n'existe aucun isomorphisme local de $\mathcal R$ qui envoie $b$ sur $d$ et qui fixe $a$ et $c$.

\vspace{1mm}

Notons par $dim_{fort}(\mathcal R)$ le nombre de composantes fortement monomorphe de la d\'ecomposition fortement monomorphe canonique de $\mathcal R$, si cette d\'ecomposition n'est pas finie, nous \'ecrivons $dim_{fort}(\mathcal R)=\infty$.

     \subsection{D\'ecomposition monomorphe}\label{subsec:decomposition monomorphe}
\begin{definition}\label{def:bloc monomorphe}
Soit $\mathcal R$ une structure relationnelle de base $E$. Un sous-ensemble $B$ de $E$ est un \emph{bloc monomorphe}\index{bloc monomorphe} de $\mathcal R$ si pour tout entier $n$ et pour toutes $n$-parties $A$ et $A'$ %de sous-ensembles \`a $n$ \'el\'ements
de $E$, les structures induites sur $A$ et $A'$ sont isomorphes d\`es que $A\setminus B=A'\setminus B$.
\end{definition}
\vspace{1mm}

Cette notion a \'et\'e introduite par Pouzet et Thi\'ery\index{Pouzet et Thi\'ery} \cite{P-T-2005,P-T-2013}.\\

    Le lemme suivant r\'eunit les propri\'et\'es les plus importantes des blocs monomorphes (voir \cite{Bou-Pouz}).

\begin{lemma}\label{lem:bloc monomorphe}
Soit $\mathcal R$ une structure relationnelle de base $E$.
\begin{enumerate}
\item[(i)] L'ensemble vide et les singletons de $E$ sont des blocs monomorphes de $\mathcal R$.
\item[(ii)] Si $B$ est un bloc monomorphe de $\mathcal R$ alors tout sous-ensemble de $B$ l'est aussi.
\item[(iii)] Si $B$ et $B'$ sont deux blocs monomorphes de $\mathcal R$ tels que $B\cap B'\neq \varnothing$ alors $B\cup B'$ est un bloc monomorphe de $\mathcal R$.
\item[(iv)] Soit $\mathscr B$ une famille de blocs monomorphes de $\mathcal R$; si $\mathscr B$ est filtrante %(c'est \`a dire que l'union de deux \'el\'ements de $\mathscr B$ est contenue dans un troisi\`eme),
    alors l'union $\mathcal B:=\underset{B_i\in\mathscr B}\bigcup B_i$ est un bloc monomorphe de $\mathcal R.$
\end{enumerate}
\end{lemma}

Un bloc monomorphe maximal pour l'inclusion est une \emph{composante monomorphe}\index{composante!monomorphe} de $\mathcal R$.

\begin{definition}\label{def:decomposition monomorphe}
Une \emph{d\'ecomposition monomorphe}\index{decomposition@d\'ecomposition!monomorphe} de $\mathcal R$  est une partition\index{partition!d'un ensemble} $\mathscr P$ de $E$ en blocs monomorphes. De mani\`ere \'equivalente, c'est une partition  $\mathscr P$ de $E$ en parties telles que pour tout entier $n$, les structures induites sur deux $n$-parties  $A$ et $A'$ de $E$ sont isomorphes d\`es que
$\left\vert A\cap B\right\vert =\left\vert A'\cap B\right\vert $ pour tout $B\in\mathscr{P}$.% les intersections $A\cap B$ et $A'\cap B$ avec toute partie $B$ de $\mathscr P$ ont la m\^eme cardinalit\'e.
\end{definition}
\medskip

Nous rappelons les r\'esultats suivants dus \`a Pouzet et Thi\'ery:
\begin{proposition}\label{prop:decomp-monomorphe}(voir \cite{P-T-2013})

Les composantes monomorphes de $\mathcal R$ forment une d\'ecomposition monomorphe qui est la moins fine de toutes les d\'ecompositions monomorphes de $\mathcal R$.
\end{proposition}

Cette d\'ecomposition monomorphe est dite \emph{minimale}\index{decomposition@d\'ecomposition!minimale} ou \emph{canonique},\index{decomposition@d\'ecomposition!canonique} elle est d\'esign\'ee par $\mathcal{P(R)}$. Nous notons par $dim_{mon}(\mathcal R)$ le nombre de blocs monomorphes de $\mathcal{P(R)}$, si $\mathcal{P(R)}$ n'est pas finie, nous \'ecrivons $dim_{mon}(\mathcal R)=\infty$. Si $dim_{mon}(\mathcal R)<\infty$ nous d\'esignons par $dim_{mon}^{\infty}(\mathcal R)$ le nombre de composantes monomorphes infinies, c'est \`a dire le nombre de blocs infinis de $\mathcal{P(R)}$. %est appel\'e  la \emph{dimension monomorphe}\index{dimension monomorphe} de $\mathcal R$.
%
%Nous rappelons est d\^u \`a Pouzet et Thi\'ery (2005) \cite{P-T-2005}.

\begin{theorem}\label{thm:chainable}(voir \cite{pouzet06, P-T-2013})

 Soit $\mathcal R$ une structure relationnelle de base $E$ et $A$ une partie de $E$. %est une composante monomorphe infinie d'une structure relationnelle $\mathcal R$ alors $\mathcal R$ est $V(\mathcal R) \setminus A$-encha\^{\i}nable.
 Consid\'erons les assertions suivantes:
 \begin{enumerate}
 \item[$(i)$] $\mathcal R$ est $E\setminus A$-encha\^{i}nable;
 \item[$(ii)$] $A$ est un bloc fortement monomorphe de $\mathcal R$;
 \item[$(iii)$] $A$ est un blocs monomorphe de $\mathcal R$.
 \end{enumerate}
 Alors $(i)\Rightarrow (ii)\Rightarrow (iii)$. Si $A$ est infinie alors $(ii)\Rightarrow (i)$. Si $A$ est une composante monomorphe infinie de $\mathcal R$ alors $(iii)\Rightarrow (ii)$.
\end{theorem}

Pour la preuve de ce Th\'eor\`eme voir la preuve du Th\'eor\`eme 2.25 de \cite{P-T-2013}.

\medskip

Comme cons\'equence nous avons:

\begin{proposition}\label{propos:monomo-fortmono}
Soit $\mathcal R$ une structure relationnelle. Toute d\'ecomposition fortement monomorphe de $\mathcal R$ est une d\'ecomposition monomorphe de $\mathcal R$. La r\'eciproque n'est pas vraie. En outre, nous avons $dim_{fort}(\mathcal R)<\infty$ %$\mathcal R$ poss\`ede une d\'ecomposition fortement monomorphe ayant un nombre fini de blocs
si et seulement si $dim_{mon}(\mathcal R)<\infty$. %$\mathcal R$ poss\`ede une d\'ecomposition monomorphe ayant un nombre fini de blocs.
De plus, le %nombre de composantes monomorphes infinies de $\mathcal R$ est \'egal au
nombre de composantes fortement monomorphes infinie de $\mathcal R$ est \'egal \`a $dim_{mon}^{\infty}(\mathcal R)$.
\end{proposition}

\begin{proof}
Le fait qu'une d\'ecomposition fortement monomorphe soit une d\'ecomposition monomorphe vient de l'implication $(ii)\Rightarrow (iii)$ du Th\'eor\`eme \ref{thm:chainable}. La r\'eciproque n'est pas vraie car un bloc monomorphe n'est pas n\'ecessairement un bloc fortement monomorphe, par exemple si $\mathcal R$ est une cha\^{i}ne de base $E$, toute partie de $E$ est un bloc monomorphe mais les blocs fortement monomorphes de $\mathcal R$ sont les intervalles de $\mathcal R$.

Par ailleurs, si $dim_{fort}(\mathcal R)<\infty$ %$\mathcal R$ poss\`ede une d\'ecomposition fortement monomorphe ayant un nombre fini de blocs
alors $dim_{mon}(\mathcal R)<\infty$. %$\mathcal R$ poss\`ede une d\'ecomposition monomorphe ayant un nombre fini de blocs.
Inversement, si $dim_{mon}(\mathcal R)<\infty$, %$\mathcal R$ poss\`ede une d\'ecomposition monomorphe ayant un nombre fini de blocs,
alors la d\'ecomposition canonique $\mathcal{P(R)}$ est finie. D'apr\`es le Th\'eor\`eme \ref{thm:chainable}, les composantes infinies dans $\mathcal{P(R)}$ sont des composantes fortement monomorphes. En partitionnant les \'el\'ements de la r\'eunion des blocs finis de $\mathcal{P(R)}$ en singletons, nous obtenons une d\'ecomposition fortement monomorphe finie. La derni\`ere partie de la proposition est due au fait que le nombre de composantes fortement monomorphes infinies est le nombre de composantes monomorphes infinies.
\end{proof}

\vspace{2mm}

Le th\'eor\`eme suivant donne une caract\'erisation d'une structure relationnelle poss\`edant une d\'ecomposition monomorphe finie.

\begin{theorem}\label{thm:charcfinitemonomorph}(voir \cite{P-T-2013})

 Une structure relationnelle $\mathcal R:=(E, (\rho_{i})_{i\in I})$ poss\`ede une d\'ecomposition monomorphe finie si et seulement s'il existe un ordre lin\'eaire $\leq$ sur $E$ et une partition finie $(E_x)_{x\in X}$ de $E$ en intervalles de $C:=(E, \leq)$ tels que tout isomorphisme local de $C$ qui pr\'eserve chaque intervalle est un isomorphisme local\index{isomorphisme!local} de $\mathcal R$. De plus, il existe une telle partition dont le nombre de blocs infinis est $dim_{mon}^{\infty}(\mathcal R)$. %le nombre de composantes infinie %la dimension monomorphe\index{dimension monomorphe} de $\mathcal R$.
\end{theorem}

En d'autres termes, la structure relationnelle $\mathcal R:=(E, (\rho_{i})_{i\in I})$ poss\`ede une d\'ecomposition monomorphe finie si et seulement si $\mathcal R$ est libre interpr\'etable\index{structure relationnelle!libre interpr\'etable} en une structure de la  forme $S:= (E, \leq, u_1,  \dots, u_l)$ o\`u $\leq$ est un ordre lin\'eaire et $u_1, \dots, u_l$ sont des relations unaires formant une partition de $E$ en intervalles\index{partition!en intervalles} de $(E, \leq)$.\\

Le th\'eor\`eme suivant d\'ecoule de l'application directe du th\'eor\`eme de compacit\'e de la logique du premier ordre: %(faire un bref rappel en annexe):

\begin{theorem} \label{compactness} Une structure relationnelle $\mathcal R$ poss\`ede une  d\'ecomposition monomorphe finie si et seulement si il existe un entier $\ell$ tel que tout \'el\'ement de $\mathcal {A(R)}$ poss\`ede une d\'ecomposition monomorphe ayant au plus  $\ell$ blocs.
\end{theorem}

La d\'ecomposition monomorphe, apr\`es son introduction dans \cite{P-T-2005}, a \'et\'e \'etudi\'ee dans plusieurs articles sur des structures particuli\`eres, les bicha\^{i}nes \cite{mont-pou}, les tournois \cite{Bou-Pouz} et les hypergraphes \cite{P-Kad}.

\subsection{D\'ecomposition en intervalles Fra\"{\i}ss\'e monomorphes}\label{subsec:decomposition Fraisse monomorphe}
\subsubsection{Notion d'intervalle Fra\"{\i}ss\'e}

\begin{definition}\label{def:intervalle fraisse}
Soit $\mathcal R:=(E, (\rho_i)_{i\in I})$ une structure relationnelle. Un sous-ensemble $J$ de $E$ est un intervalle\index{intervalle!Fra\"{\i}ss\'e} Fra\"{\i}ss\'e\footnote{La notion classique d'intervalle pour les structures binaires donn\'ee dans la D\'efinition \ref{def:intervalle} \'etant due \`a Gallai et afin d'\'eviter les confusions, nous utiliserons l'appellation ``intervalle Fra\"{\i}ss\'e'' pour l'intervalle de la D\'efinition \ref{def:intervalle fraisse} et ``intervalle Gallai'' pour l'intervalle de la D\'efinition \ref{def:intervalle}.} de $\mathcal R$ ou un $\mathcal R$-intervalle, si tout automorphisme local\index{automorphisme!local} de $\mathcal R_{\restriction_J}$ prolong\'e par l'identit\'e sur $E\setminus J$ est un automorphisme local de $\mathcal R$.
\end{definition}

Notons, d'apr\`es Fra\"{\i}ss\'e, que si $m$ est le maximum des arit\'es des $\rho_i$, alors $J$ est un $\mathcal R$-intervalle si et seulement si l'union d'un automorphisme local quelconque de $\mathcal R_{\restriction_J}$ d\'efini sur $p$ \'el\'ements $(p\leq m-1)$ et de l'identit\'e sur $m-p$ \'el\'ements de $E\setminus J$ donne un automorphisme local de $\mathcal R$. En particulier, si $\mathcal R$ est un ordre total, nous retrouvons la caract\'erisation de la notion classique d'intervalle (pour $m=2$ et $p=1$): pour tous $a,a' \in J$ et $x\notin J$;
$$a,a'<x (mod~\mathcal R)\text{ ou }a,a'>x (mod~\mathcal R).$$

Etant donn\'ee une structure relationnelle $\mathcal R$ de base $E$, l'ensemble vide, les singletons de $E$ et l'ensemble $E$ sont des intervalles Fra\"{\i}ss\'e de $\mathcal R$.

%\vspace{2mm}

\begin{fact}
Pour les structures binaires constitu\'ees de relations telles que chacune est soit r\'eflexive soit irr\'eflexive, les notions d'intervalle Fra\"{\i}ss\'e et intervalle Gallai co\"{\i}ncident. %avec la notion classique d'intervalle sur  %c'est \`a dire sur des relations ou des structures $1$-monomorphes.
\end{fact}

\vspace{1mm}

Pour une structure binaire, tout intervalle Gallai\index{intervalle!Gallai} est un intervalle Fra\"{\i}ss\'e, la r\'eciproque n'est pas vraie comme le montre l'exemple suivant. Soit $G:=(V,E)$ le graphe de la \figurename~\ref{contreexemple-intfraisse}. $J=\{x,y,z\}$ est un intervalle Fra\"{\i}ss\'e mais n'est pas un intervalle Gallai.

\begin{figure}[h]
\centering
\psset{unit=1cm}
\begin{pspicture}(-3,-2)(3,0.5)
\psdots[dotsize=5pt](-1,-1.5)(-2,0)(0,0)(2,0)(2,-1.5)
\psline[linewidth=0.3pt](-2,0)(-1,-1.5)(0,0)
\psline[linewidth=0.3pt](2,0)(2,-1.5)
\uput{0.3}[r](-1,-1.5){$a$}
\uput{0.3}[r](-2,0){$x$}
\uput{0.3}[r](0,0){$y$}
\uput{0.3}[r](2,0){$z$}
\uput{0.3}[r](2,-1.5){$b$}
\pscircle(-2,0.2){0.2}
\pscircle(0,0.2){0.2}
\end{pspicture}
\caption{\label{contreexemple-intfraisse}}
\end{figure}

\vspace{2mm}

Soit $\mathcal R$ une structure relationnelle de base $E$. Voici quelques propri\'et\'es des intervalles Fra\"{\i}ss\'e (voir \cite{fraisse}):
\begin{properties}\label{propriete:interv-fraisse}
\begin{enumerate}
\item Toute intersection de $\mathcal R$-intervalles est un $\mathcal R$-intervalle.
\item La r\'eunion d'un ensemble filtrant\index{ensemble!ordonn\'e!filtrant} pour l'inclusion de $\mathcal R$-intervalles est un $\mathcal R$-intervalle.
\item Si $J$ est un intervalle de $\mathcal{R}$ et $D\subseteq E$ alors $J\cap D$ est un intervalle de $\mathcal{R}_{\restriction_D}$.
\item Inversement, dans le cas particulier o\`u $D$ est un $\mathcal R$-intervalle, alors les intervalles de $\mathcal{R}_{\restriction_D}$ sont exactement les intersections avec $D$ des intervalles de $\mathcal R$.
\end{enumerate}
\end{properties}

La propri\'et\'e $4$ ci-dessus n'est pas vraie si $D$ est un sous-ensemble quelconque de $E$, comme le montre l'exemple suivant.

\vspace{1mm}

Soit l'ensemble ordonn\'e $P:=(V,\leq)$ dont le diagramme de Hasse\index{diagramme de Hasse} est donn\'e par la \figurename~\ref{exemple}. Posons $D=\{x,y,z\}$, l'ensemble $F=\{x,y\}$ est un intervalle de $P_{\restriction_D}$ mais il n'existe aucun intervalle $J$ de $P$ tel que $F=J\cap D$. En effet, si un tel intervalle $J$ existait, il contiendrait forc\'ement $w$ qui est comparable \`a $y$ mais incomparable \`a $x$ et donc contiendrait $z$, ce qui donnerait $J\cap D\neq F$.

\begin{figure}[h]
\centering
\psset{unit=1cm}
\begin{pspicture}(-3,-2)(3,2)
\psdots[dotsize=5pt](-1,-1.5)(-2,0)(0,0)(2,0)(2,-1.5)(1,1.5)
\psline[linewidth=0.3pt](-2,0)(-1,-1.5)(0,0)(1,1.5)(2,0)(2,-1.5)
\uput{0.3}[r](-1,-1.5){$u$}
\uput{0.3}[r](-2,0){$x$}
\uput{0.3}[r](0,0){$y$}
\uput{0.3}[r](2,0){$z$}
\uput{0.3}[r](2,-1.5){$v$}
\uput{0.3}[r](1,1.5){$w$}
\end{pspicture}
\caption{\label{exemple}}
\end{figure}

\vspace{2mm}

Remarquons \'egalement que contrairement aux intervalles Gallai des structures binaires, la r\'eunion de deux intervalles Fra\"{\i}ss\'e qui ont une intersection non vide n'est pas forc\'ement un intervalle Fra\"{\i}ss\'e comme le montre l'exemple suivant.

\vspace{1mm}

Soit $G=(V,E)$ le graphe de la \figurename~\ref{exemple-fraisse}. Les ensembles $\{a,b,c\}$ et $\{c,d,e\}$ sont des intervalles Fra\"{\i}ss\'e mais leur r\'eunion $A=\{a,b,c,d,e\}$ n'est pas un intervalle Fra\"{\i}ss\'e, car l'automorphisme local de $A$ qui envoie $a$ sur $d$ et $b$ sur $e$ prolong\'e par l'identit\'e sur $V\setminus A$ n'est pas un automorphisme local de $G$.

\begin{figure}[h]
\centering
\psset{unit=1cm}
\begin{pspicture}(-3,-2)(7,0.5)
\psdots[dotsize=5pt](-1,-1.5)(-2,0)(0,0)(2,0)(2,-1.5)(4,0)(6,0)(5,-1.5)
\psline[linewidth=0.3pt](-2,0)(-1,-1.5)(0,0)
\psline[linewidth=0.3pt](4,0)(5,-1.5)(6,0)
\psline[linewidth=0.3pt](2,0)(2,-1.5)
\uput{0.3}[r](-1,-1.5){$x$}
\uput{0.3}[r](-2,0){$a$}
\uput{0.3}[r](0,0){$b$}
\uput{0.3}[r](2,0){$c$}
\uput{0.3}[r](2,-1.5){$y$}
\uput{0.3}[r](5,-1.5){$z$}
\uput{0.3}[r](4,0){$d$}
\uput{0.3}[r](6,0){$e$}
\pscircle(-2,0.2){0.2}
\pscircle(0,0.2){0.2}
\pscircle(4,0.2){0.2}
\pscircle(6,0.2){0.2}
\end{pspicture}
\caption{\label{exemple-fraisse}}
\end{figure}

\subsubsection{D\'ecomposition Fra\"{\i}ss\'e monomorphe}
\begin{definition}
Soit $\mathcal R$ une structure relationnelle de base $E$. Un sous-ensemble $A$ de $E$ est un \emph{$\mathcal R$-intervalle monomorphe}\index{intervalle!Fra\"{\i}ss\'e!monomorphe} ou un \emph{intervalle Fra\"{\i}ss\'e monomorphe} si $A$ est un intervalle Fra\"{\i}ss\'e et la restriction, $\mathcal R_A$, de $\mathcal R$ \`a $A$ est monomorphe.
\end{definition}

Etant donn\'ee une structure relationnelle $\mathcal R$ de base $E$, l'ensemble vide et les singletons de $E$ sont des $\mathcal R$-intervalles monomorphes.

\vspace{1mm}

Voici quelques propri\'et\'es des $\mathcal R$-intervalles monomorphes qui d\'ecoulent naturellement de la d\'efinition et des Propri\'et\'es \ref{propriete:interv-fraisse}:

\begin{properties}
\begin{enumerate}
\item Toute intersection de $\mathcal R$-intervalles monomorphes est un $\mathcal R$-intervalle monomorphe.
\item La r\'eunion d'un ensemble filtrant pour l'inclusion de $\mathcal R$-intervalles monomorphes est un $\mathcal R$-intervalle monomorphe.
\item Si $J$ est un $\mathcal R$-intervalle monomorphe et $D\subseteq E$ alors $J\cap D$ est un $\mathcal{R}_{\restriction_D}$-intervalle monomorphe.
\end{enumerate}
\end{properties}

\begin{lemma}\label{lem:union inter-fraisse}
La r\'eunion de deux $\mathcal R$-intervalles monomorphes d'intersection non vide est un $\mathcal R$-intervalle monomorphe.
\end{lemma}

\begin{proof}
Vient de la d\'efinition et du fait que deux parties de m\^eme cardinalit\'e d'un $\mathcal R$-intervalle monomorphe sont toujours isomorphes. La preuve se fait de la m\^eme fa\c{c}on que celle du Lemme \ref{lem:union fortement-monomo}.
\end{proof}

\begin{definition}
Soit $\mathcal R$ une structure relationnelle de base $E$. Une \emph{d\'ecomposition Fra\"{\i}ss\'e monomorphe}\index{decomposition@d\'ecomposition!Fra\"{\i}ss\'e monomorphe} de $\mathcal R$ ou une \emph{$\mathcal R$-d\'ecomposition monomorphe} est une partition  de $E$ en $\mathcal R$-intervalles monomorphes.
\end{definition}

\begin{proposition}\label{prop:decomp-fraisse}
Soit $\mathcal R$ une structure relationnelle de base $E$. Tout $\mathcal R$-intervalle monomorphe est contenu dans un $\mathcal R$-intervalle monomorphe maximal pour l'inclusion. Les $\mathcal R$-intervalles monomorphes maximaux de $\mathcal R$ forment une $\mathcal R$-d\'ecomposition monomorphe dont toute autre $\mathcal R$-d\'ecomposition monomorphe est plus fine.
\end{proposition}

\begin{proof}
Soit un \'el\'ement $a\in E$. Comme $\{a\}$ est un $\mathcal R$-intervalle monomorphe alors, d'apr\`es le Lemme \ref{lem:union inter-fraisse}, l'union $F(a)$ de tous les $\mathcal R$-intervalles monomorphes contenant $a$ est un $\mathcal R$-intervalle monomorphe. De plus, tout $\mathcal R$-intervalle monomorphe maximal pour l'inclusion est de la forme $F(a)$.

Par ailleurs, deux $\mathcal R$-intervalles monomorphes $F(a)$ et $F(b)$ sont soit \'egaux soit disjoints, donc l'ensemble $\{F(a)\}_{a\in E}$ forme une partition de $E$ en $\mathcal R$-intervalles monomorphes et si $P'$ est une autre partition de $E$ en $\mathcal R$-intervalles monomorphes alors chaque \'el\'ement de $P'$ est contenu dans un $\mathcal R$-intervalle monomorphe maximal.
\end{proof}

\vspace{1mm}

Nous notons par $dim_{Fraisse}(\mathcal R)$ le nombre de $\mathcal R$-intervalles monomorphes maximaux et si ce nombre n'est pas fini, nous \'ecrivons $dim_{Fraisse}(\mathcal R)=\infty$. %et le nombre de ses blocs

\begin{lemma}\label{lem:R-intervalle-bloc.monomorphe}
Soit $\mathcal R$ une structure relationnelle de base $E$. Tout $\mathcal R$-intervalle monomorphe est un bloc monomorphe.
\end{lemma}

\begin{proof}
Soit $A$ un $\mathcal R$-intervalle monomorphe et soient $B$ et $B'$ deux $n$-parties de $E$ telles que $B\setminus A=B'\setminus A$. Posons $\overline{B}=B\cap A$ et $\overline{B'}=B'\cap A$. Comme $A$ est un $\mathcal R$-intervalle monomorphe, il existe un isomorphisme qui envoie $\overline{B}$ sur $\overline{B'}$ et qui, prolong\'e par l'identit\'e sur $E\setminus A$ est un automorphisme local de $\mathcal R$. %Comme $f(B)=B'$
Donc les restrictions de $\mathcal R$ \`a $B$ et $B'$ sont isomorphes.
\end{proof}

\vspace{1mm}

La r\'eciproque du Lemme \ref{lem:R-intervalle-bloc.monomorphe} est fausse. En effet, il suffit de consid\'erer une partie d'une cha\^{i}ne qui est un bloc monomorphe mais pas n\'ecessairement un intervalle Fra\"{\i}ss\'e.

\begin{lemma}\label{lem:monomo-forte.monom}
Soit $\mathcal R$ une structure relationnelle. Tout $\mathcal R$-intervalle monomorphe est un bloc fortement monomorphe de $\mathcal R$.
\end{lemma}

\begin{proof}
Si $A$ est un $\mathcal R$-intervalle monomorphe. Soient $B, B'$ deux $n$-parties de $A$ ($n\in\NN$). Comme $\mathcal R_{\restriction_A}$ est monomorphe, alors $\mathcal R_{\restriction_B}$ et $\mathcal R_{\restriction_{B'}}$ sont isomorphes. Soit $f$ cet isomorphisme, $f$ prolong\'e par l'identit\'e sur $E\setminus A$ est un automorphisme local de $\mathcal R$ car $A$ est un $\mathcal R$-intervalle. D'où $A$ est une partie fortement monomorphe.
\end{proof}

\vspace{1mm}

La r\'eciproque du Lemme \ref{lem:monomo-forte.monom} est fausse. En effet, consid\'erons le contre-exemple suivant.

\noindent Soit $\rho$ une relation ternaire d\'efinie sur $\mathbb N$ comme suit:
$$\rho=\{(x,y,z):x,y,z\in\mathbb N\setminus\{0\}\}\cup\{(0,i,j):0<i<j\}$$
alors $\mathbb N\setminus\{0\}$ est une partie fortement monomorphe, car $\rho$ est $\{0\}$-monomorphe, mais $\mathbb N\setminus\{0\}$ n'est pas un intervalle Fra\"{\i}ss\'e. %il suffit de consid\'erer un isomorphisme qui envoie x<y sur x'>y', prolong\'e par l'identit\'e sur $\{0\}$ ce n'est pas un automorphisme local de $\rho$.

\begin{corollary}\label{cor:decompfrai-decompfortmono}
Soit $\mathcal R$ une structure relationnelle de base $E$. Toute $\mathcal R$-d\'ecomposition monomorphe est une d\'ecomposition fortement monomorphe de $\mathcal R$.
\end{corollary}

La r\'eciproque du Corollaire \ref{cor:decompfrai-decompfortmono} est fausse. En effet, si nous consid\'erons le contre-exemple pr\'ec\'edent, $\{\{0\},\NN\setminus\{0\}\}$ est une d\'ecomposition fortement monomorphe mais n'est pas une $\mathcal R$-d\'ecomposition monomorphe.

\vspace{2mm}

Il s'ensuit imm\'ediatement et en utilisant la Proposition \ref{propos:monomo-fortmono}
\begin{proposition}
Soit $\mathcal R$ une structure relationnelle. Si $dim_{Fraisse}(\mathcal R)<\infty$ %$\mathcal R$ poss\`ede une d\'ecomposition Fra\"{\i}ss\'e monomorphe ayant un nombre fini d'intervalles
alors $dim_{fort}(\mathcal R)<\infty$ %$\mathcal R$ poss\`ede une d\'ecomposition fortement monomorphe ayant un nombre fini de parties
et donc $dim_{mon}(\mathcal R)<\infty$.% $\mathcal R$ poss\`ede une d\'ecomposition monomorphe ayant un nombre fini de blocs. %De plus, la dimension monomorphe de $\mathcal R$ est le plus petit nombre de $\mathcal R$-intervalles monomorphes infinis d'une $\mathcal R$-d\'ecomposition monomorphe.
\end{proposition}

Nous pouvons alors conclure que si $\mathcal R$ est une structure relationnelle et si nous d\'esignons par
$\mathscr A_1$ l'ensemble des $\mathcal R$-intervalles monomorphes, par
$\mathscr A_2$  l'ensemble des blocs fortement monomorphes de $\mathcal R$
et par $\mathscr A_3$ l'ensemble des blocs monomorphes de $\mathcal R$, alors
$$\mathscr A_1\subseteq \mathscr A_2\subseteq \mathscr A_3.$$

\subsection{Cas d'une structure ordonn\'ee: d\'ecomposition en intervalles}\label{subsec:decomposition en intervalle}

% \subsection{D\'ecomposition en intervalles pour les structures ordonn\'ees.}
        \subsubsection{Intervalle monomorphe}

\begin{definition}\label{def:intervalle monomorphe}
  Soit $\mathcal R:=(E,\leq, (\rho_i)_{i\in I})$ une structure relationnelle ordonn\'ee. Posons $C:=(E,\leq)$.
  Un sous-ensemble $A$ de $E$ est un \emph{intervalle monomorphe}\index{intervalle!monomorphe} de $\mathcal R$ si $A$ est un intervalle %\footnote{Le mot intervalle est utilis\'e ici au sens classique du terme. Dans le cas d'une cha\^{i}ne, c'est un ensemble $A$ d'\'el\'ements de $E$ tels que si $x,y\in A$ et $x\leq z\leq y$ alors $z\in A$.}
  de $C$ tel que pour toutes parties finies $F,F'$ de $E$ de m\^eme cardinalit\'e $n$ $(n\in\mathbb N)$, les structures induites sur $F$ et $F'$ sont isomorphes d\`es que $F\setminus A=F'\setminus A$.
  \end{definition}

%  \noindent L'ensemble vide et les singletons de $E$ sont des intervalles monomorphes.

%  \noindent L'ensemble $E$ est un intervalle monomorphe si et seulement si $\mathcal R$ est monomorphe.
%\vspace{1mm}

  Nous avons les propri\'et\'es suivantes qui d\'ecoulent naturellement de la d\'efinition:%imm\'ediatement le fait suivant.
  \begin{properties}
   Soit $\mathcal R:=(E,\leq, (\rho_i)_{i\in I})$ une structure relationnelle ordonn\'ee. Posons $C:=(E,\leq)$.
   \begin{enumerate}
   \item L'ensemble vide et les singletons de $E$ sont des intervalles monomorphes.
   \item L'ensemble $E$ est un intervalle monomorphe de $\mathcal R$ si et seulement si $\mathcal R$ est monomorphe.
\item  Toute partie d'un intervalle monomorphe de $\mathcal R$ qui est un intervalle pour $C$ est un intervalle monomorphe de $\mathcal R$.
\item toute intersection d'intervalles monomorphes de $\mathcal R$ est un intervalle monomorphe de $\mathcal R$.
\item La r\'eunion de deux intervalles monomorphes dont l'intersection est non vide est un intervalle monomorphe de $\mathcal R$.
\item La r\'eunion d'une famille filtrante pour l'inclusion d'intervalles monomorphes est un intervalle monomorphe.
  \end{enumerate}
  \end{properties}

  \begin{lemma}\label{lem:inter-bloc monomorphe}
  Soit $\mathcal R:=(E,\leq, (\rho_i)_{i\in I})$ une structure relationnelle ordonn\'ee, posons   $C:=(E,\leq)$.
    Tout intervalle monomorphe de $\mathcal R$ est un bloc monomorphe de $\mathcal R$. Inversement, tout bloc monomorphe de $\mathcal R$ qui est un intervalle de $C$ est un intervalle monomorphe de $\mathcal R$.%La r\'eciproque n'est pas vraie.
    \end{lemma}

     \begin{proof}
    %Soit $A$ un intervalle monomorphe de $\mathcal R$, montrons que $A$ est un bloc monomorphe de $\mathcal R$. Soit $n\in\mathbb N$ et soient $B$ et $B'$ deux $n$-parties de $E$ telles que $B\setminus A =B'\setminus A$. Montrons que
La premi\`ere partie d\'ecoule des d\'efinitions \ref{def:bloc monomorphe} et \ref{def:intervalle monomorphe}. Pour la r\'eciproque, toute partie d'une cha\^{\i}ne est un bloc monomorphe mais pas forc\'ement un intervalle de $C$, si nous exigeons que ce soit un intervalle, nous aurons un intervalle monomorphe.
    \end{proof}
\vspace{2mm}

\begin{consequence}\label{consequence}
Toute composante monomorphe de $\mathcal R$ est une union d'intervalles monomorphes de $\mathcal R$.
\end{consequence}

  \begin{lemma}\label{lem:inter-enchain} Soit $\mathcal R:=(E,\leq, (\rho_i)_{i\in I})$ une structure relationnelle ordonn\'ee et soit $A$ un sous-ensemble de $E$. Posons $C:=(E,\leq)$.
  Si $A$ est un intervalle monomorphe de $\mathcal R$ alors $\mathcal R_{\restriction_A}$ est encha\^{i}n\'ee par $C_{\restriction_A}$. Inversement, si $\mathcal R_{\restriction_A}$ est encha\^{i}n\'ee par $C_{\restriction_A}$ et $A$ un intervalle Fra\"{\i}ss\'e de $\mathcal R$ alors $A$ est un intervalle monomorphe de $\mathcal R$.
  \end{lemma}

  \begin{proof}
  Si $A$ un intervalle monomorphe de $\mathcal R$, montrons que $\mathcal R_{\restriction_A}$ est encha\^{i}n\'ee par $C_{\restriction_A}$. Soit $f$ un automorphisme local de $C_{\restriction_A}$ et notons $A_1$ et $A_2$ son domaine et son codomaine respectivement. Comme $A$ est un intervalle monomorphe de $\mathcal R$ et  $A_1\setminus A=A_2\setminus A=\varnothing$, les restrictions de $\mathcal R$ \`a $A_1$ et $A_2$ sont isomorphes. Il s'ensuit que $f$ est un automorphisme local de $\mathcal R_{\restriction_A}$. Par suite $\mathcal R_{\restriction_A}$ est encha\^{i}n\'ee par $C_{\restriction_A}$.

 \vspace{1mm}

 \noindent  Pour l'inverse, $A$ \'etant un intervalle Fra\"{\i}ss\'e de $\mathcal R$, c'est un intervalle de $C$. Soient $F,F'$ deux parties finies de $E$ ayant m\^eme cardinalit\'e $n~(n\in\mathbb N)$ telles que $F\setminus A=F'\setminus A$. Montrons que les restrictions $\mathcal R_{\restriction_F}$ et $\mathcal R_{\restriction_{F'}}$ sont isomorphes. Posons $F_A=F\cap A$ et $F'_A=F'\cap A$, les parties $F_A$ et $F'_A$ ont m\^eme cardinalit\'e, donc il existe un automorphisme local $g$ de $C_{\restriction_A}$ qui envoie $F_A$ sur $F'_A$. $\mathcal R_{\restriction_A}$ \'etant encha\^{i}n\'ee par $C_{\restriction_A}$, $g$ est un automorphisme local de  $\mathcal R_{\restriction_A}$. Comme $A$ est un intervalle Fra\"{\i}ss\'e de $\mathcal R$, en prolongeant $g$ par l'identit\'e sur $E\setminus A$, nous obtenons un automorphisme local, disons $\tilde{g}$, de $\mathcal R$. Cet automorphisme est unique, car $\mathcal R$ est ordonn\'ee et sa restriction \`a $F$ est un automorphisme local de $\mathcal R$ qui envoie $F$ sur $F'$ (car $F\setminus A=F'\setminus A$). Il s'ensuit que $A$ est un intervalle monomorphe.
  \end{proof}

\vspace{2mm}

Nous avons alors

\begin{lemma}\label{lem:equivintmonomorphes}
Soit $\mathcal R:=(E,\leq, (\rho_i)_{i\in I})$ une structure relationnelle ordonn\'ee et soit $A$ un sous-ensemble de $E$. Les assertions suivantes sont \'equivalentes:
\begin{enumerate}
\item[(i)] $A$  est un intervalle monomorphe,
\item[(ii)] $A$ est un $\mathcal R$-intervalle monomorphe,
\item[(iii)] $A$ est un bloc fortement monomorphe.
\end{enumerate}
\end{lemma}

\begin{proof}
Montrons $(i)\Rightarrow (ii)$.  Soit $J$ un intervalle monomorphe de $\mathcal R$, donc $\mathcal R_J$ est monomorphe. Montrons que $J$ est un $\mathcal R$-intervalle. Soit $f$ un automorphisme local de $\mathcal R_{\restriction_J}$ et soient $J_1$ et $J_2$ son domaine et codomaine. Montrons que $f$ prolong\'ee par l'identit\'e sur $E\setminus J$ est un automorphisme local de $\mathcal R$. Posons $A=J_1\cup (E\setminus J)$ et $A'=J_2\cup (E\setminus J)$. Les ensembles $A$ et $A'$ sont de m\^eme cardinalit\'e et $A\setminus J=A'\setminus J$, donc les restrictions de $\mathcal R$ \`a $A$ et $A'$ sont isomorphes car $J$ est un intervalle monomorphe de $\mathcal R$. Notons par $\tilde{f}$ cet isomorphisme. Comme $\mathcal R$ est ordonn\'ee $\tilde{f}$ est unique. $\tilde{f}$ est un isomorphisme local de $C$ et $J$ un intervalle de $C$ donc $\tilde{f}$ envoie les \'el\'ements de $A\setminus J_1=A\setminus J$ sur $A'\setminus J_2=A'\setminus J$, et $J_1$ sur $J_2$. Comme $A\setminus J_1= A'\setminus J_2=A\setminus J$, donc sur $A\setminus J$,  $\tilde{f}$ co\"{\i}ncide avec l'identit\'e et  $\tilde{f}_{\restriction_{J_1}}=f$. D'o\`u $\tilde{f}$ est $f$ prolong\'ee par l'identit\'e sur $E\setminus J$.

L'implication $(ii)\Rightarrow (i)$ vient du fait qu'un $\mathcal R$-intervalle est un intervalle de la cha\^{i}ne $(E,\leq)$, la monomorphie implique le reste. On conclut que les assertions $(i)$ et $(ii)$ sont \'equivalentes.

L'implication $(ii)\Rightarrow (iii)$ est donn\'ee par le Lemme \ref{lem:monomo-forte.monom}. La r\'eciproque $(iii)\Rightarrow (ii)$ est v\'erifi\'ee dans le cas des structures ordonn\'ees; en effet, si $A$ est un bloc fortement monomorphe alors $\mathcal R_{\restriction_A}$ est monomorphe. Soit $f$ un automorphisme local de $\mathcal R_{\restriction_A}$ de domaine $J$ et co-domaine $J'$. La structure $\mathcal R$ \'etant ordonn\'ee, $f$ est unique. Comme $A$ est un bloc fortement monomorphe, alors l'isomorphisme $f$ de $\mathcal R_{\restriction_J}$ sur $\mathcal R_{\restriction_{J'}}$, prolong\'e par l'identit\'e sur $E\setminus A$ est un automorphisme local de $\mathcal R$. D'où $A$ est un $\mathcal R$-intervalle monomorphe.
\end{proof}

\vspace{2mm}

Remarquons qu'un $\mathcal R$-intervalle n'est pas forc\'ement un intervalle monomorphe de $\mathcal R$. En effet, consid\'erons le contre exemple suivant:

\medskip

\noindent  Soit $\mathcal R:=(E,\leq,\rho)$ telle que $E=\{x,y,z,t\}$, $\leq$  l'ordre lin\'eaire  donn\'e par: $x\leq y\leq z\leq t$ et $\rho$ la relation binaire  donn\'ee par:
  $\rho=\{(x,y),(x,z),(x,t), (z,y)\}$. Le sous-ensemble $\{y,z,t\}$ est un intervalle Fra\"{\i}ss\'e de $\mathcal R$ mais pas un intervalle monomorphe.

\vspace{2mm}

  Dans le cas particulier des structures binaires ordonn\'ees, nous avons le r\'esultat suivant:

 \begin{proposition}\label{prop:comp-intervalle mono}
  Soit $\mathcal R:=(E,\leq ,\rho _{1},\dots,\rho _{k})$ une structure binaire ordonn\'ee de type $k$ et posons $C:=(E,\leq)$.
Toute composante monomorphe $B_{\mathcal R}$ de  $\mathcal R$ telle que $\vert B_{\mathcal R} \vert\neq 2$ est un intervalle monomorphe maximal de $\mathcal R$ et tout intervalle monomorphe maximal non trivial $I_{\mathcal R}$ de  $\mathcal R$ est une composante monomorphe de $\mathcal R$.
  \end{proposition}

  \begin{proof} Soit $B_{\mathcal R}$ une composante monomorphe de $\mathcal R$. Si $B_{\mathcal R}$ a un seul \'el\'ement, c'est un intervalle monomorphe. Supposons que $B_{\mathcal R}$ a, au moins trois \'el\'ements. Si $B_{\mathcal R}$ n'est pas un intervalle monomorphe alors d'apr\`es le Lemme \ref{lem:inter-bloc monomorphe}, $B_{\mathcal R}$ n'est pas un intervalle de $C$. Soit $J$ le plus petit intervalle de $C$ contenant $B_{\mathcal R}$.
Montrons que $J$ est un intervalle monomorphe, ce qui, en vertu du Lemme \ref{lem:inter-bloc monomorphe}, constitue une contradiction avec le fait que $B_{\mathcal R}$ est une composante monomorphe.

Soient $x,y,z\in B_{\mathcal R}$ et $a\in J\setminus B_{\mathcal R}$. Nous pouvons supposer, sans perte de g\'en\'eralit\'e, que $x<a<y<z$. En utilisant le fait que les restrictions de $\mathcal R$ \`a $\{x,a\}$, $\{y,a\}$ et $\{z,a\}$ sont isomorphes et de m\^eme les restrictions de $\mathcal R$ \`a $\{x,y,a\}$, $\{x,z,a\}$ et $\{y,z,a\}$ sont isomorphes, nous obtenons que la restriction de $\mathcal R$ \`a une partie \`a deux \'el\'ements tels que l'un soit dans $J\setminus B_{\mathcal R}$ et l'autre dans $B_{\mathcal R}$ est isomorphe \`a la restriction de $\mathcal R$ \`a une partie \`a deux \'el\'ements de $B_{\mathcal R}$. Si $b$ est un autre \'el\'ement de $J\setminus B_{\mathcal R}$, nous montrons de m\^eme que la restriction de $\mathcal R$ \`a $\{a,b\}$ est isomorphe \`a la restriction de $\mathcal R$ \`a une partie \`a deux \'el\'ements de $B_{\mathcal R}$. Donc $J$ est $2$-monomorphe.

Montrons que $J$ est un intervalle de $\mathcal R$. Soit $c\notin J$, en utilisant le fait que les restrictions de $\mathcal R$ \`a $\{x,a,c\}$ et $\{y,a,c\}$ sont isomorphes on d\'eduit que $J$ est un intervalle de $\mathcal R$. D'apr\`es le Corollaire \ref{cor:encha-monomorphe} la restriction de $\mathcal R$ \`a $J$ est ench\^{i}n\'ee par la restriction de $C$ \`a $J$. D'apr\`es le Lemme \ref{lem:inter-enchain}, $J$ est un intervalle monomorphe comme requis. On d\'eduit alors que  $B_{\mathcal R}$ est une composante monomorphe.

\medskip

Pour la r\'eciproque, soit $A$ un intervalle monomorphe maximal non trivial. Si $A$ a au moins trois \'el\'ements alors $A$ est contenu dans une composante monomorphe qui serait, d'apr\`es ce qui pr\'ec\`ede, un intervalle monomorphe maximal donc \'egal \`a $A$. De m\^eme si $A$ a deux \'el\'ements alors $A$ est une composante monomorphe car sinon, $A$ serait contenu dans une composante monomorphe d'au moins trois \'el\'ements qui serait un intervalle monomorphe maximal, ce qui contredirait la maximalit\'e de $A$.
  \end{proof}

  \vspace{2mm}

  Une composante monomorphe \`a deux \'el\'ements n'est pas forc\'ement un intervalle monomorphe comme nous pouvons le voir dans cet exemple.

  Soit $\mathcal R:=(E,\leq,\rho)$ où $E=\{x,y,z,t\}$ avec\\ $$x<y<z<t \text{ et }\rho=\{(x,t),(y,t),(z,t),(z,x),(t,t)\}.$$
  Alors $A:=\{x,z\}$ est une composante monomorphe mais n'est pas un intervalle monomorphe. Les seuls intervalles monomorphes de cette structure sont le vide et les singletons.

        \subsubsection{D\'ecomposition en intervalles}
\begin{definition}\label{def:decomposition en intervalles}
        Soit $\mathcal R:=(E,\leq, (\rho_i)_{i\in I})$ une structure relationnelle ordonn\'ee et posons  $C:=(E,\leq)$.
 Une \emph{d\'ecomposition en intervalles}\index{decomposition@d\'ecomposition!en intervalles} de $\mathcal R$ est une partition $\mathscr P$ de $E$ en intervalles monomorphes de $\mathcal R$. En d'autres termes, c'est une partition $\mathscr P$ de $E$ en intervalles $J$ de la cha\^{i}ne $C:= (E, \leq)$ telle que pour tout entier $n$ et toute paire $A,~A'$ de $n$-parties de $E$, les structures induites sur $A$ et $A'$ sont isomorphes d\`es que les traces\index{trace} $A\cap J$ et  $A'\cap J$ sont de m\^eme cardinalit\'e pour tout intervalle $J$ de $\mathscr P$.
   \end{definition}

   \vspace{1mm}

   Par  exemple, si $\mathcal R$ est  la bicha\^{i}ne $(E, \leq, \leq')$, $\mathscr P$ est une  d\'ecomposition en intervalles de $E$ si et seulement si chaque bloc $J$ de $\mathscr P$ est un intervalle pour chacun des deux ordres qui co\"{\i}ncident ou sont oppos\'es sur  $J$ (voir \cite{mont-pou}).\\

Remarquons que les intervalles monomorphes maximaux de $\mathcal R$ forment une d\'ecomposition en intervalles de $\mathcal R$ dont toute autre d\'ecomposition en intervalles est plus fine. Notons par $dim_{int}(\mathcal R)$ le nombre d'intervalles monomorphes maximaux de $\mathcal R$, si ce nombre n'est pas fini, nous \'ecrivons $dim_{int}(\mathcal R)=\infty$.\\

%Dans le cas o\`u $\mathcal R$ est une structure ordonn\'ee, nous avons:
Nous avons naturellement \`a partir du Lemme \ref{lem:equivintmonomorphes}

\begin{lemma}
Soit $\mathcal R$ une structure ordonn\'ee de base $E$ et $\mathscr P$ une partition de $E$. Les assertions suivantes sont \'equivalentes:
\begin{enumerate}
\item[(i)] $\mathscr P$  est une d\'ecomposition en intervalles de $\mathcal R$,
\item[(ii)] $\mathscr P$  est une $\mathcal R$-d\'ecomposition monomorphe,
\item[(iii)] $\mathscr P$  est une d\'ecomposition fortement monomorphe de $\mathcal R$.
\end{enumerate}
\end{lemma}

    A partir du Lemme \ref{lem:inter-bloc monomorphe}, nous avons:

\begin{lemma}\label{lem:dec intervalle=dec.mono}
    Toute d\'ecomposition en intervalles d'une structure ordonn\'ee $\mathcal R$ est une d\'ecomposition monomorphe de $\mathcal R$. La r\'eciproque n'est pas vraie.
    \end{lemma}

Cependant, en appliquant le Th\'eor\`eme \ref{thm:chainable}, nous avons si $\mathcal R:=(E, \leq, (\rho_{i})_{i\in I})$:

\begin{lemma}\label{lem:comp.infini=intervalle}
Toute composante monomorphe infinie $A$ de $\mathcal R$ est un intervalle de la cha\^{i}ne $C:=(E,\leq)$, donc un intervalle monomorphe maximal (pour l'inclusion) de $\mathcal R$.
\end{lemma}

\begin{proof}
Soit $A$ une composante monomorphe infinie de $\mathcal R$. D'apr\`es le Th\'eor\`eme \ref{thm:chainable}, $A$ est un bloc fortement monomorphe de $\mathcal R$. D'apr\`es le Lemme \ref{lem:equivintmonomorphes}, $A$ est un intervalle monomorphe.
Il s'ensuit que $A$ est un intervalle de $C$. %et donc $A$ est un intervalle monomorphe de $\mathcal R$ en vertu du Lemme \ref{lem:inter-bloc monomorphe}.
Il est \'evident que $A$ est maximal puisque c'est une composante monomorphe.
    \end{proof}

\vspace{2mm}

Ceci  entra\^{\i}ne:

\begin{theorem} \label{thm:intervalsecomposition} Si $\mathcal R$ est une structure relationnelle ordonn\'ee alors $dim_{int}(\mathcal R)<\infty$ %$\mathcal R$ poss\`ede une d\'ecomposition  en intervalles ayant un nombre fini d'intervalles
si et seulement si $dim_{mon}(\mathcal R)<\infty$; %elle poss\`ede une d\'ecomposition monomorphe ayant un nombre fini de blocs;
en outre, $dim_{mon}^{\infty}(\mathcal R)$ est %le nombre de composantes monomorphes infinies  %dimension monomorphe de $\mathcal R$ est
le plus petit nombre d'intervalles infinis d'une d\'ecomposition en intervalles de $\mathcal R$.
\end{theorem}

\begin{proof}
Si $dim_{int}(\mathcal R)<\infty$ %$\mathcal R$ poss\`ede une d\'ecomposition en intervalles ayant un nombre fini d'intervalles
alors, d'apr\`es le Lemme \ref{lem:dec intervalle=dec.mono}, $dim_{mon}(\mathcal R)<\infty$. %elle poss\`ede une d\'ecomposition monomorphe ayant un nombre fini de blocs.
Inversement, si $dim_{mon}(\mathcal R)<\infty$, %$\mathcal R$ poss\`ede une d\'ecomposition monomorphe ayant un nombre fini de blocs,
alors d'apr\`es le Lemme \ref{lem:comp.infini=intervalle}, les composantes infinies sont des intervalles. En partitionnant les \'el\'ements de la r\'eunion des composantes finies (qui sont en nombre fini) en singletons, nous obtenons une d\'ecomposition en intervalles ayant un nombre fini d'intervalles, donc $dim_{int}(\mathcal R)<\infty$. La deuxi\`eme partie vient directement du fait que le nombre de composantes monomorphes infinies est le plus petit nombre de blocs monomorphes infinis dans une d\'ecomposition monomorphe de $\mathcal R$ et du Lemme \ref{lem:comp.infini=intervalle}.
\end{proof}

\vspace{2mm}

La propri\'et\'e suivante permet de tester l'existence d'une d\'ecomposition monomorphe ayant un nombre fini de blocs sur les relations formant $\mathcal R$.
\begin{proposition}\label{prop:decomp-finie}
 Soit $k\in \mathbb N$ et soit $\mathcal R:= (E, \leq, (\rho_i)_{i<k})$ une structure relationnelle ordonn\'ee.  $\mathcal R$ admet une d\'ecomposition monomorphe finie si et seulement si toute structure $\mathcal R_i:=(E, \leq, \rho_i)$ $(i< k)$ poss\`ede une telle d\'ecomposition.
\end{proposition}

\begin{proof}
La condition n\'ecessaire est \'evidente car toute d\'ecomposition monomorphe de $\mathcal R$ est une d\'ecomposition
monomorphe de $\mathcal R_i$ pour tout $i<k$.
Pour la condition suffisante, remarquons d'abord que si $A$ est un bloc monomorphe pour toutes les structures $\mathcal R_i,~i<k$ alors
$A$ est un bloc monomorphe pour $\mathcal R$. Supposons que $\mathcal R_i$ poss\`ede une d\'ecomposition monomorphe
finie pour tout $i<k$, alors chacune de ces $k$ d\'ecompositions est une partition finie de $E$.
Comme l'ensemble de toutes les partitions d'un ensemble est un treillis complet, l'ensemble de ces $k$
d\'ecompositions poss\`ede un infimum (ie, une partition qui est la moins fine de toutes les partitions plus fines que les $k$
d\'ecompositions), notons par $P$ cette partition ($P$ est l'infimum des $k$ partitions). $P$ est forc\'ement finie car les $k$ d\'ecompositions sont finies.
Chaque bloc $A$ de la partition $P$ est inclus dans un bloc de la d\'ecomposition de $\mathcal R_i$ pour tout $i<k$, donc $A$
est un bloc monomorphe pour $\mathcal R_i,~i<k$, d'o\`u $A$ est un bloc monomorphe pour $\mathcal R$.
Il s'ensuit que $P$ est une d\'ecomposition monomorphe finie de $\mathcal R$.
\end{proof}

 \subsection{Une relation d'\'equivalence pour la d\'ecomposition monomorphe}\label{subsection:canonical}
 Nous allons retrouver la notion de d\'ecomposition monomorphe en d\'efinissant la relation suivante (voir \cite{oudrar-pouzet2014}).

\vspace{2mm}

 Soit $\mathcal R:= (V, (\rho_i)_{i\in I})$ une structure relationnelle. Soient $x$ et $y$ deux \'el\'ements de $V$ et $F$ une partie de $V\setminus \{x,y\}$. Disons que $x$ et $y$ sont \emph{$F$-\'equivalents} et notons $x\simeq_{F, \mathcal R} y$  si les restrictions de
$\mathcal R$ \`a $\{x\}\cup F$ et \`a $\{y\}\cup F$ sont isomorphes.

\vspace{2mm}

Soit $k$ un entier, disons que $x$ et $y$ sont \emph{$k$-\'equivalents}\index{relation!de $k$-\'equivalence} et notons $x\simeq_{k, \mathcal R} y$ si $x\simeq_{F, \mathcal R} y$ pour toute partie $F$ \`a $k$ \'el\'ements de $V\setminus \{x, y\}$. Nous disons qu'ils sont \emph{$(\leq k)$-\'equivalents}\index{relation!de $\leq k$-\'equivalence} et nous notons   $x\simeq_{\leq k, \mathcal R}y$ s'ils sont $k'$-\'equivalents pour tout $k'\leq k$. Enfin, nous disons qu'ils sont \emph{\'equivalents}, ce que nous notons $x\simeq_{\mathcal R} y$,  s'il sont $k$-\'equivalents pour tout entier $k$.

\vspace{1mm}

Par exemple, la $0$-\'equivalence est une relation d'\'equivalence, deux sommets $x, y$ \'etant dans une m\^eme classe  si les restrictions de $\mathcal R$ \`a $x$ et $y$ sont isomorphes.

\begin{lemma}\label{lem:equivalence}
Les relations $\simeq_{k, \mathcal R}$, $\simeq_{\leq k, \mathcal R}$ et $\simeq_{\mathcal R}$ sont des relations d'\'equivalence sur $V$.
\end{lemma}

\begin{proof}
Il suffit de voir que $\simeq_{k, \mathcal R}$ est une relation d'\'equivalence. Pour cela, il suffit de v\'erifier que cette relation est transitive. Soient $ x,y, z\in V$ avec $x\simeq_{k, \mathcal R} y$ et $y\simeq_{k,\mathcal R} z$. V\'erifions que $x\simeq_{k,\mathcal R} z$. Nous pouvons supposer ces trois \'el\'ements deux \`a deux distincts. Soit $F$ une partie \`a $k$ \'el\'ements de  $V\setminus \{x,z\}$. Nous avons deux cas.\\ \textbf{Cas 1)} $y\not \in F$. Dans ce cas, nous avons respectivement $\mathcal R_{\restriction F\cup\{x\}}$ isomorphe \`a $\mathcal R_{\restriction F\cup\{y\}}$ et $\mathcal R_{\restriction F\cup\{y\}}$ isomorphe \`a $\mathcal R_{\restriction F\cup\{z\}}$. Donc $\mathcal R_{\restriction F\cup\{x\}}$ est isomorphe \`a $\mathcal R_{\restriction F\cup\{z\}}$; ainsi $x\simeq_{k,\mathcal R} z$. \\ \textbf{Cas 2)} $y\in F$.
Posons $F':= (F\setminus \{y\})\cup \{z\}$. Comme $x\simeq_{k,\mathcal R} y$ nous avons  $\mathcal R_{\restriction F'\cup\{x\}}$ isomorphe \`a $\mathcal R_{\restriction F'\cup\{y\}}=\mathcal R_{\restriction F\cup \{z\}}$. De m\^eme, si nous posons $F'':= (F\setminus \{y\})\cup \{x\}$ alors,  comme $y\simeq_{k,\mathcal R} z$ nous avons  $\mathcal R_{\restriction F\cup \{x\}}=\mathcal R_{\restriction F''\cup\{y\}}$ isomorphe \`a $\mathcal R_{\restriction F''\cup\{z\}}$. Comme $F'\cup\{x\}=F''\cup\{z\}$ le r\'esultat s'ensuit.
\end{proof}

\vspace{2mm}

Nous avons naturellement, d'apr\`es la d\'efinition:

\begin{fact}\label{fact:classe-monomorphe}
La restriction de $\mathcal R$ \`a toute classe de $k$-\'equivalence est $k+1$-monomorphe.
\end{fact}

\begin{lemma}\label{lem:equi-blocmono}
 Chaque classe d'\'equivalence de $\simeq _{\mathcal R}$ est un bloc monomorphe de $\mathcal R$ et tout bloc monomorphe est inclus dans une classe d'\'equivalence.
\end{lemma}

\begin{proof}
Soit  $C$ une classe d'\'equivalence de $\simeq_{\mathcal R}$. Nous montrons que $C$ est un bloc monomorphe.  Soit   $A, A'\in [V]^h$  tels que  $A\setminus C=A'\setminus C$. Nous devons montrer  que $\mathcal R_{\restriction A}$ et  $\mathcal R_{\restriction A'}$  sont isomorphes. Soit $\ell := \vert A\setminus A'\vert$. Si  $\ell=0$, $A=A'$, il n'y a rien \`a montrer. Si $\ell =1$,  alors $A= \{x\}\cup (A\cap A')$ et $A'=\{y\}\cup (A\cap A')$, avec $x,y\in C$; dans ce cas $\mathcal R_{\restriction A}$ et $\mathcal R_{\restriction A'}$  sont isomorphes puisque $x\simeq_{\mathcal R} y$. Si $\ell >1$, posons $K:=A \setminus C$ et $k:= \vert K\vert$. Nous pouvons trouver une suite de $(h-k)$ sous-ensembles d'\'el\'ements de $C$, disons $A_0, \dots A_i, \dots, A_{\ell-1}$ tels que $A_0= A\cap C$,  $A_{\ell-1}= A'\cap C$ et la diff\'erence sym\'etrique %\footnote{La diff\'erence sym\'etrique de deux ensembles $A$ et $B$ est l'ensemble $(A\setminus B)\cup(B\setminus A)$.}
de $A_i$ et $A_{i+1}$ est de cardinalit\'e $2$. %$0$ or $1$. From the case $\ell=1$
Nous avons $\mathcal R_{\restriction  A_i\cup K}= \mathcal R_{\restriction A_{i+1}\cup K}$ pour $i<\ell-1$. Donc $\mathcal R_{\restriction A}$ et $\mathcal R_{\restriction A'}$ sont isomorphes.

  Puisque les classes d'\'equivalence sont des blocs monomorphes elles forment une d\'ecomposition  monomorphe. Trivialement, les \'el\'ements d'un bloc monomorphe sont deux \`a deux \'equivalents pour $\simeq_{\mathcal R}$ donc contenus dans une classe d'\'equivalence.
\end{proof}

\vspace{2mm}

D'apr\`es le Lemme  \ref{lem:equi-blocmono} nous obtenons:

\begin{proposition}\label{prop: equiv=decomp}
les classes d'\'equivalence  de $\simeq_{\mathcal R}$ forment une d\'ecomposition de  $\mathcal R$ en blocs monomorphes et toute autre d\'ecomposition en blocs monomorphes est plus fine qu'elle. Ainsi la d\'ecomposition de $\mathcal R$ en classes d'\'equivalence de $\simeq _{\mathcal R}$
co\"{\i}ncide avec la d\'ecomposition de $\mathcal R$ en composantes  monomorphes.
\end{proposition}

Certaines propri\'et\'es mentionn\'ees dans la sous-section \ref{subsec:decomposition monomorphe}, exemple Th\'eor\`eme \ref{compactness} %et Lemme \ref{lemma.minimalGrowthRate},
se retrouvent ais\'ement \`a partir de la Proposition \ref{prop: equiv=decomp}.

\vspace{2mm}

Soit $\mathcal R$ une structure relationnelle de base $E$. Soient  $x$, $y$ deux \'el\'ements non  \'equivalents de $E$ ($x\not\simeq_{\mathcal R} y$). Soit $F$ une partie finie minimale qui en t\'emoigne (en particulier $F\subset E\setminus\{x,y\}$ et les restrictions de $\mathcal R$ \`a $F\cup \{x\}$ et $F\cup \{y\}$ ne sont pas isomorphes) et soit $F':= F\cup \{x,y\}$.

\begin{lemma}\label{intervalle}
La restriction $R_{\restriction_{ F'}}$ ne contient pas d'intervalle Fra\"{\i}ss\'e $A$ contenant $\{x, y\} $ et distinct de $F'$.
\end{lemma}
\begin{proof}
Supposons que $R_{\restriction_{ F'}}$ contient un tel intervalle $A$,  posons $F_1:= A \setminus \{x,y\}$. Par minimalit\'e de $F$, $F_1\cup \{x\}$  est isomorphe \`a $F_1\cup \{y\}$.  Soit $f$ cet isomorphisme;  comme $A$ est un intervalle Fra\"{\i}ss\'e, $f$ prolong\'e par l'identit\'e sur $F\setminus A$ est un isomorphisme de  $F\cup \{x\}$  sur $F\cup \{y\}$. Contradiction.
\end{proof}

 \subsection{Hypomorphie et \'equivalence}\label{subsection: hypomorphy}
Soient  $\mathcal R$ et $\mathcal R'$ deux structures relationnelles sur le m\^eme ensemble $V$ et $k$ un entier. Nous disons, d'apr\`es Fra\"{\i}ss\'e\index{Fraisse@Fra\"{\i}ss\'e} et Lopez\index{Lopez} \cite{fraisse-lopez},  que $\mathcal R$ et $\mathcal R'$ sont \emph{$k$-hypomorphes}\index{relation!d'hypomorphie} si les restrictions $\mathcal R_{\restriction A}$ et $\mathcal R'_{\restriction A}$ sont isomorphes pour toutes $k$-partie $A$ de $V$. Soit  $\mathcal R$ une structure relationnelle de base  $V$ et $x,y\in V$. Posons $V':= V\setminus \{x,y\}$. En identifiant $x$ et $y$ \`a un \'el\'ement $z$ (exemple $z:= \{x,y\}$) nous obtenons deux structures $\mathcal R(x)$ et $\mathcal R(y)$ sur $V' \cup \{z\}$. Formellement,  $\mathcal R(x)$ est telle que l'application qui fixe $V'$ point par point et envoie $x$ sur $z$ est un isomorphisme de $\mathcal R_{\restriction  {V\setminus \{y\}}}$ sur $\mathcal R(x)$. De m\^eme pour $\mathcal R_{\restriction {V\setminus \{x\}}}$ et $\mathcal R(y)$.

Puisque $\mathcal R(x)$ et $\mathcal R(y)$ co\"{\i}ncident sur $V'$ nous avons imm\'ediatement:

\begin{lemma}\label{lem:iteration}
Soit  $x,y\in  V$ et $k\in \mathbb N$ alors  $x\simeq _{k,\mathcal  R} y$ si et seulement si $\mathcal R(x)$ et $\mathcal R(y)$ sont $(k+1)$-hypomorphes.
\end{lemma}

Il a \'et\'e montr\'e que deux structures relationnelles qui sont $(k+1)$-hypomorphes  sont $k$-hypomorphes d\`es que leur domaine poss\`ede au moins $2k+1$ \'el\'ements (voir \cite {pouzet 79} Corollaire 2.3.2 et aussi \cite{pouzet76} où ceci s'ensuit \`a partir d'une propri\'et\'e des matrices d'incidence due \`a Gottlieb\index{Gottlieb} \cite{gottlieb} et Kantor\index{Kantor} \cite{kantor}. Pour le rappel du r\'esultat de Gottlieb et Kantor voir \cite{Dam-Lo-Po-Ka}).

\begin{corollary}\label{lem:gottlieb-kantor}
Les relations d'\'equivalence $\simeq_{k,\mathcal R}$ et $\simeq_{\leq k,\mathcal  R}$ co\"{\i}ncident d\`es que $\vert V(R)\vert \geq 2k+1$.
 \end{corollary}

La notion d'hypomorphie a \'et\'e assez bien \'etudi\'ee, particuli\`erement pour les structures binaires. Disons que deux structures $\mathcal R$ and $\mathcal R'$ sur le m\^eme ensemble $V$ sont \emph{hypomorphes}\index{structure relationnelle!hypomorphe} si elles sont $k$-hypomorphes pour tout $k< \vert V\vert$ et que $\mathcal R$ est \emph{reconstructible}\index{structure relationnelle!reconstructible} si $\mathcal R$ est isomorphe \`a toute structure $\mathcal R'$ qui lui est hypomorphe.
Lopez  a montr\'e que toute structure binaire finie d\'efinie sur au moins $7$ \'el\'ements  est  reconstructible (\cite{lopez 72} pour les relations binaires, \cite{lopez 78} pour les structures binaires, voir l'expos\'e complet dans  \cite{fraisse-lopez}). De ce fait  avec le  Lemme \ref{lem:iteration}  il s'ensuit que pour une structure binaire $\mathcal R$ les relations d'\'equivalence $\simeq _{\leq 6, \mathcal R}$ et $\simeq _{\mathcal R}$ co\"{\i}ncident. Nous montrons dans la sous-section \ref{subsec:equiv-struct-binaire} qu'en fait les relations d'\'equivalence $\simeq _{\leq 3, \mathcal R}$ et $\simeq _{\mathcal R}$ co\"{\i}ncident pour une structure binaire $\mathcal R$. Ille\index{Ille} (1992) \cite{ille 92} a montr\'e que les structures ordonn\'ees finies dont le maximum des arit\'es est $m$ et de cardinalit\'e au moins un entier $s(m)$  sont reconstructibles. Ainsi, de mani\`ere similaire, il s'ensuit que pour une structure ordonn\'ee $\mathcal R$  d'arit\'e au plus $m$, il existe un entier $i(m)$  ($i(m)\leq s(m)$) tel que les relations d'\'equivalence $\simeq _{\leq i(m),\mathcal R}$ et $\simeq _{\mathcal R}$ co\"{\i}ncident.

\vspace{1mm}

Le seuil pour une structure binaire sera examin\'e dans le chapitre \ref{chap:graphe ordonne}.

\vspace{1mm}

 Pouzet\index{Pouzet} \cite {pouzet 79} a montr\'e qu'il n'y a pas de seuil de reconstructibilit\'e pour les relations ternaires. Il a construit, pour tout entier $m$ une paire de relations ternaires $\mathcal S$ et $\mathcal S'$ sur un ensemble $W$ de cardinalit\'e $2m$, qui sont hypomorphes   mais ne sont pas isomorphes. Dans son exemple $W$ contient un \'el\'ement $z$ (en fait deux) tel que $\mathcal S$ et $\mathcal S'$ co\"{\i}ncident sur $W\setminus \{z\}$. Ainsi nous pouvons dupliquer $z$ et construire une relation ternaire $\mathcal R$ sur $2m+1$ \'el\'ements ayant deux \'el\'ements $x,y$ tels que $x\not \simeq_{\mathcal R}y$ mais $x\simeq_{\leq 2m-1, \mathcal R}y$.

\begin{problem} Existence d'une structure relationnelle $\mathcal R$ de  signature finie telle que $\simeq_{\mathcal R}$ a une infinit\'e de classes d'\'equivalence et $\simeq_{k, \mathcal R}$ n'en a qu'un nombre fini pour tout entier  $k$.
 \end{problem}

 \subsection{Le r\^ole du belordre}\label{subsection:wqo}
% Dans cette partie, nous montrons une extension du Lemme \ref{lem:reduction} et de la Proposition \ref{cor:minimal} (Lemme \ref{lem:reductionbis} et Th\'eor\`eme \ref{thm:minimal1}).% et la premi\`ere partie du Th\'eor\`eme  \ref{thm:poly-expo1} (Theorem\ref{thm: polynomial-interval}).  Pour se faire, nous rappelons la notion de d\'ecomposition monomorphe et les r\'esultats principaux de \cite{P-T-2013}. %(subsection \ref{subsection:monomorphic}. One of ours tools is the notion of well-quasi-order (subsection \ref{subsection:wqo}).
 \vspace{1mm}

Nous rappelons la d\'efinition suivante (voir \cite{P-T-2013}).
\vspace{2mm}

Soit $\mathcal R:=(E,(\rho_i)_{i\in I})$ une structure relationnelle.
 Soit $\mathscr{P}:= (E_x)_{x\in X}$ une partition de l'ensemble  $E$ en parties disjointes. Le \emph{spectre}\index{ensemble!spectre d'un -}  d'un sous-ensemble fini $A$  de $E$,  relativement \`a la partition $\mathscr P$, est la famille $d(A):= (d_x(A))_{x\in X}$ o\`u $d_x(A):= \vert A\cap E_x\vert$.
 \vspace{1mm}

 %Si $\mathscr{P}$ est form\'ee des blocs d'une d\'ecomposition monomorphe de $\mathcal R$ alors
 Nous avons de mani\`ere \'evidente:

 \begin{lemma}\label{lem:spectre}
 Si deux sous-ensembles $A$ et $B$ ont le m\^eme spectre relativement \`a une m\^eme partition en blocs monomorphes de $E$, alors les restrictions de $\mathcal R$ \`a $A$ et $B$ sont isomorphes. La r\'eciproque n'est pas vraie.
 \end{lemma}
% pour la r\'eciproque penser \`a une partition de la chaine. toutes les parties de m\^eme cardinalit\'e sont isomorphes mais n'ont pas toutes le m\^eme spectre.

\begin{proposition}\label{prop-wqo-ages}
L'\^age d'une structure relationnelle admettant une d\'ecomposition monomorphe finie est h\'er\'editairement belordonn\'e.
\end{proposition}

\begin{proof}
Soient $\mathcal R$ une telle structure relationnelle, $E:= V(\mathcal R)$, $\mathscr P:= (E_x)_{x\in X}$ une d\'ecomposition monomorphe finie de $\mathcal R$, $l:= \vert X\vert$. Pour tout $x\in X$, si $E_x$ est infini nous posons $C_{x}:=(\mathbb N, \leq)$ et si $E_x$ est fini, alors $C_{x}:=(\{0, \dots, \vert E_x\vert\},\leq)$ ou $\leq$ est l'ordre naturel. Soit $C$ le produit direct  $\Pi_{x\in X} C_x$ des cha\^{i}nes $C_x,~x\in X$. Soient $A$ et $A'$ deux sous-ensembles de $E$. A partir du Lemme \ref{lem:spectre}, si $d(A)\leq d(A')$ alors $\mathcal R_{\restriction _A} \leq\mathcal R_{\restriction_ {A'}}$. D'o\`u $\mathcal A(\mathcal R)$ ordonn\'e par abritement peut \^etre consid\'er\'e comme l'image surjective de $C$ par une application croissante. D'apr\`es le lemme\footnote{Le lemme de Dickson dit la chose suivante: Pour $n\geq 1$, toute partie non vide de $\mathbb N^n$ muni de l'ordre produit,  a un nombre fini d'\'el\'ements minimaux.}  de Dickson, $C$ est belordonn\'e, d'o\`u  $\mathcal A(\mathcal R)$ est belordonn\'e.  Pour montrer que $\mathcal A(\mathcal R)$ est h\'er\'editairement belordonn\'e,  nous utilisons le Th\'eor\`eme \ref{thm:charcfinitemonomorph} et le th\'eor\`eme d'Higman sur les mots (voir Th\'eor\`eme \ref{theo:higman}) \cite{higman52}. En effet, $\mathcal R$ est libre interpr\'etable en une structure de la forme $\mathcal S:= (E, \leq, (u_1,  \dots, u_{\ell}))$ o\`u $\leq$ est un ordre lin\'eaire et $u_1, \dots, u_{\ell}$ un nombre fini de relations unaires formant une partition de $E$ en intervalles de $(E, \leq)$. Soit $P$ un ensemble belordonn\'e. Montrons que $\mathcal A(\mathcal R). P$, l'ensemble des couples $(\mathcal T,f)$ avec $\mathcal T\in \mathcal A(\mathcal R)$ et $f:V(\mathcal T)\rightarrow P$, est belordonn\'e (voir section \ref{subsection:heredit-belordonne} du chapitre \ref{chap:generalite}). Pour cela, il suffit de montrer que $\mathcal A(\mathcal S). P$ est belordonn\'e. Nous utilisons le th\'eor\`eme d'Higman. Soit $L$ l'anticha\^{i}ne \`a $\ell$ \'el\'ements $\{1, \dots, \ell\}$ et $Q$ le produit direct $P\times L$. Comme produit de deux belordres (dans ce cas une somme directe de $\ell$ copies du belordre $P$) $Q$ est belordonn\'e. Il s'ensuit, d'apr\`es le th\'eor\`eme d'Higman (Th\'eor\`eme \ref{theo:higman}), que l'ensemble $Q^{\star}$ des mots sur $Q$ est belordonn\'e. Pour conclure, montrons que $Q^{\star}$ est l'image de $\mathcal A(\mathcal S). P$ par un plongement et par cons\'equent, $\mathcal A(\mathcal S). P$ est belordonn\'e.

Consid\'erons l'application $h:\mathcal A(\mathcal S). P\rightarrow Q^{\star}$ qui, \`a tout \'el\'ement $(\mathcal S_A,f)\in \mathcal A(\mathcal S). P$, où $\mathcal S_A$ est une restriction de $\mathcal S$ \`a un sous-ensemble fini $A$ de $E$ (consid\'er\'ee \`a l'isomorphie pr\`es) et $f: A \rightarrow P$,   associe le mot $w_A:= \alpha_0\cdots \alpha_{n-1}$ sur $Q$ d\'efini comme suit $n:= \vert A\vert$, $\alpha_i:= (f(a_{i}), \chi(a_i))$ pour $i<n$, où $\{a_0, \dots, a_{n-1}\}$ est une \'enum\'eration des \'el\'ements de $A$ par rapport \`a l'ordre lin\'eaire $\leq$ sur $E$ et $\chi(a):=j$ si  $a\in u_j$, $j\in\{1,\cdots, \ell\}$.
\vspace{1mm}

$1)-$ Montrons que $h$ est injective. Soient $A$ et $A'$ deux parties de $E$ et deux applications $f_A: A \rightarrow P$ et $f_{A'}: A' \rightarrow P$ tels que $h(\mathcal S_A,f_A)=h(\mathcal S_{A'},f_{A'})$, donc $w_A=w_{A'}$. Les sous-ensembles $A$ et $A'$ sont donc de m\^eme cardinalit\'e $n$. Posons $A= \{a_0, \dots, a_{n-1}\}$ et $A'=\{a'_0, \dots, a'_{n-1}\}$, les \'el\'ements \'etant \'enum\'er\'es par rapport \`a l'ordre total $\leq$. Nous avons $(f_A(a_{i}), \chi(a_i))=(f_{A'}(a'_{i}), \chi(a'_i))$, donc $f_A(a_{i})=f_{A'}(a'_{i})$ et $\chi(a_i)=\chi(a'_i)$ pour tout $i<n$. Nous d\'eduisons que $A$ et $A'$ ont le m\^eme spectre pour la partition $u_1,\cdots,u_{\ell}$ de $E$ qui est une d\'ecomposition monomorphe de $\mathcal S$ (Lemme \ref{lem:spectre}). % que $d(A)=d(A')$. %d'apr\`es la d\'efinition de $\chi(a_i)$.
D'o\`u $\mathcal S_A$ et $\mathcal S_{A'}$ sont isomorphes. Ces restrictions \'etant prises \`a l'isomorphie pr\`es, nous d\'eduisons que $(\mathcal S_A,f_A)=(\mathcal S_{A'},f_{A'})$.
\vspace{1mm}

\vspace{1mm}

$2)-$ Reste \`a montrer que $(\mathcal S_A,f_A)\leq (\mathcal S_{A'},f_{A'})$ si et seulement si $w_A\leq w_{A'}$. Supposons $(\mathcal S_A,f_A)\leq (\mathcal S_{A'},f_{A'})$ pour $A:=\{a_0,\dots,a_{n-1}\}$ et $A':=\{a'_0,\dots,a'_{m-1}\}$ o\`u $n\leq m$. Il existe alors un abritement $k$ de $\mathcal S_A$ dans $\mathcal S_{A'}$  tel que
\begin{equation}
f_A(a_i)\leq f_{A'}(k(a_i))\; \text{ pour tout } i<n.\label{eq:1}
 \end{equation}
 Donc $k$ peut-\^etre consid\'er\'e comme un abritement de la cha\^{i}ne $0\leq 1\leq\dots\leq n-1$ dans la cha\^{i}ne $0\leq 1\leq\dots\leq m-1$ et nous pouvons \'ecrire $k(a_i)=a'_{k(i)}$. Si $w_A\nleq w_{A'}$, il existe $i<n$ tel que $\alpha_i\nleq \alpha'_{k(i)}$ ce qui signifie que $$(f_A(a_{i}), \chi(a_i))\nleq(f_{A'}(a'_{k(i)}), \chi(a'_{k(i)}))=(f_{A'}(k(a_{i})), \chi(k(a_{i}))).$$
 D'apr\`es la relation \eqref{eq:1}, nous d\'eduisons que $\chi(a_i)\neq\chi(k(a_{i}))$, c'est \`a dire qu'il existe $j\neq j'$ tel que $a_i\in u_j$ et $k(a_i)\in u_{j'}$, ce qui est absurde puisque $k$ est un abritement de $\mathcal S_A$ dans $\mathcal S_{A'}$. D'o\`u $w_A\leq w_{A'}$.
\vspace{1mm}

 Inversement, si $w_A\leq w_{A'}$, il existe un abritement $p$ de la cha\^{i}ne $0\leq 1\leq\dots\leq n-1$ dans la cha\^{i}ne $0\leq 1\leq\dots\leq m-1$ tel que $\alpha_i\leq \alpha'_{p(i)}$ autrement dit:
 \begin{equation}
 (f_A(a_{i}), \chi(a_i))\leq(f_{A'}(a'_{p(i)}), \chi(a'_{p(i)}))\;\text{ pour tout }i<n.\label{eq:2}
 \end{equation}
Soit $p'$ une application de $A$ dans $A'$ d\'efinie par $p'(a_i):=a'_{p(i)}$. Il est facile de v\'erifier que $p'$ est un abritement de $\mathcal S_A$ dans $\mathcal S_{A'}$. % Pour le montrer: par construction $p'$ est un abritement de $(A,\leq_{\restriction_A})$ dans $(A',\leq_{\restriction_{A'}})$. Reste \`a montrer que $a_i$ et $a'_{p(i)}$ sont dans le m\^eme intervalle $u_j$ ce qui est donn\'e par l'equation \ref{eq:2}. $p'$ v\'erifie $f_A(a_{i})\leq f_{A'}(p'(a_{i}))$.
\end{proof}

\bigskip

D'apr\`es le Th\'eor\`eme \ref{theo:pouzet-borne} %(\cite{pouzet-belordre72})
nous pouvons d\'eduire:

\begin{corollary}\label{cor:wqo-monomorphy}
Si la signature est finie, l'\^age d'une structure admettant une d\'ecomposition monomorphe finie a un nombre fini de bornes.
\end{corollary}

Une cons\'equence est l'extension suivante de la Proposition \ref{prop-wqo-ages}:

\begin{proposition}\label{prop:hwqo}
Soit $\mathscr C$ une classe h\'er\'editaire de structures relationnelles de signature $\mu$ finie. S'il existe un entier $\ell$ tel que $dim_{mon}(\mathcal R)\leq \ell$ pour toute structure $\mathcal R$ de $\mathscr C$ %toute structure d'une classe h\'er\'editaire $\mathscr C$ poss\`ede une d\'ecomposition monomorphe ayant au plus $\ell$ blocs
alors $\mathscr C$  est h\'er\'editairement belordonn\'ee.
\end{proposition}

\begin{proof}
D'apr\`es la Proposition \ref{prop:classe-ideal} et le Corollaire \ref{cor:wqo-monomorphy}, $\mathscr C$ est belordonn\'ee. En effet, si $\mathscr C$ n'est pas belordonn\'ee, d'apr\`es la Proposition \ref{prop:classe-ideal}, elle contient un \^age ayant une infinit\'e de bornes, ce qui contredit le Corollaire \ref{cor:wqo-monomorphy}. %sinon il existerait un \^age ayant une infinit\'e de bornes ce qui contredirait le corollaire \ref{cor:wqo-monomorphy}.
 Il s'ensuit que $\mathscr C$ est une union finie d'id\'eaux (voir Th\'eor\`eme \ref{theo:erdos-tarski}), chacun \'etant l'\^age d'une structure relationnelle. D'apr\`es le Th\'eor\`eme \ref{compactness} une telle structure  relationnelle poss\`ede une d\'ecomposition monomorphe finie. Par la Proposition \ref{prop-wqo-ages} son \^age est h\'er\'editairement belordonn\'e. Etant une union finie de classes h\'er\'editairement belordonn\'ee, $\mathscr C$ est h\'er\'editairement belordonn\'ee (voir Propri\'et\'e \ref{pro:unionhere-wqo}).
\end{proof}

\medskip

Nous avons le lemme suivant:

\begin{lemma} \label{lem:reductionbis}
Si une classe h\'er\'editaire $\mathscr C$  de structures relationnelles  finies de signature $\mu$ finie contient, pour tout entier $\ell$, une structure finie $\mathcal S$ telle que $dim_{mon}(\mathcal S)>\ell+1$, %qui ne poss\`ede pas de d\'ecomposition monomorphe ayant au plus $\ell+1$ blocs,
alors elle contient une classe h\'er\'editaire $\mathscr A$, ayant la m\^eme propri\'et\'e, qui est minimale pour l'inclusion. De plus, $\mathscr A$ est l'\^age d'une structure relationnelle $\mathcal R$ telle que $dim_{mon}(\mathcal R)=\infty$. %n'ayant pas de d\'ecomposition monomorphe finie.
\end{lemma}

\begin{proof} Disons que $\mathscr C$ v\'erifie la propri\'et\'e $(P)$ si et seulement si pour tout entier  $\ell$, $\mathscr C$ contient une structure finie $\mathcal S$ telle que $dim_{mon}(\mathcal S)>\ell+1$. %qui ne poss\`ede pas de d\'ecomposition monomorphe ayant au plus $l+1$ blocs.

Si $\mathscr C$ ne contient pas d'anticha\^{i}ne infinie (par rapport \`a l'abritement), alors la collection de ses sous-classes h\'er\'editaires est bien fond\'ee (voir Th\'eor\`eme \ref{theo:higman-equivalence}). Donc $\mathscr A$ est un \'el\'ement minimal de l'ensemble de toutes les sous-classes h\'er\'editaires de $\mathscr C$ qui v\'erifient $(P)$.

     Si $\mathscr C$ contient une anticha\^{i}ne infinie, nous appliquons la Proposition \ref{prop:classe-ideal}. Donc $\mathscr C$ contient un \^age $\mathcal A(\mathcal R)$ qui est belordonn\'e et qui a un nombre infini de bornes. D'apr\`es le Corollaire \ref{cor:wqo-monomorphy} $\mathcal R$  ne peut pas avoir une d\'ecomposition monomorphe finie. D'o\`u, par le Th\'eor\`eme \ref{compactness}, $\mathcal A(\mathcal R)$ v\'erifie $(P)$. Comme $\mathcal A(\mathcal R)$ ne poss\`ede pas d'anticha\^{i}ne infinie, nous retournons au premier cas.

     Il est clair que la classe $\mathscr A$ ne peut pas \^etre l'union de deux classes  h\'er\'editaires propres, donc elle doit \^etre filtrante pour l'abritement. Donc, d'apr\`es Fra\"{\i}ss\'e (Th\'eor\`eme \ref{theo:age=ideal}), c'est l'\^age d'une structure relationnelle $\mathcal R$. Et cette structure relationnelle n'a pas de d\'ecomposition monomorphe finie.
\end{proof}

\vspace{2mm}

Il en r\'esulte pour le cas des classes de structures ordonn\'ees:
\begin{corollary} \label{lem:reduction}
 Si une classe h\'er\'editaire  $\mathscr C$  de structures relationnelles  ordonn\'ees finies de signature finie contient, pour tout entier $\ell$, une structure finie $\mathcal S$ telle que $dim_{int}(\mathcal S)>\ell+1$, %qui ne poss\`ede pas de d\'ecomposition en intervalles ayant au plus $\ell+1$ intervalles,
 alors elle contient une classe h\'er\'editaire $\mathscr A$, ayant la m\^eme propri\'et\'e, qui est minimale pour l'inclusion. De plus, $\mathscr A$ est l'\^age d'une structure relationnelle ordonn\'ee $\mathcal R$ telle que $dim_{int}(\mathcal R)=\infty$. %et celle-ci ne poss\`ede pas de d\'ecomposition en intervalles ayant un nombre fini d'intervalles.
\end{corollary}

        \subsection{Application au profil}\label{subsection:profile}

        Les structures relationnelles poss\`edant une d\'ecomposition monomorphe en un nombre fini de blocs et les classes form\'ees de ces structures ont des propri\'et\'es int\'eressantes notamment la croissance de leurs profils. %Nous \'enon\c{c}ons, ci-dessous, quelques r\'esultats mais certaines preuves, n\'ecessitant l'introduction de nouvelles notions, seront donn\'ees plus loin. \\
\vspace{1mm}

Si une structure relationnelle infinie $\mathcal R$ poss\`ede une d\'ecomposition monomorphe finie avec $k+1$ blocs infinis (ie, $dim_{mon}^{\infty}(\mathcal R)=k+1$),
alors naturellement, le profil de $\mathcal A(\mathcal R)$ est born\'e par un polyn\^ome dont le degr\'e est au plus $k$. En effet, dans ce cas nous avons
$$\varphi_{\mathcal R}(n)\leq \underset{s\leq r}\sum\binom{r}{s}\binom{n+k-s}{k}\leq 2^r\binom{n+k}{k}$$
où $r$ est la cardinalit\'e de l'union des blocs monomorphes finis.

\vspace{1mm}

 En fait dans ce cas, le profil  %et ceci et le r\'esultat principal de  \cite{P-T-2013},
 est un quasi-polyn\^ome (c'est \`a dire une somme $a_{k}(n)n^{k}+\cdots+ a_0(n)$ dont les coefficients $a_{k}(n), \dots, a_0(n)$ sont des fonctions p\'eriodiques) dont le degr\'e est le nombre de blocs infinis de la d\'ecomposition monomorphe canonique de $\mathcal R$ moins $1$ comme stipul\'e par le th\'eor\`eme suivant:

\begin{theorem}\label{theo:Pou-Thie}(voir \cite{P-T-2005,P-T-2013})

Soit $\mathcal R$ une structure relationnelle infinie admettant une d\'ecomposition monomorphe en un nombre fini de blocs et soit $k$ sa dimension monomorphe. % nombre de blocs infinis de la d\'ecomposition canonique de $\mathcal R$,
Alors, la s\'erie g\'en\'eratrice $\mathcal H_{\varphi_{\mathcal R}}$ est une fraction rationnelle de la forme suivante, avec $P\in \mathbb Z[x]$ et $P(1)\neq 0$: $$\dfrac{P(x)}{(1-x)(1-x^2)\cdots(1-x^k)}.$$
En particulier, son profil $\varphi_{\mathcal R}$ est un quasi-polyn\^ome de degr\'e $k-1$, d'o\`u $\varphi_{\mathcal R}(n)\simeq an^{k-1}$ pour un r\'eel positif $a$.
\end{theorem}

\vspace{1mm}

 Ce profil n'est pas n\'ecessairement un polyn\^ome; un  exemple simple est le suivant. \\ Consid\'erons la somme directe de $k+1$ copies ($k\geq 1$) du graphe complet $K_{\infty}$ sur un nombre infini de sommets, le profil est la fonction partition d'entier en un nombre fini de blocs (voir exemples de profils en page \pageref{par:exemple}).

 \vspace{1mm}

Dans le cas des structures ordonn\'ees, il n'y a pas de sym\'etries: si deux structures sont isomorphes, l'isomorphisme est unique. Nous montrons que dans ce cas le profil est un polyn\^ome. %obtenons facilement la premi\`ere partie du Th\'eor\`eme \ref{thm:poly-expo1} en utilisant
 \vspace{1mm}

Nous rappelons auparavant,
un lemme de \cite{P-T-2013} valable pour la notion de d\'ecomposition monomorphe.

Soit $(E_x)_{x\in X}$ une partition d'un ensemble $E$. Soit $d\in \mathbb N$.
 Un sous-ensemble $A$ de $E$ est dit  \emph{$d$-large} si pour tout $x\in X$,
$d_x(A) \geq d~$ lorsque $A\not\supseteq E_x$.

\begin{lemma}\label{lemma.minimalGrowthRate} (voir \cite{P-T-2013})

  Soit $\mathcal R$ une structure relationnelle infinie de base $E$ ayant une d\'ecomposition monomorphe finie  et soit $\mathcal{P(R)}$ sa d\'ecomposition canonique. Alors, il existe un entier $d$ tel que pour tout sous-ensemble $A$ de $E$ qui est $d$-large\index{ensemble!$d$-large}  (par rapporte \`a  $\mathcal P(\mathcal R)$)  la partition $\mathcal P(\mathcal R)_{\restriction A}$  induite par $\mathcal P(\mathcal R)$ sur $A$ co\"{\i}ncide  avec $\mathcal P(\mathcal R_{\restriction A})$. En particulier le spectre de $A$ par rapport \`a $\mathcal{P}(\mathcal R)$ co\"{\i}ncide avec le spectre de $A$ par rapport \`a $\mathcal {P}(\mathcal R_{\restriction A})$.
   \end{lemma}

\begin{theorem}\label{thm: polynomial-interval}
Soit $\mathscr C$ une classe h\'er\'editaire de structures ordonn\'ees de signature $\mu$ finie. S'il existe un entier  $\ell$ tel que $dim_{int}(\mathcal R)\leq \ell+1$ pour tout \'el\'ement $\mathcal R$ de $\mathscr C$ % poss\`ede une d\'ecomposition en intervalles ayant au plus $\ell+1$ classes
alors le profil de $\mathscr C$  est un polyn\^ome de degr\'e au plus $\ell$.
\end{theorem}

\begin{proof}
D'apr\`es le Th\'eor\`eme \ref{thm:intervalsecomposition} et la Proposition \ref{prop:hwqo},  $\mathscr C$ est h\'er\'editairement belordonn\'ee, donc belordonn\'ee. Puisque la collection des sous-classes h\'er\'editaires de $\mathscr C$ est bien fond\'ee nous pouvons raisonner par induction. Consid\'erons une classe h\'er\'editaire $\mathscr C'\subseteq \mathscr C$ telle que $\varphi_{\mathscr {C''}}$ est un polyn\^ome pour toute sous-classe h\'er\'editaire propre $\mathscr {C''}\subset  \mathscr C'$ et montrons que  $\varphi_{\mathscr C'}$ est un polyn\^ome. Nous pouvons supposer que $\mathscr C'$ n'est pas vide. %(sinon, le profil est ind\'efini).

Si $\mathscr C'$ est une union de deux sous-classe h\'er\'editaires propres, disons $\mathscr C _1'$ et $\mathscr C _2'$ alors, d'apr\`es l'\'equation:
\begin{equation}\label{eq-profile}
\varphi_{\mathscr C'}=  \varphi_{\mathscr C_1'}+ \varphi_{\mathscr C_2'}-\varphi_{\mathscr C_1'\cap \mathscr C_2'}
\end{equation}
nous d\'eduisons  que $\varphi_{\mathscr C'}$ est un polyn\^ome. Nous pouvons donc supposer que $\mathscr C'$ n'est pas l'union de deux sous-classes h\'er\'editaires propres. Puisque $\mathscr C'$ n'est pas vide, c'est donc l'\^age d'une structure relationnelle  $\mathcal R$. D'apr\`es le Th\'eor\`eme \ref{compactness} et le Th\'eor\`eme \ref{thm:intervalsecomposition},  $\mathcal R$ a une d\'ecomposition en intervalles ayant un nombre fini de blocs. Soit $\mathcal P:= (V_{x})_{x\in  X}$ la d\'ecomposition canonique de $\mathcal R$, $X_{\infty}:= \{ x\in X:   V_{x}\;  \text{ est infini}\; \}$,  $K$ l'union des composantes monomorphes finies de $\mathcal R$,  $\l+1:= \vert X_{\infty}\vert$ et $d$ donn\'e par le Lemme \ref{lemma.minimalGrowthRate}. Soit $\mathcal F$ l'ensemble des sous-ensembles $A$ de $V(\mathcal R)$ contenant $K$ qui sont $d$-larges.  D'apr\`es le Lemme \ref{lemma.minimalGrowthRate}, $\mathcal P(\mathcal R_{\restriction A})$  est induite par  $\mathcal P(\mathcal R)$ et a le m\^eme nombre de classes que $\mathcal P(\mathcal R)$ pour tout $A\in\mathcal F$. Puisque $\mathcal R$ est ordonn\'ee, deux \'el\'ements  $A, A' \in \mathcal F$ ont le m\^eme spectre si et seulement si $\mathcal R_{\restriction A}$ et $\mathcal R_{\restriction A'}$ sont isomorphes. Donc la collection $\mathscr F:= \{\mathcal R_{\restriction A}: A \in \mathcal F\}$ est un segment final non vide de $\mathscr C'$ et l'application qui, \`a $\mathcal R_{\restriction A}$ associe   $(d_{x}(A)-d)_{x\in X_{\infty}}$  est un isomorphisme (d'ensembles ordonn\'es) de $\mathscr F$ sur $\mathbb N^{l+1}$.  Soit $\vartheta$ la fonction g\'en\'eratrice de $\mathbb N^{l+1}$ (qui compte, pour tout entier $m$ le nombre d'\'el\'ements de $\mathbb N^{l+1}$ dont la somme des composantes vaut  $m$). Nous avons $\vartheta (m)={ m+l \choose l}$ et $\varphi_{\mathscr F}(n)= \vartheta (n-(l+1)d- \vert K\vert)$;  ainsi $\varphi_{\mathscr F}$ est un polyn\^ome de degr\'e $l$. Puisque  $\mathscr C'\setminus \mathscr F$ est une classe h\'er\'editaire propre de $\mathscr C'$,  par induction son profil est un polyn\^ome. Il s'ensuit que le profil de $\mathscr C'$ est un polyn\^ome.
\end{proof}

\vspace{2mm}

La r\'eciproque du Th\'eor\`eme \ref{thm: polynomial-interval} n'est pas vraie en g\'en\'eral, des exemples dans le cas des graphes sont en page \pageref{subsubsec:les dix graphes}.  Nous montrons, dans le chapitre \ref{chap:graphe ordonne}, que la r\'eciproque  est vraie dans le cas d'une structure binaire ordonn\'ee.

\vspace{1mm}

Dans le cas o\`u le profil d'une classe de structures relationnelles ordonn\'ees n'est pas born\'e par un polyn\^ome, nous avons le r\'esultat suivant:
%Le saut dans la croissance  des profils, d'une croissance  polyn\^omiale \`a une croissance  plus rapide a \'et\'e obtenu par Pouzet dans \cite{pouzet.tr.1978} pour les \^ages. L'extension suivante de la Proposition \ref{cor:minimal} montre comment ce r\'esultat se g\'en\'eralise aux classes h\'er\'editaires.
\begin{proposition} \label{cor:minimal} Si une classe h\'er\'editaire $\mathscr C$ de structures relationnelles ordonn\'ees finies de signature $\mu$ finie a un profil non born\'e par un polyn\^ome alors elle contient une classe h\'er\'editaire $\mathscr A$ minimale pour l'inclusion ayant la m\^eme propri\'et\'e.
\end{proposition}
%\vspace{2mm}

Cette proposition donne une instance d'un r\'esultat plus g\'en\'eral valable pour une classe de structures relationnelles finies (pas n\'ecessairement ordonn\'ees).% au cas de classe de structures

\begin{theorem}\label{thm:minimal1}
%If the profile of a hereditary class $\mathscr C$ of finite relational structures (with a finite signature $\mu$) is not bounded by a polynomial then it contains a hereditary class  $\mathscr A$ with this property which is minimal w.r.t. inclusion.
Si une classe h\'er\'editaire $\mathscr C$ de structures relationnelles finies (de signature $\mu$ finie) a un profil non born\'e par un polyn\^ome alors elle contient une classe h\'er\'editaire $\mathscr A$ minimale pour l'inclusion ayant la m\^eme propri\'et\'e.
\end{theorem}

\begin{proof}
La preuve %du Th\'eor\`eme \ref{thm:minimal1}
suit, en tout point, les \'etapes de la preuve du Lemme \ref{lem:reductionbis} en utilisant, \`a la place du
Corollaire \ref{cor:wqo-monomorphy} le Th\'eor\`eme \ref{theo:pouzet-polyborne}.
\end{proof}

\section{Multi-encha\^{i}nabilit\'e et structures invariantes}\label{sec:enchainable-invariance}
\subsection{Presque multi-encha\^{\i}nabilit\'e}

   Une \emph{presque multicha\^{\i}ne}\index{presque multicha\^{i}ne} est une structure relationnelle $\mathcal R$ sur un ensemble de la forme $V:= F\cup (L\times K)$ où $F$ et  $K$ sont deux ensembles finis, pour laquelle il existe un ordre lin\'eaire sur $L$  tel que pour tout isomorphisme local $h$ de la chaîne $C:= (L, \leq) $ l'application $(h, 1_K)$
\'etendue par l'identit\'e sur $F$ est un isomorphisme local de $\mathcal R$ (l'application
$(h, 1_K)$ est d\'efinie par  $(h, 1_K)(x, y):= (h(x), y)$).
\vspace{1mm}

Soit $m$ un entier. Nous disons que $\mathcal R$ a au plus le type  $(\omega, m)$  si $C$ est de type $\omega$ et $\vert F\vert+ \vert K\vert\leq m$.  Une structure relationnelle qui est isomorphe \`a une presque multichaîne est \emph{presque multi-enchaînable}.\index{structure relationnelle!presque multi-encha\^{i}nable}

\vspace{1mm}

La notion de presque multi-enchaînabilit\'e a \'et\'e introduite par \emph{Pouzet} \cite{pouzet.tr.1978} (voir \cite{pouzet06}). Le cas particulier  $\vert K\vert =1$ est
la  notion de \emph{presque enchaînabilit\'e}\index{structure relationnelle!presque encha\^{i}nable} introduite par \emph{Fra\"{\i}ss\'e} (voir \cite{fraisse}).

% It seems to be a little bit hard to swallow, still it is the key in proving the second part of Theorem \ref{thm:poly-expo1} as well that Theorem \ref{theo:croiss-profil}.

Une structure presque multi-enchaînable peut avoir une d\'ecomposition monomorphe finie (exemple. la somme directe d'un nombre fini de copies du graphe complet d\'enombrable). Voici un simple test pour d\'eterminer si une structure presque multi-enchaînable poss\`ede une telle d\'ecomposition.
\vspace{1mm}

%D\'esignons par $\mathscr S_{\mu}$ la classe des structures relationnelles de signature $\mu$ qui n'ont pas de d\'ecomposition monomorphe finie.

\begin{lemma}\label{lem:test}  Soit $\mathcal R$ une presque multicha\^{\i}ne infinie d\'efinie sur  $V:= F\cup (L\times K)$. Alors $dim_{mon}(\mathcal R)=\infty$ si et seulement s'il existe $x,x'\in L$ et $y\in K$ tels que $(x,y)\not \simeq_{\mathcal R} (x',y)$. Ceci est particuli\`erement le cas si les restrictions de $\mathcal R$ \`a $\{(x,y)\}\cup F\cup (\{x\}\times (K \setminus \{y\}))$ et $\{(x',y)\}\cup F\cup (\{x\}\times (K \setminus \{y\}))$ ne sont pas isomorphes.
\end{lemma}

\begin{proof}
La premi\`ere partie est \'evidente. Si $dim_{mon}(\mathcal R)=\infty$ %$\mathcal R \in \mathscr S_{\mu}$
alors $V$ contient un nombre infini de paires d'\'el\'ements non \'equivalents $a_n:= (x_n, y_n)$. Puisque $K$ est fini, il existe  $n\not=n'$ tels que $y_n=y_{n'}$. Posons $x=x_n$, $x':= x_{n'}$  et $y= y_n$.

Inversement, supposons qu'il existe deux \'el\'ements $(x,y), (x',y)$ de $L\times K$ tels que $(x,y) \not \simeq _{\mathcal R} (x',y)$. Nous pouvons supposer, sans perte de g\'en\'eralit\'e, que $x<x'$. Donc il existe un entier $k$ tel que $(x,y) \not \simeq _{k,\mathcal R} (x',y)$.  Soit $A$ une $k$-partie de $V\setminus \{(x, y),(x', y)\}$ telle que les restrictions de $\mathcal R$ \`a $\{(x,y)\} \cup A$ et $\{(x', y)\} \cup A$ ne soient pas isomorphes. Soit $C(A)$ la projection de $A$ sur $L$ et posons $\ell=\vert C(A)\vert$. Comme $C$ est une cha\^{i}ne et $L$ est infini, il existe une infinit\'e de copies de $C(A)\cup\{x,x'\}$ sur $L$ (des $\ell+2$-parties de $L$ dont les \'el\'ements sont dispos\'es dans la cha\^{i}ne $C$ de la m\^eme façon que ceux de $C(A)\cup\{x,x'\}$). Soit $B$ une de ces copies et soit $h$ l'isomorphisme local de $C$ qui transforme $C(A)\cup\{x,x'\}$ en $B$.
  %Disons que deux \'el\'ements $u, v\in L$ sont $\ell$-\'equivalents si la transformation de $u$ en $v$ dans toute $\ell$-partie de $L$ en un isomorphisme local de  $C$ et de mani\`ere similaire pour la transformation de $v$ en $u$. Ceci d\'efinit une relation d'\'equivalence qui d\'ecompose $L$ en un nombre fini d'intervalles de  $C$. Parmi ces intervalles certains sont infinis. Nous pouvons choisir un sous-ensemble infini $L'$ tel que pour tous $x_1<x_1'$ dans $L'$ l'application qui transforme $x$ en $x_1$ et $x'$ en $x'_1$ s'\'etend \`a toute $\ell$-partie de $L$ en un  isomorphisme local de $C$. Soit $h$ une telle application.
  D\'esignons par $\hat h$ l'application $ (h, 1_K)$ \'etendue par l'identit\'e sur $F$. Soient $(x_1,y)=\hat h(x,y)$ et $(x'_1,y)=\hat h(x',y)$. D'apr\`es la d\'efinition de la multi-enchaînabilit\'e, l'application $\hat h$ est un isomorphisme local de $\mathcal R$.  Il s'ensuit que les restrictions de $\mathcal R$ \`a $A\cup\{(x, y),(x', y)\}$ et $\hat h(A)\cup\{(x_1, y),(x'_1, y)\}$  sont isomorphes. Comme les restrictions de $\mathcal R$ \`a $\{(x,y)\} \cup A$ et $\{(x', y)\} \cup A$ ne sont pas isomorphes, nous d\'eduisons que les restrictions de $\mathcal R$ \`a    $\hat h(A)\cup\{(x_1,y)\}$ et $\hat h(A)\cup\{(x'_1, y)\}$ ne sont pas isomorphes, ainsi $(x_1,y) \not \simeq _{\mathcal R} (x_1',y)$. Nous avons donc une infinit\'e de paires $(x,x')$ de $L^2$ avec $x<x'$ tel que $(x,y) \not \simeq _{\mathcal R} (x',y)$.% Il s'ensuit que $L'\times \{y\}$ est inclus dans une classe d'\'equivalence.
\end{proof}

\begin{question} Soit $y\in K$, est-il vrai que pour la relation $\simeq_{\mathcal R}$ ou bien tous les \'el\'ements de  $L\times \{y\}$ sont in\'equivalents ou bien ils sont tous \'equivalents?
\end{question}

\begin{proposition}\label{prop:number}
Si la signature $\mu$ est finie, il existe un nombre fini de structures presque multi-encha\^{\i}nables deux \`a deux non isomorphes de signature $\mu$  qui sont au plus de type  $(\omega, m)$.
\end{proposition}

\begin{proof}
Soit $\ell$ le maximum de la signature $\mu$. Soit $V:= F\cup (L\times K)$ où $F$ et  $K$ sont deux ensembles finis tels que $\vert F\vert +\vert K\vert\leq a$ et $"\leq"$ un ordre lin\'eaire sur  $L$ de type $\omega$. Soit $L'$ une $\ell$-partie de  $L$.
Observons que si  $\mathcal R$ et $\mathcal R'$  sont deux structures relationnelles  presque multi-enchaînables  sur $V$ (avec cette d\'ecomposition de $V$) qui  co\"{\i}ncident sur $F\cup (L'\times K)$ elles sont \'egales. Le nombre de structures relationnelles de signature $\mu$ sur un ensemble fini \'etant fini,  le nombre de telles structures  $\mathcal R$ est forc\'ement fini.
\end{proof}

\vspace{2mm}

Pour construire des structures presque multi-enchaînables, un des outils est le th\'eor\`eme de  Ramsey. Nous l'utiliserons \`a travers la notion de structures invariantes donn\'ee dans \cite{Bou-Pouz} que nous rappelons ci-dessous.

            \subsection{Structures invariantes}\label{subsec:invariance}

Soit  $C:=(L,\leq )$ une chaîne.
Pour tout entier $n$, soit $[C]^{n}$ l'ensemble des  $n$-tuples $\overrightarrow{a}:=(a_{1},...,a_{n})\in L^{n}$ tels que $a_{1}<...<a_{n}$. Cet ensemble est identifi\'e \`a l'ensemble des $n$-parties de  $L$.

Pour tout automorphisme local $h$ de  $C$ de domaine $D$, posons  $h(\overrightarrow{a}):=(h(a_{1}),...,h(a_{n}))$ pour tout $\overrightarrow{a}\in \lbrack D]^{n}$.

Soit $\mathfrak{L}:=\left\langle C,{\mathcal R},\Phi \right\rangle $ un triplet form\'e d'une chaîne $C$ sur $L$,  d'une  structure  relationnelle ${\mathcal R}:=(V,(\rho _{i})_{i\in I})$ et d'un ensemble  $\Phi $ d'applications, chacune \'etant une application $\psi $ de $[C]^{a(\psi )}$ dans $V$, où $a(\psi )$ est un entier, l'\emph{arit\'e} de $\psi$.

Nous disons que
$\mathfrak{L}$ est \emph{invariante} si:
\begin{equation}\label{eq:invariance}
\rho _{i}(\psi _{1}(\overrightarrow{\alpha }_{1}),...,\psi _{m_{i}}(%
\overrightarrow{\alpha }_{m_{i}}))=\rho _{i}(\psi _{1}(h(\overrightarrow{%
\alpha }_{1})),...,\psi _{m_{i}}(h(\overrightarrow{\alpha }_{m_{i}})))
\end{equation}
pour tout  $i\in I$ et tout  automorphisme local $h$ de $C$ dont le domaine contient $\overrightarrow{\alpha }_{1},...,\overrightarrow{\alpha }_{m_{i}},$
 où $m_{i}$ est l'arit\'e de $\rho _{i}$, $\psi _{1},...,\psi _{m_{i}}\in\Phi ,$ $\overrightarrow{\alpha }_{j}\in \lbrack C]^{a(\psi _{j})}$ pour $j=1,...,m_{i}.$
\vspace{1mm}

Cette condition exprime le fait que chaque $\rho _{i}$ est invariante pour la transformation des $m_{i}$-tuples de $V$ qui est induite sur $V$ par les automorphismes locaux de $C$. Par exemple, si ${\rho}$ est une relation binaire et $\Phi =\{\psi\} $ alors
\begin{equation*}
\rho(\psi (\overrightarrow{\alpha }),\psi (\overrightarrow{\beta }))=%
\rho(\psi (h(\overrightarrow{\alpha })),\psi (h(\overrightarrow{\beta})))
\end{equation*}%
ce qui signifie que, $\rho (\psi (\overrightarrow{\alpha }),\psi (\overrightarrow{\beta }))$ d\'epend uniquement des  positions relative de $\overrightarrow{\alpha }$ et $\overrightarrow{\beta }$ sur la chaîne $C$.
\vspace{1mm}

Si $\mathfrak{L}:=\left\langle C,{\mathcal R},\Phi \right\rangle $ et $A$ est un sous-ensemble de $L$, posons $\Phi_{\restriction_{A}}:=\{\psi _{\restriction_{\lbrack A]^{a(\psi )}}}:\psi \in\Phi \}$ et $\mathfrak{L}_{\restriction_{A}}:=\left\langle
C_{\restriction_{A}},{\mathcal R},\Phi _{\restriction_{A}}\right\rangle $ la restriction de $\mathfrak{L}$ \`a $A$.
\vspace{1mm}

Le r\'esultat suivant est une cons\'equence du th\'eor\`eme de Ramsey\index{Ramsey} (Th\'eor\`eme \ref{thm:ramsey}):

\begin{theorem}\label{thm:ramsey-invariant}(voir \cite{Bou-Pouz})

Soit $\mathfrak{L}:=\left\langle C,{\mathcal R},\Phi \right\rangle $ une structure telle que le domain $L$ de $C$ est infini, ${\mathcal R}$
form\'ee d'un nombre fini de relations et $\Phi $ est fini. Alors il existe un sous-ensemble infini $L^{\prime }$ de $L$ tel que $\mathfrak{L}_{\upharpoonleft_{L^{\prime }}}$ est invariante.
\end{theorem}

\subsection{Une application}

 D\'esignons par $\mathscr S_{\mu}$ la classe des structures relationnelles de signature $\mu$ qui n'ont pas de d\'ecomposition monomorphe finie.

 \begin{conjecture}\label{conjecture2}
Il existe un sous-ensemble fini $\mathfrak A$ form\'e de structures incomparables de $\mathscr S_{\mu}$ tel que tout \'el\'ement de $\mathscr S_{\mu}$ abrite un \'el\'ement de $\mathfrak A$.
\end{conjecture}

Notons que si nous rempla\c{c}ons  $\mathscr S_{\mu}$ par la classe  $\mathscr B$ form\'ee de bicha\^{\i}nes, $\mathfrak A$ poss\`ede vingt \'el\'ements \cite{mont-pou},
 tandis que si $\mathscr S_{\mu}$ est remplac\'ee par la classe $\mathscr T$ form\'ee de tournois, $\mathfrak A$ poss\`ede douze \'el\'ements \cite{Bou-Pouz}. Nous montrons (Th\'eor\`eme \ref{theo:base}) que la Conjecture \ref{conjecture2} est \'egalement vraie dans le cas de structures ordonn\'ees. Les cas des graphes non dirig\'es et des graphes ordonn\'es sont renvoy\'es au chapitre \ref{chap:graphe ordonne}.
 \vspace{1mm}

 Soit la classe $\mathscr D_{\mu}$  des structures relationnelles  ordonn\'ees de signature $\mu$ (finie) ne poss\`edant pas de d\'ecomposition en intervalles ayant un nombre fini de blocs.

%Nous avons le th\'eor\`eme suivant dont la preuve, n\'ecessitant l'introduction d'autres notions %bas\'ee sur le th\'eor\`eme de Ramsey (voir paragraphe \ref{}),
%sera donn\'ee plus loin:

 \begin{theorem}\label{theo:base}
Il existe un sous-ensemble fini $\mathfrak A$ form\'e de structures incomparables de $\mathscr D_{\mu}$ tel que tout \'el\'ement de $\mathscr D_{\mu}$ abrite un \'el\'ement de $\mathfrak A$.
\end{theorem}

\vspace{1mm}

La preuve de ce th\'eor\`eme n\'ecessite le r\'esultat suivant.
\vspace{1mm}

  Soit $k$ un entier. Soit $\mathscr D_{\mu, k}$ la classe des structures ordonn\'ees $\mathcal R$ de signature $\mu$ telle que $\simeq_{\leq k,\mathcal R}$ poss\`ede un nombre infini de classes.

\begin{proposition}\label{mainthm}
Soit $k$ un entier et $\mu$ finie. Tout \'el\'ement de  $\mathscr D_{\mu, k}$ abrite un \'el\'ement de $\mathscr D_{\mu,k}$ qui est presque  multi-encha\^{i}nable et est de type au plus  $(\omega, k+1)$.
\end{proposition}

\begin{proof}
Soit  $\mathcal R\in \mathscr D_{\mu, k}$ de base $V$. Alors $\simeq_{\leq k,\mathcal R}$ a un nombre infini de classes.  Prendre une suite infinie $(x_{p})_{p\in \mathbb N}$ d'\'el\'ements deux \`a deux non \'equivalents. Pour toute paire $(p,q)$, $p<q$,  nous pouvons trouver un sous-ensemble $\mathcal F(p,q)$
de $V\setminus \{x_p, x_q\}$ d'au plus $k$ \'el\'ements t\'emoignant du fait que $x_{p}$ et $x_{q}$ sont non \'equivalents. En fait,  nous pouvons trouver une famille $ \Phi$ d'applications  $f:\mathbb{N}\rightarrow V\mathbb{,}$ $g_{i}:[\mathbb{N]}^{2}\rightarrow V$ pour $i=1,\dots, k$ telles que pour tous $p<q\in \mathbb{N}$, $x_{p}=f(p)$, $\mathcal F(p,q)=\{g_{i}(p,q),~i=1,\dots, k\}$.

Soit $C:= (\mathbb N, \leq)$ et $\mathfrak {L}:=\left\langle C,{\mathcal R},\Phi \right\rangle $.  Le Th\'eo\`eme de Ramsey sous la forme du Th\'eor\`eme  \ref{thm:ramsey-invariant} assure qu'il existe un sous-ensemble infini  $X$
de $\mathbb{N}$ tel que $\mathfrak {L}_{\restriction_X}$ est  invariante.

Nous pouvons supposer, sans perte de g\'en\'eralit\'e, que $X= \mathbb N$.

\noindent\begin{claim}\label{basicclaim}
Supposons que  $\mathfrak {L}$ est  invariante.
\begin{enumerate}
\item Soit $i,j\in \{1, \dots, k\}$. Alors $g_i=g_j$ si et seulement s'il existe  $p<q$ tels que $g_i(p,q)=g_j(p,q)$.
\item Soit $i\in\{1,\dots, k\}$. Alors $g_i$ est constante si et seulement s'il existe trois entiers $p<q<r$ tels que $g_i(p,q)=g_i(q,r)$.
\item Si les restrictions $\mathcal R_{\restriction_ {\{f(p)\} \cup \mathcal F(p,q)}}$ et  $\mathcal R_{\restriction_{\{f(r)\} \cup \mathcal F(p,q)}}$ sont isomorphes pour des entiers $p<q<r$ alors les  restrictions correspondantes pour $p'<q'<r'$ sont  isomorphes.
\end{enumerate}
\end{claim}

%\noindent{\bf Preuve de l'affirmation \ref{basicclaim}}.
\begin{proofclaim}
\begin{enumerate}
\item La condition n\'ecessaire est \'evidente. Pour la condition suffisante, supposons qu'il existe  $p<q$ tels que $g_i(p,q)=g_j(p,q)$. Soient $p'<q'$ et $h$ l'isomorphisme local de $(\mathbb N, \leq)$ qui envoie $p,q$ sur $p',q'$. Appliquons deux fois l'Equation \eqref{eq:invariance} \`a l'ordre d\'efini sur $\mathcal R$, nous avons:
$$g_i(p,q)\leq g_j(p,q)\Leftrightarrow g_i(p',q')\leq g_j(p',q')$$
     et $$g_j(p,q)\leq g_i(p,q)\Leftrightarrow g_j(p',q')\leq g_i(p',q')$$ Nous obtenons $g_i(p',q')=g_j(p',q')$ pour tous $p'<q'$.
\item De m\^eme que pour $1.$, la condition n\'ecessaire est \'evidente. Pour la condition sufisante, soient $p'<q'<r'$ et $h$ l'isomorphisme local de  $(\mathbb N, \leq)$ qui envoie $p,q, r$ sur $p',q', r'$. De m\^eme que pr\'ec\'edemment, en appliquant deux fois l'Equation \eqref{eq:invariance}
\`a l'ordre d\'efini sur $\mathcal R$ nous obtenons  $g_i(p',q')=g_i(q',r')$.  En particulier, nous obtenons  $g_i(n,m)=g_i(m,l)$ pour tout $l>m$. Ainsi pour tous $n<m$, la valeur $g_i(n,m)$ est ind\'ependante de $m$.  %With the fact that $g_i(n',m')=g_i(m',k')$
Il s'ensuit que $g_i$ est constante. %\hfill $\Box$
\item Se d\'emontre de la m\^eme façon.
\end{enumerate}
\end{proofclaim}

Ainsi, si $\mathfrak L$ est invariante, la cardinalit\'e de $\mathcal F(p,q)$ est constante pour tous $p<q$. Soit $\ell$ cette cardinalit\'e.\\

Nous disons que  $\mathfrak {L}$ est de type $(I)$ si les restrictions $\mathcal R_{\restriction_{ \{f(p)\} \cup \mathcal F(p,q)}}$ et  $\mathcal R_{\restriction_{\{f(r)\} \cup \mathcal F(p,q)}}$ sont isomorphes pour des entiers  $p<q<r$. Autrement elle est de type $(II)$.

\noindent \begin{claim}\label{keyclaim}
Soient $p<q<r$.
Si $\mathfrak{L}$ est de type $(I)$ alors les restrictions $\mathcal R_{\restriction_{ \{f(q)\} \cup \mathcal F(p,q)}}$ et  $\mathcal R_{\restriction_{ \{f(r)\} \cup \mathcal F(p,q)}}$ ne sont pas isomorphes, par contre si $\mathfrak {L}$ est de type $(II)$
les restrictions $\mathcal R_{\restriction_{ \{f(p)\} \cup \mathcal F(p,q)}}$ et  $\mathcal R_{\restriction_{ \{f(r)\} \cup \mathcal F(p,q)}}$ ne sont pas isomorphes.
\end{claim}

%\noindent{\bf Proof of Claim \ref{keyclaim}}.
\begin{proofclaim}
Puisque $\mathfrak L$ est invariante, si elle est de type $(I)$ alors les restrictions $\mathcal R_{\restriction_{ \{f(p)\} \cup \mathcal F(p,q)}}$ et  $\mathcal R_{\restriction_{ \{f(r)\} \cup \mathcal F(p,q)}}$ sont isomorphes. Comme $f(p)$ et $f(q)$ ne sont pas \'equivalents et $\mathcal F(p,q)$ en t\'emoigne, les restrictions $\mathcal R_{\restriction_{ \{f(p)\} \cup \mathcal F(p,q)}}$ et  $\mathcal R_{\restriction_{ \{f(q)\} \cup \mathcal F(p,q)}}$ ne sont pas isomorphes, ainsi les restrictions $\mathcal R_{\restriction_{ \{f(q)\} \cup \mathcal F(p,q)}}$ et  $\mathcal R_{\restriction_{ \{f(r)\} \cup \mathcal F(p,q)}}$ ne sont pas isomorphes. Par d\'efinition, si $\mathfrak L$ est de type $(II)$ les restrictions $\mathcal R_{\restriction_{ \{f(p)\} \cup \mathcal F(p,q)}}$ et  $\mathcal R_{\restriction_{ \{f(r)\} \cup \mathcal F(p,q)}}$ ne sont pas isomorphes.
%\hfill  $\Box$
\end{proofclaim}
\medskip

Soit $J$ la plus grande (au sens de l'inclusion) partie de $\{1,\dots,k\}$ v\'erifiant, pour tout $j\in \{1,\dots,k\}\setminus J$ il existe $i\in J$ tel que $g_i=g_j$. Posons $\vert J\vert=\ell$. Soit $F':=\{i\in J: g_i \text{ est constante} \}$. %Si $g_i\neq g_j$ pour tous $i,j\in J$, nous posons $F'=J$, sinon $F'$ est la plus grande partie de $J$ v\'erifiant, pour tout $j\in J\setminus F'$ il existe $i\in F'$ tel que $g_i=g_j$.
Posons  $F=\underset{i\in F'}\cup g_i(\mathbb N^ 2)$ l'image de ces $g_i$ et $K':= \{0\}\cup (J \setminus F')$. Nous pouvons supposer, quitte \`a r\'eind\'exer les \'el\'ements, que  %$J$ est form\'e des derniers \'el\'ements de $\{1, \dots k\}$, ainsi
$K':= \{0,\dots k'\}$. Soit $V':= (\mathbb N\times K') \cup F'$. Nous d\'efinissons une application $\tilde \Phi$ de $V'$  dans $V$ comme suit.

Nous posons $\tilde \Phi(i)=g_i(0,1)$ pour $i\in F'$ et $\tilde \Phi(n,i):= g_i(2n, 2n+1)$ pour $i=1,\dots, k'$. Si $\mathfrak {L}$ est de type $(I)$, nous posons   $\tilde \Phi(n,0):= f(2n+1)$  et si $\mathfrak {L}$ est de type $(II)$, nous posons  $\tilde \Phi(n,0):= f(2n)$.   Soit $Im (\tilde \Phi)$ l'image de  $\tilde \Phi$.

\begin{claim}\label{lem:infiniteclasses}
La restriction de $\mathcal R$ \`a $Im(\tilde \Phi)$ est presque multi-encha\^{i}nable et appartient \`a $\mathscr D_{\mu,k}$.
\end{claim}

%{\bf Proof of claim \ref{lem:infiniteclasses}. }
\begin{proofclaim}
L'application $\tilde \Phi$ est injective d'apr\`es la d\'efinition de $F'$ et de $K'$ et l'Affirmation \ref{basicclaim}.% les images de deux \'el\'emets distincts de $F'$ sont distinctes. $\tilde \Phi(n,i)\neq \tilde \Phi(n,j)$ pour tout $n\in\mathbb N$ et tous $i,j\in\{1,\dots,k'\}$ %Plus de detail a mettre dans le basic claim

Soit $\mathcal S$ l'image inverse de  $\mathcal R$ par  $\tilde \Phi$ (c'est \`a dire la structure d\'efinie sur $V'$ et qui est isomorphe \`a la restriction de $\mathcal R$ \`a $\tilde \Phi(V')$). Nous montrons que $\mathcal S$ est une presque multicha\^{i}ne en montrant que pour tout isomorphisme local $h$ de $C:=(\mathbb N, \leq)$ l'application $\hat h:=(h, 1_{K'})$
\'etendue par l'identit\'e sur $F'$ est un isomorphisme local de $\mathcal S$ (l'application
$(h, 1_{K'})$ est d\'efinie par  $(f, 1_{K'})(x, y):= (f(x), y)$). Ceci est donn\'e par l'invariance de $\mathfrak L$.  Maintenant,  observons que
par construction, $(n,0) \not \simeq_{\leq k,S}(n',0)$ pour $n<n'$. En effet, posons  $A':=(\{n\}\times K'\setminus\{0\}) \cup F'$ et $A:=\tilde \Phi(A')$.  Nous avons que $\mathcal S_{\restriction \{(n,0)\} \cup A'}$ est isomorphe \`a  $\mathcal R_{\{\tilde \Phi(n, 0)\}\cup A}$ et $\mathcal S_{\restriction \{(n',0)\} \cup A'}$ est isomorphe \`a  $\mathcal R_{\{\tilde \Phi(n', 0)\}\cup A}$. Nous avons $A:=\{ g_i(2n, 2n+1): i=1,\dots, k'\} \cup F$ (c'est \`a dire $A=\mathcal F(2n,2n+1)$). Si  $\mathfrak L$ est de  type $(I)$, $\tilde \Phi(n, 0)=f(2n+1)$,  $\tilde \Phi(n', 0)=f(2n'+1)$    et par l'Affirmation \ref{keyclaim}, les restrictions $\mathcal R_{\restriction \{f(2n+1)\} \cup A}$ et ${\mathcal R}_{\restriction \{f(2n'+1)\} \cup A}$ ne sont pas isomorphes, autrement dit $(n,0) \not \simeq_{\leq k,\mathcal S}(n',0)$. Si  $\mathfrak {L}$ est de type $(II)$, nous avons  $\tilde \Phi(n,0)= f(2n)$ et $\tilde \Phi(n', 0)=f(2n')$.  En appliquant encore l'Affirmation \ref{keyclaim}, les restrictions $\mathcal R_{\restriction \{f(2n)\} \cup A}$ et $\mathcal R_{\restriction \{f(2n')\} \cup A}$ ne sont pas isomorphes, c'est \`a dire  $(n,0) \not \simeq_{\leq k,\mathcal S}(n',0)$.

Ainsi, d'apr\`es le Lemme \ref{lem:test}, $\mathcal S\in \mathscr D_{\mu,k}$, donc  $\mathcal R_{\restriction_{Im(\tilde \Phi)}}\in \mathscr D_{\mu,k}$.
%\hfill      $\Box$
\end{proofclaim}

Ceci termine la preuve de la Proposition \ref{mainthm}.
\end{proof}

\vspace{2mm}

\textbf{Preuve du Th\'eor\`eme \ref{theo:base}.}
 Soit $\mathcal R\in \mathscr D_{\mu}$. Soit $m$ le maximum de $\mu$. Soit $k=i(m)$ tel que $\simeq _{\mathcal R}$ co\"{\i}ncide avec $\simeq _{\leq k,\mathcal R}$ (voir le paragraphe \ref{subsection: hypomorphy}). D'apr\`es la Proposition \ref{mainthm}, $\mathcal R$ poss\`ede une restriction appartenant \`a $\mathscr D_{\mu,k}$ qui est presque  multi-enchaînable et est de type au plus $(\omega, k+1)$. D'apr\`es la Proposition \ref{prop:number} le nombre de ces structures est fini.
\hfill      $\Box$

\vspace{2mm}

Comme cons\'equence et en utilisant le Lemme \ref{lem:reduction}, nous avons le th\'eor\`eme de dichotomie pour les classes de structures ordonn\'ees:

\begin{theorem}\label{theo:dichotomie-ordon}
Soit $\mathscr C$ une classe h\'er\'editaire de structures relationnelles ordonn\'ees finies de signature finie. Alors ou bien il existe un entier $\ell$ tel que tout membre de $\mathscr C$ poss\`ede une d\'ecomposition en intervalles ayant au plus $\ell+1$ intervalles et dans ce cas le profil de $\mathscr C$ est un polyn\^ome de degr\'e au plus $\ell$, ou bien $\mathscr C$ contient l'\^age d'une structure presque multi-encha\^{i}nable appartenant \`a un ensemble fini.
\end{theorem}

Nous conjecturons que les profils des structures presque multi-encha\^{i}nables de $\mathfrak A$ dans le Th\'eor\`eme \ref{theo:base} sont au moins exponentiels, c'est \`a dire:

\begin{conjecture}\label{cojecture:ordonne}
Soit $\mathscr C$ une classe h\'er\'editaire de structures ordonn\'ees finies de signature finie. Alors, ou bien il existe un entier $\ell$ tel que tout membre de $\mathscr  C$ poss\`ede une d\'ecomposition monomorphe en au plus $\ell+1$ blocs, auquel cas $\mathscr C$ est une union finie d'\^ages de structures ordonn\'ees, chacune ayant une d\'ecomposition en intervalles ayant au plus $\ell+1$ blocs et le profil de $\mathscr C$ est un polyn\^ome, ou bien le profil de $\mathscr C$ est au moins exponentiel.
\end{conjecture}

Dans le chapitre \ref{chap:graphe ordonne}, nous d\'emontrons la Conjecture \ref{cojecture:ordonne} dans le cas particulier o\`u $\mathscr C$ est une classe h\'er\'editaire de structures binaires ordonn\'ees finies de type $k$.

\section{Une tentative de description des structures \`a profil polynomial}\label{sec:cellulaire}

Comme nous l'avons not\'e, le profil de l'\^age d'une structure relationnelle admettant une d\'ecomposition monomorphe finie est born\'e par un polyn\^ome
et en fait est  un quasi polyn\^ome. La r\'eciproque est fausse, nous verrons dans la section \ref{subsec:graph}, en page \pageref{subsubsec:les dix graphes} des exemples de graphes ayant un profil polynomial parmi les graphes sans d\'ecomposition monomorphe finie. Pouzet a montr\'e en 1978 \cite{pouzet.tr.1978} que si le profil d'une structure relationnelle est
 born\'e par un polyn\^ome, alors son \^age est celui d'une structure presque multi-encha\^{i}nable. Ce r\'esultat est loin d'une caract\'erisation, comme nous le verrons dans la section \ref{subsec:graphesordonnes}, de tels \^ages
 peuvent avoir un profil exponentiel. Par contre, Pouzet  en 2006 \cite{pouzet06} a caract\'eris\'e les  \^age de graphes (non dirig\'es et sans boucle)  dont le profil
 est born\'e par un polyn\^ome \`a l'aide d'une notion voisine de celle de d\'ecomposition monomorphe, la cellularit\'e, introduite par Schmerl en 1990 \cite{Sch}. Nous pr\'esentons cette notion ci-dessous.% ainsi qu'une conjecture (cf. citez decomposition celluler du 14-4-14).

\subsection{Structures cellulaires}
%Disons avec Schmerl qu'
Une structure  $\mathcal R:= (V,
(\rho_i)_{i \in I})$ est dite \emph{cellulaire}\index{structure relationnelle!cellulaire} s'il existe un sous-ensemble fini $F$ de $V$ et une \'enum\'eration
$(a_{x, y})_{(x,y)\in L\times K}$ des  \'el\'ements de $V\setminus F$ par un ensemble
 $L\times K$, o\`u  $K$ est fini, de sorte que pour chaque permutation
$f$ de $L$ l'application  $(f,1_K)$
\'etendue par l'identit\'e sur  $F$ est un automorphisme local de $\mathcal R$ (l'application
$(f,1_K)$ est d\'efinie  par  $(f,1_K)(x, y):= (f(x), y)$).
Notons que par d\'efinition une structure cellulaire est presque multi-encha\^{i}nable.

\vspace{1mm}

%Il est imm\'ediat de constater que
Le profil d'une relation cellulaire est born\'e par un polyn\^ome, en effet si $\vert K\vert=k$ alors $\varphi_{\mathcal R}(n)\leq \binom{n+k}{k}$. %il s'agit de prendre n elements dans la reunion de F avec $L\times {y}$ pour $y\in K$

\vspace{2mm}

D'apr\`es le Th\'eor\`eme 2.14 de (Pouzet, 2006 \cite{pouzet06}) nous avons
\begin{theorem}\label{theo:graphe cellulair}\cite{pouzet06}

Le profil de l'\^age d'un graphe (non dirig\'es et sans boucle) est born\'e par un polyn\^ome si et seulement si c'est l'\^age d'un graphe cellulaire.\index{graphe!cellulaire}
\end{theorem}

 Ce th\'eor\`eme s'ensuit des deux r\'esultats suivants  (Th\'eor\`eme \ref{profilepouzet2} et
Lemme \ref{cellularlemma})

\begin{lemma}\cite{pouzet06}\label{cellularlemma}

La croissance du profil d'un graphe non cellulaire presque multi-encha\^inable est au moins exponentielle.
\end{lemma}

\begin{theorem}\cite{pouzet06}\label{profilepouzet2}

 Soit $\mathcal R$ une structure relationnelle de signature finie. Si le profil de $\mathcal R$ est born\'e (sup\'erieurement) par un polyn\^ome alors $\mathcal R$ est une structure presque multi-encha\^{i}nable.
\end{theorem}

\subsection{D\'ecomposition cellulaire}
Les  graphes mis a part, il n'est pas vrai en g\'en\'eral que les structures relationnelles dont le profil est born\'e par un polyn\^ome soient des structures cellulaires. Comme le lecteur peut facilement le constater, si $\mathcal P$ est un ensemble ordonn\'e donn\'e par la somme directe de deux copies de la cha\^{i}ne des entiers non n\'egatifs, alors $\mathcal P$ a une d\'ecomposition monomorphe finie (chaque bloc est form\'e d'une cha\^{i}ne) donc son profil est born\'e par un polyn\^ome (il est lin\'eaire), en fait $\varphi_{\mathcal P}(n)=\lfloor\frac{n}{2}\rfloor+1$. Par contre, $\mathcal P$ n'est  pas cellulaire. Si au lieu de $\mathcal P$, on consid\`ere la somme directe de deux cliques d\'enombrables cette structure est cellulaire et a m\^eme profil.
Ceci nous am\`ene \`a la notion de d\'ecomposition cellulaire.

\vspace{2mm}

Appelons \emph{d\'ecomposition cellulaire}\index{decomposition@d\'ecomposition!cellulaire} d'une structure  $\mathcal R:= (V,
(\rho_i)_{i \in I})$ la donn\'ee d'un sous-ensemble fini $F$ de $V$ et d'une \'enum\'eration
$(a_{x, y})_{(x,y)\in L\times K}$ des  \'el\'ements de $V\setminus F$ par un ensemble
 $L\times K$  o\`u  $K$ est fini, de sorte que  les restrictions de $\mathcal R$ \`a deux parties finies $A$ et $A'$ de $V$ soient isomorphes si $A\cap F=A'\cap F$ et les \emph{vecteurs fr\'equences}\index{vecteur fr\'equence} $\chi_A$ et $\chi_{A'}$ sont \'egaux; le vecteur fr\'equence $\chi_A$ est \'egal \`a $(\chi_A(K'))_{\emptyset \not = K'\subseteq K}$ et $\chi_A(K')$ est le nombre d'\'el\'ements $y$ de $L$ tels que $A\cap (\{y\}\times K)= \{y\}\times K'$.

Clairement,  les ensembles intervenant dans la d\'efinition de structure cellulaire forment une d\'ecomposition cellulaire.
Ainsi une structure cellulaire a une d\'ecomposition cellulaire.  La r\'eciproque est fausse, si $\mathcal P$ est la somme directe de deux copies de la cha\^{i}ne des entiers non n\'egatifs, la famille $(E_n)_{n\in \NN}$ dans laquelle $E_n:= \{n\}\times \{0, 1\}$ est une d\'ecomposition cellulaire mais $\mathcal P$ n'est pas cellulaire.

Cette notion g\'en\'eralise la notion de d\'ecomposition monomorphe.
\begin{proposition}
Si une structure relationnelle $\mathcal R$ a une d\'ecomposition monomorphe finie alors une restriction du m\^eme \^age admet une d\'ecomposition cellulaire.
\end{proposition}

\begin{proof}
Soit $(V_x)_{x\in X}$ une d\'ecomposition monomorphe finie de $\mathcal R$. Soit $X'$ l'ensemble des $x\in X$ tels que $V_x$ soit infini. Pour tout $x\in X'$ soit $f_x$ une injection de $\NN$ dans $V_x$; pour chaque $n\in \NN$ soit $E_n:= \{f_x(n):x\in X'\}$. Soit $F:= \underset{x\in X\setminus X'}{\bigcup} V_x$ et soit $V':= F\cup\underset{n \in\NN} {\bigcup} E_n$ et $\mathcal R':= \mathcal R_{\restriction_{V'}}$. Alors $\mathcal R'$ est du m\^eme \^age que $\mathcal R$ et $(F, (E_n)_n)$ est une d\'ecomposition cellulaire de $\mathcal R'$.
En effet, posons $V'_x:= V_x$ si $x\in X\setminus X'$ et $V'_x:=  Im(f_x)$ si $x\in X'$. Alors, $(V'_x)_{x\in X}$ est une d\'ecomposition monomorphe de $\mathcal R'$ et donc $\mathcal R'$ a  m\^eme \^age  que $\mathcal R$. Si deux sous-ensembles finis de $V'$ ont m\^eme fr\'equence sur $(F, (E_n)_n)$ alors comme on le constate ais\'ement, leurs traces sur $V'_x$ ont m\^eme cardinalit\'e et donc puisque $(V'_x)_{x\in X}$ est une d\'ecomposition monomorphe de $\mathcal R'$, les restrictions \`a ces deux sous-ensembles sont isomorphes. Ce qui prouve que $(F, (E_n)_n)$ est une d\'ecomposition cellulaire.
\end{proof}

Nous avons facilement:
\begin{lemma}
Si une structure relationnelle a une d\'ecomposition cellulaire alors son profil est born\'e par un polyn\^ome.
\end{lemma}

\begin{problem}
Est-ce qu'une structure relationnelle  d'arit\'e born\'ee a une  d\'ecomposition cellulaire  d\`es que son  profil est born\'e par un polyn\^ome?
 \end{problem}

Comme autre probl\`emes, mentionnons les suivants:

\begin{problems}
\begin{enumerate}
\item Est-ce que le profil de  $\mathcal R$ est un  quasi-polyn\^ome  si   $\mathcal R$ a une d\'ecomposition cellulaire?
\item Est-ce que la fonction g\'en\'eratrice de $\mathcal R$ est une fraction rationnelle  d\`es que $\mathcal R$ a une d\'ecomposition cellulaire?
\end{enumerate}
\end{problems}

\clearemptydoublepage

\chapter[Cas des structures binaires]{Exemples de structures binaires sans d\'ecomposition monomorphe finie}\label{chap:graphe ordonne}

\section{Introduction}
Dans ce chapitre nous nous int\'eressons au cas des structures binaires avec un int\'er\^et particulier pour celles qui sont ordonn\'ees. %Nous commençons par donner les propri\'et\'es de la relation d\'equivalence d\'efinie dans le chapitre \ref{chap:monomorphe} dans certains cas particuliers de structures binaires.
Le r\'esultat principal de ce chapitre est le r\'esultat de dichotomie, pour les structures binaires ordonn\'ees, donn\'e par le  Th\'eor\`eme \ref{theo:dichotomie} o\`u nous montrons que le profil d'une structure binaire ordonn\'ee de type $k$ est soit polynomial soit born\'e par une exponentielle. R\'esultat qui sera \'etendu aux classes h\'er\'editaires de structures binaires ordonn\'ees.

\section{La relation d'\'equivalence sur quelques exemples}\label{sec:exemple-equivalence}
Dans cette section nous revenons sur la relation d'\'equivalence d\'efinie dans la sous-section \ref{subsection:canonical}, %du chapitre \ref{chap:monomorphe},
nous \'etudions les propri\'et\'es de cette relation sur certains cas particuliers de structurs binaires. %Nous commençons par le lemme suivant qui est vrai pour une structure relationnelle quelconque.

\subsection{Le cas des bicha\^{i}nes}

Soit "$\leq$" un ordre lin\'eaire d\'efini sur un ensemble $V$ et soient $x,y,z\in V$. Nous disons que $z$ est entre $x$ et $y$ modulo $\leq$, ou juste entre $x$ et $y$ s'il n'y a pas de risque de confusion, si ou bien $x\leq z\leq y$ ou bien $y\leq z\leq x$.

\begin{lemma}\label{lem:equiv-bichaine}
Soit   $\mathcal B:=(V, \leq_1, \leq_2)$ une bichaine.\index{bichaine@bicha\^{i}ne} Deux \'el\'ements   $x,y\in V$ sont $1$-\'equivalents si et seulement si tout \'el\'ement $z\in V$ qui est entre $x$ et $y$ pour l'un des ordres l'est aussi pour l'autre ordre.
\end{lemma}

\begin{proof} Soit $z\in V\setminus\{x,y\}$. Supposons $z$ entre $x$ et $y$ pour $\leq_1$. Si $x$ et $y$ sont $\{z\}$-\'equivalents  alors par d\'efinition, les restrictions a $\{x,z\}$ et $\{y,z\}$ sont isomorphes, donc les deux ordres coïncident sur
$\{x,z\}$ si et seulement s'ils coincident sur  $\{y,z\}$. Donc  si  $z$ est entre $x$ et $y$ pour $\leq_1$ alors ou bien les deux ordres coïncident sur $\{x,z\}$  et dans ce cas ils coïncident sur $\{y, z\}$  ou bien ils sont oppos\'es sur les deux, donc dans tous les cas  $z$ est entre $x$ et $y$ pour $\leq_2$.

Inversement, supposons que tout \'el\'ement $z$ qui est entre $x$ et $y$ pour l'un des ordres l'est pour l'autre ordre et montrons que $x$ et $y$ sont $1$-\'equivalents. Soit $z\in V\setminus\{x,y\}$. Alors $z$ est entre $x$ et $y$ modulo $\leq_1$ si et seulement si $z$ est entre $x$ et $y$ modulo $\leq_2$ et si $z$ n'est pas entre $x$ et $y$ (modulo $\leq_1$ et $\leq_2$) alors $z\leq_i x$ si et seulement si $z\leq_i y$ et $x\leq_i z$ si et seulement si $y\leq_i z$ pour $i\in\{1,2\}$. Donc les restrictions de $\mathcal B$ \`a $\{x,z\}$ et $\{y,z\}$ sont isomorphes. Il s'ensuit que $x$ et $y$ sont $1$-\'equivalents.
\end{proof}

\vspace{2mm}

Monteil et Pouzet\index{Monteil et Pouzet} dans leur article \cite{mont-pou} ont \'etudi\'e la d\'ecomposition monomorphe d'une bicha\^{i}ne. Ils ont d\'efini la relation d'\'equivalence suivante.

\vspace{2mm}

Soit $\mathcal B:= (V, \leq_1, \leq_2)$ une bicha\^{i}ne.
Deux \'el\'ements $x,y\in V$ sont \emph{\'equivalents}, fait not\'e $x\equiv y$,  si pour tout $t\in V$,
$t$ est entre $x,y$ modulo $\leq_1$ si et seulement si $t$ est entre $x,y$ modulo $\leq_2$.
Si $B_{\leq_i}(x,y)$ d\'esigne le sous-ensemble des \'el\'ements $z\in V$ qui sont entre $x$ et $y$ modulo $\leq_i$ pour $i\in\{1,2\}$, alors  $x$ et $y$ sont \'equivalents si et seulement si  $B_{\leq_1}(x,y)=B_{\leq_2}(x,y)$.

\begin{lemma}\label{lem:equivalence-bichaine} (voir \cite{mont-pou})

\begin{enumerate}
%\item The  relation $\equiv$ is an equivalence relation.
\item Sur chaque classe d'\'equivalence, les deux ordres coïncident ou sont oppos\'es.
\item \label{statement:four}Soient $x,y,z, t\in V$. Si  $x{\leq_1}z{\leq_1}y{\leq_1}t$ et si de plus $x\equiv y$ et $z\equiv t$, alors $x,y, z$ et $t$ sont \'equivalents.
\item Toute classe d'\'equivalence est une union de composantes monomorphes.
    \end{enumerate}
\end{lemma}

\vspace{1mm}

\noindent Ainsi, d'apr\`es le Lemme \ref{lem:equiv-bichaine}, pour deux \'el\'ements $x,y\in V$ nous avons
$$x\simeq_{1, \mathcal B}y \Leftrightarrow x\equiv y$$
Nous avons le r\'esultat suivant

\begin{lemma}\label{lem:bichaine partiefinie}
Soient   $\mathcal B:=(V, \leq_1, \leq_2)$ une bicha\^{\i}ne, $F$ une partie finie de $V$ d'au moins trois \'el\'ements et $x,y\in V\setminus F$. Si $\mathcal B_{\restriction_{ \{x\}\cup F'}}$ et  $\mathcal B_{\restriction_{ \{y\}\cup F'}}$ sont isomorphes pour chaque partie stricte $F'$ de $F$ alors $\mathcal B_{\restriction_{ \{x\}\cup F}}$ et  $\mathcal B_{\restriction_{ \{y\}\cup F}}$ sont isomorphes.
\end{lemma}

\begin{proof}
D'apr\`es le Lemme \ref{lem:equiv-bichaine} ci-dessus, tout \'el\'ement $z$  qui est entre $x$ et $y$ pour l'un des ordres est entre $x$ et $y$ pour l'autre. En prenant des parties $F'$ \`a deux \'el\'ements contenues dans l'intervalle commun d\'etermin\'e par $x$ et $y$, nous montrons qu'en fait, ou bien ces  deux ordres coïncident sur cet intervalle ou bien ils sont oppos\'es et ceci est suffisant.
\end{proof}

\vspace{1mm}

Du Lemme \ref{lem:bichaine partiefinie} nous pouvons d\'eduire directement

\begin{theorem}
Les relations d'\'equivalence $\simeq_{\leq 2,\mathcal B}$ et $\simeq_{\mathcal B}$ coïncident pour une bicha\^{i}ne $\mathcal B$.
\end{theorem}

        \subsection{Le cas des graphes}

Soit $G$ un  graphe non dirig\'e (sym\'etrique).
\begin{lemma}\label{lem:equiv-graphe}
Deux \'el\'ements $x, y\in V(G)$ sont \'equivalents si et seulement s'ils forment un intervalle de $G$.
 \end{lemma}
\begin{proof}
 La condition est \'evidemment suffisante; pour en voir la n\'ecessit\'e, noter que s'ils ne forment pas un intervalle alors il existe $z\in V(G)$  qui n'est pas li\'e de la m\^eme fa\c{c}on avec $x$ qu'avec $y$ ($z$ poss\`ede un seul voisin parmi $x$ et $y$) et donc $x$ et $y$ ne seraient pas \'equivalents.
\end{proof}

\begin{theorem}\label{coro:graphe}
Pour un graphe $G$ les relations $\simeq_{1,G}$ et $\simeq_{G}$ coïncident. %La $1$-\'equivalence coincide avec l'\'equivalence.
Dans ce cas, les classes d'\'equivalence sont les intervalles maximaux qui sont des cliques ou des ind\'ependants.
\end{theorem}

\subsection{Le cas des tournois}\label{subsec:cas tournois}

\begin{proposition}
Si $T$ est un tournoi,\index{tournoi}  la $\leq3$-\'equivalence est \'egale \`a l'\'equivalence.
\end{proposition}

Ce r\'esultat se d\'eduit trivialement du Lemme \ref{faibleseparation} ci-dessous.

\vspace{2mm}

Pour cela introduisons quelques  d\'efinitions et notations.
 Un {\it diamant},\index{diamant} resp. un   {\it double diamant},\index{diamant!double diamant}  est un tournoi  obtenu en remplaçant un sommet d'un tournoi \`a deux sommets, resp. d'un tournoi acyclique \`a trois sommets,   par un $3$-cycle (un cycle \`a trois \'el\'ements).  Un double diamant  est  autodual  si et seulement si le sommet du milieu du tournoi acyclique \`a trois sommets est remplac\'e par un $3$-cycle.
Soit %\footnote{Cette notation signifie $D:=1\underset{C_3}\oplus 1\underset{C_3}\oplus C_3$, la somme lexicographique des tournois $1, 1, C_3$ index\'ee par $C_3$.}
$D:= C_3(1,1, C_3)$ (on substitue  un $3$ -cycle \`a un \'el\'ement d'un $3$-cycle); soit  $B:=C_3(1, 1, 2)$ (on substitue \`a un \'el\'ement d'un $3$ cycle  un tournoi \`a deux \'el\'ements).

\vspace{2mm}

Notons que dans $D$ les deux \'el\'ements du $3$-cycle qui n'ont pas \'et\'e touch\'es par la substitution ne sont pas $3$-\'equivalents mais sont $2$-\'equivalents, de m\^eme dans le cas du double diamant autodual.

\vspace{2mm}

Boudabbous et Pouzet\index{Boudabbous et Pouzet} dans leur article \cite{Bou-Pouz} ont \'etudi\'e la d\'ecomposition monomorphe d'un tournoi.
Nous rappelons deux de leurs r\'esultats (Lemmes 9 et 13 de \cite{Bou-Pouz}).

\begin{lemma}\label{separation}(voir \cite{Bou-Pouz})

Deux sommets  $x,y$  d'un tournoi $T$ ne sont pas dans un m\^eme intervalle acyclique %\footnote{Un intervalle acyclique de $T$ est un intervalle de $T$ qui ne contient aucun cycle.}
de $T$ si et seulement si $x$ et $y$ sont distincts et ou bien:
\begin{enumerate}
\item[(i)] $x$ et $y$  appartiennent \`a un $3$-cycle, ou bien
\item[(ii)]  $x$ et $y$ appartiennent \`a un diamant, ou bien
\item[(iii)] $x$ et $y$ appartiennent \`a un double diamant  autodual.
 \end{enumerate}
\end{lemma}

\begin{lemma}\label{lem:interv-acyclique} \cite{Bou-Pouz}

Soit $T$ un tournoi et $A$ un sous-ensemble de $V(T)$.
\begin{enumerate}
\item Si $A$ est un intervalle acyclique\index{intervalle!acyclique} %\footnote{L'intervalle est d\'efini dans la section \ref{subsec:som-lex et decomp} du chapitre \ref{sect:str.rela.bin}. Un intervalle acyclique est un intervalle qui ne comporte pas de cycles.}
    alors $A$ est un bloc monomorphe.
    \item \label{statement:two}Si $A$ est une composante monomorphe alors ou bien $A$ est un intervalle qui est un $3$-cycle, ou bien $A=\{a,b\}$ et $A\cup C(\{a,b\})$\footnote{$C(\{a,b\})$ est l'ensemble des sommets $x$ tel que $\{a,b,x\}$ est un $3$-cycle de $T$.} est un intervalle, ou bien $A$ est une composante acyclique\index{composante!acyclique} de $T$ (une composante acyclique est un intervalle acyclique maximal pour l'inclusion) de $T$.
\end{enumerate}
\end{lemma}

Nous avons le lemme suivant:
\begin{lemma} \label{faibleseparation}
Deux sommets $x,y$ d'un tournoi $T$ ne sont pas dans une m\^eme classe d'equivalence si et seulement s'ils sont distincts et ou bien
\begin{enumerate}
\item[1)] $x$, $y$ sont  deux sommets d'un seul $3$-cycle d'un $B$, ou bien
\item[2)] $x$, $y$ appartiennent \`a un bord d'un diamant (non dans le cycle), ou bien
\item[3)] $x$, $y$ sont les extr\'emit\'es d'un double diamant autodual, ou bien
\item[4)] $x$, $y$ sont les sommets d'un $D$ appartenant  \`a trois de ses $3$-cycles.
\end{enumerate}
\end{lemma}

\begin{proof}
La condition suffisante se v\'erifie facilement. Pour la condition n\'ecessaire,
si $x$ et $y$ ne sont pas \'equivalents alors il existe une partie $F\subseteq V(T)\setminus\{x,y\}$ telle que les restrictions de $T$ \`a $F\cup\{x\}$ et $F\cup\{y\}$ ne sont pas isomorphes. Consid\'erons une telle partie $F$ de cardinalit\'e minimale.
Dans le cas d'un tournoi,  $F$ a au moins deux \'el\'ements (car un tournoi est $2$-monomorphe). Nous pouvons  remplacer $T$ par $T':= T_{\restriction_ {F\cup \{x, y\}}}$. Comme  $F\cup \{x\}$ et $F\cup \{y\}$ ne sont pas isomorphes, au moins un de ces deux ensembles contient un $3$-cycle.

\textbf{a)} Si $F$ a deux \'el\'ements, alors un seul de ces ensembles contient un $3$-cycle. %sinon ils seraient isomorphes
  Posons $F':=F\cup \{x,y\}$.\\
 - Si $F'$  ne contient qu'un seul $3$-cycle alors $F'$ est un diamant dont $x,y$ est un bord et on est en \emph{2)}. \\
  - Si $F'$ contient un autre $3$-cycle alors $F'$ est comme en \emph{1)}.

\textbf{b)}  Supposons que $F$ a au moins trois \'el\'ements.
\vspace{1mm}

\textbf{Cas 1:} Il existe un $3$-cycle dans $F'$ qui passe par $x$ et $y$. Soit $C$ l'ensemble des sommets $z\in F$ tels que $x,y,z$ forment un  $3$-cycle.  Soit $E:=F\setminus C$.

\begin{claim}\label{cla:claimvide} $E$ est vide.
\end{claim}
\begin{proofclaim}
Montrons que si $E$ est non vide alors $C\cup \{x,y\}$ est un intervalle de $T'$ ce qui, en raison du Lemme \ref{intervalle} constitue une contradiction. Pour cela, supposons $E$ non vide. Soit $u\in E$ et $z\in C$.

\begin{fact}  Si $\{x,y,z\}$ n'est pas un intervalle de $\{x,y,z,u\}$ alors il existe un $3$-cycle contenant $u$, $z$  et un seul des sommets $x$, $y$.
\end{fact}
\noindent {\bf  Preuve.} Si $\{x,y,z\}$ n'est pas un intervalle de $\{x,y,z,u\}$ alors $\{x,y,z,u\}$ n'est pas un diamant donc contient un $3$-cycle contenant $u$. Comme $u\not \in C$ ce $3$-cycle ne peut contenir $x$ et $y$, donc  il contient $z$  et un seul des sommets $x$, $y$.
\hfill $\blacklozenge$

\begin{fact} Ce   $3$-cycle ne peut contenir  $y$.
\end{fact}

\noindent {\bf  Preuve.}
Sinon, $\{u, z, y\}$ et $\{x, y, z\}$ \'etant deux $3$-cycles, %car $z\in C$
 $\{x,z,u\}$ ne peut-\^etre un $3$-cycle. Compte tenu du r\^ole sym\'etrique jou\'e par $x$ et $y$, ce $3$-cycle ne peut contenir $y$.
\hfill $\blacklozenge$

\vspace{2mm}

De ces deux faits, nous d\'eduisons que $x$ et $y$ ne sont pas $\{z,u\}$ \'equivalents, contredisant la minimalit\'e de $F$.
 Ainsi $\{x,y,z\}$ est un intervalle de $\{x,y,z,u\}$.   De ceci d\'ecoule que $C\cup \{x,y\}$ est un intervalle de $F'$, comme requis.
D'o\`u $E=\varnothing$, ceci termine la preuve de l'Affirmation \ref{cla:claimvide}.
\end{proofclaim}

\vspace{1mm}

Pour conclure, les seuls $3$-cycles inclus dans $F'$ passent par $x$ et $y$ ou ne les contiennent pas. S'il n'y a pas de $3$-cycle ne contenant ni $x$ ni $y$ alors $F\cup \{x\}$ et $F\cup \{y\}$  sont acycliques donc $x$ et $y$ sont \'equivalents dans $T'$, ce qui n'est pas le cas. Donc, il y a un $3$-cycle ne contenant ni $x$ ni $y$. Ce $3$-cycle \'etant inclus dans  $C$, il forme  un $D$ avec $x$ et $y$ et nous sommes dans le cas \emph{4)}.

\vspace{1mm}

\textbf{Cas 2:} Il ne passe aucun $3$-cycle ni par $x$ ni par $y$ dans $F'$. Donc tous les $3$-cycles de $T'$ ont leurs sommets dans $F$. Consid\'erons les composantes fortement connexes de $T'$, deux sommets appartiennent \`a une m\^eme composante fortement connexe s'ils sont \'egaux ou appartiennent \`a un m\^eme cycle (orient\'e). Comme $T'$ est un tournoi, deux sommets sont dans une m\^eme composante fortement connexe s'il appartiennent \`a un $3$-cycle.
Comme il n'y a aucun cycle contenant $x$ et aucun cycle contenant $y$, la composante de $x$ est r\'eduite \`a $x$ et celle de $y$ est r\'eduite \`a $y$. L'ensemble de ces composantes fortement connexes forme un ordre total dans lequel les composantes de $x$ et de $y$ ne se suivent pas (il y a des composantes entre les deux) car sinon $x$ et $y$ seraient \'equivalents dans $T'$. Les composantes fortement connexes qui sont entre celles de $x$ et de $y$ ne sont pas toutes des singletons, car sinon $x$ et $y$ seraient \'equivalents. Donc, il y a un $3$-cycle entre les deux, c'est \`a dire que $x$ et $y$ sont les extr\'emit\'es d'un double diamant autodual et on est en \emph{3)}.
\end{proof}

\begin{theorem} Soit $T$ un tournoi. Le nombre de classes de $2$-\'equivalence  est fini si et seulement si le nombre de classes d'\'equivalence est fini.
\end{theorem}

\begin{proof} Si le nombre de classes d'\'equivalence est fini, il est \'evident que le nombre de classe de $2$-\'equivalence est fini car chaque classe de $2$-\'equivalence est une union de classes d'\'equivalence. Supposons que le nombre de classes de $2$-\'equivalence soit fini et \'egal \`a $k$. Notons d'abord qu'une classe d'\'equivalence \'etant une composante monomorphe (Lemme \ref{lem:equi-blocmono}) sa forme est donn\'ee par l'assertion \ref{statement:two} du Lemme \ref{lem:interv-acyclique}.

Soit $C$  une classe de $2$-\'equivalence. Supposons que cette classe a au moins $4$ \'el\'ements. Il est facile de voir que $C$ est acyclique, en effet, s'il existe un $3$-cycle, disons $x,y,z$ dans $C$, alors tout \'el\'ement $t\in C\setminus\{x,y,z\}$ serait non $2$-\'equivalent aux sommets $x, y$ et $z$.

  Plus fortement, aucun  $3$-cycle n'intersecte $C$ en deux \'el\'ements distincts. Supposons qu'un  $3$-cycle $xy z$ intersecte $C$ en $x,y$ avec $x<y ~(mod T)$. Soit $t$ un \'el\'ement de $C$ distinct de $x,y$. Comme $t$ est $2$-\'equivalent a $x$, alors les restrictions de $T$ \`a $\{x,y,z\}$ et $\{y,z,t\}$ sont isomorphes. Mais ceci est impossible puisque la restriction \`a $\{x,z,t\}$ n'est pas un cycle et la restriction \`a $\{x,y,z\}$ est un $3$-cycle.

Si $C$ est un intervalle de $T$, alors $C$ est une classe d'\'equivalence en vertu du Lemme \ref{lem:interv-acyclique}. Supposons que $C$ ne soit pas un intervalle de $T$. Dans ce cas, il existe deux \'el\'ements distincts $x, y$ de $C$ qui ne sont pas dans une m\^eme composante acyclique. Appliquons le Lemme \ref{separation} de Boudabbous-Pouzet \cite{Bou-Pouz}. Le premier cas n'est pas possible d'apr\`es ce qu'on vient de voir. Le second cas n'est pas possible non plus, en effet si  $x, y$ est dans un diamant il est sur un bord, et si $y$ est la pointe du diamant et $x, z, t$ forment le $3$-cycle alors $x$ et $y$ ne sont pas $\{ z,t\}$-\'equivalents. Reste donc le troisi\`eme cas. Soient donc $u,v,w $ formant un $3$-cycle tels que $x$ soit la pointe d'un diamant positif  et $y$ la pointe d'un diamant n\'egatif (un diamant est positif si le sommet $x$ qui n'appartient pas au $3$-cycle est tel que $u<x~(mod T)$ pour tout sommet $u$ du $3$-cycle, si $x<u~(mod T)$ pour tout sommet $u$ du $3$-cycle, le diamant est n\'egatif). Ce $3$-cycle est forc\'ement disjoint de $C$. Soit $z$ un autre \'el\'ement de $C$, alors $z, u,v,w$ forment un diamant (en effet, comme $x$ et $z$ sont $2$-\'equivalents et que $x, u,v$ ne forment pas un $3$-cycle, alors $z, u,v$ n'est pas un cycle, de m\^eme on obtient que $z, u,w$ et $z, v,w$ non plus; ainsi $z,u,v,w$ ne contient qu'un $3$-cycle, c'est donc un diamant. La classe $C$ se partage en deux intervalles $C^{-}$ et $C^+$ form\'es chacun  des \'el\'ements qui forment avec $u,v,w$ un diamant n\'egatif, resp positif.

\noindent Disons que deux \'el\'ements $x$ et $y$ de $C$ sont \emph{fortement \'equivalents} si en supposant $x<y$ on ne peut les compl\'eter par un $3$-cycle  tels que $x$ soit la pointe d'un diamant positif et $y$ la pointe d'un diamant n\'egatif. Nous obtenons une relation d'\'equivalence dont les classes sont des intervalles de $T$ donc des classes d'\'equivalence. Donc la classe $C$ donne lieu \`a plusieurs classes d'\'equivalence d\'efinies par des $3$-cycles, mais comme ces $3$-cycles sont dans des classes diff\'erentes, il y a au plus $k-1$ tels cycles ce qui donne au plus $k$ intervalles de $C$. Donc le nombre de classes d'\'equivalence est fini.
\end{proof}

        \subsection{Le cas des graphes dirig\'es}\label{subsec:equiv-graphe-ordonne}

Soit $G=(V,E)$ un graphe dirig\'e\index{graphe!dirig\'e} sans boucle, donc $G$ est une relation binaire irr\'eflexive.

Le r\'esultat suivant se d\'eduit trivialement du Lemme \ref{lem:separationgraphe} ci-dessous.
\begin{theorem}\label{theo.equiv-graphe}
Si $G$ est un graphe dirig\'e sans boucles alors la $\leq 3$-\'equivalence est \'egale \`a l'\'equivalence.
\end{theorem}

Rappelons que ce r\'esultat a \'et\'e d\'emontr\'e ind\'ependemment par Boudabbous (2013) \cite{boudabbous}.

\vspace{2mm}

Nous introduisons quelques notations. Soit $W:= \{0,1\}^2$. Pour tout \'el\'ement $u:=(a,b)\in W$, posons $\overline u:= (b,a)$. Posons \'egalement $0:=(0,0)$.

Soit $G:= (V,E)$ un graphe dirig\'e sans boucle. Soit $\rho$ la fonction caract\'eristique de $E$ et $d$ l'application de $V^2$ dans $W$ d\'efinie par
%Soit $d$  l'application de $V^2$ dans $\{0,1\}^2$ d\'efinie par
$$d(x,y):=(\rho(x,y),\rho(y,x))\;\text{ pour tous }x,y\in V.$$
Il est clair que $d(x,y)= \overline {d(y,x)}$ pour tous $x,y\in V$. Il est \'egalement clair que
 $$d(x,y)=(1,0)\Leftrightarrow d(y,x)=(0,1)$$ et que $$d(x,y)=d(y,x)\Leftrightarrow \rho(x,y)=\rho(y,x).$$

Nous disons que la paire $\{x,y\}$ est \emph{vide} si $d(x,y)=0$, \emph{asym\'etrique} si $d(x,y)\in\{(1,0),(0,1)\}$ et \emph{sym\'etrique} si $d(x,y)=(1,1)$.

\smallskip

Soient $V'=\{x,y,z,t\}\subseteq V$, supposons qu'il existe dans $G$ un chemin de longueur $2$, que nous noterons $P_3$, reliant les sommets $x,y,z$ de $V'$ avec $d(x,z)=d(z,y)=(1,0)$.

\smallskip

\begin{figure}[h]
\centering
\psset{unit=1cm}
\begin{pspicture}(-7.5,-4)(7,2)
\psdots[dotsize=5pt](-6,0.5)(-4,0.5)(-1,0.5)(1,0.5)(4,0.5)(6,0.5)(-6,1.5)(-4,1.5)(-1,1.5)(1,1.5)(4,1.5)(6,1.5)
\psdots[dotsize=5pt](-3.5,-1.5)(-1.5,-1.5)(1.5,-1.5)(3.5,-1.5)(-3.5,-2.5)(-1.5,-2.5)(1.5,-2.5)(3.5,-2.5)
\psline[linewidth=0.4pt, arrowlength=4, arrowsize=1mm]{->}(-6,1.5)(-6,0.8) %fleche de x \`a y dans $Q_1$
\psline[linewidth=0.4pt](-6,0.8)(-6,0.5)
\psline[linewidth=0.4pt, arrowlength=4, arrowsize=1mm]{->}(-6,0.5)(-4.8,0.5)    %fleche de y \`a z dans $Q_1$
\psline[linewidth=0.4pt](-4.8,0.5)(-4,0.5)
\psline[linewidth=0.4pt, arrowlength=4, arrowsize=1mm]{->}(-4,0.5)(-4,1.2)    %fleche de z \`a t dans $Q_1$
\psline[linewidth=0.4pt](-4,1.2)(-4,1.5)
\psline[linewidth=0.4pt, arrowlength=4, arrowsize=1mm]{->}(-4,1.5)(-5.2,1.5)  %fleche de t \`a x dans $Q_1$
\psline[linewidth=0.4pt](-6,1.5)(-5.2,1.5)
\psline[linewidth=0.4pt, arrowlength=4, arrowsize=1mm]{->}(-1,1.5)(-1,0.8)    %fleche de x \`a y dans $Q_2$
\psline[linewidth=0.4pt](-1,0.8)(-1,0.5)
\psline[linewidth=0.4pt, arrowlength=4, arrowsize=1mm]{->}(-1,0.5)(0.2,0.5)   %fleche de y \`a z dans $Q_2$
\psline[linewidth=0.4pt](1,0.5)(0.2,0.5)
\psline[linewidth=0.4pt](-1,1.5)(1,1.5)(1,0.5)                %fleche de z \`a t et de t \`a x dans $Q_2$
\psline[linewidth=0.4pt, arrowlength=4, arrowsize=1mm]{->}(4,1.5)(4,0.8)      %fleche de x \`a y dans $Q_3$
\psline[linewidth=0.4pt](4,0.8)(4,0.5)
\psline[linewidth=0.4pt, arrowlength=4, arrowsize=1mm]{->}(4,0.5)(5.2,0.5)   %fleche de y \`a z dans $Q_3$
\psline[linewidth=0.4pt](5.2,0.5)(6,0.5)
\psline[linewidth=0.4pt, arrowlength=4, arrowsize=1mm]{->}(-3.5,-1.5)(-3.5,-2.2) %fleche de x \`a y dans $Q_4$
\psline[linewidth=0.4pt](-3.5,-2.2)(-3.5,-2.5)
\psline[linewidth=0.4pt, arrowlength=4, arrowsize=1mm]{->}(-3.5,-2.5)(-2.3,-2.5)    %fleche de y \`a z dans $Q_4$
\psline[linewidth=0.4pt](-2.3,-2.5)(-1.5,-2.5)
\psline[linewidth=0.4pt, arrowlength=4, arrowsize=1mm]{->}(-3.5,-1.5)(-2.3,-1.5) %fleche de x \`a t dans $Q_4$
\psline[linewidth=0.4pt](-2.3,-1.5)(-1.5,-1.5)
\psline[linewidth=0.4pt, arrowlength=4, arrowsize=1mm]{->}(-1.5,-1.5)(-1.5,-2.2)    %fleche de t \`a z dans $Q_4$
\psline[linewidth=0.4pt](-1.5,-2.2)(-1.5,-2.5)
\psline[linewidth=0.4pt, arrowlength=4, arrowsize=1mm]{->}(1.5,-1.5)(1.5,-2.2) %fleche de x \`a y dans $Q_5$
\psline[linewidth=0.4pt](1.5,-2.2)(1.5,-2.5)
\psline[linewidth=0.4pt, arrowlength=4, arrowsize=1mm]{->}(1.5,-2.5)(2.7,-2.5)    %fleche de y \`a z dans $Q_5$
\psline[linewidth=0.4pt](2.7,-2.5)(3.5,-2.5)
\psline[linewidth=0.4pt, arrowlength=4, arrowsize=1mm]{->}(1.5,-1.5)(2.7,-1.5) %fleche de x \`a t dans $Q_5$
\psline[linewidth=0.4pt](2.7,-1.5)(3.5,-1.5)
\psline[linewidth=0.4pt, arrowlength=4, arrowsize=1mm]{->}(3.5,-2.5)(3.5,-1.8)    %fleche de z \`a t dans $Q_5$
\psline[linewidth=0.4pt](3.5,-1.8)(3.5,-1.5)
%\psline[linewidth=0.3pt](-2,0)(-1,-1.5)(0,0)
%\psline[linewidth=0.3pt](2,0)(2,-1.5)
\uput{0.3}[u](-6,1.5){$x$}
\uput{0.3}[u](-1,1.5){$x$}
\uput{0.3}[u](4,1.5){$x$}
\uput{0.3}[u](-3.5,-1.5){$x$}
\uput{0.3}[u](1.5,-1.5){$x$}
\uput{0.7}[d](-5,0.5){$Q_1$}
\uput{0.3}[d](-6,0.5){$z$}
\uput{0.3}[d](-1,0.5){$z$}
\uput{0.3}[d](4,0.5){$z$}
\uput{0.3}[d](-3.5,-2.5){$z$}
\uput{0.3}[d](1.5,-2.5){$z$}
\uput{0.7}[d](0,0.5){$Q_2$}
\uput{0.3}[u](-4,1.5){$t$}
\uput{0.3}[u](1,1.5){$t$}
\uput{0.3}[u](6,1.5){$t$}
\uput{0.3}[u](-1.5,-1.5){$t$}
\uput{0.3}[u](3.5,-1.5){$t$}
\uput{0.7}[d](5,0.5){$Q_3$}
\uput{0.7}[d](-2.5,-2.5){$Q_4$}
\uput{0.7}[d](2.5,-2.5){$Q_5$}
\uput{0.3}[d](-4,0.5){$y$}
\uput{0.3}[d](1,0.5){$y$}
\uput{0.3}[d](6,0.5){$y$}
\uput{0.3}[d](-1.5,-2.5){$y$}
\uput{0.3}[d](3.5,-2.5){$y$}
\end{pspicture}
\caption{\label{type-carre} Les carr\'es de type $Q_i$ pour $1\leq i\leq 5$.}
\end{figure}
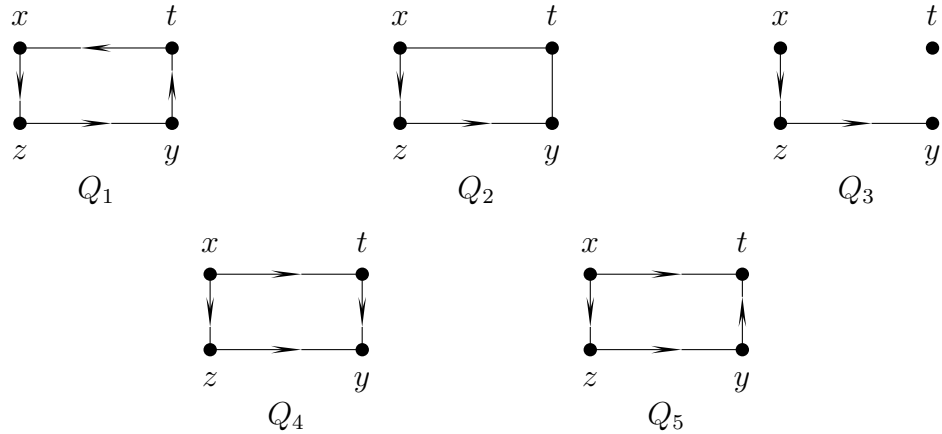

Nous disons que $V'$ est un \emph{carr\'e}\index{carr\'e} de type:
\begin{enumerate}
\item  $Q_1$ si $d(y,t)=d(t,x)=(1,0)$ (les sommets de $V'$ forment un circuit de longueur $4$).
\item  $Q_2$ si $d(y,t)=d(t,x)=(1,1)$.
\item  $Q_3$ si $d(y,t)=d(t,x)=(0,0)$.
\item  $Q_4$ si $d(y,t)=d(t,x)=(0,1)$.
\item  $Q_5$ si $d(y,t)=d(x,t)=(1,0)$.
\end{enumerate}
Les paires de sommets $\{x,y\}$ et $\{z,t\}$ sont les \emph{diagonales}\index{carr\'e!diagonale d'un -} de ces carr\'es, la paire $\{x,y\}$ \'etant la \emph{diagonale principale}\index{carr\'e!diagonale principale d'un -} pour chacun des carr\'es $Q_i~(2\leq i\leq 5)$. Les valeurs de $d$ en ces diagonales ne sont pas n\'ecessairement \'egales. Ces carr\'es sont repr\'esent\'es sur la \figurename~\ref{type-carre}.

\medskip

\begin{figure}[h]
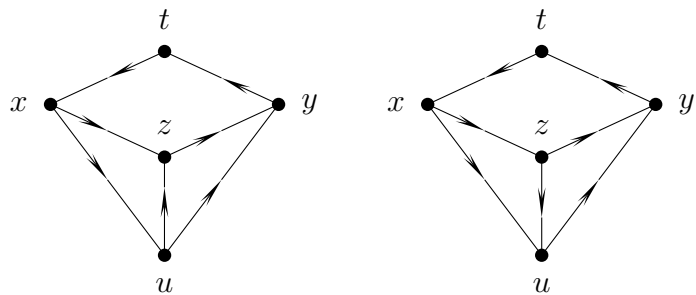

\centering
\input{prisme1}\qquad \input{prisme2}
\caption{\label{prisme} Repr\'esentation des prismes.}
\end{figure}

De m\^eme, consid\'erons $V''=\{x,y,z,t,u\}$ un sous-ensemble de $V$. Nous disons que les  sommets de $V''$ forment un \emph{prisme}\index{prisme} $P$ si $\{x,y,z,t\}$ et $\{x,u,y,t\}$ sont deux carr\'es de type $Q_1$, avec $d(x,z)=d(z,y)=d(y,t)=d(t,x)=(1,0)$, et la paire $\{z,u\}$ est asym\'etrique. Les diagonales des deux carr\'es $\{x,y,z,t\}$ et $\{x,u,y,t\}$ sont les diagonales du prisme (dans notre cas se sont les paires $\{x,y\}$, $\{z,t\}$ et $\{t,u\}$), la diagonale commune (ie, $\{x,y\}$) \'etant la diagonale principale de $P$ (voir \figurename~\ref{prisme}).

Nous avons de mani\`ere \'evidente le fait suivant:
\begin{fact}\label{fact:$1$-equiv}
Deux sommets $x,y$ de $V$ ne sont pas $1$-\'equivalents si et seulement s'il existe $z\in V\setminus\{x,y\}$ tel que les paires $\{x,z\}$ et $\{y,z\}$ ne sont pas toutes les deux vides, sym\'etriques ou asym\'etriques.
\end{fact}

%\vspace{2mm}
%%%%%%%%%%%%%%%%%%%%%%%%%%%%%%%%%%%%%%%%%%%%%%%%%%%%%%%%%%%%%%%

Nous avons:
\begin{lemma}\label{lem:$2$-equiv}
Soit $G:=(V,E)$ un graphe dirig\'e ayant au moins quatre sommets. Deux sommets $x,y$ de $V$ qui sont $1$-\'equivalents ne sont pas $2$-\'equivalents si et seulement s'ils sont distincts et ou bien
\begin{enumerate}
\item[(1)] $\{x,y\}$ est une diagonale d'un carr\'e de type $Q_1$ dont l'autre diagonale est asym\'etrique, ou bien
\item[(2)] $\{x,y\}$ est la diagonale principale d'un carr\'e de type $Q_2$ dont l'autre diagonale n'est pas sym\'etrique, ou bien
\item[(3)] $\{x,y\}$ est la diagonale principale d'un carr\'e de type $Q_3$ dont l'autre diagonale n'est pas vide, ou bien
\item[(4)] $\{x,y\}$ est la diagonale principale d'un carr\'e de type $Q_4$ dont l'autre diagonale n'est pas asym\'etrique, ou bien
\item[(5)] $\{x,y\}$ est une diagonale d'un carr\'e de type $Q_5$ tel que l'autre diagonale est vide, sym\'etrique ou asym\'etrique et dans ce dernier cas elle est dans un $3$-cycle.
\end{enumerate}
\end{lemma}

\begin{proof}
La condition suffisante se v\'erifie facilement. Pour la condition n\'ecessaire, soit $F=\{z,z'\}$ une partie \`a deux \'el\'ements de $V\setminus\{x,y\}$ telle que les restrictions de $G$ \`a $F\cup\{x\}$ et $F\cup\{y\}$ ne sont pas isomorphes. Posons $V'=F\cup\{x,y\}$ et $G'=G_{\restriction_{V'}}$.

Comme $x$ et $y$ sont $1$-\'equivalents alors d'apr\'es le Fait \ref{fact:$1$-equiv}, les paires $\{x,z\}$ et $\{y,z\}$ (resp. $\{x,z'\}$ et $\{y,z'\}$) sont toutes les deux soit vides, soit sym\'etriques soit asym\'etriques.

Si $d(x,z)=d(y,z)$ et $d(x,z')=d(y,z')$ les restrictions de $G$ \`a $F\cup\{x\}$ et $F\cup\{y\}$ sont isomorphes, donc nous avons

\qquad - ou bien $\{x,z\}$ et $\{y,z\}$ sont asym\'etriques et $d(x,z)=d(z,y)$,

\qquad - ou bien $\{x,z'\}$ et $\{y,z'\}$ sont asym\'etriques et $d(x,z')=d(z',y)$.

Supposons, sans perte de g\'en\'eralit\'e, que $d(x,z)=d(z,y)\in\{(1,0),(0,1)\}$. Il existe alors, dans $G'$, un chemin $P_3$ d'extr\'emit\'es $x,y$. Nous pouvons encore supposer, sans perte de g\'en\'eralit\'e, que $d(x,z)=d(z,y)=(1,0)$, nous avons alors deux cas:
\vspace{1mm}

\textbf{Cas 1:} $d(x,z')=d(z',y)\in\{(1,0),(0,1)\}$, alors

\qquad \textbf{1-a)} si $d(x,z')=d(z',y)=(1,0)$ alors la paire $\{z,z'\}$ ne peut pas \^etre asym\'etrique car sinon $x$ et $y$ seraient $\{z,z'\}$-\'equivalents et on est dans le cas \emph{(4)} du lemme,

\qquad \textbf{1-b)} si $d(x,z')=d(z',y)=(0,1)$, alors $\{z,z'\}$ ne peut pas \^etre vide ou sym\'etrique car sinon $x$ et $y$ seraient $\{z,z'\}$-\'equivalents. Donc $\{z,z'\}$ est asym\'etrique et on est dans le cas \emph{(1)} du lemme.
\vspace{1mm}

\textbf{Cas 2:} $d(x,z')=d(y,z')$, nous avons alors quatre cas:

\qquad \textbf{2-a)} $d(x,z')=d(y,z')=(1,0)$ alors $d(z,z')\neq (1,0)$ car sinon $x$ et $y$ seraient $\{z,z'\}$-\'equivalents et on est dans le cas \emph{(5)} du lemme.

\qquad \textbf{2-b)} $d(x,z')=d(y,z')=(0,1)$, nous avons la m\^eme conclusion que dans le cas \textbf{2-a)} ci-dessus.

\qquad \textbf{2-c)} $d(x,z')=d(y,z')=(0,0)$ alors $d(z,z')\neq (0,0)$ car sinon $x$ et $y$ seraient $\{z,z'\}$-\'equivalents et on est dans le cas \emph{(3)} du lemme.

\qquad \textbf{2-d)} $d(x,z')=d(y,z')=(1,1)$ alors $d(z,z')\neq (1,1)$ car sinon $x$ et $y$ seraient $\{z,z'\}$-\'equivalents et on est dans le cas \emph{(2)} du lemme.
\end{proof}

\begin{remark}
Les cas 1) et 2) du Lemme \ref{faibleseparation} se retrouvent respectivement dans le cas (1) pour $\{x,y\}$ assym\'etrique et le cas (5) pour $d(x,y)=(1,0)$ et $\{z,z'\}$ assym\'etrique du Lemme \ref{lem:$2$-equiv}.
\end{remark}

\begin{lemma}\label{lem:separationgraphe}
Soit $G:=(V,E)$ un graphe dirig\'e. Deux sommets $x,y$ de $V$ qui sont $\leq 2$-\'equivalents ne sont pas \'equivalents si et seulement s'ils sont distincts et ou bien
\begin{enumerate}
\item[(1)] il existe trois sommets de $V$ qui forment avec l'un des sommets $x,y$ un diamant positif et avec l'autre un diamant n\'egatif, ou bien
\item[(2)] $\{x,y\}$ est la diagonales principale d'un prisme dont chacune des autres diagonales est soit vide soit sym\'etrique.
\end{enumerate}
\end{lemma}

\begin{proof}
La condition suffisante est facile \`a v\'erifier. Pour la condition n\'ecessaire, soit $F$ une partie minimale de $V\setminus\{x,y\}$ telle que les restrictions de $G$ \`a $F\cup\{x\}$ et $F\cup\{y\}$ ne soient pas isomorphes. Les sommets $x$ et $y$ \'etant $\leq 2$-\'equivalents, $F$ a au moins trois \'el\'ements. D'apr\`es le Fait \ref{fact:$1$-equiv}, pour tout $z\in F$, les paires $\{x,z\}$ et $\{y,z\}$  sont toutes les deux soit vides, soit sym\'etriques soit asym\'etriques. Si $d(x,z)=d(y,z)$ pour tout $z\in F$ les restrictions de $G$ \`a $F\cup\{x\}$ et $F\cup\{y\}$ sont isomorphes, ce qui contredit le choix de $F$. Donc il existe $z\in F$ tel que $\{x,z\}$ et $\{y,z\}$ sont asym\'etriques et $d(x,z)=d(z,y)$.
Posons $$\begin{array}{l}
D=\{z\in F : d(x,z)\neq d(y,z)\},\\
%H=\{z\in F/d(x,z)=d(y,z)=(1,1)\},\\
%T=\{z\in F/d(x,z)=d(y,z)=(0,0)\},\\
%\text{et }U=\{z\in F/ \{x,z\} \text{ et }\{y,z\}\text{ sont asym\'etriques et }d(x,z)=d(y,z)\}.\\
\text{et }D'=F\setminus D.\\
\end{array}$$
Il est clair que $D\neq\varnothing$.% et que $\{D,H,T,U\}$ est une partition de $F$. %Nous avons trois cas.

\begin{claim}
$D'=\varnothing$.
\end{claim}
\begin{proofclaim}
Supposons $D'\neq\varnothing$. Soit $z\in D$ et $h\in D'$. Alors comme $x$ et $y$ sont  $2$-\'equivalents, les restrictions de $G$ \`a $\{x,z,h\}$ et \`a $\{y,z,h\}$ sont isomorphes. On constate facilement qu'il n'y a qu'un isomorphisme, celui qui fixe $h$, envoie $x$ sur $z$ et $z$ sur $y$. Par cons\'equent $d(z,h)=d(y,h)=d(x,h)$. Il s'ensuit que $D\cap \{x,y\}$ est un intervalle Fra\"{\i}ss\'e de $G_{\restriction_{F\cup \{x,y\}}}$. Ceci contredit le Lemme \ref{intervalle}.
%Si $H\cup T\cup U\neq\varnothing$ alors la $2$-\'equivalence de $x,y$ impose $d(z,h)=d(x,h)=d(y,h)$ pour tout $z\in D$ et $h\in H\cup T\cup U$ et donc $D\cup\{x,y\}$ est un intervalle de $G_{\restriction_{F\cup\{x,y\}}}$, ce qui contredit le Lemme \ref{intervalle}.
\end{proofclaim}

Donc $F=D$, % ce qui implique que $\vert D\vert\geq 3$,
nous avons alors deux cas.

\vspace{1mm}

\textbf{Cas 1:} $d(x,z)$ est constante pour tout $z\in F$. Dans ce cas, pour tous $z,z'\in F$, les sommets $x,y,z,z'$ forment un carr\'e de type $Q_4$ et la $2$-\'equivalence de $x,y$ impose \`a $\{z,z'\}$ d'\^etre asym\'etrique d'apr\`es le Lemme \ref{lem:$2$-equiv}, donc $G_{\restriction_F}$ est un tournoi.% De m\^eme les paires $\{z,h\}$ sont sym\'etriques pour tout $z\in D$ et tout $h\in H$ et les paires $\{z,t\}$ sont
\begin{claim}
$F$ poss\`ede un $3$-cycle.
\end{claim}
\begin{proofclaim}
Si $F$ n'a pas de $3$-cycle alors $G_{\restriction_F}$ est un ordre total strict. Comme $d(x,z)$ est constante, $G_{\restriction_{F\cup \{x\}}}$ est un ordre total strict. De m\^eme  $G_{\restriction_{F\cup \{y\}}}$. Et donc les restrictions de $G$ \`a $F\cup\{x\}$ et $F\cup\{y\}$ sont isomorphe, ce qui est impossible de par la d\'efinition de $F$.
\end{proofclaim}

Soient alors $a,b,c$ les sommets du $3$-cycle de $F$, alors les restrictions de $G$ \`a $\{a,b,c,x\}$ et $\{a,b,c,y\}$ sont des diamants, l'un positif et l'autre n\'egatif et on est dans le cas \emph{(1)} du Lemme \ref{lem:separationgraphe}.
\vspace{1mm}

\textbf{Cas 2:} il existe deux sommets $a,b\in F$ tels que $d(x,a)\neq d(x,b)$, ce qui signifie que $d(x,a)= d(b,x)$. Supposons, sans perte de g\'en\'eralit\'e que $d(x,a)= d(b,x)=(1,0)$. Donc les sommets $x,y,a,b$ forment un carr\'e de type $Q_1$ et ceci, en vertu du Lemme \ref{lem:$2$-equiv}, impose \`a $\{a,b\}$ d'\^etre vide ou sym\'etrique.

Comme $\vert F\vert\geq 3$, il existe $c\in F\setminus\{a,b\}$ tel que l'un parmi les deux sous-ensembles $\{x,y,a,c\}$ et $\{x,y,b,c\}$ forme un carr\'e de type $Q_1$ et l'autre un carr\'e de type $Q_4$. Plus pr\'ecis\'ement, $c$ v\'erifie ou bien $d(x,c)=d(x,a)$ ou bien $d(x,c)=d(x,b)$. Supposons, sans perte de g\'en\'eralit\'e, que $d(x,c)=d(x,a)$. Dans ce cas le sous-ensemble $\{x,y,a,c\}$ forme un carr\'e de type $Q_4$ et le sous-ensemble $\{x,y,b,c\}$ forme un carr\'e de type $Q_1$. Comme $x$ et $y$ sont $2$-\'equivalents, le Lemme \ref{lem:$2$-equiv} impose \`a $\{a,c\}$ d'\^etre asym\'etrique et \`a $\{b,c\}$ d'\^etre vide ou sym\'etrique. Les sommets $\{x,y,a,b,c\}$ forment alors un prisme de diagonale principale $\{x,y\}$ et chacune des deux autres diagonales $\{a,b\}$ et $\{b,c\}$ est vide ou sym\'etrique. On est dans le cas \emph{(2)} du Lemme \ref{lem:separationgraphe}.
%Si $\{x,y,a,c\}$ (resp. $\{x,y,b,c\}$) est un carr\'e de type $Q_1$ alors la $2$-\'equivalence de $x,y$ impose \`a $\{a,c\}$ (resp. $\{b,c\}$) d'\^etre vide ou sym\'etrique et \`a $\{b,c\}$ (resp. $\{a,c\}$) d'\^etre asym\'etrique. Dans les deux cas les sommets $\{x,y,a,b,c\}$ forment un prisme de diagonale principale $\{x,y\}$ et chacune des autres diagonales ($\{a,b\}$ et $\{a,c\}$ dans le premier cas et $\{a,b\}$ et $\{b,c\}$ dans le second) \'etant vide ou sym\'etrique. On est dans le cas \emph{(2)} du Lemme \ref{lem:separationgraphe}.
\end{proof}

\begin{remark}
Les cas 3) et 4) du Lemme \ref{faibleseparation} se retrouvent dans le cas (1) du Lemme \ref{lem:separationgraphe}.
\end{remark}
%%%%%%%%%%%%%%%%%%%%%%%%%%%%%%%%%%%%%%%%%%         Rajout\'e apr\`es la soutenance %%%%%%%%%%%%%%%%%%%%%%%%%%%%%%%%%%%%%%%%%%%%

 \subsection{Le cas des structures binaires}\label{subsec:equiv-struct-binaire}

Dans cette sous-section nous utilisons les d\'efinitions et notations introduites dans la sous-section \ref{subsec:equiv-graphe-ordonne}.

\begin{theorem}
Si $\mathcal R$ est une structure binaire de type $k$, $k\in\mathbb N$, alors la $\leq 3$-\'equivalence est \'egale \`a l'\'equivalence.
\end{theorem}

Ce r\'esultat se d\'eduit trivialement du Lemme \ref{lem:separbinaire} ci dessous.

\vspace{1mm}

Soit $\mathcal R:=(E,\rho _{1},\dots,\rho _{k})$ une structure binaire de type $k$, $k\geq 2$. Pour tout $i\in\{1,\dots,k\}$, identifions la relation $\rho_i$ \`a sa fonction caract\'eristique et posons $$d_i(x,y):=(\rho_i(x,y),\rho_i(y,x)) \text{ pour tous } x,y\in E.$$
Pour $u:= (u_i)_{i= 1,\dots, k}$ posons $\overline u:= (\overline{u}_i)_{i= 1,\dots, k}$. Soit $d$ l'application de $V^2$ dans $W^k$ d\'efinie par $d(x,y):= (d_i(x,y))_{i= 1,\dots, k}$. Nous avons $d(x,y)= \overline{d(y,x)}$.
\vspace{2mm}

Nous avons de mani\`ere \'evidente le lemme suivant:

\begin{lemma}\label{lem:1equivbinaire}
Soit $\mathcal R:=(E,\rho _{1},\dots,\rho _{k})$ une structure binaire de type $k$ o\`u $E$ poss\`ede au moins trois \'el\'ements. Deux \'el\'ements  $x$ et $y$ de $E$ sont $1$-\'equivalents si et seulement si $d(x,z)=d(y,z)$ ou $d(x,z)=d(z,y)$ pour tout $z\in E\setminus\{x,y\}$.
\end{lemma}

\vspace{1mm}

Posons pour tout $i\in\{1,\dots,k\}$, $\mathcal R_i:=(E,\rho_i)$. Pour tout $i$, $\mathcal R_i$ est un graphe dirig\'e.  Remarquons que si $x\simeq_{n,\mathcal R} y$, $n\in\NN^{\star}$, alors pour tout $i\in\{1,\dots,k\}$ on a  $x\simeq_{n,\mathcal R_i} y$. La r\'eciproque est fausse.

\medskip

Nous avons:

\begin{lemma}\label{lemm:2equiv-binaire}
Soit $\mathcal R:=(E,\rho _{1},\dots,\rho _{k})$ une structure binaire de type $k$ o\`u $E$ poss\`ede au moins quatre \'el\'ements. Deux sommets $x,y$ de $E$ qui sont $1$-\'equivalents ne sont pas $2$-\'equivalents si et seulement si ils sont distincts et ou bien
\begin{enumerate}
\item il existe $i\in\{1,\dots,k\}$ tel que $x\not\simeq_{2,\mathcal R_i} y$, ou bien
\item $x\simeq_{2,\mathcal R_i} y,~\forall i\in\{1,\dots,k\}$ et il existe $i\neq j$ et $z,z'\in E\setminus\{x,y\}$ tels que l'on ait une seule des situations suivantes.
    \begin{enumerate}
    \item Dans $\mathcal R_i$, $\{x,y,z,z'\}$ est un carr\'e de type $Q_1$ de diagonale principale $\{x,y\}$ avec $\{z,z'\}$ vide ou sym\'etrique et dans $\mathcal R_j$, $\{x,y,z,z'\}$ est un carr\'e de type $Q_4$ de diagonale principale $\{x,y\}$ avec $\{z,z'\}$ asym\'etrique.
    \item Pour chacune des structures $\mathcal R_i$ et $\mathcal R_j$, $F=\{x,y,z,z'\}$ est un carr\'e de type $Q_4$ de diagonale principale $\{x,y\}$ avec $\{z,z'\}$ asym\'etrique tel que l'isomorphisme $f$ qui envoie ${\mathcal R_i}_{\restriction_F}$ sur ${\mathcal R_j}_{\restriction_F}$ fixe une seule des paires $\{x,y\}$ ou $\{z,z'\}$.
    \end{enumerate}
\end{enumerate}
\end{lemma}

\begin{proof}
On v\'erifie facilement que ces conditions sont suffisante. Voyons qu'elles sont n\'ecessaires. Soient donc $x$ et $y$ deux \'el\'ements $1$-\'equivalents qui ne sont pas $2$-\'equivalents. Nous avons deux cas.

\vspace{1mm}

\textbf{Cas 1:} Il existe $i\in\{1,\dots,k\}$ tel que $x\not\simeq_{2,\mathcal R_i} y$. C'est la condition \emph{1.} du Lemme \ref{lemm:2equiv-binaire}. %Donc $x$ et $y$ v\'erifient  le Lemme \ref{lem:$2$-equiv} pour la relation $\mathcal R_i$.
\vspace{1mm}

\textbf{Cas 2:} Pour tout $i\in\{1,\dots,k\}$, on a $x\simeq_{2,\mathcal R_i} y$. Comme $x$ et $y$ sont $1$-\'equivalents, d'apr\`es le Lemme \ref{lem:1equivbinaire} nous avons trois cas:
\medskip

\qquad \textbf{2.1:} Pour tout $z\in E\setminus\{x,y\}$, $d_i(x,z)=d_i(y,z),~\forall i\in\{1,\dots,k\}$. Dans ce cas, pour toute partie $A$ \`a deux \'el\'ements de $E\setminus\{x,y\}$, l'application $f$ qui envoie $x$ sur $y$ et fixe les \'el\'ements de $A$ est un isomorphisme de $\mathcal R_{\restriction_{A\cup\{x\}}}$ sur $\mathcal R_{\restriction_{A\cup\{y\}}}$. Il s'ensuit que $x$ et $y$ sont $2$-\'equivalents, ce qui contredit l'hypoth\`ese, donc ce cas ne peut pas se pr\'esenter.
\medskip

\qquad \textbf{2.2:} Il existe $z\neq z'$ tels que

\qquad $d_i(x,z)=d_i(y,z),~\forall i\in\{1,\dots,k\}$,

\qquad et $d_i(x,z')=d_i(z',y)\neq d_i(y,z'),~\forall i\in\{1,\dots,k\}$.

Donc les paires $\{x,z'\}$ et $\{y,z'\}$ sont asym\'etriques et pour tout $i\in\{1,\dots,k\}$, la restriction de $\mathcal R_i$ \`a $\{x,y,z,z'\}$ a l'une des formes donn\'ees dans la \figurename~\ref{cas binaire-preuve}.

\begin{figure}[h]
\centering
\psset{unit=1cm}
\begin{pspicture}(-7.5,-3)(7,3)
\psdots[dotsize=5pt](-6,0.5)(-2,0.5)(2,0.5)(6,0.5)(-7,1.5)(-5,1.5)(-3,1.5)(-1,1.5)(1,1.5)(3,1.5)(5,1.5)(7,1.5)(-6,2.5)(-2,2.5)(2,2.5)(6,2.5)
\psdots[dotsize=5pt](-6,-0.5)(-2,-0.5)(2,-0.5)(6,-0.5)(-7,-1.5)(-5,-1.5)(-3,-1.5)(-1,-1.5)(1,-1.5)(3,-1.5)(5,-1.5)(7,-1.5)(-6,-2.5)(-2,-2.5)(2,-2.5)(6,-2.5)
\psline[linewidth=0.4pt, arrowlength=4, arrowsize=1mm]{->}(-7,1.5)(-6.5,2) %fleche de x \`a z' dans (i)
\psline[linewidth=0.4pt](-6.5,2)(-6,2.5)
\psline[linewidth=0.4pt, arrowlength=4, arrowsize=1mm]{->}(-6,2.5)(-5.5,2)    %fleche de z' \`a y dans (i)
\psline[linewidth=0.4pt](-5.5,2)(-5,1.5)
\psline[linewidth=0.4pt, arrowlength=4, arrowsize=1mm]{->}(-3,1.5)(-2.5,2)    %fleche de x \`a z' dans (ii)
\psline[linewidth=0.4pt](-2.5,2)(-2,2.5)
\psline[linewidth=0.4pt, arrowlength=4, arrowsize=1mm]{->}(-2,2.5)(-1.5,2)    %fleche de x \`a z' dans (ii)
\psline[linewidth=0.4pt](-1.5,2)(-1,1.5)
\psline[linewidth=0.4pt](-3,1.5)(-2,0.5)(-1,1.5)
\psline[linewidth=0.4pt, arrowlength=4, arrowsize=1mm]{->}(1,1.5)(1.5,2)    %fleche de x \`a z' dans (iii)
\psline[linewidth=0.4pt](1.5,2)(2,2.5)
\psline[linewidth=0.4pt, arrowlength=4, arrowsize=1mm]{->}(2,2.5)(2.5,2)    %fleche de z' \`a y dans (iii)
\psline[linewidth=0.4pt](2.5,2)(3,1.5)
\psline[linewidth=0.4pt, arrowlength=4, arrowsize=1mm]{->}(1,1.5)(1.5,1)    %fleche de x \`a z dans (iii)
\psline[linewidth=0.4pt](1.5,1)(2,0.5)
\psline[linewidth=0.4pt, arrowlength=4, arrowsize=1mm]{->}(3,1.5)(2.5,1)    %fleche de y \`a z dans (iii)
\psline[linewidth=0.4pt](2.5,1)(2,0.5)
\psline[linewidth=0.4pt, arrowlength=4, arrowsize=1mm]{->}(5,1.5)(5.5,2)    %fleche de x \`a z' dans (iv)
\psline[linewidth=0.4pt](5.5,2)(6,2.5)
\psline[linewidth=0.4pt, arrowlength=4, arrowsize=1mm]{->}(6,2.5)(6.5,2)    %fleche de z' \`a y dans (iv)
\psline[linewidth=0.4pt](6.5,2)(7,1.5)
\psline[linewidth=0.4pt, arrowlength=4, arrowsize=1mm]{->}(6,0.5)(5.5,1)    %fleche de x \`a z dans (iv)
\psline[linewidth=0.4pt](5.5,1)(5,1.5)
\psline[linewidth=0.4pt, arrowlength=4, arrowsize=1mm]{->}(6,0.5)(6.5,1)    %fleche de y \`a z dans (iii)
\psline[linewidth=0.4pt](6.5,1)(7,1.5)
\uput{0.3}[l](-7,1.5){$x$}
\uput{0.3}[l](-3,1.5){$x$}
\uput{0.3}[l](1,1.5){$x$}
\uput{0.3}[l](5,1.5){$x$}
\uput{0.7}[l](-6,2.5){\emph{(i)}}
\uput{0.3}[r](-5,1.5){$y$}
\uput{0.3}[r](-1,1.5){$y$}
\uput{0.3}[r](3,1.5){$y$}
\uput{0.3}[r](7,1.5){$y$}
\uput{0.7}[l](-2,2.5){\emph{(ii)}}
\uput{0.2}[u](-6,2.5){$z'$}
\uput{0.2}[u](-2,2.5){$z'$}
\uput{0.2}[u](2,2.5){$z'$}
\uput{0.2}[u](6,2.5){$z'$}
\uput{0.7}[l](2,2.5){\emph{(iii)}}
\uput{0.7}[l](6,2.5){\emph{(iv)}}
\uput{0.2}[l](-6,0.5){$z$}
\uput{0.2}[l](-2,0.5){$z$}
\uput{0.2}[l](2,0.5){$z$}
\uput{0.2}[l](6,0.5){$z$}
\psline[linewidth=0.4pt, arrowlength=4, arrowsize=1mm]{->}(-6,-0.5)(-6.5,-1) %fleche de z' \`a x dans (i')
\psline[linewidth=0.4pt](-6.5,-1)(-7,-1.5)
\psline[linewidth=0.4pt, arrowlength=4, arrowsize=1mm]{->}(-5,-1.5)(-5.5,-1)    %fleche de y \`a z' dans (i')
\psline[linewidth=0.4pt](-5.5,-1)(-6,-0.5)
\psline[linewidth=0.4pt, arrowlength=4, arrowsize=1mm]{->}(-1,-1.5)(-1.5,-1)    %fleche de z' \`a x dans (ii')
\psline[linewidth=0.4pt](-1.5,-1)(-2,-0.5)
\psline[linewidth=0.4pt, arrowlength=4, arrowsize=1mm]{->}(-2,-0.5)(-2.5,-1)    %fleche de y \`a z' dans (ii')
\psline[linewidth=0.4pt](-2.5,-1)(-3,-1.5)
\psline[linewidth=0.4pt](-3,-1.5)(-2,-2.5)(-1,-1.5)
\psline[linewidth=0.4pt, arrowlength=4, arrowsize=1mm]{->}(2,-0.5)(1.5,-1)    %fleche de x \`a z' dans (iii')
\psline[linewidth=0.4pt](1.5,-1)(1,-1.5)
\psline[linewidth=0.4pt, arrowlength=4, arrowsize=1mm]{->}(3,-1.5)(2.5,-1)    %fleche de z' \`a y dans (iii')
\psline[linewidth=0.4pt](2.5,-1)(2,-0.5)
\psline[linewidth=0.4pt, arrowlength=4, arrowsize=1mm]{->}(3,-1.5)(2.5,-2)    %fleche de x \`a z dans (iii')
\psline[linewidth=0.4pt](2.5,-2)(2,-2.5)
\psline[linewidth=0.4pt, arrowlength=4, arrowsize=1mm]{->}(1,-1.5)(1.5,-2)    %fleche de y \`a z dans (iii')
\psline[linewidth=0.4pt](1.5,-2)(2,-2.5)
\psline[linewidth=0.4pt, arrowlength=4, arrowsize=1mm]{->}(6,-0.5)(5.5,-1)    %fleche de x \`a z' dans (iv')
\psline[linewidth=0.4pt](5.5,-1)(5,-1.5)
\psline[linewidth=0.4pt, arrowlength=4, arrowsize=1mm]{->}(7,-1.5)(6.5,-1)    %fleche de z' \`a y dans (iv')
\psline[linewidth=0.4pt](6.5,-1)(6,-0.5)
\psline[linewidth=0.4pt, arrowlength=4, arrowsize=1mm]{->}(6,-2.5)(5.5,-2)    %fleche de x \`a z dans (iv')
\psline[linewidth=0.4pt](5.5,-2)(5,-1.5)
\psline[linewidth=0.4pt, arrowlength=4, arrowsize=1mm]{->}(6,-2.5)(6.5,-2)    %fleche de y \`a z dans (iv')
\psline[linewidth=0.4pt](6.5,-2)(7,-1.5)
\uput{0.3}[l](-7,-1.5){$x$}
\uput{0.3}[l](-3,-1.5){$x$}
\uput{0.3}[l](1,-1.5){$x$}
\uput{0.3}[l](5,-1.5){$x$}
\uput{0.7}[l](-6,-0.5){\emph{(i')}}
\uput{0.3}[r](-5,-1.5){$y$}
\uput{0.3}[r](-1,-1.5){$y$}
\uput{0.3}[r](3,-1.5){$y$}
\uput{0.3}[r](7,-1.5){$y$}
\uput{0.7}[l](-2,-0.5){\emph{(ii')}}
\uput{0.3}[r](-6,-0.5){$z'$}
\uput{0.3}[r](-2,-0.5){$z'$}
\uput{0.3}[r](2,-0.5){$z'$}
\uput{0.3}[r](6,-0.5){$z'$}
\uput{0.7}[l](2,-0.5){\emph{(iii')}}
\uput{0.7}[l](6,-0.5){\emph{(iv')}}
\uput{0.3}[d](-6,-2.5){$z$}
\uput{0.3}[d](-2,-2.5){$z$}
\uput{0.3}[d](2,-2.5){$z$}
\uput{0.3}[d](6,-2.5){$z$}
\end{pspicture}
\caption{\label{cas binaire-preuve}}
\end{figure}

Posons $F=\{x,y,z,z'\}$. Dans les graphes \emph{(i)} et \emph{(i')} de la \figurename~\ref{cas binaire-preuve}, $F$ est un carr\'e de type $Q_3$. Dans \emph{(ii)} et \emph{(ii')}, $F$ est un carr\'e de type $Q_2$ et dans \emph{(iii)}, \emph{(iv)}, \emph{(iii')} et \emph{(iv')}, $F$ est un carr\'e de type $Q_5$.

Comme $x\simeq_{2,\mathcal R_i}y$ pour tout $i\in\{1,\dots,k\}$, le Lemme \ref{lem:$2$-equiv} permet d'avoir la relation entre $z$ et $z'$ dans chacun des cas. On obtient:

$\bullet~~\{z,z'\}$ est vide dans \emph{(i)} et \emph{(i')}.

$\bullet~~\{z,z'\}$ est sym\'etrique dans \emph{(ii)} et \emph{(ii')}.

$\bullet~~\{z,z'\}$ est asym\'etrique et n'appartient pas \`a un $3$-cycle dans les autres cas.

Dans chacun des cas ci-dessus, l'application $f$ qui fixe $z$ et envoie $x$ sur $z'$ et $z'$ sur $y$ est un isomorphisme de $\mathcal R_{\restriction_{\{x,z,z'\}}}$ sur $\mathcal R_{\restriction_{\{y,z,z'\}}}$. Ceci signifie que $x$ et $y$ sont $2$-\'equivalents, ce qui contredit l'hypoth\`ese. Ce cas ne peut donc pas se pr\'esenter.
\medskip

\qquad \textbf{2.3:} Pour tout $z\in E\setminus\{x,y\}$, $d_i(x,z)=d_i(z,y)\neq d_i(y,z),~\forall i\in\{1,\dots,k\}$. Donc, pour tous $z,z'\in E\setminus\{x,y\}$, $z\neq z'$, la restriction de $\mathcal R_i$ \`a $F=\{x,y,z,z'\}$ a l'une des formes donn\'ees dans la \figurename~\ref{cas2 binaire-preuve} suivante.

\begin{figure}[h]
\centering
\psset{unit=1cm}
\begin{pspicture}(-7.5,0)(7,3)
\psdots[dotsize=5pt](-6,0.5)(-2,0.5)(2,0.5)(6,0.5)(-7,1.5)(-5,1.5)(-3,1.5)(-1,1.5)(1,1.5)(3,1.5)(5,1.5)(7,1.5)(-6,2.5)(-2,2.5)(2,2.5)(6,2.5)
\psline[linewidth=0.4pt, arrowlength=4, arrowsize=1mm]{->}(-7,1.5)(-6.5,2) %fleche de x \`a z' dans (j_1)
\psline[linewidth=0.4pt](-6.5,2)(-6,2.5)
\psline[linewidth=0.4pt, arrowlength=4, arrowsize=1mm]{->}(-6,2.5)(-5.5,2)    %fleche de z' \`a y dans (j_1)
\psline[linewidth=0.4pt](-5.5,2)(-5,1.5)
\psline[linewidth=0.4pt, arrowlength=4, arrowsize=1mm]{->}(-5,1.5)(-5.5,1)    %fleche de y \`a z dans (j_1)
\psline[linewidth=0.4pt](-5.5,1)(-6,0.5)
\psline[linewidth=0.4pt, arrowlength=4, arrowsize=1mm]{->}(-6,0.5)(-6.5,1)    %fleche de z \`a x dans (j_1)
\psline[linewidth=0.4pt](-6.5,1)(-7,1.5)
\psline[linewidth=0.4pt, arrowlength=4, arrowsize=1mm]{->}(-3,1.5)(-2.5,2)    %fleche de x \`a z' dans (j_2)
\psline[linewidth=0.4pt](-2.5,2)(-2,2.5)
\psline[linewidth=0.4pt, arrowlength=4, arrowsize=1mm]{->}(-2,2.5)(-1.5,2)    %fleche de z' \`a y dans (j_2)
\psline[linewidth=0.4pt](-1.5,2)(-1,1.5)
\psline[linewidth=0.4pt, arrowlength=4, arrowsize=1mm]{->}(-3,1.5)(-2.5,1)    %fleche de x \`a z dans (j_2)
\psline[linewidth=0.4pt](-2.5,1)(-2,0.5)
\psline[linewidth=0.4pt, arrowlength=4, arrowsize=1mm]{->}(-2,0.5)(-1.5,1)    %fleche de z \`a y dans (j_2)
\psline[linewidth=0.4pt](-1.5,1)(-1,1.5)
\psline[linewidth=0.4pt, arrowlength=4, arrowsize=1mm]{->}(2,2.5)(1.5,2)    %fleche de z' \`a x dans (j_3)
\psline[linewidth=0.4pt](1.5,2)(1,1.5)
\psline[linewidth=0.4pt, arrowlength=4, arrowsize=1mm]{->}(3,1.5)(2.5,2)    %fleche de y \`a z' dans (j_3)
\psline[linewidth=0.4pt](2.5,2)(2,2.5)
\psline[linewidth=0.4pt, arrowlength=4, arrowsize=1mm]{->}(3,1.5)(2.5,1)    %fleche de y \`a z dans (j_3)
\psline[linewidth=0.4pt](2.5,1)(2,0.5)
\psline[linewidth=0.4pt, arrowlength=4, arrowsize=1mm]{->}(2,0.5)(1.5,1)    %fleche de z \`a x dans (j_3)
\psline[linewidth=0.4pt](1.5,1)(1,1.5)
\psline[linewidth=0.4pt, arrowlength=4, arrowsize=1mm]{->}(6,2.5)(5.5,2)    %fleche de z' \`a x dans (j_4)
\psline[linewidth=0.4pt](5.5,2)(5,1.5)
\psline[linewidth=0.4pt, arrowlength=4, arrowsize=1mm]{->}(7,1.5)(6.5,2)    %fleche de y \`a z' dans (j_4)
\psline[linewidth=0.4pt](6.5,2)(6,2.5)
\psline[linewidth=0.4pt, arrowlength=4, arrowsize=1mm]{->}(5,1.5)(5.5,1)    %fleche de x \`a z dans (j_4)
\psline[linewidth=0.4pt](5.5,1)(6,0.5)
\psline[linewidth=0.4pt, arrowlength=4, arrowsize=1mm]{->}(6,0.5)(6.5,1)    %fleche de z \`a y dans (j_4)
\psline[linewidth=0.4pt](6.5,1)(7,1.5)
\uput{0.3}[l](-7,1.5){$x$}
\uput{0.3}[l](-3,1.5){$x$}
\uput{0.3}[l](1,1.5){$x$}
\uput{0.3}[l](5,1.5){$x$}
\uput{0.7}[l](-6,2.5){$(j_1)$}
\uput{0.3}[r](-5,1.5){$y$}
\uput{0.3}[r](-1,1.5){$y$}
\uput{0.3}[r](3,1.5){$y$}
\uput{0.3}[r](7,1.5){$y$}
\uput{0.7}[l](-2,2.5){$(j_2)$}
\uput{0.2}[u](-6,2.5){$z'$}
\uput{0.2}[u](-2,2.5){$z'$}
\uput{0.2}[u](2,2.5){$z'$}
\uput{0.2}[u](6,2.5){$z'$}
\uput{0.7}[l](2,2.5){$(j_3)$}
\uput{0.7}[l](6,2.5){$(j_4)$}
\uput{0.2}[d](-6,0.5){$z$}
\uput{0.2}[d](-2,0.5){$z$}
\uput{0.2}[d](2,0.5){$z$}
\uput{0.2}[d](6,0.5){$z$}
\end{pspicture}
\caption{\label{cas2 binaire-preuve}}
\end{figure}

Dans les graphes $(j_1)$ et $(j_4)$ de la \figurename~\ref{cas2 binaire-preuve}, $F$ est un carr\'e de type $Q_1$ et dans $(j_2)$ et $(j_3)$, $F$ est un carr\'e de type $Q_4$.

Comme $x\simeq_{2,\mathcal R_i}y$ pour tout $i\in\{1,\dots,k\}$, le Lemme \ref{lem:$2$-equiv} donne:

$\bullet~~\{z,z'\}$ est vide ou sym\'etrique dans $(j_1)$ et $(j_4)$.

$\bullet~~\{z,z'\}$ est asym\'etrique dans $(j_2)$ et $(j_3)$.

Les configurations de la \figurename~\ref{cas2 binaire-preuve} donnent alors lieu aux huit configurations de la \figurename~\ref{cas3 binaire-preuve} suivante.

\begin{figure}[h]
\centering
\psset{unit=1cm}
\begin{pspicture}(-7.5,-3)(7,3)
\psdots[dotsize=5pt](-6,0.5)(-2,0.5)(2,0.5)(6,0.5)(-7,1.5)(-5,1.5)(-3,1.5)(-1,1.5)(1,1.5)(3,1.5)(5,1.5)(7,1.5)(-6,2.5)(-2,2.5)(2,2.5)(6,2.5)
\psdots[dotsize=5pt](-6,-0.5)(-2,-0.5)(2,-0.5)(6,-0.5)(-7,-1.5)(-5,-1.5)(-3,-1.5)(-1,-1.5)(1,-1.5)(3,-1.5)(5,-1.5)(7,-1.5)(-6,-2.5)(-2,-2.5)(2,-2.5)(6,-2.5)
\psline[linewidth=0.4pt, arrowlength=4, arrowsize=1mm]{->}(-7,1.5)(-6.5,2) %fleche de x \`a z' dans (j_1)
\psline[linewidth=0.4pt](-6.5,2)(-6,2.5)
\psline[linewidth=0.4pt, arrowlength=4, arrowsize=1mm]{->}(-6,2.5)(-5.5,2)    %fleche de z' \`a y dans (j_1)
\psline[linewidth=0.4pt](-5.5,2)(-5,1.5)
\psline[linewidth=0.4pt, arrowlength=4, arrowsize=1mm]{->}(-5,1.5)(-5.5,1)    %fleche de y \`a z dans (j_1)
\psline[linewidth=0.4pt](-5.5,1)(-6,0.5)
\psline[linewidth=0.4pt, arrowlength=4, arrowsize=1mm]{->}(-6,0.5)(-6.5,1)    %fleche de z \`a x dans (j_1)
\psline[linewidth=0.4pt](-6.5,1)(-7,1.5)
\psline[linewidth=0.4pt, arrowlength=4, arrowsize=1mm]{->}(-3,1.5)(-2.5,2)    %fleche de x \`a z' dans (j_2)
\psline[linewidth=0.4pt](-2.5,2)(-2,2.5)
\psline[linewidth=0.4pt, arrowlength=4, arrowsize=1mm]{->}(-2,2.5)(-1.5,2)    %fleche de z' \`a y dans (j_2)
\psline[linewidth=0.4pt](-1.5,2)(-1,1.5)
\psline[linewidth=0.4pt, arrowlength=4, arrowsize=1mm]{->}(-3,1.5)(-2.5,1)    %fleche de x \`a z dans (j_2)
\psline[linewidth=0.4pt](-2.5,1)(-2,0.5)
\psline[linewidth=0.4pt, arrowlength=4, arrowsize=1mm]{->}(-2,0.5)(-1.5,1)    %fleche de z \`a y dans (j_2)
\psline[linewidth=0.4pt](-1.5,1)(-1,1.5)
\psline[linewidth=0.4pt, arrowlength=4, arrowsize=1mm]{->}(-2,0.5)(-2,1.5)    %fleche de z \`a z' dans (j_2)
\psline[linewidth=0.4pt](-2,1.5)(-2,2.5)
\psline[linewidth=0.4pt, arrowlength=4, arrowsize=1mm]{->}(2,2.5)(1.5,2)    %fleche de z' \`a x dans (j_3)
\psline[linewidth=0.4pt](1.5,2)(1,1.5)
\psline[linewidth=0.4pt, arrowlength=4, arrowsize=1mm]{->}(3,1.5)(2.5,2)    %fleche de y \`a z' dans (j_3)
\psline[linewidth=0.4pt](2.5,2)(2,2.5)
\psline[linewidth=0.4pt, arrowlength=4, arrowsize=1mm]{->}(3,1.5)(2.5,1)    %fleche de y \`a z dans (j_3)
\psline[linewidth=0.4pt](2.5,1)(2,0.5)
\psline[linewidth=0.4pt, arrowlength=4, arrowsize=1mm]{->}(2,0.5)(1.5,1)    %fleche de z \`a x dans (j_3)
\psline[linewidth=0.4pt](1.5,1)(1,1.5)
\psline[linewidth=0.4pt, arrowlength=4, arrowsize=1mm]{->}(2,0.5)(2,1.5)    %fleche de z \`a z' dans (j_3)
\psline[linewidth=0.4pt](2,1.5)(2,2.5)
\psline[linewidth=0.4pt, arrowlength=4, arrowsize=1mm]{->}(6,2.5)(5.5,2)    %fleche de z' \`a x dans (j_4)
\psline[linewidth=0.4pt](5.5,2)(5,1.5)
\psline[linewidth=0.4pt, arrowlength=4, arrowsize=1mm]{->}(7,1.5)(6.5,2)    %fleche de y \`a z' dans (j_4)
\psline[linewidth=0.4pt](6.5,2)(6,2.5)
\psline[linewidth=0.4pt, arrowlength=4, arrowsize=1mm]{->}(5,1.5)(5.5,1)    %fleche de x \`a z dans (j_4)
\psline[linewidth=0.4pt](5.5,1)(6,0.5)
\psline[linewidth=0.4pt, arrowlength=4, arrowsize=1mm]{->}(6,0.5)(6.5,1)    %fleche de z \`a y dans (j_4)
\psline[linewidth=0.4pt](6.5,1)(7,1.5)
\uput{0.3}[l](-7,1.5){$x$}
\uput{0.3}[l](-3,1.5){$x$}
\uput{0.3}[l](1,1.5){$x$}
\uput{0.3}[l](5,1.5){$x$}
\uput{0.7}[l](-6,2.5){$(j_1')$}
\uput{0.3}[r](-5,1.5){$y$}
\uput{0.3}[r](-1,1.5){$y$}
\uput{0.3}[r](3,1.5){$y$}
\uput{0.3}[r](7,1.5){$y$}
\uput{0.7}[l](-2,2.5){$(j_2')$}
\uput{0.2}[u](-6,2.5){$z'$}
\uput{0.2}[u](-2,2.5){$z'$}
\uput{0.2}[u](2,2.5){$z'$}
\uput{0.2}[u](6,2.5){$z'$}
\uput{0.7}[l](2,2.5){$(j_3')$}
\uput{0.7}[l](6,2.5){$(j_4')$}
\uput{0.2}[l](-6,0.5){$z$}
\uput{0.2}[l](-2,0.5){$z$}
\uput{0.2}[l](2,0.5){$z$}
\uput{0.2}[l](6,0.5){$z$}
\psline[linewidth=0.4pt, arrowlength=4, arrowsize=1mm]{->}(-7,-1.5)(-6.5,-1) %fleche de x \`a z' dans (j1'')
\psline[linewidth=0.4pt](-6.5,-1)(-6,-0.5)
\psline[linewidth=0.4pt, arrowlength=4, arrowsize=1mm]{->}(-6,-0.5)(-5.5,-1)    %fleche de z' \`a y dans (j1'')
\psline[linewidth=0.4pt](-5.5,-1)(-5,-1.5)
\psline[linewidth=0.4pt, arrowlength=4, arrowsize=1mm]{->}(-5,-1.5)(-5.5,-2) %fleche de y \`a z dans (j1'')
\psline[linewidth=0.4pt](-5.5,-2)(-6,-2.5)
\psline[linewidth=0.4pt, arrowlength=4, arrowsize=1mm]{->}(-6,-2.5)(-6.5,-2)    %fleche de z \`a x dans (j1'')
\psline[linewidth=0.4pt](-6.5,-2)(-7,-1.5)
\psline[linewidth=0.4pt](-6,-0.5)(-6,-2.5)
\psline[linewidth=0.4pt, arrowlength=4, arrowsize=1mm]{->}(-3,-1.5)(-2.5,-1)    %fleche de x \`a z' dans (j2'')
\psline[linewidth=0.4pt](-2.5,-1)(-2,-0.5)
\psline[linewidth=0.4pt, arrowlength=4, arrowsize=1mm]{->}(-2,-0.5)(-1.5,-1)    %fleche de z' \`a y dans (j2'')
\psline[linewidth=0.4pt](-1.5,-1)(-1,-1.5)
\psline[linewidth=0.4pt, arrowlength=4, arrowsize=1mm]{->}(-3,-1.5)(-2.5,-2)    %fleche de x \`a z dans (j2'')
\psline[linewidth=0.4pt](-2.5,-2)(-2,-2.5)
\psline[linewidth=0.4pt, arrowlength=4, arrowsize=1mm]{->}(-2,-2.5)(-1.5,-2)    %fleche de z \`a y dans (j2'')
\psline[linewidth=0.4pt](-1.5,-2)(-1,-1.5)
\psline[linewidth=0.4pt, arrowlength=4, arrowsize=1mm]{->}(-2,-0.5)(-2,-1.5)    %fleche de z' \`a z dans (j2'')
\psline[linewidth=0.4pt](-2,-1.5)(-2,-2.5)
\psline[linewidth=0.4pt, arrowlength=4, arrowsize=1mm]{->}(2,-0.5)(1.5,-1)    %fleche de z' \`a x dans (j3'')
\psline[linewidth=0.4pt](1.5,-1)(1,-1.5)
\psline[linewidth=0.4pt, arrowlength=4, arrowsize=1mm]{->}(3,-1.5)(2.5,-1)    %fleche de y \`a z' dans (j3'')
\psline[linewidth=0.4pt](2.5,-1)(2,-0.5)
\psline[linewidth=0.4pt, arrowlength=4, arrowsize=1mm]{->}(2,-2.5)(1.5,-2)    %fleche de z \`a x dans (j3'')
\psline[linewidth=0.4pt](1.5,-2)(1,-1.5)
\psline[linewidth=0.4pt, arrowlength=4, arrowsize=1mm]{->}(3,-1.5)(2.5,-2)    %fleche de y \`a z dans (j3'')
\psline[linewidth=0.4pt](2.5,-2)(2,-2.5)
\psline[linewidth=0.4pt, arrowlength=4, arrowsize=1mm]{->}(2,-0.5)(2,-1.5)    %fleche de z' \`a z dans (j3'')
\psline[linewidth=0.4pt](2,-1.5)(2,-2.5)
\psline[linewidth=0.4pt, arrowlength=4, arrowsize=1mm]{->}(7,-1.5)(6.5,-1)    %fleche de y \`a z' dans (j4'')
\psline[linewidth=0.4pt](6.5,-1)(6,-0.5)
\psline[linewidth=0.4pt, arrowlength=4, arrowsize=1mm]{->}(6,-0.5)(5.5,-1)    %fleche de z' \`a x dans (j4'')
\psline[linewidth=0.4pt](5.5,-1)(5,-1.5)
\psline[linewidth=0.4pt, arrowlength=4, arrowsize=1mm]{->}(5,-1.5)(5.5,-2)    %fleche de x \`a z dans (j4'')
\psline[linewidth=0.4pt](5.5,-2)(6,-2.5)
\psline[linewidth=0.4pt, arrowlength=4, arrowsize=1mm]{->}(6,-2.5)(6.5,-2)    %fleche de z \`a y dans (j4'')
\psline[linewidth=0.4pt](6.5,-2)(7,-1.5)
\psline[linewidth=0.4pt](6,-2.5)(6,-0.5)
\uput{0.3}[l](-7,-1.5){$x$}
\uput{0.3}[l](-3,-1.5){$x$}
\uput{0.3}[l](1,-1.5){$x$}
\uput{0.3}[l](5,-1.5){$x$}
\uput{0.7}[l](-6,-0.5){$(j_1'')$}
\uput{0.3}[r](-5,-1.5){$y$}
\uput{0.3}[r](-1,-1.5){$y$}
\uput{0.3}[r](3,-1.5){$y$}
\uput{0.3}[r](7,-1.5){$y$}
\uput{0.7}[l](-2,-0.5){$(j_2'')$}
\uput{0.3}[r](-6,-0.5){$z'$}
\uput{0.3}[r](-2,-0.5){$z'$}
\uput{0.3}[r](2,-0.5){$z'$}
\uput{0.3}[r](6,-0.5){$z'$}
\uput{0.7}[l](2,-0.5){$(j_3'')$}
\uput{0.7}[l](6,-0.5){$(j_4'')$}
\uput{0.3}[d](-6,-2.5){$z$}
\uput{0.3}[d](-2,-2.5){$z$}
\uput{0.3}[d](2,-2.5){$z$}
\uput{0.3}[d](6,-2.5){$z$}
\end{pspicture}
\caption{\label{cas3 binaire-preuve}}
\end{figure}

Remarquons que les configurations $(j_1')$ et $(j_4')$ sont isomorphes ainsi que $(j_1'')$ et $(j_4'')$ et l'isomorphisme dans ces deux situations est le m\^eme, il permute $x$ et $y$ et fixe $z$ et $z'$. Si la restriction de $\mathcal R_i$ \`a $F$ a l'une des formes $(j_1')$, $(j_1'')$, $(j_4')$ et $(j_4'')$ pour tout $i\in\{1,\dots,k\}$ alors $x$ et $y$ sont $2$-\'equivalents, donc ceci ne peut pas avoir lieu.

Remarquons \'egalement que les configurations $(j_2')$, $(j_3')$, $(j_2'')$ et $(j_3'')$ sont isomorphes. Les isomorphismes les reliant deux \`a deux sont de trois types:

\qquad $-$ Soit il permute $x$ et $y$ et fixe $z$ et $z'$ comme c'est le cas pour l'isomorphisme reliant les graphes $(j_2')$ et $(j_3')$ et celui reliant les graphes $(j_2'')$ et $(j_3'')$.

\qquad $-$ Soit il fixe $x$ et $y$ et permute $z$ et $z'$ comme c'est le cas pour l'isomorphisme reliant les graphes $(j_2')$ et $(j_2'')$ et celui reliant les graphes $(j_3')$ et $(j_3'')$.

\qquad $-$ Soit il permute $x$ et $y$ et permute $z$ et $z'$ comme c'est le cas pour l'isomorphisme reliant les graphes $(j_2')$ et $(j_3'')$ et celui reliant les graphes $(j_2'')$ et $(j_3')$.
\medskip

On v\'erifie que si la restriction de $\mathcal R_i$ \`a $F$ est soit dans $\{(j_2'),(j_3'')\}$ pour tout $i\in\{1,\dots,k\}$, soit dans $\{(j_2''),(j_3')\}$ pour tout $i\in\{1,\dots,k\}$, alors $x$ et $y$ sont $2$-\'equivalents, donc ceci ne peut pas avoir lieu.
\medskip

Il s'ensuit qu'il existe $i\neq j$ et $z,z'\in E\setminus\{x,y\}$ tels que l'on ait l'une des deux situations suivantes:

\qquad $\centerdot$ ou bien la restriction de $\mathcal R_i$ \`a $F=\{x,y,z,z'\}$ est dans l'ensemble $\{(j_1'),(j_1''),(j_4'),(j_4'')\}$ et celle de $\mathcal R_j$ \`a $F$ dans l'ensemble $\{(j_2'),(j_2''),(j_3'),(j_3'')\}$ et on est dans le cas $2.(a)$ du Lemme \ref{lemm:2equiv-binaire}.

\qquad $\centerdot$ ou bien la paire de restrictions $\{{\mathcal R_i}_{\restriction_F},{\mathcal R_j}_{\restriction_F}\}$ est dans l'ensemble
$$\Big\{\{(j_2'),(j_3')\};\{(j_2''),(j_3'')\};\{(j_2'),(j_2'')\};\{(j_3'),(j_3'')\}\Big\}$$ et on est dans le cas $2.(b)$ du Lemme \ref{lemm:2equiv-binaire}. Ceci termine la preuve.
\end{proof}

\bigskip

Dans ce qui suit, nous allons g\'en\'eraliser le Lemme \ref{lem:separationgraphe} au cas des structures binaires de type $k$. Rappelons que $W=\{0,1\}^2$ et consid\'erons $W^k$ comme un ensemble muni d'une involution not\'ee $\overline{(.)}$. Donc une structure binaire de type $k$, $\mathcal R:=(E,\rho_1,\dots,\rho_k)$ peut-\^etre, simplement, repr\'esent\'ee par $R:= (E, d)$ o\`u $d$ est une application de $E^2$ dans $W^k$ telle que $d(x,y)= \overline {d(y,x)}$ pour tout $x,y\in E$. %Ceci est invent\'e pour exprimer que la connaissance de $d(x,y)$ d\'etermine celle de $d(y,x)$.

\medskip

Nous commen\c{c}ons par \'etendre quelques d\'efinitions et notations introduites dans la sous-section \ref{subsec:equiv-graphe-ordonne}. Soit $\mathcal R$ une structure binaire de type $k$ d\'efinie sur un ensemble $E$ et soient $x,y\in E$. Nous disons que la paire $\{x,y\}$ est \emph{sym\'etrique} si $d(x,y)=d(y,x)$ et \emph{asym\'etrique} si $d(x,y)\neq d(y,x)$.
\smallskip

Un \emph{diamant} est une structure binaire de type $k$ d\'efinie sur un ensemble \`a quatre \'el\'ements $\{x,a,b,c\}$ tels que:
        \begin{itemize}
        \item $\{u,v\}$ est asym\'etrique pour tout $\{u,v\}\subset\{x,a,b,c\}$,
        \item $d(x,a)=d(x,b)=d(x,c)=\alpha$ et
        \item $d(a,b)=d(b,c)=d(c,a)\in\{\alpha,\overline{\alpha}\}$.
        \end{itemize}
L'ensemble $\{a,b,c\}$ est la \emph{base} du diamant et $x$ est son sommet.

\noindent Deux diamants de m\^eme base $\{a,b,c\}$ et de sommet $x$ et $y$ respectivement sont dits \emph{oppos\'es} si $d(x,z)=\overline{d(y,z)}$ pour tout $z\in\{a,b,c\}$.

\smallskip
Un \emph{prisme} est une structure binaire de type $k$ d\'efinie sur un ensemble \`a cinq \'el\'ements $\{x,y,a,b,c\}$ tels que:
        \begin{itemize}
        \item $\{x,z\}$ et $\{y,z\}$ sont asym\'etrique pour tout $z\in\{a,b,c\}$,
        \item $d(x,a)=d(a,y)=d(y,b)=d(b,x)=d(x,c)=d(c,y)=\alpha$,
        \item $\{a,c\}$ est asym\'etrique et $d(a,c)\in\{\alpha,\overline{\alpha}\}$.
        \end{itemize}
Les paires $\{x,y\}$, $\{a,b\}$ et $\{b,c\}$ sont les \emph{diagonales} du prisme, $\{x,y\}$ \'etant la diagonale principale.

Pour $k=1$, nous retrouvons les d\'efinitions de diamant et prisme donn\'ees dans les sous-sections \ref{subsec:cas tournois} et \ref{subsec:equiv-graphe-ordonne} respectivement.

\begin{lemma}\label{lem:separbinaire}
Soit $\mathcal R$ une structure binaire de type $k$ d\'efinie sur un ensemble $E$ ayant au moins cinq \'el\'ements. Deux \'el\'ements $x,y$ de $E$ qui sont $\leq2$-\'equivalents ne sont pas \'equivalents si et seulement si ils sont distincts et ou bien
    \begin{enumerate}
    \item[(1)] il existe trois \'el\'ements de $E$ qui forment avec $x$ et $y$ deux diamants oppos\'es, ou  bien
    \item[(2)] $\{x,y\}$ est la diagonale principale d'un prisme dont les autres diagonales sont sym\'etriques.
    \end{enumerate}

\end{lemma}

\begin{proof}
La condition suffisante se v\'erifie facilement. Montrons qu'elle est n\'ecessaire. Soient $x,y\in E$ deux \'el\'ements qui sont $\leq2$-\'equivalents mais non \'equivalents. Soit $F$ une partie finie minimale de $E\setminus\{x,y\}$ telle que les restrictions de $\mathcal R$ \`a $F\cup\{x\}$ et $F\cup\{y\}$ ne soient pas isomorphes. Les \'el\'ements $x$ et $y$ \'etant $\leq 2$-\'equivalents, $F$ a au moins trois \'el\'ements. Soit $C=\{z\in F:d(x,z)= d(y,z)\}$.

    \begin{claim}\label{claim1-preuve lemme equiv}
    $C\neq F$.
    \end{claim}
 \begin{proofclaim}
 Si $C=F$ alors les restrictions de $\mathcal R$ \`a $F\cup\{x\}$ et $F\cup\{y\}$ sont isomorphes, ce qui constitue une contradiction avec le choix de $F$.
 \end{proofclaim}

  \begin{claim}\label{claim2-preuve lemme equiv}
  $C=\varnothing$.
    \end{claim}
 \begin{proofclaim}
 Supposons $C\neq\varnothing$ et soit $h\in C$ et $z\in F\setminus C$. Donc $d(x,h)=d(y,h)$ et $d(x,z)\neq d(y,z)$. les \'el\'ements $x$ et $y$ \'etant $2$-\'equivalents, les restrictions de $\mathcal R$ \`a $\{x,h,z\}$ et $\{y,h,z\}$ sont isomorphes. Le sommet $z$ ne peut-\^etre fix\'e par cet isomorphisme car sinon $x$ est envoy\'e sur $y$ et donc $d(x,z)=d(y,z)$ ce qui n'est pas le cas. De m\^eme $z$ ne peut-\^etre envoy\'e sur $h$ car on airait toujours $d(x,z)=d(y,z)$. Ainsi $z$ est envoy\'e sur $y$ et $x$ sur $z$ ce qui donne $d(h,z)= d(h,y)$. Avec l'egalit\'e $d(h,x)=d(h,y)$ nous d\'eduisons que $d(h,z)=d(h,x)= d(h,y)$ pour tout $z\in F\setminus C$. Il s'ensuit que $F\setminus C)\cup\{x,y\}$ est un intervalle propre de $\mathcal R_{\restriction_{F\cup\{x,y\}}}$. Comme $F\setminus C\varsubsetneq F$ alors les restrictions de $\mathcal R$ \`a $(F\setminus C)\cup\{x\}$ et $(F\setminus C)\cup\{y\}$ sont isomorphes et cet isomorphisme s'\'etend, par l'identit\'e sur $C$, en un isomorphisme de $\mathcal R_{\restriction_{F\cup\{x\}}}$ sur $\mathcal R_{\restriction_{F\cup\{y\}}}$, ce qui constitue une contradiction avec le choix de $F$.
 \end{proofclaim}

 \smallskip

 En conclusion, nous avons $d(x,z)\neq d(y,z)$ pour tout $z\in F$. Mais, en raison de la $1$-\'equivalence de $x$ et $y$, le Lemme \ref{lem:1equivbinaire} donne $d(x,z)= d(z,y)=\overline{d(y,z)}$ pour tout $z\in F$. Donc $d(y,z)\neq d(z,y)$ pour tout $z\in F$. Il s'ensuit que $\{x,z\}$ et $\{y,z\}$ sont forc\'ement asym\'etriques pour tout $z\in F$. Nous avons alors deux cas:

\medskip

\textbf{Cas 1:} $d(x,z)$ est constante pour $z\in F$. Posons alors $d(x,z)=\alpha$ pour $z\in F$. Donc $d(z,y)=\alpha$ pour tout $z\in F$ avec $\alpha\neq\overline{\alpha}$.

    \begin{fact}\label{fact:preuve lemme equiv}
    Si $d(x,u)=d(x,v)=\alpha$ pour $u,v\in F$ alors $\{u,v\}$ est asym\'etrique et $d(u,v)\in\{\alpha,\overline{\alpha}\}$.
    \end{fact}
\begin{proofclaim}
Si $d(x,u)=d(x,v)=\alpha$ alors $d(u,y)=d(v,y)=\alpha$. Les restrictions de $\mathcal R$ \`a $\{x,u,v\}$ et $\{y,u,v\}$ \'etant isomorphes, il n'y a que deux isomorphisme possibles, l'un envoie $x$ sur $v$, $v$ sur $u$ et $u$ sur $y$ et l'autre envoie $x$ sur $u$, $u$ sur $v$ et $v$ sur $y$. Il s'ensuit que $d(u,v)\in\{\alpha,\overline{\alpha}\}$ et donc $\{u,v\}$ est asym\'etrique.
 \end{proofclaim}

 \smallskip

 Il s'ensuit que si la valeur de $d(x,z)$ est constante pour tout $z\in F$ et vaut $\alpha$ alors pour tous $u,v\in F$, $u\neq v$, on a $d(u,v)\in\{\alpha,\overline{\alpha}\}$. Comme $\vert F\vert\geq 3$, nous avons deux situations:
  \begin{enumerate}
  \item ou bien il existe trois \'el\'ements $a,b,c\in F$ tels que $d(a,b)=d(b,c)=d(c,a)\in\{\alpha,\overline{\alpha}\}$ (donc $\{a,b,c\}$ est un $3$-cycle pour $d$). Dans ce cas, $\{a,b,c\}$ forme avec $x$ et $y$ deux $k$-diamants oppos\'es et on est dans le cas \emph{(1)} du Lemme \ref{lem:separbinaire}. En particulier, $F$ a seulement $3$ \'el\'ements.
  \item ou bien il n'existe pas dans $F$ trois \'el\'ements comme ci-dessus et dans ce cas, les \'el\'ements de $F$ forment une cha\^{i}ne pour $d$ et donc les restrictions de $\mathcal R$ \`a $F\cup\{x\}$ et $F\cup\{y\}$ sont isomorphes, ce qui constitue une contradiction.
\end{enumerate}
\medskip

\textbf{Cas 2:} $d(x,z)$ n'est pas constante pour $z\in F$. Donc il existe deux \'el\'ements $a,b\in F$ tels que $d(x,a)=\alpha$ et $d(x,b)=\beta$ avec $\alpha\neq\beta$, $\alpha\neq\overline{\alpha}$ et $\beta\neq\overline{\beta}$. Donc $d(a,y)=\alpha$ et $d(b,y)=\beta$.

\begin{claim}\label{claim3-preuve lemme equiv}
$\{a,b\}$ est sym\'etrique et $\beta=\overline{\alpha}$.
\end{claim}
\begin{proofclaim}
Consid\'erons l'isomorphisme qui envoie $\{x,a,b\}$ sur $\{y,a,b\}$. Cet isomorphisme ne peut pas envoyer $x$ sur $a$ car sinon, si $a$ est envoy\'e sur $y$ on aurait $d(x,b)=d(a,b)=d(y,b)$, donc $\beta=\overline{\beta}$, ce qui est impossible et si $a$ est envoy\'e sur $b$ on aurait $d(x,b)=d(a,y)$ ce qui est impossible puisque $\alpha\neq\beta$. De m\^eme  $x$ ne peut-\^etre envoy\'e sur $b$ car sinon $b$ est envoy\'e sur $y$ et on aurait $d(x,a)=d(b,a)=d(y,a)$ ce qui n'est pas le cas. Donc cet isomorphisme envoie $x$ sur $y$ et permute $a$ et $b$. Nous avons alors:
\begin{equation}
\begin{array}{l}
\alpha=d(x,a)=d(y,b)=\overline{\beta}\\
\beta=d(x,b)=d(y,a)=\overline{\alpha}\\
d(a,b)=d(b,a)%
\end{array}\label{eq:preuveclaim}
\end{equation}
Il s'ensuit que $\{a,b\}$ est sym\'etrique et $\beta=\overline{\alpha}$.
\end{proofclaim}

\smallskip

Comme $\vert F\vert\geq 3$, il existe $c\in F\setminus\{a,b\}$. D'apr\`es ce qui pr\'ec\`ede, $d(x,c)\in\{\alpha,\overline{\alpha}\}$. Supposons, sans perte de g\'en\'eralit\'e que $d(x,c)=\alpha$. Donc $d(c,y)=\alpha$. En utilisant les relations \eqref{eq:preuveclaim}, nous avons
$d(x,a)=d(a,y)=d(y,b)=d(b,x)=d(x,c)=d(c,y)=\alpha$ avec $\{a,b\}$ et $\{b,c\}$ qui sont sym\'etriques. D'apr\`es le Fait \ref{fact:preuve lemme equiv}, $\{a,c\}$ est asym\'etrique et $d(a,c)\in\{\alpha,\overline{\alpha}\}$.
Donc les sommets $\{x,y,a,b,c\}$ forment un $k$-prisme de diagonale principale $\{x,y\}$, les autres diagonales qui sont $\{a,b\}$ et $\{b,c\}$ \'etant sym\'etriques et on est dans le cas \emph{(2)} du Lemme \ref{lem:separbinaire}.
Ceci termine la preuve.
\end{proof}

Il ressort de la preuve du Lemme \ref{lem:separbinaire}.
\begin{corollary}
Soit $\mathcal R:=(E,\rho _{1},\dots,\rho _{k})$ une structure binaire de type $k$. Deux \'el\'ements $x,y$ de $E$ qui sont $\leq2$-\'equivalents ne sont pas \'equivalents si et seulement si il existe $i\in\{1,\dots,k\}$ tel que $x$ et $y$ ne sont pas \'equivalents dans $\mathcal R_i$.
\end{corollary}

    \subsection{Le cas des structures binaires ordonn\'ees}\label{subsection:structures binaires ordo}

Soit $\mathcal R:=(E,\leq ,\rho _{1},\dots,\rho _{k})$ une structure binaire ordonn\'ee de type $k$, o\`u chacune des relations $\rho_i$ est identifi\'ee \`a sa fonction caract\'eristique. Posons $C:=(E,\leq)$. %la relation d'\'equivalence $\simeq_{\mathcal R}$ a les propri\'etes suivantes:

Pour tous $x,y\in E$ et tout $1\leq i\leq k$,  posons $d_i(x,y):=(\rho_i(x,y),\rho_i(y,x))$ et $d(x,y):= (d_i(x,y))_{i= 1,\dots, k}$. Notons par $I_{\leq}(x,y)$ le sous-ensemble des \'el\'ements $z\in E$ qui sont entre $x$ et $y$ modulo $\leq$, c'est \`a dire le plus petit intervalle de $C$ contenant $x,y$.

\vspace{2mm}

Nous avons de mani\`ere \'evidente

\begin{fact}\label{fact:1-equiv-binaire}
Deux sommets $x$ et $y$ de $E$ sont $1$-\'equivalents si et seulement si %pour tout $1\leq i\leq k$
on a
$$\left\{\begin{array}{ll}
d(x,z)=d(z,y)& \text{ si }z\in I_{\leq}(x,y),\\
d(x,z)=d(y,z)& \text{ si }z\notin I_{\leq}(x,y)\end{array}\right.$$
\end{fact}

\begin{lemma}\label{lem;2-equiv=inter}
Soit $\mathcal R$ une structure binaire ordonn\'ee de type $k$ d\'efinie sur $E$ ayant au moins trois \'el\'ements.
%De m\^eme que la $2$-monomorphie implique la $1$-monomorphie pour une structure binaire d'au moins trois \'el\'ements,
\begin{enumerate}
\item[(1)]  Si deux \'el\'ements $x$ et $y$ sont $1$-\'equivalents alors tous les \'el\'ements de $I_{\leq}(x,y)$ sont $0$-\'equivalents.
\item[(2)] Si deux \'el\'ements $x$ et $y$ sont $2$-\'equivalents alors $I_{\leq}(x,y)$ est un intervalle de $\mathcal R$.
\item[(3)] Si $E$ a au moins cinq \'el\'ements et $x$ et $y$ sont $2$-\'equivalents alors les restrictions de $\mathcal R$ \`a $I_{\leq}(x,y)\setminus\{x\}$ et $I_{\leq}(x,y)\setminus\{y\}$ sont $2$-monomorphes donc isomorphes.
\end{enumerate}
\end{lemma}

\begin{proof}
Le (1) du lemme vient directement du fait que si deux restrictions de $\mathcal R$ sont isomorphes, l'isomorphisme est unique.
\vspace{1mm}

Pour le (2) si $I_{\leq}(x,y)=\{x,y\}$ le r\'esultat est \'evident d'apr\`es le Fait \ref{fact:1-equiv-binaire}. Si $I_{\leq}(x,y)\neq\{x,y\}$, soit $z\in I_{\leq}(x,y)\setminus\{x,y\}$ et $a\notin I_{\leq}(x,y)$. Il est facile de voir que $d(a,z)=d(a,x)=d(a,y)$, en effet $x$ et $y$ \'etant $2$-\'equivalents, les restrictions de $\mathcal R$ \`a $\{x,a,z\}$ et $\{y,a,z\}$ sont isomorphes et  l'isomorphisme est unique.
\vspace{1mm}

Pour la preuve de l'assertion (3), d'apr\`es le Corollaire \ref{lem:gottlieb-kantor}, $x$ et $y$ sont $1$-\'equivalents. Le r\'esultat est \'evident si  $I_{\leq}(x,y)$ a deux \'el\'ements et si $I_{\leq}(x,y)$ a trois \'el\'ements, il se d\'eduit de la $1$-\'equivalence de $x$ et $y$. Si $I_{\leq}(x,y)$ a au moins quatre \'el\'ements, il suffit de montrer que pour tous $z,z'\in I_{\leq}(x,y)\setminus\{x,y\}$, les restrictions \`a $\{x,z,z'\}$ et $\{y,z,z'\}$ sont $2$-monomorphes. Les restrictions \`a $\{x,z,z'\}$ et $\{y,z,z'\}$ sont isomorphes et l'isomorphisme est unique. Donc si on suppose $x<z<z'<y$, cet isomorphisme envoie $x$ sur $z$, $z$ sur $z'$ et $z'$ sur $y$. On a donc, $$d(x,z)=d(z,z')=d(z',y)$$ et $$d(x,z')=d(z,y).$$ Nous devons montrer que,
$$d(x,z)=d(z,z')=d(z',y)=d(x,z')=d(z,y).$$ Ceci vient du fait que $E$ a au moins cinq \'el\'ements, donc il existe $a\in E\setminus\{x,y,z,z'\}$. Si $a\notin I_{\leq}(x,y)$ on utilise le fait que les restrictions de $\mathcal R$ \`a $\{a,x,z\}$ et $\{a,y,z\}$ sont isomorphes et si $a\in I_{\leq}(x,y)$ on utilise le fait que les restrictions de $\mathcal R$ \`a $\{a,x,z\}$ et $\{a,y,z\}$  resp. \`a $\{a,x,z'\}$ et $\{a,y,z'\}$ sont isomorphes.
\end{proof}

\vspace{2mm}

Du Lemme \ref{lem;2-equiv=inter} nous d\'eduisons

\begin{theorem}\label{lemme-ordered}
Les relations d'\'equivalence  $\simeq_{\leq 2,{\mathcal R}}$ et $\simeq _{\mathcal R}$ co\"{\i}ncident.
\end{theorem}

\begin{proof}
Supposons que deux \'el\'ements $x$ et $y$ soient $\leq 2$-\'equivalents mais non \'equivalents. Il existe alors un sous-ensemble $F$ d'au moins trois \'el\'ements tel que les restrictions de $\mathcal R$ \`a $F\cup\{x\}$ et $F\cup\{y\}$ ne sont pas isomorphes. Il est clair que $F\subseteq I_{\leq}(x,y)$ car sinon $I_{\leq}(x,y)$ serait un intervalle de $\mathcal R_ {\restriction_{F\cup\{x,y\}}}$ ce qui contredirait le Lemme \ref{intervalle}. D'apr\`es (3) du Lemme \ref{lem;2-equiv=inter}, on d\'eduit que les restrictions de $\mathcal R$ \`a $F\cup\{x\}$ et $F\cup\{y\}$ sont isomorphes, ce qui contredit l'hypoth\`ese. Donc $x$ et $y$ sont \'equivalents.
\end{proof}

\vspace{2mm}

De la Proposition \ref{prop:comp-intervalle mono}, le Lemme \ref{lem:equi-blocmono} et le Lemme \ref{lem:binaire2-monomorphe} nous avons

\begin{lemma}\label{lem:formeclasseequiv} Soit $\mathcal R:=(E,\leq ,\rho _{1},\dots,\rho _{k})$ une structure binaire ordonn\'ee de type $k$.
Toute classe d'\'equivalence d'au moins trois \'el\'ements est un intervalle de $\mathcal R$, donc un intervalle de $\leq$, sur lequel toute relation $\rho_i$ qui est r\'eflexive est soit une cha\^{\i}ne qui  co\"{\i}ncide ou est oppos\'ee \`a $\leq$, soit une  clique r\'eflexive, soit une anticha\^{\i}ne et toute relation $\rho_i$ qui est irr\'eflexive est soit un tournoi acyclique qui co\"{\i}ncide ou est oppos\'e \`a l'ordre strict $<$, soit une clique ou une relation vide. Inversement, un intervalle de $\mathcal R$ n'est pas n\'ecessairement une classe de $1$-\'equivalence, mais tout intervalle qui est contenu dans une classe de $1$-\'equivalence est contenu dans une classe d'\'equivalence.
\end{lemma}

\begin{proof}
Du Lemme \ref{lem:equi-blocmono} nous avons que toute classe d'\'equivalence est une composante monomorphe. De la Proposition \ref{prop:comp-intervalle mono} nous avons que toute composante monomorphe ayant au moins trois \'el\'ements est un intervalle monomorphe maximal, donc un intervalle de $\mathcal R$. La forme des relations $\rho_i$ vient du (3) du Lemme \ref{lem;2-equiv=inter} et du Lemme \ref{lem:binaire2-monomorphe}.

Pour l'inverse, le fait que tout intervalle qui est contenu dans une classe de $1$-\'equivalence est contenu dans une classe d'\'equivalence est d\^u \`a la forme des classes de $1$-\'equivalence donn\'ee par le Fait \ref{fact:classe-monomorphe} (en page \pageref{fact:classe-monomorphe}) et le Corollaire \ref{cor:encha-monomorphe}. Pour voir qu'un intervalle de $\mathcal R$ n'est pas n\'ecessairement une classe de $1$-\'equivalence, voici un contre-exemple. Soit $\mathcal R:=(E,\leq,\rho)$ o\`u $E=\{x,y,z,t\}$ avec\\ $$x<y<z<t \text{ et }\rho=\{(x,y),(x,z),(x,t),(z,t),(t,z)\}.$$
   $A=\{y,z,t\}$ est un intervalle de $\mathcal R$ mais n'est pas une classe de $1$-\'equivalence car $y$ et $z$ ne sont pas $\{t\}$-\'equivalents. Les classes de $1$-\'equivalence sont $\{x\}$, $\{y\}$ et $\{z,t\}$. Ces classes sont \'egalement des classes d'\'equivalence.
\end{proof}

\begin{lemma}\label{lem:adjacent}
Si deux classes d'\'equivalence sont telles que leur union est un intervalle de $\leq$  alors elles ne sont pas contenues dans une m\^eme classe de $1$-\'equivalence.
\end{lemma}

\begin{proof}
Ceci est d\^u au fait que si deux \'el\'ements $x$ et $y$ sont tels que $I_{\leq}(x,y)=\{x,y\}$ alors $x$ et $y$ sont $1$-\'equivalents si et seulement si ils sont \'equivalents car dans ce cas $\{x,y\}$ est un intervalle de $\mathcal R$ d'apr\`es le Fait \ref{fact:1-equiv-binaire} (en page \pageref{fact:1-equiv-binaire}).
\end{proof}

\vspace{2mm}

Comme cons\'equence du Fait \ref{fact:1-equiv-binaire}, du Th\'eor\`eme \ref{lemme-ordered} et du Lemme \ref{lem;2-equiv=inter} nous avons:

\begin{lemma}\label{lem:separation-binaire}  Soit $\mathcal R:=(E,\leq ,\rho _{1},\dots,\rho _{k})$ une structure binaire ordonn\'ee de type $k$.
Deux sommets  $x,y$ de $E$ ne sont pas \'equivalents si et seulement si ils sont distincts et ou bien
\begin{enumerate}
\item $x,y$ ne sont pas $1$-\'equivalent, ou bien
\item $x,y$ sont $1$-\'equivalents et une seule parmi les deux conditions suivantes a lieu;
    \begin{enumerate}
    \item Il existe $z,z'\in I_{\leq}(x,y)\setminus\{x,y\}$ tels que ou bien  $d_i(x,z)\neq d_i(x,z')$ ou bien $d_i(x,z)\neq d_i(z,z')$ pour un $1\leq i\leq k$.
    \item Il existe $z\in I_{\leq}(x,y)\setminus\{x,y\}$ et $z'\notin I_{\leq}(x,y)$ tels que $d_i(z',z)\neq d_i(z',x)=d_i(z',y)$ pour un $1\leq i\leq k$.
    \end{enumerate}
\end{enumerate}
\end{lemma}

\medskip

Le Lemme \ref{lem:separation-binaire} stipule que si deux sommets $x, y$ sont $1$-\'equivalents mais non \'equivalents, alors les classes d'\'equivalence $C_x$ et $C_y$ de $x$ et $y$ respectivement  sont s\'epar\'ees par au moins une troisi\`eme classe $C_z$ ($C_z\subseteq I_{\leq}(x,y)$) qui n'appartient pas \`a la m\^eme classe de $1$-\'equivalence que $C_x$ et $C_y$.

            \section{Structures binaires ayant une infinit\'e de classes d'\'equivalence. Application au profil}\label{section:graphs-ordered}
Dans cette section, nous montrons que la Conjecture \ref{conjecture2} (page \pageref{conjecture2}) est vraie dans le cas des graphes non dirig\'es et des graphes dirig\'es ordonn\'es et que la r\'eciproque du Th\'eor\`eme \ref{thm: polynomial-interval} est vraie dans le cas des  structures binaires ordonn\'ees.

\vspace{2mm}

%Nous avons les r\'esultats suivants:
Dans le cas des graphes non dirig\'es (sans boucle) nous avons

\begin{theorem}\label{prop: graphs}
Il existe un ensemble $\mathfrak A$ form\'e de dix graphes infinis tel que tout graphe qui ne se d\'ecompose pas en une somme lexicographique finie de cliques et d'ind\'ependants abrite une copie de l'un de ces graphes.% There are ten infinite graphs such that a graph does not decompose into a finite lexicographic sum of cliques and independent sets iff it contains a copy of one of these ten graphs.
\end{theorem}

%\vspace{1mm}

Dans le cas des graphes dirig\'es ordonn\'es, nous avons le r\'esultat suivant si nous remplaçons  $\mathscr S_{\mu}$, dans la Conjecture \ref{conjecture2},  par la classe $\mathscr D$ form\'ee de graphes dirig\'es ordonn\'es qui n'ont pas de d\'ecomposition monomorphe finie:

 \begin{theorem}\label{thm:graph-ordonne}
 Il existe un ensemble $\mathfrak A$ form\'e de mille deux cent quarante six graphes (dirig\'es) ordonn\'es tel que tout \'el\'ement de  $\mathscr D$ abrite un \'el\'ement de $\mathfrak A$.
\end{theorem}

%\vspace{1mm}

Nous avons montr\'e, dans le Th\'eor\`eme \ref{thm: polynomial-interval} que toute structure ordonn\'ee qui poss\`ede une d\'ecomposition monomorphe finie et ainsi toute classe form\'ee de telles structures, a un profil polyn\^omial. De mani\`ere g\'en\'erale, une structure relationnelle infinie ayant un profil born\'e par un polyn\^ome n'a pas forc\'ement une d\'ecomposition monomorphe finie, sauf si la structure est binaire ordonn\'ee (des exemples de graphes sans d\'ecomposition monomorphe finie ayant un profil polyn\^omial sont donn\'es en page \pageref{subsubsec:les dix graphes}). Ceci est une cons\'equence du r\'esultat suivant dont la preuve s'appuit, en partie, sur le Th\'eor\`eme \ref{thm:graph-ordonne}:

\begin{theorem}\label{theo:dichotomie}
 Si une structure binaire ordonn\'ee $\mathcal R$, de type $k$, a une d\'ecomposition monomorphe finie alors le profil de son \^age $Age(\mathcal R)$ est un polyn\^ome,   autrement  ce profil est au moins exponentiel.
\end{theorem}

Ce r\'esultat de dichotomie s'\'etend aux classes h\'er\'editaires de structures binaires ordonn\'ees.

 \begin{theorem}\label{theo:dichotomie-classe} Soit $\mathscr C$ une classe h\'er\'editaire de structures binaires ordonn\'ees finies de type $k$. Alors, ou bien il existe un entier $\ell$ tel que tout membre de $\mathscr  C$ poss\`ede une d\'ecomposition monomorphe en au plus $\ell+1$ blocs, auquel cas $\mathscr C$ est une union finie d'\^ages de structures binaires ordonn\'ees, chacune ayant une d\'ecomposition monomorphe en au plus $\ell+1$ blocs et le profil de $\mathscr C$ est un polyn\^ome, ou bien le profil de $\mathscr C$ est au moins exponentiel.
\end{theorem}

\begin{proof}
S'il existe un entier $\ell$ tel que tout membre de $\mathscr  C$ poss\`ede une d\'ecomposition monomorphe en au plus $\ell+1$ blocs, alors d'apr\`es la Proposition \ref{prop:hwqo},  $\mathscr  C$ est h\'er\'editairement belordonn\'ee, c'est donc une union finie d'\^ages de structures binaires ordonn\'ees, chacune ayant une d\'ecomposition monomorphe en au plus $\ell+1$ blocs d'apr\`es le Th\'eor\`eme \ref{compactness} et le profil de $\mathscr C$ est un polyn\^ome d'apr\`es le Th\'eor\`eme \ref{thm: polynomial-interval}.

S'il n'existe aucun entier $\ell$ tel que tout membre de $\mathscr  C$ poss\`ede une d\'ecomposition monomorphe en au plus $\ell+1$ blocs, alors d'apr\`es le Lemme \ref{lem:reduction}, $\mathscr  C$ contient une classe h\'er\'editaire $\mathscr A$, ayant la m\^eme propri\'et\'e, qui est minimale pour l'inclusion. Cette classe est l'\^age d'une structure relationnelle binaire ordonn\'ee $\mathcal R$ qui ne poss\`ede pas de d\'ecomposition monomorphe finie (d'apr\`es le Th\'eor\`eme \ref{compactness}). Donc d'apr\`es le Th\'eor\`eme \ref{theo:dichotomie}, le profil de $\mathcal R$, et donc de $\mathscr A$, est au moins exponentiel. Il s'ensuit que le profil de $\mathscr C$ est au moins exponentiel.
\end{proof}

\vspace{2mm}

La preuve du Th\'eor\`eme \ref{theo:dichotomie} repose sur la description des structures de $\mathfrak A$ du Th\'eor\`eme \ref{thm:graph-ordonne}. Pour les preuves des Th\'eor\`emes \ref{prop: graphs} et \ref{thm:graph-ordonne}, nous utilisons le th\'eor\`eme de Ramsey\index{Ramsey} sous la forme du Th\'eor\`eme  \ref{thm:ramsey-invariant} pour trouver les \'el\'ements de l'ensemble $\mathfrak A$ dans chacun de ces th\'eor\`emes. Ces membres, comme il a \'et\'e montr\'e dans la Proposition \ref{mainthm},  sont presque multi-enchaînables et sont d\'efinis sur $F\cup(L\times K)$ avec $L:=\mathbb N$, $\vert K\vert=2$ et $F=\varnothing$  dans le cas des graphes non dirig\'es et  $\vert F\vert\leq 1$) dans le cas des graphes dirig\'es ordonn\'es. Avant de faire une description de ces graphes (resp. graphes dirig\'es ordonn\'es), donnons d'abord la mani\`ere de les construire.

        \subsection{Le cas des graphes. Preuve du Th\'eor\`eme \ref{prop: graphs}}\label{subsec:graph}

Les graphes consid\'er\'es dans ce paragraphe sont non dirig\'es.

\vspace{2mm}

Soit  $G:= (V, E)$ un graphe non dirig\'e.
Nous rappelons que si un graphe $G$ est une somme lexicographique $\underset{i\in H}\sum G_i$ d'une famille de graphes $G_i$, index\'ee par un  graphe $H$, alors, si les ensembles $V(G_i)$ sont deux \`a deux disjoints, ils forment une partition de $V$ en intervalles. Inversement, si l'ensemble $V$ est partitionn\'e en intervalles, alors $G$ est la somme lexicographique des  graphes induits par les blocs de cette partition.

\medskip

 La preuve du Th\'eor\`eme \ref{prop: graphs} se fait comme suit. Soit $G$ un graphe qui ne se d\'ecompose pas en une somme lexicographique finie $\underset{i\in H}\sum G_i$ d'une famille de graphes  $G_i$,  chacun \'etant une clique ou un ind\'ependant, index\'ee par un graphe fini $H$. Alors, d'apr\`es le Th\'eor\`eme \ref{coro:graphe},  $G$ poss\`ede un nombre infini de classes d'\'equivalence. Donc, il existe une application injective  $f:\NN\longrightarrow V(G)$ telle que les images de deux  \'el\'ements distincts sont non \'equivalentes. Nous pouvons alors d\'efinir une application $g:[\NN]^2\longrightarrow V(G)$ telle que $g(n,m)$ t\'emoigne du fait que $f(n)$ et $f(m)$ sont non \'equivalents pour tous $n<m$. C'est \`a dire,
 \begin{equation}
 \{f(n),g(n,m)\}\in E(G)\Leftrightarrow\{f(m),g(n,m)\}\notin E(G).\label{eq:equivalence-graphe}
 \end{equation}

Soit $\omega:=(\NN,\leq)$, $\Phi:=\{f,g\}$ et $\mathfrak{L}:=\left\langle \omega,G,\Phi \right\rangle $. Le Th\'eor\`eme \ref{thm:ramsey-invariant} permet de trouver un sous-ensemble infini $X$ de $\NN$ tel que $\mathfrak L_{\restriction_X}$ est invariante.

\vspace{1mm}

En indexant les \'el\'ements de l'ensemble $X$ par des entiers, nous pouvons supposer que $X=\NN$ et donc $\mathfrak{L}$ est invariante.

\begin{claim}\label{graphclaim}
\begin{enumerate}
\item $\{f(n),f(m)\}\in E(G)\Leftrightarrow\{f(n'),f(m')\}\in E(G),~\forall n<m,~n'<m'$.
\item $\{g(n,m),g(n',m')\}\in E(G)\Leftrightarrow\{g(k,l),g(k',l')\}\in E(G),~\forall n<m<n'<m',~k<l<k'<l'$.
\item $\{f(n),g(n,m)\}\in E(G)\Leftrightarrow\{f(k),g(k,l)\}\in E(G),~\forall n<m,~k<l$.
\item $\{f(m),g(n,m)\}\in E(G)\Leftrightarrow\{f(l),g(k,l)\}\in E(G),~\forall n<m,~k<l$.
\item $\{f(k),g(n,m)\}\in E(G)\Leftrightarrow\{f(l),g(p,q)\}\in E(G),~\forall n<m<k,~p<q<l$.
\item $\{f(n),g(m,k)\}\in E(G)\Leftrightarrow\{f(p),g(q,l)\}\in E(G),~\forall n<m<k,~p<q<l$.
\item Si  $\{f(n),g(n,m)\}$ et  $\{f(k),g(n,m)\}$ sont tous les deux des ar\^etes ou des non ar\^etes pour trois entiers $n<m<k$, il en est de m\^eme pour $\{f(n'),g(n',m')\}$ et  $\{f(k'),g(n',m')\}$  pour tous $n'<m'<k'$.
\item $g(n,m)$ et $f(k)$ sont distincts pour tous entiers  $n<m$ et $k$ .
\item $g(n,m)$ et $g(n',m')$ sont distincts pour tous $n<m<n'<m'$.
\end{enumerate}
\end{claim}

\begin{proofclaim}
%\textbf{Proof of Claim \ref{graphclaim}}
%$(1)$, $(2)$, $(3)$, $(4)$ and $(5)$ follow from invariance.
Les six premieres assertions d\'ecoulent de l'invariance de $\mathfrak{L}$. L'assertion \emph{7} est la m\^eme que l'assertion \emph{3} de l'Affirmation \ref{basicclaim} (dans la preuve de la Proposition \ref{mainthm} en page \pageref{basicclaim}), elle d\'ecoule des assertions \emph{3} et \emph{5} ci-dessus. Pour montrer l'assertion \emph{8}, supposons, sans perte de g\'en\'eralit\'e, qu'il existe des entiers  $n<m<k$ tels que $g(n,m)=f(k)$. Par construction des fonctions $f$ et $g$ nous avons $$\{f(n),g(n,m)\}\in E(G)\Leftrightarrow\{f(m),g(n,m)\}\notin E(G).$$ D'apr\`es  l'assertion \emph{1}, $\{f(n),f(k)\}$ et $\{f(m),f(k)\}$ sont tous deux des ar\^etes ou des non ar\^etes, en remplaçant $f(k)$ par $g(n,m)$ nous arrivons \`a une contradiction avec l'\'equivalence \eqref{eq:equivalence-graphe}.
      Pour l'assertion \emph{9}, supposons qu'il existe des entiers $n<m<n'<m'$ tels que  $g(n,m)=g(n',m')$. D'apr\`es l'assertion \emph{5},  $\{f(n'),g(n,m)\}$ et  $\{f(m'),g(n,m)\}$ sont tous deux des ar\^etes ou des non ar\^etes. En remplaçant $g(n,m)$ par $g(n',m')$, nous obtenons  $\{f(n'),g(n',m')\}$ et $\{f(m'),g(n',m')\}$ sont tous deux des ar\^etes ou des non ar\^etes. Ce qui contredit l'\'equivalence \eqref{eq:equivalence-graphe}.
%\hfill       $\Box$
\end{proofclaim}

\vspace{2mm}

 Soit $F:\NN\times\{0,1\}\longrightarrow V(G)$. Posons $F(n,1):=g(2n,2n+1)$. A partir de l'assertion \emph{7} de l'Affirmation \ref{graphclaim} nous avons deux cas. Si $\{f(n),g(n,m)\}$ et  $\{f(k),g(n,m)\}$ sont tous deux des ar\^etes ou des non ar\^etes pour des entiers  $n<m<k$, alors nous posons  $F(n,0):=f(2n+1)$. Autrement, nous posons  $F(n,0):=f(2n)$.

\vspace{1mm}

Soit $G'$ le graph ayant pour sommets l'ensemble  $V(G'):=\NN\times\{0,1\}$ tel que
 $$\{x,y\}\in E(G')\Leftrightarrow\{F(x),F(y)\}\in E(G)$$
pour tous $x,y$ de $\NN\times\{0,1\}$.\\ Par construction de $G'$ nous avons:
\begin{equation}\label{eq1}
\quad\{(n,0),(n,1)\}\in E(G')\Leftrightarrow\{(m,0),(n,1)\}\notin E(G')
\end{equation}
pour tous entiers $n<m$, alors $(n,0)$ et $(m,0)$ sont non \'equivalents pour tous $n<m$ et donc  $G'$ poss\`ede un nombre infini de classes  d'\'equivalence.
%Since $\mathfrak L$ is invariant, we have:
 Par l'invariance de $\mathfrak L$ le graphe $G'$ est une presque multichaîne et donc $G_{\restriction_{Im(F)}}$ est presque multi-enchaînable. %following claim

\vspace{2mm}

 Ainsi, pour construire $G'$ il suffit de d\'ecider des relations reliant les quatre sommets $(0,0)$, $(0,1)$, $(1,0)$ et $(1,1)$, les relations restantes  seront d\'eduites \`a partir des isomorphismes locaux de  $C:=(\NN,\leq)$ et de l'Affirmation \ref{graphclaim}.\\ Par exemple si $\{(0,0),(0,1)\}\in E(G')$, alors $\{(0,1),(1,0)\}\notin E(G')$ et l'invariance donne $\{(1,0),(1,1)\}\in E(G')$ nous avons alors \`a d\'ecider sur les relations entre les sommets des ensembles $\{(0,0),(1,0)\}$, $\{(0,0),(1,1)\}$  et $\{(0,1),(1,1)\}$, ceci donne huit ($2^3$) graphes, mais deux des graphes obtenus sont isomorphes, donc nous avons seulement sept graphes. Si  $\{(0,0),(0,1)\}\notin E(G')$, alors $\{(0,1),(1,0)\}\in E(G')$ et par invariance $\{(1,0),(1,1)\}\notin E(G')$, nous avons huit autres cas mais cinq ont d\'ej\`a \'et\'e obtenus dans le cas pr\'ec\'edent, ainsi nous avons seulement trois nouveaux graphes. Le total donne les dix graphes annonc\'es.

\vspace{2mm}

 Si $\{(0,0),(1,0)\}\in E(G')$ alors $\{(n,0),(m,0)\}\in E(G')$ pour tous $n<m$ et donc $\{(n,0), ~n\in\NN\}$ est une clique, elle est no\'ee $K_{\mathbb N}$ et si  $\{(0,0),(1,0)\}\notin E(G')$ alors $\{(n,0),(m,0)\}\notin E(G')$ pour tous $n<m$ et donc $\{(n,0), ~n\in\NN\}$ est un stable et est not\'e $I_{\mathbb N}$.

\subsubsection{Description des dix graphes}

Notons par  $G_i, 1\leq i\leq 10$, les dix graphes. Ils ont le m\^eme ensemble de sommets $V(G_i):=\mathbb N\times\{0,1\}$. Posons  $A:=\mathbb N\times\{0\}$ et $B:=\mathbb N\times\{1\}$.

\vspace{2mm}

 Pour $~i:=1,2,3$, les sous-ensembles $A$ et $B$ sont des stables et une paire $\{(n,0),(m,1)\}$ pour $n,m\in\mathbb N$ est une ar\^ete de $G_1$ si $n=m$, une ar\^ete de $G_2$ si $n\leq m$ et une ar\^ete de $G_3$ si $n\neq m$.\\ Ainsi,  $G_1$ est la somme directe d'une infinit\'e de copies de $K_2$ (le graphe complet \`a deux sommets) et $G_2$ est le biparti demi-complet de Schmerl et Trotter \cite{S-T} (les graphes $G_2$ et $G_3$ ont \'et\'e \'etudi\'e dans le chapitre \ref{sec:age de graphe}).

\vspace{2mm}

 Pour $4\leq i\leq 7$, l'un des sous-ensembles $A$, $B$ est une clique et l'autre est un stable. Ces graphes ont pour ensembles d'ar\^etes
 $$E(G_4)=E(G_1)\cup\{\{(n,0),(m,0)\}, n\neq m\in\mathbb N\},$$
  $$E(G_5)=E(G_2)\cup\{\{(n,0),(m,0)\}, n\neq m\in\mathbb N\},$$
  $$E(G_6)=E(G_2)\cup\{\{(n,1),(m,1)\}, n\neq m\in\mathbb N\}$$
  et $$E(G_7)=E(G_3)\cup\{\{(n,0),(m,0)\}, n\neq m\in\mathbb N\}.$$
   Les graphes $G_i$  pour $i=8,9,10$, sont tels que les sous-ensembles $A$ et $B$ sont tous les deux des cliques avec $E(G_8)\cap E(G_1)=E(G_1)$, $E(G_9)\cap E(G_2)=E(G_2)$ et $E(G_{10})\cap E(G_3)=E(G_3)$. Les graphes $G_7$, $G_8$ et $G_{10}$ sont les compl\'ementaires des graphes $G_4$, $G_3$, et  $G_1$ respectivement. %le graphe $G_7$ le compl\'ementaire de $G_4$ et le graphe $G_{10}$ le compl\'ementaire de $G_1$,
   Chacun des graphes $G_5$ et $G_6$ est \'equimorphe, mais non isomorphe, \`a son graphe compl\'ementaire. Toute restriction finie de $G_5$ (resp. $G_6$) s'abrite dans $G_6$ (resp. $G_5$). Donc, les graphes  $G_5$ et $G_6$ ont le m\^eme \^age.
   Chacun des graphes $G_2$ et $G_9$ s'abrite dans le graphe compl\'ementaire de l'autre.
   Ces dix  graphes sont donn\'es dans la \figurename~\ref{tab:graphes}.

\begin{figure}
\begin{center}
\small
\begin{tabular}[l]{ccc}
\input{grapheG1}&&\input{grapheG2}\\
    $\begin{array}{c}
    G_1:=(\mathbb N\times\{0,1\},E_1)\\
    E_1:=\{\{(n,0),(n,1)\},~n\in\mathbb N\}%
    \end{array}$&&$\begin{array}{c}
                                G_2:=(\mathbb N\times\{0,1\},E_2)\\
                                    E_2:=\{\{(n,0),(m,1)\},~n\leq m\in\mathbb N\}%
                                    \end{array}$\\
                                    &&\\
\input{grapheG3}&&\input{grapheG4}\\
$\begin{array}{c}
G_3:=(\mathbb N\times\{0,1\},E_3)\\
E_3:=\{\{(n,0),(m,1)\},~n\neq m\in\mathbb N\}%
\end{array}$&&$\begin{array}{c}
G_4:=(\mathbb N\times\{0,1\},E_4)\\
E_4:=E_1\cup \{\{(n,0),(m,0)\},~n\neq m\in\mathbb N\}%
\end{array}$\\
                 &&\\
\input{grapheG5}&&\input{grapheG6}\\
$\begin{array}{c}
G_5:=(\mathbb N\times\{0,1\},E_5)\\
E_5:=E_2\cup \{\{(n,0),(m,0)\},~n\neq m\in\mathbb N\}%
\end{array}$&&$\begin{array}{c}
G_6:=(\mathbb N\times\{0,1\},E_6)\\
E_6:=E_2\cup\{\{(n,1),(m,1)\},~n\neq m\in\mathbb N\}%
\end{array}$\\
                        &&\\
\input{grapheG7}&&\input{grapheG8}\\
$\begin{array}{c}
G_7:=(\mathbb N\times\{0,1\},E_7)\\
E_7:=E_3\cup \{\{(n,0),(m,0)\},~n\neq m\in\mathbb N\}%
\end{array}$&&$\begin{array}{c}
G_8:=(\mathbb N\times\{0,1\},E_8)\\
E_8:=E_4\cup  \{\{(n,0),(m,0)\},~n\neq m\in\mathbb N\}%
\end{array}$\\
                            &&\\
\input{grapheG9}&&\input{grapheG10}\\
$\begin{array}{c}
G_9:=(\mathbb N\times\{0,1\},E_9)\\
E_9:=E_5\cup \{\{(n,1),(m,1)\},~n\neq m\in\mathbb N\}\\
\text{ou }E_9:=E_6\cup \{\{(n,0),(m,0)\},~n\neq m\in\mathbb N\}%
\end{array}$&&$\begin{array}{c}
G_{10}:=(\mathbb N\times\{0,1\},E_{10})\\
E_{10}:=E_7\cup  \{\{(n,1),(m,1)\},~n\neq m\in\mathbb N\}%
\end{array}$\\
\end{tabular}
\caption{Les graphes minimaux n'admettant pas de d\'ecomposition monomorphe finie.}
\label{tab:graphes}
\end{center}
\end{figure}

\subsubsection{Profils des dix graphes}\label{subsubsec:les dix graphes}

Si $G$ est un graphe, nous d\'esignons par $\varphi_G$ son profil. Clairement, le profil de $G$ est \'egal \`a celui de son compl\'ementaire $\overline{G}$ et si $G$ abrite un graphe
 $H$ alors $\varphi_H\leq \varphi_G$. Donc, d'apr\`es la description ci-dessus, les graphes $G_1$, $G_2$, $G_3$, $G_4$ et $G_5$ ont les m\^emes  profils que $G_{10}$, $G_9$, $G_8$, $G_7$ et $G_6$ respectivement. Comme nous le verrons ci-dessous, les profils $\varphi_{G_1}$, $\varphi_{G_3}$ et $\varphi_{G_4}$ ont des croissances polynomiales tandis que  $\varphi_{G_2}$ et $\varphi_{G_5}$ ont des croissances exponentielles. Ceci montre en particulier qu'une structure ayant un profil polynomialement born\'e n'a pas forc\'ement une d\'ecomposition monomorphe finie.
\vspace{2mm}

 \textbf{Profil de $G_1$:} nous avons $\varphi_{G_1}(n)=\lfloor\frac{n}{2}\rfloor +1,\forall n\in\mathbb N$, sa croissance est polynomiale, $\varphi_{G_1}\simeq\dfrac{n}{2}$. En effet, tout sous-graphe de $G_1$ d\'efini sur $n$ sommets est isomorphe, pour deux entiers $p, q$, \`a $pK_2\oplus q$,  qui est la somme directe de $p$ copies de $K_2$, la clique \`a deux sommets, et d'un ind\'ependant d'ordre $q$, tels que  $p\leq \frac{n}{2}$ et $q=n-2p$. Donc, nous pouvons repr\'esenter tout sous-graphe d'ordre  $n$ par un couple $(n,p)$ d'entiers avec $p\leq \frac{n}{2}$ et inversement. Deux sous-graphes repr\'esent\'es par $(n,p)$ et $(n',p')$ sont isomorphes si et seulement si $n=n'$ et $p=p'$. Sa fonction g\'en\'eratrice est rationnelle (voir dans la sous-section \ref{subsec:profil} en page \pageref{par:exemple}), elle est donn\'ee par:
$$F_{G_1}(x)=\dfrac{1}{(1-x)(1-x^2)}.$$  %Then $\varphi_{G_1}(n)=\lfloor\frac{n}{2}\rfloor +1,\forall n\in\mathbb N$ and thus is bounded above by a polynomial.
%$$\varphi'_3(n)=\vert\{(n,p)/ p\leq \dfrac{n}{2}, p\in\mathbb N\}\vert=\lfloor\frac{n}{2}\rfloor +1,\forall n\in\mathbb N.$$
\vspace{2mm}

\textbf{Profil de $G_2$:} Les premi\`eres valeurs du profil pour $n=0, 1, 2, 3, 4, 5, 6, 7, 8, 9$ sont $1$, $1$, $2$, $3$, $6$, $10$, $20$, $36$, $72$, $136$. Il v\'erifie la r\'ecurrence suivante: $\varphi_{G_2}(n)=\varphi_{G_2}(n-1)+2^{n-3}$ pour $n\geq 3$ et $n$  impair; $\varphi_{G_2}(n)=\varphi_{G_2}(n-1)+2^{n-3}+2^{\frac{n-4}{2}}$ pour $n\geq 4$ et $n$ pair. Ainsi $\varphi_{G_2}(n)\geq 2^{n-3}$ pour $n\geq 3$ (voir paragraphe \ref{subsec:graphe $G_1$} du chapitre \ref{sec:age de graphe}). Sa fonction g\'en\'eratrice est rationnelle, elle est donn\'ee par:
$$F_{G_2}(x)=\dfrac{1-x-2x^2+x^3}{(1-2x)(1-2x^2)}.$$  %and then is bounded bellow by an exponential.
\vspace{2mm}

\textbf{Profil de $G_3$:} Les premi\`eres valeurs du profil pour $n=0, 1, 2, 3, 4, 5, 6, 7$ sont $1,1$, $2,3$, $6,6$, $10,10$. Il est donn\'e par $\varphi_{G_3}(n)=\underset{k=0}{\overset{\lfloor n/2\rfloor}{\sum}}(k+1)$ pour tout $n\neq 2$, sa croissance est polynomiale  de puissance $2$, $\varphi_0(n)\simeq\dfrac{n^2}{8}$. La fonction g\'en\'eratrice est la fonction rationnelle (voir paragraphe \ref{subsec:profilG0} du chapitre \ref{sec:age de graphe}): $$F_{G_3}(x)=\dfrac{1-x^2+x^3+2x^4-2x^5-x^6+x^7}{(1-x)(1-x^2)^2}.$$ %Thus, it is bounded above by a polynomial.

\textbf{Profil de $G_4$:} Les premi\`eres valeurs du profil pour $n=0, 1, 2, 3, 4, 5, 6, 7, 8, 9$ sont $1$, $1$, $2$, $4$, $7$, $10$, $14$, $18$, $23$, $28$. Il v\'erifie la relation suivante: $\varphi_{G_4}(n)=\dfrac{1}{4}(n-1)(n+5)$ pour $n\geq 3$ et $n$ impair; $\varphi_{G_4}(n)=\frac{1}{4}n(n+4)-1$ pour $n\geq 4$ et $n$ pair et nous avons $\varphi_4\simeq\dfrac{n^2}{4}$. Sa fonction g\'en\'eratrice est la fonction rationnelle (voir paragraphe \ref{subsec:profilG4} du chapitre \ref{sec:age de graphe}) $$F_{G_4}(x)=\dfrac{1-x+2x^3-x^5}{(1-x)^3(1+x)}.$$

 \textbf{Profil de $G_5$:} $\varphi_{G_5}(n)=2^{n-1}$ pour $n\geq 1$. Les premi\`eres valeurs du profil pour $n=0, 1, 2, 3, 4, 5, 6$ sont $1,1, 2, 4,8, 16, 32$ (voir paragraphe \ref{subsec:graphe $G_5$} du chapitre \ref{sec:age de graphe}). Sa fonction g\'en\'eratrice est la fonction rationnelle, elle est donn\'ee par: $$F_{G_5}(x)=\dfrac{1+x-2x^2}{(1-2x)}.$$

        \subsection{Graphes dirig\'es ordonn\'es. Preuve  du Th\'eor\`eme \ref{thm:graph-ordonne}}\label{subsec:graphesordonnes}

Soit $\mathcal G:=(V,\leq,\rho)$ un graphe dirig\'e ordonn\'e o\`u %$V$ is the set of vertices, $\leq$ a linear order on $V$ and $\rho$ a binary relation on $V$ which is
$\rho$ est identifi\'ee \`a sa fonction caract\'eristique. %Pour \'eviter des complications inutiles, nous supposons que $\rho$ est soit r\'eflexive soit irr\'eflexive. %that is a map from $V\times V$ in $\{0,1\}$. We suppose that $\rho$ is reflexive.
Le sous-ensemble $E$ de $V^2$ tel que $(x,y)\in E$ si et seulement si $\rho(x,y)=1$ est l'ensemble des arcs de $\mathcal G$ et $G:=(V,E)$ est le graphe dirig\'e associ\'e \`a $\mathcal G$.   Pour tous $x,y\in V$,  $d(x,y):=(\rho(x,y),\rho(y,x))$ et $I_{\leq}(x,y)$ sont % $d:V^2\longrightarrow \{0,1\}^2$
d\'efinis comme dans la section \ref{subsection:structures binaires ordo} ci-dessus.

Un graphes dirig\'es ordonn\'es \'etant une structure binaire ordonn\'ee (de type 1), tous les r\'esultats de la section \ref{subsection:structures binaires ordo} sont vrais dans ce cas.

\vspace{2mm}

\begin{lemma}\label{lem:separation2}
Si un graphe dirig\'e  ordonn\'e $\mathcal G:=(V,\leq,\rho)$ poss\`ede une infinit\'e de classes d'\'equivalence alors il v\'erifie une seule des deux assertions suivantes:
\begin{enumerate}
 \item Il existe un sous-ensemble infini $A\subseteq V$ tel que deux \'el\'ements  distincts de $A$ sont $0$-\'equivalents mais non $1$-\'equivalents.
 \item Il existe deux sous-ensembles infinis disjoints $A_1$ et $A_2$ de $V$ tels que deux  \'el\'ements distincts de $A_i$, $i=1,2$, sont $1$-\'equivalents mais non \'equivalents et pour tous $x,y\in A_i$,  la trace de l'intervalle  $I_{\leq}(x,y)$ sur $A_j$, $j\neq i\in\{1,2\}$, est non vide.
  \end{enumerate}
\end{lemma}

\begin{proof}
Si $\mathcal G$ poss\`ede une infinit\'e de classes d'\'equivalence alors nous avons deux cas:

    \textbf{Cas 1:} $\mathcal G$ poss\`ede une infinit\'e de classes de $1$-\'equivalence. Ces classes sont r\'eparties suivant deux classes de $0$-\'equivalence. D'apr\`es le th\'eor\`eme de Ramsey\index{Ramsey} (Th\'er\`eme \ref{thm:ramsey}), il existe une classe de $0$-\'equivalence qui contient une infinit\'e de classes de $1$-\'equivalence. D\'esignons cette classe par $Cl$. Comme toute classe d'\'equivalence est contenue dans une classe de $1$-\'equivalence, nous pouvons piquer un \'el\'ement de chaque classe de $1$-\'equivalence appartenant \`a $Cl$ pour former un sous-ensemble $A$ de $V$. L'ensemble $A$ v\'erifie l'assertion 1 du lemme.

 \vspace{1mm}

    \textbf{Cas 2:} $\mathcal G$ poss\`ede un nombre fini de classes de $1$-\'equivalence. Supposons un ensemble $F$ form\'e d'une infinit\'e d'\'el\'ements deux \`a deux non \'equivalents. Comme $\mathcal G$ poss\`ede un nombre fini de classes de $1$-\'equivalence, on peut supposer que tous les \'el\'ements de $F$ sont $1$-\'equivalents (en utilisant le th\'eor\`eme de Ramsey), donc $F\subseteq C$ o\`u $C$ est une classe de $1$-\'equivalence. Soient $a, b\in F$. Les \'el\'ements de l'intervalle $I_{\leq}(a,b)$ ne peuvent pas \^etre tous $1$-\'equivalents car sinon $I_{\leq}(a,b)$ serait un intervalle monomorphe d'apr\`es le Lemme \ref{lem:formeclasseequiv}. Donc, il existe $c\in I_{\leq}(a,b)$ qui est dans une classe de $1$-\'equivalence $C'\neq C$. Nous pouvons extraire de $F$ une suite monotone $(a_i)_{i\geq 0}$. D'apr\`es ce qui pr\'ec\`ede, pour tout $i\geq 0$, il existe $c_i\in I_{\leq}(a_i,a_{i+1})$ avec $c_i$ appartenant \`a une classe de $1$-\'equivalence diff\'erente de $C$. Nous pouvons donc trouver une sous-suite infinie $(c_i')_{i\geq 0}$ de $(c_i)_{i\geq 0}$ dont les \'el\'ements sont tous dans une m\^eme classe de $1$-\'equivalence (toujours en utilisant le th\'eor\`eme de Ramsey). Soit alors $(a_i')_{i\geq 0}$ une sous-suite de $(a_i)_{i\geq 0}$ tel que $c_i'\in I_{\leq}(a_i',a_{i+1}')$. Posons $A_1=\{a_i', i\in\mathbb N\}$ et $A_2=\{c_i', i\in\mathbb N\}$.
    Les ensembles $A_1$ et $A_2$ v\'erifient l'assertion 2 du lemme.
\end{proof}

\vspace{2mm}

La preuve du Th\'eor\`eme \ref{thm:graph-ordonne} se fait comme suit. Soit $\mathcal G:=(V,\leq,\rho)$ un graphe dirig\'e ordonn\'e qui a une infinit\'e de classes d'\'equivalence, alors d'apr\`es le Lemme \ref{lem:separation2}, nous avons deux cas.
%\medskip

\subsubsection{Etude du premier cas}

%\textbf{Cas 1:}
Si $\mathcal G$ v\'erifie l'assertion 1 du Lemme \ref{lem:separation2}, alors, comme il a \'et\'e fait, ci-dessus, pour les graphes non dirig\'es, nous pouvons trouver deux applications  $f:\NN\longrightarrow V$ et $g:[\NN]^2\longrightarrow V$ telles que $f(\NN)=A$ ($A$ donn\'e par le Lemme \ref{lem:separation2}). Ainsi $f(n)$ et $f(m)$ ne sont pas $1$-\'equivalents pour tous $n<m$ et $g(n,m)$ t\'emoigne de ce fait, c'est \`a dire que les restrictions de $\mathcal G$ \`a $\{f(n),g(n,m)\}$ et $\{f(m),g(n,m)\}$ ne sont pas isomorphes pour tous $n<m$.

\vspace{1mm}

Soit $\Phi:=\{f,g\}$ et $\mathfrak{L}:=\left\langle \omega,\mathcal G,\Phi \right\rangle $.
%Then, the restriction of $\mathcal G$ to the union of images of these two maps has infinitely many equivalence classes.
Le Th\'eor\`eme \ref{thm:ramsey-invariant} permet de trouver un sous-ensemble infini $X\subseteq\NN$ tel que  $\mathfrak L_{\restriction_X}$ est invariante. %the maps are invariant for $\mathcal G$ ($\mathfrak{L}:=\left\langle \omega,\mathcal G,\Phi \right\rangle $ with $\Phi:=\{f,g\}$, $\mathfrak L\restriction_X$ is invariant).
En indexant les \'el\'ements de  $X$ par des entiers, nous pouvons supposer $X=\NN$ et ainsi $\mathfrak{L}$ est invariante.

\begin{claim}\label{order-graphclaim}
\begin{enumerate}
\item $f(n)\leq f(m)\Leftrightarrow f(n')\leq f(m'),~\forall n<m,~n'<m'$.
\item $d(f(n), f(m))=d(f(n'), f(m')),~\forall n<m,~n'<m'$.
\item $g(n,m)\leq g(k,l)\Leftrightarrow g(n',m')\leq g(k',l'),~\forall n<m\leq k<l,~n'<m'\leq k'<l'$.
\item $d(g(n,m),g(k,l))= d(g(n',m'),g(k',l')),~\forall n<m\leq k<l,~n'<m'\leq k'<l'$.
%\item $d(g(n,m),g(k,l))= d(g(n',m'),g(k',l')),~\forall n<m<k<l,~n'<m'<k'<l'$.
\item $g(n,m)\in I_{\leq}(f(n),f(m))\Leftrightarrow g(n',m')\in I_{\leq}(f(n'),f(m')),~\forall n<m,~n'<m'$.
\item $g(n,m)\leq f(k)$ pour certains entiers $n<m<k \Leftrightarrow g(n,m)\leq f(l)$ pour tout $l>m$.
\item $d(f(n),g(n,m))= d(f(k),g(k,l)),~\forall n<m,~k<l$.
\item $d(g(n,m),f(k))= d(g(p,q),f(l)),~\forall n<m<k,~p<q<l$.
\item Si les restrictions   $\mathcal G_{\restriction_{\{f(n),g(n,m)\}}}$ et  $\mathcal G_{\restriction_{\{f(k),g(n,m)\}}}$  sont isomorphes pour des entiers $n<m<k$ alors   $\mathcal G_{\restriction_{\{f(n'),g(n',m')\}}}$ et  $\mathcal G_{\restriction_{\{f(k'),g(n',m')\}}}$ sont isomorphes pour tous $n'<m'<k'$.
\item $g(n,m)$ et $f(k)$ sont distincts pour tous entiers $n<m$ et $k$ .
\item $g(n,m)\neq g(n',m')$ pour tous $n<m<n'<m'$.
\end{enumerate}
\end{claim}

\begin{proofclaim}
La preuve se fait de la m\^eme façon que pour l'Affirmation \ref{graphclaim}.
%\hfill       $\Box$
\end{proofclaim}
\medskip

 Soit $F:\NN\times\{0,1\}\longrightarrow V(\mathcal G)$. Posons $F(n,1):=g(2n,2n+1)$. A partir de l'assertion  9 de l'Affirmation \ref{order-graphclaim} nous avons deux cas. Si $\mathcal G_{\restriction_{\{f(n),g(n,m)\}}}$ et  $\mathcal G_{\restriction_{\{f(k),g(n,m)\}}}$ sont isomorphes pour certains entiers $n<m<k$, alors nous posons $F(n,0):=f(2n+1)$. Autrement nous posons $F(n,0):=f(2n)$.

Nous d\'efinissons un graphe dirig\'e ordonn\'e $\mathcal G_1:=(V_1,\leq_1,\rho_1)$ ayant pour ensemble de sommets $V_1:=\NN\times\{0,1\}$ tel que
$$\left\{\begin{array}{l} x\leq_1 y\Leftrightarrow F(x)\leq F(y)\\ d_1(x,y)=d(F(x),F(y))
\end{array}\right.$$
pour tous $x,y\in \NN\times\{0,1\}$, o\`u $d_1(x,y):=(\rho_1(x,y),\rho_1(y,x))$. % $\rho'$ being the associated binary relation.

Par construction de $\mathcal G_1$ nous avons:
\begin{equation}\label{eq2}
{\mathcal G_1}_{\restriction_{\{(n,0),(n,1)\}}} \text{ et }{\mathcal G_1}_{\restriction_{\{(n,1),(m,0)\}}} \text{ ne sont pas isomorphes.}
\end{equation}
pour tous entiers $n<m$. Donc $\mathcal G_1$ poss\`ede une infinit\'e de classes d'\'equivalence  et ainsi  $\mathcal G_1\in\mathscr D$.
\medskip

 Par l'invariance de $\mathfrak L$ le graphe $\mathcal G_1$ est une presque  multicha\^{\i}ne et donc $\mathcal G_{\restriction_{Im(F)}}$ est presque  multi-encha\^{\i}nable. %almost multichainable. %the following claim

\vspace{2mm}

Pour construire $\mathcal G_1$, comme pour les graphes non dirig\'es, il suffit de d\'ecider sur les relations entre les quatre sommets $(0,0)$, $(0,1)$, $(1,0)$ et $(1,1)$, les relations restantes seront d\'eduites \`a partir des isomorphismes locaux de $C:=(\NN,\leq_1)$ et de l'Affirmation \ref{order-graphclaim}.
\smallskip

\noindent  Par exemple si nous supposons que $(0,1)$ est entre $(0,0)$ et $(1,0)$ pour l'ordre $\leq_1$, alors
\begin{enumerate}
\item si $(0,0)<_1(1,0)$, n\'ecessairement nous avons $(0,0)<_1(0,1)<_1(1,0)$. L'invariance de la structure donne $(1,0)<_1(1,1)$. Par la  transitivit\'e de l'ordre nous obtenons $(0,1)<_1(1,1)$ et donc $(0,0)<_1(1,1)$. Soit $n<m$, par l'invariance, les in\'equations pr\'ec\'edentes sont v\'erifi\'ees si nous  remplaçons $0$ et $1$ dans la premi\`ere composante de chacune des paires ci-dessus par $n$ et $m$. Alors $(n,0)<_1(n,1)<_1(m,0)<_1(m,1)$ pour tous $n<m$. Donc $\leq_1$ est ordonn\'e comme $\omega$, la cha\^{i}ne des entiers naturels. Nous consid\'erons alors toutes les situations pour la relation $\rho_1$ en prenant en consid\'eration la condition \eqref{eq2} qui est dans ce cas;
    \begin{equation}
d_1((n,0),(n,1))\neq d_1((n,1),(m,0)),~\forall n<m\in\NN.\label{eq:3}
\end{equation}
Ceci donne  $3.(4^4)$ cas, mais certains graphes peuvent s'abriter dans d'autres. %we obtain, in this case two hundred and eighty eight graphs.
\item Si $(1,0)<_1(0,0)$, les m\^emes consid\'erations que pr\'ec\'edemment donnent $\leq_1$  ordonn\'e comme $\omega^*$. Nous obtenons les m\^emes cas pour $\rho_1$.
\end{enumerate}
 Et si nous supposons que $(0,1)$ n'est pas entre  $(0,0)$ et $(1,0)$ pour $\leq_1$, alors, avec les m\^emes arguments, $\leq_1$ est isomorphe \`a un ordre de l'ensemble  $\{\underline{2}^*.\omega, \underline{2}.\omega^*, \omega+\omega, \omega^*+\omega, \omega+\omega^*, \omega^*+\omega^*\}$, o\`u $\underline{2}$ est la cha\^{i}ne \`a deux \'el\'ements $\{0,1\}$ ordonn\'es naturellement, $\underline{2}^*$ son ordre dual, $\underline{2}^*.\omega$ est le produit lexicographique de $\underline{2}^*$ et $\omega$, c'est \`a dire la cha\^{i}ne $\omega$ dans laquelle on substitut chaque \'el\'ement par la cha\^{i}ne $\underline{2}^*$.
  La relation $\rho_1$ est d\'efinie comme ci-dessus, la condition \eqref{eq2} est dans ce cas
 \begin{equation} d_1((n,0),(n,1))\neq d_1((m,0),(n,1)),~\forall n<m\in\NN.\label{eq:4}
\end{equation}

Soit $\mathfrak A_1$ l'ensemble des graphes dirig\'es ordonn\'es minimaux (pour l'abritement) induits par ce cas.

%\bigskip

 \subsubsection{Etude du second cas}

 Si $\mathcal G$ v\'erifie l'assertion 2 du Lemme \ref{lem:separation2}, alors, par construction de $A_1$ et $A_2$ du Lemme \ref{lem:separation2} nous avons:

 - $d$ est constante sur toute paire $(x,y)$   d'\'el\'ements distincts de $A_i,~i\in\{1,2\}$ telle que $x<y$ %relation $\rho$ is constant on $A_1$ and on $A_2$
et

- tout sommet de $A_i$ se trouve entre deux sommets de $A_j$ (relativement \`a l'ordre $\leq$) pour $i,j\in\{1,2\}, i\neq j$.

\noindent Alors $(A_1\cup A_2,\leq_{\restriction_{A_1\cup A_2}})$ est ordonn\'e comme $\omega$ ou $\omega^*$. %the relation between $A_1$ and $A_2$ are also constant.
Les ensembles $A_1$ et $A_2$ \'etant, chacun, contenu dans une classe de $1$-\'equivalence, nous avons alors deux situations:

\vspace{1mm}

 \paragraph{Premi\`ere situation.} Il existe $x_0\in V\setminus (A_1\cup A_2)$ qui t\'emoigne du fait que $A_1$ et $A_2$ sont dans deux  classes de $1$-\'equivalence diff\'erentes. %To simplify the enumeration of different cases, we can define two functions;
Comme $A_1\cup A_2$ est ordonn\'e comme $\omega$ ou $\omega^{\star}$, nous pouvons supposer que nous avons ou bien $x_0\leq a$ ou bien $x_0\geq a$, pour tout $a\in A_1\cup A_2$, car sinon, nous pouvons trouver un sous-ensemble infini de $A_1\cup A_2$ pour lequel c'est v\'erifi\'e.
\smallskip

    Nous pouvons alors trouver des applications
    $f',g',g'':\NN\longrightarrow V$ telles que, $f'(\NN)=A_1$, $g'(\NN)=A_2$,  $g''(\NN):=\{x_0\}$, avec  $g'(i)\in I_{\leq}(f'(i),f'(i+1)),~\forall i\in\NN$ et les restrictions de $\mathcal G$ \`a $\{ f'(i), g'(i), g''(0)\}$ et $\{ f'(j), g'(i), g''(0)\}$ ne sont pas isomorphes pour tous $i,j\in\NN$.
\smallskip

    Posons $F':=\{a\}\cup(\NN\times\{0,1\})\longrightarrow V$, avec $a\notin\NN\times\{0,1\}$,  telle que $F'(a):=g''(0)$, $F'(n,0):=f'(n)$ et $F'(n,1):=g'(n)$. Soit $\mathcal G_2:=(V_2,\leq_2,\rho_2)$ un graphe ordonn\'e ayant pour sommets l'ensemble $V_2=\{a\}\cup(\NN\times\{0,1\})$ tel que
$$\left\{\begin{array}{l} x\leq_2 y\Leftrightarrow F'(x)\leq F'(y)\\ d_2(x,y)=d(F'(x),F'(y))
\end{array}\right.$$
pour tous $x,y\in V_2$, o\`u $d_2(x,y):=(\rho_2(x,y),\rho_2(y,x))$. % $\rho_2$ being the associated binary relation.
\smallskip

Par construction, $\mathcal G_2$ v\'erifie la relation \eqref{eq3} ci-dessous, ainsi $\mathcal G_2$ poss\`ede une infinit\'e de classes d'\'equivalence.
\begin{equation}\label{eq3}
\forall n<m\in\NN,~{\mathcal G_2}_{\restriction_{\{a,(n,0),(n,1)\}}} \text{ et }{\mathcal G_2}_{\restriction_{\{a,(n,1),(m,0)\}}} \text{ ne sont pas isomorphes.}
\end{equation}

$\mathcal G_2$ v\'erifie \'egalement l'Observation \ref{observation1} ci-dessous qui se d\'eduit directement du fait que les \'el\'ements de $A_i$ pour $i\in\{0,1\}$ sont $1$-\'equivalents.
\begin{observation}\label{observation1}
\begin{enumerate}
\item $a\leq_2 (n,i)\Leftrightarrow a\leq_2 (m,i),~\forall n<m,~ i\in\{0,1\}$.
\item $a\leq_2 (n,0)\Leftrightarrow a\leq_2 (n,1),~\forall n\in\NN$.
\item $d_2(a,(n,i))=d_2(a,(m,i)), ~\forall n<m,~ i\in\{0,1\}$.
%\item $d''(0,(n,1))=d''(0,(m,1)), ~\forall n<m,~n,m\in\NN^*$.
%\item $(n,0)\leq'' (m,0)\Leftrightarrow (n',0)\leq'' (m',0),~\forall n<m,~n'<m'$.
\item $d_2((n,i),(m,i))=d_2((n',i),(m',i)), \forall n<m,~n'<m',~ i\in\{0,1\}$.
%\item $(n,1)\leq'' (m,1)\Leftrightarrow (n',1)\leq'' (m',1),~\forall n<m,~n'<m'$.
%\item $d''((n,1),(m,1))=d''((n',1),(m',1)), \forall n<m,n'<m'; n,m,n',m'\in\NN^*$.
%\item $(n,0)\leq'' (n,1)\Leftrightarrow (m,0)\leq'' (m,1),~\forall n<m$.
\item $d_2((n,0),(n,1))=d_2((m,0),(m,1)), ~\forall n<m$.
%\item $(n,0)\leq'' (n,1)\Leftrightarrow (n,1)\leq'' (m,0),~\forall n<m$.
\item $d_2((n,0),(n,1))=d_2((n,1),(m,0)), ~\forall n<m$.
%\item $(n,0)\leq'' (m,1)\Leftrightarrow (m,1)\leq'' (m+1,0),~\forall n<m$.
\item $d_2((n,0),(m,1))=d_2((m,1),(m+1,0)), ~\forall n<m$.
\end{enumerate}
\end{observation}

Il s'ensuit que $\mathcal G_2$ est une presque multicha\^{\i}ne.
Avec ceci, nous obtenons un sous-ensemble fini $\mathfrak B_1$ de graphes ordonn\'es ayant le m\^eme ensemble de sommets $\{a\}\cup(\NN\times\{0,1\})$. %if $\leq$ is ordered as $\omega$ and the same number if it is ordered as $\omega^*$.
\medskip

\paragraph{Deuxi\`eme situation.} Il n'exite pas de sommet $x_0$ comme ci-dessus, alors, d'apr\`es le Lemme \ref{lem:separation-binaire}, deux sommets $x,y$ de $A_i$ sont s\'epar\'es par deux sommets $z,z'$ avec $z\in I_{\leq}(x,y)\cap A_j$, $j\neq i$ et $z'\notin I_{\leq}(x,y)$. Dans ce cas et \`a partir des Lemmes \ref{lem:separation-binaire} et \ref{lem:separation2}, la relation entre deux \'el\'ements de $A_i$, pour au moins un $i=1,2$, est diff\'erente de la relation entre deux \'el\'ements qui se trouvent l'un dans $A_1$ et l'autre dans $A_2$. % $x,y$, avec $x\in A_1$ et $y\in A_2$.
Nous pouvons alors d\'efinir deux applications $f_1, g_1:\NN\longrightarrow V$ telles que, $f_1(\NN)=A_1$, $g_1(\NN)=A_2$, $g_1(i)\in I_{\leq}(f_1(i),f_1(i+1)),~\forall i\in\NN$. Posons $F'':\NN\times\{0,1\}\longrightarrow V$ telle que $F''(n,0):=f(n)$ et $F''(n,1):=g(n)$.\\
Nous pouvons d\'efinir un graphe dirig\'e ordonn\'e $\mathcal G_3:=(V_3,\leq_3,\rho_3)$ avec l'ensemble de sommets $V_3=\NN\times\{0,1\}$ tel que
$$\left\{\begin{array}{l} x\leq_3 y\Leftrightarrow F''(x)\leq F''(y)\\ d_3(x,y)=d(F''(x),F''(y))
\end{array}\right.$$
pour tous $x,y\in \NN\times\{0,1\}$, o\`u $d_3(x,y):=(\rho_3(x,y),\rho_3(y,x))$. Comme nous l'avons signal\'e auparavant, par construction de $\mathcal G_3$, l'ordre $\leq_3$ est isomorphe \`a $\omega$ ou $\omega^*$ avec $(n,1)\in I_{\leq_3}((n,0),(n+1,0))$ et
pour tous $n<m$ nous avons l'un des trois cas $$d_3((n,0),(m,0))\neq d_3((n,0),(m,1)),$$ ou  $$d_3((n,0),(m,0))\neq d_3((n,1),(m,1)),$$ ou  $$d_3((n,0),(n,1))\neq d_3((n,1),(m,1)).$$
Alors, avec le fait que $A_1$ et $A_2$ sont chacun dans une classe de $1$-\'equivalence, $\mathcal G_3$ v\'erifie l'Observation \ref{observation2} ci-dessous:
\begin{observation}\label{observation2}

\begin{enumerate}
\item $d_3((n,i),(m,i))=d_3((n',i),(m',i)), ~\forall n<m,~n'<m',~~i\in\{0,1\}$.
\item $d_3((n,0),(n,1))=d_3((m,0),(m,1)), ~\forall n<m$.
\item $d_3((n,0),(n,1))=d_3((n,1),(m,0)), ~\forall n<m$.
\item $d_3((n,0),(m,1))=d_3((m,1),(m+1,0)), ~\forall n<m$.
\end{enumerate}
\end{observation}

Il est alors clair que $\mathcal G_3$ est une presque multicha\^{\i}ne d\'efinie sur $\NN\times\{0,1\}$ et ayant une infinit\'e de classes d'\'equivalence, donc $\mathcal G_3\in\mathscr D$.
Soit $\mathfrak B_2$ l'ensemble des graphes dirig\'es ordonn\'es obtenus dans ce cas. Il s'av\`ere que des membres de $\mathfrak A_1\cup\mathfrak B_1$ peuvent abriter un membre de $\mathfrak B_2$. Soit alors  $\mathfrak A$ l'ensemble des  graphes  minimaux (pour l'abritement) de $\mathfrak A_1\cup\mathfrak B_1\cup\mathfrak B_2$. Cet ensemble contient mille deux cent quarante six graphes. Ils sont d\'ecris ci-dessous.
%\vspace{2mm}

        \subsubsection{Description des graphes}
 Si $G$ est un graphe dirig\'e, le graphe $G'$, obtenu \`a partir de $G$ en rajoutant tout arc $u$ tel que  $u^{-1}$ est un arc de $G$ est le \emph{sym\'etris\'e}\index{graphe!sym\'etris\'e d'un -} de $G$. Ainsi, $G'$ est un graphe dirig\'e dans lequel la relation $\rho$ est sym\'etrique et donc $G$ peut \^etre consid\'er\'e comme un graphe non dirig\'e.

\vspace{2mm}

  D\'esignons par $\mathcal G_{\ell,k}^{(p)}:=(V_{\ell,k}^{(p)},\leq_{\ell,k}^{(p)},\rho_{\ell,k}^{(p)})$ les graphes dirig\'es ordonn\'es de $\mathfrak A$ (du Th\'eor\`eme \ref{thm:graph-ordonne}), o\`u $p, \ell$ et $k$ sont des entiers tels que $1\leq p\leq 10$, $1\leq \ell\leq 8$. L'ensemble des sommets  $V_{\ell,k}^{(p)}$ est soit $\NN\times\{0,1\}$ (si $\mathcal G_{\ell,k}^{(p)}$ est dans $\mathfrak A_1\cup \mathfrak B_2$) soit $\{a\}\cup\NN\times\{0,1\}$ (si $\mathcal G_{\ell,k}^{(p)}$ est dans $\mathfrak B_1$).

\vspace{1mm}

  \noindent Les  graphes ordonn\'es avec la m\^eme valeur de $p$ sont dits de \emph{classe} $p$, tous ont leurs restrictions \`a $A:=\NN\times\{0\}$ (resp. $B:=\NN\times\{1\}$) qui sont isomorphes. S'il ont la m\^eme valeur de  $\ell$, leurs ordres lin\'eaires sont de m\^eme type, $\ell$ prend les valeurs de $1$ \`a $8$ si l'ordre lin\'eaire est isomorphe \`a $\omega$, $\omega^*$, $\underline{2}^*.\omega$, $\underline{2}.\omega^*$, $\omega+\omega$, $\omega^*+\omega$, $\omega+\omega^*$, $\omega^*+\omega^*$ respectivement.

  \medskip

   L'entier $k$ \'enum\`ere les graphes pour toutes valeurs de  $p$ et $\ell$. Diff\'erentes classes n'ont pas n\'ecessairement les m\^emes cardinalit\'es.

\medskip

   Pour $p=1$ si $\ell=1,2$ nous avons $1\leq k\leq 18$ et si $5\leq \ell\leq 8$ nous avons $1\leq k\leq 15$, il n'y a pas de graphe si $\ell=3,4$. Pour $p=2,3$, si $\ell=1,2$ nous avons $1\leq k\leq 21$, si $\ell=3,4$, alors $k=1$ et si $5\leq \ell\leq 8$ nous avons  $1\leq k\leq 15$. Pour $p=4$ nous avons les m\^emes valeurs de  $k$ que pour  $p=3$ except\'e pour $\ell=3,4$ o\`u nous n'avons aucun graphe. Pour $5\leq p\leq 10$, si $\ell=1,2$ nous avons  $1\leq k\leq 22$ et si $5\leq \ell\leq 8$ nous avons  $1\leq k\leq 24$, pas de graphes si  $\ell=3,4$. Le total est de mille deux cent quarante six graphes comme annonc\'e. %(1246).%c'est ce que je trouve en faisant le d\'ecompte

  \vspace{1mm}

  Dans chaque classe $p$, lorsque $\ell=1$, (l'ordre lin\'eaire  $\leq_{\ell,k}^{(p)}$ est isomorphe \`a $\omega$), nous avons, %for all $n\in\NN$,
    $(0,0)<_{\ell,k}^{(p)}(0,1)<_{\ell,k}^{(p)}(1,0)$  quand l'ensemble des sommets est $\NN\times\{0,1\}$  et $a<_{\ell,k}^{(p)}(0,0)<_{\ell,k}^{(p)}(0,1)<_{\ell,k}^{(p)}(1,0)$  lorsqu'il est $\{a\}\cup\NN\times\{0,1\}$. Les relations sont invers\'ees lorsque  $\ell=2$.
Si $\ell\geq 3$, l'ensemble des sommets est $\NN\times\{0,1\}$.
\vspace{1mm}

Nous donnons les repr\'esentations des graphes pour certaines classes. Toutes ces repr\'esentations  sont faites sur les six sommets $\{(0,0),(1,0),(2,0),(0,1),(1,1),(2,1)\}$ pour les graphes de $\mathfrak A_1\cup\mathfrak B_2$ et sur les sept sommets $\{a,(0,0),(1,0),(2,0),(0,1),(1,1),(2,1)\}$ pour ceux de $\mathfrak B_1$ (les boucles ne sont pas repr\'esent\'ees), l'ordre $\leq_{\ell,k}^{(p)}$ \'etant isomorphe \`a $\omega$.
\vspace{1mm}

  Pour \'eviter des complications inutiles, nous supposons que dans chacun de ces graphes la relation $\rho_{\ell,k}^{(p)}$ est r\'eflexive (nous pouvons la supposer irr\'eflexive, les graphes obtenus vont coïncider sur les paires de sommets distincts avec les graphes obtenus pour $\rho$ r\'eflexive) (les autres cas peuvent se d\'eduire facilement). Nous pouvons alors d\'ecrire les graphes dirig\'es associ\'es $G_{\ell,k}^{(p)}=(V_{\ell,k}^{(p)},E_{\ell,k}^{(p)})$. \\

 Pour $n\in\NN$, soit $a_n:=((n,0),(n,1))$.

  \textbf{\underline{Classe $p=1$}:} Les restrictions de $G_{\ell,k}^{(1)}$ aux ensembles  $A:=\NN\times\{0\}$ et $B:=\NN\times\{1\}$  sont toutes les deux des anticha\^{\i}nes.
\medskip

  \textbf{I)} Si $\ell=1,2$ alors $1\leq k\leq 18$. %Consider first that the linear order is isomorphic to $\omega$, in this case $1\leq k\leq 18$.
   Les graphes $\mathcal G_{\ell,k}^{(1)}$ pour $1\leq k\leq 9$ sont dans $\mathfrak A_1$, ils sont dans $\mathfrak B_2$ pour $10\leq k\leq 12$ et dans $\mathfrak B_1$ pour $13\leq k\leq 18$.
\vspace{1mm}

  $\bullet$  Pour $1\leq k\leq 12$. %the vertex set of $\mathcal G_k^{(1)}$ is $A\cup B$. %These graphs are in $\mathfrak A_1$ for $1\leq k\leq 9$ and in $\mathfrak B_2$ for $10\leq k\leq 12$.
   Une paire $(x,x')$ de sommets, o\`u $x=(n,i), x'=(n',i')$, est
\begin{itemize}
\item un arc de $G_{\ell,1}^{(1)}$ si $n=n'$ et $i<i'$;
\item un arc de $G_{\ell,2}^{(1)}$ si $(x',x)$ est un arc de $G_{\ell,1}^{(1)}$. Ainsi $G_{\ell,2}^{(1)}$ est le dual de $G_{\ell,1}^{(1)}$;
\item un arc de $G_{\ell,3}^{(1)}$ si $n=n'$ et $i\neq i'$. Le graphe $G_{\ell,3}^{(1)}$ est auto-dual;
\item un arc de $G_{\ell,4}^{(1)}$ si $n\leq n'$ et $i<i'$;
\item un arc de $G_{\ell,5}^{(1)}$ si $(x',x)$ est un arc de $G_{\ell,4}^{(1)}$. Ainsi $G_{\ell,5}^{(1)}$ est le dual de $G_{\ell,4}^{(1)}$;
\item un arc de $G_{\ell,6}^{(1)}$ s'il est soit un arc de $G_{\ell,4}^{(1)}$ soit un arc de $G_{\ell,5}^{(1)}$. Ainsi  $G_{\ell,6}^{(1)}$ est le sym\'etris\'e de $G_{\ell,4}^{(1)}$ (et de $G_{\ell,5}^{(1)}$), il est auto-dual;
\item un arc de $G_{\ell,7}^{(1)}$ si  $i<i'$. Le graphe $G_{\ell,7}^{(1)}$ est \'equimorphe \`a son dual. %et le graphe ordonn\'e $\mathcal G_{\ell,7}^{(1)}$ est \'equimorphe \`a son dual;
 \item un arc de $G_{\ell,8}^{(1)}$ si ou bien ($n\leq n'$ et $i<i'$) ou bien ($n>n'$ et $i<i'$) ou bien ($n<n'$ et $i>i'$). %Le graphe $G_{l,8}^{(1)}$ est auto-dual;
\item un arc de $G_{\ell,9}^{(1)}$ si $(x',x)$ est un arc de $G_{\ell,8}^{(1)}$. Ainsi  $G_{\ell,9}^{(1)}$ est le dual de $G_{\ell,8}^{(1)}$;
\item un arc de $G_{\ell,10}^{(1)}$ si ou bien ($n\leq n'$ et $i<i'$) ou bien ($n<n'$ et $i>i'$);
\item un arc de $G_{\ell,11}^{(1)}$ si $(x',x)$ est un arc de $G_{\ell,10}^{(1)}$. Ainsi $G_{\ell,11}^{(1)}$ est le dual de $G_{\ell,10}^{(1)}$;
\item un arc de $G_{\ell,12}^{(1)}$ si  $i\neq i'$. Le graphe $G_{\ell,12}^{(1)}$ est le sym\'etris\'e de $G_{\ell,10}^{(1)}$ (et de $G_{\ell,11}^{(1)}$).
\end{itemize}

\vspace{2mm}

D\'esignons par $\underline{2}$ l'ensemble ordonn\'e de base $2:=\{0,1\}$  tel que $0<1$. L'ensemble ordonn\'e $\underline{2}^*$ est son dual. D\'esignons par  $K_2$ le graphe complet r\'eflexif \`a deux sommets et par $\Delta_{\NN}$ l'anticha\^{\i}ne ayant $\NN$ comme ensemble de sommets. Alors $G_{\ell,1}^{(1)}$ est isomorphe \`a $\underline{2}.\Delta_{\NN}$, le produit lexicographique de $\underline{2}$ par $\Delta_{\NN}$. Le graphe $G_{\ell,2}^{(1)}$ est isomorphe \`a $\underline{2}^*.\Delta_{\NN}$, $G_{\ell,3}^{(1)}$ est isomorphe \`a $K_2.\Delta_{\NN}$, le graphe $G_{\ell,6}^{(1)}$ est le graphe biparti demi-complet de Shmerl- Trotter et le graphe $G_{\ell,7}^{(1)}$ est isomorphe \`a la somme ordinale %\footnote{Voir la d\'efinition dans le paragraphe \ref{subsec:somme et produit} du chapitre \ref{chap:generalite}.}
$\Delta_{\NN}+\Delta_{\NN}$.

\vspace{2mm}

$\bullet$  Pour $13\leq k\leq 18$, l'ensemble des sommets est $\{a\}\cup A\cup B$.  %the vertex set of $\mathcal G_k^{(1)}$ is $\{a\}\cup A\times B$ with  $A$ and $B$  both antichains. These graphs are in $\mathfrak B_1$.
Une paire $(x,x')$ de sommets est
\begin{itemize}
\item un arc de $G_{\ell,13}^{(1)}$ si $x=a, x'=(n,1)$;
\item un arc de $G_{\ell,14}^{(1)}$ si $(x',x)$ est un arc de $G_{\ell,13}^{(1)}$. Ainsi $G_{\ell,14}^{(1)}$ est le dual de $G_{\ell,13}^{(1)}$;
\item un arc de $G_{\ell,15}^{(1)}$ s'il est soit un arc de $G_{\ell,13}^{(1)}$ soit un arc de $G_{\ell,14}^{(1)}$. Ainsi  $G_{\ell,15}^{(1)}$ est le sym\'etris\'e de $G_{\ell,13}^{(1)}$ (et de $G_{\ell,14}^{(1)}$);
\item un arc de $G_{\ell,16}^{(1)}$ si ou bien $x=a$ et $x'=(n,0)$ ou bien $x=(n,1)$ et $x'=a$; ce graphe est auto-dual;
\item un arc de $G_{\ell,17}^{(1)}$ s'il est soit un arc de $G_{\ell,16}^{(1)}$ soit un arc de $G_{\ell,13}^{(1)}$;
\item un arc de $G_{\ell,18}^{(1)}$ si $(x',x)$ est un arc de $G_{\ell,17}^{(1)}$. Ainsi $G_{\ell,18}^{(1)}$ est le dual de $G_{\ell,17}^{(1)}$;
\end{itemize}
\medskip

\begin{figure}[!hbp]
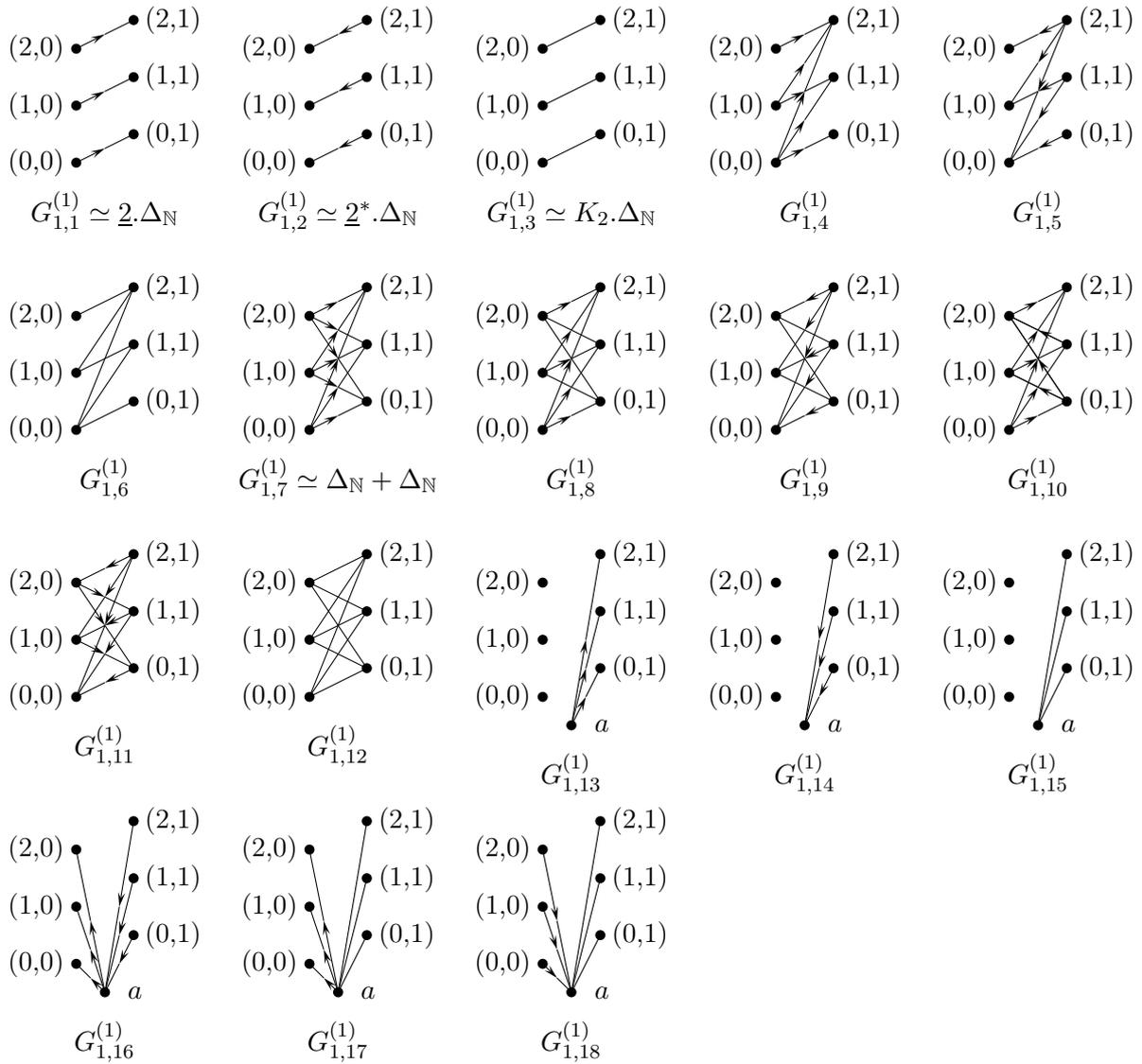

\begin{center}
\small
\begin{tabular}[l]{ccccccccc}
\input{graph1-class1}&&\input{graph2-class1}&&\input{graph3-class1}&&\input{graph4-class1}&&\input{graph5-class1}\\
\input{graph6-class1}&&\input{graph7-class1}&&\input{graph8-class1}&&\input{graph9-class1}&&\input{graph10-class1}\\
\input{graph11-class1}&&\input{graph12-class1}&&\input{graph13-class1}&&\input{graph14-class1}&&\input{graph15-class1}\\
\input{graph16-class1}&&\input{graph17-class1}&&\input{graph18-class1}&&&&\\
&&&&&&&&\\
\end{tabular}
\caption{Les graphes minimaux de classe $p=1$ pour $\ell=1$.}
\label{repre:graphe-classe1}
\end{center}
\end{figure}

Ces graphe sont repr\'esent\'es dans la \figurename ~\ref{repre:graphe-classe1}.

\vspace{3mm}

\textbf{II)} Si $\ell=3,4$ il n'y a aucun graphe.

\vspace{3mm}

\textbf{III)} Si $5\leq \ell\leq 8$, nous avons les m\^emes exemples pour chaque valeurs de $\ell$ et leur nombre est $15$, dans l'ordre $\leq_{\ell,k}^{(1)}$, les \'el\'ements de $A$ sont plac\'es avant ceux de $B$.

$\centerdot$ $G_{\ell,k}^{(1)}=G_{1,k}^{(1)}$ pour tout $1\leq k\leq 6$.

$\centerdot$ $G_{\ell,k}^{(1)}=G_{1,k+1}^{(1)}$ pour tout $7\leq k\leq 10$.

$\centerdot$ $G_{\ell,11}^{(1)}$ est obtenu \`a partir de $G_{1,7}^{(1)}$ en rajoutant tous les arcs $((n,1),(m,0))$ pour $n\geq m$.

$\centerdot$ $G_{\ell,12}^{(1)}$ est obtenu \`a partir de $G_{1,10}^{(1)}$ en rajoutant tous les arcs $((n,1),(m,0))$ pour $n\geq m$, $G_{\ell,12}^{(1)}$ est le dual de $G_{\ell,11}^{(1)}$

$\centerdot$ $G_{\ell,13}^{(1)}$ est non dirig\'e, c'est le graphe $G_3$ de la sous-section \ref{subsec:graph} (voir \figurename~\ref{tab:graphes} en page \pageref{tab:graphes}).

$\centerdot$ $G_{\ell,14}^{(1)}$ est obtenu \`a partir de $G_{\ell,13}^{(1)}$ en rajoutant tous les arcs $a_n$ pour $n\in\NN$.

$\centerdot$ $G_{\ell,15}^{(1)}$ est obtenu \`a partir de $G_{\ell,13}^{(1)}$ en rajoutant tous les arcs $a_n^{-1}$ pour $n\in\NN$.

\vspace{2mm}

\vspace{3mm}

 \textbf{\underline{Classe $p=2$}:} Dans ce cas, les restrictions de $G_{\ell,k}^{(2)}$ \`a $A$ et $B$ sont des cha\^{\i}nes isomorphes \`a la cha\^{i}ne $\omega$.
\medskip

 \textbf{I)} Si $\ell=1,2$ nous avons $1\leq k\leq 21$, les graphes $\mathcal G_{\ell,k}^{(2)}$ pour $1\leq k\leq 12$ sont dans $\mathfrak A_1$, ils sont dans $\mathfrak B_1$ pour $13\leq k\leq 18$ et dans $\mathfrak B_2$ pour $19\leq k\leq 21$.

\vspace{1mm}

 $\bullet$  Pour $1\leq k\leq 9$ le graphe $G_{\ell,k}^{(2)}$ %has vertex set equal to $A\cup B$ and
 co\"{\i}ncide avec $G_{\ell,k}^{(1)}$ sur les paires de sommets de $A\times B$. Ainsi $G_{\ell,7}^{(2)}$ est une cha\^{\i}ne isomorphe \`a $\omega+\omega$.

\vspace{1mm}

 $\bullet$ Pour $10\leq k\leq 12$, le graphe $G_{\ell,k}^{(2)}$ %is $A\cup B$, it
 co\"{\i}ncide avec $G_{\ell,10}^{(1)}$ sur les paires de sommets de $A\times B$ avec
 \begin{enumerate}
 \item la suppression de l'arc $a_n, ~n\in\NN$ si $k=10$;  alors $G_{\ell,10}^{(2)}$ est isomorphe \`a $\Delta_2.\omega$, le produit lexicographique de l'anticha\^{\i}ne sur deux sommets $\Delta_2$ avec $\omega$, la cha\^{\i}ne des entiers de $\NN$.
 \item le remplacement de $a_n$ par $a_n^{-1},~n\in\NN$ si $k=11$; alors $G_{\ell,11}^{(2)}$ est isomorphe \`a $\underline{2}^*.\omega$.
 \item le rajout des arcs $a_n^{-1},~n\in\NN$ si $k=12$; alors $G_{\ell,12}^{(2)}$ est isomorphe \`a $K_{2}.\omega$.
 \end{enumerate}

\vspace{1mm}

 $\bullet$ Pour $13\leq k\leq 18$, l'ensemble des arcs sur $A\times B$ du graphe $G_{\ell,k}^{(2)}$ est l'union des ensembles d'arcs de $G_{\ell,10}^{(1)}$ et $G_{\ell,k}^{(1)}$.

 \vspace{1mm}

 $\bullet$ Les ensembles d'arcs des graphes $G_{\ell,19}^{(2)}$ et $G_{\ell,20}^{(2)}$ sur $A\times B$ co\"{\i}ncident avec ceux de  $G_{\ell,11}^{(1)}$ et $G_{\ell,12}^{(1)}$ respectivement.

 \vspace{1mm}

 $\bullet$ L'ensemble des arcs de $G_{\ell,21}^{(2)}$  sur $A\times B$ est vide.
 Alors $G_{\ell,21}^{(2)}$ est isomorphe \`a $\omega\oplus\omega$, la somme directe %\footnote{Voir la d\'efinition de la somme directe en \ref{subsec:somme et produit} du chapitre \ref{chap:generalite}.}
 de deux cha\^{\i}nes isomorphes \`a $\omega$.

\medskip

\begin{figure}[!hbp]
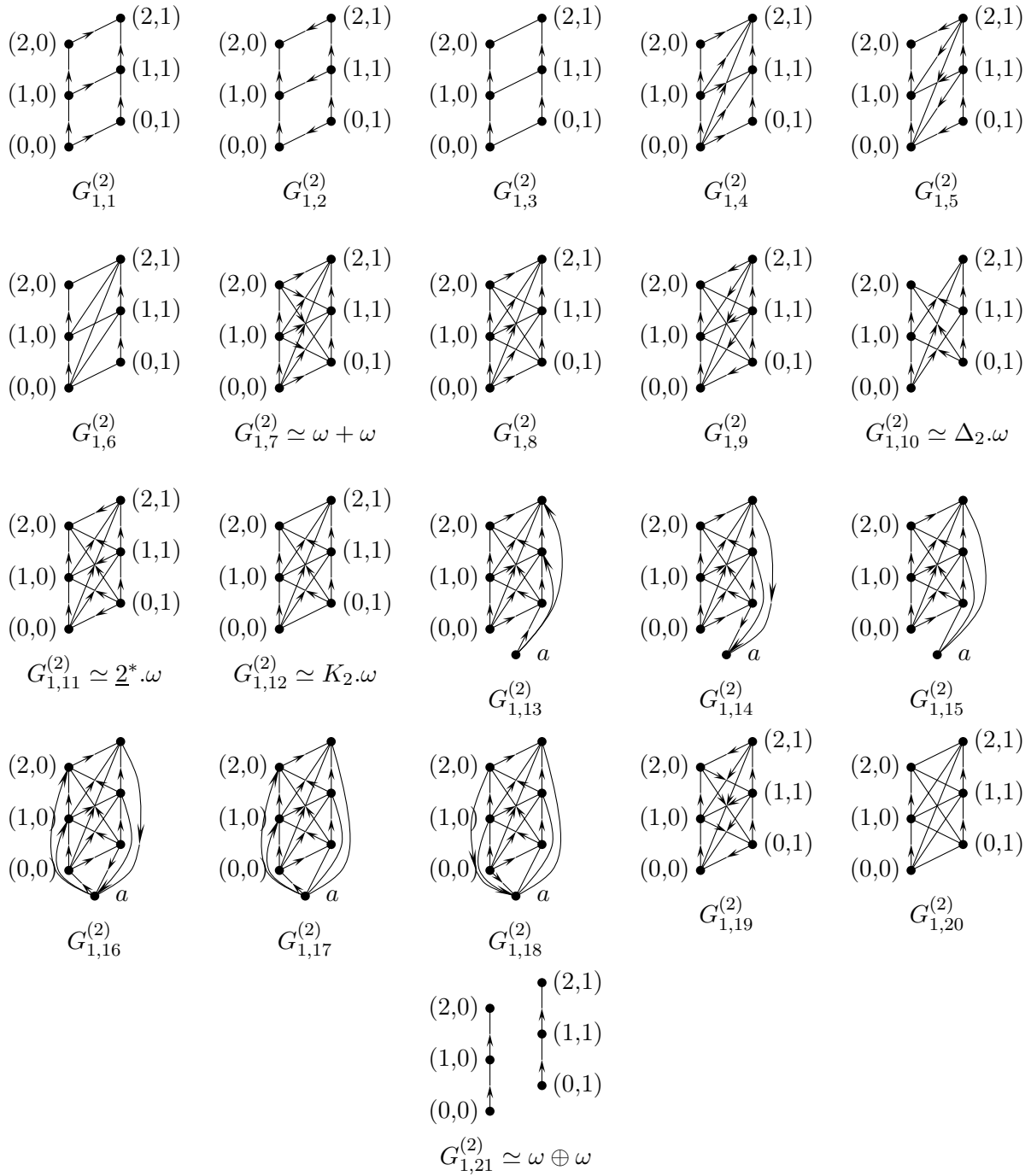

\begin{center}
\small
\begin{tabular}[l]{ccccccccc}
\input{graph1-class2}&&\input{graph2-class2}&&\input{graph3-class2}&&\input{graph4-class2}&&\input{graph5-class2}\\
\input{graph6-class2}&&\input{graph7-class2}&&\input{graph8-class2}&&\input{graph9-class2}&&\input{graph10-class2}\\
\input{graph11-class2}&&\input{graph12-class2}&&\input{graph13-class2}&&\input{graph14-class2}&&\input{graph15-class2}\\
\input{graph16-class2}&&\input{graph17-class2}&&\input{graph18-class2}&&\input{graph19-class2}&&\input{graph20-class2}\\
&&&&\input{graph21-class2}&&&&\\
\end{tabular}
\caption{Les graphes minimaux de classe $p=2$ pour $\ell=1$.}
\label{repre:graphe-classe2}
\end{center}
\end{figure}

Les graphes $G_{1,k}^{(2)}, 1\leq k\leq 21$ sont repr\'esent\'es dans la \figurename ~\ref{repre:graphe-classe2}.

 \vspace{3mm}

 \textbf{II)} Si $\ell=3,4$ dans ce cas $k=1$, nous avons un seul graphe %dans lequel  $\rho_{l,k}^{(2)}$ is a linear order isomorphic to $\omega$.
 qui est un ordre lin\'eaire isomorphe \`a $\omega$.

 \vspace{2mm}

 \textbf{III)} Si $5\leq \ell\leq 8$, nous avons les m\^emes exemples pour chaque valeur de $\ell$ et leur nombre est $15$.

 $\centerdot$ $G_{\ell,k}^{(2)}=G_{1,k}^{(2)}$ pour tout $1\leq k\leq 6$.

 $\centerdot$ $G_{\ell,7}^{(2)}$ est un ordre  lin\'eaire isomorphe \`a $\omega$.

$\centerdot$ $G_{\ell,k}^{(2)}=G_{1,k}^{(2)}$ pour tout $8\leq k\leq 9$.

$\centerdot$ $G_{\ell,10}^{(2)}=G_{1,19}^{(2)}$.

 $\centerdot$ Pour $11\leq k\leq 15$, le graphe $G_{\ell,k}^{(2)}$ co\"{\i}ncide sur $A\times B$ avec $G_{\ell,k}^{(1)}$.

 \vspace{2mm}

\vspace{3mm}

  \textbf{\underline{Classe $p=3$}:} Dans ce cas, les restrictions de $G_{\ell,k}^{(3)}$ \`a $A$ et $B$ sont des cha\^{\i}nes isomorphes \`a $\omega^*$.
\medskip

  \textbf{I)} Si $\ell=1,2$ nous avons $1\leq k\leq 21$, les graphes $\mathcal G_{\ell,k}^{(3)}$ pour $1\leq k\leq 12$ sont dans $\mathfrak A_1$, ils sont dans $\mathfrak B_1$ pour $13\leq k\leq 18$ et dans $\mathfrak B_2$ pour $19\leq k\leq 21$.

\vspace{1mm}

 $\bullet$ Si $1\leq k\leq 9$, le  graphe $G_{\ell,k}^{(3)}$ co\"{\i}ncide sur $A\times B$ avec $G_{\ell,k}^{(1)}$. Alors $G_{\ell,7}^{(3)}$ est une cha\^{\i}ne isomorphe \`a $\omega^*+\omega^*$.

\vspace{2mm}

$\bullet$ Si $10\leq k\leq 21$, le graphe $G_{\ell,k}^{(3)}$ est le dual de $G_{\ell,k}^{(2)}$. Ainsi le graphe $G_{\ell,10}^{(3)}$ est isomorphe \`a $\Delta_2.\omega^*$, le graphe $G_{\ell,11}^{(3)}$ est isomorphe \`a $\underline{2}.\omega^*$, le graphe $G_{\ell,12}^{(3)}$ est isomorphe \`a $K_{2}.\omega^*$ et le graphe  $G_{\ell,21}^{(3)}$ est isomorphe \`a $\omega^*\oplus\omega^*$.

\vspace{3mm}

 \textbf{II)} Si $\ell=3,4$, nous avons un seul graphe qui est un ordre lin\'eaire isomorphe \`a $\omega^*$. %we have only one example in which $\rho_{l,k}^{(3)}$ is a linear order isomorphic to $\omega^*$.

 \vspace{2mm}

 \textbf{III)} Si $5\leq \ell\leq 8$, nous avons les m\^emes exemples pour chaque valeur de $\ell$ et leur nombre est $15$.

 $\centerdot$ $G_{\ell,k}^{(3)}=G_{1,k}^{(3)}$ pour tout $1\leq k\leq 6$.

 $\centerdot$ $G_{\ell,7}^{(3)}$ est un ordre lin\'eaire isomorphe \`a $\omega^*$.

 $\centerdot$ $G_{\ell,k}^{(3)}=G_{1,k}^{(3)}$ pour tout $8\leq k\leq 9$.

$\centerdot$ $G_{\ell,10}^{(3)}=G_{1,19}^{(3)}$.

 $\centerdot$ Pour $11\leq k\leq 15$, le graphe $G_{\ell,k}^{(3)}$ co\"{\i}ncide sur $A\times B$ avec $G_{\ell,k}^{(1)}$.

\vspace{3mm}

%Set $H:=F\cup F'$.
 \textbf{\underline{Classe $p=4$}:} Dans ce cas, $A$ et $B$  sont tous les deux des  cliques r\'eflexives.

\medskip

  \textbf{I)} Si $\ell=1,2$ nous avons $1\leq k\leq 21$, les graphes pour $1\leq k\leq 12$ sont dans $\mathfrak A_1$, ils sont dans $\mathfrak B_1$ pour $13\leq k\leq 18$ et dans $\mathfrak B_2$ pour $19\leq k\leq 21$.

\vspace{1mm}

 $\centerdot$ Si $1\leq k\leq 9$, le graphe $G_{\ell,k}^{(4)}$ co\"{\i}ncide avec $G_{\ell,k}^{(1)}$ sur les paires de $A\times B$. Alors $G_{\ell,7}^{(4)}$ est isomorphe \`a $K_{\NN}+K_{\NN}$, la somme ordinale de deux cliques reflexives ayant chacune pour ensemble de sommets $\mathbb N$. %All these graphs have the same vertex set $\NN\times 2$.

%\vspace{1mm}

$\centerdot$ Le graphe $G_{\ell,10}^{(4)}$ est le sym\'etris\'e de $G_{\ell,10}^{(2)}$.

%\vspace{2mm}

$\centerdot$ Le graphe $G_{\ell,11}^{(4)}$ (respectivement $G_{\ell,12}^{(4)}$) est obtenu \`a partir de $G_{\ell,10}^{(4)}$ en rajoutant les arcs $a_n, n\in\NN$ (respectivement $a_n^{-1}, n\in\NN$). Le graphe $G_{\ell,12}^{(4)}$ est le dual de $G_{\ell,11}^{(4)}$.

\vspace{2mm}

$\centerdot$ Pour $13\leq k\leq 18$, le graphe $G_{\ell,k}^{(4)}$ est obtenu \`a partir de $G_{\ell,k}^{(2)}$ en prenant le sym\'etris\'e sur $A\cup B$, les arcs restants (c'est \`a dire ceux pour lesquels un des sommets est $a$) \'etant les m\^emes que dans $G_{\ell,k}^{(2)}$. %All these graphs have the same vertex set  $\{0\}\cup\NN^*\times 2$.

\vspace{2mm}
$\centerdot$ Le graphe $G_{\ell,19}^{(4)}$ co\"{\i}ncide avec $G_{\ell,11}^{(1)}$ sur les paires de $A\times B$.

\vspace{2mm}
$\centerdot$ Le graphe $G_{\ell,20}^{(4)}$ est le dual de $G_{\ell,19}^{(4)}$.

\vspace{2mm}
 $\centerdot$ Le graphe $G_{\ell,21}^{(4)}$ est le sym\'etris\'e de $G_{\ell,21}^{(2)}$, il est isomorphe \`a $K_{\NN}\oplus K_{\NN}$.

\vspace{3mm}

\textbf{II)} Si $\ell=3,4$ nous n'avons aucun graphe. %the linear order is isomorphic to an order of the following set $\{ \omega+\omega, \omega^*+\omega, \omega+\omega^*, \omega^*+\omega^*\}$, we have the same examples which are

\vspace{3mm}

\textbf{III)} Si $5\leq \ell\leq 8$, nous avons les m\^emes exemples pour chaque valeur de $\ell$ et leur nombre est  $15$.

 $\centerdot$ $G_{\ell,k}^{(4)}=G_{1,k}^{(4)}$ pour tout $1\leq k\leq 6$.

 $\centerdot$ $G_{\ell,k}^{(4)}=G_{1,k+1}^{(4)}$ pour tout $7\leq k\leq 8$.

$\centerdot$ $G_{\ell,9}^{(4)}=G_{1,19}^{(4)}$.

$\centerdot$ $G_{\ell,10}^{(4)}=G_{1,20}^{(4)}$.

 $\centerdot$ Pour $11\leq k\leq 15$, le graphe $G_{\ell,k}^{(4)}$ co\"{\i}ncide sur $A\times B$ avec $G_{\ell,k}^{(1)}$.

\vspace{3mm}

%Set $F_i:=\{(e,e'):e=(x,i),e'=(x',i), x< x'\}$, where $i\in 2$.

 \textbf{\underline{Classe $p=5$}:}  Dans ce cas tous les graphes ont le m\^eme ensemble de sommets qui est $A\cup B$ tel que l'un des deux ensembles  $A$ ou bien $B$ forme une cha\^{\i}ne isomorphe \`a $\omega$ l'autre est une anticha\^{\i}ne.

\medskip

 \textbf{I)} Si $\ell=1,2$ nous avons $1\leq k\leq 22$,

 $\centerdot$ Si $1\leq k\leq 12$, le graphe $G_{\ell,k}^{(5)}$ est tel que $A$ forme une cha\^{\i}ne, $B$ une anticha\^{\i}ne, le reste des arcs \'etant le m\^eme que dans $G_{\ell,k}^{(1)}$. % the set $F_0$ of edges.
\smallskip

 $\centerdot$ Si $13\leq k\leq 21$, le graphe $G_{\ell,k}^{(5)}$ est tel que $B$ forme une cha\^{\i}ne, $A$ une anticha\^{\i}ne, le reste des arcs \'etant le m\^eme que dans $G_{\ell,k-12}^{(1)}$.
\smallskip

 $\centerdot$ Dans $G_{\ell,22}^{(5)}$, $A$ est ordonn\'e lin\'eairement comme $\omega$ et $B$ est une anticha\^{\i}ne, il n'y a pas d'autre arc. Ainsi $G_{\ell,22}^{(5)}$ est isomorphe \`a $\omega\oplus \Delta_{\NN}$.
\medskip

\begin{figure}[!hbp]
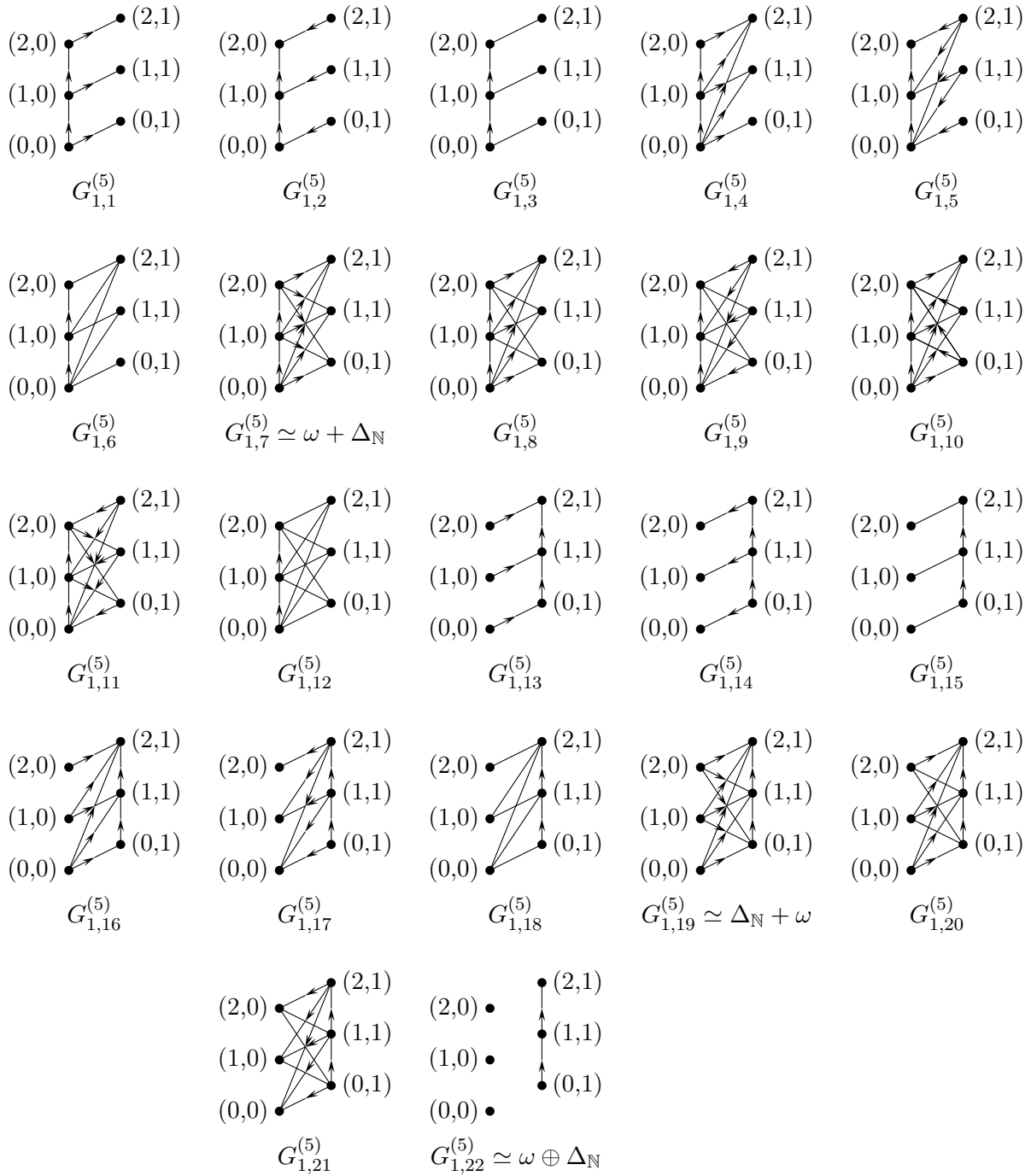

\begin{center}
\small
\begin{tabular}[l]{ccccccccc}
\input{graph1-class5}&&\input{graph2-class5}&&\input{graph3-class5}&&\input{graph4-class5}&&\input{graph5-class5}\\
\input{graph6-class5}&&\input{graph7-class5}&&\input{graph8-class5}&&\input{graph9-class5}&&\input{graph10-class5}\\
\input{graph11-class5}&&\input{graph12-class5}&&\input{graph13-class5}&&\input{graph14-class5}&&\input{graph15-class5}\\
\input{graph16-class5}&&\input{graph17-class5}&&\input{graph18-class5}&&\input{graph19-class5}&&\input{graph20-class5}\\
&&\input{graph21-class5}&&\input{graph22-class5}&&&&\\
\end{tabular}
\caption{Les graphes minimaux de classe $p=5$ pour $\ell=1$.}
\label{repre:graphe-classe5}
\end{center}
\end{figure}

Les graphe $G_{1,k}^{(5)},~1\leq k\leq 22$ sont repr\'esent\'es dans la \figurename ~\ref{repre:graphe-classe5}.

\vspace{3mm}

\textbf{II)} Si $\ell=3,4$ il n'y a aucun graphe. %the linear order is isomorphic to an order of the following set $\{ \omega+\omega, \omega^*+\omega, \omega+\omega^*, \omega^*+\omega^*\}$, we have the same examples which are

\vspace{2mm}

 \textbf{III)} Si $5\leq \ell\leq 8$, nous avons les m\^emes exemples pour chaque valeur de $\ell$ et leur nombre est  $24$.

 $\centerdot$ $G_{\ell,k}^{(5)}=G_{1,k}^{(5)}$ pour tout $1\leq k\leq 6$.

 $\centerdot$ $G_{\ell,k}^{(5)}=G_{1,k+1}^{(5)}$ pour tout $7\leq k\leq 10$.

  $\centerdot$ $G_{\ell,k}^{(5)}=G_{1,k+2}^{(5)}$ pour tout $11\leq k\leq 16$.

   $\centerdot$ $G_{\ell,k}^{(5)}=G_{1,k+3}^{(5)}$ pour tout $17\leq k\leq 18$.

 $\centerdot$ Les graphes $G_{\ell,19}^{(5)}$, $G_{\ell,20}^{(5)}$ et $G_{\ell,21}^{(5)}$ co\"{\i}ncident sur $A\times B$ avec respectivement $G_{\ell,13}^{(1)}$, $G_{\ell,14}^{(1)}$ et $G_{\ell,15}^{(1)}$ tels que $A$ est ordonn\'e comme $\omega$ et $B$ est une anticha\^{\i}ne.

 $\centerdot$ Les graphes $G_{\ell,22}^{(5)}$, $G_{\ell,23}^{(5)}$ et $G_{\ell,24}^{(5)}$ co\"{\i}ncident sur $A\times B$ avec respectivement $G_{\ell,13}^{(1)}$, $G_{\ell,14}^{(1)}$ et $G_{\ell,15}^{(1)}$ tel que  $B$ est ordonn\'e comme $\omega$ et $A$ est une anticha\^{\i}ne.

 \vspace{2mm}

 \vspace{3mm}

\textbf{\underline{Classe $p=6$}:} Dans ce cas tous les graphes ont le m\^eme ensemble de sommets qui est $A\cup B$ tel que l'un des deux ensembles  $A$ ou  $B$ est une cha\^{\i}ne isomorphe \`a $\omega^*$ l'autre \'etant une anticha\^{\i}ne. Sur les paires de $A\times B$, les graphes $G_{\ell,k}^{(6)}$  sont obtenus de la m\^eme façon que dans le cas $p=5$. % with the only difference that we add the set $F'_i$ instead of the set $F_i$.
Tout graphe de cette classe est le dual d'un des graphes de la classe pr\'ec\'edente. %Every graph in this case is the dual of one graph of the previous case.

\vspace{3mm}

\textbf{\underline{Classes $7\leq p\leq 10$}:}  Les graphes $G_{\ell,k}^{(p)}$  sont obtenus de la m\^eme façon que dans le cas $p=6$ avec la diff\'erence que  $A$ ou $B$ est une  clique r\'eflexive au lieu d'une cha\^{\i}ne si $p=7$, $A$ et $B$ sont tous les deux ordonn\'es lin\'eairement, l'un comme $\omega$ et l'autre comme $\omega^*$ si $p=8$, l'un des ensembles $A$ ou $B$ est ordonn\'e lin\'eairement comme $\omega$ et l'autre est une clique r\'eflexive si $p=9$ et l'un des deux ensembles $A$ ou $B$ est ordonn\'e lin\'eairement comme $\omega^*$ et l'autre est une clique r\'eflexive dans le cas $p=10$. %In this
\medskip

Nous obtenons, parmi tous les graphes cit\'es, les vingt bicha\^{\i}nes $\mathcal B:=(V,\leq,\leq')$ de Monteil et Pouzet \cite{mont-pou}. %with $\leq$ isomorphic to $\omega$ or $\omega^*$ (in the eight remaining bichains, the order $\leq$ is of order type $\alpha+\beta$ with $\alpha, \beta \in \{\omega, \omega^{*}\}$).

%If $\rho$ is irreflexive, we obtain the same examples without the loops.
\bigskip

Nous avons le r\'esultat suivant.
\begin{lemma}
Aucun des graphes $\mathcal G_{\ell,k}^{(p)}$ ne s'abrite dans un autre.
\end{lemma}
\begin{proof}
Il suffit de remarquer que pour tout triplet $(p,\ell,k)$ la restriction de $\mathcal G_{\ell,k}^{(p)}$ \`a l'ensemble $\{(0,0),(1,0),(0,1),(1,1)\}$, si  $\mathcal G_{\ell,k}^{(p)}\in\mathfrak A_1\cup\mathfrak B_2$ et \`a l'ensemble $\{a,(0,0),(1,0),(0,1),(1,1)\}$, si $\mathcal G_{\ell,k}^{(p)}\in\mathfrak B_1$,  ne s'abrite dans aucun autre graphe $\mathcal G_{\ell',k'}^{(p')}$ pour $(p',\ell',k')\neq (p,\ell,k)$.
\end{proof}

\subsubsection{Profils des \'el\'ements de $\mathfrak A$}

Pour toutes valeurs des entiers $p,\ell$ et $k$, d\'esignons par $\varphi_{\ell,k}^{(p)}$ le profil du graphe $\mathcal G_{\ell,k}^{(p)}$.
\begin{proposition}\label{prop:profil-ordonne}
Le profil d'un \'el\'ement de $\mathfrak A$ est ou bien donn\'e par l'une des cinq fonctions:  $\varphi_1(n):=2^n-1$, $\varphi_2(n):=2^n-n$,  $\varphi_3(n):=2^{n-1}$, $\varphi_4(n):=2^{n-1}+1$ et la fonction de Fibonacci ou bien  born\'e inf\'erieurement par l'une d'entre elles.
\end{proposition}

\vspace{2mm}

La preuve de cette proposition d\'ecoule des lemmes suivants.

\begin{lemma}\label{lemm-profil1}
 Les graphes ordonn\'es $\mathcal G_{\ell,k}^{(p)}$ pour $1\leq\ell\leq 2$ et ($p=1$ et $1\leq k\leq 3$) ou bien ($2\leq p\leq 4$ et $10\leq k\leq 12$) ont un profil qui a la croissance de la fonction de Fibonacci.
\end{lemma}
\begin{proof}
Soit $\mathcal G_{\ell,k}^{(p)}=(V_{\ell,k}^{(p)},\leq_{\ell,k}^{(p)},\rho_{\ell,k}^{(p)})$ un des graphes cit\'es dans le Lemme \ref{lemm-profil1}. Comme $\ell\in\{1, 2\}$, alors $\leq_{\ell,k}^{(p)}$ est isomorphe \`a $\omega$ ou $\omega^{\star}$. D'apr\`es la description des graphes faite pr\'ec\'edemment, les restrictions de $\mathcal G_{\ell,k}^{(p)}$ aux paires de sommets diff\'erentes de $\{(n,0),(n,1)\}$  sont toutes isomorphes. %D\'esignons par $\varphi_{\ell,k}^{(p)}$ le profil de $\mathcal G_{\ell,k}^{(p)}$.
Pour calculer $\varphi_{\ell,k}^{(p)}(r)$ pour $r\in\mathbb N$,  soit $r$ sommets distincts ordonn\'es suivant $\leq_{\ell,k}^{(p)}$. Alors, ou bien cet ordre se termine par une paire $\{(n,0),(n,1)\}$ avec $n\in\mathbb N$ et dans ce cas le nombre de sous-graphes non isomorphes d'ordre $r$ est $\varphi_{\ell,k}^{(p)}(r-2)$, ou bien il ne se termine pas par une telle paire et dans ce cas le nombre de sous-graphes d'ordre $r$ est $\varphi_{\ell,k}^{(p)}(r-1)$. On obtient alors

$$\left\{\begin{array}{l}
\varphi_{\ell,k}^{(p)}(0)=\varphi_{\ell,k}^{(p)}(1)=1\\
\varphi_{\ell,k}^{(p)}(r)=\varphi_{\ell,k}^{(p)}(r-2)+\varphi_{\ell,k}^{(p)}(r-1)\;\text{ pour }r\geq 2
\end{array}
\right.$$
\end{proof}

\begin{lemma}\label{lemme:profilexpo}
Le profil du graphe ordonn\'e $\mathcal G_{\ell,k}^{(p)}$ pour $1\leq\ell\leq 2$ avec $(p=1\text{ et } 7\leq k\leq 12)$ ou bien $(p=2\text{ et } k\in\{5,6,9\})$ ou bien $(p=3\text{ et } k\in\{4,6,8\})$ ou bien $(p=4\text{ et }k\in\{4,5,7\})$ est donn\'e par $\varphi_{\ell,k}^{(p)}(r)=2^{r}-1,~~r\geq 1$.
\end{lemma}
\begin{proof}
Dans ces cas, nous pouvons coder chaque sous-graphe d'ordre $r$ par un mot de longueur $r$ sur un alphabet \`a deux lettres $\{0,1\}$. Soit $r$ sommets distincts ordonn\'es suivant $\leq_{\ell,k}^{(p)}$. A chaque sommet dans cet ordre, nous associons $0$ s'il est dans $\mathbb N\times\{0\}$ et $1$ s'il est dans $\mathbb N\times\{1\}$, les lettres \'etant plac\'ees de gauche \`a droite. Sachant que les mots compos\'es d'une seule lettre donnent des sous-graphes isomorphes. Le nombre de sous-graphes non isomorphes d'ordre $r$ est donc le nombre de mots diff\'erents de longueur $r$ moins $1$, d'o\`u le r\'esultat.
\end{proof}

\begin{lemma}
Le profil du graphe ordonn\'e $\mathcal G_{\ell,k}^{(p)}$ pour $1\leq\ell\leq 2$ avec $(p=1, ~k\in\{4,5,6\})$ ou bien $(p=2,~k\in\{4,7,8\})$ ou $(p=3,~ k\in\{5,7,9\})$ ou bien $(p=4,~ k\in\{6,8,9\})$ est donn\'e par: %$\varphi_{\ell,k}^{(p)}(r)=2^{r}-r,~~r\geq 3$, les premi\`ere valeur \'etant $1$, $1$, $2$ pour $r=0$, $1$ et $2$.
$$\left\{\begin{array}{l}
\varphi_{\ell,k}^{(p)}(0)=\varphi_{\ell,k}^{(p)}(1)=1;~\varphi_{\ell,k}^{(p)}(2)=2,\\
\varphi_{\ell,k}^{(p)}(r)=2^{r}-r,~~r\geq 3.
\end{array}
\right.$$
\end{lemma}

\begin{proof}
Dans ces cas, nous pouvons coder chaque sous-graphe d'ordre $r$ par un mot de longueur $r$ sur un alphabet \`a deux lettres $\{0,1\}$ comme dans la preuve du Lemme \ref{lemme:profilexpo}.
Par exemple, pour $p=1$, nous associons \`a chaque sommet de $\mathbb N\times\{0\}$ la lettre $0$ et \`a chaque sommet de $\mathbb N\times\{1\}$ la lettre $1$. Si $p=2$ nous faisons l'inverse. Alors tous les mots de m\^eme longueur $r$ de la forme $\underset{q}{\underbrace{1\dots1}}\underset{r-q}{\underbrace{0\dots 0}}$ donnent des sous-graphes qui sont isomorphes entre eux pour toute valeur de $q~(0\leq q\leq r)$.  Le nombre de sous-graphes non isomorphes d'ordre $r$ est donc le nombre de mots diff\'erents de longueur $r$ moins $r$, d'o\`u le r\'esultat.
\end{proof}

\begin{lemma}\label{lemm-profil2}
Le profil du graphe ordonn\'e $\mathcal G_{\ell,k}^{(p)}$ pour $1\leq\ell\leq 2$  avec $(p=1,~ k\in\{13,14,15\})$ ou bien $(2\leq p\leq 3,~ k\in\{13,16,17,19,20,21\})$ ou bien $(p=4,~ k\in\{15,17,18,19,20,21\})$ est donn\'e par $\varphi_{\ell,k}^{(p)}(r)=2^{r-1},~~r\geq 1$.
\end{lemma}

\begin{proof}
Pour $k\neq 19, 20, 21$, tous les autres graphes donn\'es dans ce lemme ont pour ensemble de sommets $\{a\}\cup\mathbb N\times\{0,1\}$. Ils ont pour particularit\'e d'\^etre monomorphes sur $\mathbb N\times\{0,1\}$. Nous pouvons coder chaque sous-graphe d'ordre $r$ par un mot sur un alphabet \`a trois lettres $\{a,0,1\}$, en associant, de la m\^eme mani\`ere que pr\'ec\'edemment, $0$ \`a tout sommet de $\mathbb N\times\{0\}$ et $1$ \`a tout sommet de $\mathbb N\times\{1\}$ et en rajoutant $a$ en d\'ebut du mot si le sous-graphe contient le sommet $a$. Sachant que pour tout $r\in \mathbb N$, tous les mots de longueur $r$ form\'es uniquement des lettres $0$ et $1$ correspondent \`a des sous-graphes isomorphes entre eux et isomorphes au mot de longueur $r$ commençant par $a$ et form\'e d'une seule lettre ($0$ pour certains graphes comme par exemple pour $\mathcal G_{\ell,13}^{(1)}$) et $1$ pour d'autres comme par exemple pour $\mathcal G_{\ell,15}^{(4)}$). % correspondent \`a des sous-graphes d'ordre $r$ isomorphes \`a ceux  ne contenant pas le sommet $a$.
Donc le nombre de sous-graphes non isomorphes d'ordre $r$ est le nombre de mots de longueur $r$ commençant par $a$, ce qui donne $2^{r-1}$.
\smallskip

\noindent Pour $k\in\{ 19, 20, 21\}$, les graphes ont pour ensemble de sommets $\mathbb N\times\{0,1\}$. Ils sont tels que tout sous-graphe ayant $q$ sommets de $\mathbb N\times\{0\}$ et $r-q$ sommets de $\mathbb N\times\{1\}$ est isomorphe \`a un des sous-graphes ayant $r-q$ sommets de $\mathbb N\times\{0\}$ et $q$ sommets de $\mathbb N\times\{1\}$. D'o\`u le r\'esultat. %le nombre de tous les sous-graphes est la moiti\'e du nombre de mots \`a deux lettres.
\end{proof}

\begin{lemma}
Le profil du graphe ordonn\'e $\mathcal G_{\ell,k}^{(p)}$ pour $1\leq\ell\leq 2$ et $(p=1,~ k\in\{16,17,18\})$ ou bien $(2\leq p\leq 3,~ k\in\{14,15,18\})$ ou bien $(p=4,~ k\in\{13,14,16\})$ est donn\'e par $\varphi_{\ell,k}^{(p)}(r)=2^{r-1}+1,~~r\geq 2$.
\end{lemma}

\begin{proof}
Tous les graphes donn\'es dans ce lemme ont pour ensemble de sommets $\{a\}\cup\mathbb N\times\{0,1\}$. Ils ont pour particularit\'e d'\^etre monomorphes sur $\mathbb N\times\{0,1\}$. Nous pouvons, comme pour la preuve du Lemme \ref{lemm-profil2} pr\'ec\'edent, coder chaque sous-graphe d'ordre $r$ par un mot sur un alphabet \`a trois lettres $\{a,0,1\}$, avec la diff\'erence qu'ici, les mots de longueur $r$ commençant par $a$ et form\'es d'une seule lettre ne donnent pas des sous-graphes isomorphes \`a ceux ne contenant pas le sommet $a$. Donc le nombre de sous-graphes non isomorphes d'ordre $r$ est le nombre de mots de longueur $r$ commençant par $a$ auquel on rajoute un, ce qui donne $2^{r-1}+1$.
\end{proof}

\begin{lemma}
$\varphi_{\ell,k}^{(p)}(r)\geq 2^{r},~~r\geq 3$ pour $1\leq\ell\leq 2$ et $(2\leq p\leq 4,~k\in\{1,2,3\})$, les premi\`eres valeurs \'etant $1$, $1$, $3$ pour $r=0$, $1$ et $2$.
\end{lemma}

\begin{proof}
A tout sous-graphe \`a $r$ sommets, ordonn\'es suivant $\leq_{\ell,k}^{(p)}$, nous pouvons associer un mot de longueur $r$ en $0, 1$. Tous les mots diff\'erents  donnent des sous-graphes non isomorphes (\`a l'exception des deux mots form\'es par une seule lettre qui donnent deux sous-graphes isomorphes), donc $\varphi_{\ell,k}^{(p)}(r)\geq 2^{r}-1$. Pour $r\geq 3$, chaque mot qui contient le facteur %\footnote{Un facteur d'une suite de lettres $a_1a_2\dots a_n$ est une sous-suite de la forme $a_ia_{i+1}\dots a_j$ pour $1\leq i<j\leq n$.}
  $01$ qui correspond \`a deux sommets $(n,0)$ et $(m,1)$ successifs par rapport \`a $\leq_{\ell,k}^{(p)}$, code deux sous-graphes non isomorphes, car nous avons deux cas, le cas o\`u $m=n$ et le cas o\`u $m>n$ qui donnent des sous-graphes non isomorphes. Il s'ensuit qu'\`a partir de $r\geq 3$ nous avons $\varphi_{\ell,k}^{(p)}(r)\geq 2^{r}$.
\end{proof}

\begin{lemma}
Pour $\ell\geq 3$ ou $p\geq 5$ le profil du graphe $\mathcal G_{\ell,k}^{(p)}$ est sup\'erieur ou \'egal \`a l'une des cinq fonctions: $\varphi_1(n):=2^n-1$, $\varphi_2(n):=2^n-n$,  $\varphi_3(n):=2^{n-1}$, $\varphi_4(n):=2^{n-1}+1$ et la fonction de Fibonacci.
\end{lemma}

\begin{proof}
Tous les graphes pour $p\geq 5$ se d\'eduisent des graphes pour $p\leq 4$ avec les restrictions aux ensembles $\mathbb N\times\{0\}$ et $\mathbb N\times\{1\}$ qui ne sont pas isomorphes, donc le nombre de sous-graphes d'ordre $r$ sera plus grand que dans le cas $p\leq 4$ dont les profils sont donn\'es par l'une des cinq fonctions ci-dessus. Pour $\ell\geq 3$, les graphes s'obtiennent \`a partir de ceux pour lesquels $\ell\leq 2$, soit en permutant $\leq_{\ell,k}^{(p)}$ et $\rho_{\ell,k}^{(p)}$ comme c'est le cas pour $\ell=3,4$ ou en modifiant l'ordre $\leq_{\ell,k}^{(p)}$ (pour $5\leq\ell\leq 8$), les raisonnement adopt\'es dans les preuves des lemmes pr\'ec\'edents restent valables.
\end{proof}

\vspace{4mm}

\textbf{Preuve du Th\'eor\`eme \ref{theo:dichotomie}.} Soit $\mathcal R:=(E,\leq,\rho_1,\dots,\rho_k)$ une structure binaire ordonn\'ee de type $k$. Si $\mathcal R$ a une d\'ecomposition monomorphe finie alors, d'apr\`es le Th\'eor\`eme \ref{compactness}, il existe un entier $\ell$ tel que tout \'el\'ement de $\mathcal {A(R)}$, l'\^age de $\mathcal R$, poss\`ede une d\'ecomposition monomorphe ayant au plus $\ell$ blocs. D'apr\`es le Th\'eor\`eme \ref{thm: polynomial-interval}, le profil de $\mathcal {A(R)}$, et donc celui de $\mathcal R$, est un polyn\^ome. Si $\mathcal R$ n'a pas de d\'ecomposition monomorphe finie, alors d'apr\`es la Proposition \ref{prop:decomp-finie}, il existe $1\leq i\leq k$ tel que $\mathcal R_i:=(E, \leq, \rho_i)$ n'a pas de d\'ecomposition monomorphe finie, donc $\mathcal R_i\in\mathscr D$ et d'apr\`es le Th\'eor\`eme \ref{thm:graph-ordonne}, $\mathcal R_i$ abrite un \'el\'ement de $\mathfrak A$. Les \'el\'ements de $\mathfrak A$, d'apr\`es la Proposition \ref{prop:profil-ordonne}, ont un profil au moins exponentiel. Donc $\mathcal R_i$ a un profil au moins exponentiel ce qui entra\^{i}ne le m\^eme r\'esultat pour $\mathcal R$.
\hfill  $\Box$ 
\clearemptydoublepage
\clearpage
\addcontentsline{toc}{chapter}{Conclusion}
\chapter*{Conclusion g\'en\'erale}
\markboth{\slshape{Conclusion g\'en\'erale}} {\slshape{conclusion g\'en\'erale}}

Dans ce travail, trois notions importantes ont \'et\'e abord\'ees: l'algébricité de la série \'enumérative d'une classe de structures finies, la notion de minimalit\'e et le phénomène de saut dans le comportement des profils. Chacune des trois notions a constitu\'e une partie de ce document. 
Cette conclusion, est consacr\'ee \`a l'expos\'e de quelques perspectives qui feraient, en plus des questions posées et des conjectures formulées dans les diff\'erentes parties,  directement suite \`a ce travail.\\

\vspace{2mm}

Au chapitre \ref{chap:permut-bichaine} de la première partie nous avons explicité le lien existant entre les permutations et les bichaînes, ce qui justifie que l'on considère les permutations comme des structures relationnelles. A ce titre, nous avons généralisé des résultats sur les permutations à des structures binaires ordonnées réflexives, en particulier le th\'eor\`eme d'Albert-Atkinson que nous avons g\'en\'eralis\'e dans le Th\'eor\`eme \ref{theo:algebraic} du chapitre \ref{sect:str.rela.bin}.

\vspace{1mm}

Certes, la littérature sur les permutations est vaste et comporte \'enormément de résultats et d'outils; il serait intéressant de les traduire en terme de la théorie des relations, notamment la notion de classe grille-g\'eom\'etrique (voir Annexe \ref{chapitre:grille-geo}) qui est, en vu des Th\'eor\`emes \ref{theorem:annexe1} et \ref{theoreme1}, un outil puissant pour montrer le caractère h\'er\'editairement rationnel et h\'er\'editairement alg\'ebrique d'une classe de permutations.

\vspace{1mm}

Nous avons ensuite illustr\'e, dans le chapitre \ref{chap:exemple-conjecture}, cette g\'en\'eralisation par la construction, en guise d'exemple, d'une classe de structures binaires ordonn\'ees dont les \'el\'ements ind\'ecomposables sont de taille au plus deux. Dans le cas où les structures de cette classe sont form\'ees d'une seule relation en plus de l'ordre lin\'eaire, la classe est caract\'eris\'ee par des bornes dont la taille est \'egale \`a trois ou quatre. Il serait int\'eressant de donner une caract\'erisation de cette classe dans le cas g\'en\'eral.

\vspace{1mm}

Nous avons \'egalement \'etabli une bijection entre la classe des structures binaires ordonn\'ees s\'eparables r\'eflexive de type $1$ et celle des $3$-permutations s\'eparables en passant par la classe des arbres binaires \'etiquet\'es par un ensemble \`a quatre \'el\'ements. Il serait int\'eressant de trouver une relation "plus combinatoire" entre les objets de ces deux classes. %une structure binaire s\'eparable r\'eflexive de type $1$ et une $3$-permutation s\'eparable.

\vspace{1mm}

Nous pensons que le Th\'eor\`eme \ref{theo:algebraic} peut encore \^etre g\'en\'eralis\'e et, \`a cet effet, nous avons propos\'e la Conjecture \ref{conjec} (en page \pageref{conjec}). Cette conjecture est v\'erifi\'ee par au moins une classe de structures, la cl\^oture par sommes de la classe des bicha\^{i}nes critiques de Schmerl et Trotter. Il serait int\'eressant de donner une r\'eponse, positive ou n\'egative, \`a cette conjecture ou au moins d'\'elargir l'ensemble d'exemples qui la v\'erifient.

\vspace{2mm}

Dans la deuxi\`eme partie, nous nous sommes int\'eress\'es \`a la notion de minimalit\'e. Nous avons montr\'e que les classes ind-minimales sont les \^ages de structures ind\'ecomposables infinies, qu'elles sont belordonn\'ees mais pas n\'ecessairement h\'er\'editairement belordonn\'ees et qu'elles
sont en nombre contin\^upotent. Il serait int\'eressant d'\'etudier d'autres propri\'et\'es, par exemple:
\begin{itemize}
\item Que doit v\'erifier une relation d'\^age ind-minimal pour que cet \^age soit h\'er\'editairement belordonn\'e?
%\item Quelle est la relation entre le noyau d'une structure relationnelle et le caract\`ere ind-minimal de son \^age? Si ce noyau n'est pas vide, est-il fini? Son cardinal est-il born\'e?
\item Il existe un nombre contin\^upotent de classes form\'ees de graphes dirig\'es sans boucles ou non dirig\'es avec des boucles. Qu'en est-il des classes de graphes non dirig\'es sans boucles?
\end{itemize}

\vspace{4mm}

Enfin, la troisi\`eme partie a port\'e sur le ph\'enomène de sauts observé dans le comportement des profils des classes héréditaires de structures finies. Nous avons montr\'e que pour des classes de structures binaires ordonn\'ees, la croissance du profil passe d'une croissance polynomiale \`a une croissance exponentielle. Nous pensons que ce r\'esultat reste vrai dans le cas des structures ordonn\'ees non n\'ecessairement binaires. %, il serait intéressant de le montrer.
La preuve de notre r\'esultat  s'appuie sur  la notion de d\'ecomposition monomorphe d'une structure relationnelle et le th\'eoreme de Ramsey.

\vspace{2mm}

Nous avons d\'ecrit, dans le cas de deux relations binaires dont un ordre lin\'eaire, une liste de structures minimales parmi celles dont le profil n'est pas polynomialement born\'e et nous avons montr\'e que les profils des structures de cette liste sont exponentiels. La cardinalit\'e de cette liste \'etant très grande, il serait int\'eressant de confirmer ce r\'esultat par un algorithme et un programme informatique.

\vspace{2mm}

Nous avons \'egalement d\'ecrit, dans le cas des graphes non dirig\'es, une liste form\'ee de dix graphes minimaux parmi ceux qui n'ont pas une d\'ecomposition monomorphe finie. Parmi ces graphes, certains ont un profil polynomialement borné. Ceci montre en particulier qu'un graphe, et donc par suite une structure binaire non ordonn\'ee, ayant un profil polynomialement born\'e n'a pas forc\'ement une d\'ecomposition monomorphe finie. Nous savons qu'un graphe (non dirig\'e sans boucle) a un profil born\'e par un polyn\^ome si et seulement si il est cellulaire au sens de Schmerl \cite{pouzet06}. Nous mentionnons \`a la section \ref{sec:cellulaire} une extension de la  notion de Schmerl susceptible de caract\'eriser les structures relationnelles a profil born\'e par un polyn\^ome.

\clearemptydoublepage
\appendix
\addcontentsline{toc}{chapter}{Annexe}%\appendix
\chapter[Codage des structures s\'eparables par des arbres binaires] {Codage des structures relationnelles binaires s\'eparables de type $1$ par des arbres binaires \'etiquet\'es}\label{codage}

\section{Introduction}

Dans la Proposition \ref{prop:generatriceseparable} de la section \ref{subsection:classe S1} (page \pageref{prop:generatriceseparable}), nous avons d\'etermin\'e la fonction g\'en\'eratrice de $\mathscr S^{re}_1$, la classe des structures binaires ordonn\'ees s\'eparables r\'eflexives de type $1$. Il s'av\`ere (voir \cite{Asin-Mans}) que cette fonction est celle des $3$-permutations s\'eparables, elle est \'egalement la m\^eme que celle des partitions guillotine en dimension $4$ (voir \cite{Ack-Ba-Pin-Rom}).

\vspace{1mm}

En effet, dans \cite{Ack-Ba-Pin-Rom}, les auteurs ont \'etudi\'e les partitions guillotine en dimension $d$ et ont donn\'e la s\'erie g\'en\'eratrice du nombre de ces partitions, elle v\'erifie
\begin{equation}
f=1+xf+(d-1)xf^2.\label{eq:part-guill}
\end{equation}

 Ils ont \'egalement \'etabli une correspondance bijective entre les partitions guillotine en dimension $d$ et les arbres binaires \'etiquet\'es par l'ensemble $\{1,\cdots,d\}$.

\vspace{2mm}

Dans \cite{Asin-Mans}, les auteurs ont donn\'e la s\'erie g\'en\'eratrice des $d$-permutations s\'eparables, elle v\'erifie
\begin{equation}
f=1+xf+(2^{d-1}-1)xf^2.\label{eq:d-permutation}
 \end{equation}
  Il ont \'egalement \'etabli une correspondance bijective entre les $d$-permutations s\'eparables et les partitions guillotine en dimension $2^{d-1}$. Il se trouve que notre s\'erie $S^{r}_1$ v\'erifie la relation \eqref{eq:part-guill} pour $d=4$ et la relation \eqref{eq:d-permutation} pour $d=3$.

\vspace{2mm}

  Dans cette partie, nous construisons  une bijection entre les structures binaires ordonn\'ees s\'eparables de type $1$ et les arbres binaires \'etiquet\'es par l'ensemble $\{1,2,3,4\}$, ceci entraine que l'ensemble $\mathscr S^{re}_1$ et l'ensemble des partitions guillotine en dimension $4$ sont isomorphes. Cela entraine \'egalement que $\mathscr S^{re}_1$ est isomorphe \`a l'ensemble des $3-$permutations s\'eparables.

\section{Structures binaires ordonn\'ees s\'eparables et arbres binaires}

Un \emph{arbre binaire}\index{arbre!binaire} est un graphe dirig\'e connexe sans cycle tel que le degr\'e de chaque sommet (n{\oe}ud) est au plus $3$. C'est un arbre avec une racine dans lequel tout sommet poss\`ede au plus deux successeurs. La \emph{racine}\index{arbre!racine d'un -} d'un arbre binaire est le sommet de degr\'e au plus $2$ qui ne poss\`ede aucun pr\'ed\'ecesseur. Avec cette racine, chaque autre n{\oe}ud a un unique pr\'ed\'ecesseur, appel\'e \emph{parent} et chaque n{\oe}ud, diff\'erent d'un sommet pendant, poss\`ede au plus deux successeurs, les \emph{fils}. Les deux fils sont souvent d\'esign\'es par \emph{fils droit} et \emph{fils gauche}. Un arbre \emph{\'etiquet\'e} est un arbre dans lequel chaque sommet poss\`ede une \emph{\'etiquette}, appel\'ee aussi \emph{couleur}. Ceci d\'efinit une application entre les sommets de l'arbre et l'ensemble des \'etiquettes.

\vspace{2mm}

Rappelons que $\mathscr S^{re}_1$ est l'ensemble des structures relationnelles binaires ordonn\'ees s\'eparables r\'eflexives de type $1$ et que $\mathscr D^1_{(2)}$, l'ensemble des structures ordonn\'ees de type $1$ d\'efinies sur des ensembles \`a deux \'el\'ements,  poss\`ede, \`a l'isomorphie pr\`es, quatre structures. Ecrivons chaque $r$ de $\mathscr D^1_{(2)}$ comme $r:=(\{0,1\},\leq, \rho)$ avec $0<1$, et
posons $\mathscr D^1_{(2)}=\{r_1,r_2,r_3,r_4\}$.
\vspace{3mm}

A une structure de $\mathscr S^{re}_1$, nous associons un arbre binaire dont tous les sommets sont \'etiquet\'es par l'ensemble $\{1,2,3,4\}$ en imposant la restriction suivante: le fils gauche n'a pas la m\^eme couleur que son parent. Cet arbre est construit de mani\`ere r\'ecursive comme suit.\\
L'arbre vide correspond à la structure \`a un \'el\'ement et l'arbre \'etiquet\'e \`a un sommet correspond \`a la structure (consid\'er\'ee comme op\'erateur) \`a deux \'el\'ements (il y a quatre couleurs).
\vspace{3mm}

Soit $\mathcal R\in\mathscr S^{re}_1$ une structure \`a $n$ \'el\'ements. Donc $\mathcal R$ est la somme lexicographique de deux structures $\mathcal R_1$ et $\mathcal R_2$, ayant respectivement $n_1$ et $n_2$ \'el\'ements avec $n_1+n_2=n$, ind\'ex\'ee par une structure $r_i$ de $\mathscr D^1_{(2)}$, c'est \`a dire que $\mathcal R=\mathcal R_1\underset{r_i}\oplus \mathcal R_2$.

 \noindent Nous associons \`a $\mathcal R$ un sommet \'etiquet\'e $i$ d'un arbre \`a $n-1$ sommets, ce sommet sera la racine de l'arbre. Ce sommet a  deux fils, le fils gauche repr\'esente $\mathcal R_1$ et le fils droit repr\'esente $\mathcal R_2$.  Si $\mathcal R_1$ ou $\mathcal R_2$ a un sommet ($n_1=1$ ou $n_2=1$) alors le sommet $i$ de l'arbre (la racine) aura un seul fils, le droit ou le gauche. Comme, par construction, nous avons
\begin{equation}
\mathscr S^{re}_1=\{\bf{1}\}\cup \underset{1\leq i\leq 4}\bigcup (\mathscr S^{re}_1(r_i)\underset{r_i}\oplus \mathscr S^{re}_1),\label{unicité-codage}
\end{equation}
la structure $\mathcal R_1$ ne peut pas se d\'ecomposer en une somme lexicographique ind\'ex\'ee par $r_i$. Donc le sommet qui la repr\'esente (c'est \`a dire le fils gauche du sommet racine) ne peut pas avoir $i$ pour \'etiquette. Nous continuons le m\^eme processus pour construire un arbre binaire (puisque chaque structure est $2$-d\'ecomposable) \`a $n-1$ sommets dont les sommets pendants ont pour couleurs les indices des structures \`a deux \'el\'ements, qui s'abritent dans $\mathcal R$ et qui interviennent dans la somme lexicographique.

\vspace{3mm}

Soient maintenant deux structures ordonn\'ees s\'eparables $\mathcal S_1$ et $\mathcal S_2$ de $\mathscr S^{re}_1$ et soient $T_1$ et $T_2$ les arbres binaires qui leurs sont associ\'es respectivement. Si $\mathcal S_1$ est $r_i$-ind\'ecomposable et $\mathcal S=\mathcal S_1\underset{r_i}\oplus \mathcal S_2$ alors l'arbre $T$ correspondant \`a $\mathcal S$ est obtenu \`a partir des deux arbres $T_1$ et $T_2$ de la mani\`ere suivante: un sommet de couleur $i$ est associ\'e \`a $\mathcal S$, ce sommet qui est la racine de $T$) aura pour fils gauche la racine de $T_1$ et pour fils droit la racine de $T_2$ (voir \figurename ~\ref{arbre}). La racine de $T_1$ n'a pas l'\'etiquette $i$ puisque $\mathcal S_1$ est $r_i$-ind\'ecomposable, donc l'arbre $T$ code bien la structure $\mathcal S$.

\begin{figure}
\begin{center}
\psset{unit=1cm}
\begin{pspicture}(-3,-2)(3,2)
\psdots[dotsize=5pt](0,1)(-1,0.5)(1,0.5)
\psline[linewidth=0.3pt](0,1)(-1,0.5)
\psline[linewidth=0.3pt](0,1)(1,0.5)
\uput{0.3}[ur](0,1){$i$}
\psellipse(-1,-0.5)(0.5,1)
\psellipse(1,-0.5)(0.5,1)
\rput(-1,-0.5){$T_1$}
\rput(1,-0.5){$T_2$}
\uput{0.3}[d](0,-1.5){$T$}
\end{pspicture}
\caption{Construction de l'arbre de $\mathcal S=\mathcal S_1\underset{r_i}\oplus \mathcal S_2$ \`a partir des arbres de $\mathcal S_1$ et $\mathcal S_2$}
\label{arbre}
\end{center}
\end{figure}

\vspace{3mm}

Il est clair, d'apr\`es la construction des structures s\'eparables, qu'\`a deux structures diff\'erentes correspondent deux arbres diff\'erents construits de la mani\`ere d\'ecrite ci-dessus et qu'avec la restriction impos\'ee sur les \'etiquettes et l'Equation \eqref{unicité-codage}, deux arbres diff\'erents correspondent \`a deux structures diff\'erentes. Donc, il existe une correspondance bijective entre les arbres binaires \`a $n$ sommets, dont les sommets sont \'etiquet\'es par les \'el\'ements de l'ensemble $\{1,2,3,4\}$ avec la restriction que le fils gauche de tout sommet porte une \'etiquette diff\'erente de celle du parent, et les structures binaires ordonn\'ees s\'eparables (r\'efl\'exives ou irr\'eflexives) de type $1$ ayant $n+1$ sommets.

\vspace{3mm}

Il se trouve que cette m\^eme classe d'arbres binaires permet de coder les partitions guillotine d'une boite de dimension $4$ \cite{Ack-Ba-Pin-Rom}, voir aussi \cite{Asin-Mans}.
\vspace{3mm}

Le probl\`eme de partitions guillotine d'une boite de dimension $d$ se d\'efinit comme suit. Etant donn\'ee une boite $B$ de dimension $d$ dans $\mathbb R^d$, une \emph{partition guillotine} de $B$ est une subdivision de $B$ en plus petites boites de $\mathbb R^d$ obtenues en coupant d'abord $B$ en deux boites de $\mathbb R^d$ par un hyperplan parall\`ele \`a ses axes, puis couper les boites  obtenues de la m\^eme mani\`ere (les directions des coupes peuvent changer). Clairement, il existe une infinit\'e de partitions guillotine avec un nombre donn\'e d'hyperplans, mais si l'on s'int\'eresse aux directions des coupes plut\^ot qu'\`a leurs positions exactes, alors le nombre de partitions (structurellement) diff\'erentes obtenues avec $n$ hyperplans de $\mathbb R^d$ est fini.

%%%%%%%%%%%%%%%%%%%%%%%%%%%%%%%%%%%%%%%%%%%%%%%%%%%%%%%%%%%%%%%%%%%%%%%%%%%%%%%%%%%%%%%%%%%%%%%%%%%%%%%%%%%

\clearemptydoublepage
\chapter{Les classes grille-g\'eom\'etriques}\label{chapitre:grille-geo}

Dans cette annexe, nous rappelons la d\'efinition des classes de permutations appel\'ees \emph{classes grille-g\'eom\'etriques}\index{classe!grille-g\'eom\'etrique} pour voir que les permutations exceptionnelles sont contenues dans des classes grille-g\'eom\'etriques. % et justifier l'applications des Th\'eor\`emes \ref{theo:classegrill} et \ref{theo:grillgeoalgebrique}.
Les propri\'et\'es de ces classes nous ont permis, dans la section \ref{sec:conjecture} en page \pageref{theo:classegrill}, de d\'eduire que la cl\^oture par sommes de la classe des bicha\^{i}nes critiques est h\'er\'editairement alg\'ebrique. Pour plus de d\'etails sur ces classes consulter \cite{albert-vatter1, albert-vatter2, huc-vatter, vatter4}.

\vspace{2mm}

Dans l'\'etude des permutations et des classes de permutations, les classes grille-g\'eom\'etriques sont consid\'er\'ees comme l'un des outils les plus puissants. Ces classes sont d\'efinies comme suit.

\vspace{2mm}

Soit $M$ une matrice dont les coefficients appartiennent \`a $\{-1,0,1\}$. Contrairement aux notations conventionnelles, les coefficients de $M$ sont index\'es de gauche \`a droite pour les colonnes et de bas en haut pour les lignes et les indices de ligne et de colonne sont invers\'es. Ainsi, $M_{k,l}$ d\'esigne le coefficient de la $k^{\text{i\`eme}}$ colonne \`a partir de la gauche et de la $l^{\text{i\`eme}}$ ligne \`a partir du bas comme indiqu\'e ci-dessous:

$$\left(\begin{array}{ccc}
(1,2)&(2,2)&(3,2)\\
(1,1)&(2,1)&(3,1)%
\end{array}
\right).$$

A la matrice $M$ est associ\'ee une grille ayant un nombre de blocs \'egal au nombre de coefficients de $M$, dont la r\'epartition en lignes et en colonnes est dict\'ee par $M$.

\vspace{2mm}

La \emph{figure standard} de $M$ est l'ensemble des points de $\mathbb R^2$ pouvant-\^etre plac\'es dans cette grille  de la mani\`ere suivante:

\qquad - sur le segment de droite reliant les points $(k-1,l-1)$ et $(k,l)$ si $M_{k,l}=1$ ou

\qquad - sur le segment de droite reliant les points $(k-1,l)$ et $(k,l-1)$ si $M_{k,l}=-1$.

\vspace{1mm}

\textbf{Exemple:} La figure standard de $M_1=\left(\begin{array}{cc}
1&-1\\
-1&1
\end{array}
\right)$
et celle de $M_2=\left(\begin{array}{cc}
-1&1\\
1&-1
\end{array}
\right)$
sont repr\'esent\'ees sur la  \figurename ~\ref{figstandar}):

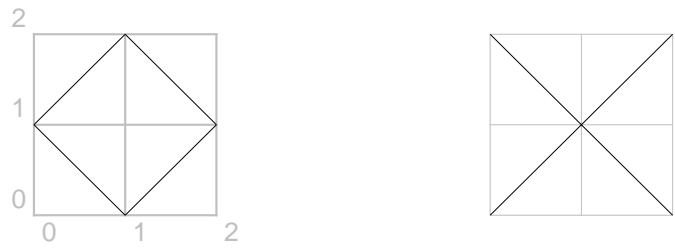
\begin{figure}[t]
\begin{center}
\psset{unit=1.2cm}
\begin{pspicture}(0,-0.5)(7,2)
\psgrid[subgriddiv=0,gridcolor=lightgray,gridlabelcolor=lightgray](0,0)(2,2)
\psline[linewidth=0.3pt](0,1)(1,2)
\psline[linewidth=0.3pt](0,1)(1,0)
\psline[linewidth=0.3pt](1,0)(2,1)
\psline[linewidth=0.3pt](1,2)(2,1)
\psline[linewidth=0.4pt,linecolor=lightgray](5,0)(5,2)
\psline[linewidth=0.4pt,linecolor=lightgray](5,0)(7,0)
\psline[linewidth=0.3pt,linecolor=lightgray](5,2)(7,2)
\psline[linewidth=0.3pt,linecolor=lightgray](7,0)(7,2)
\psline[linewidth=0.3pt,linecolor=lightgray](6,0)(6,2)
\psline[linewidth=0.3pt,linecolor=lightgray](5,1)(7,1)
%\psgrid[subgriddiv=0,gridcolor=lightgray,gridlabel=0](5,0)(7,2)
\psline[linewidth=0.3pt](5,0)(6,1)(7,0)
%\psline[linewidth=0.3pt](6,1)(7,0)
\psline[linewidth=0.3pt](5,2)(6,1)(7,2)
\end{pspicture}
\caption{Figures standards de $M_1$ \`a gauche et $M_2$ \`a droite}
\label{figstandar}
\end{center}
\end{figure}

\vspace{2mm}

\begin{figure}[t]
\begin{center}
\psset{unit=1.2cm}
\begin{pspicture}(0,-0.5)(4,2)
\psgrid[subgriddiv=0,gridcolor=lightgray,gridlabelcolor=lightgray](0,0)(3,2)
\psline[linewidth=0.3pt](0,2)(1,1)
\psline[linewidth=0.3pt](1,1)(2,2)
\psline[linewidth=0.3pt](1,1)(2,0)
\psline[linewidth=0.3pt](3,0)(2,1)(3,2)
\psdots[dotsize=4pt](0.9,1.1)(0.2,1.8)(1.9,1.9)(1.4,0.6)(2.2,1.2)(2.7,1.7)(2.6,0.4)
\uput{0.1}[ur](2.6,0.4){\footnotesize{1}}
\uput{0.1}[dr](2.7,1.7){\footnotesize{5}}
\uput{0.1}[r](2.2,1.2){\footnotesize{4}}
\uput{0.1}[ur](1.4,0.6){\footnotesize{2}}
\uput{0.1}[l](1.9,1.9){\footnotesize{7}}
\uput{0.1}[d](0.2,1.8){\footnotesize{6}}
\uput{0.1}[l](0.9,1.1){\footnotesize{3}}
\end{pspicture}
\caption{Repr\'esentation de la permutation $6327415$ sur la figure standard de $M$.}
\label{figstandarpermut}
\end{center}
\end{figure}
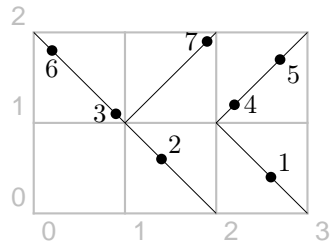

La \emph{classe grille-g\'eom\'etrique} de $M$, not\'ee $Geom(M)$ est l'ensemble de toutes les permutations qui peuvent-\^etre r\'eparties en blocs tels que la sous-suite contenue dans chaque bloc est monotone et peuvent-\^etre repr\'esent\'ees sur la figure standard de $M$. Ces permutations sont d\'efinies de la mani\`ere suivante. Choisir $n$ points sur la figure tels que deux points donn\'es ne soient pas sur une m\^eme ligne horizontale ou verticale. Puis, num\'eroter les points de $1$ \`a $n$ en allant de bas en haut et r\'e\'ecrire ces nombres en allant de gauche \`a droite. Un exemple qui montre que la permutation $6327415$ se trouve dans $Geom(M)$ pour $M=\left(\begin{array}{ccc}
-1&1&1\\
0&-1&-1
\end{array}
\right)$
est donn\'e dans la \figurename ~\ref{figstandarpermut}.

\vspace{2mm}

La notion de classe grille-g\'eom\'etrique permet de faire une description structurelle des classes de permutations. Ces classes ont des propri\'et\'es int\'eressantes.

\begin{theorem}\label{theorem:annexe1}\cite{albert-vatter1}

Toute classe grille-g\'eom\'etrique est h\'er\'editairement rationnelle et a un nombre fini de bornes.
 \end{theorem}

Une classe de permutations est \emph{g\'eom\'etriquement en grille}\index{classe!g\'eom\'etriquement en grille} si elle est contenue dans une classe grille-g\'eom\'etrique d'une matrice $M$.

\begin{theorem}\label{theoreme1}\cite{albert-vatter2}

 La cl\^oture par sommes de toute classe  g\'eom\'etriquement en grille est h\'er\'editairement alg\'ebrique.
\end{theorem}

\medskip

Nous rappelons que, pour tout $n\in\NN,~n\geq 4$, nous avons quatre sortes de bicha\^{\i}nes critiques (voir d\'etails en page \pageref{prmut-critique}),
$(V_n,L_{n,1},L_{n,2})$, $(V_n,L_{n,2},L_{n,1})$, $(V_n,L_{n,1},(L_{n,2})^{-1})$ et $(V_n,(L_{n,2})^{-1},L_{n,1})$ o\`u $V_n:=\{0, \dots, n-1\}\times \{0,1\}$, $L_{n,1}$ et $L_{n,2}$ sont deux ordres lin\'eaires sur $V_n$ donn\'es par: $L_{n,1}:=(0,0)<(0,1)<\cdots<(i,0)<(i,1)\cdots <(n-1,0)<(n-1,1)$ et $L_{n,2}:=(n-1,0) <\cdots<(n-i,0) <\cdots <(0, 0)< (n-1,1)<\cdots  <(n-i,1)\cdots <(0, 1)$.

\vspace{2mm}

Les types d'isomorphie des bicha\^{\i}nes critiques sont les permutations exceptionnelles, elles sont d\'ecrites en terme de permutations de $1, \dots, 2m$ pour $m\geq 2$ dans \cite{A-A}:
\begin{enumerate}
\item[(i)] $\sigma_m^1:=2.4.6....2m.1.3.5....2m-1.$
\item[(ii)] $\sigma_m^2:=2m-1.2m-3....1.2m.2m-2....2.$
\item[(iii)] $\sigma_m^3:=m+1.1.m+2.2....2m.m.$
\item[(vi)] $\sigma_m^4:=m.2m.m-1.2m-1....1.m+1.$
\end{enumerate}

\vspace{2mm}

\noindent Remarquons que $\sigma_m^3=(\sigma_m^1)^{-1}$ et $\sigma_m^4=(\sigma_m^2)^{-1}$.

\vspace{3mm}

Ces quatre classes de permutations sont des classes g\'eom\'etriquement en grille. En effet, nous avons:

\begin{itemize}
\item $\{\sigma_m^1/~m\in\NN,~m\geq 4\}\subseteq Geom(M_1)$ avec $M_1=(1~~~1)$.
\item $\{\sigma_m^2/~m\in\NN,~m\geq 4\}\subseteq Geom(M_2)$ avec $M_2=(-1~~~-1)$.
\item $\{\sigma_m^3/~m\in\NN,~m\geq 4\}\subseteq Geom(M_3)$ avec $M_3=\left(\begin{array}{c}1\\
                                                                                1
                                                                                \end{array}\right)$.
\item $\{\sigma_m^4/~m\in\NN,~m\geq 4\}\subseteq Geom(M_3)$ avec $M_3=\left(\begin{array}{c}-1\\
                                                                                -1
                                                                                \end{array}\right)$.
\end{itemize}

\vspace{2mm}

Il s'ensuit, d'apr\`es le Th\'eor\`eme \ref{theoreme1}, que la cl\^oture par sommes de la classe des permutations exceptionnelles, et donc la cl\^oture par sommes de la classe des bicha\^{i}nes critiques,  est h\'er\'editairement alg\'ebrique. 
\clearemptydoublepage
\addcontentsline{toc}{chapter}{Notations}\label{notation}

%\addcontentsline{toc}{chapter}{Notations}\label{notation}
\markboth{\slshape{Notations}} {\slshape{Notations}}

\chapter*{Notations}

Dans ce document sont utilis\'ees les notations suivantes.\\
%On désignera par :

$\centerdot$ $\varnothing$ d\'esigne l'ensemble vide, $\subset$, $\subseteq$, $\cup$,  $\cap$ d\'esignent les relations ensemblistes: l'inclusion au sens strict, l'inclusion au sens large, la r\'eunion et l'intersection des ensembles respectivement.\\

$\centerdot$ $A\setminus B$ : compl\'ementaire de $B$ dans $A$.\\

$\centerdot$ $\vert E\vert$ : cardinalit\'e de l'ensemble $E$.\\

$\centerdot$ $\mathscr P(E)$ : ensemble des parties de $E$.\\

$\centerdot$ $E^n$ : ensemble des $n$-uples d'\'el\'ements de $E$.\\

$\centerdot$ $[E]^n$ : ensemble des parties \`a $n$ \'el\'ements de $E$.\\

$\centerdot$ $\Delta_E:=\{(x,x):x\in E\}$.\\

$\centerdot$ $n$ un entier, $[n]:=\{1,2,\dots,n\}$, $[0]:=\varnothing$.\\

$\centerdot$ $\mathfrak S_n$ : ensemble des permutations de $\{1,2,\dots,n\}$ et $\mathfrak S:=\underset{n\in\mathbb N}\cup\mathfrak S_n$.\\

$\centerdot$ Si $G$ est un graphe, $\overline{G}$ est son graphe complémentaire.\\

$\centerdot$ $K_n$ : graphe complet à $n$ sommets.\\

$\centerdot$ $I_n$ : ind\'ependant à $n$ sommets ou graphe vide à $n$ sommets.\\

$\centerdot$ $K_{\mathbb N}$ : graphe complet ayant pour ensemble de sommets l'ensemble $\mathbb N$ des entiers non n\'egatifs.\\

$\centerdot$ $I_{\mathbb N}$ : l'ind\'ependant ayant pour ensemble de sommets l'ensemble $\mathbb N$ des entiers non n\'egatifs.\\

$\centerdot$ $\Delta_{\mathbb N}$ : l'anticha\^{i}ne ayant pour ensemble de sommets l'ensemble $\mathbb N$ des entiers non n\'egatifs.\\

$\centerdot$ $\omega$ : chaîne des entiers naturels ($\mathbb N$ muni de l'ordre naturel).\\

%$\centerdot$ $\rho,\delta,\dots$ désignent des relations $n$-aire, $n\in\mathbb N^{\star}$.\\

$\centerdot$  Si $\rho$ est une relation $n$-aire, $\rho^c$ désigne sa relation complémentaire.\\

$\centerdot$ Si $\rho$ est une relation binaire sur $E$, $\rho^{-1}:=\{(x,y)\in E^2:(y,x)\in\rho\}$ est sa relation inverse ou duale not\'ee \'egalement $\rho^{\star}$.\\

$\centerdot$ Si $\mathcal R$ est une structure relationnelle alors

        \hspace{2cm} - $V(\mathcal R)$ (ou $dom(\mathcal R)$) d\'esigne son domaine ou sa base.

        \hspace{2cm} - $\tau(\mathcal R)$ d\'esigne son type d'isomorphie.

        \hspace{2cm} - $\mathcal R_{\restriction_A}$ d\'esigne sa restriction  \`a $A\subseteq dom(\mathcal R)$, elle est parfois not\'ee $\mathcal R_A$.

        \hspace{2cm} - $\varphi_{\mathcal R}$ désigne le profil de $\mathcal R$.\\

$\centerdot$ $\Omega_\mu$ : classe des structures relationnelles finies de signature $\mu$.\\

$\centerdot$ $\textit{T}_{\mu}$ : Collection des types d'isomorphie des structures de $\Omega_\mu$.\\

$\centerdot$  $\mathcal {A(R)}$: \^age de $\mathcal R$ est la collection des restrictions finies de $\mathcal R$ consid\'er\'ees à l'isomorphie pr\`es ou collection des types d'isomorphie des restrictions finies de $\mathcal R$.\\

$\centerdot$ Pour $\mathcal R$ et $\mathcal R'$ dans $\Omega_\mu,\;\mathcal R\leq \mathcal R'$ signifie $\mathcal R$ s'abrite dans $\mathcal R'$.\\

$\centerdot$ Si $\mathcal P$ est un ensemble pr\'eordonn\'e et $A\subseteq \mathcal P$ alors :

 \hspace{2cm} - $Forb(A):=\{x\in\mathcal P\;:\; y\nleq x, \forall y\in A\}$.

\hspace{2cm} - $\downarrow A:=\{x\in\mathcal P\; : x\leq y \text{ pour au moins un } y\in A\}$.

\hspace{2cm} - $\uparrow A:=\{x\in\mathcal P\; : x\geq y \text{ pour au moins un } y\in A\}$.\\

$\centerdot$ $\varphi_{\mathscr C}$ : profil de la classe $\mathscr C$.\\

$\centerdot$ $\mathcal R:=(E,\rho_1,\dots,\rho_k)$, o\`u $\rho_i$ est une relation binaire pour tout $1\leq i\leq k$, d\'esigne une structure binaire de type $k$.\\

$\centerdot$ $\mathcal R:=(E,\leq,\rho_1,\dots,\rho_k)$, o\`u $\leq$ est un ordre total et $\rho_i$ une relation binaire pour tout $1\leq i\leq k$, d\'esigne une structure binaire ordonn\'ee de type $k$.\\

$\centerdot$ $\Omega_k$ : classe des structures binaires finies de type $k$.\\

$\centerdot$ $\Theta_k$ : classe des structures binaires ordonn\'ees finies de type $k$, $\Theta_k\subseteq\Omega_{k+1}$.\\

$\centerdot$ $\Gamma_k:=\{\mathcal R\in\Theta_k: \mathcal R \text{ réflexive}\}$ : classe des structures binaires ordonn\'ees finies r\'eflexives de type $k$, $\Gamma_k\subseteq\Theta_k$.\\

$\centerdot$ $Ind(\mathscr C)$ : classes des membres  ind\'ecomposables de la classe $\mathscr C$.\\

$\centerdot$ $Ind(\mathcal R)$ : collection des sous-structures finies ind\'ecomposables de la structure binaire $\mathcal R$.\\

$\centerdot$ Pour $\mathscr I\subseteq Ind(\Omega_k)$, $\sum\mathscr I:=\{\mathcal R\in\Omega_k\; : \;Ind(\mathcal R)\subseteq\mathscr I\}$.\\

$\centerdot$ $\mathscr D_{(2)}^k:=\{\mathcal R\in\Gamma_k:\vert dom(\mathcal R)\vert=2\}$ : ensemble des structures binaires ordonn\'ees, r\'eflexives de type $k$ d\'efinies sur deux \'el\'ements, $\mathscr D_{(2)}^k\subseteq\Gamma_k$.\\

$\centerdot$ $\mathscr S_k$ : classe de structures binaires ordonn\'ees séparables de type $k$, $\mathscr S_k\subseteq\Theta_k$.\\

$\centerdot$ $\mathscr T_k^s:=\{\tau(\mathcal R):\mathcal R\in\mathscr S_k\}$ : classe des types d'isomorphie des structures de $\mathscr S_k$.\\

$\centerdot$ $\mathscr S^{re}_k$ (resp. $\mathscr S^{ir}_k$) : sous-classe form\'ee des structures r\'eflexives (resp. irr\'eflexives) de $\mathscr S_k$,
 $\mathscr S^{re}_k:=\mathscr S_k\cap\Gamma_k$.\\

$\centerdot$ $\mathscr T^{re}_k$ (resp. $\mathscr T^{ir}_k$) : sous-classe form\'ee des types d'isomorphie des structures de $\mathscr S^{re}_k$ (resp. $\mathscr S^{ir}_k)$; $\mathscr T^{re}_k\subseteq\mathscr T_k^s$.\\

$\centerdot$ Si $i$ est un entier, $\tau$  une structure de $\Theta_k$ et $\mathscr A$ une sous-classe de $\Theta_k$ alors:

       \qquad \qquad - $\mathscr A(\tau)$ d\'esigne l'ensemble des membres de $\mathscr A$ qui sont $\tau$-ind\'ecomposables.

       \qquad \qquad - $\mathscr A_{(i)}$ d\'esigne l'ensemble des membres de $\mathscr A$ qui ont $i$ \'el\'ements.

       \qquad \qquad - $\mathscr A_{(\geq i)}$ d\'esigne l'ensemble des membres de $\mathscr A$ qui ont au moins $i$ \'el\'ements.\\

$\centerdot$ $x\simeq_{F,\mathcal R}y$ : signifie $x$ et $y$ sont $F$-\'equivalents dans $\mathcal R$, c'est \`a dire les restrictions de $\mathcal R$ \`a $\{x\}\cup F$ et \`a $\{y\}\cup F$ sont isomorphes.\\

$\centerdot$ $x\simeq_{k,\mathcal R}y$ : signifie $x$ et $y$ sont $k$-\'equivalents dans $\mathcal R$, c'est \`a dire que $x\simeq_{F,\mathcal R}y$ pour toute partie $F$ \`a $k$ \'el\'ements qui ne contient pas $x$ et $y$.\\

$\centerdot$ $x\simeq_{\leq k,\mathcal R}y$ : signifie $x$ et $y$ sont $\leq k$-\'equivalents dans $\mathcal R$, c'est \`a dire que $x\simeq_{k',\mathcal R}y$ pour tout $k'\leq k$.\\

$\centerdot$ $x\simeq_{\mathcal R}y$ : signifie $x$ et $y$ sont \'equivalents dans $\mathcal R$, c'est \`a dire que $x\simeq_{k,\mathcal R}y$ pour tout entier $k$.\\

$\centerdot$  $\mathscr D_{\mu,k}$ : la classe des structures ordonn\'ees $\mathcal R$ de signature $\mu$ telle que $\simeq_{\leq k,\mathcal R}$ poss\`ede un nombre infini de classes.\\

$\centerdot$ $dim_{fort}(\mathcal R)$ : nombre de composantes fortement monomorphes de $\mathcal R$. Si ce nombre n'est pas fini, on \'ecrit $dim_{fort}(\mathcal R)=\infty$.\\

$\centerdot$ $dim_{mon}(\mathcal R)$ : nombre de composantes monomorphes de $\mathcal{P(R)}$, la d\'ecomposition monomorphe canonique de $\mathcal R$. Si ce nombre n'est pas fini, on \'ecrit $dim_{mon}(\mathcal R)=\infty$.\\

$\centerdot$ $dim_{mon}^{\infty}(\mathcal R)$ : nombre de composantes monomorphes infinies de $\mathcal{P(R)}$ lorsque $dim_{mon}(\mathcal R)<\infty$.\\

$\centerdot$ $dim_{Fraisse}(\mathcal R)$ : nombre d'intervalles Fra\"{\i}ss\'e monomorphes maximaux de $\mathcal R$. Si ce nombre n'est pas fini, on \'ecrit $dim_{Fraisse}(\mathcal R)=\infty$.\\

$\centerdot$ $dim_{int}(\mathcal R)$ : nombre d'intervalles monomorphes maximaux de $\mathcal R$. Si ce nombre n'est pas fini, on \'ecrit $dim_{int}(\mathcal R)=\infty$.\\

$\centerdot$ $\mathscr S_{\mu}$ : classe de structures relationnelles de signature $\mu$ (finie)  ne poss\`edant pas de d\'ecomposition monomorphe ayant un nombre fini de blocs.\\

$\centerdot$ $\mathscr D_{\mu}$ : classe de structures relationnelles  ordonn\'ees de signature $\mu$ (finie) ne poss\`edant pas de d\'ecomposition en intervalles ayant un nombre fini de blocs.\\

\clearemptydoublepage

\clearemptydoublepage

\addcontentsline{toc}{chapter}{index}
\printindex

\end{document}